# Statistical mechanics of mean-field disordered systems
*A Hamilton-Jacobi approach*

Tomas Dominguez and Jean-Christophe Mourrat

# Contents

# Overview

The central goal of this book is to present new techniques for studying the behaviour of mean-field systems with disordered interactions. We will focus in particular on certain problems of statistical inference, and on spin glasses. We will also consider simpler models of statistical mechanics for which the interactions are not disordered, as a useful training ground. We will mostly be interested in computing the asymptotic behaviour of a fundamental quantity called the free energy, in the limit of large system size. In the context of statistical inference, the free energy is essentially the mutual information between the observations and the signal that we wish to recover; knowing its limit behaviour allows us to identify how much of the signal can be recovered from the observations. The purpose of this work is to present a way to approach this class of problems using techniques from partial differential equations, and specifically, the theory of Hamilton-Jacobi equations. We strove to make this book self-contained, and in particular, no prior knowledge of Hamilton-Jacobi equations or other partial differential equations is assumed.

To get a sense of what spin glasses are, let us start by describing the simplest such model, namely the Sherrington-Kirkpatrick (SK) model [237]. This model is often motivated by considering the problem of splitting $N$ individuals into two groups. An assignment of individuals into two groups can be encoded by a vector $\sigma \in \Sigma_N := \{-1, +1\}^N$, with the understanding that $\sigma_i$ represents the group to which the individual indexed by $i$ is assigned. For each pair of individuals $(i, j)$, we are given a number $g_{ij}$ that encodes the quality of the interaction between individuals $i$ and $j$, with $g_{ij}$ being large if $i$ and $j$ get along very well, and $g_{ij}$ being very negative if they cannot stand each other. We would like to find an assignment into groups that maximizes the sum total of the interactions, that is, we want to maximize, over the set $\Sigma_N$, the "comfort function"

$$c_N(\sigma) := \sum_{i,j=1}^{N} g_{ij} \sigma_i \sigma_j. \tag{0.1}$$

Given the problem we are trying to encode, it would probably have felt more natural to write $\mathbf{1}_{\{\sigma_i = \sigma_j\}}$ in place of $\sigma_i \sigma_j$ above. However, the form $\sigma_i \sigma_j$ is more standard, partly for historical reasons having to do with the modelling of magnetic materials. Since $\sigma_i \sigma_j = 2\mathbf{1}_{\{\sigma_i = \sigma_j\}} - 1$, writing one or the other is inconsequential. It will be convenient not to impose any symmetry on the matrix $(g_{ij})_{1 \leqslant i, j \leqslant N}$, so the interaction between $i$ and $j$, with $i \neq j$, is actually encoded by $g_{ij} + g_{ji}$. It will also be convenient to include some diagonal terms $(g_{ii})_{1 \leqslant i \leqslant N}$; these terms only change the comfort function by an additive constant, and their contribution will in any case be of lower order.

In order to gain some insight into a *typical* instance of this optimization problem, we assume a particular structure on the interaction parameters $(g_{ij})_{1 \leqslant i, j \leqslant N}$. Namely, we assume that the $(g_{ij})_{1 \leqslant i, j \leqslant N}$ are independent standard Gaussian random variables. By "standard Gaussian", we always mean Gaussian random variables with zero mean and unit variance.

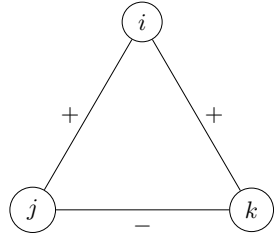

**Figure 0.1** A simple example of frustration.

One quickly realizes that the optimization of the comfort function (0.1) is not going to be easy because the coefficients $(g_{ij})_{1 \leqslant i, j \leqslant N}$ do not have a sign. A glimpse of the difficulty can already be perceived when we consider situations in which three individuals $i$, $j$ and $k$ are as depicted in Figure 0.1, that is, with $g_{ij} > 0$ and $g_{ik} > 0$, but with $g_{jk} < 0$. In view of the sign of $g_{ij}$, we would rather want $i$ and $j$ to belong to the same group, and similarly with $i$ and $k$. But the sign of $g_{ik}$ suggests to rather have $i$ and $k$ in different groups, and all of these "local preferences" cannot be reconciled at once. Physicists use the word *frustration* to describe situations of this nature, and, a *glass* is a system that is subject to such frustrations. More generally, the presence of these frustrations suggests that it is likely to be difficult to optimize the function (0.1). For instance, a naive method that would try to optimize each *spin* $\sigma_i$ one at a time in order to decrease the value of (0.1) is unlikely to reach the global optimum. The word "glass" points to the fact that this phenomenology is also present for regular window glass, as well as a number of other materials. Indeed, the making of a glass requires the very rapid cooling of an initially liquid material. A slow cooling would have allowed the material to find its preferred crystal structure, but the fast cooling has blocked the particles in a highly disordered configuration, and the complex geometric arrangement of the particles is believed to create jams that render the finding of the preferred crystalline state very difficult to achieve.

In order to gain insight into the problem of optimizing the comfort function (0.1), it is useful to first focus on the evaluation of the asymptotic behaviour of the maximum itself, and more generally of quantities of the form

$$\frac{1}{N}\mathbb{E}\log \sum_{\sigma \in \Sigma_N} \exp\left(\frac{\beta}{\sqrt{N}} \sum_{i,j=1}^{N} g_{ij}\sigma_i\sigma_j\right), \tag{0.2}$$

where $\beta \geq 0$ is a tunable parameter, and $\mathbb{E}$ denotes the expectation with respect to the randomness of the $(g_{ij})$. The quantity in (0.2) is called the *free energy*, and it is closely related to the random probability measure, called the *Gibbs measure*, that assigns to each configuration $\sigma \in \Sigma_N$ a probability proportional to

$$\exp\left(\frac{\beta}{\sqrt{N}} \sum_{i,j=1}^{N} g_{ij}\sigma_i\sigma_j\right). \tag{0.3}$$

To understand the presence of the factor of $N^{-1/2}$ in the exponential above, one needs to realize that the maximum of the comfort function (0.1) will typically be of order $N^{3/2}$, see Exercise 6.1. With the factor of $N^{-1/2}$ in the exponential in (0.2), we thus ensure that the terms that have the largest contribution to the sum in (0.2) are exponential in $N$. Since we are then summing over $2^N$ terms, this allows us to interpolate between a situation in which the entropy dominates, for small $\beta$, where the Gibbs measure resembles the uniform law over $\Sigma_N$, and a situation in which the energy dominates, for large $\beta$, where the Gibbs measure concentrates on the configurations that essentially realize the maximal value of the comfort function (0.1). In particular, a good approximation of $\beta N^{-3/2}$ times the maximum of the comfort function (0.1) is obtained by choosing $\beta$ sufficiently large, not depending on $N$, in (0.2). This is made precise in Exercise 6.3.

The problem of identifying the large-$N$ limit of the free energy (0.2) turns out to be surprisingly difficult. Starting in the late 1970's, Giorgio Parisi and collaborators [177, 178, 179, 217, 218, 219, 220, 221] proposed a way to solve this problem using sophisticated non-rigorous techniques, and the limit of (0.2) is now known as the *Parisi formula*. A rich phenomenology progressively emerged concerning the structure of the associated Gibbs measure, which turns out to organize itself along an ultrametric structure. Some key elements of this picture that were uncovered in the physics literature were then progressively put on a rigorous mathematical footing [132, 210, 211, 250, 253, 254]. Although many interesting questions remain, the mathematical understanding of the SK model, and of some generalizations thereof, is by now very substantial [41, 48, 61, 134, 210, 211, 213, 214, 215, 252].

There are however many seemingly innocent generalizations of the SK model that mostly remain mathematically mysterious. Motivated in part by considerations that relate to artificial neural networks (see e.g. [257]), we would like for instance to consider generalizations of the SK model in which the spins are organized over two (or more) layers. To encode this precisely, we can represent a configuration as a pair $\sigma = (\sigma_1, \sigma_2) \in \Sigma_N^2$, where $\sigma_1 = (\sigma_{1,i})_{1 \leq i \leq N}$ and $\sigma_2 = (\sigma_{2,i})_{1 \leq i \leq N}$ represent the spins in each of the two layers, and where we choose the layers to be of the same size for convenience of notation. The energy function can then be written as

$$H_N(\sigma) := \frac{1}{\sqrt{N}} \sum_{i,j=1}^{N} g_{ij} \sigma_{1,i} \sigma_{2,j}. \tag{0.4}$$

We will refer to this as the *bipartite model*. An illustration of the graph of interactions between the coordinates of $\sigma$ is in Figure 0.2. Perhaps surprisingly, this model is much less understood than the SK model. In particular, the limit free energy of this model has not yet been identified rigorously.

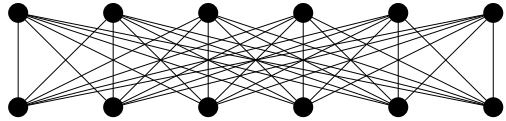

**Figure 0.2** Illustration of the bipartite model for $N = 6$. Elementary units are organized in two layers and only interact across layers.

Inspired by [42, 43, 64, 131, 200], we propose to approach this problem using a point of view based on partial differential equations [194, 196]. It turns out that one can identify the limit of the free energy (0.2) in the SK model as

$$\lim_{N \to +\infty} -\frac{1}{N} \mathbb{E} \log \sum_{\sigma \Sigma_N} \exp\left( \sqrt{\frac{2t}{N}} \sum_{i,j=1}^{N} g_{ij} \sigma_i \sigma_j - Nt \right) = f(t, 0), \tag{0.5}$$

where $f = f(t, \mathsf{q}) : \mathbb{R}_{\geq 0} \times \mathcal{Q}_2(\mathbb{R}_{\geq 0}) \to \mathbb{R}$ is the solution to the infinite-dimensional Hamilton-Jacobi equation

$$\partial_t f(t, \mathsf{q}) - \int_0^1 \partial_\mathsf{q} f(t, \mathsf{q}, u)^2 \, du = 0 \quad \text{on} \quad \mathbb{R}_{>0} \times \mathcal{Q}_2(\mathbb{R}_{\geq 0}), \tag{0.6}$$

with $\mathcal{Q}_2(\mathbb{R}_{\geq 0})$ being the space of square-integrable non-decreasing paths from $[0, 1)$ to $\mathbb{R}_{\geq 0}$, and $\partial_\mathsf{q} f(t, \mathsf{q}, \cdot)$ denoting the Gateaux derivative of $f$ at $(t, \mathsf{q})$, see (6.74). The initial condition $f(0, \cdot)$ to this equation is described by a functional transform of the Bernoulli measure $\delta_{-1} + \delta_{+1}$ encoding the "reference" law of one spin. Theorems 6.7 and 6.8 state this result more precisely. One of the most interesting aspects of this statement is that it suggests a very natural candidate for the limit free energy of the bipartite model and its generalizations, also phrased in the language of Hamilton-Jacobi equations, see Questions 6.9 and 6.11 for precise statements. An inequality between the limit free energy and this candidate limit is proved in [195, 197], but the converse bound remains an open problem.

The same approach has already led to more complete results concerning certain problems of statistical inference, as will be explained in detail in this book. For clarity of exposition, we will focus on the following relatively simple situation. We observe a noisy version of a rank-one matrix of the form $\bar{x}\bar{x}^*$, where $\bar{x} = (\bar{x}_1,\ldots,\bar{x}_N)$ is a vector of independent and identically distributed (i.i.d.) random variables, and the superscript $*$ denotes the transposition operator. Precisely, we assume that we observe the matrix

$$Y := \sqrt{\frac{2t}{N}}\bar{x}\bar{x}^* + W, \qquad (0.7)$$

where $t \geq 0$ is a free parameter that allows us to vary the signal-to-noise ratio, and $W = (W_{ij})_{1 \leq i,j \leq N}$ is a matrix of independent standard Gaussian random variables. We will see in Chapter 4 that this problem is closely related to that of community detection in random networks of Erdős-Rényi type with diverging average degree. We wish to answer the following question: given the observation of $Y$, can we recover meaningful information about the signal $\bar{x}\bar{x}^*$? One way to assess this is to monitor the minimal mean-square error

$$\mathsf{mmse}_N(t) := \frac{1}{N^2} \inf_g \mathbb{E}|\bar{x}\bar{x}^* - g(Y)|^2 = \frac{1}{N^2}\mathbb{E}|\bar{x}\bar{x}^* - \mathbb{E}[\bar{x}\bar{x}^* \mid Y]|^2 \qquad (0.8)$$

between the signal $\bar{x}\bar{x}^*$ and its noisy observation $Y$. Here, the infimum is taken over the set of measurable functions $g$, and for a matrix $a$, we write $|a| := \sqrt{\mathrm{tr}(aa^*)}$, with tr denoting the trace operator. Assuming that $\mathbb{E}\bar{x}_1 = 0$ for simplicity, we can compare this minimal mean-square error with the error one would make by simply using the null estimator, defining

$$\mathsf{var}_N(t) := \frac{1}{N^2}\mathbb{E}|\bar{x}\bar{x}^*|^2. \qquad (0.9)$$

It turns out that there exists a critical parameter $t_c \in (0,+\infty)$ such that the following holds: for $t < t_c$, essentially no information can be obtained about $\bar{x}\bar{x}^*$, in the sense that the difference between $\mathsf{mmse}_N(t)$ and $\mathsf{var}_N(t)$ becomes vanishingly small as $N$ becomes large; while for $t > t_c$, this difference remains bounded away from zero as $N$ tends to infinity. In fact, we will be able to characterize exactly the large-$N$ limit of $\mathsf{mmse}_N(t)$. This result can be obtained using a variety of methods [32, 33, 34, 96, 109, 161, 162]. We will use a Hamilton-Jacobi approach developed in [69, 71, 72, 75, 193, 194], and identify the large-$N$ limit of the minimal mean-square error (0.8) to be

$$\lim_{N\to+\infty} \mathsf{mmse}_N(t) = \left(\mathbb{E}|\bar{x}_1|^2\right)^2 - \partial_t f(t,0), \qquad (0.10)$$

where $f : \mathbb{R}_{\geq 0} \times \mathbb{R}_{\geq 0} \to \mathbb{R}$ is the solution to the Hamilton-Jacobi equation

$$\partial_t f(t,h) - \left(\partial_h f(t,h)\right)^2 = 0 \quad \text{on} \quad \mathbb{R}_{>0} \times \mathbb{R}_{>0}, \qquad (0.11)$$

subject to some explicit initial condition. This Hamilton-Jacobi approach has been extended to a very large class of other models of statistical inference in [75] (alternative approaches have not yet reached this level of generality).

Before going into this, we will train our skills on the analysis of the much simpler Curie-Weiss model and its generalizations. The Curie-Weiss model is the non-disordered version of the SK model: its energy function is obtained by setting all the interaction parameters $g_{ij}$ to be equal to 1 in (0.1). The identification of the limit free energy of this model is a straightforward consequence of Cramér's theorem on large deviations for sums of i.i.d. random variables. This shortcut through large deviations does not seem to be available for problems in statistical inference or spin glasses. Our first goal will be to recover the limit free energy of the Curie-Weiss model and its variants using the Hamilton-Jacobi approach. This will give us the opportunity to develop this approach and its associated toolbox in the simplest possible context. This toolbox will then be ready for us to use, with no further addition necessary, when we turn to the problems of statistical inference.

**Organization of the book**

This book is organized as follows. We start in Chapter 1 with an introduction to the basics of statistical mechanics, where we motivate the notion of a Gibbs measure on physical grounds, and introduce the Curie-Weiss model as well as its generalizations. In Chapter 2, we develop the fundamentals of convex analysis and large deviation principles, and use these to compute the limit of the free energy in the Curie-Weiss model and its generalizations. This large-deviation approach is not applicable to the disordered models we want to consider next, so we need to develop an alternative approach.

In Chapter 3, we define the notion of viscosity solution to a Hamilton-Jacobi equation, and use it to recover the limit free energy of the Curie-Weiss model. We discover technical challenges to applying the same method to generalized versions of the Curie-Weiss model, and develop a new "selection principle" based on convexity to overcome these. We then turn to statistical inference in Chapter 4, focusing on the problem of recovering a large symmetric rank-one matrix from a noisy observation. We discover that the tools developed in the previous chapter apply to this setting as well, and allow us to give a closed-form description of its phase transitions. Chapter 5 is preparatory work for a discussion of the more challenging case of spin glasses. The first half of this chapter is a self-contained introduction to Poisson point processes, including limit theorems on extreme values of independent and identically distributed random variables, which we believe to be of wide interest. We finally turn to the setting of spin glasses in Chapter 6. For the Sherrington-Kirkpatrick model, we show how to relate the Parisi formula with the Hamilton-Jacobi approach. We conclude with a more informal discussion on the status of current research for more challenging models. Appendix A is a self-contained presentation of many of the basic results in real analysis and probability theory that are used throughout the book. Solutions to the exercises are also provided.

**Related works**

Our main goal with this book is to provide the reader with an inviting presentation of some results and open problems in mean-field disordered systems which we find fascinating, following the common thread of the Hamilton-Jacobi approach. We wish to stress however that we do not claim this to be the one-and-only definitive way to approach these problems.

Trying to survey the range of works and topics related to those discussed in this book is a daunting task. We will attempt to fulfill it to the best of our abilities, hoping that the references listed here will at least provide the reader with good entry points to explore further according to their interests. Apologies to all those whose work is not included here.

Chapters 1 and 2 cover very classical material on statistical mechanics, convex analysis, and large deviations. References on these topics include [117], [53, 108, 138, 187, 230], and [92, 93, 258] respectively.

Concerning the theory of viscosity solutions to Hamilton-Jacobi equations discussed in Chapter 3, the most classical reference on the topic is probably [84]; we also mention [38, 39, 115] for more accessible presentations that focus on first-order equations. The notion of viscosity solution was introduced in [83, 85, 114], and the Hopf-Lax and Hopf formulas were introduced in [140, 159] and proved to be viscosity solutions in [165] and [37, 166] respectively. The infinite-dimensional Hamilton-Jacobi equations discussed in Chapter 6 are developed in [73, 74, 195].

The problems of statistical inference discussed in Chapter 4 are also explored in the surveys [185, 268], along with several aspects not covered here. We also mention the monograph [175] presenting the interplay between statistical mechanics, information theory, and combinatorial optimization, with extensive coverage of a class of methods called message-passing or belief-propagation algorithms that have shown their usefulness in a wide variety of circumstances.

The calculation of the limit free energy (or mutual information) for the inference problem discussed in Chapter 4 and its generalizations has been approached using a large variety of techniques. Together with an interpolation in the spirit of [132, 133, 134], these may involve algorithmic approaches in [32, 95, 96, 148, 162], a cavity method in the spirit of [13] in [161, 163, 172, 182, 228], concentration of measure in [109], or adaptive interpolation in [33, 34, 35, 168, 169, 227]. The Hamilton-Jacobi approach presented in Chapter 4 has been developed in [69, 71, 72, 75, 193, 194]. Examples of models that are covered by this approach but currently not by other techniques are discussed in Section 7 of [168].

In Chapter 4 on statistical inference, we are mostly concerned with determining whether or not one can recover meaningful information about the signal from a noisy observation, and if so how much. In practice, it is also fundamental to know whether or not one can recover this information in a reasonable amount of time. We barely scratch the surface on this point in Section 4.4. The example in Proposition 4.21 shows a situation in which for some regime of parameters, it is theoretically possible to recover meaningful information about the signal, but the standard (and rapid) PCA method completely fails. In fact, one expects that in this regime, there is no polynomial-time algorithm that recovers meaningful information about the signal. This difference between the moment when one can theoretically infer non-trivial information about the signal, and the one when one can do so efficiently, is often called the statistical-to-computational gap. It is believed to exist in a wide variety of situations, and we refer to [31, 123, 125, 268] for surveys on this.

In Section 4.5, we discuss the relationship between the inference problem (0.7) and the problem of detecting communities from the observation of a random graph where the probability for an edge to be drawn between two nodes depends on the community of each node. A survey on this problem of community detection is [1]; we also refer to [91] where many predictions were put forward at the physics level of rigour. In the setting in which the average degree of the graph of connections remains bounded, the identification of the regime of parameters for which one can reconstruct non-trivial information about the community structure from the observation of the random graph was achieved in [171, 188, 191] for models with two communities. To the best of our knowledge, the identification of a closed-form description of this regime of parameters remains open in settings with more than two communities, and the difficulties seem to be of a similar nature as those encountered in the analysis of the bipartite spin-glass model (0.4). We refer in particular to [5] for positive results in this direction, and again to [91] for several predictions. A related problem concerns the determination of the asymptotic mutual information between the community structure and the graph of connections. This problem was successfully resolved in "convex" cases [4, 80] and even in some "non-convex" cases [3, 129, 149, 189, 192, 267], but a fully general solution has still not been identified [129]. (What counts as a "convex" or a "non-convex" model should hopefully become clear upon reading Chapter 6.) This problem seems very similar to those presented in Questions 6.9 and 6.11 in the context of spin glasses, as discussed further in [105, 153].

Books with a focus on spin glasses include [55, 60, 63, 67, 81, 88, 179, 203, 205, 211, 239, 253, 254]. The construction of Poisson-Dirichlet cascades in Sections 5.5 and 5.6 essentially follows [211]. We also refer to [211] for a complete proof of the Parisi formula (6.8).

The idea that the limit free energy of the SK model should be described by the Parisi formula emerged in the series of works [217, 218, 219] based on the non-rigorous replica method. The physical understanding of the model progressed further in many contributions including [178, 177, 179, 220]. A mathematical proof of the Parisi formula was then obtained in [132, 250]. A more robust proof was later found in [210, 212] using the idea of ultrametricity discussed in Section 5.7, and also inspired by [11, 13, 21, 57, 128, 207, 208, 209]. We refer to [211] for a more thorough presentation of historical developments on this topic. The Parisi formula for spherical models was obtained in [77, 249], and takes a form that was predicted in [86].

Early contributions in the physics literature on the bipartite spin-glass model (0.4) include [119, 120, 152, 156]. On the mathematical side, the limit free energy of spin glasses with multiple types has been identified in a number of cases. For Hamiltonians such that (6.163) holds, the limit free energy has been identified in [41, 206, 213, 214, 215] when $\xi$ is convex over $\mathbb{R}^{D \times D}$; see also [47, 48, 155, 216] in the spherical case. The more general case when $\xi$ is convex over the space of positive semi-definite matrices is obtained in [70, 73, 195, 197] using some of the ideas presented here. For fully general models in the form of (6.163), Theorem 6.12 from [70] imposes strong constraints on what the limit free energy can be; it is at present unclear whether this statement is a complete characterization of the limit free energy or not (even though there are choices of $(t, \mathsf{q})$ such that (6.171) is satisfied for multiple pairs $(\mathsf{q}', \mathsf{p})$, it may still be the case that, for instance, there is only one continuous function $f$ that satisfies the properties listed in Theorem 6.12). Spin glasses with multiple types were already present in the first proof of the Parisi formula from [250, 251, 253, 254], in the form of two copies of the original system coupled together through a constraint on their overlap.

That there exist connections between limit free energies of statistical-mechanics systems and solutions to Hamilton-Jacobi equations dates back at least to [64, 200]; see also [49] for a survey of related contemporary research topics. In the context of spin glasses, heuristic connections between limit free energies and Hamilton-Jacobi equations were first pointed out in [8, 42, 43, 131], under a replica-symmetric or one-step replica symmetry breaking assumption. The possibility to rephrase the Parisi formula in terms of a Hamilton-Jacobi equation as in Theorems 6.7 and 6.8 is from [196, 198]. Theorems 6.10 and 6.12 are from [195] (generalized in [197]) and [70] respectively. A high-temperature version of Theorem 6.12 is in [97].

For spherical models with multiple types, the limit free energy is identified in some non-convex cases in [243, 244, 245, 246] using an approach inspired by the early work [255], and in [29] for the bipartite spherical SK model. We also mention [51, 145, 154, 174] for a geometric analysis of the energy landscape of spin glasses with multiple types; see also [23, 24, 118, 241] in the single-type case.

By taking a low-temperature limit, the determination of the limit free energy of spin glasses allows one to infer the asymptotic value of the maximum of the Hamiltonian; see Exercise 6.3 and [26, 78, 102]. In the context of the Sherrington-Kirkpatrick model, with $H_N$ as in (6.3), this would amount to determining the asymptotic behaviour of the maximum of $H_N(\sigma)/N$ over $\sigma \in \Sigma_N$; we denote the limit by OPT. One may ask whether there exists an efficient algorithm for actually finding a configuration $\sigma \in \Sigma_N$ such that $H_N(\sigma)/N$ is close to OPT. The answer turns out to depend on the specifics of the model, and relates to the overlap gap property [124]. A specific value ALG described by an analogue of the Parisi formula was identified such that there exists an efficient algorithm that can identify a configuration $\sigma$ with $H_N(\sigma)/N$ approximately equal to ALG with high probability [110, 184, 236, 242]. Moreover, it was argued in [142] that no algorithm within a broad class can exist which improves upon this value. Remarkably, these results on algorithmic thresholds have been extended to non-convex spherical models such as the bipartite model in (0.4), despite the fact that we do not know of a characterization of OPT or of the free energy in this case [143, 144]. The reference [27] surveys a number of results related to optimization algorithms for spin glasses.

The study of mean-field spin glasses has inspired developments in many other contexts. Works that explore connections between spin glasses and neural networks include [9, 40, 44, 45, 46, 89, 122, 136, 170, 257]; see [7, 18, 90, 113, 121, 137, 226, 256] for overviews. For a suitable choice of the reference measure denoted by $P_N$ in Chapter 6, the law of the spins of a given type in the bipartite model (0.4) is the same as the law of the spins in the Hopfield model [40]. Some versions of the perceptron model can be obtained similarly. The Hopfield and perceptron models were introduced in [17, 141, 167, 173, 181] as toy models for memory storage and retrieval or classification tasks. Early works on the statistical mechanics of these models include [19, 20, 126, 127, 222, 223, 224]; recent rigorous works include [56, 101, 254, 265]. Among many other topics related to spin glasses, we mention random constraint satisfaction problems [98, 100, 158, 175, 180, 183], the random assignment and travelling salesman problems [14, 15, 176], error correcting codes in information theory [229], and combinatorial problems such as graph colouring [79, 99, 199]. A recent book on "spin glass theory and far beyond" is [67].

## Acknowledgements

This book grew out of lectures given by JCM at the 2021 CRM-PIMS summer school in Montreal, at ETH Zurich in 2022, and to graduate students at New York University and ENS Lyon. The two of us (TD and JCM) met on the occasion of the summer school in Montreal and have been working together since then. We are extremely grateful to the organizers of these events for their support, including financial support to TD from ETH Zurich to take part in the writing of this book. We also very warmly thank the participants of these two events and of the graduate courses for their interest and very useful feedback.

**Course outline**

We describe here how the series of lectures in Zurich were organized. We do not believe that this is the optimal organization, and some of the results discussed in the book were not known at the time (in particular Theorem 6.12), but we hope that it can still be of some interest. This series of lectures comprised 12 sessions of $2 \times 45$ min. The first session was an overview of motivations and a description of the plan for the course. The second session started with a refresher on large deviations concluding with Theorem 2.16, followed by the contents of Section 1.1. The third session covered the rest of Chapter 1 and finished with Theorem 2.19. The fourth session covered Sections 3.1 and 3.2. The fifth session contained the proof of the comparison principle in Theorem 3.5 in the simpler case where the space domain $\mathbb{R}^d$ is replaced by the torus, the statements of the variational formulas in Theorems 3.8 and 3.13 without proofs, and concluded with Corollary 2.20 recovered using the Hamilton-Jacobi approach. The plan for the sixth session was to cover Section 3.6, with a direct proof of the convex selection principle in Theorem 3.21 that bypasses Lemma 3.23, and also to present some material from Subsection 2.1.3 on subdifferentials. This actually spilled over to the seventh session. The rest of this session was spent discussing the setup of the statistical inference problem in Chapter 4, and motivating it informally in relation with the problem of community detection discussed in Section 4.5. The eighth session covered Section 4.1 and Section 4.3 up to (4.83), mostly skipping Section 4.2 except for the simple Gaussian integration by parts in (4.29) and some simple generalization of it. The ninth session was meant to cover the rest of Section 4.3, taking the concentration of the free energy for granted (any estimate stating that (4.101) tends to zero with $N$ will do, and one can show an upper bound of the order of $N^{-1/3}$ on this quantity using simple arguments based on the Efron-Stein and Gaussian Poincaré inequalities in Exercises 4.7 and 4.8). This actually spilled over a bit to the tenth session. Also, the proof of Theorem 4.9 was only obtained under the additional assumption that $\mathbb{E}\bar{x}_1 = 0$, since the convex selection principle was only shown in the form of Theorem 3.21, but the more refined Lemma 3.22 is needed to conclude in general here. The major part of the tenth session was spent covering Section 6.2 and its analogue for the bipartite model, and discussing the new difficulties that show up. The eleventh session started with a discussion of the random energy model from Section 6.3, and then presented the main results of Section 5, in particular Proposition 5.13, as well as elementary properties of the Poisson-Dirichlet process from Section 5.5. In the last session, the Poisson-Dirichlet cascades from Section 5.6 were defined, the contents of Section 6.4 were covered without proving everything, and with some hand-waving for the remainder of Chapter 6. A recording of these lectures can be found at `https://tinyurl.com/HJ-ETHZ`.

Tomas Dominguez                      Jean-Christophe Mourrat
University of Toronto                ENS Lyon and CNRS

# Chapter 1
# Introduction to statistical mechanics

In this chapter, we introduce the basic objects from statistical mechanics that we will explore throughout the rest of the book. In Section 1.1, we introduce and motivate the notion of a Gibbs measure from a physical point of view, as well as the notion of free energy. The rest of the book will be devoted to computing this fundamental quantity for models of increasing complexity. In Section 1.2, we consider the classical Ising model, and give a brief historical overview of its developments. Finally, in Section 1.3, we trivialize the geometry of the Ising model to obtain the Curie-Weiss model. The Curie-Weiss model is simple enough that its free energy can be computed exactly using the classical large deviation principles discussed in Chapter 2, but it is rich enough to capture many of the challenges present in more sophisticated models. For this reason, we will then use it as a test-bed for the development of the Hamilton-Jacobi approach.

## 1.1  Gibbs measures

We start by motivating the notion of a Gibbs measure on physical grounds. In this section we will avoid delving into too many technicalities, so we will allow ourselves to not be as rigorous as in the rest of the book. Consider a system of $N$ units, which can each be in one of $K$ states $\{1, 2, \ldots, K\}$. We can think of these units as being particles. We think of $K$ as being fixed, while $N$ is very large and will be sent to infinity. Each state $k \in \{1, \ldots, K\}$ has an associated energy $e_k > 0$, and we assume that the system is "isolated" in the sense that its total energy is fixed, say at $N\bar{e}$ for some constant $\bar{e}$ in the convex hull of $\{e_1, \ldots, e_K\}$. The state of the system can be represented by a vector $x = (x_1, \ldots, x_N) \in \{1, \ldots, K\}^N$, and we denote by

$$N_k(x) := \sum_{n=1}^{N} \mathbf{1}_{\{x_n = k\}}. \tag{1.1}$$

the number of units in state $k$. The constraint on the total energy can be encoded as

$$\sum_{k=1}^{K} N_k(x) e_k = N\bar{e}. \tag{1.2}$$

We assume that the system is uniformly distributed among all possible configurations $x$ that satisfy the constraint (1.2). This assumption can be related to Hamiltonian dynamics and the Boltzmann hypothesis, but we simply take it for granted here.

A natural question is to determine the probability of finding a particular unit in state $k$,

$$\mathbb{P}\{x_i = k\} = \frac{1}{N} \sum_{n=1}^{N} \mathbb{P}\{x_n = k\} = \frac{1}{N} \mathbb{E} N_k(x), \tag{1.3}$$

when the number $N$ of units in the system is very large. We will argue, without full rigour nor attention to taking integer parts at the right places, that $N_k(x)/N$ converges in probability to a constant as $N$ tends to infinity, and we will identify this constant. In order to do so, we fix a probability vector

$$p = (p_1, \ldots, p_K) \in \mathbb{R}_{\geq 0}^{K} \quad \text{with} \quad \sum_{k=1}^{K} p_k = 1, \tag{1.4}$$

and we determine the asymptotic behaviour of the number of configurations $x$ such that $(N_1(x), \ldots, N_K(x)) = (Np_1, \ldots, Np_K)$. For a fixed $N$, this number is given by the multinomial coefficient

$$\binom{N}{Np_1}\binom{N - Np_1}{Np_2} \cdots \binom{N - \sum_{k=1}^{K-1} Np_k}{Np_K} = \frac{N!}{(Np_1)! \cdots (Np_K)!}. \tag{1.5}$$

Using that $N! = \exp(N \log N - N + \mathcal{O}(\log N))$, we can rewrite this as

$$\exp\left(N \log N - N + \sum_{k=1}^{K} \left[(Np_k) \log(Np_k) - Np_k\right] + \mathcal{O}(\log N)\right)$$

$$= \exp\left(-N \sum_{k=1}^{K} p_k \log p_k + \mathcal{O}(\log N)\right). \tag{1.6}$$

This motivates the introduction of the quantity

$$S(p) := -\sum_{k=1}^{K} p_k \log(p_k), \tag{1.7}$$

which is called the *entropy* of $p$. The calculation leading to (1.6) shows that, if $p$ and $p'$ are probability vectors with $S(p) > S(p')$, then there are exponentially more configurations with $(N_1(x), \ldots, N_K(x)) = (Np_1, \ldots, Np_K)$ than there are configurations with $(N_1(x), \ldots, N_K(x)) = (Np'_1, \ldots, Np'_K)$. In other words, if $S(p) > S(p')$, then the configurations consistent with $p$ vastly outnumber those that are consistent with $p'$, and so if we pick one configuration uniformly at random, it is overwhelmingly more likely that we observe $p$ rather than $p'$. To be more precise, with overwhelming probability, we will observe $(N_1(x), \ldots, N_K(x)) \simeq (Np_1^*, \ldots, Np_K^*)$, for the vector $p^*$ that maximizes the entropy (1.7) within the set

$$\Lambda := \left\{ p \in [0,1]^K \mid \sum_{k=1}^K p_k = 1 \text{ and } \sum_{k=1}^K p_k e_k = \overline{e} \right\}. \tag{1.8}$$

Using the strict concavity of the entropy, one can show that this maximizer $p^*$ exists and is unique. Our initial question regarding the probability of finding a particular unit in a given state has therefore been answered: asymptotically as $N$ tends to infinity, the proportion of particles in state $k \in \{1, \ldots, K\}$ is about $p_k^*$. We can in fact go a bit further and make $p^*$ more explicit. By the Lagrange multiplier theorem, there exist constants $\alpha, \beta \in \mathbb{R}$ such that for every $k \in \{1, \ldots, K\}$,

$$\partial_{p_k} S(p^*) = \alpha + \beta e_k. \tag{1.9}$$

Rearranging gives a constant $Z$ such that for every $k \in \{1, \ldots, K\}$,

$$p_k^* = \frac{\exp(-\beta e_k)}{Z}, \tag{1.10}$$

and since $\sum_{k=1}^K p_k^* = 1$, the constant $Z$ can be rewritten as

$$Z = \sum_{k=1}^K \exp(-\beta e_k). \tag{1.11}$$

The probability measure on $\{1, \ldots, K\}$ with probability vector $p^*$ given by (1.10) is called the *Gibbs measure* at inverse temperature $\beta$. When one builds the theory of thermodynamics from first principles, this Lagrange multiplier $\beta$ is *defined* to be the inverse temperature of the physical system (up to multiplication by the Boltzmann constant which we set equal to 1 here). The parameter $\beta$ is fixed in such a way that the average energy of the system with respect to the Gibbs measure is $\overline{e}$,

$$\frac{\sum_{k=1}^K e_k \exp(-\beta e_k)}{\sum_{k=1}^K \exp(-\beta e_k)} = \sum_{k=1}^K p_k^* e_k = \overline{e}. \tag{1.12}$$

In most physical systems, there is room for very high energy levels, so we typically have $\beta > 0$; however, this is not guaranteed in our current context, and will depend on how $\overline{e}$ compares with the arithmetic average of $(e_1, \ldots, e_K)$. A direct calculation allows us to relate the entropy, the average energy and the inverse temperature of the system,

$$S(p^*) = \sum_{k=1}^{K} p_k^* \log\big(Z\exp(\beta e_k)\big) = \log(Z) + \beta \bar{e}. \tag{1.13}$$

The quantity $\log(Z)$ should be called the *free entropy*; however, throughout the book we will call any quantity of this form the *free energy* of the system. In the language of physics, the free energy should be minus the free entropy divided by $\beta$, so that it is indeed expressed in the same units as $\bar{e}$. But since the interpretation of the inverse temperature will become less transparent as we proceed, and since the term free energy is much more commonly used, we will stick to the less proper terminology. The quantity $Z$ itself is often called the *partition function*. We mention in passing that the identity (1.13) is a shadow of a convex-duality relationship between the entropy and the free energy of a system, if we think of them as functions of $\bar{e}$ and $-\beta$ respectively; in particular,

$$-S(p^*) = \sup_{\beta' \in \mathbb{R}} \left(-\beta'\bar{e} - \log \sum_{k=1}^{K} \exp(-\beta' e_k)\right). \tag{1.14}$$

We also refer to Corollaries 4.14 and 4.15 in [59] for a more general view on this.

We have therefore been able to show that the probability of finding a microscopic unit of the system in a particular state is given by the Gibbs measure at inverse temperature $\beta$, and that the inverse temperature $\beta$ can be measured from macroscopic quantities of the system such as the free energy, the entropy and the average energy. Although we have considered a setting in which the microscopic and macroscopic parts of the system are made of identical units, a refined version of the argument allows to lift this restriction. To frame our findings in concrete terms, suppose we have a piece of material at equilibrium in a room at inverse temperature $\beta$, and assume the accessible energy levels of the material are $(e_1, \ldots, e_K)$. We should then expect to find the piece of material in state $k$ with probability proportional to $\exp(-\beta e_k)$. For further motivation on the physical terminology and the notion of a Gibbs measure, we refer the interested reader to Chapter 1 in [117].

**Exercise 1.1.** Let $(a_k)_{k \geq 1}$ be a sequence of real numbers.

(i) Show that, for each $K \geq 1$, we have

$$\lim_{N \to +\infty} \left| \frac{1}{N} \log\left(\sum_{k=1}^{K} \exp(Na_k)\right) - \max_{1 \leq k \leq K} a_k \right| = 0. \tag{1.15}$$

(ii) How rapidly can $K$ grow with $N$ for (1.15) to still remain valid?

## 1.2    The Ising model

The Ising model is undoubtedly the most famous model in the field of statistical mechanics. It was introduced by Willhelm Lenz and his student Ernst Ising in 1920 to gain a theoretical understanding of the phase transition from ferromagnetic to paramagnetic behaviour first observed in 1895 by Pierre Curie. In the lab, Curie found that, when subjected to a magnetic field, certain materials have a "memory" of the magnetic field they have been exposed to, while others do not. Materials of the former kind are said to exhibit *ferromagnetic* behaviour, while materials of the latter kind are said to exhibit *paramagnetic* behaviour. Perhaps surprisingly, Curie observed that materials can be made to transition from ferromagnetic to paramagnetic behaviour by increasing their temperature [87]. Shortly after these empirical observations, Auguste Piccard and Pierre Weiss gave a partly phenomenological explanation of what they called the *magnetocaloric phenomenon*, and proposed some formulas for the magnetic susceptibility of materials near the critical temperature at which the phase transition occurs [263]. It was following this line of investigation that Wilhelm Lenz and Ernst Ising introduced the Ising model in 1920.

In the Ising model, atoms are arranged on a finite box in the integer lattice $\mathbb{Z}^d$. Each atom carries a magnetic moment called a *spin* which can either be in the $-1$ or in the $+1$ orientation. The atoms interact with their nearest neighbours in such a way as to favour the alignment of spins. To formalize this model, let

$$B_N := \{-N, \ldots, N\}^d \tag{1.16}$$

denote the box of side-length $2N$ in the integer lattice. Given two lattice points $i, j \in B_N$, we write $i \sim j$ to mean that $i$ and $j$ are nearest-neighbours in the sense that they differ by one unit in at most one coordinate. A configuration of spins may be encoded by a vector $\sigma \in \{\pm 1\}^{B_N} := \{-1, +1\}^{B_N}$, and the interactions between the spins can be described by the energy function, or *Hamiltonian*, of the system,

$$H_N(\sigma) := \sum_{\substack{i,j \in B_N \\ i \sim j}} \sigma_i \sigma_j. \tag{1.17}$$

Physicists would typically add a minus sign to the Hamiltonian, but here we omit this minus sign for convenience of notation. With our sign convention, the system has a preference for larger values of the energy function as opposed to smaller ones. The discussion in Section 1.1 implies that, at inverse temperature $\beta > 0$, the probability of finding the system in the configuration $\sigma$ is the Gibbs weight

$$G_N(\beta, \sigma) := \frac{\exp(\beta H_N(\sigma))}{Z_N(\beta)} \tag{1.18}$$

with partition function

$$Z_N(\beta) := \sum_{\sigma \in \{\pm 1\}^{B_N}} \exp(\beta H_N(\sigma)). \tag{1.19}$$

In view of the Hamiltonian in (1.17), this means that the system will have a preference for configurations in which neighbouring spins are equal. In the presence of an external magnetic field of intensity $h \in \mathbb{R}$, the Hamiltonian at inverse temperature $\beta > 0$ is enriched with the term $h \sum_{i \in B_N} \sigma_i$ to obtain

$$H_N(\beta, h, \sigma) := \beta \sum_{\substack{i,j \in B_N \\ i \sim j}} \sigma_i \sigma_j + h \sum_{i \in B_N} \sigma_i \tag{1.20}$$

whose associated Gibbs measure is

$$G_N(\beta, h, \sigma) := \frac{\exp H_N(\beta, h, \sigma)}{Z_N(\beta, h)} \tag{1.21}$$

for the enriched partition function

$$Z_N(\beta, h) := \sum_{\sigma \in \{\pm 1\}^{B_N}} \exp H_N(\beta, h, \sigma). \tag{1.22}$$

From a physical perspective, it would be more natural that the inverse temperature $\beta$ multiplies both terms in (1.20), in other words to replace $h$ by $\beta h$ in (1.20), but the parametrization chosen in (1.20) is much more convenient to work with; and this is why we include it into the definition of $H_N$ as opposed to encoding the dependency in $\beta$ as in (1.18). To study the empirical observations of Curie from a theoretical point of view, we need to understand the mean magnetization

$$m_N(\beta, h) := \sum_{\sigma \in \{\pm 1\}^{B_N}} \left( \frac{1}{|B_N|} \sum_{i \in B_N} \sigma_i \right) G_N(\beta, h, \sigma), \tag{1.23}$$

where $|B_N|$ denotes the cardinality of the box $B_N$. One possible approach to understand the behaviour of the mean magnetization (1.23) when the number of particles $N$ tends to infinity is to start by studying the free energy

$$F_N(\beta, h) := \frac{1}{|B_N|} \log Z_N(\beta, h). \tag{1.24}$$

Indeed, a direct computation shows that $m_N(\beta, h) = \partial_h F_N(\beta, h)$. Moreover, using a box-decomposition argument, it is possible to show that the sequence $(F_N)_{N \geqslant 1}$ admits a limit $f : \mathbb{R}_{\geqslant 0} \times \mathbb{R} \to \mathbb{R}$, and that if $(\beta, h)$ is a point of differentiability of $f$, then

$$m(\beta, h) := \lim_{N \to +\infty} m_N(\beta, h) = \partial_h f(\beta, h). \tag{1.25}$$

The limit free energy $f$ therefore encodes much of the physical information about the system. It also has the advantage of being more "robust" to perturbations of the model than the mean magnetization $m_N(\beta, h)$ — for instance, it does not depend on boundary conditions one might introduce to encode some interaction between the spins on the inner boundary of the box $B_N$ and the outside environment. In contrast to the empirical evidence found by Curie, Ising [146] was able to show that in the one-dimensional setting, $d = 1$, the limit magnetization curve $h \mapsto m(\beta, h)$ remains smooth for every fixed $\beta > 0$. He also suggested that the same would happen in higher dimensions, $d \geq 2$. However, in 1936, Rudolf Peierls was able to prove that in any dimension $d \geq 2$, the function $h \mapsto m(\beta, h)$ does in fact have a discontinuity at $h = 0$ when $\beta$ is sufficiently large [225]. This is depicted in Figure 1.1 and is in agreement with Curie's empirical observations. Indeed, it implies that at low temperature and for $h > 0$, the difference between the fraction of spins in the $+1$ and in the $-1$ directions remains bounded away from zero, no matter how small $h$ is. In other words, the system retains some global ordering even in the limit of $h > 0$ going to zero. In contrast, in the high temperature regime, the random fluctuations dominate and no ordering can be preserved as $h > 0$ is sent to zero.

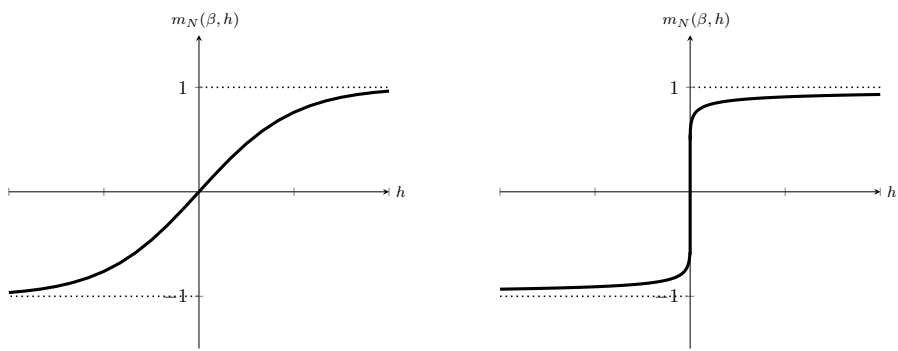

**Figure 1.1** Dependence of the magnetization (1.23) on the magnetic intensity $h$ in the two-dimensional Ising model for $N$ large and $\beta > 0$ small (on the left) and large (on the right).

The one-dimensional and two-dimensional Ising models are explored in Exercises 1.2 and 1.3. In Exercise 1.2, we compute the limit of the free energy (1.24) in the one-dimensional Ising model explicitly. In higher dimensions, $d \geq 2$, it is by no means guaranteed that a self-contained description of the limit free energy should exist. Surprisingly, in the 1940's, Lars Onsager and others managed to give an explicit description of the limit free energy of the two-dimensional Ising model [150, 151, 204, 247, 266]. This allowed them to determine several critical exponents which were completely unexpected. By opening a first window into the study of critical phenomena and discovering many unexpected results, their work became extremely influential. The approach developed in this series of works rested on an ambitious generalization of the method used in Exercise 1.2 for the one-dimensional Ising model based on "transfer matrices". Among other hurdles, it required the asymptotic analysis of the determinants of large Toeplitz matrices, which itself appealed to considerations involving orthogonal polynomials. Other approaches to determine the limit free energy were developed later, for instance by leveraging connections with domino tilings [68, 186]. The calculation of the limit free energy for the two-dimensional Ising model is, in some sense, miraculous, and generalizations to higher dimensions, $d \geqslant 3$, are unlikely to be found. In contrast, the models we consider in the rest of the book are of *mean-field* type. This means that the rich geometry of the integer lattice is replaced by a much simpler geometry, typically that of the complete graph, making the models invariant under any permutation of the index set. This very high degree of symmetry gives more hope that the identification of the limit free energy will be amenable to robust analytic techniques. For much more about the Ising model, we refer the interested reader to [117].

**Exercise 1.2.** For every inverse temperature $\beta > 0$ and magnetic intensity $h \in \mathbb{R}$, consider the one-dimensional free energy (1.24),

$$F_N(\beta, h) := \frac{1}{N} \log \sum_{\sigma \in \{\pm 1\}^N} \exp\left(\beta \sum_{i=1}^{N-1} \sigma_i \sigma_{i+1} + \beta \sigma_N \sigma_1 + h \sum_{i=1}^{N} \sigma_i\right), \quad (1.26)$$

(i) By rewriting $F_N$ in terms of the $N^{\text{th}}$ power of a 2-by-2 matrix, show that

$$\lim_{N \to +\infty} F_N(\beta, h) = \log\left(e^\beta \cosh(h) + \left(e^{2\beta} \cosh^2(h) - 2 \sinh(2\beta)\right)^{\frac{1}{2}}\right). \quad (1.27)$$

(ii) Deduce that the limit (1.25) of the mean magnetization is continuous in $h$. It may be helpful to refer to Propositions 2.11 and 2.15 as well as Exercise 2.6.

**Exercise 1.3.** Through a slight abuse of notation, consider the Hamiltonian

$$H_N(\sigma) := \sum_{\substack{i,j \in B_{N+1} \\ i \sim j}} \sigma_i \sigma_j \quad (1.28)$$

on $\{\pm 1\}^{B_N}$ with the boundary condition $\sigma_i := 1$ for $i \notin B_N$ and the convention that the sum (1.28) is taken over the edge set of the lattice $B_{N+1}$ as opposed to its vertex set. For each pair of neighbours $i \sim j \in B_{N+1}$ with $\sigma_i \sigma_j = -1$, draw an edge centred at the midpoint $\frac{i+j}{2}$ in the direction orthogonal to the edge $(i, j)$. These edges create contours that delimit regions in which the sign of the spins is constant. Denote by $\Gamma(\sigma)$ the set of these edges and consider the Gibbs measure defined for a bounded and measurable $f = f(\sigma)$ by

$$\langle f(\sigma) \rangle := \frac{\sum_{\sigma \in \{\pm 1\}^{B_N}} f(\sigma) \exp(\beta H_N(\sigma))}{\sum_{\sigma \in \{\pm 1\}^{B_N}} \exp(\beta H_N(\sigma))}. \quad (1.29)$$

(i) Show that $H_N(\sigma) + 2|\Gamma(\sigma)|$ does not depend on $\sigma \in \{\pm 1\}^{B_N}$.

(ii) Show that if $\sigma_0 = -1$, then the origin must be surrounded by at least one contour in $\Gamma(\sigma)$.

(iii) For any contour $\gamma$, show that $\langle \mathbf{1}_{\{\gamma \subseteq \Gamma(\sigma)\}} \rangle \leq \exp(-2\beta|\gamma|)$.

(iv) Show that for $\beta$ sufficiently large $\liminf_{N \to +\infty} \langle \sigma_0 \rangle > 0$.

(v) Conclude that for large $\beta$ the Ising model exhibits ferromagnetic behaviour.

## 1.3  The Curie-Weiss model

The Curie-Weiss model is the non-disordered mean-field model obtained by trivializing the geometry of the Ising model. In the Curie-Weiss model, the energy associated with each configuration $\sigma \in \{\pm 1\}^N := \{-1, +1\}^N$ at inverse temperature $t \geq 0$ and with external magnetic field of intensity $h \in \mathbb{R}$ is

$$H_N(t, h, \sigma) := \frac{t}{N} \sum_{i,j=1}^{N} \sigma_i \sigma_j + h \sum_{i=1}^{N} \sigma_i. \tag{1.30}$$

Up to relabelling the spin coordinates, this corresponds to the one-dimensional Ising model on the box $B_{N/2}$ where the restriction that a spin only interacts with its nearest neighbours is dropped, and a spin can instead interact with all other spins. To ensure that the Hamiltonian (1.30) remains of order $N$ despite these additional terms, the interaction term in (1.30) is divided by $N$. We have chosen to denote the inverse temperature parameter by $t \geq 0$ as opposed to $\beta \geq 0$ as this parameter can be interpreted as a time variable in the partial differential equation that will appear later. Just like in the Ising model, the probability of finding the Curie-Weiss model in the configuration $\sigma \in \{\pm 1\}^N$ is the Gibbs weight

$$G_N(t, h, \sigma) := \frac{\exp H_N(t, h, \sigma)}{Z_N(t, h)}, \tag{1.31}$$

with partition function

$$Z_N(t, h) := \sum_{\sigma \in \{\pm 1\}^N} \exp H_N(t, h, \sigma). \tag{1.32}$$

As in Exercise 1.3, we denote the Gibbs average of any bounded and measurable function $f : \{\pm 1\}^N \to \mathbb{R}$ by

$$\langle f(\sigma) \rangle := \frac{\sum_{\sigma \in \{\pm 1\}^N} f(\sigma) \exp H_N(t, h, \sigma)}{\sum_{\sigma \in \{\pm 1\}^N} \exp H_N(t, h, \sigma)}. \tag{1.33}$$

Notice that the bracket $\langle \cdot \rangle$ depends on $N$, $t$ and $h$, although this dependence is kept implicit in the notation. As discussed in Section 1.2, a key goal of ours is to understand the mean magnetization

$$m_N(t,h) := \left\langle \frac{1}{N} \sum_{i=1}^{N} \sigma_i \right\rangle \tag{1.34}$$

as a function of the model parameters, in the regime of large $N$. To understand this object in the setting of the Curie-Weiss model, we will first study the free energy

$$F_N(t,h) := \frac{1}{N} \log \frac{1}{2^N} \sum_{\sigma \in \{\pm 1\}^N} \exp H_N(t,h,\sigma). \tag{1.35}$$

The factor of $2^{-N}$ in (1.35) has been introduced so that we can interpret the normalized sum $2^{-N} \sum_{\sigma \in \{\pm 1\}^N}$ as integration against the uniform probability measure on $\{\pm 1\}^N$; its presence is inconsequential as it amounts to subtracting $\log 2$ from $F_N$.

In Chapter 2 we will use the classical theory of large deviations to explicitly compute the limit of the free energy (1.35). This will be possible as the Hamiltonian (1.30) depends on a spin configuration $\sigma \in \{\pm 1\}^N$ only through the sample average

$$S_N(\sigma) := \frac{1}{N} \sum_{i=1}^{N} \sigma_i. \tag{1.36}$$

Indeed, the Hamiltonian (1.30) can be written as

$$H_N(t,h,\sigma) = N\big(t S_N(\sigma)^2 + h S_N(\sigma)\big). \tag{1.37}$$

From the formula for the limit free energy, we will be able to leverage the envelope theorem discussed in Section 2.4 to understand the magnetization (1.34).

Unfortunately, for more complex mean-field models, it is not possible to determine the limit of the free energy by expressing the Hamiltonian as a function of a simple quantity such as $S_N(\sigma)$ whose large deviations can be figured out separately. In Chapter 3 we will therefore develop a Hamilton-Jacobi approach independent of large deviation principles to determine the limit of the free energy (1.35). To face as many of the difficulties presented by more complex models as possible in a simple setting where large deviation principles can be used to verify our results, we will also consider variants of the Curie-Weiss model which generalize it in two ways:

(i) we will replace the term $S_N(\sigma)^2$ in the Hamiltonian (1.37) by $\xi(S_N(\sigma))$ for an arbitrary smooth function $\xi \in C^\infty(\mathbb{R}; \mathbb{R})$;

(ii) we will replace the uniform probability measure on $\{\pm 1\}^N$ by an arbitrary probability measure $P_N$ on $\mathbb{R}^N$ whose support is contained in the closed Euclidean ball of radius $\sqrt{N}$ centred at the origin.

The free energy of the generalized Curie-Weiss model is therefore

$$F_N(t,h) := \frac{1}{N} \log \int_{\mathbb{R}^N} \exp N\big(t\xi(S_N(\sigma)) + hS_N(\sigma)\big) \, dP_N(\sigma). \qquad (1.38)$$

The free energy (1.35) can be recovered by setting $\xi(x) = x^2$ and by choosing $P_N$ to be the uniform probability measure on $\{\pm 1\}^N$. Once we have understood how to use the Hamilton-Jacobi approach to determine the limit of the free energy in the generalized Curie-Weiss model, we will see that the techniques developed along the way can be used to study many interesting models such as rank-one matrix estimation, community detection, and, at least partially, spin glass models.

# Chapter 2
# Convex analysis and large deviation principles

In this chapter we use the theory of large deviations to establish a variational formula for the limit of the free energy in the generalized Curie-Weiss model. In Section 2.1, we give a brief introduction to convex analysis. In particular, we prove the Fenchel-Moreau duality theorem and the local Lipschitz continuity of convex functions. While not really pertaining to convex analysis, we also show the Rademacher theorem on the almost everywhere differentiability of locally Lipschitz continuous functions. We conclude the section by introducing and studying the basic properties of subdifferentials of convex functions. In Section 2.2, we prove large deviation principles, which we then use in Section 2.3 to compute the limit of the free energy in the generalized Curie-Weiss model. Section 2.4 is devoted to the envelope theorem, which allows us to differentiate the variational formula for the limit free energy and study the limit mean magnetization in the Curie-Weiss model.

## 2.1 Convex analysis

Convex analysis is the branch of real analysis devoted to the study of convex sets and convex functions. Although this theory can be developed on very general vector spaces [108], we will restrict our attention to the Euclidean space $\mathbb{R}^d$. A set $C \subseteq \mathbb{R}^d$ is *convex* if for every $x, y \in C$ and $\alpha \in (0,1)$,

$$\alpha x + (1-\alpha)y \in C. \tag{2.1}$$

A function $f : \mathbb{R}^d \to \mathbb{R} \cup \{+\infty\}$ is *convex* if for every $x, y \in \mathbb{R}^d$ and $\alpha \in (0,1)$,

$$f(\alpha x + (1-\alpha)y) \leq \alpha f(x) + (1-\alpha)f(y). \tag{2.2}$$

This definition generalizes the classical definition of convexity of a real-valued function $f$ defined on a convex set $C \subseteq \mathbb{R}^d$. Indeed, if a function $f : C \to \mathbb{R}$ is convex in the classical sense, then extending it to be $+\infty$ on $\mathbb{R}^d \smallsetminus C$ gives a function which is convex in the sense (2.2). Conversely, if a function $f : \mathbb{R}^d \to \mathbb{R} \cup \{+\infty\}$ is convex in the sense (2.2), then it is convex in the classical sense on its *effective domain*

$$\mathrm{dom}\, f := \{x \in \mathbb{R}^d \mid f(x) < +\infty\}, \tag{2.3}$$

which is a convex set. We say that a function $f : \mathbb{R}^d \to \mathbb{R} \cup \{+\infty\}$ is *proper* if its effective domain is not empty. The convexity of functions and sets are intimately related. A function $f : \mathbb{R}^d \to \mathbb{R} \cup \{+\infty\}$ is convex if and only if its *epigraph*

$$\mathrm{epi}\, f := \{(x, \lambda) \in \mathbb{R}^d \times \mathbb{R} \mid f(x) \leq \lambda\} \tag{2.4}$$

is convex, while a set $C \subseteq \mathbb{R}^d$ is convex if and only if its *indicator function*

$$I_C(x) = \begin{cases} 0 & \text{if } x \in C \\ +\infty & \text{otherwise} \end{cases} \tag{2.5}$$

is convex. The notion of lower semi-continuity will also play an important role. We say that a function $f : \mathbb{R}^d \to \mathbb{R} \cup \{+\infty\}$ is *lower semi-continuous* if its epigraph is closed. It is readily verified that a function $f : \mathbb{R}^d \to \mathbb{R} \cup \{+\infty\}$ is lower semi-continuous if and only if, for every sequence $(x_n)_{n \geq 1} \subseteq \mathbb{R}^d$ converging to some point $x \in \mathbb{R}^d$, we have $f(x) \leq \liminf_{n \to +\infty} f(x_n)$.

**Exercise 2.1.** Let $C_1, C_2 \subseteq \mathbb{R}^d$ be convex sets. Show that the set difference

$$C_2 - C_1 := \{x_2 - x_1 \mid x_1 \in C_1 \text{ and } x_2 \in C_2\} \tag{2.6}$$

is also convex.

**Exercise 2.2.** Let $C$ be a convex set whose interior is not empty. Show that

$$\mathrm{int}(\overline{C}) = \mathrm{int}(C), \tag{2.7}$$

where we use the notation $\mathrm{int}(A)$ to denote the interior of a set $A$. Deduce that $\partial C = \partial \overline{C}$.

**Exercise 2.3.** Let $C$ be a closed convex set whose interior is not empty. Prove that $C$ is the closure of its interior,

$$C = \overline{\mathrm{int}(C)} \tag{2.8}$$

**Exercise 2.4.** Let $A$ be a subset of $\mathbb{R}^d$, and denote by

$$\mathrm{conv}(A) := \Big\{ \sum_{i=1}^n \alpha_i x_i \,\Big|\, x_i \in A, \lambda_i \in [0, 1] \text{ with } \sum_{i=1}^n \lambda_i = 1, \text{ and } n \geq 1 \Big\} \tag{2.9}$$

its convex hull. Show that any $x \in \mathrm{conv}(A)$ is a convex combination of at most $d + 1$ points of $A$.

**Exercise 2.5.** Let $A$ be a compact subset of $\mathbb{R}^d$. Prove that its convex hull $\mathrm{conv}(A)$ is also compact.

**Exercise 2.6.** Let $X$ be a bounded random variable taking values in $\mathbb{R}^d$. Show that the log-Laplace transform $\psi : \mathbb{R}^d \to \mathbb{R}$ defined by

$$\psi(\lambda) := \log \mathbb{E} \exp(\lambda \cdot X) \tag{2.10}$$

is convex.

**Exercise 2.7.** Fix an index set $I$, and for each $\alpha \in I$, let $f_\alpha : \mathbb{R}^d \to \mathbb{R} \cup \{+\infty\}$ be a convex function. Show that the supremum function

$$\sup_{\alpha \in I} f_\alpha \tag{2.11}$$

is convex. Prove the same statement with convex replaced by lower semi-continuous.

**Exercise 2.8.** Let $f = f(x,y) : \mathbb{R}^d \times \mathbb{R}^k \to \mathbb{R} \cup \{+\infty\}$ be a convex function (of the pair $(x,y)$). For every $x \in \mathbb{R}^d$, we define

$$g(x) := \inf_{y \in \mathbb{R}^k} f(x,y).$$

Assuming that $g$ takes values in $\mathbb{R} \cup \{+\infty\}$, show that $g$ is convex.

**Exercise 2.9.** We fix a function $f : \mathbb{R}^d \to \mathbb{R} \cup \{+\infty\}$ and define

$$\mathcal{L} := \left\{ g : \mathbb{R}^d \to \mathbb{R} \cup \{+\infty\} \mid g \leqslant f \text{ and } g \text{ is lower semi-continuous} \right\}. \tag{2.12}$$

Assuming that the set $\mathcal{L}$ is not empty, we define the lower-semicontinuous envelope $\underline{f} : \mathbb{R}^d \to \mathbb{R} \cup \{+\infty\}$ of the function $f$ as the supremum of all functions in $\mathcal{L}$, that is $\underline{f} = \sup_{g \in \mathcal{L}} g$. For every $x \in \mathbb{R}^d$, we write

$$\liminf_{y \to x} f(y) := \lim_{r \searrow 0} \inf\{f(y) \mid |y-x| \leqslant r\}. \tag{2.13}$$

(i) Show that $\underline{f}$ is lower semi-continuous.

(ii) Show that for every $x \in \mathbb{R}^d$, we have $\underline{f}(x) = \liminf_{y \to x} f(y)$.

(iii) Show that if $f$ is convex, then $\underline{f}$ is convex.

### 2.1.1   The Fenchel-Moreau duality theorem

The main result in convex analysis that we will discuss is the Fenchel-Moreau theorem. This result ensures that a convex and lower semi-continuous function is equal to its convex bi-dual. It will play an important role in obtaining variational formulas for limit free energies from the Hamilton-Jacobi approach. For instance, it will be used in Chapter 3 to recover the formula for the limit free energy in the Curie-Weiss model that we will prove in this chapter. The *convex dual* of a proper function $f : \mathbb{R}^d \to \mathbb{R} \cup \{+\infty\}$ is the function $f^* : \mathbb{R}^d \to \mathbb{R} \cup \{+\infty\}$ defined by

$$f^*(\lambda) := \sup_{x \in \mathbb{R}^d} (\lambda \cdot x - f(x)), \tag{2.14}$$

If $f^*$ is proper, then we can iterate this operation and obtain the *convex bi-dual* of $f$, which is the function $f^{**} : \mathbb{R}^d \to \mathbb{R} \cup \{+\infty\}$ defined by

$$f^{**}(x) := (f^*)^*(x) = \sup_{\lambda \in \mathbb{R}^d} (x \cdot \lambda - f^*(\lambda)). \tag{2.15}$$

For instance, if $f(x) := p \cdot x + a$ is the *affine function* with normal vector $p \in \mathbb{R}^d$ and intercept $a \in \mathbb{R}$, then its convex dual is given by

$$f^*(\lambda) := \begin{cases} -a & \text{if } \lambda = p \\ +\infty & \text{otherwise} \end{cases} \tag{2.16}$$

while its convex bi-dual is $f$ itself,

$$f^{**}(x) = p \cdot x + a = f(x). \tag{2.17}$$

This is a very important special case of the Fenchel-Moreau theorem. Indeed, to prove the Fenchel-Moreau theorem we will first show that a convex and lower semi-continuous function $f : \mathbb{R}^d \to \mathbb{R} \cup \{+\infty\}$ is the supremum of its affine minorants, and we will then leverage this observation as well as the basic properties of the convex dual established in Exercise 2.10. A function $g : \mathbb{R}^d \to \mathbb{R}$ is an *affine minorant* of $f$ if $g$ is an affine function and $g \leq f$. This representation result for convex functions will be deduced from the supporting hyperplane theorem, which itself is a special case of the Hahn-Banach separation theorem. The Hahn-Banach separation theorem is the geometrically intuitive statement that any two disjoint convex sets $C_1, C_2 \subseteq \mathbb{R}^d$ may be separated by a hyperplane. Although this result is classical, to keep the book as self-contained as possible, we will provide a full proof. We will first use projection operators onto closed and convex sets to show that convex sets can be linearly separated from singletons that they do not contain, and that this separation is strict if the convex set is closed. We will then deduce the Hahn-Banach separation theorem from this special case.

**Lemma 2.1** (Projection). *If $C \subseteq \mathbb{R}^d$ is a non-empty closed convex set and $x \in \mathbb{R}^d$, then there exists a unique point $P_C(x) \in C$ with*

$$|x - P_C(x)|^2 = \inf\{|x - y|^2 \mid y \in C\}. \tag{2.18}$$

*Moreover, a point $y \in C$ is the projection $P_C(x)$ of $x$ onto $C$ if and only if for all $z \in C$,*

$$(x - y) \cdot (z - y) \leq 0. \tag{2.19}$$

*Proof.* Fix $z \in C$, and observe that

$$\inf\{|x - y|^2 \mid y \in C\} = \inf\{|x - y|^2 \mid y \in C \text{ and } |x - y| \leq |x - z|\}.$$

The infimum on the right side is achieved at some point $P_C(x)$ as it consists in minimizing the continuous function $y \mapsto |x - y|^2$ over a compact set. To prove the uniqueness of this minimizer, suppose that $y_1, y_2$ are two minimizers, and let $\bar{y} = \frac{y_1 + y_2}{2}$ be their average. Introduce the differences $x_1 := y_1 - x$ and $x_2 := y_2 - x$, and observe that

$$|x_1|^2 + |x_2|^2 = \frac{1}{2}|x_1 + x_2|^2 + \frac{1}{2}|x_2 - x_1|^2.$$

This implies that

$$|y_1 - x|^2 + |y_2 - x|^2 = 2|\bar{y} - x|^2 + |y_2 - y_1|^2 \geq |y_1 - x|^2 + |y_2 - x|^2 + |y_2 - y_1|^2,$$

where the second inequality uses the fact that $|\bar{y} - x|^2 \geq |y_i - x|^2$ for $i = 1$ and $i = 2$. Rearranging shows that $y_2 = y_1$ and establishes the uniqueness of the projection $P_C(x)$. We now show that the projection $y = P_C(x)$ satisfies (2.19). Fix $z \in C$ and $\alpha \in (0, 1)$. The convexity of $C$ implies that $P_C(x) + \alpha(z - P_C(x)) \in C$, so the definition of $P_C(x)$ reveals that

$$|P_C(x) - x|^2 \leq |P_C(x) - x + \alpha(z - P_C(x))|^2$$
$$= |P_C(x) - x|^2 + 2\alpha(P_C(x) - x) \cdot (z - P_C(x)) + \alpha^2|z - P_C(x)|^2.$$

Rearranging, dividing by $\alpha$ and letting $\alpha$ tend to zero shows that $y = P_C(x)$ satisfies (2.19). Conversely, suppose that $y \in C$ satisfies (2.19), and fix $z \in C$. If $y = x$, it is clear that $y = P_C(x)$, so let us assume that $y \neq x$. Expanding (2.19) reveals that

$$0 \geq (x - y) \cdot (z - y) = (x - y) \cdot (z - x + x - y) = (x - y) \cdot (z - x) + |x - y|^2.$$

It follows by the Cauchy-Schwarz inequality that

$$|x - y|^2 \leq |x - y||z - x|.$$

Dividing by $|x - y| \neq 0$ and leveraging the uniqueness of $P_C(x)$ shows that $y = P_C(x)$, and completes the proof. ∎

**Theorem 2.2** (Supporting hyperplane). *If $C \subseteq \mathbb{R}^d$ is a non-empty convex set and $x \notin C$, then there exists a non-zero vector $v \in \mathbb{R}^d$ with*

$$v \cdot x \leqslant \inf\{v \cdot y \mid y \in C\}. \tag{2.20}$$

*Moreover, if $C$ is closed, then the vector $v$ can be chosen so that the inequality in (2.20) is strict.*

*Proof.* We decompose the proof into two steps.

*Step 1: closed $C$.* The vector $v := P_C(x) - x$ is well-defined and non-zero by Lemma 2.1. Moreover, its characterization (2.19) ensures that for any $z \in C$,

$$0 \geqslant (x - P_C(x)) \cdot (z - P_C(x)) = -v \cdot (z - v - x) = -v \cdot z + v \cdot x + |v|^2.$$

Rearranging shows that $v \cdot z \geqslant v \cdot x + |v|^2$ which establishes (2.20) with strict inequality.

*Step 2: general $C$.* Without loss of generality, we may assume that the interior of the set $C$ is not empty, that is, $\operatorname{int}(C) \neq \varnothing$. Indeed, if $\operatorname{int}(C)$ is empty, then $C$ must lie in an affine set of dimension less than $d$, and the normal vector $v \in \mathbb{R}^d$ to any hyperplane containing this affine set satisfies (2.20), provided that we adjust the sign of $v$ so that $v \cdot x \leqslant 0$.

We now distinguish two cases. First, if $x \notin \overline{C}$, then applying Step 1 to $\overline{C}$ gives (2.20) for $C$. On the other hand, if $x \in \overline{C} \setminus C \subseteq \partial C$, then Exercise 2.2 yields a sequence $(x_n)_{n \geqslant 1} \subseteq \mathbb{R}^d \setminus \overline{C}$ converging to $x$. Applying Step 1 to $\overline{C}$ and $x_n$ gives a sequence $(v_n)_{n \geqslant 1}$ of non-zero vectors with $v_n \cdot x_n \leqslant v_n \cdot z$ for all $z \in \overline{C}$. Since $v_n \neq 0$, the normalized vector

$$\bar{v}_n := \frac{v_n}{|v_n|}$$

is well-defined and lies in the unit ball. It follows by compactness of the unit ball that the sequence $(\bar{v}_n)_{n \geqslant 1}$ admits a subsequential limit $v \neq 0$ with $v \cdot x \leqslant v \cdot z$ for all $z \in \overline{C}$. This completes the proof. ∎

**Theorem 2.3** (Hahn-Banach). *If $C_1, C_2 \subseteq \mathbb{R}^d$ are disjoint non-empty convex sets, then there exists a non-zero vector $v \in \mathbb{R}^d$ with*

$$\sup\{v \cdot x_1 \mid x_1 \in C_1\} \leqslant \inf\{v \cdot x_2 \mid x_2 \in C_2\}. \tag{2.21}$$

*Proof.* Consider the set

$$C := C_2 - C_1 = \{x_2 - x_1 \mid x_1 \in C_1 \text{ and } x_2 \in C_2\},$$

and recall that it is convex by Exercise 2.1. Moreover, as $C_1$ and $C_2$ are disjoint, we have $0 \notin C$. It follows by the supporting hyperplane theorem that there exists a non-zero vector $v \in \mathbb{R}^d$ with $0 \leqslant v \cdot x$ for every $x \in C$, or equivalently $v \cdot x_1 \leqslant v \cdot x_2$ for all $x_1 \in C_1$ and $x_2 \in C_2$. This completes the proof. ∎

**Proposition 2.4** (Envelope of affine minorants). *A function $f : \mathbb{R}^d \to \mathbb{R} \cup \{+\infty\}$ is convex and lower semi-continuous if and only if it is the supremum of its affine minorants.*

*Proof.* Every affine function is convex and lower continuous, and the supremum of a family of convex and lower semi-continuous functions is also convex and lower semi-continuous by Exercise 2.7. This establishes the converse implication. Notice also that the equivalence is clearly valid when $f$ is identically $+\infty$. From now on, we therefore assume that $f$ is proper, convex and lower semi-continuous, we fix $x \in \mathbb{R}^d$ as well as $\mu < f(x)$, and we aim to find an affine function $g$ such that $g \leq f$ and $g(x) \geq \mu$. Since $f$ is convex and lower semi-continuous, its epigraph (2.4) is convex and closed. Moreover, it does not contain the pair $(x, \mu)$. It follows by the supporting hyperplane theorem for closed sets that there exist $p \in \mathbb{R}^d$, $a \in \mathbb{R}$, and $c_1 < c_2$ such that for every $(y, \lambda) \in \operatorname{epi} f$,

$$p \cdot x + a\mu = c_1 < c_2 \leq p \cdot y + a\lambda. \tag{2.22}$$

Recall that we assume that $f$ is proper; that is, there exists $x_0 \in \mathbb{R}^d$ such that $f(x_0) < +\infty$. Since $\{x_0\} \times [f(x_0), +\infty)$ is a subset of $\operatorname{epi} f$, we must have that $a \geq 0$. We decompose the rest of the proof into two steps. First we treat the case $a > 0$, and then the case $a = 0$.

*Step 1: $a > 0$.* Dividing (2.22) by $a > 0$, we find that $f(y) \geq a^{-1} p \cdot (x - y) + \mu$ for every $y \in \operatorname{dom} f$, and this inequality clearly extends to every $y \in \mathbb{R}^d$. The mapping $y \mapsto a^{-1} p \cdot (x - y) + \mu$ is thus an affine minorant of $f$, and takes the value $\mu$ at $x$, as desired.

*Step 2: $a = 0$.* As a preliminary step we show that the set of affine minorants of $f$ is not empty. Applying the same reasoning as that leading to (2.22) with $x$ replaced by $x_0$, and, say, with $\mu$ replaced by $f(x_0) - 1 < +\infty$, we get the existence of $p_0 \in \mathbb{R}^d$ and $a_0 \geq 0$ such that

$$p_0 \cdot x_0 + a_0(f(x_0) - 1) < p_0 \cdot y + a_0 \lambda$$

for every $(y, \lambda) \in \operatorname{epi} f$. Since $(x_0, f(x_0)) \in \operatorname{epi} f$, it cannot be that $a_0 = 0$. Arguing as in Step 1, we thus deduce that the mapping $g_0(y) := a_0^{-1} p_0 \cdot (x_0 - y) + f(x_0) - 1$ is an affine minorant of $f$. Coming back to the main thread of the argument, in the case when $a = 0$, we obtain from (2.22) that for every $y \in \operatorname{dom} f$,

$$c_2 - c_1 + p \cdot (x - y) \leq 0.$$

Using the affine minorant $g_0$ just constructed, we deduce that for every $M > 0$, and $y \in \operatorname{dom} f$,

$$f(y) \geq g_0(y) + M\bigl(c_2 - c_1 + p \cdot (x - y)\bigr), \tag{2.23}$$

and this inequality is in fact valid for every $y \in \mathbb{R}^d$. Since $c_2 - c_1 > 0$, we can make sure that the value on the right side of (2.23) at $y = x$ is as large as desired, in particular larger than $\mu$, by choosing $M$ sufficiently large. Once $M$ is thus chosen, the right side of (2.23) defines a suitable affine minorant of $f$. This completes the proof. ∎

The Fenchel-Moreau theorem extends this result by identifying an explicit set of affine minorants of $f$ whose supremum is $f$. In essence, for each slope $p$, we optimize the intercept to build the affine minorant with slope $p$ that just touches $f$. In Exercise 2.10, it is shown that the convex bi-dual of a proper function $f : \mathbb{R}^d \to \mathbb{R} \cup \{+\infty\}$ is a convex and lower semi-continuous minorant of $f$. The Fenchel-Moreau theorem states that as soon as $f$ is convex and lower semi-continuous, the inequality $f \geqslant f^{**}$ is in fact an equality.

**Theorem 2.5** (Fenchel-Moreau). *A proper function $f : \mathbb{R}^d \to \mathbb{R} \cup \{+\infty\}$ is convex and lower semi-continuous if and only if it is equal to its convex bi-dual,*

$$f = f^{**}. \tag{2.24}$$

**Remark 2.6.** Throughout this section, we chose to avoid considering functions that may also take the value $-\infty$; as a result, we required that $f^*$ be proper in order to define the bi-dual $f^{**}$. When stating the identity (2.24), we understand that it implicitly implies that $f^{**}$ is well-defined, or in other words, that $f^*$ is proper. If $f$ is not proper, that is if $f$ is constant and equal to $+\infty$, then we could decide that $f^* = -\infty$ and $f^{**} = +\infty$, which is indeed equal to $f$ in this case as well.

*Proof of Theorem 2.5.* Since $f^{**}$ is convex and lower semi-continuous by (i) in Exercise 2.10, the converse implication is clear. To show the direct implication, we fix a proper convex and lower semi-continuous function $f : \mathbb{R}^d \to \mathbb{R} \cup \{+\infty\}$, and aim to show that $f^{**} = f$. We start by showing that $f^{**}$ is well-defined and that $f^{**} \geqslant f$. Let $x \in \mathbb{R}^d$ and $\varepsilon > 0$. Invoking Proposition 2.4 gives an affine minorant $g$ of $f$ with $g(x) \geqslant f(x) - \varepsilon$. Since $g \leqslant f$, we can use (iii) in Exercise 2.10 to obtain that $g^* \geqslant f^*$. Recalling from (2.16) that the dual of an affine function is proper, we obtain that $f^*$ is proper as well. It thus follows that $f^{**}$ is well-defined, and another application of (iii) in Exercise 2.10 yields that $g^{**} \leqslant f^{**}$. Remembering that an affine function satisfies the Fenchel-Moreau theorem by (2.17) shows that in fact

$$f(x) - \varepsilon \leqslant g(x) \leqslant f^{**}(x).$$

Letting $\varepsilon$ tend to zero completes the proof of the inequality $f \leqslant f^{**}$. The converse inequality is immediate, see part (ii) of Exercise 2.10. ∎

**Remark 2.7.** As was pointed out in the previous remark, we chose throughout to exclude the possibility that convex functions take the value $-\infty$. If we were to allow for this possibility, we would then modify the definition and require that (2.2) holds whenever it is unambiguous; or equivalently, that it holds for every $x, y \in \mathbb{R}^d$ such that $f(x) < +\infty$ and $f(y) < +\infty$. If a convex function $f$ takes the value $-\infty$ at some point $x \in \mathbb{R}^d$, then on any half-line emanating from $x$, the function $f$ will either be $-\infty$ all the way; or else it must take the value $-\infty$ for a while, and then jump to $+\infty$, with the value of $f$ at the transition point being some arbitrary element of $\mathbb{R} \cup \{-\infty, +\infty\}$. If $f$ is also lower semi-continuous, then it must be $-\infty$ on some closed convex set, and $+\infty$ outside of it. There cannot be any affine minorant to a function that takes the value $-\infty$ somewhere, and in particular, if we were to define the bi-dual for such a function, it could only be the constant function equal to $-\infty$.

**Exercise 2.10.** Let $f : \mathbb{R}^d \to \mathbb{R} \cup \{+\infty\}$ be a proper function. Show that its convex dual satisfies the following properties.

(i) $f^*$ is convex and lower semi-continuous.

(ii) $f^{**} \leq f$.

(iii) If $f \leq g$ for some proper function $g : \mathbb{R}^d \to \mathbb{R} \cup \{+\infty\}$, then $f^* \geq g^*$.

**Exercise 2.11.** Consider the norm function $f(x) := \frac{1}{2}|x|^2$. Show that $f^* = f$.

**Exercise 2.12.** Let $f : \mathbb{R}^d \to \mathbb{R}$ be a real-valued Lipschitz continuous convex function with Lipschitz constant $L$. Show that $f^*(x) = +\infty$ for all $x \in \mathbb{R}^d$ with $|x| > L$.

**Exercise 2.13.** Let $f : \mathbb{R}^d \to \mathbb{R} \cup \{+\infty\}$ be a convex function, and denote by $\underline{f}$ the lower semi-continuous envelope of $f$, as defined in Exercise 2.9. Show that we have $f^* = (\underline{f})^*$ and $f^{**} = \underline{f}$.

**Exercise 2.14.** Let $K$ be a closed convex cone — this means that $K$ is a closed convex set with the additional property that, for all $x \in K$ and $\lambda > 0$, we have $\lambda x \in K$. We also assume that $K$ is not empty. The *polar* of $K$ is the closed convex cone

$$K^\circ := \{v \in \mathbb{R}^d \mid v \cdot x \leq 0 \text{ for all } x \in K\}. \tag{2.25}$$

Show that $K = K^{\circ\circ}$.

**Exercise 2.15.** We use the notion of cone and of the polar of a set $K$ as defined in Exercise 2.14. Let $C$ be a closed convex set. The *normal* to $C$ at a point $x \in C$ is the closed convex cone

$$\mathbf{n}_C(x) := \{v \in \mathbb{R}^d \mid v \cdot (x' - x) \leq 0 \text{ for all } x' \in C\}. \tag{2.26}$$

The *tangent* to $C$ at a point $x \in C$ is the closed convex cone

$$\mathbf{T}_C(x) := \overline{\{\lambda (x' - x) \mid x' \in C \text{ and } \lambda \geq 0\}}. \tag{2.27}$$

(i) Show that for all $x \in \text{int}(C)$, we have $\mathbf{n}_C(x) = \{0\}$ and $\mathbf{T}_C(x) = \mathbb{R}^d$.

(ii) Show that for all $x \in C$, we have $\mathbf{n}_C(x) = \mathbf{T}_C(x)^\circ$ and $\mathbf{T}_C(x) = \mathbf{n}_C(x)^\circ$.

(iii) Fix $x \in C$. Prove that $v \in \mathbf{T}_C(x)$ if and only if there is a sequence $(x_i)_{i \geq 1} \subseteq C$ converging to $x$ and a sequence $(t_i)_{i \geq 1} \subseteq \mathbb{R}_{>0}$ decreasing to 0 with $t_i^{-1}(x_i - x) \to v$.

**Exercise 2.16.** The purpose of this exercise is to generalize Lemma 2.1 to the Hilbert space setting. Let $H$ be a Hilbert space and let $C \subseteq H$ be a closed convex subset of $H$.

(i) Show that for every $x \in H$, there is a unique point $P_C(x) \in C$ with

$$\|x - P_C(x)\|^2 = \inf\{\|x - y\|^2 \mid y \in C\}. \tag{2.28}$$

(ii) Prove that a point $y \in C$ is the projection $P_C(x)$ of $x$ onto $C$ if and only if for all $z \in C$,

$$(x - y) \cdot (z - y) \leq 0. \tag{2.29}$$

**Exercise 2.17.** The purpose of this exercise is to prove the Riesz representation theorem on a Hilbert space $H$.

(i) Let $C$ be a closed subspace of $H$, and denote by

$$C^\perp := \{x \in H \mid x \cdot y = 0 \text{ for all } y \in C\} \tag{2.30}$$

its orthogonal complement. Use Exercise 2.16 to show that $H = C \oplus C^\perp$.

(ii) Let $f : H \to \mathbb{R}$ be a continuous linear functional on $H$. Show that there exists a unique $y \in H$ with $f(x) = x \cdot y$ for all $x \in H$.

### 2.1.2    Continuity and almost everywhere differentiability

In this section we will show that a convex function is locally Lipschitz continuous on the interior of its effective domain. Although this does not pertain to convex analysis, we will also show Rademacher's theorem stating that a locally Lipschitz continuous function is almost everywhere differentiable. The multi-dimensional version of Rademacher's theorem relies on its one-dimensional counterpart whose proof is the content of Exercise 2.18.

We start by showing that a convex function that is locally bounded from above is in fact locally bounded and locally Lipschitz continuous. For every $r > 0$ and $x \in \mathbb{R}^d$, we denote by $B_r(x)$ the open Euclidean ball of radius $r$ centred at $x$.

**Lemma 2.8.** *Let $f : \mathbb{R}^d \to \mathbb{R} \cup \{+\infty\}$ be a proper and convex function. If there exist $\delta > 0$ and $x^* \in \mathbb{R}^d$ such that $f$ is bounded from above by $M > 0$ on the open ball $B_{2\delta}(x^*)$, then the following assertions hold.*

(i) *There exists $m \in \mathbb{R}$ such that for every $x \in B_{2\delta}(x^*)$, we have $m \leq f(x) \leq M$.*

(ii) *For every $x, y \in B_\delta(x^*)$, we have $|f(x) - f(y)| \leq \frac{M-m}{\delta}|x-y|$.*

*Proof.* The convexity of $f$ reveals that for any $x \in B_{2\delta}(x^*)$,

$$f(x^*) = f\left(\frac{1}{2}(2x^* - x) + \frac{1}{2}x\right) \leq \frac{1}{2}f(x^* + (x^* - x)) + \frac{1}{2}f(x) \leq \frac{1}{2}M + \frac{1}{2}f(x).$$

This implies that $f(x) \geq 2f(x^*) - M$, so the first assertion holds with $m := 2f(x^*) - M$. To prove the second assertion, fix $x, y \in B_\delta(x^*)$ and consider the point $z := y + \delta \frac{y-x}{|y-x|}$ in $B_{2\delta}(x^*)$. Since $y$ lies on the line segment joining $x$ and $z$,

$$y = \frac{|y-x|}{\delta + |y-x|} z + \frac{\delta}{\delta + |y-x|} x,$$

the convexity of $f$ implies that

$$f(y) - f(x) \leq \frac{|y-x|}{\delta + |y-x|}(f(z) - f(x)) \leq \frac{M-m}{\delta}|y-x|.$$

Interchanging the roles of $x$ and $y$ completes the proof. ∎

**Proposition 2.9** (continuity). *Let $f : \mathbb{R}^d \to \mathbb{R} \cup \{+\infty\}$ be a proper and convex function. If $C \subseteq \text{int}(\text{dom}(f))$ is a compact set, then there exists $L < +\infty$ such that for every $x, y \in C$,*

$$|f(x) - f(y)| \leq L|x - y|. \tag{2.31}$$

*Proof.* For each $x \in \mathbb{R}^d$ and $r > 0$, denote by $C_r(x) := \prod_{i=1}^{d}[x_i - r, x_i + r]$ the cube of side-length $2r$ centred at $x$, or equivalently the closed $\ell^\infty$-ball of radius $r$ centred at $x$. We fix $x^* \in C \subseteq \text{int}(\text{dom} f)$, and pick $\delta = \delta(x^*) > 0$ and $r > 0$ such that $x^* \in B_{2\delta}(x^*) \subseteq C_r(x^*) \subseteq \text{dom}(f)$. Denote by $(v_i)_{i \leq 2^d}$ the vertices of $C_r(x^*)$, and observe that any point in the cube $C_r(x^*)$ may be written as a convex combination of the vertices $(v_i)_{i \leq 2^d}$. In particular, for any $y \in B_{2\delta}(x^*)$, there exist non-negative scalars $(\alpha_i)_{i \leq 2^d}$ with

$$y = \sum_{i=1}^{2^d} \alpha_i v_i \quad \text{and} \quad \sum_{i=1}^{2^d} \alpha_i = 1.$$

It follows by convexity of $f$ that $f(y) \leq \max\{f(v_0), \ldots, f(v_n)\}$. Invoking Lemma 2.8 shows that $f$ is Lipschitz continuous on $B_\delta(x^*)$, and remembering that $C$ is compact implies that $f$ is in fact Lipschitz continuous on $C$. This completes the proof. ∎

Together with Rademacher's theorem on the almost everywhere differentiability of Lipschitz functions, this result implies that a real-valued convex function is almost everywhere differentiable.

**Theorem 2.10** (Rademacher). *If $U \subseteq \mathbb{R}^d$ is open and $f : U \to \mathbb{R}$ is Lipschitz continuous, then $f$ is differentiable almost everywhere on $U$.*

*Proof.* The proof proceeds in three steps. First, we use the one-dimensional Rademacher theorem to show that for every direction vector $v \in \mathbb{R}^d$, the directional derivative

$$D_v f(x) := \frac{d}{dt}\bigg|_{t=0} f(x+tv) \tag{2.32}$$

exists almost everywhere on $U$, then we show that it is almost everywhere given by $\nabla f(x) \cdot v$, and finally we use the separability and compactness of the unit sphere and the Lipschitz continuity of $f$ to conclude that $f$ is differentiable almost everywhere.

*Step 1: almost everywhere existence of the directional derivative.* Fix a direction vector $v \in \mathbb{R}^d \smallsetminus \{0\}$, and denote by

$$A_v := \{x \in U \mid D_v f(x) \text{ exists}\}$$

the set of points where the derivative of $f$ in the direction $v$ exists. Denote by $v^\perp$ the orthogonal complement of the subspace $\mathbb{R}v$, so that by Exercise 2.17, we have the decomposition $\mathbb{R}^d = \mathbb{R}v \oplus v^\perp$. By the one-dimensional Rademacher theorem established in Exercise 2.18, for each $y \in v^\perp$, the intersection of the set $U \smallsetminus A_v$ with the line parallel to $\mathbb{R}v$ through the point $y$ is of zero measure. It follows by the Fubini-Tonelli theorem that $U \smallsetminus A_v$ is of zero measure, and therefore that the directional derivative $D_v f$ exists almost everywhere on $U$.

*Step 2: $D_v f(x) = \nabla f(x) \cdot v$ almost everywhere.* Fix a direction vector $v \in \mathbb{R}^d$ as well as a smooth function $\phi \in C_c^\infty(U; \mathbb{R})$ of compact support. A linear change of variables reveals that for any $t > 0$ small enough,

$$\int_U \left(\frac{f(x+tv) - f(x)}{t}\right) \phi(x) \, dx = -\int_U f(x) \left(\frac{\phi(x-tv) - \phi(x)}{-t}\right) dx.$$

Since $f$ is Lipschitz continuous and $\phi$ is compactly supported, both integrands in this equality are uniformly bounded and vanish outside a compact set. Using the dominated convergence theorem to let $t$ tend to zero shows that

$$\int_U D_v f(x) \phi(x) \, dx = -\int_U f(x) (\nabla \phi(x) \cdot v) \, dx.$$

In particular, choosing $v = e_i$ for the canonical basis vector $e_i \in \mathbb{R}^d$ reveals that

$$\int_U \partial_{x_i} f(x) \phi(x) \, dx = -\int_U f(x) \partial_{x_i} \phi(x) \, dx.$$

It follows that

$$\int_U D_v f(x) \phi(x) \, dx = -\sum_{i=1}^d v_i \int_U f(x) \partial_{x_i} \phi(x) \, dx = \int_U (\nabla f(x) \cdot v) \phi(x) \, dx, \tag{2.33}$$

where here we just use $\nabla f(x)$ as a notation for

$$\nabla f(x) = (\partial_{x_1} f(x), \ldots, \partial_{x_d} f(x)) = (D_{e_1} f(x), \ldots, D_{e_d} f(x)),$$

which is well-defined almost everywhere on $U$ by the previous step. Since the function $D_v f - \nabla f \cdot v$ is uniformly bounded by the Lipschitz continuity of $f$, applying the equality (2.33) along a sequence $(\phi_n)_{n \geq 1}$ of smooth functions of compact support that approximate $D_v f - \nabla f \cdot v$ in $L^1(\mathbb{R}^d; \mathbb{R})$ shows that $D_v f = \nabla f \cdot v$ almost everywhere on $U$.

*Step 3: concluding that $f$ is differentiable almost everywhere.* Let $(v_i)_{i \geq 1} \subseteq S^{d-1}$ be a countable set of direction vectors that is dense in the unit sphere $S^{d-1} \subseteq \mathbb{R}^d$. By the additivity of measure as well as Steps 1 and 2, the complement of the set

$$A := \bigcap_{i=1}^{\infty} \{x \in U \mid \nabla f(x) \text{ exists}, D_{v_i} f(x) \text{ exists and } D_{v_i} f(x) = \nabla f(x) \cdot v_i\}$$

has zero measure. Here we understand the phrase "$\nabla f(x)$ exists" to mean that $\partial_{x_i} f(x) = D_{e_i} f(x)$ exists for every $i \in \{1, \ldots, d\}$.

We now fix $x \in A$ and show that $f$ is differentiable at $x$. If we introduce the error function

$$R(x, v, t) := \frac{f(x+tv) - f(x)}{t} - \nabla f(x) \cdot v,$$

this comes down to proving that, given $\varepsilon > 0$, it is possible to find $\delta = \delta(\varepsilon, x) > 0$ such that for all $v \in S^{d-1}$ and $t > 0$ with $|t| < \delta$, we have

$$|R(x, v, t)| \leq \varepsilon. \tag{2.34}$$

Writing $L$ for the Lipschitz constant of $f$, we have that $|\partial_{x_i} f(x)| \leq L$ and $|\nabla f(x)| \leq \sqrt{d} L$, so the Cauchy-Schwarz inequality implies that $R$ is Lipschitz continuous in $v$. Indeed, for any $v, v' \in S^{d-1}$, we have

$$|R(x, v, t) - R(x, v', t)| \leq (\sqrt{d} + 1) L |v - v'|.$$

By the compactness of $S^{d-1}$, there is a finite sub-collection $(v_i)_{i \leq n}$ of direction vectors with the property that the open balls $(B_\varepsilon(v_i))_{i \leq n}$ cover the unit sphere $S^{d-1}$. Since each error term $R(x, v_i, t)$ vanishes as $t$ tends to zero, it is possible to find $\delta > 0$ such that $|R(x, v_i, t)| \leq \varepsilon$ whenever $|t| < \delta$ and $1 \leq i \leq n$. Now, if $v \in S^{d-1}$ belongs to the open ball $B_\varepsilon(v_i)$ and $|t| < \delta$, then

$$|R(x, v, t)| \leq |R(x, v_i, t)| + |R(x, v_i, t) - R(x, v, t)| \leq \varepsilon + (\sqrt{d} + 1)\varepsilon.$$

Redefining $\varepsilon$ establishes (2.34) and completes the proof. ∎

**Exercise 2.18.** The purpose of this exercise is to prove the one-dimensional Rademacher theorem for a Lipschitz function $F : [a, b] \to \mathbb{R}$.

(i) For any step function of the form $\phi := \sum_{i=1}^{n} \phi_i \mathbf{1}_{[a_{i-1}, a_i)}$, with $\phi_1, \ldots, \phi_n \in \mathbb{R}$ and $a \leqslant a_0 \leqslant \cdots \leqslant a_n \leqslant b$, define

$$T(\phi) := \sum_{i=1}^{n} \phi_i \big( F(a_i) - F(a_{i-1}) \big). \tag{2.35}$$

Prove that the operator $T$ admits a unique extension to a continuous functional on $L^2([a,b]; \mathbb{R})$.

(ii) Use Exercise 2.17 to find $f \in L^2([a,b]; \mathbb{R})$ such that, for all $x \in [a,b]$, we have

$$F(x) = F(a) + \int_a^x f(t) \, dt. \tag{2.36}$$

(iii) Using the Lebesgue differentiation theorem in Theorem A.16, conclude that the function $F$ is almost everywhere differentiable on $[a,b]$.

**Exercise 2.19.** Let $A$ be a $d$-by-$k$ matrix, $b \in \mathbb{R}^d$, and $c \in \mathbb{R}^k$. We denote by $A^*$ the transpose of the matrix $A$, and by $\leqslant$ the partial order on $\mathbb{R}^d$ given by its product structure; in other words, for every $x, y \in \mathbb{R}^d$, we write $x \leqslant y$ to mean that the inequality holds coordinate by coordinate. The goal of this exercise is to identify a general condition under which the "linear programs"

$$\inf \{ b \cdot x \mid x \in \mathbb{R}^d, \ x \geqslant 0, \ Ax \geqslant c \} \quad \text{and} \quad \sup \{ c \cdot y \mid y \in \mathbb{R}^k, \ y \geqslant 0, \ A^* y \leqslant b \} \tag{2.37}$$

are equal.

(i) Slightly adapting the notation in (2.5), for every set $A$, we may write $I_{\{x \in A\}}$ instead of $I_A(x)$; for instance, we understand that $I_{\{x \geqslant 0\}}$ is equal to 0 if $x \geqslant 0$, and is $+\infty$ otherwise. What is the convex dual of the mapping $x \mapsto I_{\{x \geqslant 0\}}$?

(ii) For every $x, z \in \mathbb{R}^d$, we write $\phi(x) := I_{\{x \geqslant 0\}} + I_{\{Ax \geqslant c\}}$ and

$$\psi(z) := \inf_{y \in \mathbb{R}^k} \big( c \cdot y + I_{\{y \leqslant 0\}} + I_{\{A^* y \geqslant z\}} \big). \tag{2.38}$$

Show that $\phi = \psi^*$.

(iii) For each $z \in \mathbb{R}^d$, we define the set

$$K_z := \{ y \in \mathbb{R}^k \mid y \geqslant 0 \text{ and } A^* y \leqslant z \}. \tag{2.39}$$

Assuming that $K_z$ is a compact set for every $z$ in a neighbourhood of $b$, show that the two variational problems in (2.37) are equal.

### 2.1.3 The subdifferential of a convex function

The local Lipschitz continuity of a real-valued convex function and the Rademacher theorem imply that a real-valued convex function is differentiable almost everywhere. It is at times useful to also be able to study the differentiability of a convex function at every single point of its effective domain using the notion of subdifferential. The *subdifferential* of a convex function $f : \mathbb{R}^d \to \mathbb{R} \cup \{+\infty\}$ at a point $x \in \mathrm{dom}\, f$ is the set

$$\partial f(x) := \{p \in \mathbb{R}^d \mid f(y) \geq f(x) + p \cdot (y - x) \text{ for all } y \in \mathbb{R}^d\}. \tag{2.40}$$

We start by showing that the subdifferential of $f$ at a point in the interior of its effective domain is not empty.

**Proposition 2.11.** *Let $f : \mathbb{R}^d \to \mathbb{R} \cup \{+\infty\}$ be a convex function. For every point $x \in \mathrm{int}(\mathrm{dom}\, f)$, the subdifferential $\partial f(x)$ is not empty.*

*Proof.* We consider the convex set $\mathcal{C} := \{(y, \mu) \in \mathrm{dom}(f) \times \mathbb{R} \mid f(y) < \mu\}$, and fix $x \in \mathrm{int}(\mathrm{dom}\, f)$. Since $(x, f(x)) \notin \mathcal{C}$ the supporting hyperplane theorem gives a non-zero vector $(v, b) \in \mathbb{R}^d \times \mathbb{R}$ with

$$0 \leq v \cdot (y - x) + b(\mu - f(x))$$

for every $(y, \mu) \in \mathcal{C}$. Since $\mu$ can be arbitrarily large, we must have $b \geq 0$. If we had $b = 0$, then we would have $0 \leq v \cdot (y - x)$ for all $y$ in a neighbourhood of $x \in \mathrm{int}(\mathrm{dom}\, f)$, which is not possible since $(v, b)$ is non-zero. This means that $b > 0$ so the vector $p := -v/b$ is well-defined and satisfies

$$\mu \geq f(x) + p \cdot (y - x)$$

for all $(y, \mu) \in \mathcal{C}$. Letting $\mu$ tend to $f(y)$ reveals that $p \in \partial f(x)$, which means that the subdifferential $\partial f(x)$ is not empty. ∎

It turns out that a convex function is differentiable at a point $x$ in the interior of its effective domain if and only its subdifferential at $x$ consists of a singleton. To prove this, we will leverage the fact that a convex function is differentiable at a point in the interior of its effective domain if and only if its directional derivative is a linear function of the direction. The *directional derivative* of a real-valued convex function $f : \mathbb{R}^d \to \mathbb{R} \cup \{+\infty\}$ at $x \in \mathrm{int}(\mathrm{dom}\, f)$ in the direction of $v \in \mathbb{R}^d$ is the function

$$D_v f(x) := \lim_{t \searrow 0} \frac{f(x + tv) - f(x)}{t}. \tag{2.41}$$

The convexity of $f$ implies that the difference quotient defining the directional derivative is a decreasing function of $t$, so the directional derivative is well-defined as a monotone limit. In general, the linearity of the directional derivative does not suffice to characterize the differentiability of a function; however, it does suffice for locally Lipschitz continuous functions, a fact that we already used at least implicitly in the proof of the Rademacher theorem. Combining this with the local Lipschitz continuity of convex functions yields the following result.

**Lemma 2.12.** *A convex function $f : \mathbb{R}^d \to \mathbb{R} \cup \{+\infty\}$ is differentiable at $x \in \text{int}(\text{dom } f)$ if and only if the map $v \mapsto D_v f(x)$ is linear. In this case, we have $D_v f(x) = \nabla f(x) \cdot v$.*

*Proof.* On the one hand, if $f$ is differentiable, then for any direction vector $v \in \mathbb{R}^d$, we have $D_v f(x) = \nabla f(x) \cdot v$. Conversely, suppose that $v \mapsto D_v f(x)$ is linear, that is, there exists $a \in \mathbb{R}^d$ such that $D_v f(x) = a \cdot v$ for every $v \in \mathbb{R}^d$. We assume for the sake of contradiction that $f$ is not differentiable at $x$. Let $(v_n)_{n \geqslant 1}$ with $|v_n| = 1$ for all $n \geqslant 1$ and $(t_n)_{n \geqslant 1} \subseteq \mathbb{R}_{>0}$ be a sequence converging to 0 such that the error term

$$R(x, v_n, t_n) := \left| \frac{f(x + t_n v_n) - f(x)}{t_n} - a \cdot v_n \right|$$

does not converge to zero. Up to passing to a subsequence, assume that $(v_n)_{n \geqslant 1}$ converges to some $v_0$ in the unit sphere. Remembering Proposition 2.9, denote by $L < +\infty$ the Lipschitz constant of $f$ in a neighbourhood of $x$. The triangle inequality implies that

$$R(x, v_n, t_n) \leqslant R(x, v_0, t_n) + (L + |a|) |v_n - v_0|.$$

Leveraging the assumption that $D_{v_0} f(x) = a \cdot v_0$ to let $n$ tend to infinity contradicts the absurd hypothesis that the sequence $(R(x, v_n, t_n))_{n \geqslant 1}$ does not converge to zero. This completes the proof. ∎

**Theorem 2.13.** *A convex function $f : \mathbb{R}^d \to \mathbb{R} \cup \{+\infty\}$ is differentiable at a point $x \in \text{int}(\text{dom } f)$ if and only if $\partial f(x)$ consists of a singleton. In this case, we have $\partial f(x) = \{\nabla f(x)\}$.*

*Proof.* We first show the direct implication. Recall from Proposition 2.11 that the subdifferential $\partial f(x)$ is not empty. We fix $p \in \partial f(x)$. By definition of the subdifferential, for every $v \in \mathbb{R}^d$ and $\lambda > 0$,

$$f(x + \lambda v) - f(x) \geqslant \lambda v \cdot p.$$

Dividing by $\lambda$ and letting $\lambda$ tend to zero shows that $(\nabla f(x) - p) \cdot v \geqslant 0$. Choosing $v = p - \nabla f(x)$ reveals that $p = \nabla f(x)$, so $\partial f(x) = \{\nabla f(x)\}$.

Conversely, suppose that the subdifferential is a singleton, $\partial f(x) = \{p\}$, and fix a direction vector $v \in \mathbb{R}^d$. The convexity of $f$ and the definition of the directional derivative imply that for all $\lambda \in \mathbb{R}$,

$$f(x) + \lambda D_v f(x) \leqslant f(x + \lambda v).$$

This implies that the convex sets

$$C := \{(x + \lambda v, f(x) + \lambda D_v f(x)) \mid \lambda \in \mathbb{R}\}$$

and

$$C' := \{(y, \mu) \in \operatorname{dom} f \times \mathbb{R} \mid f(y) < \mu\}$$

are disjoint. It follows by the Hahn-Banach separation theorem that there exists a non-zero vector $(a, b) \in \mathbb{R} \times \mathbb{R}^d$ with

$$a(f(x) + \lambda D_v f(x)) + b \cdot (x + \lambda v) \leq a\mu + b \cdot y \tag{2.42}$$

for all $(y, \mu) \in C'$ and $\lambda \in \mathbb{R}$. Taking $\lambda = 0$ shows that

$$af(x) + b \cdot x \leq a\mu + b \cdot y \tag{2.43}$$

for all $(y, \mu) \in C'$. Since $\mu$ can be arbitrarily large, we must have $a \geq 0$. If we had $a = 0$, then we would have $0 \leq b \cdot (y - x)$ for all $y$ in a neighbourhood of $x \in \operatorname{int}(\operatorname{dom} f)$, which is not possible since $(a, b)$ is non-zero. Dividing through by $a$ and letting $\mu$ tend to $f(y)$ in (2.43) shows that $-b/a \in \partial f(x)$, and therefore $b/a = -p$. Combining this with (2.42) and letting $\mu$ tend to $f(y)$ in the resulting bound gives

$$f(x) + \lambda D_v f(x) - p \cdot (x + \lambda v) \leq f(y) - p \cdot y$$

for all $y \in \operatorname{dom} f$ and $\lambda \in \mathbb{R}$. Taking $y = x$ reveals that $\lambda(D_v f(x) - p \cdot v) \leq 0$ for all $\lambda \in \mathbb{R}$, which implies that $D_v f(x) = p \cdot v$. In particular, the map $v \mapsto D_v f(x)$ is linear. Invoking Lemma 2.12 completes the proof. ∎

The next proposition gives a sort of continuity property of the subdifferential as we vary the base point.

**Proposition 2.14.** *Let $f : \mathbb{R}^d \to \mathbb{R} \cup \{+\infty\}$ be a convex function, and let $(x_n, p_n)_{n \geq 1}$ be a sequence of points in $\operatorname{dom} f \times \mathbb{R}^d$ with $p_n \in \partial f(x_n)$ for each $n \geq 1$ that converges to some point $(x, p) \in \operatorname{dom} f \times \mathbb{R}^d$. If $f$ is continuous at $x \in \operatorname{dom} f$, then $p \in \partial f(x)$.*

*Proof.* Fix $y \in \mathbb{R}^d$ as well as $n \geq 1$. Since $p_n \in \partial f(x_n)$, we have

$$f(y) \geq f(x_n) + p_n \cdot (y - x_n)$$

Letting $n$ tend to infinity and using the continuity of $f$ at $x$ completes the proof. ∎

We recall that if the point $x$ in Proposition 2.14 belongs to $\operatorname{int}(\operatorname{dom} f)$, then the continuity of $f$ at $x$ is automatically satisfied by Proposition 2.9. Although we will not use this fact, we also mention that one can use Proposition 2.14 to show that if a convex function $f$ is uniformly Lipschitz continuous on its effective domain, then the subdifferential of $f$ at every point $x \in \operatorname{dom} f$ is not empty, including at points $x$ that sit on the boundary of $\operatorname{dom} f$.

**Proposition 2.15.** *Let $x \in \mathbb{R}^d$, and for each integer $n \geqslant 1$, let $f_n : \mathbb{R}^d \to \mathbb{R} \cup \{+\infty\}$ be a convex function such that $x \in \text{int}(\text{dom} f_n)$ and with $f_n$ differentiable at $x \in \mathbb{R}^d$. If $(f_n)_{n \geqslant 1}$ converges pointwise to some function $f : \mathbb{R}^d \to \mathbb{R} \cup \{+\infty\}$ with $x \in \text{dom} f$, and if the sequence of derivatives $(\nabla f_n(x))_{n \geqslant 1}$ converges to some vector $p \in \mathbb{R}^d$, then $p \in \partial f(x)$.*

*Proof.* Since $f_n$ is differentiable at $x \in \text{int}(\text{dom} f)$, we have $\nabla f_n(x) \in \partial f_n(x)$ by Theorem 2.13. It follows by definition of the subdifferential that, for every $y \in \mathbb{R}^d$,

$$f_n(y) \geqslant f_n(x) + \nabla f_n(x)(y - x).$$

Letting $n$ tend to infinity and using the pointwise convergence of $f_n$ to $f$ completes the proof. ∎

For a more in-depth discussion of convex analysis, we refer the interested reader to [53, 108, 138, 187, 230].

**Exercise 2.20.** We use the notation $\text{conv}(C)$ from Exercise 2.4 for the convex hull of a set $C$, and the notation $\mathbf{n}_C(x)$ from Exercise 2.15 for the normal cone at $x \in C$. Let $f : \mathbb{R}^d \to \mathbb{R} \cup \{+\infty\}$ be a convex function such that $\text{dom} f$ has non-empty interior. For each $x \in \text{dom} f$, let

$$S(x) := \Big\{ \lim_{i \to +\infty} \nabla f(x_i) \,\Big|\, (x_i)_{i \geqslant 1} \subseteq \text{int}\,\text{dom}\,f \text{ converges to } x \text{ and}$$

$$\text{for every } i \geqslant 1, \text{ the function } f \text{ is differentiable at } x_i \Big\} \quad (2.44)$$

be the set of limits of gradients of $f$ along sequences that converge to $x$. Suppose that $f$ is Lipschitz on $\text{dom} f$, and show that

$$\partial f(x) = \mathbf{n}_{\text{dom} f}(x) + \text{conv}(S(x)). \quad (2.45)$$

Although we do not prove this, let us mention that the Lipschitz assumption is not necessary. Indeed, if $f : \mathbb{R}^d \to \mathbb{R} \cup \{+\infty\}$ is a convex and lower semicontinuous function such that $\text{dom} f$ has non-empty interior, then

$$\partial f(x) = \mathbf{n}_{\text{dom} f}(x) + \overline{\text{conv}}(S(x)). \quad (2.46)$$

We refer the interested reader to Theorem 25.6 in [230] or Theorem 4.5 in [50] for a proof of this more general formula.

## 2.2 Large deviation principles

One of the first quantities that one tries to understand in any basic course in probability theory is the sample average of a sequence of i.i.d. coin tosses. More specifically, given a collection of independent Bernoulli random variables $(X_n)_{n \geq 1}$ with mean $p \in [0,1]$, one is often presented with two main results about the sample average

$$S_N := \frac{1}{N} \sum_{n=1}^{N} X_n \qquad (2.47)$$

of $N$ of these random variables. The first is the law of large numbers, which states that the sample average $S_N$ converges almost surely to the mean $p$ as the number of trials $N$ goes to infinity. The second is the central limit theorem, which is the assertion that the typical deviations of the sample average from the mean $p$ are of order $1/\sqrt{N}$ and are normally distributed. Large deviation principles are concerned with the rare deviations, of order one, of the sample average $S_N$ from its mean $p$.

In the setting of i.i.d. coin tosses, studying large deviations comes down to analyzing the asymptotic behaviour of the probabilities

$$\mathbb{P}\left\{S_N = \frac{k}{N}\right\} = \binom{N}{k} p^k (1-p)^{N-k} \qquad (2.48)$$

for $k \in \{0,\ldots,N\}$. Stirling's formula implies that for a fixed $x = \frac{k}{N} \in (0,1)$, we have

$$\log \binom{N}{k} = -N\big(x \log(x) + (1-x) \log(1-x)\big) + \mathcal{O}(\log(N)). \qquad (2.49)$$

This means that

$$\log \mathbb{P}\left\{S_N = \frac{k}{N}\right\} = -N I(x) + \mathcal{O}(\log(N)) \qquad (2.50)$$

for the rate function

$$I(x) := x \log\left(\frac{x}{p}\right) + (1-x) \log\left(\frac{1-x}{1-p}\right). \qquad (2.51)$$

The probability that $S_N$ lies above $x \in [p,1)$ is therefore bounded from below by

$$\mathbb{P}\{S_N \geq x\} \geq \mathbb{P}\left\{S_N = \frac{\lfloor Nx \rfloor}{N}\right\} = \exp(-N I(x) + \mathcal{O}(\log(N))), \qquad (2.52)$$

and bounded from above by

$$\mathbb{P}\{S_N \geq x\} \leq \sum_{k=\lfloor xN \rfloor}^{N} \mathbb{P}\left\{S_N = \frac{k}{N}\right\} \leq \exp\big(-N I(x) + \mathcal{O}(\log(N))\big), \qquad (2.53)$$

where we have used that $I$ is increasing on the interval $[p,1)$ in the second inequality. Combining these bounds and letting $N$ tend to infinity reveals that for every $x \in [p,1)$,

$$\lim_{N \to +\infty} \frac{1}{N} \log \mathbb{P}\{S_N \geq x\} = -I(x). \qquad (2.54)$$

Similar arguments can be used to show that for $x \in (0, p]$.

$$\lim_{N \to +\infty} \frac{1}{N} \log \mathbb{P}\{S_N \leqslant x\} = -I(x). \tag{2.55}$$

In particular, the large deviations of $(S_N)_{N \geqslant 1}$ from its mean $p$ are entirely described by the rate function $I$. As shown in Exercise 2.21, the equalities (2.54) and (2.55) are equivalent to the fact that for every Borel set $A \subseteq \mathbb{R}$,

$$-\inf_{x \in \mathrm{int}(A)} I(x) \leqslant \liminf_{N \to +\infty} \frac{1}{N} \log \mathbb{P}\{S_N \in A\}$$
$$\leqslant \limsup_{N \to +\infty} \frac{1}{N} \log \mathbb{P}\{S_N \in A\} \leqslant -\inf_{x \in \overline{A}} I(x), \tag{2.56}$$

where we recall that we write $\mathrm{int}(A)$ to denote the interior of the set $A$. This string of inequalities is more amenable to generalization than (2.54) and (2.55) as it requires neither that we stipulate the mean $p$ of the random variables $(S_N)_{N \geqslant 1}$, nor that the underlying space be ordered. We will therefore say that a sequence $(S_N)_{N \geqslant 1}$ of random variables taking values in a topological space $\mathcal{S}$ satisfies a *large deviation principle* with *rate function* $I : \mathcal{S} \to \mathbb{R}$ if, for every Borel set $A \subseteq \mathcal{S}$, we have

$$\limsup_{N \to +\infty} \frac{1}{N} \log \mathbb{P}\{S_N \in A\} \leqslant -\inf_{x \in \overline{A}} I(x), \tag{2.57}$$

and

$$\liminf_{N \to +\infty} \frac{1}{N} \log \mathbb{P}\{S_N \in A\} \geqslant -\inf_{x \in \mathrm{int}(A)} I(x). \tag{2.58}$$

Since $\mathrm{int}(A) \subseteq A \subseteq \overline{A}$, there is no loss in generality in restricting our attention to closed sets in the upper bound (2.57), and to open sets in the lower bound (2.58). Heuristically, if $(S_N)_{N \geqslant 1}$ satisfies a large deviation principle, then for every $x \in \mathcal{S}$,

$$\mathbb{P}\{S_N \simeq x\} \simeq \exp(-NI(x)), \tag{2.59}$$

so $S_N$ should concentrate around the minima of its rate function. This is consistent with the fact that the rate function (2.51) is minimized at the mean $p$. Although the theory of large deviations can be developed for rather general spaces $\mathcal{S}$, we will restrict our attention to the case $\mathcal{S} = \mathbb{R}$ here, and derive a general criterion involving log-Laplace transforms to determine when a sequence $(S_N)_{N \geqslant 1}$ satisfies a large deviation principle. Extending this result to the case of $\mathcal{S} = \mathbb{R}^d$ is an interesting exercise which we leave to the curious reader.

The log-Laplace transform of a random variable $X$ is the function $\psi : \mathbb{R} \to \mathbb{R} \cup \{+\infty\}$ defined by

$$\psi(\lambda) := \log \mathbb{E} \exp(\lambda X). \tag{2.60}$$

To motivate the general large deviation principle that we will prove, and see how the log-Laplace transform is related to it, consider a sequence $(X_n)_{n \geq 1}$ of bounded and centred i.i.d. random variables and denote by

$$S_N := \frac{1}{N} \sum_{n=1}^{N} X_n \qquad (2.61)$$

the sample average of $N$ of these random variables. Chebyshev's inequality implies that for any $x \in \mathbb{R}$ and $\lambda \geq 0$, we have

$$\mathbb{P}\{S_N \geq x\} \leq \exp(-\lambda N x) \mathbb{E} \exp(\lambda N S_N) = \exp\left(-N\left(\lambda x - \frac{1}{N}\psi_N(\lambda N)\right)\right), \qquad (2.62)$$

where $\psi_N$ denotes the log-Laplace transform of $S_N$. If we write $\psi = \psi_1$ for the log-Laplace transform of $X_1$, then the independence of the random variables $(X_n)_{n \geq 1}$ reveals that

$$\frac{1}{N}\psi_N(\lambda N) = \frac{1}{N}\log \mathbb{E}\exp(\lambda N S_N) = \log \mathbb{E}\exp(\lambda X_1) = \psi(\lambda). \qquad (2.63)$$

Substituting this into (2.62) and taking the infimum over all $\lambda \geq 0$ in the resulting bound shows that

$$\frac{1}{N}\log \mathbb{P}\{S_N \geq x\} \leq -\sup_{\lambda \geq 0}(\lambda x - \psi(\lambda)). \qquad (2.64)$$

This suggests that the sample mean $S_N$ might satisfy a large deviation principle with rate function given by the convex dual $\psi^*$. This is indeed correct and known as Cramér's theorem. It is a special case of the general large deviation principle that we now prove.

**Theorem 2.16** (General LDP on $\mathbb{R}$). *Let $(S_N)_{N \geq 1}$ be a sequence of real random variables, and for each integer $N \geq 1$ denote by $\psi_N : \mathbb{R} \to \mathbb{R}$ the log-Laplace transform of $S_N$. If there exists $\psi \in C^1(\mathbb{R}; \mathbb{R})$ such that for all $\lambda \in \mathbb{R}$, we have*

$$\lim_{N \to +\infty} \frac{1}{N}\psi_N(\lambda N) = \psi(\lambda), \qquad (2.65)$$

*then $(S_N)_{N \geq 1}$ satisfies a large deviation principle with rate function $\psi^*$.*

The knowledgeable reader will recognize in Theorem 2.16 a version of the Gärtner-Ellis theorem (see for instance Theorem 2.3.6 of [92]), with a slightly less general but significantly simpler assumption. Although we will not prove this here, we point out that Theorem 2.16 is also valid as stated for random variables taking values in $\mathbb{R}^d$.

*Proof of Theorem 2.16.* We start by deriving a handful of convenient properties of the convex dual

$$\psi^*(x) = \sup_{\lambda \in \mathbb{R}} (\lambda x - \psi(\lambda)). \tag{2.66}$$

Evaluating (2.66) at $\lambda = 0$, and using that $\psi(0) = 0$, we see that $\psi^*$ is non-negative,

$$\psi^* \geq 0. \tag{2.67}$$

We set $m := \psi'(0)$. Since $\psi$ is convex by Exercise 2.6 and the fact that a limit of convex functions is convex, we have that $\psi$ is bounded from below by the linear function $\lambda \mapsto m\lambda$ (this can be seen as a consequence of Theorem 2.13). It follows that

$$\psi^*(m) = 0 \quad \text{and} \quad \sup_{\lambda \leq 0}(\lambda x - \psi(\lambda)) \leq \sup_{\lambda \leq 0}(\lambda x - \lambda m) = 0 \tag{2.68}$$

for every $x \geq m$. In particular, for every $x \geq m$, the supremum in (2.66) can be restricted to $\mathbb{R}_{\geq 0}$,

$$\psi^*(x) = \sup_{\lambda \geq 0}(\lambda x - \psi(\lambda)). \tag{2.69}$$

We decompose the rest of the proof into four steps. First we derive an upper bound for the probability that $S_N$ is greater than some $x \geq m$, and we use this bound to establish the large deviation upper bound. We then obtain a lower bound for the probability that $S_N$ is an a small neighbourhood of some value $x \in \mathbb{R}$, and finally we use this lower bound to establish the large deviation lower bound.

*Step 1: upper bound on $\mathbb{P}\{S_N \geq x\}$.* In this step, we show that for every $x \geq m$,

$$\limsup_{N \to +\infty} \frac{1}{N} \log \mathbb{P}\{S_N \geq x\} \leq -\psi^*(x). \tag{2.70}$$

Given $\lambda \geq 0$, Chebyshev's inequality implies that

$$\mathbb{P}\{S_N \geq x\} \leq \exp(-\lambda N x) \mathbb{E} \exp(\lambda N S_N) = \exp\bigl(-\lambda N x + \psi_N(\lambda N)\bigr).$$

It follows by the assumption (2.65) that

$$\limsup_{N \to +\infty} \frac{1}{N} \log \mathbb{P}\{S_N \geq x\} \leq \psi(\lambda) - \lambda x.$$

Taking the supremum over $\lambda \geq 0$ and remembering (2.69) establishes (2.70).

*Step 2: large deviation upper bound.* We now show the large deviation upper bound; that is, we prove that for every closed set $A \subseteq \mathbb{R}$,

$$\limsup_{N \to +\infty} \frac{1}{N} \log \mathbb{P}\{S_N \in A\} \leq -\inf_{x \in A} \psi^*(x). \tag{2.71}$$

First, we observe that by the same reasoning as in Step 1, for every $x \leqslant m$,

$$\limsup_{N \to +\infty} \frac{1}{N} \log \mathbb{P}\{S_N \leqslant x\} \leqslant -\psi^*(x). \tag{2.72}$$

Given a closed set $A \subseteq \mathbb{R}$, we use the value $m = \psi'(0)$ defined before Step 1 to partition the probability that $S_N$ lies in $A$,

$$\mathbb{P}\{S_N \in A\} \leqslant \mathbb{P}\{S_N \in A \cap (-\infty, m]\} + \mathbb{P}\{S_N \in A \cap [m, +\infty)\}. \tag{2.73}$$

If $A \cap [m, +\infty)$ is not empty, then we denote by $x \geqslant m$ its infimum. This infimum belongs to $A$ since $A$ is closed, and

$$\mathbb{P}\{S_N \in A \cap [m, +\infty)\} \leqslant \mathbb{P}\{S_N \geqslant x\}.$$

It follows by (2.70) that

$$\limsup_{N \to +\infty} \frac{1}{N} \log \mathbb{P}\{S_N \in A \cap [m, +\infty)\} \leqslant -\psi^*(x) \leqslant -\inf_{x \in A} \psi^*(x). \tag{2.74}$$

The last inequality remains valid if $A \cap [m, +\infty)$ is empty provided that we understand that $\log 0 = -\infty$, and that if $A$ itself is empty, then the infimum over $A$ is $+\infty$. Similar observations based on (2.72) as opposed to (2.70) reveal that

$$\limsup_{N \to +\infty} \frac{1}{N} \log \mathbb{P}\{S_N \in A \cap (-\infty, m]\} \leqslant -\inf_{x \in A} \psi^*(x).$$

Combining this with (2.74) and (2.73), we obtain (2.71).

*Step 3: lower bound for* $\mathbb{P}\{S_N \simeq x\}$. Since $\psi^*$ is convex, non-negative and vanishes at $m$ by Exercise 2.10 and equations (2.67) and (2.68), it must be non-increasing on $(-\infty, m]$ and non-decreasing on $[m, +\infty)$. Moreover, as $\psi'$ is continuous, its image must be an interval of $\mathbb{R}$ whose endpoints we denote by $d_- \leqslant d_+ \in [-\infty, +\infty]$. In this step, we fix $x \in (d_-, d_+)$ and show that for every $\varepsilon > 0$,

$$\liminf_{N \to +\infty} \frac{1}{N} \log \mathbb{P}\{S_N \in [x - \varepsilon, x + \varepsilon]\} \geqslant -\psi^*(x). \tag{2.75}$$

Looking at the optimality condition for the supremum in (2.69), it is natural to define $\lambda \in \mathbb{R}$ to be such that $x = \psi'(\lambda)$. Such a $\lambda$ exists since $x$ belongs to $(d_-, d_+)$. For definiteness, we assume that $\lambda \geqslant 0$ and observe that

$$\mathbb{P}\{S_N \in [x - \varepsilon, x + \varepsilon]\} \geqslant \exp(-\lambda N(x + \varepsilon)) \mathbb{E} \exp(\lambda N S_N) \mathbf{1}_{\{S_N \in [x-\varepsilon, x+\varepsilon]\}}. \tag{2.76}$$

In order to conclude, we would like to remove the indicator function in the expectation above. Loosely speaking, we would like to verify that our choice of $\lambda$ ensures that

$$\frac{\mathbb{E} \exp(\lambda N S_N) \mathbf{1}_{\{S_N \notin [x-\varepsilon, x+\varepsilon]\}}}{\mathbb{E} \exp(\lambda N S_N)} \ll 1.$$

This is a question about a probability upper bound very much like the one studied in Step 1. For every $\mu \geq 0$, we have

$$\mathbb{E}\exp(\lambda NS_N)\mathbf{1}_{\{S_N \geq x+\varepsilon\}} \leq \exp(-N\mu(x+\varepsilon))\mathbb{E}\exp((\lambda+\mu)NS_N)$$
$$= \exp\bigl(-N\mu(x+\varepsilon) + \psi_N(N\lambda+N\mu)\bigr). \quad (2.77)$$

Similarly, for every $\mu \geq 0$, we have

$$\mathbb{E}\exp(\lambda NS_N)\mathbf{1}_{\{S_N \leq x-\varepsilon\}} \leq \exp(N\mu(x-\varepsilon))\mathbb{E}\exp((\lambda-\mu)NS_N)$$
$$= \exp\bigl(N\mu(x-\varepsilon) + \psi_N(N\lambda-N\mu)\bigr). \quad (2.78)$$

Since $\psi \in C^1(\mathbb{R};\mathbb{R})$, any $\mu > 0$ sufficiently small will be such that

$$|\psi(\lambda+\mu) - \psi(\lambda) - \mu\psi'(\lambda)| + |\psi(\lambda-\mu) - \psi(\lambda) + \mu\psi'(\lambda)| \leq \frac{\varepsilon\mu}{4}.$$

Fixing such a $\mu > 0$, recalling the convergence (2.65) and remembering that we chose $\psi'(\lambda) = x$, we obtain that for every $N$ sufficiently large,

$$|\psi_N(N\lambda+N\mu) - \psi_N(N\lambda) - N\mu x| + |\psi_N(N\lambda-N\mu) - \psi_N(N\lambda) + N\mu x| \leq \frac{N\varepsilon\mu}{2}.$$

Combining this with (2.77) and (2.78) reveals that for every $N$ sufficiently large,

$$\mathbb{E}\exp(\lambda NS_N)\mathbf{1}_{\{S_N \notin [x-\varepsilon,x+\varepsilon]\}} \leq 2\exp\left(-\frac{N\varepsilon\mu}{2}\right)\mathbb{E}\exp(\lambda NS_N),$$

and in particular,

$$\mathbb{E}\exp(\lambda NS_N)\mathbf{1}_{\{S_N \in [x-\varepsilon,x+\varepsilon]\}} \geq \frac{1}{2}\mathbb{E}\exp(\lambda NS_N).$$

Substituting this into (2.76), we get that, for every $N$ sufficiently large,

$$\mathbb{P}\{S_N \in [x-\varepsilon,x+\varepsilon]\} \geq \frac{1}{2}\exp\bigl(-\lambda N(x+\varepsilon)\bigr)\mathbb{E}\exp(\lambda NS_N).$$

Observing that the left side of (2.75) only gets smaller as $\varepsilon$ decreases, letting $N$ tend to infinity and then $\varepsilon$ tend to zero establishes (2.75).

*Step 4: large deviation lower bound.* Finally, we show the large deviation lower bound; that is, we prove that for every open set $A \subseteq \mathbb{R}$,

$$\limsup_{N\to+\infty} \frac{1}{N}\log\mathbb{P}\{S_N \in A\} \geq -\inf_{x\in A}\psi^*(x). \quad (2.79)$$

To start with, we extend (2.75) to every $x \in \mathbb{R}$. If $x \notin [d_-, d_+]$ this is immediate since in this case $\psi^*(x) = +\infty$, by a small variant of Exercise 2.12. Suppose that $d_+ < +\infty$ and that $x = d_+$. Since $d_- \leq d_+$ are the endpoints of the interval spanned by $\psi'$, and since $\psi'(0) = m$, we must have that $m \leq d_+$. Since $\psi^*$ is non-decreasing on $[m, +\infty)$, it follows that
$$\psi^*(d_+) \geq \lim_{\substack{x \to d_+ \\ x < d_+}} \psi^*(x),$$
so again (2.75) holds. The same argument also applies to the case when $d_- < +\infty$ and $x = d_-$. We now pick an open set $A \subseteq \mathbb{R}$. For every $x \in A$, we can find $\varepsilon > 0$ sufficiently small that $[x-\varepsilon, x+\varepsilon] \subseteq A$. We deduce from (2.75) that
$$\liminf_{N \to +\infty} \frac{1}{N} \log \mathbb{P}\{S_N \in A\} \geq -\psi^*(x).$$
Taking the supremum over $x \in A$ establishes (2.79) and completes the proof. ∎

**Theorem 2.17** (Cramér). *If $(X_n)_{n \geq 1}$ is a sequence of i.i.d. random variables with finite log-Laplace transform $\psi : \mathbb{R} \to \mathbb{R}$, then the sequence $(S_N)_{N \geq 1}$ of sample averages (2.61) satisfies a large deviation principle with rate function $\psi^*$.*

*Proof.* Using that $\psi$ takes only finite values, one can check that $\psi$ must be continuously differentiable. The claim is then an immediate consequence of the general large deviation result in Theorem 2.16 and the equality (2.63) which relates the log-Laplace transform of $S_N$ to that of $X_1$. ∎

To finish our discussion of large deviation principles, we explore the necessity of assuming that the rate function in Theorem 2.16 be continuously differentiable. This possibly surprising assumption will also appear in the Hamilton-Jacobi approach when we prove the convex selection principle in Section 3.6.

**Example 2.18.** In this example we show that the assumption $\psi \in C^1(\mathbb{R}; \mathbb{R})$ in Theorem 2.16 is necessary. Consider a sequence $(S_N)_{N \geq 1}$ of Rademacher random variables,
$$\mathbb{P}\{S_N = 1\} = \mathbb{P}\{S_N = -1\} = \frac{1}{2},$$
and denote by $\psi_N$ the log-Laplace transform of $S_N$ (which here does not actually depend on $N$). A direct computation shows that $\psi_N(\lambda) = \log \cosh(\lambda)$, and by Exercise 1.1,
$$\psi(\lambda) := \lim_{N \to +\infty} \frac{1}{N} \psi_N(\lambda N) = \lim_{N \to +\infty} \frac{1}{N} \log\left(\frac{e^{N\lambda}}{2} + \frac{e^{-N\lambda}}{2}\right) = |\lambda|.$$

Since $\psi$ is Lipschitz with Lipschitz constant equal to one, Exercise 2.12 implies that
$$\psi^*(x) = \begin{cases} 0 & x \in [-1, 1], \\ +\infty & \text{otherwise.} \end{cases}$$

If Theorem 2.16 were true, then for all $N$ large enough, we would have

$$\frac{1}{N}\log \mathbb{P}\{S_N \in (-1/2, 1/2)\} \geq - \inf_{x \in (-1/2, 1/2)} \psi^*(x) - \frac{1}{2} = -\frac{1}{2},$$

and therefore

$$0 = \mathbb{P}\{S_N \in (-1/2, 1/2)\} \geq \exp(-N/2).$$

This contradiction establishes the necessity of the assumption $\psi \in C^1(\mathbb{R}; \mathbb{R})$.

We briefly discuss further the phenomenon at play in this example, without going into details — the reader can safely skip this imprecise paragraph. Suppose that a sequence $(S_N)_{N \geq 1}$ of random variables with log-Laplace transform $\psi_N$ satisfies a large deviation principle with rate function $I$ that may or may not be convex. Informally, this means that (2.59) holds. In Example 2.18 we would need to have $I(x) = 0$ for $x \in \{-1, +1\}$ and $I(x) = +\infty$ otherwise. In general, when (2.59) holds, we should expect that

$$\lim_{N \to +\infty} \frac{1}{N} \psi_N(\lambda N) = I^*(\lambda), \qquad (2.80)$$

and Theorem 2.16 states that if $\psi = I^*$ is $C^1$, then $I^{**}$ is the correct rate function. Using the Fenchel-Moreau theorem and Exercise 2.10, one can show that in general, the convex bi-dual $I^{**}$ is the convex envelope, or largest convex minorant, of the function $I$. The assumption in Theorem 2.16 that the function $\psi$ is continuously differentiable can be seen to be equivalent to the statement that the convex dual $\psi^*$ has no flat piece, a flat piece being an open interval over which $\psi^*$ is affine (see also Exercise 2.23 for a related statement). This means that whenever (2.59) holds, the assumption that $\psi$ is continuously differentiable implies that the convex envelope of $I$ has no flat piece. This can only occur if $I$ was convex in the first place, and then indeed $I^{**} = I$ by the Fenchel-Moreau theorem.

For a more general and in-depth discussion of large deviation principles, we refer the interested reader to [92, 93, 258].

**Exercise 2.21.** Show that (2.54) and (2.55) are equivalent to the string of inequalities (2.56).

**Exercise 2.22.** Let $X$ be a Bernoulli random variable with parameter $p \in [0, 1]$, and denote by $\psi$ its log-Laplace transform. Show that the convex dual of $\psi$ is the rate function (2.51),

$$\psi^*(x) = x \log\left(\frac{x}{p}\right) + (1-x) \log\left(\frac{1-x}{1-p}\right). \qquad (2.81)$$

## 2.3 Analyzing the Curie-Weiss model

The large deviation principle in Theorem 2.16 allows us to compute the limit free energy in the generalized Curie-Weiss model introduced in Section 1.3. We recall that this model is defined in terms of a probability measure $P_N$ on $\mathbb{R}^N$ with the property that $|\sigma| \leq \sqrt{N}$ for $P_N$-a.e. sample $\sigma$, as well as a smooth function $\xi \in C^\infty(\mathbb{R};\mathbb{R})$ on the real line. For every $\sigma \in \mathbb{R}^N$, we denote by

$$S_N(\sigma) := \frac{1}{N} \sum_{i=1}^{N} \sigma_i \tag{2.82}$$

its sample average, or magnetization. For every $\sigma$ in the support of $P_N$, we have that

$$|S_N(\sigma)|^2 = \left|\frac{1}{N}\sum_{i=1}^N \sigma_i\right|^2 \leq \frac{1}{N}\sum_{i=1}^N \sigma_i^2 = \frac{|\sigma|^2}{N} \leq 1. \tag{2.83}$$

This means that the magnetization lies in the interval $[-1,1]$,

$$S_N(\sigma) \in [-1,1]. \tag{2.84}$$

The free energy in the generalized Curie-Weiss model is given by

$$F_N(t,h) := \frac{1}{N} \log \int_{\mathbb{R}^N} \exp N\bigl(t\xi(S_N(\sigma)) + hS_N(\sigma)\bigr)\, dP_N(\sigma). \tag{2.85}$$

To see how a large deviation principle for the sequence $(S_N)_{N\geq 1}$ can shed light on the asymptotic behaviour of the free energy (2.85), we denote by $\psi_N$ the log-Laplace transform of $S_N$ and suppose that the limit

$$\psi(h) := \lim_{N\to+\infty} \frac{1}{N} \psi_N(hN) = \lim_{N\to+\infty} F_N(0,h) \tag{2.86}$$

exists and defines a continuously differentiable function $\psi \in C^1(\mathbb{R};\mathbb{R})$. The general large deviation principle in Theorem 2.16 then implies that the sequence $(S_N)_{N\geq 1}$ satisfies a large deviation principle with rate function $\psi^*$. Heuristically, this means that for every $m \in [-1,1]$,

$$\mathbb{P}\{S_N \simeq m\} \simeq \exp(-N\psi^*(m)). \tag{2.87}$$

If we fix an integer $K \geq 1$ and denote by $-1 = m_0 < m_1 < \ldots < m_K = 1$ a partition of the interval $[-1,1]$ into sub-intervals of width $2/K$, then we can localize the free energy according to the approximate magnetization of its spin configurations,

$$F_N(t,h) \simeq \frac{1}{N} \log \sum_{k=0}^{K} \int_{\mathbb{R}^N} \mathbf{1}_{\{S_N(\sigma)\in[m_k,m_{k+1})\}} \exp N\bigl(t\xi(m_k) + hm_k\bigr)\, dP_N(\sigma). \tag{2.88}$$

Here $m_{K+1}$ can be any value larger than one, say $m_{K+1} = +\infty$. The large deviation principle (2.87) implies that for $K$ large enough and $1 \leqslant k \leqslant K$,

$$\int_{\mathbb{R}^N} \mathbf{1}_{\{S_N(\sigma) \in [m_k, m_{k+1})\}} \, dP_N(\sigma) \simeq \exp(-N\psi^*(m_k)). \tag{2.89}$$

It follows by Exercise 1.1 that for large enough values of $K$,

$$F_N(t, h) \simeq \max_{0 \leqslant k \leqslant K} \left( \xi(m_k) + hm_k - \psi^*(m_k) \right). \tag{2.90}$$

We should therefore be able to leverage a large deviation principle for the magnetization (2.82) to show that the limit free energy (2.85) in the generalized Curie-Weiss model is given by

$$\lim_{N \to +\infty} F_N(t, h) = \sup_{m \in [-1, 1]} \left( \xi(m) + hm - \psi^*(m) \right). \tag{2.91}$$

By choosing $\xi(x) = x^2$ and $P_N$ to be uniform on $\{\pm 1\}^N := \{-1, +1\}^N$, we will be able to deduce a formula for the limit of the free energy (1.35) in the Curie-Weiss model.

**Theorem 2.19.** *For each integer $N \geqslant 1$, we denote by $F_N : \mathbb{R}_{\geqslant 0} \times \mathbb{R} \to \mathbb{R}$ the free energy (2.85) in the generalized Curie-Weiss model, and suppose that for every $h \in \mathbb{R}$ the limit*

$$\psi(h) := \lim_{N \to +\infty} F_N(0, h) \tag{2.92}$$

*exists. If $\psi \in C^1(\mathbb{R}; \mathbb{R})$, then the limit free energy $f : \mathbb{R}_{\geqslant 0} \times \mathbb{R} \to \mathbb{R}$ in the generalized Curie-Weiss model is given by*

$$f(t, h) := \lim_{N \to +\infty} F_N(t, h) = \sup_{m \in [-1, 1]} \left( t\xi(m) + hm - \psi^*(m) \right). \tag{2.93}$$

*Proof.* For each integer $N \geqslant 1$, we denote by $\psi_N$ the log-Laplace transform of the magnetization $S_N$ defined in (2.82). The function $\psi \in C^1(\mathbb{R}; \mathbb{R})$ can be written as

$$\psi(h) = \lim_{N \to +\infty} F_N(0, h) = \lim_{N \to +\infty} \frac{1}{N} \psi_N(hN), \tag{2.94}$$

so the general large deviation principle in Theorem 2.16 implies that the sequence $(S_N)_{N \geqslant 1}$ satisfies a large deviation principle with rate function $\psi^*$. This means that for every Borel set $A \subseteq \mathbb{R}$, we have

$$\limsup_{N \to +\infty} \frac{1}{N} \log P_N\{S_N \in A\} \leqslant \inf_{x \in \overline{A}} \psi^*(x) \tag{2.95}$$

and

$$\liminf_{N \to +\infty} \frac{1}{N} \log P_N\{S_N \in A\} \geqslant - \inf_{x \in \text{int}(A)} \psi^*(x). \tag{2.96}$$

The rest of the proof proceeds in two steps. First we establish the upper bound in (2.93), and then we obtain the matching lower bound.

*Step 1: the upper bound.* We fix an integer $K \geq 1$, and for each $k \in \{-K, \ldots, K-1\}$, we introduce the set
$$A_k := \left[\frac{k}{K}, \frac{k+1}{K}\right].$$

The bound (2.84) implies that
$$S_N(\sigma) \in [-1, 1] = \bigcup_{k=-K}^{K-1} A_k$$

for every $\sigma$ in the support of $P_N$. It follows by monotonicity of the logarithm that

$$\begin{aligned}
F_N(t,h) &\leq \frac{1}{N} \log \sum_{k=-K}^{K-1} \int_{\mathbb{R}^N} \mathbf{1}_{\{S_N \in A_k\}} \exp N(t\xi(S_N) + hS_N) \, dP_N \\
&\leq \frac{1}{N} \log \sum_{k=-K}^{K-1} \max_{m \in A_k} \left(\exp N(t\xi(m) + hm)\right) P_N\{S_N \in A_k\} \\
&\leq \frac{\log(2K)}{N} + \max_{-K \leq k < K} \left(\max_{m \in A_k}(t\xi(m) + hm) + \frac{1}{N} \log P_N\{S_N \in A_k\}\right).
\end{aligned}$$

Using the large deviation upper bound (2.95) to let $N$ tend to infinity reveals that

$$\limsup_{N \to +\infty} F_N(t,h) \leq \max_{-K \leq k < K} \left(\max_{m \in A_k}(t\xi(m) + hm) - \inf_{m \in A_k} \psi^*(m)\right). \tag{2.97}$$

We next observe that since $\xi \in C^\infty(\mathbb{R}; \mathbb{R})$, there exists a constant $C < +\infty$ such that for every $m, m' \in [-1, 1]$,
$$|\xi(m) - \xi(m')| \leq C|m - m'|.$$

For each $k \in \{-K, \ldots, K-1\}$, and $m \in A_k$, we therefore have that
$$\max_{m' \in A_k}\left(t\xi(m') + hm'\right) - \psi^*(m) \leq t\xi(m) + hm - \psi^*(m) + \frac{C+h}{K}.$$

Taking the supremum over $m \in A_k$ and over $k \in \{-K, \ldots, K-1\}$, and then combining this with (2.97), we obtain that
$$\limsup_{N \to +\infty} F_N(t,h) \leq \sup_{m \in [-1,1]} \left(t\xi(m) + hm - \psi^*(m)\right) + \frac{C+h}{K}.$$

Letting $K$ tend to infinity then yields the upper bound
$$\limsup_{N \to +\infty} F_N(t,h) \leq \sup_{m \in [-1,1]} \left(t\xi(m) + hm - \psi^*(m)\right). \tag{2.98}$$

*Step 2: the lower bound.* We fix $\varepsilon > 0$, $m_0 \in [-1,1]$, $\delta > 0$, and consider the open set $U_\delta := (m_0 - \delta, m_0 + \delta)$. The monotonicity of the logarithm implies that

$$F_N(t,h) \geq \frac{1}{N} \log \int_{\mathbb{R}^N} \mathbf{1}_{\{S_N \in U_\delta\}} \exp N\big(t\xi(S_N) + hS_N\big) \, dP_N$$

$$\geq \inf_{m \in U_\delta} \big(t\xi(m) + hm\big) + \frac{1}{N} \log P_N\{S_N \in U_\delta\}.$$

Using the large deviation lower bound (2.96) to let $N$ tend to infinity reveals that

$$\liminf_{N \to +\infty} F_N(t,h) \geq \inf_{m \in U_\delta} \big(t\xi(m) + hm\big) - \psi^*(m_0).$$

Letting $\delta$ tend to zero and leveraging the continuity of the map $m \mapsto t\xi(m) + hm$ gives the lower bound

$$\liminf_{N \to +\infty} F_N(t,h) \geq t\xi(m_0) + hm_0 - \psi^*(m_0).$$

Since $m_0 \in [-1,1]$ was arbitrary, this completes the proof. ∎

**Corollary 2.20.** *The limit of the free energy* (1.35) *in the (standard) Curie-Weiss model is given by*

$$f(t,h) := \lim_{N \to +\infty} F_N(t,h) = \sup_{m \in [-1,1]} \big(tm^2 + hm - \psi^*(m)\big) \qquad (2.99)$$

*for the function* $\psi^* : \mathbb{R} \to \mathbb{R}$ *defined by*

$$\psi^*(m) := \frac{1+m}{2} \log(1+m) + \frac{1-m}{2} \log(1-m), \quad \text{for } m \in [-1,1], \qquad (2.100)$$

*with the understanding that* $0 \log(0) = 0$, *and* $\psi^* = +\infty$ *on* $\mathbb{R} \setminus [-1,1]$.

*Proof.* The Curie-Weiss model corresponds to the generalized Curie-Weiss model with $\xi(x) = x^2$ and $P_N$ uniform on the set $\{\pm 1\}^N$. This means that for every $h \in \mathbb{R}$,

$$F_N(0,h) = \frac{1}{N} \log \frac{1}{2^N} \sum_{\sigma \in \{\pm 1\}^N} \exp\left(h \sum_{i=1}^N \sigma_i\right) = \log \cosh(h).$$

In particular, the limit

$$\psi(h) = \lim_{N \to +\infty} F_N(0,h) = \log \cosh(h)$$

exists and defines a continuously differentiable function $\psi \in C^1(\mathbb{R}; \mathbb{R})$. It follows by Theorem 2.19 that the limit free energy is

$$f(t,h) = \sup_{m \in [-1,1]} \big(tm^2 + hm - \psi^*(m)\big)$$

for the convex dual

$$\psi^*(m) = \sup_{h\in\mathbb{R}}\left(hm - \psi(h)\right) = \sup_{h\in\mathbb{R}}\left(hm - \log\cosh(h)\right).$$

The symmetry of $\psi$ implies that of $\psi^*$, so it suffices to compute the convex dual for $m \in [0,1]$. The derivative of the function $g(h) := hm - \psi(h)$ is $g'(h) = m - \tanh(h)$. Since $\tanh(h) \in (-1,1)$, if $m = 1$ the function $g$ is increasing and

$$\psi^*(1) = \lim_{h\to+\infty}(h - \log\cosh(h)) = \log(2).$$

On the other hand, if $m \in [0,1)$, then $g$ admits a unique critical point $h^*$ which satisfies

$$\tanh(h^*) = m.$$

Since $g$ is concave, we must have $\psi^*(m) = g(h^*)$. A direct computation shows that

$$e^{2h^*} = \frac{1-m}{1+m}$$

from which it is readily verified that

$$\psi^*(m) := \frac{1+m}{2}\log(1+m) + \frac{1-m}{2}\log(1-m). \tag{2.101}$$

(The attentive reader will have observed that this must coincide with the rate function (2.51), up to a change of variables.) This completes the proof. ∎

## 2.4 The envelope theorem and the Curie-Weiss magnetization

The original motivation for introducing the free energy (1.35) was to understand the asymptotic behaviour of the mean magnetization

$$m_N(t,h) := \left\langle \frac{1}{N}\sum_{i=1}^N \sigma_i \right\rangle \tag{2.102}$$

in the Curie-Weiss model, and to show that for large enough values of the inverse temperature parameter $t$, it exhibits ferromagnetic behaviour. Recall that $\langle\cdot\rangle$ denotes the average under the Gibbs measure (1.31), see (1.33). A direct computation reveals that the mean magnetization is the derivative of the free energy,

$$m_N(t,h) = \partial_h F_N(t,h). \tag{2.103}$$

The sequence $(F_N)_{N\geq 1}$ is convex in $h$ by Exercise 2.6 and, by Corollary 2.20, it converges to the limit free energy

$$f(t,h) := \sup_{m\in[-1,1]}\left(tm^2 + hm - \psi^*(m)\right), \tag{2.104}$$

where $\psi^* : \mathbb{R} \to \mathbb{R}$ is given by

$$\psi^*(m) := \frac{1+m}{2}\log(1+m) + \frac{1-m}{2}\log(1-m). \qquad (2.105)$$

We would like to pass to the limit in (2.103), and to study the properties of $\partial_h f$. This should hopefully allow us to identify the phase transition of the model from paramagnetic to ferromagnetic behaviour as we increase $t$, similarly to (v) of Exercise 1.3.

The problem of differentiating a function which, like the function $f$ in (2.104), is defined as the envelope, or supremum, of a family of functions, is sufficiently important that we will discuss it in general. Given a continuous function $g \in C(\mathbb{R}^d \times \mathbb{R}^d; \mathbb{R})$ that is continuously differentiable in its first variable, we will show that the function

$$f(x) := \sup_{y \in \mathbb{R}^d} g(x,y) \qquad (2.106)$$

is differentiable at a point $x \in \mathbb{R}^d$ if and only if the function $y \mapsto \nabla_x g(x,y)$ is constant on the set of optimizers

$$\mathcal{O}_x := \{y \in \mathbb{R}^d \mid f(x) = g(x,y)\}. \qquad (2.107)$$

This holds in particular whenever $\mathcal{O}_x$ is a singleton. In this case, the derivative of $f$ can be computed by simply bringing the derivative into the supremum.

**Theorem 2.21** (Envelope). *Let $g \in C(\mathbb{R}^d \times \mathbb{R}^d; \mathbb{R})$ be a continuous function that is continuously differentiable in its first variable, let $f : \mathbb{R}^d \to \mathbb{R}$ be its envelope (2.106), and let $\mathcal{O}_x$ be as in (2.107). We fix $x \in \mathbb{R}^d$, and assume that there exists a compact set $K$ which contains $\mathcal{O}_{x'} \neq \emptyset$ for every $x'$ in a small enough neighbourhood of $x$. As $x' \in \mathbb{R}^d$ tends to $x$, we have the expansion*

$$f(x') = f(x) + \sup_{y \in \mathcal{O}_x} \{(x'-x) \cdot \nabla_x g(x,y)\} + o(|x - x'|). \qquad (2.108)$$

*In particular, the envelope function $f : \mathbb{R}^d \to \mathbb{R}$ is differentiable at $x \in \mathbb{R}^d$ if and only if the set*

$$\mathcal{D}_x := \{\nabla_x g(x,y) \mid y \in \mathcal{O}_x\} \qquad (2.109)$$

*is a singleton. In this case, for any $y \in \mathcal{O}_x$, we have $\nabla f(x) = \nabla_x g(x,y)$.*

*Proof.* We decompose the proof into three steps.

*Step 1: lower bound for* (2.108). Since $g$ is continuously differentiable, we have for every $y \in \mathcal{O}_x$ that, as $x'$ tends to $x$,

$$g(x',y) = g(x,y) + (x'-x) \cdot \nabla_x g(x,y) + o(|x - x'|). \qquad (2.110)$$

Moreover, since $\mathcal{O}_x$ is contained in the compact set $K$, the expansion (2.110) holds uniformly over $y \in \mathcal{O}_x$. By the definition of $f$ as a supremum and the fact that $y \in \mathcal{O}_x$, we obtain that, as $x'$ tends to $x$,

$$f(x') \geqslant f(x) + (x' - x) \cdot \nabla_x g(x, y) + o(|x - x'|). \tag{2.111}$$

Taking the supremum over $y \in \mathcal{O}_x$ yields the lower bound for (2.108).

*Step 2: upper bound for* (2.108). We now turn to the converse bound. Let $(x_n)_{n \geqslant 1} \subseteq \mathbb{R}^d \smallsetminus \{x\}$ be a sequence converging to $x$. Our goal is to show that

$$\limsup_{n \to +\infty} \frac{1}{|x_n - x|} \left( f(x_n) - f(x) - \sup_{y \in \mathcal{O}_x} \{(x_n - x) \cdot \nabla_x g(x, y)\} \right) \leqslant 0. \tag{2.112}$$

For every $n \geqslant 1$ sufficiently large, we can pick $y_n \in \mathcal{O}_{x_n}$. Since $\mathcal{O}_{x_n} \subseteq K$ for every $n$ sufficiently large, we can, after the extraction of a subsequence, assume that $(y_n)_{n \geqslant 1}$ converges to some point $\bar{y} \in \mathbb{R}^d$. In fact, it suffices that we show (2.112) for every sequence $(x_n, y_n)_{n \geqslant 1}$ with $x_n \neq x$ converging to $x$, with $y_n \in \mathcal{O}_{x_n}$, and with $y_n$ converging to some $\bar{y} \in \mathbb{R}^d$. This can be seen by arguing by contradiction: if (2.112) fails to hold for some sequence $(\bar{x}_n)_{n \geqslant 1}$, then we can construct another sequence $(x_n, y_n)_{n \geqslant 1}$ with $y_n \in \mathcal{O}_{x_n}$, with (2.112) still failing, and with $y_n$ converging to some $\bar{y}$. So we assume that these conditions hold from now on. For every $y' \in \mathbb{R}^d$, we have

$$g(x_n, y') \leqslant \sup_{y \in \mathbb{R}^d} g(x_n, y) = g(x_n, y_n),$$

and letting $n$ tend to infinity, we obtain that $g(x, y') \leqslant g(x, \bar{y})$ — in other words, we have that $\bar{y} \in \mathcal{O}_x$. By the definition of $f$ as a supremum and the fundamental theorem of calculus, we have

$$g(x_n, y_n) - f(x) \leqslant g(x_n, y_n) - g(x, y_n) = \int_0^1 (x_n - x) \cdot \nabla_x g(t x_n + (1-t) x, y_n) \, dt.$$

Since $y_n \in \mathcal{O}_{x_n}$ and $\bar{y} \in \mathcal{O}_x$, we obtain that

$$f(x_n) \leqslant f(x) + \sup_{y \in \mathcal{O}_x} \{(x_n - x) \cdot \nabla_x g(x, y)\}$$

$$+ |x_n - x| \left| \int_0^1 \nabla_x g(t x_n + (1-t) x, y_n) \, dt - \nabla_x g(x, \bar{y}) \right|. \tag{2.113}$$

Since $\nabla_x g$ is continuous, we have

$$\lim_{n \to +\infty} \int_0^1 \nabla_x g(t x_n + (1-t) x, y_n) \, dt = \nabla_x g(x, \bar{y}),$$

so letting $n$ tend to infinity in (2.113) yields (2.112).

*Step 3: differentiability condition.* By (2.108), it is clear that if the set $\mathcal{D}_x$ in (2.109) is a singleton, then $f$ is differentiable at $x$, and $\nabla f(x) = \nabla_x g(x,y)$ for every $y \in \mathcal{O}_x$. Conversely, suppose that the function $f$ is differentiable at $x$. Using (2.111) and the differentiability assumption, we obtain that, for every $y \in \mathcal{O}_x$ and as $x'$ tends to $x$,

$$(x'-x) \cdot \big(\nabla f(x) - \nabla_x g(x,y)\big) \geq o(|x'-x|).$$

This means that $\nabla_x g(x,y) = \nabla f(x)$, so $\mathcal{D}_x$ is a singleton as announced. ∎

To apply this result to determine the limit of the mean magnetization (2.102) as $N$ tends to infinity, given $t \in \mathbb{R}_{\geq 0}$ and $h \in \mathbb{R}$, we need to study the set of maximizers of the function $g_{t,h} : [-1,1] \to \mathbb{R}$ defined by

$$g_{t,h}(m) := tm^2 + hm - \psi^*(m). \tag{2.114}$$

We will focus in particular on the case when $h$ is close to zero. Figure 2.1 displays the graph of this function for $(t,h) = (0.1,0)$ and $(t,h) = (0.6,0)$, and it suggests that there should be a transition at some critical inverse temperature $t_c \in (0.1, 0.6)$, where the number of maximizers jumps from one to two.

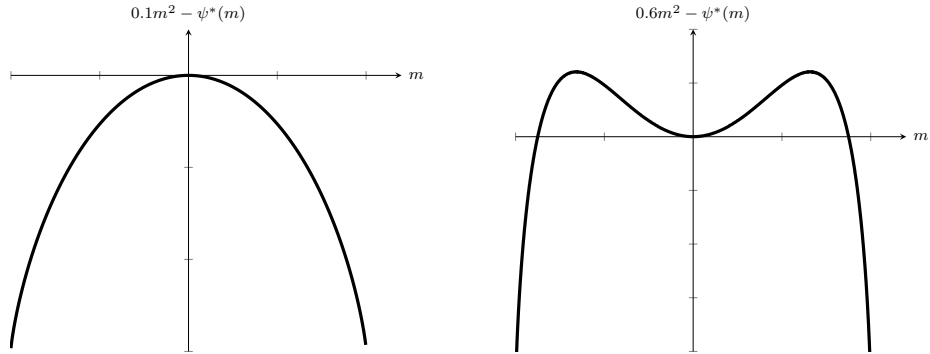

**Figure 2.1** Graphs of $m \mapsto tm^2 - \psi^*(m)$ for $t = 0.1$ and $t = 0.6$ in the Curie-Weiss model.

**Proposition 2.22.** *If $t \leq \frac{1}{2}$, then for every $h \in \mathbb{R}$, the function (2.114) has a unique maximizer $m_h(t) \in [-1,1]$. On the other hand, if $t > \frac{1}{2}$, then for all $|h|$ small enough, the function (2.114) has exactly two local maximizers $m_h^-(t) \leq m_h^+(t) \in [-1,1]$. In all cases, these local maximizers are solutions to the fixed point equation*

$$m = \tanh(h + 2tm). \tag{2.115}$$

*For $t > \frac{1}{2}$, there exists $m^*(t) > 0$ (not depending on $h$) such that for all $|h|$ small enough,*

$$m_h^-(t) < -m^*(t) < m^*(t) < m_h^+(t). \tag{2.116}$$

*Proof.* Fix $t \in \mathbb{R}_{\geq 0}$ and $h \in \mathbb{R}$. A direct computation gives that for every $m \in [-1, 1]$,

$$\partial_m g_{t,h}(m) = 2tm + h + \frac{1}{2}\log\left(\frac{1-m}{1+m}\right).$$

It follows that the critical points $m$ of $g_{t,h}$ are the solutions to the critical point equation

$$m = \frac{e^{2(h+2tm)} - 1}{e^{2(h+2tm)} + 1} = \tanh(h + 2tm).$$

Differentiating one more time, wee see that, for every $m \in [-1, 1]$,

$$\partial_m^2 g_{t,h}(m) = 2t - \frac{1}{1-m^2}.$$

We now distinguish two cases. First, if $t \leq \frac{1}{2}$, then $\partial_m^2 g_{t,h}(m) \leq 0$, with the inequality being strict for every non-zero $m \in [-1, 1]$. This means that $\partial_m g_{t,h}$ is a decreasing function that tends to $+\infty$ as $m$ approaches $-1$ and tends to $-\infty$ as $m$ approaches $1$. It follows by the intermediate value theorem that $g_{t,h}$ has a unique critical point $m_h(t)$, and the sign of $\partial_m g_{t,h}$ on each side of this critical point reveals that it is a maximizer. In the case $t > \frac{1}{2}$, the function $\partial_m^2 g_{t,h}$ has exactly two zeros,

$$\pm m^*(t) := \pm\sqrt{1 - \frac{1}{2t}},$$

and it is positive on the interval $(-m^*(t), m^*(t))$ and negative on its complement. Since $\partial_m g_{t,0}(0) = 0$, the continuity of $h \mapsto \partial_m g_{t,h}$ implies that for $|h|$ small enough, the function $\partial_m g_{t,h}$ must have a zero on the interval $(-m^*(t), m^*(t))$. Since $\partial_m g_{t,h}$ tends to $+\infty$ as $m$ approaches $-1$, decreases for $m \leq -m^*(t)$, increases on $(-m^*(t), m^*(t))$ and has a zero in this interval before decreasing to $-\infty$ as $m$ approaches $1$, it must have two zeros $m_h^\pm(t) \in [-1, 1]$ outside the interval $(-m^*(t), m^*(t))$. Analyzing the sign of $\partial_m g_{t,h}$ reveals that these are the only local maximizers of $g_{t,h}$; the critical point in the interval $(-m^*(t), m^*(t))$ is a local minimizer. This completes the proof. ∎

If $t > \frac{1}{2}$ and $h > 0$, then $g_{t,h}(m_h^+(t)) > g_{t,h}(m_h^-(t))$, so the envelope theorem and Proposition 2.22 imply that $m(t, h) = \partial_h f(t, h) = m_h^+(t)$. The sequence $(m_h^+(t))_{h>0}$ is uniformly bounded by one, and by (2.115) any of its subsequential limits as $h \searrow 0$ must satisfy the fixed point equation $m = \tanh(2tm)$. Together with the separation bound (2.116), this means that it must be $m_0^+(t)$. That is,

$$\lim_{h \searrow 0} m(t, h) = m_0^+(t) > 0. \tag{2.117}$$

A similar argument shows that

$$\lim_{h \nearrow 0} m(t, h) = m_0^-(t) = -m_0^+(t) < 0, \tag{2.118}$$

where the last equality uses the symmetry of the fixed point equation (2.115) when $h = 0$. The qualitative behaviour of the mapping $h \mapsto m(t,h)$ is therefore similar to that in the two-dimensional Ising model depicted in Figure 1.1. The phase transition at the critical temperature $t_c = \frac{1}{2}$ is explored further in Exercise 2.24. By varying the choice of the function $\xi$ and the sequence of measures $P_N$ in the generalized Curie-Weiss model, one can devise models with a wide variety of behaviours, such as multiple phase transitions [107].

**Exercise 2.23.** Let $f : \mathbb{R}^d \to \mathbb{R}$ be a function such that

$$\lim_{|x| \to +\infty} \frac{f(x)}{|x|} = +\infty. \tag{2.119}$$

We also assume that the function $f$ is strictly convex, that is, for every $x \neq y \in \mathbb{R}^d$ and $\alpha \in (0,1)$, we have

$$f(\alpha x + (1-\alpha)y) < \alpha f(x) + (1-\alpha) f(y). \tag{2.120}$$

Show that the effective domain of $f^*$ is $\mathbb{R}^d$ and that $f^*$ is continuously differentiable.

**Exercise 2.24.** Denote by $f : \mathbb{R}_{\geq 0} \times \mathbb{R} \to \mathbb{R}$ the limit free energy (2.104) in the Curie-Weiss model and by $t_c := \frac{1}{2}$ its critical inverse temperature parameter. We recall from Proposition 2.22 that if $t \leq t_c$, then the function (2.114) for $h = 0$ has a unique maximizer $m_0(t) = 0$, while it has two maximizers $\pm m_0(t)$ if $t > t_c$.

(i) Identify $\partial_t f(t,0)$, and show that it is a continuous function of $t$.

(ii) Identify an exponent $\delta \geq 0$ such that, as $h$ tends to zero and in a sense to be specified,
$$f(t_c, h) \simeq |h|^{1+1/\delta}.$$

(iii) Identify an exponent $\beta \geq 0$ such that, as $t$ tends to $t_c$ and in a sense to be specified,
$$m_0(t) \simeq (t - t_c)_+^\beta,$$
where $x_+ := \max(0,x)$ denotes the positive part.

The notation used here for the exponents $\delta$ and $\beta$ follows the convention of the physics literature. In dimensions $d \geq 4$, the Ising model on $\mathbb{Z}^d$ has the same critical exponents as the Curie-Weiss model, with a logarithmic correction when $d = 4$ [10, 12]. For the two-dimensional Ising model, we have $\delta = 15$ and $\beta = \frac{1}{8}$ [150, 151, 204, 247, 266]. There is no known or conjectured closed-form description of these exponents in three dimensions; sophisticated numerical estimates are derived in [111, 112, 157].

# Chapter 3
# Hamilton-Jacobi equations

In Chapter 2, we computed the limit of the free energy in the Curie-Weiss model via large deviation principles. Unfortunately, this approach seems inapplicable for the models coming from statistical inference or spin glass theory which we are ultimately interested in. In this chapter, we introduce the Hamilton-Jacobi approach, and we use it to compute the limit of the free energy in the Curie-Weiss model and its generalization. The approach centres around the observation that the finite-volume free energy satisfies an approximate Hamilton-Jacobi equation, up to a small error term. Section 3.1 presents this computation in the context of the Curie-Weiss and generalized Curie-Weiss models, and justifies the need to introduce a precise notion of solution to a Hamilton-Jacobi equation. The notion of a viscosity solution is defined in Section 3.2, and we observe that any subsequential limit of the free energy in the Curie-Weiss model must be a viscosity solution. We then show in Section 3.3 that there can be at most one viscosity solution to a Hamilton-Jacobi equation with a prescribed initial condition. This already guarantees that we have identified the limit free energy of the Curie-Weiss model uniquely. In Section 3.4, in order to recover the explicit formula for the limit free energy found in the previous chapter, we develop variational representations of viscosity solutions whenever the initial condition or the non-linearity in the equation is convex (or concave). In Section 3.5, we explore whether any aspect of this variational structure is preserved in the absence of convexity or concavity. Finally, in Section 3.6, we focus on a new difficulty that emerges when trying to analyze the generalized Curie-Weiss model. To overcome this, we present a convenient criterion for deciding when a given convex function is the viscosity solution to a Hamilton-Jacobi equation.

## 3.1 A Hamilton-Jacobi approach to Curie-Weiss

Recall from Section 1.3 that, given a smooth function $\xi \in C^\infty(\mathbb{R}; \mathbb{R})$ and a probability measure $P_N$ on $\mathbb{R}^N$ with the property that $|\sigma| \leq \sqrt{N}$ for $P_N$-a.e. $\sigma$, the Hamiltonian of the generalized Curie-Weiss model at the point $(t, h, \sigma) \in \mathbb{R}_{\geq 0} \times \mathbb{R} \times \mathbb{R}^N$ is defined by

$$H_N(t, h, \sigma) := Nt\xi\left(\frac{1}{N}\sum_{i=1}^N \sigma_i\right) + h\sum_{i=1}^N \sigma_i. \tag{3.1}$$

The free energy in the generalized Curie-Weiss model is given by

$$F_N(t, h) := \frac{1}{N}\log\int_{\mathbb{R}^N} \exp H_N(t, h, \sigma)\,dP_N(\sigma), \tag{3.2}$$

and the average of any bounded and measurable function $f : \mathbb{R}^N \to \mathbb{R}$ with respect to the Gibbs measure (1.31) is denoted by

$$\langle f(\sigma)\rangle := \frac{\int_{\mathbb{R}^N} f(\sigma)\exp H_N(t, h, \sigma)\,dP_N(\sigma)}{\int_{\mathbb{R}^N} \exp H_N(t, h, \sigma)\,dP_N(\sigma)}. \tag{3.3}$$

We recall that, although this is kept implicit in the notation, the bracket $\langle \cdot \rangle$ depends on the choice of parameters $t$ and $h$. The (standard) Curie-Weiss model corresponds to the choice $\xi(x) = x^2$ and $P_N$ uniform on the hypercube $\{\pm 1\}^N := \{-1, +1\}^N$. To illustrate the Hamilton-Jacobi approach in the simplest possible setting, let us first focus on this situation.

In the context of the Curie-Weiss model, a direct computation shows that

$$\partial_t F_N(t, h) = \left\langle \left(\frac{1}{N}\sum_{i=1}^N \sigma_i\right)^2 \right\rangle \quad \text{and} \quad \partial_h F_N(t, h) = \left\langle \frac{1}{N}\sum_{i=1}^N \sigma_i \right\rangle. \tag{3.4}$$

It follows that

$$\partial_t F_N(t, h) - \left(\partial_h F_N(t, h)\right)^2 = \left\langle\left(\frac{1}{N}\sum_{i=1}^N \sigma_i\right)^2\right\rangle - \left\langle\frac{1}{N}\sum_{i=1}^N \sigma_i\right\rangle^2 = \mathrm{Var}\left(\frac{1}{N}\sum_{i=1}^N \sigma_i\right). \tag{3.5}$$

Another direct computation reveals that

$$\partial_h^2 F_N(t, h) = \left\langle\frac{1}{N}\left(\sum_{i=1}^N \sigma_i\right)^2\right\rangle - \frac{1}{N}\left\langle\sum_{i=1}^N \sigma_i\right\rangle^2 = N\,\mathrm{Var}\left(\frac{1}{N}\sum_{i=1}^N \sigma_i\right). \tag{3.6}$$

Together with (3.5), this implies that the finite-volume free energy (3.2) in the Curie-Weiss model satisfies

$$\partial_t F_N(t, h) - \left(\partial_h F_N(t, h)\right)^2 = \frac{1}{N}\partial_h^2 F_N(t, h). \tag{3.7}$$

Since $P_N$ is uniform on the hypercube $\{\pm 1\}^N$, the free energy at "time" 0 is independent of $N$. Indeed,

$$F_N(0, h) = \frac{1}{N}\log\frac{1}{2^N}\sum_{\sigma\in\{\pm 1\}^N}\exp\left(h\sum_{i=1}^N \sigma_i\right) = \frac{1}{N}\log\left(\frac{1}{2}(e^h + e^{-h})\right)^N = \log\cosh(h).$$

We introduce the notation

$$\psi(h) := \log\cosh(h) = F_1(0,h) = F_N(0,h). \tag{3.8}$$

These observations suggest that the limit free energy $f : \mathbb{R}_{\geqslant 0} \times \mathbb{R} \to \mathbb{R}$ in the Curie-Weiss model should satisfy the Hamilton-Jacobi equation

$$\partial_t f - (\partial_h f)^2 = 0 \quad \text{on} \quad \mathbb{R}_{>0} \times \mathbb{R}, \tag{3.9}$$

subject to the initial condition $f(0,\cdot) = \psi$. In Corollary 2.20, we showed that

$$f(t,h) = \sup_{m \in [-1,1]} \left( tm^2 + hm - \psi^*(m) \right). \tag{3.10}$$

By the envelope theorem (Theorem 2.21), one can show that at every point $(t,h)$ of differentiability of $f$, we have

$$\partial_t f(t,h) = m_0^2(t,h) \quad \text{and} \quad \partial_h f(t,h) = m_0(t,h), \tag{3.11}$$

for a maximizer $m_0(t,h)$ of the right side of (3.10). In particular, the limit free energy $f$ does indeed satisfy the Hamilton-Jacobi equation (3.9) at all its points of differentiability. Together with Rademacher's theorem (Theorem 2.10), this implies that $f$ satisfies the equation (3.9) almost everywhere. Notice that $f$ is Lipschitz continuous because each $F_N$ is, by (3.4). A natural question arises at this point. Could we identify the limit free energy $f$ as the *unique* Lipschitz continuous function with $f(0,\cdot) = \psi$ which satisfies the equation (3.9) almost everywhere? Unfortunately, the answer is *no*. In fact, the following construction shows that this Hamilton-Jacobi equation admits infinitely many almost-everywhere solutions which are Lipschitz continuous.

**Example 3.1.** For simplicity, we will construct infinitely many Lipschitz functions that satisfy the Hamilton-Jacobi equation (3.9) almost everywhere and are constant equal to zero at the initial time. Similar constructions can be performed for more general initial conditions. Temporarily disregarding the question of the initial condition, we start by noticing that the functions $(t,h) \mapsto 0$, $(t,h) \mapsto t + h$ and $(t,h) \mapsto t - h$ are all solutions to (3.9). It thus follows that the Lipschitz function

$$\widetilde{f}(t,h) := \begin{cases} t + h & \text{if } h \in [-t,0] \\ t - h & \text{if } h \in [0,t] \\ 0 & \text{otherwise} \end{cases} \tag{3.12}$$

displayed in Figure 3.1 is an almost everywhere solution to (3.9), as it is obtained by "gluing together" these different solutions, and we can disregard the measure-zero set of points where they are joined together. We also have $\widetilde{f}(0,\cdot) = 0$.

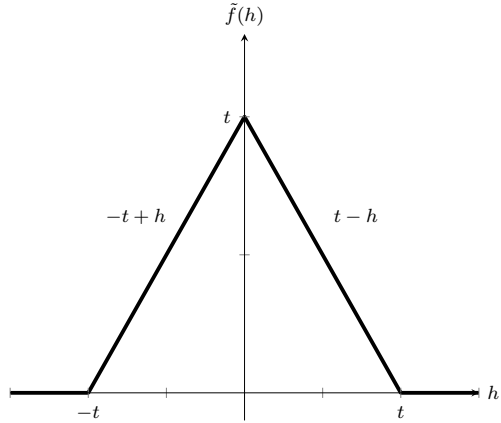

**Figure 3.1** Graph of the function $h \mapsto \widetilde{f}(t,h)$, for a fixed value of $t > 0$.

Evidently, the null solution also satisfies these properties, so we already have two almost-everywhere solutions with the same initial condition. Moreover, any translation in space of $\widetilde{f}$ also satisfies this property; and we can also delay the "emergence of the corner" to some arbitrary time. So in fact there are uncountably many Lipschitz functions that solve the equation (3.9) almost everywhere and vanish at the initial time.

We could try to impose uniqueness of solutions by strengthening the regularity assumptions; for instance, we could impose that a solution $f$ to (3.9) be in $C^1(\mathbb{R}_{\geqslant 0} \times \mathbb{R}; \mathbb{R})$ and satisfy the equation everywhere. The problem with this idea is that in this case, the set of solutions can be empty, and indeed, we have already seen in (2.117) and (2.118) that our candidate solution $f$ in (3.10) is not everywhere differentiable. We therefore need to identify a notion of solution that is more stringent than the "almost-everywhere solutions" explored above, but less stringent than asking the solution to be in $C^1(\mathbb{R}_{\geqslant 0} \times \mathbb{R}; \mathbb{R})$. In a nutshell, we will also require the function $f$ to satisfy a certain form of the maximum principle. We observe that whenever two functions $f$ and $g$ satisfy the Hamilton-Jacobi equation (3.10) with "viscosity" parameter $\varepsilon > 0$,

$$\partial_t f - (\partial_h f)^2 = \varepsilon \Delta f \quad \text{on} \quad \mathbb{R}_{>0} \times \mathbb{R},$$

if their initial conditions are ordered, say $f(0,\cdot) \leqslant g(0,\cdot)$, then this ordering is preserved at all later times $t \geqslant 0$ as well, $f(t,\cdot) \leqslant g(t,\cdot)$. While we will obtain this monotonicity property as a consequence of the definition of solution explained below, we point out that one can essentially also go the other way around [16].

Before introducing this appropriate notion of solution, let us see what the computations leading to (3.7) would look like in the context of the generalized

Curie-Weiss model. In this setting, one can readily see that

$$\partial_t F_N(t,h) = \left\langle \xi\left(\frac{1}{N}\sum_{i=1}^N \sigma_i\right)\right\rangle \quad \text{and} \quad \partial_h F_N(t,h) = \left\langle \frac{1}{N}\sum_{i=1}^N \sigma_i\right\rangle, \qquad (3.13)$$

and that (3.6) still holds. It follows that

$$\partial_t F_N(t,h) - \xi\left(\partial_h F_N(t,h)\right) = \left\langle \xi\left(\frac{1}{N}\sum_{i=1}^N \sigma_i\right)\right\rangle - \xi\left(\left\langle \frac{1}{N}\sum_{i=1}^N \sigma_i\right\rangle\right). \qquad (3.14)$$

If the variance of the magnetization is small, we therefore expect $\partial_t F_N - \xi(\partial_h F_N)$ to be close to zero. Recall that for every $\sigma$ in the support of $P_N$, the magnetization is bounded by one,

$$\left|\frac{1}{N}\sum_{i=1}^N \sigma_i\right| \leq \left(\frac{1}{N}\sum_{i=1}^N \sigma_i^2\right)^{\frac{1}{2}} = \frac{|\sigma|}{\sqrt{N}} \leq 1, \qquad (3.15)$$

and so this bound also holds with probability one under the Gibbs measure. Moreover, as $\xi$ is smooth, there exists a constant $C < +\infty$ such that for every $x, y \in [-1, 1]$,

$$|\xi(y) - \xi(x) - (y-x)\xi'(x)| \leq C(x-y)^2.$$

This implies that

$$\left|\xi\left(\frac{1}{N}\sum_{i=1}^N \sigma_i\right) - \xi\left(\left\langle \frac{1}{N}\sum_{i=1}^N \sigma_i\right\rangle\right) - \xi'\left(\left\langle \frac{1}{N}\sum_{i=1}^N \sigma_i\right\rangle\right)\left(\frac{1}{N}\sum_{i=1}^N \sigma_i - \left\langle \frac{1}{N}\sum_{i=1}^N \sigma_i\right\rangle\right)\right|$$

$$\leq C\left|\frac{1}{N}\sum_{i=1}^N \sigma_i - \left\langle \frac{1}{N}\sum_{i=1}^N \sigma_i\right\rangle\right|^2. \qquad (3.16)$$

Averaging with respect to the Gibbs measure (3.3), substituting the resulting bound into (3.14), and recalling (3.6) yields that

$$\left|\partial_t F_N(t,h) - \xi\left(\partial_h F_N(t,h)\right)\right| \leq C\mathrm{Var}\left(\frac{1}{N}\sum_{i=1}^N \sigma_i\right) = \frac{C}{N}\partial_h^2 F_N(t,h). \qquad (3.17)$$

If we assume that the limit

$$\psi(h) := \lim_{N \to +\infty} F_N(0,h) \qquad (3.18)$$

of the initial conditions exists, which is for instance the case when $P_N$ is a product measure, then we expect the limit free energy $f : \mathbb{R}_{\geq 0} \times \mathbb{R} \to \mathbb{R}$ in the generalized Curie-Weiss model to satisfy the Hamilton-Jacobi equation

$$\partial_t f - \xi(\partial_h f) = 0 \quad \text{on} \quad \mathbb{R}_{>0} \times \mathbb{R} \qquad (3.19)$$

subject to the initial condition $f(0,\cdot) = \psi$. This will however be harder to prove rigorously than for the Curie-Weiss model due to the presence of the absolute value in the inequality (3.17). The problem is that functions that are only known to satisfy the inequality (3.17) do not necessarily satisfy the maximum principle, while the notion of viscosity solutions takes this property as its foundation. We will explain how to circumvent this difficulty by leveraging the convexity of the solutions.

## 3.2  Viscosity solutions to Hamilton-Jacobi equations

In this section, given a locally Lipschitz continuous non-linearity $\mathsf{H} : \mathbb{R}^d \to \mathbb{R}$ and a Lipschitz continuous initial condition $\psi : \mathbb{R}^d \to \mathbb{R}$, we define a notion of solution $f = f(t,x) : \mathbb{R}_{\geq 0} \times \mathbb{R}^d \to \mathbb{R}$ to the Hamilton-Jacobi equation

$$\partial_t f - \mathsf{H}(\nabla f) = 0 \quad \text{on} \quad \mathbb{R}_{>0} \times \mathbb{R}^d, \tag{3.20}$$

subject to the initial condition $f(0, \cdot) = \psi$. We understand the gradient to be taken with respect to the $x$ variable, leaving the $t$ variable aside: $\nabla f(x) = (\partial_{x_1} f, \ldots, \partial_{x_d} f)$. The equations arising in the context of the Curie-Weiss and generalized Curie-Weiss models correspond to the one-dimensional case, $d = 1$ and, respectively, the non-linearities $\mathsf{H}(p) = p^2$ and $\mathsf{H}(p) = \xi(p)$.

A natural way to define a solution to the equation (3.20) is to add a small "viscosity" parameter $\varepsilon > 0$, and to consider the second-order parabolic equation

$$\partial_t f_\varepsilon - \mathsf{H}(\nabla f_\varepsilon) = \varepsilon \Delta f_\varepsilon \quad \text{on} \quad \mathbb{R}_{>0} \times \mathbb{R}^d \tag{3.21}$$

subject to the initial condition $f_\varepsilon(0, \cdot) = \psi$. That smooth solutions exist for (3.21) can be shown by classical techniques, because the Laplacian term is dominant on very small scales, being of higher order than the non-linear term in the equation. The solution to the Hamilton-Jacobi equation (3.20) can then be defined as the limit of the solutions to (3.21) as the viscosity parameter $\varepsilon$ tends to zero. Although we will not pursue this route rigorously here, let us suppose for a moment that for each $\varepsilon > 0$, we have been able to define a smooth solution $f_\varepsilon : \mathbb{R}_{\geq 0} \times \mathbb{R}^d \to \mathbb{R}$ to (3.21), and that we have been able to show that, as $\varepsilon$ tends to zero, the sequence $(f_\varepsilon)_{\varepsilon > 0}$ converges to some function $f : \mathbb{R}_{\geq 0} \times \mathbb{R}^d \to \mathbb{R}$ in the topology of local uniform convergence.

We would like to say, in some sense to be discovered, that the limit thus obtained is a solution to (3.20). As was seen in the previous section, the main difficulty is that we do not want to impose a solution to (3.19) to be differentiable everywhere; but we do not obtain unique solutions if we simply ignore the points of non-differentiability. In analogy with the notion of weak solutions, we would like to introduce smooth test functions and somehow move the derivatives of $f$ onto the test functions, so that we could introduce some constraint for what a solution is allowed to do at points of non-differentiability. This transfer of the derivatives onto the test functions will not be obtained by some integration by parts here. Rather, we seek a strategy that accords well with the fact that our approximations $(f_\varepsilon)_{\varepsilon > 0}$ satisfy the maximum principle. We pick a test function $\phi \in C^\infty(\mathbb{R}_{>0} \times \mathbb{R}^d; \mathbb{R})$, and assume that $f - \phi$ achieves a strict local maximum at the point $(t^*, x^*) \in \mathbb{R}_{>0} \times \mathbb{R}^d$. If $f$ is smooth at $(t^*, x^*)$, then we clearly have

$$\bigl(\partial_t \phi(t^*, x^*), \nabla \phi(t^*, x^*)\bigr) = \bigl(\partial_t f(t^*, x^*), \nabla f(t^*, x^*)\bigr), \tag{3.22}$$

so we expect that
$$\partial_t \phi(t^*, x^*) - \mathsf{H}(\nabla \phi(t^*, x^*)) = 0. \qquad (3.23)$$

This identity may no longer be valid when $f$ is not differentiable at $(t^*, x^*)$, but we will now argue that the quantity on the left side of (3.23) must always be non-positive.

Using Exercise 3.1, we can find a sequence $(t_\varepsilon^*, x_\varepsilon^*)_{\varepsilon > 0}$ converging to $(t^*, x^*)$ with the property that, for every $\varepsilon > 0$ sufficiently small, the function $f_\varepsilon - \phi$ achieves a local maximum at $(t_\varepsilon^*, x_\varepsilon^*) \in \mathbb{R}_{>0} \times \mathbb{R}^d$, and therefore

$$\partial_t (f_\varepsilon - \phi)(t_\varepsilon^*, x_\varepsilon^*) = 0, \ \nabla(f_\varepsilon - \phi)(t_\varepsilon^*, x_\varepsilon^*) = 0, \text{ and } \Delta(f_\varepsilon - \phi)(t_\varepsilon^*, x_\varepsilon^*) \leq 0. \qquad (3.24)$$

It follows by (3.21) that

$$\left(\partial_t \phi - \mathsf{H}(\nabla \phi)\right)(t_\varepsilon^*, x_\varepsilon^*) = \left(\partial_t f - \mathsf{H}(\nabla f)\right)(t_\varepsilon^*, x_\varepsilon^*) = \varepsilon \Delta f_\varepsilon(t_\varepsilon^*, x_\varepsilon^*) \leq \varepsilon \Delta \phi(t_\varepsilon^*, x_\varepsilon^*).$$

Letting $\varepsilon$ tend to zero and leveraging the smoothness of $\phi$ shows that

$$\left(\partial_t \phi - \mathsf{H}(\nabla \phi)\right)(t^*, x^*) \leq 0. \qquad (3.25)$$

We have thus shown that whenever a smooth function $\phi$ is such that $f - \phi$ has a strict local maximum at $(t^*, x^*) \in \mathbb{R}_{>0} \times \mathbb{R}^d$, the inequality (3.25) holds at $(t^*, x^*)$. An analogous argument shows that whenever a smooth function $\phi$ is such that $f - \phi$ has a strict local minimum at $(t^*, x^*) \in \mathbb{R}_{>0} \times \mathbb{R}^d$, we have

$$\left(\partial_t \phi - \mathsf{H}(\nabla \phi)\right)(t^*, x^*) \geq 0. \qquad (3.26)$$

As we will see in the next section, these properties of $f$ we just derived are sufficient to determine it uniquely, once the initial condition $f(0, \cdot)$ is fixed. We will therefore take them as the definition of being a solution to the Hamilton-Jacobi equation (3.20).

When $f - \phi$ has a local maximum at $(t^*, x^*)$, we will often say that "$\phi$ touches $f$ from above at $(t^*, x^*)$". The reason is that, for our purposes, adding a constant to the test function $\phi$ is irrelevant to the discussion, so we may as well assume indeed that $(f - \phi)(t^*, x^*) = 0$. For this reason, we may also call the point $(t^*, x^*)$ the "contact point". We feel that the wording has some intuitive appeal, as we can imagine taking some arbitrary smooth function $\phi$ that is way above $f$, and then progressively "sliding it down" until the graphs of $\phi$ and $f$ touch — at the contact point. This is illustrated in Figure 3.2. Notice that not every point can be "touched" in this way. For instance, there is no smooth function that touches the absolute value function from above at the origin.

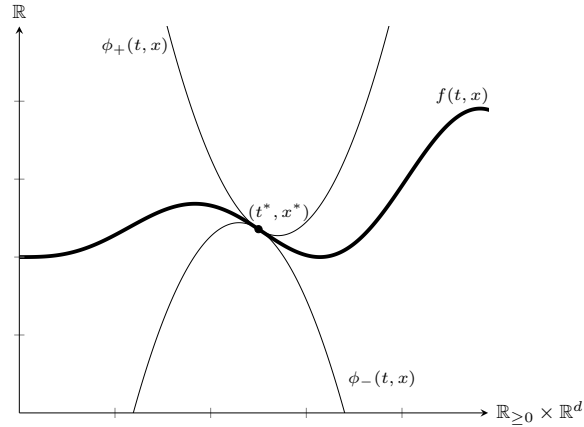

**Figure 3.2** The function $\phi_+$ touches the function $f$ from above at the point $(t^*, x^*)$ while the function $\phi_-$ touches the function $f$ from below at the point $(t^*, x^*)$.

**Definition 3.2.** A continuous function $u : \mathbb{R}_{\geq 0} \times \mathbb{R}^d \to \mathbb{R}$ is a *viscosity subsolution* to the Hamilton-Jacobi equation (3.20) if, for every $(t^*, x^*) \in \mathbb{R}_{>0} \times \mathbb{R}^d$ and $\phi \in C^\infty(\mathbb{R}_{>0} \times \mathbb{R}^d; \mathbb{R})$ with the property that $u - \phi$ has a local maximum at the point $(t^*, x^*) \in \mathbb{R}_{>0} \times \mathbb{R}^d$, we have

$$\left(\partial_t \phi - \mathsf{H}(\nabla \phi)\right)(t^*, x^*) \leq 0. \tag{3.27}$$

A continuous function $v : \mathbb{R}_{\geq 0} \times \mathbb{R}^d \to \mathbb{R}$ is a *viscosity supersolution* to the Hamilton-Jacobi equation (3.20) if, for every $(t^*, x^*) \in \mathbb{R}_{>0} \times \mathbb{R}^d$ and $\phi \in C^\infty(\mathbb{R}_{>0} \times \mathbb{R}^d; \mathbb{R})$ with the property that $v - \phi$ has a local minimum at the point $(t^*, x^*) \in \mathbb{R}_{>0} \times \mathbb{R}^d$, we have

$$\left(\partial_t \phi - \mathsf{H}(\nabla \phi)\right)(t^*, x^*) \geq 0. \tag{3.28}$$

A continuous function $f : \mathbb{R}_{\geq 0} \times \mathbb{R}^d \to \mathbb{R}$ is a *viscosity solution* to the Hamilton-Jacobi equation (3.20) if it is both a viscosity subsolution and a viscosity supersolution to (3.20).

As shown in Exercises 3.2–3.4, in the definition of a viscosity subsolution or viscosity supersolution, replacing "local maximum" by "strict local maximum" or by "global maximum", or replacing the requirement that $\phi \in C^\infty(\mathbb{R}_{>0} \times \mathbb{R}^d; \mathbb{R})$ by the requirement that $\phi \in C^1(\mathbb{R}_{>0} \times \mathbb{R}^d; \mathbb{R})$ lead to equivalent definitions. In Exercise 3.5, we verify that a $C^1$ function that satisfies the equation (3.20) everywhere is indeed a viscosity solution, so the notion of viscosity solution is more "permissive" than prescribing the function to be $C^1(\mathbb{R}_{\geq 0} \times \mathbb{R}^d; \mathbb{R})$ and to solve the equation everywhere. One can also show that a viscosity solution must satisfy the equation (3.20) at every point of differentiability, see for instance Theorem 10.1.2.1 in [115]. By the Rademacher theorem (Theorem 2.10), a Lipschitz viscosity solution must therefore satisfy the equation (3.20) almost everywhere. In other words, the notion of viscosity solution is indeed more stringent than that of "almost-everywhere solution" explored in Section 3.1. This is corroborated by the fact that the solution constructed in Example 3.1 is not a viscosity solution to the Hamilton-Jacobi equation (3.9).

**Example 3.3.** Recall the definition of the function $\widetilde{f}$ in Example 3.1, fix $t_0 > 0$, and consider the smooth function $\phi(t,x) = t$. Observe that $\widetilde{f} - \phi \leq 0$, and $f(t_0, 0) - \phi(t_0, 0) = 0$, so $\widetilde{f} - \phi$ has a local maximum at $(t_0, 0)$. If $\widetilde{f}$ were a viscosity subsolution to the Hamilton-Jacobi equation (3.9), we would have $(\partial_t \phi - (\partial_x \phi)^2)(t_0, 0) \leq 0$. However, we see that $(\partial_t \phi - (\partial_x \phi)^2)(t_0, 0) = 1$. This shows that $\widetilde{f}$ is not a viscosity solution to (3.9).

We now verify that any subsequential limit of the sequence $(F_N)_{N \geq 1}$ of free energies in the Curie-Weiss model is a viscosity solution to the Hamilton-Jacobi equation (3.9). Notice first that the sequence of free energies $(F_N)_{N \geq 1}$ is precompact in the topology of local uniform convergence, by the Arzelà-Ascoli theorem. Indeed, the function $F_N(0, \cdot)$ does not depend on $N$ by (3.8), and $\partial_h F_N$ and $\partial_t F_N$ take values in the compact interval $[-1, 1]$ by (3.4).

**Proposition 3.4.** *For each $N \geq 1$, let $F_N$ denote the free energy (1.35) of the (standard) Curie-Weiss model, and let $f : \mathbb{R}_{\geq 0} \times \mathbb{R} \to \mathbb{R}$ be a subsequential limit of $(F_N)_{N \geq 1}$ in the topology of local uniform convergence. The function $f$ is a viscosity solution to the Hamilton-Jacobi equation (3.9) with initial condition $f(0, h) = \psi(h) := \log \cosh(h)$.*

*Proof.* This is essentially the argument that led us to (3.25). Recall from (3.7) that

$$\partial_t F_N - (\partial_h F_N)^2 = \frac{1}{N} \partial_h^2 F_N \quad \text{on} \quad \mathbb{R}_{\geq 0} \times \mathbb{R}. \tag{3.29}$$

Fix a smooth function $\phi \in C^\infty(\mathbb{R}_{>0} \times \mathbb{R}; \mathbb{R})$ with the property that $f - \phi$ has a strict local maximum at $(t^*, h^*) \in \mathbb{R}_{>0} \times \mathbb{R}$. Since $(F_N)_{N \geq 1}$ converges locally uniformly to $f$, using Exercise 3.1 it is possible to find a sequence $(t_N, h_N)_{N \geq 1} \subseteq \mathbb{R}_{>0} \times \mathbb{R}$ converging to $(t^*, h^*)$ with the property that, for every $N$ large enough, $F_N - \phi$ has a local maximum at $(t_N, h_N)$. Since $t^* > 0$, we have that $t_N > 0$ for every $N$ sufficiently large, so the derivatives in $t$ and $h$ of $F_N - \phi$ must vanish at $(t_N, h_N)$, and we must also have that $\partial_h^2 (F_N - \phi)(t_N, h_N) \leq 0$. Using also the identity (3.29), we obtain

$$\left(\partial_t \phi - (\partial_h \phi)^2\right)(t_N, h_N) = \frac{1}{N} \partial_h^2 F_N(t_N, h_N) \leq \frac{1}{N} \partial_h^2 \phi(t_N, h_N).$$

Since $\phi$ is smooth, letting $N$ tend to infinity shows that $f$ is a viscosity subsolution to the Hamilton-Jacobi equation (3.9). A similar argument shows that it is also a viscosity supersolution to this equation. By (3.8), the initial condition is satisfied. This completes the proof. ∎

**Exercise 3.1.** Let $(f_N)_{N \geqslant 1}$ be a sequence of continuous functions on $\mathbb{R}^d$ converging locally uniformly to a function $f : \mathbb{R}^d \to \mathbb{R}$. If $f$ has a strict local maximum at $x \in \mathbb{R}^d$, show that there exists $r > 0$ and a sequence of points $(x_N)_{N \geqslant 1} \subseteq \mathbb{R}^d$ converging to $x$ such that, for every $N$ sufficiently large, we have $f_N(x_N) = \sup_{B_r(x)} f_N$.

**Exercise 3.2.** Show that in the definition of a viscosity subsolution, replacing "local maximum" by "strict local maximum" yields an equivalent definition.

**Exercise 3.3.** Show that in the definition of a Lipschitz viscosity subsolution, replacing "local maximum" by "global maximum" or "strict global maximum" yields an equivalent definition.

**Exercise 3.4.** Show that in the definition of a viscosity subsolution, replacing the regularity of test functions "$\phi \in C^\infty(\mathbb{R}_{>0} \times \mathbb{R}^d; \mathbb{R})$" by "$\phi \in C^1(\mathbb{R}_{>0} \times \mathbb{R}^d; \mathbb{R})$" yields an equivalent definition.

**Exercise 3.5.** Show that if $f \in C^1(\mathbb{R}_{\geqslant 0} \times \mathbb{R}^d; \mathbb{R})$ satisfies (3.20) everywhere, then $f$ is a viscosity solution.

## 3.3 Uniqueness of solutions via the comparison principle

The result in Proposition 3.4 ensures that any subsequential limit of the free energy in the Curie-Weiss model is a viscosity solution to (3.9). In particular, this shows that a viscosity solution to this equation with initial condition $\psi$ exists. To show that the finite-volume free energies indeed converge to *the* solution to (3.9), we need to assert the uniqueness of solutions to Hamilton-Jacobi equations with a prescribed initial condition. We will in fact show a more general result known as the comparison principle. The comparison principle formalizes the idea that the Hamilton-Jacobi equation (3.20) should preserve the ordering of initial conditions. In fact, since the equation is invariant under the addition of a constant to the solution, we will show that if $u$ is a subsolution and $v$ is a supersolution to (3.20), then the function $u - v$ achieves its supremum at time zero. In particular, if $u(0, \cdot) \leqslant v(0, \cdot)$, then this ordering is preserved by the evolution.

If we assume for a moment that $u$ and $v$ are smooth functions, then we can argue heuristically to get a sense of why this result might be true. For smooth functions $u$ and $v$, saying that $u$ and $v$ are a subsolution and a supersolution to the Hamilton-Jacobi equation (3.20) amounts to saying that

$$\partial_t u - \mathsf{H}(\nabla u) \leqslant 0 \quad \text{and} \quad \partial_t v - \mathsf{H}(\nabla v) \geqslant 0. \tag{3.30}$$

Arguing by contradiction, suppose that

$$\sup_{\mathbb{R}_{\geqslant 0} \times \mathbb{R}^d} (u - v) > \sup_{\{0\} \times \mathbb{R}^d} (u - v). \tag{3.31}$$

Up to subtracting a small increasing function of $t$ to $u$, such as $\varepsilon t$ for some sufficiently small $\varepsilon > 0$, we can make sure that (3.31) is still valid, and also that we have improved (3.30) into

$$\partial_t u - \mathsf{H}(\nabla u) < 0 \quad \text{and} \quad \partial_t v - \mathsf{H}(\nabla v) \geq 0. \tag{3.32}$$

Now, if we assume that the supremum on the left side of (3.31) is achieved at some point $(t^*, x^*)$, then we must have $t^* > 0$, so the first order derivatives of $u$ and $v$ must coincide at this point. But this contradicts (3.32), so (3.31) cannot be true. Of course, there is much left to be desired with this argument, since we assumed that $u$ and $v$ are smooth, and also that the supremum on the left side of (3.31) is achieved at some point. We start by discussing how one can go around the second problem.

To simplify the discussion, we will temporarily assume that the space $\mathbb{R}^d$ is replaced by the unit torus $\mathbb{T}^d = \mathbb{R}^d/\mathbb{Z}^d$, so that the variable $x$ lives in a compact space without boundary. We take $u, v : \mathbb{R}_{\geq 0} \times \mathbb{T}^d \to \mathbb{R}$ to be a viscosity subsolution and supersolution to (3.20) respectively, and we assume that they are smooth. Our starting point is again to argue by contradiction, assuming that there exists some time $T > 0$ with

$$\sup_{[0,T] \times \mathbb{T}^d} (u - v) > \sup_{\{0\} \times \mathbb{T}^d} (u - v). \tag{3.33}$$

Denoting $\chi(t) := \frac{\varepsilon}{T-t}$, we can select $\varepsilon > 0$ sufficiently small that

$$\sup_{[0,T] \times \mathbb{T}^d} (u - v - \chi) > \sup_{\{0\} \times \mathbb{T}^d} (u - v - \chi). \tag{3.34}$$

Since $u - v$ is uniformly bounded over $[0, T] \times \mathbb{T}^d$, it is clear that approximate optimizers of the left side of (3.34) will remain away from the final time $T$. Using also that $[0, T] \times \mathbb{T}^d$ is compact, we can construct an optimizer $(t^*, x^*) \in [0, T] \times \mathbb{T}^d$ for this supremum, and it is clear that we must in fact have $t^* \in (0, T)$. Once this is verified, we can use the differential condition at the maximum to ascertain that

$$\partial_t(u-v)(t^*, x^*) - \frac{\varepsilon}{(T-t^*)^2} = \partial_t(u - v - \chi)(t^*, x^*) = 0 \tag{3.35}$$

and

$$\nabla(u-v)(t^*, x^*) = 0. \tag{3.36}$$

Since $v$ is a smooth supersolution to the equation, we can combine this with the second inequality in (3.30) to obtain that

$$\big(\partial_t u - \mathsf{H}(\nabla u)\big)(t^*, x^*) - \frac{\varepsilon}{(T-t^*)^2} = \big(\partial_t v - \mathsf{H}(\nabla v)\big)(t^*, x^*) \geq 0. \tag{3.37}$$

This contradicts the assumption (3.30) that $u$ is a smooth subsolution to the Hamilton-Jacobi equation (3.20). The role of the perturbation function $\chi$ is therefore two-fold.

First, it ensures that the optimum $t^*$ is detached from the right endpoint of the interval $[0,T]$, that is, $t^* < T$. The function $\chi$ also allows us to strengthen the inequalities (3.30) into those in (3.32), if we understand that the function $u$ is redefined to be $u - \chi$.

In order to establish the comparison principle rigorously, we need to resolve two main problems. The first, and most fundamental, is of course that $u$ and $v$ cannot be assumed to be differentiable everywhere. The second is that the spatial variable takes values in $\mathbb{R}^d$ rather than the torus. To tackle the first of these problems, we will "double the variables", and rather optimize a function that involves $u(t,x) - v(t',x')$ plus a smooth penalty term that strongly encourages $(t,x)$ and $(t',x')$ to stay close together. This will naturally provide us with smooth test functions that touch $u$ and $v$ from above and below respectively. The second problem, that the variable $x$ lives in an unbounded space, will be tackled by introducing another "cutoff" function, similar to the function $\chi$ used above, but in the space variable. This additional spatial cutoff will in fact allow us to prove a somewhat stronger result than the comparison principle. To state this concisely, we denote the Lipschitz semi-norm of a function $h : \mathbb{R}^d \to \mathbb{R}$ by

$$\|h\|_{\mathrm{Lip}} = \sup_{x \neq x'} \frac{|h(x') - h(x)|}{|x' - x|}, \tag{3.38}$$

and introduce the space of uniformly Lipschitz functions in space

$$\mathfrak{L} = \left\{ u : \mathbb{R}_{\geq 0} \times \mathbb{R}^d \to \mathbb{R} \mid u \text{ is continuous with } \sup_{t \geq 0} \|u(t,\cdot)\|_{\mathrm{Lip}} < +\infty \right\}. \tag{3.39}$$

**Theorem 3.5** (Comparison principle). *Let $u, v \in \mathfrak{L}$ be a viscosity subsolution and a viscosity supersolution to the Hamilton-Jacobi equation (3.20) respectively. Introduce the Lipschitz constant $L := \max(\sup_{t \geq 0} \|u(t,\cdot)\|_{\mathrm{Lip}}, \sup_{t \geq 0} \|v(t,\cdot)\|_{\mathrm{Lip}})$ as well as the local Lipschitz constant*

$$V := \sup \left\{ \frac{|\mathsf{H}(p') - \mathsf{H}(p)|}{|p' - p|} \;\middle|\; |p|, |p'| \leq L \right\}. \tag{3.40}$$

*For every $R, M \in \mathbb{R}$ with $M > 2L$, the mapping*

$$(t,x) \mapsto u(t,x) - v(t,x) - M(|x| + Vt - R)_+ \tag{3.41}$$

*achieves its supremum at a point in $\{0\} \times \mathbb{R}^d$.*

In the statement above, we use the notation $r_+ := \max(0, r)$ to denote the positive part of a real number $r \in \mathbb{R}$. Since the proof of Theorem 3.5 is a bit long, we suggest that the reader ignore any term related to the cutoff in space $M(|x| + Vt - R)_+$ and its smoothed variants on first reading, in effect showing the comparison principle with the unbounded space domain $\mathbb{R}^d$ replaced by the compact space $\mathbb{T}^d$.

*Proof of Theorem 3.5.* Suppose for the sake of contradiction that there exists $T > 0$ with
$$\sup_{[0,T] \times \mathbb{R}^d} (u - v - \varphi) > \sup_{\{0\} \times \mathbb{R}^d} (u - v - \varphi), \tag{3.42}$$
where $\varphi(t,x) := M(|x| + Vt - R)_+$. The proof proceeds in three steps. First we smoothen and perturb (3.42), then we use a variable doubling argument to obtain a system of inequalities, and finally we contradict this system of inequalities.

*Step 1: smoothing and perturbing.* Let $\varepsilon_0 \in (0,1)$ be a parameter to be determined, and let $\theta \in C^\infty(\mathbb{R};\mathbb{R})$ be a non-decreasing function such that, for every $r \in \mathbb{R}$,
$$(r - \varepsilon_0)_+ \leqslant \theta(r) \leqslant r_+.$$

We introduce the function
$$\Phi(t,x) := M\theta\left(\left(\varepsilon_0 + |x|^2\right)^{1/2} + Vt - R\right)$$
defined on $\mathbb{R}_{\geqslant 0} \times \mathbb{R}^d$. The choice of $\theta$ and the bound $(a+b)_+ \leqslant a_+ + b_+$ imply that
$$\varphi(t,x) \leqslant \Phi(t,x) + M\varepsilon_0 \leqslant \varphi(t,x) + M\varepsilon_0^{1/2} + M\varepsilon_0,$$
where we have used that $(a+b)^{\frac{1}{2}} \leqslant a^{\frac{1}{2}} + b^{\frac{1}{2}}$ for $a,b > 0$. It follows by (3.42) that
$$\sup_{\{0\} \times \mathbb{R}^d} (u - v - \Phi) < \sup_{[0,T] \times \mathbb{R}^d} (u - v - \Phi) + M\varepsilon_0 + M\varepsilon_0^{1/2},$$
so choosing $\varepsilon_0 > 0$ small enough guarantees that
$$\sup_{[0,T] \times \mathbb{R}^d} (u - v - \Phi) > \sup_{\{0\} \times \mathbb{R}^d} (u - v - \Phi). \tag{3.43}$$

This is a smoothed version of the hypothesis (3.42) we aim to contradict. We now also add a cutoff function in time. For a small parameter $\varepsilon > 0$ to be determined, we introduce the function
$$\chi(t,x) := \Phi(t,x) + \frac{\varepsilon}{T-t},$$
defined on $\mathbb{R}_{\geqslant 0} \times \mathbb{R}^d$. Choosing $\varepsilon > 0$ small enough ensures that
$$\sup_{[0,T] \times \mathbb{R}^d} (u - v - \chi) > \sup_{\{0\} \times \mathbb{R}^d} (u - v - \chi). \tag{3.44}$$

*Step 2: system of inequalities.* For each $\alpha \geqslant 1$, we define the function $\Psi_\alpha : [0,T] \times \mathbb{R}^d \times [0,T] \times \mathbb{R}^d \to \mathbb{R} \cup \{-\infty\}$ by
$$\Psi_\alpha(t,x,t',x') := u(t,x) - v(t',x') - \frac{\alpha}{2}\left(|t-t'|^2 + |x-x'|^2\right) - \chi(t,x). \tag{3.45}$$

We now argue that the function $\Psi_\alpha$ achieves its supremum at a point $(t_\alpha, x_\alpha, t'_\alpha, x'_\alpha)$ which remains bounded as $\alpha$ tends to infinity. In order to do so, we write $C < +\infty$ to denote a constant whose value might change as we proceed through the argument, and which may depend on $T$, $M$, $R$, $V$, $\sup_{t \leq T}|u(t,0)|$, $\sup_{t \leq T}|v(t,0)|$ and $L$. For every $x \in \mathbb{R}^d$ with $|x| > R+1$, the bound $\Phi(t,x) \geq M(|x| + Vt - R - 1)_+$ reveals that

$$\Psi_\alpha(t,x,t',x') \leq L(|x|+|x'|) + u(t,0) - v(t',0) - \frac{\alpha}{2}|x-x'|^2 - \Phi(t,x)$$
$$\leq L(|x|+|x'|) - \frac{\alpha}{2}|x-x'|^2 - M|x| + C$$
$$\leq (2L-M)|x| + L|x-x'| - \frac{\alpha}{2}|x-x'|^2 + C$$
$$\leq (2L-M)|x| + \frac{L^2}{2\alpha} + C.$$

We have used that the function $y \mapsto Ly - \frac{\alpha}{2}y^2$ achieves its maximum at $y = \frac{L}{\alpha}$. Observe also that the supremum of (3.45) is bounded from below by $\Psi_\alpha(0,0,0,0)$, which does not depend on $\alpha$, and that $M > 2L$. This implies that $x_\alpha$ remains bounded as $\alpha$ tends to infinity, and that

$$\alpha\left(|t_\alpha - t'_\alpha|^2 + |x_\alpha - x'_\alpha|^2\right) + \chi(t_\alpha, x_\alpha) \leq C. \tag{3.46}$$

It follows that, up to the extraction of a subsequence, there exist $t_0 \in [0,T]$ and $x_0 \in \mathbb{R}^d$ such that $t_\alpha \to t_0$, $t'_\alpha \to t_0$, $x_\alpha \to x_0$ and $x'_\alpha \to x_0$ as $\alpha \to +\infty$. By (3.46), we have that $\chi(t_\alpha, x_\alpha)$ remains bounded, and thus in particular, it must be that $t_0 \in [0,T)$. On the other hand, the continuity of $u$, $v$ and $\chi$ together with the bounds

$$\sup_{[0,T] \times \mathbb{R}^d} (u - v - \chi) \leq \Psi_\alpha(t_\alpha, x_\alpha, t'_\alpha, x'_\alpha) \leq u(t_\alpha, x_\alpha) - v(t'_\alpha, x'_\alpha) - \chi(t_\alpha, x_\alpha)$$

imply that

$$(u - v - \chi)(t_0, x_0) = \sup_{[0,T] \times \mathbb{R}^d}(u - v - \chi).$$

By (3.44), we thus deduce that $t_0 \in (0,T)$. This means that $t_\alpha, t'_\alpha \in (0,T)$ for all $\alpha$ large enough. We have therefore found a quadruple $(t_\alpha, x_\alpha, t'_\alpha, x'_\alpha)$ such that $\Psi_\alpha$ achieves its supremum at $(t_\alpha, x_\alpha, t'_\alpha, x'_\alpha)$, and with $t_\alpha, t'_\alpha \in (0,T)$ for $\alpha$ large enough. With this in mind, we fix $\alpha \geq 1$ large enough, and introduce the smooth functions $\phi, \phi' \in C^\infty((0,T) \times \mathbb{R}^d; \mathbb{R})$ defined by

$$\phi(t,x) := v(t'_\alpha, x'_\alpha) + \frac{\alpha}{2}\left(|t - t'_\alpha|^2 + |x - x'_\alpha|^2\right) + \chi(t,x),$$
$$\phi'(t',x') := u(t_\alpha, x_\alpha) - \frac{\alpha}{2}\left(|t' - t_\alpha|^2 + |x' - x_\alpha|^2\right) - \chi(t_\alpha, x_\alpha).$$

Since $(t_\alpha, x_\alpha, t'_\alpha, x'_\alpha)$ maximizes $\Psi_\alpha$, the function $u - \phi$ achieves its maximum at the point $(t_\alpha, x_\alpha) \in \mathbb{R}_{>0} \times \mathbb{R}^d$, while the function $v - \phi'$ achieves its minimum at $(t'_\alpha, x'_\alpha) \in \mathbb{R}_{>0} \times \mathbb{R}^d$. It follows by the definition of a viscosity subsolution and supersolution that

$$(\partial_t \phi - \mathsf{H}(\nabla \phi))(t_\alpha, x_\alpha) \leq 0 \quad \text{and} \quad (\partial_t \phi' - \mathsf{H}(\nabla \phi'))(t'_\alpha, x'_\alpha) \geq 0. \tag{3.47}$$

This is the system of inequalities that we now strive to contradict.

*Step 3: reaching a contradiction.* A direct computation shows that

$$(\partial_t \phi - \mathsf{H}(\nabla \phi))(t_\alpha, x_\alpha) = \alpha(t_\alpha - t'_\alpha) + \partial_t \Phi(t_\alpha, x_\alpha) + \frac{\varepsilon}{(T - t_\alpha)^2} - \mathsf{H}(\nabla \phi(t_\alpha, x_\alpha)), \tag{3.48}$$

and

$$(\partial_t \phi' - \mathsf{H}(\nabla \phi'))(t'_\alpha, x'_\alpha) = \alpha(t_\alpha - t'_\alpha) - \mathsf{H}(\nabla \phi'(t'_\alpha, x'_\alpha)). \tag{3.49}$$

To compare these two quantities, we would like to use the local Lipschitz continuity of the non-linearity H. With the definition of $V$ in mind, we must therefore verify that

$$|\nabla \phi(t_\alpha, x_\alpha)| \leq L \quad \text{and} \quad |\nabla \phi'(t'_\alpha, x'_\alpha)| \leq L. \tag{3.50}$$

Fix $z \in \mathbb{R}^d$ and $\varepsilon > 0$. Since $u - \phi$ achieves a local maximum at $(t_\alpha, x_\alpha) \in \mathbb{R}_{>0} \times \mathbb{R}^d$, and $u$ is uniformly Lipschitz continuous with Lipschitz constant $L$,

$$\phi(t_\alpha, x_\alpha + \varepsilon z) - \phi(t_\alpha, x_\alpha) \geq u(t_\alpha, x_\alpha + \varepsilon z) - u(t_\alpha, x_\alpha) \geq -\varepsilon L |z|.$$

Dividing by $\varepsilon$ and letting $\varepsilon$ tend to zero reveals that $\nabla \phi(t_\alpha, x_\alpha) \cdot z \geq -L|z|$. Choosing $z = -\nabla \phi(t_\alpha, x_\alpha)$ gives the first inequality in (3.50); the second inequality is obtained in an identical manner. Together with (3.49) and the definition of $V$, this implies that (3.49) is bounded from above by

$$\alpha(t_\alpha - t'_\alpha) - \mathsf{H}(\nabla \phi(t_\alpha, x_\alpha)) + V|\nabla \phi(t_\alpha, x_\alpha) - \nabla \phi'(t'_\alpha, x'_\alpha)|$$
$$= \alpha(t_\alpha - t'_\alpha) - \mathsf{H}(\nabla \phi(t_\alpha, x_\alpha)) + V|\nabla \Phi(t_\alpha, x_\alpha)|.$$

A direct computation shows that $V|\nabla \Phi(t_\alpha, x_\alpha)| \leq \partial_t \Phi(t_\alpha, x_\alpha)$, so in fact

$$(\partial_t \phi' - \mathsf{H}(\nabla \phi'))(t'_\alpha, x'_\alpha) \leq \alpha(t_\alpha - t'_\alpha) + \partial_t \Phi(t_\alpha, x_\alpha) - \mathsf{H}(\nabla \phi(t_\alpha, x_\alpha))$$
$$< (\partial_t \phi - \mathsf{H}(\nabla \phi))(t_\alpha, x_\alpha) \leq 0,$$

where the strict inequality is due to the term $\frac{\varepsilon}{(T-t)^2}$ in (3.48), and the final inequality leverages the first inequality in (3.47). This contradicts the second inequality in (3.47) and completes the proof. ∎

**Corollary 3.6.** *If $u, v \in \mathcal{L}$ are a viscosity subsolution and a viscosity supersolution to the Hamilton-Jacobi equation (3.20) respectively, then*

$$\sup_{\mathbb{R}_{\geq 0} \times \mathbb{R}^d} (u - v) = \sup_{\{0\} \times \mathbb{R}^d} (u - v). \tag{3.51}$$

*Proof.* Suppose for the sake of contradiction that there is a point $(t^*, x^*) \in \mathbb{R}_{>0} \times \mathbb{R}^d$ such that
$$(u-v)(t^*, x^*) > \sup_{\{0\} \times \mathbb{R}^d} (u-v). \tag{3.52}$$

In the notation of Theorem 3.5, we choose $M := 2L+1$ and $R := |x^*| + Vt^*$, so that
$$u(t^*, x^*) - v(t^*, x^*) - M(|x^*| + Vt^* - R)_+ = (u-v)(t^*, x^*).$$

By the assumption (3.52), this is strictly greater than
$$\sup_{x \in \mathbb{R}^d} (u(0,x) - v(0,x)) \geq \sup_{x \in \mathbb{R}^d} (u(0,x) - v(0,x) - M(|x| - R)_+).$$

This contradicts Theorem 3.5 and thus completes the proof. ∎

Finally, we verify that the comparison principle implies the uniqueness of a viscosity solution with a given initial condition.

**Corollary 3.7.** *If $u, v \in \mathcal{L}$ are two viscosity solutions to the Hamilton-Jacobi equation (3.20) with the same initial condition $u(0, \cdot) = v(0, \cdot)$, then $u = v$.*

*Proof.* Since $u$ is a viscosity subsolution and $v$ is a viscosity supersolution to the Hamilton-Jacobi equation (3.20), the comparison principle in Corollary 3.6 implies that $u \leq v$. A symmetric argument reveals that $v \leq u$, and completes the proof. ∎

To relate this uniqueness result back to the limit free energy in the Curie-Weiss model, we recall what we have done so far. Using the Arzelà-Ascoli theorem, we argued that the sequence $(F_N)_{N \geq 1}$ of free energies in the Curie-Weiss model is precompact. In Proposition 3.4, leveraging the approximate Hamilton-Jacobi equation (3.7), we showed that any subsequential limit must be a viscosity solution to the Hamilton-Jacobi equation (3.9) with initial condition $\psi(h) := \log \cosh(h)$. Together with the uniqueness result in Corollary 3.7, this implies that $(F_N)_{N \geq 1}$ must itself converge to the unique solution $f$ to the Hamilton-Jacobi equation (3.7).

We have thus succeeded in showing the convergence of the free energy in the Curie-Weiss model, and we have found an intrinsic characterization of its limit in terms of Hamilton-Jacobi equations. We could go about and study properties of this limit using PDE techniques, most importantly the comparison principle we have just seen. For instance, we can observe that the initial condition $\psi = \log \cosh$ for the Curie-Weiss model is non-negative and bounded from above by a parabola: there exists $C < +\infty$ such that $\psi(h) \leq Ch^2$. Using the comparison principle, we deduce that the solution $f$ will remain non-negative and bounded from above by the solution with the initial condition $h \mapsto Ch^2$. This solution can be identified explicitly as $(t,h) \mapsto \frac{Ch^2}{1-4Ct}$, for $t < 1/(4C)$. So we deduce in particular that $f$ will remain differentiable at $h = 0$, with null derivatives, at least up until the time $t = 1/(4C)$.

Of course, we also know from the arguments based on large deviation principles that the limit free energy can be written as a variational formula, see Corollary 2.20. As we will now see, this variational formula can be recovered using general results for solutions to Hamilton-Jacobi equations of the form (3.20), under suitable convexity assumptions.

**Exercise 3.6.** Let $\mathsf{H} \colon \mathbb{R}^d \to \mathbb{R}$ be a uniformly Lipschitz non-linearity with Lipschitz constant $V > 0$, let $f$ be a viscosity subsolution to (3.20), and let $\Phi \in C^\infty(\mathbb{R}_{>0} \times \mathbb{R}^d; \mathbb{R})$ be a smooth function with $\partial_t \Phi \geqslant V|\nabla \Phi|$. Show that $f - \Phi$ is a subsolution to (3.20).

**Exercise 3.7.** Let $f$ be a Lipschitz viscosity solution to (3.20), and let $L := \|f(0,\cdot)\|_{\mathrm{Lip}}$ denote the Lipschitz constant of its initial condition. Show that $\sup_{t \geqslant 0} \|f(t,\cdot)\|_{\mathrm{Lip}} \leqslant L$.

**Exercise 3.8.** Fix $L \in \mathbb{R}$. Let $f$ and $g$ be two Lipschitz viscosity solutions to (3.20) with $\sup_{t \geqslant 0} \|f(t,\cdot)\|_{\mathrm{Lip}} \leqslant L$, $\sup_{t \geqslant 0} \|g(t,\cdot)\|_{\mathrm{Lip}} \leqslant L$, and let $V$ be defined by (3.40). Suppose that there exist $x \in \mathbb{R}^d$ and $T > 0$ with $f(0,\cdot) = g(0,\cdot)$ on $B_{VT}(x)$. Prove that $f(T,x) = g(T,x)$.

**Exercise 3.9.** Let $\mathsf{K} \colon \mathbb{R}^d \to \mathbb{R}$ be a Lipschitz function, let $f$ be a Lipschitz viscosity solution to (3.20), and let $g$ be a Lipschitz viscosity solution to $\partial_t g - \mathsf{K}(\nabla g) = 0$ on $\mathbb{R}_{>0} \times \mathbb{R}^d$. Suppose that $f(0,\cdot) = g(0,\cdot)$, that $L$ denotes a Lipschitz constant of $f$ and $g$, and that $\mathsf{H}$ and $\mathsf{K}$ coincide on $B_L(0)$. Prove that $f = g$.

## 3.4 Variational representations of viscosity solutions

The Hamilton-Jacobi approach we have developed so far is able to identify the limit free energy in the Curie-Weiss model as the unique viscosity solution to the Hamilton-Jacobi equation (3.9). In this section, we develop variational representations for the solution to the Hamilton-Jacobi equation (3.20), under appropriate convexity assumptions. This will allow us to recover the variational formula for the limit free energy in the Curie-Weiss model obtained in Corollary 2.20 using the Hamilton-Jacobi approach. We will prove two variational formulas, the Hopf-Lax formula in the setting when the non-linearity $\mathsf{H}$ is convex, and the Hopf formula in the setting when the initial condition $\psi$ is convex. This implies in particular that solutions to the Hamilton-Jacobi equation (3.20) do exist under the stated convexity assumptions. A general result of existence of solutions for these equations can be obtained using the classical Perron method, but we will not develop this point here, and only refer the interested reader to [36, 84, 104] for more on this. In the context of the problems of statistical mechanics that we aim to study, the existence of solutions should come for free since the limit free energy is expected to be such a solution.

As we will see, under different assumptions, the Hopf-Lax and the Hopf formulas allow us to write the solution to the Hamilton-Jacobi equation (3.20) with initial condition $\psi$ as a saddle-point problem for the functional defined, for each $(t,x) \in \mathbb{R}_{\geq 0} \times \mathbb{R}^d$ and $(y,p) \in \mathbb{R}^d \times \mathbb{R}^d$, by

$$\mathcal{J}_{t,x}(y,p) := \psi(y) + p \cdot (x-y) + t\mathsf{H}(p). \tag{3.53}$$

More precisely, the Hopf-Lax and Hopf formulas are representations of the value of the solution to the Hamilton-Jacobi equation (3.20) at the point $(t,x) \in \mathbb{R}_{\geq 0} \times \mathbb{R}^d$ as saddle-point problems over the variables $(y,p) \in \mathbb{R}^d \times \mathbb{R}^d$. Notice that for each fixed $(y,p) \in \mathbb{R}^d \times \mathbb{R}^d$, the mapping $(t,x) \mapsto \mathcal{J}_{t,x}(y,p)$ is a solution to (3.20), so we can think of the Hopf-Lax and Hopf formulas as representations of the solution with initial condition $\psi$ as envelopes of this family of solutions. At a point of differentiability of this envelope function, the function will be tangent to one particular solution in this set indexed by $(y,p) \in \mathbb{R}^d \times \mathbb{R}^d$, and since the equation is of first order, a function that is tangent to a smooth solution at a point must solve the equation at that point. Of course, this is only an informal discussion, and the purpose of the next two sections is to prove the Hopf-Lax and Hopf formulas rigorously.

Since the proofs in this section are a bit long and the ideas developed there will not reappear later, the reader may consider skipping these proofs on first reading.

### 3.4.1   Hopf-Lax formula

The Hopf-Lax formula is a variational representation for the unique solution to the Hamilton-Jacobi equation (3.20) in the setting when the non-linearity $\mathsf{H}$ is convex. It reads as follows.

**Theorem 3.8** (Hopf-Lax formula [165]). *If $\psi : \mathbb{R}^d \to \mathbb{R}$ is a Lipschitz continuous initial condition and $\mathsf{H} : \mathbb{R}^d \to \mathbb{R}$ is a convex and locally Lipschitz continuous non-linearity, then the Hopf-Lax function*

$$f(t,x) := \sup_{y \in \mathbb{R}^d} \inf_{p \in \mathbb{R}^d} \mathcal{J}_{t,x}(y,p) = \sup_{y \in \mathbb{R}^d} \left( \psi(y) - t\mathsf{H}^*\left(\frac{y-x}{t}\right) \right) \tag{3.54}$$

*is the unique viscosity solution in $\mathcal{L}$ to the Hamilton-Jacobi equation (3.20).*

We recall from (2.14) that $\mathsf{H}^*$ denotes the convex dual of $\mathsf{H}$. If we assume that the non-linearity $\mathsf{H}$ is concave instead of convex, then we can appeal to Theorem 3.8 to obtain a variational representation of the viscosity solution to (3.20) as well. Indeed, this follows from the observation that a function $f$ is a viscosity solution to (3.20) if and only if the function $g := -f$ is a viscosity solution to

$$\partial_t g + \mathsf{H}(-\nabla g) = 0. \tag{3.55}$$

In the case when H is concave, we thus obtain that the solution $f$ to (3.20) with initial condition $\psi$ is given by

$$f(t,x) = \inf_{y \in \mathbb{R}^d} \sup_{p \in \mathbb{R}^d} \mathcal{J}_{t,x}(y,p). \tag{3.56}$$

To prove Theorem 3.8, we first verify that the function (3.54) satisfies the right initial condition, and that the supremum in its definition is attained. We then show that it satisfies a semigroup property from which we deduce that it belongs to the space $\mathcal{L}$. Finally, we show that it is the unique viscosity solution in $\mathcal{L}$ to the Hamilton-Jacobi equation (3.20) with initial condition $\psi$.

**Lemma 3.9.** *Under the assumptions of Theorem 3.8, the Hopf-Lax function (3.54) satisfies the right initial condition,*

$$f(0,\cdot) = \psi. \tag{3.57}$$

*Proof.* We interpret the formula in (3.54) at $t = 0$ as

$$f(0,x) = \sup_{y \in \mathbb{R}^d} \inf_{p \in \mathbb{R}^d} \big(\psi(y) + p \cdot (x-y)\big).$$

Taking $y = x$ on the right side of this expression gives the lower bound

$$f(0,x) \geq \psi(x).$$

On the other hand, given $y \in \mathbb{R}^d$, choosing

$$p = p(y) := \begin{cases} -\|\psi\|_{\mathrm{Lip}} \frac{x-y}{|x-y|} & \text{if } y \neq x \\ 0 & \text{otherwise} \end{cases}$$

gives the upper bound

$$f(0,x) \leq \sup_{y \in \mathbb{R}^d} \big(\psi(y) + p(y) \cdot (x-y)\big) = \sup_{y \in \mathbb{R}^d} \big(\psi(y) - \|\psi\|_{\mathrm{Lip}} |x-y|\big) \leq \psi(x),$$

where we have used the Lipschitz continuity of the initial condition $\psi$. Combining these lower and upper bounds completes the proof. ∎

**Lemma 3.10.** *Under the assumptions of Theorem 3.8, for any $(t,x) \in \mathbb{R}_{\geq 0} \times \mathbb{R}^d$, there exists $y \in \mathbb{R}^d$ with*

$$f(t,x) = \psi(y) - t\mathrm{H}^*\Big(\frac{y-x}{t}\Big). \tag{3.58}$$

*Proof.* Fix $\lambda > 0$ and $y \in \mathbb{R}^d$, and observe that

$$\mathsf{H}^*(y) = \sup_{p \in \mathbb{R}^d} \left( p \cdot y - \mathsf{H}(p) \right) \geq \lambda |y| - \sup_{|z| \leq \lambda} |\mathsf{H}(z)|,$$

where the inequality is obtained by taking $p = \lambda \frac{y}{|y|}$. Since H is locally Lipschitz continuous, and therefore locally bounded, the supremum on the right side of this expression is finite. Dividing by $|y|$, letting $\lambda$ tend to infinity, and then letting $|y|$ tend to infinity reveals that

$$\liminf_{|y| \to +\infty} \frac{\mathsf{H}^*(y)}{|y|} = +\infty. \tag{3.59}$$

This confirms that $t\mathsf{H}^*(\frac{y-x}{t})$ should be interpreted as $+\infty$ when $t = 0$ and $y \neq x$. It also implies that, given a point $(t,x) \in \mathbb{R}_{\geq 0} \times \mathbb{R}^d$, it is possible to find $R > 0$ large enough such that for all $y \in \mathbb{R}^d$ with $|y - x| > tR$, we have

$$t\mathsf{H}^*\left(\frac{y-x}{t}\right) \geq (\|\psi\|_{\mathrm{Lip}} + 1)|y - x|. \tag{3.60}$$

It follows by the Lipschitz continuity of $\psi$ that for all $y \in \mathbb{R}^d$ with $|y - x| > tR$, we have

$$\psi(y) - t\mathsf{H}^*\left(\frac{y-x}{t}\right) \leq \psi(x) - |y - x|.$$

This means that the supremum defining the Hopf-Lax function in (3.54) may be restricted to a bounded set. Together with the fact that the function $y \mapsto \psi(y) - t\mathsf{H}^*(\frac{y-x}{t})$ is upper semicontinuous and locally bounded from above, as $\mathsf{H}^*(z) \geq -\mathsf{H}(0)$, this implies that the supremum on the right side of (3.54) is achieved at some point $y \in \mathbb{R}^d$. This completes the proof. ∎

**Lemma 3.11** (Semigroup property). *Under the assumptions of Theorem 3.8, for every pair $t > s \geq 0$ and $x \in \mathbb{R}^d$,*

$$f(t,x) = \sup_{y \in \mathbb{R}^d} \left( f(s,y) - (t-s)\mathsf{H}^*\left(\frac{y-x}{t-s}\right) \right). \tag{3.61}$$

*Proof.* Fix $y, z \in \mathbb{R}^d$. Since $\mathsf{H}^*$ is convex by Exercise 2.10, we have

$$\mathsf{H}^*\left(\frac{y-x}{t}\right) \leq \frac{s}{t}\mathsf{H}^*\left(\frac{y-z}{s}\right) + \frac{t-s}{t}\mathsf{H}^*\left(\frac{z-x}{t-s}\right).$$

Substituting this bound into (3.54) yields

$$f(t,x) \geq \psi(y) - s\mathsf{H}^*\left(\frac{y-z}{s}\right) - (t-s)\mathsf{H}^*\left(\frac{z-x}{t-s}\right).$$

Taking the supremum over all $y \in \mathbb{R}^d$ gives

$$f(t,x) \geq f(s,z) - (t-s)\mathsf{H}^*\left(\frac{z-x}{t-s}\right),$$

and taking the supremum over all $z \in \mathbb{R}^d$ establishes the lower bound

$$f(t,x) \geq \sup_{y \in \mathbb{R}^d}\left(f(s,y) - (t-s)\mathsf{H}^*\left(\frac{y-x}{t-s}\right)\right).$$

To obtain the matching upper bound, we invoke Lemma 3.10 to find $y \in \mathbb{R}^d$ with

$$f(t,x) = \psi(y) - t\mathsf{H}^*\left(\frac{y-x}{t}\right).$$

Defining $z := \frac{s}{t}x + \frac{t-s}{t}y$, we observe that

$$\frac{z-x}{t-s} = \frac{y-x}{t} = \frac{y-z}{s}.$$

In particular, taking $y \in \mathbb{R}^d$ in (3.54) gives

$$f(s,z) - (t-s)\mathsf{H}^*\left(\frac{z-x}{t-s}\right) \geq \psi(y) - s\mathsf{H}^*\left(\frac{y-z}{s}\right) - (t-s)\mathsf{H}^*\left(\frac{z-x}{t-s}\right)$$
$$= \psi(y) - t\mathsf{H}^*\left(\frac{y-x}{t}\right)$$
$$= f(t,x).$$

Taking the supremum over $z \in \mathbb{R}^d$ establishes the matching upper bound and completes the proof. ∎

**Lemma 3.12.** *Under the assumptions of Theorem 3.8, we have $f \in \mathcal{L}$ with*

$$\sup_{t \geq 0}\|f(t,\cdot)\|_{\mathrm{Lip}} = \|\psi\|_{\mathrm{Lip}}. \tag{3.62}$$

*Proof.* Fix $(t,x,x') \in \mathbb{R}_{\geq 0} \times \mathbb{R}^d \times \mathbb{R}^d$ and invoke Lemma 3.10 to find $y \in \mathbb{R}^d$ with

$$f(t,x) = \psi(y) - t\mathsf{H}^*\left(\frac{y-x}{t}\right).$$

Taking $y - x + x' \in \mathbb{R}^d$ in (3.54) gives the lower bound

$$f(t,x') \geq \psi(y-x+x') - t\mathsf{H}^*\left(\frac{y-x}{t}\right).$$

It follows that

$$f(t,x) - f(t,x') \leq \psi(y) - \psi(y-x+x') \leq \|\psi\|_{\mathrm{Lip}}|x-x'|.$$

Reversing the roles of $x$ and $x'$ gives $y' \in \mathbb{R}^d$ with

$$f(t,x') - f(t,x) \leq \psi(y') - \psi(y' - x' + x) \leq \|\psi\|_{\mathrm{Lip}}|x - x'|.$$

Combining these two bounds establishes the spatial Lipschitz continuity (3.62) of the Hopf-Lax function. To conclude that $f \in \mathfrak{L}$, there remains to show that $f$ is also continuous in time; we show that it is in fact Lipschitz continuous in time as well. Fix $x \in \mathbb{R}^d$ and $t > s \geq 0$. The semigroup property in Lemma 3.11 implies that

$$f(t,x) \geq f(s,x) - (t-s)\mathsf{H}^*(0) \geq f(s,x) - (t-s)\mathsf{H}(0). \tag{3.63}$$

Combining the semigroup property in Lemma 3.11 with the spatial Lipschitz continuity (3.62) reveals that

$$\begin{aligned} f(t,x) &\leq f(s,x) + \sup_{y \in \mathbb{R}^d}\left(\|\psi\|_{\mathrm{Lip}}|y-x| - (t-s)\mathsf{H}^*\left(\frac{y-x}{t-s}\right)\right) \\ &= f(s,x) + (t-s)\sup_{z \in \mathbb{R}^d}\left(\|\psi\|_{\mathrm{Lip}}|z| - \mathsf{H}^*(z)\right) \\ &\leq f(s,x) + (t-s)\sup_{|p| \leq \|\psi\|_{\mathrm{Lip}}}\sup_{z \in \mathbb{R}^d}\left(z \cdot p - \mathsf{H}^*(z)\right), \end{aligned}$$

where the final inequality uses that $\|\psi\|_{\mathrm{Lip}}|z| = z \cdot \frac{\|\psi\|_{\mathrm{Lip}} z}{|z|}$. Invoking the Fenchel-Moreau theorem (Theorem 2.5) and remembering (3.63) establishes the temporal Lipschitz continuity of the Hopf-Lax function,

$$|f(t,x) - f(s,x)| \leq |t-s|\sup_{|p| \leq \|\psi\|_{\mathrm{Lip}}}|\mathsf{H}(p)|. \tag{3.64}$$

The convexity of the non-linearity $\mathsf{H}$ has played its part. This completes the proof. ∎

*Proof of Theorem 3.8.* The proof proceeds in two steps. First we show that the Hopf-Lax function (3.54) is a viscosity supersolution to the Hamilton-Jacobi equation (3.20), and then that it is also a viscosity subsolution to this equation. Together with Lemma 3.12 and the uniqueness result in Corollary 3.7, this proves that the Hopf-Lax function (3.54) is the unique viscosity solution in $\mathfrak{L}$ to the Hamilton-Jacobi equation (3.20).

*Step 1: viscosity supersolution.* Consider a smooth function $\phi \in C^\infty(\mathbb{R}_{>0} \times \mathbb{R}^d; \mathbb{R})$ with the property that $f - \phi$ has a local minimum at the point $(t^*, x^*) \in \mathbb{R}_{>0} \times \mathbb{R}^d$. By definition of a local minimum, for every $s \in (0, t^*)$ sufficiently small and $y \in \mathbb{R}^d$, we have

$$\phi(t^*, x^*) - \phi(t^* - s, x^* + sy) \geq f(t^*, x^*) - f(t^* - s, x^* + sy).$$

It follows by the semigroup property in Lemma 3.11 that

$$\phi(t^*,x^*) - \phi(t^* - s, x^* + sy) \geq s\mathsf{H}^*\left(\frac{sy}{s}\right) = s\mathsf{H}^*(y).$$

Dividing by $s$ and letting $s$ tend to zero reveals that

$$\partial_t \phi(t^*,x^*) - y \cdot \nabla \phi(t^*,x^*) + \mathsf{H}^*(y) \geq 0.$$

Taking the infimum over $y \in \mathbb{R}^d$ and invoking the Fenchel-Moreau theorem (Theorem 2.5) shows that

$$\bigl(\partial_t \phi - \mathsf{H}(\nabla \phi)\bigr)(t^*,x^*) = \partial_t \phi(t^*,x^*) - \mathsf{H}^{**}\bigl(\nabla \phi(t^*,x^*)\bigr) \geq 0.$$

The convexity of the non-linearity $\mathsf{H}$ has played its part. This verifies the supersolution criterion for the Hopf-Lax function $f$.

*Step 2: viscosity subsolution.* Consider a smooth function $\phi \in C^\infty(\mathbb{R}_{>0} \times \mathbb{R}^d; \mathbb{R})$ with the property that $f - \phi$ has a local maximum at the point $(t^*,x^*) \in \mathbb{R}_{>0} \times \mathbb{R}^d$, and suppose for the sake of contradiction that there is $\delta > 0$ such that for all $(t',x')$ sufficiently close to $(t^*,x^*)$, we have

$$\bigl(\partial_t \phi - \mathsf{H}(\nabla \phi)\bigr)(t',x') \geq \delta > 0.$$

Using the convexity of $\mathsf{H}$ and the Fenchel-Moreau theorem, this may be recast as the assumption that for all $(t',x')$ sufficiently close to $(t^*,x^*)$ and all $y \in \mathbb{R}^d$, we have

$$\partial_t \phi(t',x') - y \cdot \nabla \phi(t',x') + \mathsf{H}^*(y) \geq \delta. \tag{3.65}$$

Leveraging the semigroup property in Lemma 3.11 and arguing as in the proof of Lemma 3.10, it is possible to find $R > 0$ with the property that for all $s > 0$ sufficiently small there is $y_s \in \mathbb{R}^d$ with $|x^* - y_s| \leq Rs$ and

$$f(t^*,x^*) = f(t^* - s, y_s) - s\mathsf{H}^*\left(\frac{y_s - x}{s}\right).$$

Writing $u(r) := (rt^* + (1-r)(t^* - s), rx^* + (1-r)y_s)$, it follows by the fundamental theorem of calculus and the absurd assumption (3.65) with $y = \frac{y_s - x^*}{s}$ that

$$\begin{aligned}
\phi(t^*,x^*) - \phi(t^* - s, y_s) &= \int_0^1 \frac{\mathrm{d}}{\mathrm{d}r} \phi(u(r)) \, \mathrm{d}r \\
&= \int_0^1 \bigl(s\partial_t \phi + (x^* - y_s) \cdot \nabla \phi\bigr)(u(r)) \, \mathrm{d}r \\
&\geq s\delta - s\mathsf{H}^*\left(\frac{y_s - x^*}{s}\right) \\
&= s\delta + f(t^*,x^*) - f(t^* - s, y_s).
\end{aligned}$$

Rearranging shows that for $s$ sufficiently small,

$$f(t^* - s, y_s) - \phi(t^* - s, y_s) \geq s\delta + f(t^*, x^*) - \phi(t^*, x^*).$$

This contradicts the local maximality of $f - \phi$ at $(t^*, x^*)$ and completes the proof. ∎

### 3.4.2   Hopf formula

The Hopf formula is a variational representation for the unique solution to the Hamilton-Jacobi equation (3.20) in the setting when the initial condition $\psi$ is convex. It reads as follows.

**Theorem 3.13** (Hopf formula [37, 166]). *If $\psi : \mathbb{R}^d \to \mathbb{R}$ is a Lipschitz continuous and convex initial condition, and $\mathsf{H} : \mathbb{R}^d \to \mathbb{R}$ is a locally Lipschitz continuous non-linearity, then the Hopf function*

$$f(t,x) := \sup_{p \in \mathbb{R}^d} \inf_{y \in \mathbb{R}^d} \mathcal{J}_{t,x}(y,p) = \sup_{p \in \mathbb{R}^d} \inf_{y \in \mathbb{R}^d} \left( \psi(y) + p \cdot (x-y) + t\mathsf{H}(p) \right) \quad (3.66)$$

*is the unique viscosity solution in $\mathfrak{L}$ to the Hamilton-Jacobi equation* (3.20).

Compared with the Hopf-Lax formula in (3.54), notice that the supremum is now taken over $p$ while the infimum is taken over $y$. Assuming that the infimum over $y$ is achieved and that $\psi$ is smooth, a differentiation in $y$ shows that $p = \nabla\psi(y)$ at the optimum. For each fixed $y$, the mapping $(t,x) \mapsto \psi(y) + \nabla\psi(y) \cdot (x-y) + t\mathsf{H}(\nabla\psi(y))$ is the solution to the Hamilton-Jacobi equation (3.20) with initial condition given by the tangent of $\psi$ at $y$. Since $\psi$ is convex, the comparison principle tells us that the solution to (3.20) should indeed lie above this family of solutions.

If we assume instead that the initial condition $\psi$ is concave, then we can again use the observation around (3.55) and appeal to Theorem 3.13 to obtain that the viscosity solution to (3.20) is given by

$$f(t,x) = \inf_{p \in \mathbb{R}^d} \sup_{y \in \mathbb{R}^d} \mathcal{J}_{t,x}(y,p). \quad (3.67)$$

We now turn to the proof of Theorem 3.13. In order to do so, we first verify that the function (3.66) satisfies the right initial condition, and that the supremum in its definition is attained. We next argue that this function is jointly convex and belongs to the space $\mathfrak{L}$. Finally, we show that it satisfies a semigroup property, and deduce that it is the unique viscosity solution in $\mathfrak{L}$ to the Hamilton-Jacobi equation (3.20). It will be convenient to notice that the Hopf function (3.66) can be written as

$$f(t,x) = \sup_{p \in \mathbb{R}^d} \left( p \cdot x + t\mathsf{H}(p) - \psi^*(p) \right) = \left( \psi^* - t\mathsf{H} \right)^*(x). \quad (3.68)$$

**Lemma 3.14.** *Under the assumptions of Theorem 3.13, the Hopf function* (3.66) *satisfies the right initial condition,*

$$f(0,\cdot) = \psi. \tag{3.69}$$

*Proof.* Taking $t = 0$ in the representation (3.68) and combining the convexity of $\psi$ with the Fenchel-Moreau theorem (Theorem 2.5) completes the proof. ∎

**Lemma 3.15.** *Under the assumptions of Theorem 3.13, for every $(t,x) \in \mathbb{R}_{\geqslant 0} \times \mathbb{R}^d$, there exists $p \in \mathbb{R}^d$ with $|p| \leqslant \|\psi\|_{\mathrm{Lip}}$ such that*

$$f(t,x) = p \cdot x + t\mathsf{H}(p) - \psi^*(p). \tag{3.70}$$

*Proof.* Observe that $\psi^*(p) = +\infty$ whenever $|p| > \|\psi\|_{\mathrm{Lip}}$ by Exercise 2.12. This means that the supremum in (3.68) may be restricted to the compact ball of radius $\|\psi\|_{\mathrm{Lip}}$ about the origin. Using also that the function $\psi^*$ is lower semicontinuous, by part (i) of Exercise 2.10, this implies that the supremum in (3.68) is achieved, as announced. ∎

**Lemma 3.16.** *Under the assumptions of Theorem 3.13, the Hopf function $f$ is jointly convex and belongs to the space $\mathfrak{L}$ with*

$$\sup_{t \geqslant 0} \|f(t,\cdot)\|_{\mathrm{Lip}} = \|\psi\|_{\mathrm{Lip}}. \tag{3.71}$$

*Proof.* The first equality of (3.68) is a representation of $f$ as a supremum of affine functions of $(t,x)$, so $f$ is jointly convex by Exercise 2.7. To prove Lipschitz continuity, we fix $t,t' \geqslant 0$ as well as $x,x' \in \mathbb{R}^d$, and invoke Lemma 3.15 to find $p \in \mathbb{R}^d$ with $|p| \leqslant \|\psi\|_{\mathrm{Lip}}$ and $f(t,x) = p \cdot x + t\mathsf{H}(p) - \psi^*(p)$. The Cauchy-Schwarz inequality implies that

$$\begin{aligned} f(t,x) - f(t',x') &\leqslant p \cdot (x - x') + \mathsf{H}(p)(t - t') \\ &\leqslant |p||x - x'| + |\mathsf{H}(p)||t - t'| \\ &\leqslant \|\psi\|_{\mathrm{Lip}}|x - x'| + \sup_{|q| \leqslant \|\psi\|_{\mathrm{Lip}}} |\mathsf{H}(q)||t - t'|. \end{aligned}$$

An identical argument with the roles of $(x,t)$ and $(x',t')$ reversed establishes the Lipschitz continuity (3.71) of the Hopf function. It also shows that $f$ is continuous in $(t,x)$, and thus that $f \in \mathfrak{L}$. ∎

**Lemma 3.17** (Semigroup property). *Under the assumptions of Theorem 3.13, for every $t,s \geqslant 0$ and $x \in \mathbb{R}^d$,*

$$f(t+s,\cdot) = (f^*(t,\cdot) - s\mathsf{H})^*. \tag{3.72}$$

*Proof.* We will use properties (i)-(iii) of the convex dual established in Exercise 2.10. Property (ii) implies that $(\psi^* - t\mathsf{H})^{**} \leq \psi^* - t\mathsf{H}$. Combining this with the representation (3.68) of the Hopf function reveals that

$$f^*(t, \cdot) - s\mathsf{H} = (\psi^* - t\mathsf{H})^{**} - s\mathsf{H} \leq \psi^* - (t+s)\mathsf{H},$$

and leveraging property (iii) gives the upper bound

$$(f^*(t, \cdot) - s\mathsf{H})^* \geq (\psi^* - (t+s)\mathsf{H})^* = f(t+s, \cdot). \tag{3.73}$$

Similarly, property (ii) implies that

$$\psi^* - (t+s)\mathsf{H} \geq (\psi^* - (t+s)\mathsf{H})^{**}.$$

Multiplying through by $\frac{t}{t+s} = 1 - \frac{s}{t+s}$ reveals that

$$\psi^* - t\mathsf{H} \geq \frac{s}{t+s}\psi^* + \frac{t}{t+s}(\psi^* - (t+s)\mathsf{H})^{**}. \tag{3.74}$$

The functions $\psi^*$ and $(\psi^* - (t+s)\mathsf{H})^{**}$ are both convex by property (i), so, as a positive linear combination of convex functions, the right side of (3.74) defines a convex function. It follows by the Fenchel-Moreau theorem and two applications of property (iii) that

$$(\psi^* - t\mathsf{H})^{**} \geq \frac{s}{t+s}\psi^* + \frac{t}{t+s}(\psi^* - (t+s)\mathsf{H})^{**},$$

and therefore,

$$(\psi^* - (t+s)\mathsf{H})^{**} - (\psi^* - t\mathsf{H})^{**} \leq \frac{s}{t+s}\left((\psi^* - (t+s)\mathsf{H})^{**} - \psi^*\right)$$
$$\leq \frac{s}{t+s}\left(\psi^* - (t+s)\mathsf{H} - \psi^*\right)$$
$$= -s\mathsf{H},$$

where the second inequality uses property (ii). Leveraging the dual representation (3.68), this inequality may be rewritten as

$$f^*(t+s, \cdot) \leq f^*(t, \cdot) - s\mathsf{H}.$$

Remembering that $f$ is convex by Lemma 3.16, appealing to the Fenchel-Moreau theorem, and using property (iii) once more gives the lower bound

$$f(t+s, \cdot) = f^{**}(t+s, \cdot) \geq (f^*(t, \cdot) - s\mathsf{H})^*.$$

Together with the upper bound (3.73), this completes the proof. ∎

*Proof of Theorem 3.13.* The proof proceeds in two steps. First we show that the Hopf function (3.66) is a viscosity subsolution to the Hamilton-Jacobi equation (3.20), and then that it is also a viscosity supersolution to this equation. Together with Lemma 3.16 and the uniqueness result in Corollary 3.7, this proves that the Hopf function (3.66) is the unique viscosity solution in $\mathfrak{L}$ to the Hamilton-Jacobi equation (3.20).

*Step 1: viscosity subsolution.* Consider a smooth function $\phi \in C^\infty(\mathbb{R}_{>0} \times \mathbb{R}^d; \mathbb{R})$ with the property that $f - \phi$ has a local maximum at the point $(t^*, x^*) \in \mathbb{R}_{>0} \times \mathbb{R}^d$. Invoking Lemma 3.15, we may find $p \in \mathbb{R}^d$ with $|p| \leq \|\psi\|_{\mathrm{Lip}}$ and

$$f(t^*, x^*) = p \cdot x^* + t^* \mathsf{H}(p) - \psi^*(p).$$

Together with the definition of a local maximum, this implies that for $s \geq 0$ small enough, $z \in \mathbb{R}^d$ and $\varepsilon \geq 0$ small enough,

$$\phi(t^*, x^*) - \phi(t^* - s, x^* + \varepsilon z) \leq f(t^*, x^*) - f(t^* - s, x^* + \varepsilon z)$$
$$\leq -p \cdot \varepsilon z + s \mathsf{H}(p), \quad (3.75)$$

where we used the expression (3.68) for $f(t^* - s, x^* + \varepsilon z)$ in the second step. Taking $s = 0$, dividing by $\varepsilon$, and letting $\varepsilon$ tend to zero shows that $(\nabla \phi(t^*, x^*) - p) \cdot z \geq 0$. Since $z \in \mathbb{R}^d$ is arbitrary, this means that

$$p = \nabla \phi(t^*, x^*). \quad (3.76)$$

On the other hand, taking $\varepsilon = 0$, dividing by $s$, and letting $s$ tend to zero in (3.75) reveals that

$$\partial_t \phi(t^*, x^*) \leq \mathsf{H}(p).$$

Rearranging and remembering (3.76) shows that $(\partial_t \phi - \mathsf{H}(\nabla \phi))(t^*, x^*) \leq 0$, which is the required subsolution criterion.

*Step 2: viscosity supersolution.* Consider a smooth function $\phi \in C^\infty(\mathbb{R}_{>0} \times \mathbb{R}^d; \mathbb{R})$ with the property that $f - \phi$ has a local minimum at the point $(t^*, x^*) \in \mathbb{R}_{>0} \times \mathbb{R}^d$. We fix $\lambda \in (0, 1]$ and $(t', x') \in \mathbb{R}_{\geq 0} \times \mathbb{R}^d$. The joint convexity of $f$ established in Lemma 3.16 and the definition of a local minimum imply that, for $\lambda$ small enough,

$$f(t', x') - f(t^*, x^*) \geq \frac{1}{\lambda}\big(f(t^* + \lambda(t' - t^*), x^* + \lambda(x' - x^*)) - f(t^*, x^*)\big)$$
$$\geq \frac{1}{\lambda}\big(\phi(t^* + \lambda(t' - t^*), x^* + \lambda(x' - x^*)) - \phi(t^*, x^*)\big).$$

Letting $\lambda$ tend to zero shows that

$$f(t', x') - f(t^*, x^*) \geq (t' - t^*)\partial_t \phi(t^*, x^*) + \nabla \phi(t^*, x^*) \cdot (x' - x^*).$$

Together with the semigroup property in Lemma 3.17 and two applications of property (ii) in Exercise 2.10, this implies that

$$f(t^*, x^*) = \left(f^*(t^* - s, \cdot) - s\mathsf{H}\right)^*(x^*)$$
$$\geqslant \left(\left(f(t^*, x^*) - s\partial_t \phi(t^*, x^*) + \nabla \phi(t^*, x^*) \cdot (\cdot - x^*)\right)^* - s\mathsf{H}\right)^*(x^*).$$

To leverage this bound, observe that

$$\left(f(t^*, x^*) - s\partial_t \phi(t^*, x^*) + \nabla \phi(t^*, x^*) \cdot (\cdot - x^*)\right)^*(q)$$
$$= s\partial_t \phi(t^*, x^*) - f(t^*, x^*) + \nabla \phi(t^*, x^*) \cdot x^* + \sup_{y \in \mathbb{R}^d} y \cdot \left(q - \nabla \phi(t^*, x^*)\right),$$

which is infinite unless $q = \nabla \phi(t^*, x^*)$. It thus follows from the definition of the convex dual that

$$f(t^*, x^*) \geqslant x^* \cdot \nabla \phi(t^*, x^*)$$
$$- \left(s\partial_t \phi(t^*, x^*) - f(t^*, x^*) + \nabla \phi(t^*, x^*) \cdot x^* - s\mathsf{H}(\nabla \phi(t^*, x^*))\right).$$

Rearranging and dividing by $s > 0$ shows that $(\partial_t \phi - \mathsf{H}(\nabla \phi))(t^*, x^*) \geqslant 0$ and completes the proof. ∎

Of course, when both the non-linearity $\mathsf{H}$ and the initial condition $\psi$ are convex, the uniqueness result in Corollary 3.7 ensures that the Hopf-Lax and the Hopf variational representations coincide. It is also possible to show this directly using the Fenchel-Moreau theorem.

**Proposition 3.18.** *If $f, g : \mathbb{R}^d \to \mathbb{R}$ are two convex functions, then*

$$\sup_{x \in \mathbb{R}^d} \inf_{y \in \mathbb{R}^d} \left(f(x) + g(y) - x \cdot y\right) = \sup_{y \in \mathbb{R}^d} \inf_{x \in \mathbb{R}^d} \left(f(x) + g(y) - x \cdot y\right). \quad (3.77)$$

*In particular, when $\psi : \mathbb{R}^d \to \mathbb{R}$ and $\mathsf{H} : \mathbb{R}^d \to \mathbb{R}$ are both convex, the Hopf-Lax function (3.54) and the Hopf function (3.66) coincide.*

*Proof.* Using the definition of the convex dual $g^*$ and then applying the Fenchel-Moreau theorem to the convex function $f$ shows that

$$\sup_{x \in \mathbb{R}^d} \inf_{y \in \mathbb{R}^d} \left(f(x) + g(y) - x \cdot y\right) = \sup_{x \in \mathbb{R}^d} \left(f(x) - g^*(x)\right)$$
$$= \sup_{x \in \mathbb{R}^d} \sup_{y \in \mathbb{R}^d} \left(x \cdot y - f^*(y) - g^*(x)\right).$$

Since this expression is symmetric in the pair $(f, g)$, this establishes (3.77). More explicitly, we can write

$$\sup_{y \in \mathbb{R}^d} \sup_{x \in \mathbb{R}^d} \left(x \cdot y - f^*(y) - g^*(x)\right) = \sup_{y \in \mathbb{R}^d} \left(g(y) - f^*(y)\right)$$
$$= \sup_{y \in \mathbb{R}^d} \inf_{x \in \mathbb{R}^d} \left(f(x) + g(y) - x \cdot y\right).$$

When both H and $\psi$ are convex, we can make the change of variables $z = x - y$ in (3.54) and (3.66) to rewrite the statement that the Hopf-Lax and Hopf formulas coincide as

$$\sup_{z \in \mathbb{R}^d} \inf_{p \in \mathbb{R}^d} \left( \psi(x+z) + t\mathsf{H}(p) - p \cdot z \right) = \sup_{p \in \mathbb{R}^d} \inf_{z \in \mathbb{R}^d} \left( \psi(x+z) + t\mathsf{H}(p) - p \cdot z \right).$$

The fact that this identity is valid follows from (3.77). This completes the proof. ∎

Applying the Hopf formula to the non-linearity $\mathsf{H}(p) = p^2$ and the convex initial condition $\psi(h) := \log \cosh(h)$ in the context of Proposition 3.4, we find that the limit free energy in the Curie-Weiss model is given by

$$f(t,h) = \sup_{m \in \mathbb{R}^d} \left( mh + t\mathsf{H}(m) - \psi^*(m) \right) = \sup_{m \in [-1,1]} \left( mh + tm^2 - \psi^*(m) \right), \quad (3.78)$$

where the second equality combines the fact that $\psi$ is Lipschitz continuous with Lipschitz constant one with Exercise 2.12. We have therefore given a proof of Corollary 2.20 on the limit of the free energy in the Curie-Weiss model using only the Hamilton-Jacobi approach.

## 3.5 Variational representations in the absence of convexity

In this section, we explore whether any aspect of the variational structures discussed in the previous section remains in the absence of any convexity or concavity assumption on the non-linearity or the initial condition. The considerations discussed here will not reappear until Section 6.6, which is the last section of the book, so the reader may consider skipping this section on first reading.

For the generalized Curie-Weiss model and the model from statistical inference studied in Chapter 4, as well as for the broader class of models from statistical inference considered in [71] and discussed briefly at the end of Section 4.3, the free energy is always convex. If we can show that it converges to the solution to a Hamilton-Jacobi equation, then we can represent its limit variationally using the Hopf variational formula. This convexity property of the free energy is however lost in the realm of spin glasses. In some cases, the non-linearity appearing in the relevant Hamilton-Jacobi equation is convex. This will allow us to appeal to a Hopf-Lax formula instead and to still represent the solution variationally; we will see in Chapter 6 that this variational representation is closely related to the Parisi formula. However, there are many spin-glass models of interest for which the non-linearity in the equation is neither convex nor concave. It is therefore interesting to wonder whether, for general solutions to Hamilton-Jacobi equations, any aspect of the variational structure displayed in the Hopf and Hopf-Lax formulas is preserved in the absence of any convexity or concavity assumption on the initial condition or the non-linearity. Recall that the Hopf-Lax and Hopf formulas are expressed as variational problems for the function $\mathcal{J}_{t,x}$ defined in (3.53). The specific property we want to discuss is whether we can always represent the solution $f(t,x)$ to the Hamilton-Jacobi equation (3.20) as a critical value of the function $\mathcal{J}_{t,x}$; in other words, we ask whether we can find a point $(y^*, p^*)$ such that $\nabla \mathcal{J}_{t,x}(y^*, p^*) = 0$ and $f(t,x) = \mathcal{J}_{t,x}(y^*, p^*)$. The next result shows that this property is indeed valid as soon as the initial condition or the non-linearity in the equation is convex (or concave). The gradient $\nabla \mathcal{J}_{t,x}$ is understood as being jointly in the two arguments of $\mathcal{J}_{t,x}$, that is, $\nabla \mathcal{J}_{t,x} = (\nabla_y \mathcal{J}_{t,x}, \nabla_p \mathcal{J}_{t,x})$.

**Theorem 3.19.** *Let $\psi \in C^1(\mathbb{R}^d; \mathbb{R})$ be Lipschitz continuous and let $\mathsf{H} \in C^1(\mathbb{R}^d; \mathbb{R})$. Suppose that among the functions $\psi$ and $\mathsf{H}$, at least one of them is convex or concave, and let $f : \mathbb{R}_{\geqslant 0} \times \mathbb{R}^d \to \mathbb{R}$ be the viscosity solution to the Hamilton-Jacobi equation (3.20) with initial condition $f(0, \cdot) = \psi$. Then*

$$\text{there exists } (y^*, p^*) \in \mathbb{R}^d \times \mathbb{R}^d \text{ with } \nabla \mathcal{J}_{t,x}(y^*, p^*) = 0$$
$$\text{such that } f(t,x) = \mathcal{J}_{t,x}(y^*, p^*). \quad (3.79)$$

We recall that under the assumptions of Theorem 3.19, the viscosity solution to (3.20) with initial condition $\psi$ exists and admits a variational representation in terms of $\mathcal{J}_{t,x}$, by Theorems 3.8 and 3.13. For instance, if $\mathsf{H}$ is convex, then by Theorem 3.8, we have

$$f(t,x) = \sup_{y \in \mathbb{R}^d} \inf_{p \in \mathbb{R}^d} \mathcal{J}_{t,x}(y, p). \quad (3.80)$$

We stress that the proof of Theorem 3.19 must leverage the convexity assumption, say on $\mathsf{H}$ here, beyond the validity of the representation (3.80). In other words, in general, it is not true that the mapping

$$(t,x) \mapsto \sup_{y \in \mathbb{R}^d} \inf_{p \in \mathbb{R}^d} \mathcal{J}_{t,x}(y, p) \quad (3.81)$$

will satisfy the property in (3.79), even in cases when the supremum and infimum are achieved. To give an example, consider the one-dimensional case with $\psi(y) = \log \cosh(y)$ and $\mathsf{H}(p) = (p^2 - 1)^2$. One can then verify that for every $x \in \mathbb{R}$ and $t \geqslant 0$,

$$\psi(x) = \sup_{y \in \mathbb{R}} \inf_{p \in \mathbb{R}} \mathcal{J}_{t,x}(y, p), \quad (3.82)$$

and that the supremum and infimum are achieved. Yet, there will typically be no $(y^*, p^*)$ with $\nabla \mathcal{J}_{t,x}(y^*, p^*) = 0$ such that $\psi(x) = \mathcal{J}_{t,x}(y^*, p^*)$. For small $t \geqslant 0$, such a critical point $(y^*, p^*) = (y^*(t,x), p^*(t,x))$ is unique, and the mapping $(t,x) \mapsto \mathcal{J}_{t,x}(y^*(t,x), p^*(t,x))$ is in fact the solution to (3.20), as discussed further below in relation with the method of characteristics.

*Proof of Theorem 3.19.* We show the result assuming that H is convex. The other cases are identical, up to a change of sign or an interchange of the roles of $\psi$ and H. We recall that since H is convex, the function $f$ satisfies (3.80). We give ourselves a small parameter $\delta \in (0,1]$, and define the perturbed non-linearity

$$H_\delta(p) := H(p) + \delta |p|^2,$$

so that the mapping $H_\delta$ is strictly convex, see (2.120) for a definition of strict convexity. We denote by $\mathcal{J}_{t,x}^{(\delta)}$ the functional $\mathcal{J}_{t,x}$ with the non-linearity H replaced by the perturbed non-linearity $H_\delta$. We decompose the rest of the proof into four steps. First, we show that the saddle-point problem (3.81) associated with the perturbed functional $\mathcal{J}_{t,x}^{(\delta)}$ admits an optimizer. We then prove that the family $(H_\delta^*)_{\delta>0}$ satisfies a lower semicontinuity property as $\delta$ tends to zero, and use this to show that the saddle-point problem (3.81) associated with the perturbed functional $\mathcal{J}_{t,x}^{(\delta)}$ is continuous at zero with respect to the perturbation parameter. We then combine these results to establish (3.79).

*Step 1: identifying an optimizer to the perturbed saddle-point problem.* In this step, we show that for each $\delta \in (0,1]$,

$$\text{there exists } (y_\delta^*, p_\delta^*) \in \mathbb{R}^d \times \mathbb{R}^d \text{ with } \nabla \mathcal{J}_{t,x}^{(\delta)}(y_\delta^*, p_\delta^*) = 0$$
$$\text{such that } \sup_{y \in \mathbb{R}^d} \inf_{p \in \mathbb{R}^d} \mathcal{J}_{t,x}^{(\delta)}(y,p) = \mathcal{J}_{t,x}^{(\delta)}(y_\delta^*, p_\delta^*). \quad (3.83)$$

We will also make sure that $(y_\delta^*, p_\delta^*)$ is chosen in a fixed compact set not depending on $\delta > 0$, so that we can then extract a converging subsequence as $\delta \to 0$.

It follows from (3.59)-(3.60) that we can find a compact set $K$ such that for every $\delta \in [0,1]$,

$$\sup_{y \in \mathbb{R}^d} \inf_{p \in \mathbb{R}^d} \mathcal{J}_{t,x}^{(\delta)}(y,p) = \sup_{y \in K} \inf_{p \in \mathbb{R}^d} \mathcal{J}_{t,x}^{(\delta)}(y,p). \quad (3.84)$$

Moreover, as was argued below (3.60), the supremum on the right side of (3.84) is achieved, say at some $y_\delta^* \in K$. Using that $H_\delta$ grows at least quadratically, one can also make sure that, for every $y \in \mathbb{R}^d$, there exists $p_\delta(y)$ such that

$$\inf_{p \in \mathbb{R}^d} \mathcal{J}_{t,x}^{(\delta)}(y,p) = \mathcal{J}_{t,x}^{(\delta)}(y, p_\delta(y)). \quad (3.85)$$

Since H is differentiable, we must have that $\nabla_p \mathcal{J}_{t,x}^{(\delta)}(y, p_\delta(y)) = 0$. Since $H_\delta$ is strictly convex, there is in fact exactly one $p_\delta(y)$ that realizes (3.85) (see the solution to Exercise 2.23 for a detailed argument). The envelope theorem (Theorem 2.21) therefore ensures that the mapping

$$y \mapsto \inf_{p \in \mathbb{R}^d} \mathcal{J}_{t,x}^{(\delta)}(y,p)$$

is differentiable, and that its gradient is $\nabla_y \mathcal{J}_{t,x}^{(\delta)}(y, p_\delta(y))$. By the optimality condition for $y_\delta^*$, we conclude that

$$\nabla_y \mathcal{J}_{t,x}^{(\delta)}(y_\delta, p_\delta(y_\delta)) = 0.$$

This means that $\nabla \psi(y_\delta) = p_\delta(y_\delta)$. In particular, we must have that $|p_\delta(y)| \leqslant \|\Psi\|_{\text{Lip}}$. We have therefore obtained (3.83), with $y_\delta^* \in K$ and $p_\delta^* := p_\delta(y_\delta^*)$ such that $|p_\delta^*| \leqslant \|\Psi\|_{\text{Lip}}$.

*Step 2: establishing a lower semicontinuity property of* $(H_\delta^*)$. Let $(q_\delta)_{\delta > 0}$ be a family of points in $\mathbb{R}^d$ that converge to $q \in \mathbb{R}^d$. In this step, we show that

$$\liminf_{\delta \to 0} H_\delta^*(q_\delta) \geqslant H^*(q). \tag{3.86}$$

More explicitly, the claim is that

$$\liminf_{\delta \to 0} \sup_{p \in \mathbb{R}^d} \left( p \cdot q_\delta - H(p) - \delta |p|^2 \right) \geqslant \sup_{p \in \mathbb{R}^d} \left( p \cdot q - H(p) \right). \tag{3.87}$$

Assuming first that $H^*(q)$ is finite, we can find $p \in \mathbb{R}^d$ such that

$$p \cdot q - H(p) \geqslant H^*(q) - \varepsilon.$$

Using $p$ as a candidate in the variational problem on the left side of (3.87) yields the desired inequality, up to an error of $\varepsilon$. Since $\varepsilon > 0$ was arbitrary, this completes the argument. If instead $H^*(q)$ is infinite, then we can argue similarly, since for each $M \in \mathbb{R}$, we can find $p \in \mathbb{R}^d$ with $p \cdot q - H(p) \geqslant M$. This completes the proof of (3.86).

*Step 3: showing the perturbed saddle-point problem is continuous at zero.* In this step, we show that

$$\lim_{\delta \to 0} \sup_{y \in \mathbb{R}^d} \inf_{p \in \mathbb{R}^d} \mathcal{J}_{t,x}^{(\delta)}(y, p) = \sup_{y \in \mathbb{R}^d} \inf_{p \in \mathbb{R}^d} \mathcal{J}_{t,x}(y, p). \tag{3.88}$$

We can rewrite (3.88) as

$$\lim_{\delta \to 0} \sup_{y \in \mathbb{R}^d} \left( \psi(y) - H_\delta^*\left(\frac{y - x}{t}\right) \right) = \sup_{y \in \mathbb{R}^d} \left( \psi(y) - H^*\left(\frac{y - x}{t}\right) \right). \tag{3.89}$$

Since for each $y \in \mathbb{R}^d$, the family $(H_\delta^*(y))_{\delta > 0}$ is non-increasing, the limit on the left side of (3.89) exists, and since $H_\delta^* \leqslant H^*$, the statement of (3.89) with the equality sign replaced by $\geqslant$ is valid. For the converse inequality, we recall that the supremum on the left side of (3.89) is achieved at some $y_\delta^* \in K$, where $K$ is a fixed compact set.

Up to the extraction of a subsequence, we may assume that $y^*_\delta$ converges to some $y^* \in K$. By the result of the previous step, we obtain that

$$\limsup_{\delta \to 0} \sup_{y \in \mathbb{R}^d} \left( \psi(y) - \mathsf{H}^*_\delta\left(\frac{y-x}{t}\right) \right) = \lim_{\delta \to 0} \left( \psi(y^*_\delta) - \mathsf{H}^*_\delta\left(\frac{y^*_\delta - x}{t}\right) \right)$$
$$\leqslant \psi(y^*) - \mathsf{H}^*\left(\frac{y^* - x}{t}\right)$$
$$\leqslant \sup_{y \in \mathbb{R}^d} \left( \psi(y) - \mathsf{H}^*\left(\frac{y-x}{t}\right) \right).$$

This completes the proof of (3.89).

*Step 4: concluding.* We are now ready to conclude the proof. Up to the extraction of a subsequence, we may assume that the points $y^*_\delta$ and $p^*_\delta$ constructed in Step 1 converge as $\delta > 0$ tends to zero, say to $y^* \in K$ and $|p^*| \leqslant \|\Psi\|_{\mathrm{Lip}}$. By the result of the previous step, we have

$$f(t,x) = \sup_{y \in \mathbb{R}^d} \inf_{p \in \mathbb{R}^d} \mathcal{J}_{t,x}(y,p) = \limsup_{\delta \to 0} \sup_{y \in \mathbb{R}^d} \inf_{p \in \mathbb{R}^d} \mathcal{J}^{(\delta)}_{t,x}(y,p)$$
$$= \lim_{\delta \to 0} \mathcal{J}^{(\delta)}_{t,x}(y^*_\delta, p^*_\delta)$$
$$= \mathcal{J}_{t,x}(y^*, p^*).$$

Since the functions $\mathsf{H}$ and $\psi$ are continuously differentiable, we can also verify by continuity that $\nabla \mathcal{J}_{t,x}(y^*, p^*) = 0$, as desired. ∎

In fact, for smooth initial condition $\psi$, say $\psi \in C^2(\mathbb{R}^d; \mathbb{R})$, the property in (3.79) remains valid for small $t \geqslant 0$ even in the absence of any convexity or concavity assumption on $\psi$ or $\mathsf{H}$. To explain this best, we start by explaining a classical method for solving the Hamilton-Jacobi equation (3.20) for a short time called the method of characteristics.

Suppose for a moment that, for some $T \in (0, +\infty]$, we have found a smooth solution $f \in C^2([0,T) \times \mathbb{R}^d; \mathbb{R})$ to the Hamilton-Jacobi equation (3.20) with initial condition $\psi \in C^2(\mathbb{R}^d; \mathbb{R})$, and that the non-linearity $\mathsf{H}$ is also in $C^2(\mathbb{R}^d; \mathbb{R})$. We fix $i \in \{1, \dots, d\}$ and write $u_i := \partial_{x_i} f$. Differentiating (3.20), we find that

$$\partial_t u_i(t,x) - \nabla \mathsf{H}(\nabla f(t,x)) \cdot \nabla u_i(t,x) = 0. \tag{3.90}$$

For each $t \in [0,T)$ and $x \in \mathbb{R}^d$, we define $X(t,x)$ such that $X(0,x) = x$ and $X(\cdot, x)$ solves the ordinary differential equation $\partial_t X(t,x) = -\nabla \mathsf{H}(\nabla f(t, X(t,x)))$. The curve $t \mapsto X(t,x)$ (or also $t \mapsto (t, X(t,x))$) is called the *characteristic* starting from $x \in \mathbb{R}^d$. We find using (3.90) that

$$\frac{\mathrm{d}}{\mathrm{d}t}\left(u_i(t, X(t,x))\right)$$
$$= \partial_t u_i(t, X(t,x)) - \nabla \mathsf{H}(\nabla f(t, X(t,x))) \cdot \nabla u_i(t, X(t,x)) = 0. \tag{3.91}$$

In other words, the derivatives of $f$ are constant along the curve $t \mapsto (t, X(t,x))$; and in fact, this curve is therefore a straight line, with $\partial_t X(t,x) = -\nabla \mathsf{H}(\nabla \psi(x))$. Finally, since

$$\frac{d}{dt}\bigl(f(t,X(t,x))\bigr) = \partial_t f(t,X(t,x)) - \nabla \mathsf{H}(\nabla \psi(x)) \cdot \nabla f(t,X(t,x))$$
$$= \mathsf{H}(\nabla f(t,X(t,x))) - \nabla \mathsf{H}(\nabla \psi(x)) \cdot \nabla f(t,X(t,x))$$
$$= \mathsf{H}(\nabla \psi(x)) - \nabla \mathsf{H}(\nabla \psi(x)) \cdot \nabla \psi(x), \tag{3.92}$$

we conclude that

$$f(t,X(t,x)) = \psi(x) + t\bigl(\mathsf{H}(\nabla \psi(x)) - \nabla \mathsf{H}(\nabla \psi(x)) \cdot \nabla \psi(x)\bigr). \tag{3.93}$$

To sum up, we have argued that if a solution to the Hamilton-Jacobi equation (3.20) is sufficiently smooth, then it must satisfy (3.93). Let us now try to use this observation to actually construct the function $f$. Given the non-linearity $\mathsf{H}$ and the initial condition $\psi$, both assumed to be in $C^2(\mathbb{R}^d;\mathbb{R})$ and with $\nabla \psi$ and $\nabla^2 \psi$ bounded, we set

$$X(t,x) := x - t\nabla \mathsf{H}(\nabla \psi(x)). \tag{3.94}$$

Using a fixed-point argument, one can check that for every $t \geq 0$ sufficiently small, the mapping $x \mapsto X(t,x)$ is bijective. This allows us to define the function $f$ according to the formula (3.93). One can then verify that $f$ is a classical solution to (3.20) on $[0,T) \times \mathbb{R}^d$ for a sufficiently small $T > 0$, as will be explored in more detail in Exercise 3.10 below. This perspective thus allows us to understand the possible low regularity of the viscosity solution in terms of characteristic lines intersecting each other.

The statement (3.79) we derived in Theorem 3.19 also has a very natural interpretation in terms of characteristic lines. Indeed, the condition that the pair $(y,p) \in \mathbb{R}^d \times \mathbb{R}^d$ satisfies $\nabla \mathcal{J}_{t,x}(y,p) = 0$ can be rewritten as

$$p = \nabla \psi(y) \quad \text{and} \quad x = y - t\nabla \mathsf{H}(p). \tag{3.95}$$

A pair $(y,p)$ therefore satisfies (3.95) if and only if $p = \nabla \psi(y)$ and the characteristic line $t \mapsto X(t,y) = y - t\nabla \mathsf{H}(\nabla \psi(y))$ satisfies $X(t,y) = x$. In short, there is a one-to-one correspondence between critical points of $\mathcal{J}_{t,x}$ and characteristic lines that pass through the point $(t,x)$. Moreover, for a pair $(y,p)$ satisfying (3.95), we have

$$\mathcal{J}_{t,x}(y,p) = \psi(y) + p \cdot (x-y) + t\mathsf{H}(p) \tag{3.96}$$
$$= \psi(y) + t\bigl(\mathsf{H}(\nabla \psi(y)) - \nabla \mathsf{H}(\nabla \psi(y)) \cdot \nabla \psi(y)\bigr), \tag{3.97}$$

in agreement with (3.93).

This connection between the method of characteristics and critical points of $\mathcal{J}_{t,x}$ has at least two interesting consequences. First, we learn that if $\psi$ and $\mathsf{H}$ are in $C^2(\mathbb{R}^d;\mathbb{R})$ with $\nabla\psi$ and $\nabla^2\psi$ bounded, then the property (3.79) is valid for every $t \geqslant 0$ sufficiently small, even in the absence of any convexity of concavity assumption on H or $\psi$. Second, we can also reinterpret Theorem 3.19 as saying that if at least one function among $\psi$ or H is convex or concave, then at every $(t,x) \in \mathbb{R}_{\geqslant 0} \times \mathbb{R}^d$, including possibly large $t \geqslant 0$, there exists a characteristic line that prescribes the correct value for $f(t,x)$. To state this more precisely, we may introduce the *wavefront*

$$W := \Big\{\big(t,\, x - t\nabla\mathsf{H}(\nabla\psi(x)),$$
$$\psi(x) + t\big(\mathsf{H}(\nabla\psi(x)) - \nabla\mathsf{H}(\nabla\psi(x))\cdot\nabla\psi(x)\big)\big) \,\big|\, t \geqslant 0,\, x \in \mathbb{R}^d\Big\}. \quad (3.98)$$

Theorem 3.19 states that whenever H or $\psi$ is convex or concave, letting $f$ denote the viscosity solution to (3.20), we have that the graph of $f$ is a subset of the wavefront.

Since we have just observed that this property is also valid for short times without any convexity assumption, one may expect that the graph of the viscosity solution to (3.20) is in fact always a subset of the wavefront. This is however not the case. Indeed, it was shown in [232] that for any fixed non-linearity $\mathsf{H} \in C^2(\mathbb{R}^d;\mathbb{R})$ which is neither convex nor concave, there exists a smooth and Lipschitz continuous function $\psi : \mathbb{R}^d \to \mathbb{R}$ such that the graph of the viscosity solution to (3.20) with initial condition $\psi$ is not a subset of the wavefront.

On the other hand, for any given $\mathsf{H} \in C^2(\mathbb{R}^d;\mathbb{R})$ and Lipschitz continuous $\psi : \mathbb{R}^d \to \mathbb{R}$, one can show that there always exists a Lipschitz function $f : \mathbb{R}_{\geqslant 0} \times \mathbb{R}^d \to \mathbb{R}$ whose graph belongs to the wavefront [231]. In particular, such a function $f$ satisfies the Hamilton-Jacobi equation (3.20) almost everywhere, and $f(0,\cdot) = \psi$. Such a function is however not unique in general.

In the context of the general non-convex spin-glass models discussed in Section 6.6, one can show that any subsequential limit of the free energy must satisfy a natural Hamilton-Jacobi equation "almost everywhere". (The quotes are due because the state space is infinite-dimensional.) It is not known whether the limit free energy is the viscosity solution to the equation. But surprisingly, it is shown in [70] that, assuming that the limit free energy exists, its graph must belong to the wavefront. We do not know whether the graph of the viscosity solution belongs to the wavefront in this case. We refer to Section 6.6 for a more precise discussion on this point.

**Exercise 3.10.** The goal of this exercise is to justify the method of characteristics for solving the Hamilton-Jacobi equation (3.20) for a short time. We assume that the non-linearity H and the initial condition $\psi$ both belong to $C^2(\mathbb{R}^d;\mathbb{R})$. We also assume that the first and second-order derivatives of $\psi$ are uniformly bounded.

(i) Show that there exists $T > 0$ such that for every $t \in [0,T)$, the mapping

$$\phi_t : x \mapsto x - t\nabla\mathsf{H}(\nabla\psi(x)) \quad (3.99)$$

is a $C^1$ diffeomorphism of $\mathbb{R}^d$; in other words, it is bijective and both $\phi_t$ and its inverse belong to $C^1(\mathbb{R}^d;\mathbb{R}^d)$.

(ii) We define the function $u : [0,T) \times \mathbb{R}^d \to \mathbb{R}$ implicitly by setting, for every $t \in [0,T)$ and $x \in \mathbb{R}^d$,

$$u(t,\phi_t(x)) := \psi(x) + t\big(\mathsf{H}(\nabla\psi(x)) - \nabla\mathsf{H}(\nabla\psi(x))\cdot\nabla\psi(x)\big). \qquad (3.100)$$

Show that $u \in C^1([0,T) \times \mathbb{R}^d;\mathbb{R})$ and that for every $t \in [0,T)$ and $x \in \mathbb{R}^d$, the gradient of $u$ is constant along the characteristic line $t \mapsto (t,\phi_t(x))$.

(iii) Deduce that $u$ satisfies the Hamilton-Jacobi equation (3.20) at every point in $[0,T) \times \mathbb{R}^d$.

## 3.6   Leveraging convexity to identify viscosity solutions

The Hopf formula and Proposition 3.4 allowed us to recover the formula obtained in Corollary 2.20 for the limit free energy in the Curie-Weiss model using the Hamilton-Jacobi approach. However, if we now try to apply the same reasoning to recover the formula obtained in Theorem 2.19 for the limit free energy in the generalized Curie-Weiss model, then we run into a problem. The upper bound

$$\partial_t F_N - \xi(\partial_h F_N) \leqslant \frac{C}{N}\partial_h^2 F_N \qquad (3.101)$$

implied by (3.17) should suffice to show that the limit free energy is a viscosity subsolution to the Hamilton-Jacobi equation (3.19); however the lower bound

$$\partial_t F_N - \xi(\partial_h F_N) \geqslant -\frac{C}{N}\partial_h^2 F_N \qquad (3.102)$$

implied by (3.17) should ring alarm bells since the *sign* in front of the Laplacian seems to be wrong. As discussed in Section 3.2, the notion of viscosity solutions is entirely built around the validity of a maximum principle; but changing the sign in front of the Laplacian destroys this property. To see this clearly, one can think about solutions of the backwards heat equation, $\partial_t f = -\Delta f$. Solutions to this equation may blow up in finite time, and it is not difficult to construct examples of initial conditions that are ordered, $f(0,\cdot) \leqslant g(0,\cdot)$, and such that this ordering is no longer valid at some subsequent time.

One may wonder whether the negative sign in front of the Laplacian in (3.102) is just an artefact; maybe a more refined analysis would allow us to obtain a better estimate. In the case when the function $\xi$ is convex, we have

$$\left\langle \xi\left(\frac{1}{N}\sum_{i=1}^N \sigma_i\right)\right\rangle \geqslant \xi\left(\left\langle \frac{1}{N}\sum_{i=1}^N \sigma_i\right\rangle\right), \qquad (3.103)$$

by Jensen's inequality. This gives the improved lower bound

$$\partial_t F_N(t,h) - \xi\left(\partial_h F_N(t,h)\right) \geq 0. \tag{3.104}$$

Arguing exactly as in the proof of Proposition 3.4, we thereby obtain the following extension of Theorem 2.19 in the convex setting. Notice that for convex $\xi$, this is indeed an extension since we no longer require that $\psi \in C^1(\mathbb{R};\mathbb{R})$.

**Proposition 3.20.** *For each integer $N \geq 1$, denote by $F_N : \mathbb{R}_{\geq 0} \times \mathbb{R} \to \mathbb{R}$ the free energy* (3.2) *in the generalized Curie-Weiss model, and suppose that for every $h \in \mathbb{R}$ the limit*

$$\psi(h) := \lim_{N \to +\infty} F_N(0,h) \tag{3.105}$$

*exists. If $\xi$ is convex, then the limit free energy $f : \mathbb{R}_{\geq 0} \times \mathbb{R} \to \mathbb{R}$ in the generalized Curie-Weiss model is the unique viscosity solution to the Hamilton-Jacobi equation* (3.19)*, and it admits the Hopf representation*

$$f(t,h) := \lim_{N \to +\infty} F_N(t,h) = \sup_{m \in [-1,1]} \left(t\xi(m) + hm - \psi^*(m)\right). \tag{3.106}$$

*Proof.* The sequence $(F_N)_{N \geq 1}$ of free energies is precompact by the Arzelà-Ascoli theorem and the Lipschitz bounds (3.13). Using the inequalities (3.101) and (3.104), we can follow the proof of Proposition 3.4 to show that any subsequential limit of $(F_N)_{N \geq 1}$ must be a viscosity solution to the Hamilton-Jacobi equation (3.19). Invoking the uniqueness result in Corollary 3.7, we conclude that $(F_N)_{N \geq 1}$ converges to the unique solution $f$ to the Hamilton-Jacobi equation (3.19). Since each of the initial conditions $h \mapsto F_N(0,h)$ is convex, their limit $\psi$ must also be convex. Invoking the Hopf formula in Theorem 3.13 yields that

$$f(t,h) = \sup_{m \in \mathbb{R}} \left(t\xi(m) + hm - \psi^*(m)\right). \tag{3.107}$$

To restrict this supremum to the interval $[-1,1]$, observe that $\partial_h F_N$ is uniformly bounded by one by the derivative computation (3.13), so the initial condition $\psi$ must be Lipschitz continuous with Lipschitz constant one. It follows by Exercise 2.12 that $\psi^*$ is infinite outside the interval $[-1,1]$. This completes the proof. ∎

When $\xi$ is convex, we were therefore able to improve upon the bound (3.102) and conclude. One may wonder if this could be achieved for general $\xi$. This is however not possible, as is most clearly demonstrated by considering the case of $\xi(p) = -p^2$. Indeed, we obtain in this case that

$$\partial_t F_N + \left(\partial_h F_N\right)^2 = -\frac{1}{N}\partial_h^2 F_N. \tag{3.108}$$

So there really is a fundamental problem with the sign here; it is not just that the bounds we have derived are too crude. And yet, we know from Theorem 2.19 and the Hopf formula that the sequence $(F_N)_{N \geqslant 1}$ of free energies does indeed converge to the viscosity solution to the Hamilton-Jacobi equation (3.19) even in this case, at least when the initial condition $\psi$ is in $C^1(\mathbb{R};\mathbb{R})$! To prove this using the Hamilton-Jacobi approach, we now introduce a new tool, which we call the convex selection principle, to identify when a *convex* function $f$ is a viscosity solution to the Hamilton-Jacobi equation (3.20).

In Example 3.1, we described a function that satisfies a Hamilton-Jacobi equation almost everywhere but is not a viscosity solution to this equation. Notice that this counterexample had corner singularities "in both directions"; formally, the second derivative was neither bounded from above nor from below. A convex function cannot look like this, since its Hessian must be non-negative. Roughly speaking, the convex selection principle states that imposing the function to be convex completely rules out the emergence of non-viscosity-type singularities, and thus restores uniqueness. We will establish this result in the setting of the Hamilton-Jacobi equation (3.20). Recall that the non-linearity $\mathsf{H} : \mathbb{R}^d \to \mathbb{R}$ is locally Lipschitz continuous.

**Theorem 3.21** (Convex selection principle [75]). *If $f : \mathbb{R}_{\geqslant 0} \times \mathbb{R}^d \to \mathbb{R}$ is a jointly convex and jointly Lipschitz continuous function that satisfies the Hamilton-Jacobi equation* (3.20) *on a dense subset of $\mathbb{R}_{\geqslant 0} \times \mathbb{R}^d$ and $f(0,\cdot) \in C^1(\mathbb{R}^d;\mathbb{R})$, then $f$ is a viscosity solution to the Hamilton-Jacobi equation* (3.20).

When we say that $f$ satisfies the Hamilton-Jacobi equation (3.20) on a dense subset of $\mathbb{R}_{\geqslant 0} \times \mathbb{R}^d$, we mean that the set

$$\left\{(t,x) \in \mathbb{R}_{>0} \times \mathbb{R}^d \mid f \text{ differentiable at } (t,x) \text{ and } \left(\partial_t f - \mathsf{H}(\nabla f)\right)(t,x) = 0\right\} \quad (3.109)$$

is dense in $\mathbb{R}_{\geqslant 0} \times \mathbb{R}^d$. Naturally, Theorem 3.21 also holds if we replace "convex" by "concave" in the statement. We already stressed that, in some sense, the notion of viscosity solution is tailored to approximations in which a small *positive* term times the Laplacian of $f$ appears on the right side of (3.20). In particular, the notion of viscosity solution is sensitive to the orientation of time; in general, it is *not* the case that the time-reversed viscosity solution to some equation will be the viscosity solution to the time-reversed equation. However, superficially, the statement of Theorem 3.21 looks invariant under time reversal. The only assumption that breaks this symmetry is that $f(0,\cdot) \in C^1(\mathbb{R}^d;\mathbb{R})$. This already hints at the fact that this assumption is necessary. In other words, Theorem 3.21 states that under the convexity assumption, pathological solutions cannot spontaneously emerge if we start from a smooth initial condition. However, if we start from a Lipschitz function that is not in $C^1(\mathbb{R}^d;\mathbb{R})$, then we may be able to exploit the singularities of the initial condition to create solutions that differ from the viscosity solution.

To prove the convex selection principle, we first show that the function $f$ in its statement must actually satisfy (3.20) at all its points of differentiability.

We then prove a convex selection principle with a somewhat weaker version of the assumption that $f(0,\cdot) \in C^1(\mathbb{R}^d;\mathbb{R})$ which will be convenient in Chapter 4 to analyze the model from statistical inference in full generality. On first reading, the reader may want to work out a simpler proof using directly the assumption that $f(0,\cdot) \in C^1(\mathbb{R}^d;\mathbb{R})$ — this stronger assumption will be valid in the most interesting cases discussed in Chapter 4. We recall that for a convex function $f : \mathbb{R}_{\geq 0} \times \mathbb{R}^d \to \mathbb{R}$, we denote by $\partial f(t,x) \subseteq \mathbb{R} \times \mathbb{R}^d$ the subdifferential of $f$ at the point $(t,x) \in \mathbb{R}_{\geq 0} \times \mathbb{R}^d$, as defined in (2.40).

**Lemma 3.22.** *If $f : \mathbb{R}_{\geq 0} \times \mathbb{R}^d \to \mathbb{R}$ is a jointly convex and jointly Lipschitz continuous function that satisfies the Hamilton-Jacobi equation (3.20) on a dense subset of $\mathbb{R}_{\geq 0} \times \mathbb{R}^d$, then it satisfies the Hamilton-Jacobi equation (3.20) at all its points of differentiability in $\mathbb{R}_{>0} \times \mathbb{R}^d$. Moreover, for every $(t,x) \in \mathbb{R}_{\geq 0} \times \mathbb{R}^d$, there exists $(a,p) \in \partial f(t,x)$ such that $a - \mathsf{H}(p) = 0$.*

*Proof.* Fix $(t,x) \in \mathbb{R}_{\geq 0} \times \mathbb{R}^d$, and let $(t_n, x_n)_{n \geq 1} \subseteq \mathbb{R}_{>0} \times \mathbb{R}^d$ be a sequence of points of differentiability of $f$ converging to $(t,x)$ at which

$$\big(\partial_t f - \mathsf{H}(\nabla f)\big)(t_n, x_n) = 0. \tag{3.110}$$

Since $f$ is differentiable at the interior point $(t_n, x_n) \in \mathbb{R}_{>0} \times \mathbb{R}^d$, Theorem 2.13 implies that $\partial f(t_n, x_n) = \{(\partial_t f(t_n, x_n), \nabla f(t_n, x_n))\}$. The joint Lipschitz continuity of $f$ implies that, up to the extraction of a subsequence, we may assume that the sequence of gradients $(\partial_t f(t_n, x_n), \nabla f(t_n, x_n))_{n \geq 1}$ converges to some vector $(a,p) \in \mathbb{R} \times \mathbb{R}^d$. By Proposition 2.14, we have $(a,p) \in \partial f(t,x)$, and by (3.110) and the continuity of H, the pair $(a,p) \in \partial f(t,x)$ is such that $a - \mathsf{H}(p) = 0$. This establishes the second part of the statement. If $t > 0$ and $f$ is differentiable at $(t,x)$, then Theorem 2.13 implies that $\partial f(t,x) = \{(\partial_t f(t,x), \nabla f(t,x))\}$. It must therefore be the case that $(a,p) = (\partial_t f(t,x), \nabla f(t,x))$, and thus that $f$ satisfies the Hamilton-Jacobi equation (3.20) at the point $(t,x)$. This completes the proof. ∎

We now state and show the refined version of Theorem 3.21.

**Lemma 3.23.** *Let $f : \mathbb{R}_{\geq 0} \times \mathbb{R}^d \to \mathbb{R}$ be a jointly convex and jointly Lipschitz continuous function that satisfies the Hamilton-Jacobi equation (3.20) on a dense subset of $\mathbb{R}_{\geq 0} \times \mathbb{R}^d$. Suppose that the initial condition $\psi := f(0,\cdot)$ is such that, for every $x \in \mathbb{R}^d$ and $p \in \partial \psi(x)$, there exists $b \in \mathbb{R}$ with $(b,p) \in \partial f(0,x)$ and $b - \mathsf{H}(p) \geq 0$. Then $f$ is a viscosity solution to the Hamilton-Jacobi equation (3.20).*

*Proof.* We decompose the proof into two steps. First we show that $f$ is a viscosity subsolution to the Hamilton-Jacobi equation (3.20), and then that it is a supersolution to this equation. The assumption on the initial condition will only play a role in showing that $f$ is a viscosity supersolution to (3.20).

*Step 1: viscosity subsolution.* Consider a smooth function $\phi \in C^\infty(\mathbb{R}_{>0} \times \mathbb{R}^d; \mathbb{R})$ with the property that $f - \phi$ has a local maximum at the point $(t^*, x^*) \in \mathbb{R}_{>0} \times \mathbb{R}^d$. We are going to show that $f$ is differentiable at $(t^*, x^*)$. By Proposition 2.11, the subdifferential $\partial f(t^*, x^*)$ contains at least one element, say $(a, p) \in \mathbb{R} \times \mathbb{R}^d$. By the definition of the subdifferential and of a local maximum, for every $(t', x')$ sufficiently close to $(t^*, x^*)$, we have

$$a(t' - t^*) + p \cdot (x' - x^*) \leqslant f(t', x') - f(t^*, x^*) \leqslant \phi(t', x') - \phi(t^*, x^*). \qquad (3.111)$$

It follows from the smoothness of $\phi$ that, as $(t', x')$ tends to $(t^*, x^*)$,

$$(t' - t^*)\big(a - \partial_t \phi(t^*, x^*)\big) + (x' - x^*)\big(p - \nabla \phi(t^*, x^*)\big) \leqslant o\big(|t' - t^*| + |x' - x^*|\big).$$

This implies that $(a, p) = (\partial_t \phi, \nabla \phi)(t^*, x^*)$. Using (3.111) once more, we obtain that $f$ is differentiable at $(t^*, x^*)$, and that $(\partial_t f, \nabla f)(t^*, x^*) = (\partial_t \phi, \nabla \phi)(t^*, x^*)$. It follows by Lemma 3.22 that

$$\big(\partial_t \phi - \mathsf{H}(\nabla \phi)\big)(t^*, x^*) = \big(\partial_t f - \mathsf{H}(\nabla f)\big)(t^*, x^*) = 0.$$

This completes the verification that $f$ is a viscosity subsolution to the Hamilton-Jacobi equation (3.20).

*Step 2: viscosity supersolution.* Consider a smooth function $\phi \in C^\infty(\mathbb{R}_{>0} \times \mathbb{R}^d; \mathbb{R})$ with the property that $f - \phi$ has a local minimum at the point $(t^*, x^*) \in \mathbb{R}_{>0} \times \mathbb{R}^d$. Together with the convexity of $f$, this implies that for every $(t', x') \in \mathbb{R}_{\geqslant 0} \times \mathbb{R}^d$ and $\varepsilon > 0$ small enough,

$$f(t', x') - f(t^*, x^*) \geqslant \varepsilon^{-1}\big(f\big((t^*, x^*) + \varepsilon(t' - t^*, x' - x^*)\big) - f(t^*, x^*)\big)$$
$$\geqslant \varepsilon^{-1}\big(\phi\big((t^*, x^*) + \varepsilon(t' - t^*, x' - x^*)\big) - \phi(t^*, x^*)\big).$$

Letting $\varepsilon$ tend to zero shows that $(\partial_t \phi, \nabla \phi)(t^*, x^*) \in \partial f(t^*, x^*)$. It therefore suffices to fix $(a, p) \in \partial f(t^*, x^*)$ and prove that

$$a - \mathsf{H}(p) \geqslant 0. \qquad (3.112)$$

Since $(a, p) \in \partial f(t^*, x^*)$, we have $f(0, y) \geqslant f(t^*, x^*) - at^* + (y - x^*) \cdot p$ for every $y \in \mathbb{R}^d$. Rearranging shows that for every $y \in \mathbb{R}^d$,

$$f(0, y) - y \cdot p \geqslant f(t^*, x^*) - at^* - p \cdot x^*. \qquad (3.113)$$

Inspired by the Hopf formula, we would like to consider an optimizer of the minimization problem

$$\inf_{y \in \mathbb{R}^d} \big(f(0, y) - y \cdot p\big). \qquad (3.114)$$

Since we cannot guarantee the existence of this optimizer, we introduce a small parameter $\varepsilon > 0$ and consider instead the perturbed minimization problem

$$\inf_{y \in \mathbb{R}^d} \left( f(0,y) - y \cdot p + \varepsilon \sqrt{1 + |y|^2} \right). \tag{3.115}$$

We see from (3.113) that we may restrict the infimum in (3.115) to values of $y$ that range in a compact set (which depends on $\varepsilon$). Together with the continuity of the functions involved, this shows that the infimum in (3.115) is achieved, say at $y_\varepsilon \in \mathbb{R}^d$. We observe that

$$\lim_{\varepsilon \to 0} \inf_{y \in \mathbb{R}^d} \left( f(0,y) - y \cdot p + \varepsilon \sqrt{1 + |y|^2} \right) = \inf_{y \in \mathbb{R}^d} \left( f(0,y) - y \cdot p \right). \tag{3.116}$$

Indeed, the existence of the limit on the left side of (3.116) and the fact that it is lower bounded by the right side are immediate. Conversely, for each $\delta > 0$, we can find $y_\delta^* \in \mathbb{R}^d$ such that

$$f(0, y_\delta^*) - y_\delta^* \cdot p \leq \inf_{y \in \mathbb{R}^d} \left( f(0,y) - y \cdot p \right) + \delta,$$

and we obtain the upper bound in (3.116) up to an error of $\delta > 0$ by using $y_\delta^*$ as a candidate in the infimum on the left side of (3.116). Since $\delta > 0$ was arbitrary, this shows (3.116). Since

$$\inf_{y \in \mathbb{R}^d} \left( f(0,y) - y \cdot p \right) + \varepsilon \sqrt{1 + |y_\varepsilon|^2} \leq f(0, y_\varepsilon) - y_\varepsilon \cdot p + \varepsilon \sqrt{1 + |y_\varepsilon|^2}$$

$$= \inf_{y \in \mathbb{R}^d} \left( f(0,y) - y \cdot p + \varepsilon \sqrt{1 + |y|^2} \right),$$

we deduce from (3.116) that

$$\lim_{\varepsilon \to 0} \varepsilon |y_\varepsilon| = 0. \tag{3.117}$$

Using the convexity of $f$ and that $y_\varepsilon$ is a minimizer of (3.115), we have that, for every $y \in \mathbb{R}^d$ and $\lambda \in (0, 1]$,

$$f(0, y_\varepsilon + y) - f(0, y_\varepsilon) \geq \lambda^{-1} \left( f(0, y_\varepsilon + \lambda y) - f(0, y_\varepsilon) \right)$$

$$\geq \lambda^{-1} \left( \lambda y \cdot p + \varepsilon \left( \sqrt{1 + |y_\varepsilon|^2} - \sqrt{1 + |y_\varepsilon + \lambda y|^2} \right) \right).$$

Letting $\lambda$ tend to zero shows that for every $y \in \mathbb{R}^d$,

$$f(0, y_\varepsilon + y) - f(0, y_\varepsilon) \geq y \cdot \left( p - \varepsilon \frac{y_\varepsilon}{\sqrt{1 + |y_\varepsilon|^2}} \right),$$

which means that

$$p_\varepsilon := p - \varepsilon \frac{y_\varepsilon}{\sqrt{1 + |y_\varepsilon|^2}} \in \partial \psi(y_\varepsilon). \tag{3.118}$$

Invoking the assumption on the initial condition gives $b_\varepsilon \in \mathbb{R}$ with $(b_\varepsilon, p_\varepsilon) \in \partial f(0, y_\varepsilon)$ and
$$b_\varepsilon - \mathsf{H}(p_\varepsilon) \geq 0. \tag{3.119}$$
Since $(b_\varepsilon, p_\varepsilon) \in \partial f(0, y_\varepsilon)$ and $(a, p) \in \partial f(t^*, x^*)$, we have
$$f(t^*, x^*) \geq f(0, y_\varepsilon) + b_\varepsilon t^* + p_\varepsilon \cdot (x^* - y_\varepsilon),$$
$$f(0, y_\varepsilon) \geq f(t^*, x^*) - at^* + p \cdot (y_\varepsilon - x^*).$$
Combining these two inequalities reveals that
$$b_\varepsilon t^* + p_\varepsilon \cdot (x^* - y_\varepsilon) \leq at^* + p \cdot (x^* - y_\varepsilon).$$
Using (3.117) and that $|p_\varepsilon - p| \leq \varepsilon$ yields that $b_\varepsilon \leq a + o_\varepsilon(1)$, and using also (3.119), we obtain that $a - \mathsf{H}(p_\varepsilon) \geq o_\varepsilon(1)$. Since $\mathsf{H}$ is continuous, this establishes (3.112) and completes the proof. ∎

*Proof of Theorem 3.21.* We write $\psi := f(0, \cdot)$. By Lemma 3.23, it suffices to show that for every $x \in \mathbb{R}^d$ and $p \in \partial \psi(x)$, there exists $b \in \mathbb{R}$ such that $(b, p) \in \partial f(0, x)$ and $b - \mathsf{H}(p) \geq 0$. Invoking Lemma 3.22 gives $(b, q) \in \partial f(0, x)$ with $b - \mathsf{H}(q) = 0$. Since $\psi \in C^1(\mathbb{R}^d; \mathbb{R})$, Theorem 2.13 implies that $p = q = \nabla \psi(x)$. This completes the proof. ∎

To apply the convex selection principle to the generalized Curie-Weiss model, it will be convenient to observe that the set of points $(t^*, x^*) \in \mathbb{R}_{\geq 0} \times \mathbb{R}^d$ for which there exists a smooth function that touches $f$ from above at $(t^*, x^*)$ is dense. This will allow us to replace the assumption that $f$ satisfies the Hamilton-Jacobi equation (3.20) on a dense set by the assumption that whenever a smooth function $\phi$ touches the function $f$ from above at a point $(t^*, x^*)$, the function $\phi$ should satisfy the Hamilton-Jacobi equation (3.20) at the point of contact $(t^*, x^*)$. This version of the convex selection principle will be well-suited for the generalized Curie-Weiss model due to the bound (3.17); it is also useful in other contexts such as statistical inference.

**Corollary 3.24.** *Let $f : \mathbb{R}_{\geq 0} \times \mathbb{R}^d \to \mathbb{R}$ be a jointly convex and jointly Lipschitz continuous function. Suppose that for any point $(t^*, x^*) \in \mathbb{R}_{>0} \times \mathbb{R}^d$ and smooth function $\phi \in C^\infty(\mathbb{R}_{\geq 0} \times \mathbb{R}^d; \mathbb{R})$ such that $f - \phi$ has a strict local maximum at $(t^*, x^*)$, we have $(\partial_t \phi - \mathsf{H}(\nabla \phi))(t^*, x^*) = 0$. If moreover $f(0, \cdot) \in C^1(\mathbb{R}^d; \mathbb{R})$, then $f$ is a viscosity solution to the Hamilton-Jacobi equation (3.20).*

*Proof.* Consider a smooth $\phi \in C^\infty(\mathbb{R}_{>0} \times \mathbb{R}^d; \mathbb{R})$ with the property that $f - \phi$ admits a strict local maximum at the point $(t^*, x^*) \in \mathbb{R}_{>0} \times \mathbb{R}^d$. Following Step 1 of the proof of Lemma 3.23, we see that $f$ is differentiable at the contact point $(t^*, x^*)$ with $(\partial_t f, \nabla f)(t^*, x^*) = (\partial_t \phi, \nabla \phi)(t^*, x^*)$. By the convex selection principle in Theorem 3.21, it therefore suffices to show that the set

$$\mathcal{A} := \{(t^*, x^*) \in \mathbb{R}_{>0} \times \mathbb{R}^d \mid \text{there exists } \phi \in C^\infty(\mathbb{R}_{>0} \times \mathbb{R}^d; \mathbb{R}) \text{ such that}$$
$$f - \phi \text{ has a strict local maximum at } (t^*, x^*)\}$$

is dense in $\mathbb{R}_{>0} \times \mathbb{R}^d$. Fix $(t_0, x_0) \in \mathbb{R}_{>0} \times \mathbb{R}^d$, and let $V$ be a compact neighbourhood of $(t_0, x_0)$. For each $\alpha \geq 1$, consider the mapping

$$\phi_\alpha(t, x) := f(t, x) - \frac{\alpha}{2}\left(|t - t_0|^2 + |x - x_0|^2\right),$$

and denote by $(t_\alpha, x_\alpha)$ a maximizer of $\phi_\alpha$ on $V$. Writing $L$ for the joint Lipschitz constant of $f$ and rearranging the bound

$$f(t_\alpha, x_\alpha) - \frac{\alpha}{2}\left(|t_\alpha - t_0|^2 + |x_\alpha - x_0|^2\right) = \phi_\alpha(t_\alpha, x_\alpha) \geq \phi_\alpha(t_0, x_0) = f(t_0, x_0)$$

reveals that

$$\left(|t_\alpha - t_0| + |x_\alpha - x_0|\right)^2 \leq \frac{2}{\alpha}|f(t_\alpha, x_\alpha) - f(t_0, x_0)| \leq \frac{2L}{\alpha}\left(|t_\alpha - t_0| + |x_\alpha - x_0|\right).$$

It follows that

$$|t_\alpha - t_0| + |x_\alpha - x_0| \leq \frac{2L}{\alpha},$$

so the sequence $(t_\alpha, x_\alpha)_{\alpha \geq 1}$ tends to $(t_0, x_0)$ as $\alpha$ tends to infinity. In particular, for $\alpha$ sufficiently large, the point $(t_\alpha, x_\alpha)$ is in the interior of $V$, and so is a local maximum of $\phi_\alpha$ as a function on $\mathbb{R}_{\geq 0} \times \mathbb{R}^d$. It is a strict local maximum for the mapping

$$(t, x) \mapsto \phi_\alpha(t, x) - \left(|t - t_\alpha|^2 + |x - x_\alpha|^2\right).$$

We have thus shown that for every $\alpha$ sufficiently large, the point $(t_\alpha, x_\alpha)$ belongs to $\mathcal{A}$. Recalling that $(t_\alpha, x_\alpha)$ tends to $(t_0, x_0)$ as $\alpha$ tends to infinity, we conclude that the set $\mathcal{A}$ is dense in $\mathbb{R}_{>0} \times \mathbb{R}^d$, as desired. ∎

We are finally in a position to use the Hamilton-Jacobi approach to prove Theorem 2.19 on the limit free energy in the generalized Curie-Weiss model.

**Theorem 3.25.** *For each integer $N \geq 1$, let $F_N : \mathbb{R}_{\geq 0} \times \mathbb{R} \to \mathbb{R}$ denote the free energy* (3.2) *in the generalized Curie-Weiss model. Suppose that for every $h \in \mathbb{R}$, the limit*

$$\psi(h) := \lim_{N \to +\infty} F_N(0, h) \qquad (3.120)$$

*exists, and that $\psi \in C^1(\mathbb{R}; \mathbb{R})$. Then, for every $t \geq 0$ and $h \in \mathbb{R}$, we have that $F_N(t, h)$ converges to $f(t, h)$ as $N$ tends to infinity, where $f : \mathbb{R}_{\geq 0} \times \mathbb{R} \to \mathbb{R}$ is the viscosity solution to the Hamilton-Jacobi equation* (3.19) *with initial condition $\psi$. Moreover, this function admits the Hopf representation*

$$f(t, h) = \lim_{N \to +\infty} F_N(t, h) = \sup_{m \in [-1, 1]} \left(t\xi(m) + hm - \psi^*(m)\right). \qquad (3.121)$$

*Proof.* The sequence $(F_N)_{N \geq 1}$ of free energies is precompact by the Arzelà-Ascoli theorem and the Lipschitz bounds (3.13). We denote by $f$ any subsequential limit, and observe that it must be jointly Lipschitz continuous by (3.13), and jointly convex by the joint convexity of each free energy $F_N$. We may therefore use the convex selection principle in Corollary 3.24 to prove that $f$ is a viscosity solution to the Hamilton-Jacobi equation (3.19). Fix $\phi \in C^\infty(\mathbb{R}_{>0} \times \mathbb{R}^d; \mathbb{R})$ with the property that $f - \phi$ has a strict local maximum at the point $(t^*, h^*) \in \mathbb{R}_{>0} \times \mathbb{R}$. Since $(F_N)_{N \geq 1}$ converges to $f$ with respect to the topology of local uniform convergence, Exercise 3.1 allows us to find a sequence $(t_N, h_N)_{N \geq 1}$ converging to $(t^*, h^*)$ such that $F_N - \phi$ has a local maximum at $(t_N, h_N)$ for each integer $N \geq 1$. The bound (3.17) gives a constant $C < +\infty$ with

$$\left|\partial_t F_N(t_N, x_N) - \xi\left(\partial_h F_N(t_N, x_N)\right)\right| \leq \frac{C}{N} \partial_h^2 F_N(t_N, x_N).$$

At the local maximum $(t_N, h_N)$ we have

$$\partial_t(F_N - \phi)(t_N, x_N) = 0, \quad \partial_h(F_N - \phi)(t_N, x_N) = 0, \quad \partial_h^2(F_N - \phi)(t_N, x_N) \leq 0,$$

from which it follows that

$$\left|\partial_t \phi(t_N, x_N) - \xi\left(\partial_h \phi(t_N, x_N)\right)\right| \leq \frac{C}{N} \partial_h^2 \phi(t_N, h_N).$$

Letting $N$ tend to infinity reveals that $\left(\partial_t \phi - \xi(\partial_h \phi)\right)(t^*, h^*) = 0$, where we have used the smoothness of $\phi$. Invoking the convex selection principle in Corollary 3.24 shows that $f$ is a viscosity solution to the Hamilton-Jacobi equation (3.19). Together with the uniqueness result in Corollary 3.7, this implies that $(F_N)_{N \geq 1}$ converges to the unique viscosity solution $f$ to the Hamilton-Jacobi equation (3.19), and arguing as in Proposition 3.20 gives the Hopf variational representation (3.121) of $f$. This completes the proof. ∎

We have thus been able to recover the key result from Section 2.3 obtained using large deviation principles. It is worth mentioning that Theorem 3.25 essentially allows us to deduce a large-deviation principle as a consequence. Indeed, if we consider a smooth function $\xi$ that is essentially zero near a point $m_0 \in \mathbb{R}$, and essentially $-\infty$ elsewhere, then $F_N(t, 0)$ is essentially computing the probability that the mean magnetization is around $m_0$, and the formula (3.121) states that this is essentially $\exp(-N\psi^*(m_0))$.

We have thus presented two methods, one based on large deviation principles and one based on the Hamilton-Jacobi approach, to identify the limit free energy in the generalized Curie-Weiss model. We do not know of another approach for doing this.

We close this chapter with a counterexample that shows that the differentiability assumption $\psi \in C^1(\mathbb{R}^d; \mathbb{R})$ in the convex selection principle cannot be dropped in general, despite the result in Proposition 3.20.

**Example 3.26.** Consider the probability measure

$$P_N := \frac{1}{2}\delta_{(1,\ldots,1)} + \frac{1}{2}\delta_{(-1,\ldots,-1)}$$

and the non-linearity $\xi(p) := -p^2$. The generalized Curie-Weiss free energy (3.2) associated with this non-linearity and this probability measure is given by

$$F_N(t,h) = \frac{1}{N}\log\left(\frac{1}{2}e^{-tN+hN} + \frac{1}{2}e^{-tN-hN}\right).$$

By Exercise 1.1, we have

$$\lim_{N\to+\infty} F_N(t,h) = |h| - t. \tag{3.122}$$

In particular, the limiting initial condition $\psi(h) := |h|$ is not in $C^1(\mathbb{R};\mathbb{R})$. On the other hand, the unique viscosity solution to the Hamilton-Jacobi equation

$$\partial_t f + \left(\partial_h f\right)^2 = 0 \tag{3.123}$$

with initial condition $\psi$ is given by the Hopf formula

$$f(t,h) = \sup_{m\in\mathbb{R}}\left(-tm^2 + hm - \psi^*(m)\right).$$

A direct computation shows that

$$\psi^*(m) = \sup_{p\in\mathbb{R}}(pm - |p|) = \begin{cases} 0 & \text{if } |m| \leq 1, \\ +\infty & \text{if } |m| > 1. \end{cases}$$

This means that

$$f(t,h) = \sup_{|m|\leq 1}\left(-tm^2 + hm\right) = \begin{cases} \frac{h^2}{4t} & \text{if } |h| \leq 2t \\ |h| - t & \text{if } |h| > 2t. \end{cases}$$

Since $(t,h) \mapsto |h| - t$ also satisfies the equation (3.123) at every point of differentiability and has the same initial condition, we have thus obtained a counter-example to the convex selection principle when the differentiability assumption $\psi \in C^1(\mathbb{R}^d;\mathbb{R})$ is dropped. We also see that the limit free energy (3.122) is not the viscosity solution in this case. Interestingly, we can identify another reference measure $P_N$ such that the associated free energy does converge to $f$. Indeed, consider now the probability measure

$$P_N := \frac{1}{2}\int_{-1}^{1}\delta_{(u,\ldots,u)}\,du.$$

In words, if we pick a uniform random variable $U$ over $[-1,1]$, then the law of the vector $(U,\ldots,U) \in \mathbb{R}^N$ is $P_N$. In this case, and still with $\xi(p) = -p^2$, we have

$$F_N(t,h) = \frac{1}{N} \log \frac{1}{2} \int_{-1}^{1} \exp N\left(-tu^2 + hu\right) du,$$

which, as $N$ tends to infinity, converges to

$$\sup_{u \in [-1,1]} \left(-tu^2 + hu\right) = f(t,h).$$

In particular, we see with these examples that, in situations in which the initial condition $\psi$ is not in $C^1(\mathbb{R};\mathbb{R})$, it is not possible in general to identify the limit free energy from the knowledge of $\psi := \lim_{N \to +\infty} F_N(0,\cdot)$ only.

# Chapter 4
# Statistical inference

So far, we have focused on models coming from statistical mechanics, and we have developed the Hamilton-Jacobi approach to compute the limit of their associated free energy. It turns out that certain problems of statistical inference share a very similar structure. In this chapter we apply the Hamilton-Jacobi approach to the problem of recovering information about a symmetric rank-one matrix given a noisy observation of it. To be more precise, we assume that, for some vector $\bar{x} = (\bar{x}_1, \ldots, \bar{x}_N)$ with independent and identically distributed coordinates, we observe $\bar{x}\bar{x}^* + W$, where the noise term $W$ is made of independent Gaussian random variables, independent from $\bar{x}$. Here and throughout this chapter, the superscript $*$ denotes the transposition operator. Our goal will be to study whether it is possible to recover meaningful information about the signal $\bar{x}\bar{x}^*$ from this noisy observation, in the regime of large $N$. The answer to this question will depend on the magnitude of the signal-to-noise ratio. When $\mathbb{E}\bar{x}_1 = 0$, it turns out that there exists a strictly positive and finite threshold such that if the signal-to-noise ratio is below this value, then essentially no information can be reconstructed about the signal from the observation, while if it is above this value, then there exists an estimator that has a non-trivial correlation with the signal itself (see Exercise 4.9 for a precise statement).

This critical threshold is analogous to the critical inverse temperature at which phase transitions occur in the Ising and Curie-Weiss models, and we will determine it using the techniques presented in the previous chapters and progressively developed in [69, 71, 72, 75, 193, 194]. In particular, we will appeal to the convex selection principle introduced in Section 3.6 to identify the limit free energy of this problem. We point out that, for this particular problem of statistical inference, a more standard approach similar to the proof of Proposition 3.20 would also work. We prefer to present an approach based on the convex selection principle since it generalizes to a much broader class of problems [75]. Alternative methods that so far have not reached the level of generality of [75] include [32, 33, 34, 35, 95, 96, 109, 148, 161, 162, 163, 168, 169, 172, 182, 227, 228].

In Section 4.1, we link the problem of statistical inference with the framework developed in earlier chapters, using the language of Gibbs measures and free energies to encode the relevant information-theoretic quantities. Section 4.2 contains a general discussion on Gaussian integration by parts and concentration results. We then combine these results with the Hamilton-Jacobi approach in Section 4.3 to identify the limit free energy. The critical signal-to-noise ratio is determined by analyzing the derivatives of the limit free energy as in Section 2.4. In Section 4.4, we test the performance of the classical principal-component analysis estimator for the symmetric rank-one matrix estimation problem, by comparing its associated mean-square error with the minimal mean-square error determined in Section 4.3. Finally, in Section 4.5 we see how, from an information-theoretic perspective, the symmetric rank-one matrix estimation problem is equivalent to the problem of detecting community structures in a dense random graph for which the probability of two nodes being connected depends only on the community of each node.

## 4.1  From statistical inference to statistical mechanics

We give ourselves a random vector $\bar{x} := (\bar{x}_1, \ldots, \bar{x}_N) \in \mathbb{R}^N$ with independent coordinates sampled from a bounded probability measure $P_1$ on the real line. We write $P_N := (P_1)^{\otimes N}$ to denote the law of $\bar{x}$. The problem of symmetric rank-one matrix estimation that we consider consists in observing

$$Y := \sqrt{\frac{2t}{N}} \bar{x}\bar{x}^* + W \in \mathbb{R}^{N \times N}, \qquad (4.1)$$

where $W = (W_{ij})_{1 \leq i,j \leq N} \in \mathbb{R}^{N \times N}$ is made of independent standard Gaussian random variables independent of the vector $\bar{x}$, and we think of it as noise. We call the parameter $t \geq 0$ the *signal-to-noise ratio*. Our main focus will be on understanding how much information we can recover about the symmetric rank-one matrix $\bar{x}\bar{x}^*$ from its noisy observation $Y$. In particular, we will aim to determine the asymptotic behaviour, as $N$ tends to infinity, of the *minimal mean-square error*,

$$\mathsf{mmse}_N(t) := \frac{1}{N^2} \inf_g \mathbb{E}|\bar{x}\bar{x}^* - g(Y)|^2 = \frac{1}{N^2} \mathbb{E}|\bar{x}\bar{x}^* - \mathbb{E}[\bar{x}\bar{x}^* \mid Y]|^2, \qquad (4.2)$$

between the signal $\bar{x}\bar{x}^*$ and its noisy observation $Y$. Here, the infimum is taken over the set of measurable functions $g$, and for any two matrices $a, b$ of the same size, we write

$$a \cdot b := \mathrm{tr}(ab^*) \quad \text{and} \quad |a| := (a \cdot a)^{\frac{1}{2}} \qquad (4.3)$$

for their entry-wise scalar product and their Euclidean norm, respectively. In the $L^2$ sense of (4.2), the best estimator for recovering $\bar{x}\bar{x}^*$ is the conditional expectation of this quantity given $Y$. Notice that in order to compute this conditional expectation, one needs to know the laws of $\bar{x}$ and $W$, so we implicitly assume that these laws are known to the observer. We will also be interested in computing the *mutual information*,

$$I_N(t) := \frac{1}{N}\mathbb{E}\int_{\mathbb{R}^N} \log\left(\frac{dP_{\bar{x}|Y}}{dP_N}(x)\right) dP_{\bar{x}|Y}(x), \tag{4.4}$$

between the vector $\bar{x}$ and the observation $Y$. Here $P_{\bar{x}|Y}$ denotes the conditional law of $\bar{x}$ given $Y$ (this is a random object that depends on the realization of $Y$; in other words, it is the evaluation at $Y$ of the mapping which sends each $y$ to the conditional law of $\bar{x}$ given that $Y = y$). The quantity $I_N(t)$ is a precise measure of the amount of information obtained about the vector $\bar{x}$ when we observe $Y$; see for instance [82] for more on this interpretation, and Exercise 4.3.

In order to analyze quantities such as the minimal mean-square error (4.2) or the mutual information (4.4), it will be useful to get a better understanding of the conditional law of the signal $\bar{x}$ given the observation $Y$. For every $x \in \mathbb{R}^N$ and $y \in \mathbb{R}^{N\times N}$, a formal calculation suggests that

$$\mathbb{P}\{\bar{x} = x \mid Y = y\} = \frac{\mathbb{P}\{\bar{x} = x \text{ and } Y = y\}}{\mathbb{P}\{Y = y\}} = \frac{\exp\left(-\frac{1}{2}\left|y - \sqrt{\frac{2t}{N}}xx^*\right|^2\right) dP_N(x)}{\int_{\mathbb{R}^N} \exp\left(-\frac{1}{2}\left|y - \sqrt{\frac{2t}{N}}x'x'^*\right|^2\right) dP_N(x')}. \tag{4.5}$$

More precisely, denoting

$$H_N^\circ(t,x) := \sqrt{\frac{2t}{N}} x \cdot Yx - \frac{t}{N}|x|^4 = \sqrt{\frac{2t}{N}} x \cdot Wx + \frac{2t}{N}(x \cdot \bar{x})^2 - \frac{t}{N}|x|^4, \tag{4.6}$$

we have for every bounded measurable function $f : \mathbb{R}^N \to \mathbb{R}$ that

$$\mathbb{E}[f(\bar{x}) \mid Y] = \frac{\int_{\mathbb{R}^N} f(x) \exp H_N^\circ(t,x) \, dP_N(x)}{\int_{\mathbb{R}^N} \exp H_N^\circ(t,x) \, dP_N(x)}. \tag{4.7}$$

This is verified rigorously in Exercise 4.2. In other words, the conditional law of $\bar{x}$ given $Y$ is the Gibbs measure associated with the Hamiltonian $H_N^\circ(t,\cdot)$. Notice also that the term $x \cdot Wx$ in (4.6) is highly reminiscent of the quantity (0.1) appearing in the definition of the Sherrington-Kirkpatrick spin-glass model. As in previous chapters, it will be convenient to focus first on studying the asymptotic behaviour of the free energy

$$F_N^\circ(t) := \frac{1}{N} \log \int_{\mathbb{R}^N} \exp H_N^\circ(t,x) \, dP_N(x). \tag{4.8}$$

Compared with the Curie-Weiss model, one important novelty in the symmetric rank-one matrix estimation problem is that the Hamiltonian $H_N^\circ(t,\cdot)$ is itself a random quantity, as it depends on $\bar{x}$ and $W$ through $Y$. In particular, the free energy (4.8) is also random, and we denote its average by

$$\overline{F}_N^\circ(t) := \mathbb{E} F_N^\circ(t) = \frac{1}{N} \mathbb{E} \log \int_{\mathbb{R}^N} \exp H_N^\circ(t,x) \, dP_N(x). \tag{4.9}$$

We will often slightly abuse terminology and refer to the average free energy (4.9) as simply the free energy. Just like in the Curie-Weiss model, we introduce notation for a random variable whose law is the random Gibbs measure (4.7). For every bounded and measurable function $f : \mathbb{R}^N \to \mathbb{R}$, we write

$$\langle f(x) \rangle := \int_{\mathbb{R}^N} f(x) \, dP_{\bar{x}|Y}(x) = \frac{\int_{\mathbb{R}^N} f(x) \exp H_N^\circ(t,x) \, dP_N(x)}{\int_{\mathbb{R}^N} \exp H_N^\circ(t,x) \, dP_N(x)}. \tag{4.10}$$

Although this is kept implicit in the notation, we stress that the bracket $\langle \cdot \rangle$ is a random quantity that depends on $t$ and on $Y$. Whenever we write expressions such as $\langle g(x,\bar{x}) \rangle$, we understand that the variable $x$ is integrated against the conditional probability measure $P_{\bar{x}|Y}$, while keeping $\bar{x}$ fixed. In more explicit notation,

$$\langle g(x,\bar{x}) \rangle = \int_{\mathbb{R}^N} g(x,\bar{x}) \, dP_{\bar{x}|Y}(x) = \frac{\int_{\mathbb{R}^N} g(x,\bar{x}) \exp H_N^\circ(t,x) \, dP_N(x)}{\int_{\mathbb{R}^N} \exp H_N^\circ(t,x) \, dP_N(x)}, \tag{4.11}$$

not to be confused with $\mathbb{E}[g(\bar{x},\bar{x}) \mid Y]$ for instance. If we have no information about the signal $\bar{x}$, when $t = 0$, then $x$ is simply an independent copy of $\bar{x}$. On the other hand, if we have perfect information about $\bar{x}$, for instance if we had instead observed $\bar{x}$ with no noise, then $x$ would be equal to $\bar{x}$.

It will also be convenient to introduce independent copies of $x$ under the Gibbs average $\langle \cdot \rangle$, often called *replicas*, which we denote by $x'$, $x''$, or also $x^1$, $x^2$, $x^3$, and so on if an arbitrary number of replicas needs to be considered. Explicitly, for every bounded and measurable function $f : \mathbb{R}^N \times \mathbb{R}^N \to \mathbb{R}$,

$$\langle f(x,x') \rangle = \frac{\int_{\mathbb{R}^N} \int_{\mathbb{R}^N} f(x,x') \exp \left( H_N^\circ(t,x) + H_N^\circ(t,x') \right) dP_N(x) \, dP_N(x')}{\left( \int_{\mathbb{R}^N} \exp H_N^\circ(t,x) \, dP_N(x) \right)^2}, \tag{4.12}$$

with the natural generalization of this expression in the case of more replicas. Compared with the setting of spin glasses explored in Chapter 6, the fact that the Gibbs measure (4.7) is a conditional expectation will fundamentally simplify the analysis. In the language of physics, the symmetric rank-one matrix estimation problem is always *replica symmetric* — this terminology will be clarified in Chapter 6. The replica symmetry will be derived from the Nishimori identity (named after [202]). This identity allows us to replace one replica $x$ by the ground-truth signal $\bar{x}$, provided that we average over all sources of randomness.

**Proposition 4.1** (Nishimori identity). *For all bounded and measurable functions $f : \mathbb{R}^N \times \mathbb{R}^{N \times N} \to \mathbb{R}$ and $g : \mathbb{R}^N \times \mathbb{R}^N \times \mathbb{R}^{N \times N} \to \mathbb{R}$, we have*

$$\mathbb{E}\langle f(x, Y)\rangle = \mathbb{E}\langle f(\bar{x}, Y)\rangle, \tag{4.13}$$

$$\mathbb{E}\langle g(x, x', Y)\rangle = \mathbb{E}\langle g(x, \bar{x}, Y)\rangle, \tag{4.14}$$

*and so on with more replicas, that is, for every integer $\ell \geq 1$ and bounded measurable function $h : (\mathbb{R}^N)^\ell \times \mathbb{R}^{N \times N} \to \mathbb{R}$,*

$$\mathbb{E}\langle h(x^1, x^2, \ldots, x^\ell, Y)\rangle = \mathbb{E}\langle h(x^1, x^2, \ldots, x^{\ell-1}, \bar{x}, Y)\rangle. \tag{4.15}$$

**Remark 4.2.** The replacement of a replica by the ground-truth signal $\bar{x}$ can only be done once: the terms in (4.14) are in general different from $\mathbb{E}\langle g(\bar{x}, \bar{x}, Y)\rangle$. The double average is also necessary, as the terms $\langle g(x, x', Y)\rangle$ and $\langle g(x, \bar{x}, Y)\rangle$ are not equal in general.

*Proof of Proposition 4.1.* It clearly suffices to prove (4.15). By Dynkin's $\pi$-$\lambda$ theorem (see Theorem A.5 and Exercise A.3), it suffices to verify this identity for functions that factorize over the variables. Precisely, we assume that the function $h$ can be written in the form

$$h(x^1, \ldots, x^\ell, Y) = h_1(x^1, \ldots, x^{\ell-1}) h_2(x^\ell) h_3(Y) \tag{4.16}$$

for some bounded measurable functions $h_1, h_2, h_3$. We can then write

$$\mathbb{E}\langle h(x^1, \ldots, x^\ell, Y)\rangle = \mathbb{E}\big(\langle h_1(x^1, \ldots, x^{\ell-1}, Y)\rangle \langle h_2(x^\ell)\rangle h_3(Y)\big)$$
$$= \mathbb{E}\big(\langle h_1(x^1, \ldots, x^{\ell-1}, Y)\rangle \mathbb{E}[h_2(\bar{x}) \mid Y] h_3(Y)\big),$$

where we used that the Gibbs measure (4.10) is the conditional law of $\bar{x}$ given the observation $Y$ in the last identity. Recalling that the measure $\langle \cdot \rangle$ depends on the randomness only through $Y$, and thus that $\langle h_1(x^1, \ldots, x^{\ell-1}, Y)\rangle$ is $Y$-measurable, we obtain that

$$\mathbb{E}\langle h(x^1, \ldots, x^\ell, Y)\rangle = \mathbb{E}\big(\mathbb{E}[\langle h_1(x^1, \ldots, x^{\ell-1}, Y)\rangle h_2(\bar{x}) h_3(Y) \mid Y]\big)$$
$$= \mathbb{E}\big(\langle h_1(x^1, \ldots, x^{\ell-1}, Y)\rangle h_2(\bar{x}) h_3(Y)\big)$$
$$= \mathbb{E}\langle h_1(x^1, \ldots, x^{\ell-1}, Y) h_2(\bar{x}) h_3(Y)\rangle.$$

Combining this with the definition (4.16) of the function $h$ completes the proof. ∎

Using this identity, we can rewrite the minimal mean-square error (4.2) and the mutual information (4.4) in terms of the Gibbs average (4.10) and the free energy (4.9).

**Proposition 4.3.** *For every $t \geq 0$ and integer $N \geq 1$, the minimal mean-square error* (4.2) *is given by*

$$\mathsf{mmse}_N(t) = \frac{1}{N^2}\mathbb{E}|\bar{x}|^4 - \frac{1}{N^2}\mathbb{E}\langle (x\cdot\bar{x})^2\rangle. \tag{4.17}$$

*Proof.* We expand the square to find that

$$N^2\mathsf{mmse}_N(t) = \mathbb{E}|\bar{x}\bar{x}^* - \langle xx^*\rangle|^2 = \mathbb{E}|\bar{x}|^4 - 2\mathbb{E}(\bar{x}\bar{x}^* \cdot \langle xx^*\rangle) + \mathbb{E}(\langle xx^*\rangle \cdot \langle xx^*\rangle),$$

and we can rewrite

$$\langle xx^*\rangle \cdot \langle xx^*\rangle = \langle x'x'^* \cdot xx^*\rangle.$$

The Nishimori identity (4.15) yields that

$$N^2\mathsf{mmse}_N(t) = \mathbb{E}|\bar{x}|^4 - \mathbb{E}\langle \bar{x}\bar{x}^* \cdot xx^*\rangle = \mathbb{E}|\bar{x}|^4 - \mathbb{E}\langle (x\cdot\bar{x})^2\rangle.$$

This completes the proof. ∎

**Proposition 4.4.** *For every $t \geq 0$ and integer $N \geq 1$, the mutual information* (4.4) *is given by*

$$\mathsf{I}_N(t) = \frac{t}{N^2}\mathbb{E}|\bar{x}|^4 - \overline{F}^\circ_N(t). \tag{4.18}$$

*Proof.* Remembering that the conditional distribution of $\bar{x}$ given $Y$ is the Gibbs measure (4.7), we can write the mutual information in terms of the Gibbs average (4.10) and the free energy (4.9) as

$$\mathsf{I}_N(t) = \frac{1}{N}\mathbb{E}\left\langle \log\left(\frac{\exp H^\circ_N(t,x)}{\int_{\mathbb{R}^N} \exp H^\circ_N(t,x')\,dP_N(x')}\right)\right\rangle = \frac{1}{N}\mathbb{E}\langle H^\circ_N(t,x)\rangle - \overline{F}^\circ_N(t).$$

Invoking the Nishimori identity (4.15) reveals that

$$\mathbb{E}\langle H^\circ_N(t,x)\rangle = \mathbb{E}H^\circ_N(t,\bar{x}) = \sqrt{\frac{2t}{N}}\mathbb{E}\bar{x}\cdot Y\bar{x} - \frac{t}{N}\mathbb{E}|\bar{x}|^4 = \frac{2t}{N}\mathbb{E}|\bar{x}|^4 - \frac{t}{N}\mathbb{E}|\bar{x}|^4 = \frac{t}{N}\mathbb{E}|\bar{x}|^4.$$

This completes the proof. ∎

Together with the fact that the coordinates of $\bar{x}$ are independent and identically distributed, these results imply that understanding the minimal mean-square error (4.2) and the mutual information (4.4) comes down to understanding the asymptotic behaviour of the quantity $\mathbb{E}\langle (x\cdot\bar{x})^2\rangle$ and of the free energy (4.9). Indeed, we have

$$\frac{1}{N^2}\mathbb{E}|\bar{x}|^4 = \left(\mathbb{E}|\bar{x}_1|^2\right)^2 + \frac{1}{N}\left(\mathbb{E}|\bar{x}_1|^4 - \left(\mathbb{E}|\bar{x}_1|^2\right)^2\right) = \left(\mathbb{E}|\bar{x}_1|^2\right)^2 + o(1). \tag{4.19}$$

To understand the limit of the free energy (4.9), we will use the Hamilton-Jacobi approach. There are two issues that we need to resolve to be able to do this. The first is that we need to be able to simplify the time derivative of the free energy (4.9),

$$\partial_t \overline{F}_N^\circ(t) = \sqrt{\frac{1}{2tN^3}} \mathbb{E}\langle x \cdot Wx \rangle + \frac{2}{N^2}\mathbb{E}\langle (x \cdot \bar{x})^2 \rangle - \frac{1}{N^2}\mathbb{E}\langle |x|^4 \rangle, \quad (4.20)$$

by integrating out the Gaussian noise $W$; the second is that we need to "enrich" the free energy (4.9) so that it depends on an additional variable that allows us to close the equation. The first of these issues will be addressed through the Gaussian integration by parts formula discussed in the next section, where we will show that

$$\partial_t \overline{F}_N^\circ(t) = \frac{1}{N^2}\mathbb{E}\langle (x \cdot \bar{x})^2 \rangle. \quad (4.21)$$

The second issue will be tackled in Section 4.3 by enriching the observation $Y$. The equality (4.21) reveals that determining the limit free energy is in fact sufficient to understand both the asymptotic mutual information and the asymptotic minimal mean-square error [135].

**Exercise 4.1.** Let $A \in \mathbb{R}^{d \times d}$ be a symmetric rank-one matrix. Show that there exists a unit vector $\bar{u} \in \mathbb{R}^d$, unique up to a sign, and a unique scalar $\lambda \in \mathbb{R}$ with $A = \lambda \bar{u}\bar{u}^*$.

**Exercise 4.2.** Show that in the context of the symmetric rank-one matrix estimation problem, the conditional law of $\bar{x}$ given $Y$ is the Gibbs measure (4.7) associated with the Hamiltonian (4.6).

**Exercise 4.3.** Fix two probability measures $\mathbb{P}$ and $\mathbb{Q}$ on a separable metric space $S$, and write $\mathbb{P} \ll \mathbb{Q}$ to mean that $\mathbb{P}$ is absolutely continuous with respect to $\mathbb{Q}$. The *relative entropy* of $\mathbb{P}$ with respect to $\mathbb{Q}$ is defined by

$$\mathrm{H}(\mathbb{P} \mid \mathbb{Q}) := \begin{cases} \int_S \frac{d\mathbb{P}}{d\mathbb{Q}} \log\left(\frac{d\mathbb{P}}{d\mathbb{Q}}\right) d\mathbb{Q} & \text{if } \mathbb{P} \ll \mathbb{Q}, \\ +\infty & \text{otherwise.} \end{cases} \quad (4.22)$$

The *mutual information* between two random variables $X$ and $Z$ on $S$ is defined by

$$\mathrm{I}(X; Z) := \mathrm{H}(\mathbb{P}_{X,Z} \mid \mathbb{P}_X \otimes \mathbb{P}_Z), \quad (4.23)$$

where $\mathbb{P}_V$ denotes the law of a random variable $V$. We now place ourselves in the context of the symmetric rank-one matrix estimation problem.

(i) Prove that the mutual information (4.4) is given by

$$\mathrm{I}_N(t) = \frac{1}{N}\mathrm{I}(\bar{x}; Y). \quad (4.24)$$

(ii) Prove that
$$\mathrm{H}(\mathbb{P}_{\bar{x},Y} \mid \mathbb{P}_{\bar{x},W}) = \mathrm{H}(\mathbb{P}_Y \mid \mathbb{P}_W) + \mathrm{I}(\bar{x};Y), \tag{4.25}$$

(iii) Show that
$$\overline{F}_N^\circ(t) = \frac{1}{N}\mathrm{H}(\mathbb{P}_Y \mid \mathbb{P}_W) \quad \text{and} \quad \mathrm{H}(\mathbb{P}_{\bar{x},Y} \mid \mathbb{P}_{\bar{x},W}) = \frac{t}{N}\mathbb{E}|\bar{x}|^4. \tag{4.26}$$

(iv) Deduce the relationship (4.18) between the mutual information (4.4) and the free energy (4.9).

## 4.2 Gaussian integration by parts and concentration inequalities

In this section, we temporarily move away from the symmetric rank-one matrix estimation problem to discuss two classical results about Gaussian random variables, namely integration by parts and concentration inequalities. In the context of statistical inference and spin glasses, the Gaussian integration by parts formula is very helpful to compute derivatives of the free energy, while the Gaussian concentration inequality ensures that the free energy concentrates around its average. In the next section we will see how to leverage these results to tackle the symmetric rank-one matrix estimation problem using the Hamilton-Jacobi approach.

We start by presenting the Gaussian integration by parts formula. In its simplest form, let $g$ be a centred Gaussian random variable with variance $v^2$, so with density

$$\varphi_v(x) := \frac{1}{\sqrt{2\pi v^2}} \exp\left(-\frac{x^2}{2v^2}\right). \tag{4.27}$$

For every differentiable and bounded function $F \in C^1(\mathbb{R};\mathbb{R})$, we have

$$\mathbb{E}gF(g) = \int_\mathbb{R} xF(x)\varphi_v(x)\,\mathrm{d}x = -v^2 \int_\mathbb{R} F(x)\varphi_v'(x)\,\mathrm{d}x, \tag{4.28}$$

so integrating by parts reveals that

$$\mathbb{E}gF(g) = v^2 \int_\mathbb{R} F'(x)\varphi_v(x)\,\mathrm{d}x = v^2\mathbb{E}F'(g) = \mathbb{E}g^2\,\mathbb{E}F'(g). \tag{4.29}$$

To see what this formula yields in the context of statistical inference, let us fix a pair $1 \leq i,j \leq N$, and denote $F(W,\bar{x}) := \langle x_i x_j \rangle$. By (4.29),

$$\mathbb{E}\langle W_{ij} x_i x_j \rangle = \mathbb{E}W_{ij}F(W,\bar{x}) = \mathbb{E}\partial_{W_{ij}}F(W,\bar{x}). \tag{4.30}$$

A direct computation shows that $\partial_{W_{ij}} F(W,\bar{x})$ is given by

$$\partial_{W_{ij}} \left( \frac{\int_{\mathbb{R}^N} x_i x_j \exp\left(\sqrt{\frac{2t}{N}} \sum_{\ell,\ell'=1}^N W_{\ell,\ell'} x_\ell x_{\ell'}\right) d\mu(x)}{\int_{\mathbb{R}^N} \exp\left(\sqrt{\frac{2t}{N}} \sum_{\ell,\ell'=1}^N W_{\ell,\ell'} x_\ell x_{\ell'}\right) d\mu(x)} \right)$$

$$= \sqrt{\frac{2t}{N}} \left( \langle x_i^2 x_j^2 \rangle - \langle x_i x_j \rangle^2 \right), \quad (4.31)$$

where we used the shorthand notation

$$d\mu(x) := \exp\left( \frac{2t}{N} (\bar{x} \cdot x)^2 - \frac{t}{N} |x|^4 \right) dP_N(x). \quad (4.32)$$

Summing over all $i$ and $j$, we thus obtain that

$$\mathbb{E}\langle x \cdot W x \rangle = \sqrt{\frac{2t}{N}} \left( \mathbb{E}\langle |x|^4 \rangle - \mathbb{E}\langle (x \cdot x')^2 \rangle \right), \quad (4.33)$$

which combined with (4.20) and the Nishimori identity yields (4.21).

When it comes to integration by parts, the identity (4.29) is all we need for the purposes of this chapter. With a view towards subsequent chapters on spin glasses, we now present a multivariate version of this Gaussian integration by parts formula, also allowing for more general functions $F$; and then also another variant that is tailored to the setting of Gibbs measures. Readers who are not interested in subsequent chapters on spin glasses can safely skip to the end of the proof of Theorem 4.6, at which place we turn to discussing concentration inequalities.

**Theorem 4.5** (Gaussian integration by parts). *Let $g = (g_1, \ldots, g_d)$ be a centred Gaussian vector in $\mathbb{R}^d$ and $F \in C^1(\mathbb{R}^d; \mathbb{R})$ be such that $\mathbb{E}|g_1 F(g)| + \mathbb{E}|\nabla F(g)| < +\infty$. We have*

$$\mathbb{E} g_1 F(g) = \sum_{i=1}^d \mathbb{E} g_1 g_i \, \mathbb{E} \partial_{x_i} F(g). \quad (4.34)$$

*Proof.* Denote by $C$ the covariance matrix of the Gaussian vector $g$. The proof proceeds in three steps. First we prove that when $C = I_d$, we have for any bounded function $G \in C^1(\mathbb{R}^d; \mathbb{R}^d)$ with $\mathbb{E}|\nabla G(g)| < +\infty$ that

$$\mathbb{E} g \cdot G(g) = \mathbb{E} \nabla \cdot C G(g). \quad (4.35)$$

Then we use a change of variables to deduce (4.35) for arbitrary $C$, and finally we prove (4.34) by choosing an appropriate $G$ and using an approximation argument.

*Step 1: proving* (4.35) *for $C = I_d$.* In the case when the covariance matrix $C$ is the identity, the Gaussian vector $g$ has density

$$\varphi(x) := \frac{1}{(2\pi)^d} \exp\left(-\frac{1}{2}|x|^2\right).$$

This implies that for any bounded function $G \in C^1(\mathbb{R}^d; \mathbb{R}^d)$ with $\mathbb{E}|\nabla G(g)| < +\infty$, we have

$$\mathbb{E} g \cdot G(g) = \int_{\mathbb{R}^d} x \cdot G(x) \varphi(x) \, dx = -\int_{\mathbb{R}^d} \nabla \varphi(x) \cdot G(x) \, dx,$$

where we used that $\nabla \varphi(x) = -x\varphi(x)$. Integrating by parts reveals that

$$\mathbb{E} g \cdot G(g) = \int_{\mathbb{R}^d} \varphi(x) \nabla \cdot G(x) \, dx = \mathbb{E} \nabla \cdot G(g).$$

The boundedness of $G$ ensures that the boundary term in this integration by parts vanishes. This establishes (4.35) when $C = I_d$.

*Step 2: deducing* (4.35) *for general* $C$. Letting $z$ be a standard Gaussian random vector in $\mathbb{R}^d$ and $g := C^{1/2}z$, the vector $g$ is centred Gaussian with covariance $C$. It follows by symmetry of $C^{1/2}$ and the previous step that

$$\mathbb{E} g \cdot G(g) = \mathbb{E} z \cdot C^{1/2} G(C^{1/2} z) = \mathbb{E} z \cdot \widetilde{G}(z)$$

for the function $\widetilde{G}(z) := C^{1/2} G(C^{1/2} z)$. A direct computation reveals that

$$\nabla \cdot \widetilde{G}(z) = \sum_{i,j,k=1}^d C^{1/2}_{ij} C^{1/2}_{ki} \partial_{x_k} G_j(g) = \sum_{j=1}^d \partial_{x_j}(CG)_j(g) = \nabla \cdot CG(g)$$

which gives (4.35) for general $C$.

*Step 3: obtaining* (4.34). Applying (4.35) to the function $G(g) := F(g)e_1$, where $e_1 \in \mathbb{R}^d$ denotes the first canonical basis vector, and observing that $(CG)_i = C_{i,1} F(g)$ for $1 \leq i \leq d$ establishes (4.34) for $F$ bounded. To lift this assumption, given $\varepsilon > 0$, define the function $\phi_\varepsilon : \mathbb{R} \to \mathbb{R}$ by

$$\phi_\varepsilon(x) := \frac{x}{\sqrt{1 + \varepsilon x^2}}.$$

Observe that $|\phi_\varepsilon|$ is bounded by $\varepsilon^{-1/2}$, so $\phi_\varepsilon \circ F$ is differentiable, bounded, and satisfies $\mathbb{E}|\nabla(\phi_\varepsilon \circ F)| < +\infty$ since $\|\nabla \phi_\varepsilon\|_\infty \leq 1$. It follows that

$$\mathbb{E} g_1 \phi_\varepsilon(F(g)) = \sum_{i=1}^d \mathbb{E} g_1 g_i \mathbb{E} \phi'_\varepsilon(F(g)) \partial_{x_i} F(g) = \sum_{i=1}^d \mathbb{E} g_1 g_i \mathbb{E} \frac{\partial_{x_i} F(g)}{(1 + \varepsilon F(g)^2)^{3/2}}.$$

Invoking the dominated convergence theorem to let $\varepsilon$ tend to zero completes the proof. ∎

As is suggested by the example in (4.33), the Gaussian integration by parts formula admits a convenient reformulation that is adapted to the context of Gibbs measures. With a view towards using it for spin glasses in subsequent chapters, we now state and prove it in a general context. We reiterate that the reader only interested in statistical inference models can safely skip this generalization and directly move to right after the proof of Theorem 4.6.

Consider a centred Gaussian process made of the two families $(y(\sigma))_{\sigma \in \Sigma}$ and $(z(\sigma))_{\sigma \in \Sigma}$ indexed by some complete and separable metric space of indices $\Sigma$, and let $\mu$ be a finite measure on $\Sigma$. Suppose that the covariances of $(y(\sigma))_{\sigma \in \Sigma}$ and $(z(\sigma))_{\sigma \in \Sigma}$ are given by

$$C_y(\sigma^1, \sigma^2) := \mathbb{E} y(\sigma^1) y(\sigma^2) \quad \text{and} \quad C_z(\sigma^1, \sigma^2) := \mathbb{E} z(\sigma^1) z(\sigma^2) \qquad (4.36)$$

for some bounded and continuous functions $C_y : \Sigma^2 \to \mathbb{R}$ and $C_z : \Sigma^2 \to \mathbb{R}$. The process $(y(\sigma))_{\sigma \in \Sigma}$ defines the Gibbs measure

$$dG(\sigma) := \frac{\exp y(\sigma) \, d\mu(\sigma)}{\int_\Sigma \exp y(\tau) \, d\mu(\tau)} \qquad (4.37)$$

on $\Sigma$, and we denote by $\langle \cdot \rangle$ the Gibbs average associated with the Gibbs measure $G$, with canonical random variable $\sigma$. We also let $(\sigma^\ell)_{\ell \geq 1}$ be independent copies of $\sigma$ under $\langle \cdot \rangle$. This means that for any integer $\ell \geq 1$ and any bounded and measurable function $f : \Sigma^\ell \to \mathbb{R}$,

$$\langle f(\sigma^1, \ldots, \sigma^\ell) \rangle := \frac{\int_{\Sigma^\ell} f(\sigma^1, \ldots, \sigma^\ell) \prod_{i=1}^\ell \exp y(\sigma^i) \, d\mu(\sigma^i)}{\left( \int_\Sigma \exp y(\sigma) \, d\mu(\sigma) \right)^\ell}. \qquad (4.38)$$

The regularity assumptions on the covariance functions (4.36) ensure that the integral defining the Gibbs average is well-defined. Indeed, the continuity of $C_y$ can be used to show that the process $(\exp y(\sigma))_{\sigma \in \Sigma}$ is continuous in mean and therefore in probability. In the context of the symmetric rank-one matrix estimation problem, we would have

$$\Sigma := \mathbb{R}^N, \ y(x) := \sqrt{\frac{2t}{N}} x \cdot W x, \text{ and } d\mu(x) = \exp\left( \frac{2t}{N} (\bar{x} \cdot x)^2 - \frac{t}{N} |x|^4 \right) dP_N(x), \quad (4.39)$$

and for the calculation leading to (4.33), we would choose $z(x) := x \cdot W x$. In general, the Gaussian integration by parts formula in Theorem 4.5 can be used to compute Gibbs averages of Gaussian processes as follows.

**Theorem 4.6** (Gibbs Gaussian integration by parts). *Let $C(\sigma^1, \sigma^2) := \mathbb{E} z(\sigma^1) y(\sigma^2)$ denote the covariance between the Gaussian processes $(z(\sigma))_{\sigma \in \Sigma}$ and $(y(\sigma))_{\sigma \in \Sigma}$. For every $n \geq 1$ and bounded measurable $\Phi : \Sigma^n \to \mathbb{R}$, we have*

$$\mathbb{E} \langle \Phi(\sigma^1, \ldots, \sigma^n) z(\sigma^1) \rangle = \mathbb{E} \Big\langle \Phi(\sigma^1, \ldots, \sigma^n) \Big( \sum_{i=1}^n C(\sigma^1, \sigma^i) - n C(\sigma^1, \sigma^{n+1}) \Big) \Big\rangle.$$
(4.40)

*Proof.* The definition of the Gibbs average (4.38) implies that

$$\mathbb{E}\langle \Phi(\sigma^1,\ldots,\sigma^n)z(\sigma^1)\rangle$$
$$= \int_{\Sigma^n} \Phi(\sigma^1,\ldots,\sigma^n)\mathbb{E}z(\sigma^1)\prod_{\ell=1}^n \frac{\exp y(\sigma^\ell)}{\int_\Sigma \exp y(\tau)\,d\mu(\tau)}\,d\mu(\sigma^\ell). \quad (4.41)$$

To simplify this expression, we fix $\sigma^1,\ldots,\sigma^n \in \Sigma$ and introduce the Gaussian process $(y'(\sigma))_{\sigma\in\Sigma}$ defined by

$$y'(\sigma) := y(\sigma) - \lambda_{\sigma^1,\sigma}z(\sigma^1) \quad \text{for} \quad \lambda_{\sigma^1,\sigma} := \frac{\mathbb{E}y(\sigma)z(\sigma^1)}{\mathbb{E}z(\sigma^1)^2}.$$

Since $\mathbb{E}y'(\sigma)z(\sigma^1) = 0$ and $y$ and $z$ are jointly Gaussian, the process $(y'(\sigma))_{\sigma\in\Sigma}$ is independent of $z(\sigma^1)$. If we condition on $(y'(\sigma))_{\sigma\in\Sigma}$ and define the function $F:\mathbb{R}\to\mathbb{R}$ by

$$F(x) := \prod_{\ell=1}^n \frac{\exp\bigl(y'(\sigma^\ell)+\lambda_{\sigma^1,\sigma^\ell}x\bigr)}{\int_\Sigma \exp\bigl(y'(\tau)+\lambda_{\sigma^1,\tau}x\bigr)\,d\mu(\tau)},$$

then the Gaussian integration by parts formula in Theorem 4.5 implies that

$$\mathbb{E}z(\sigma^1)\prod_{\ell=1}^n \frac{\exp y(\sigma^\ell)}{\int_\Sigma \exp y(\eta)\,d\mu(\eta)} = \mathbb{E}z(\sigma^1)^2 \mathbb{E}F'(z(\sigma^1)).$$

A direct computation shows that

$$F'(x) = \sum_{\ell=1}^n \lambda_{\sigma^1,\sigma^\ell}F(x) - nF(x)\frac{\int_\Sigma \lambda_{\sigma^1,\sigma^{n+1}}\exp(y'(\sigma^{n+1})+\lambda_{\sigma^1,\sigma^{n+1}}x)\,d\mu(\sigma^{n+1})}{\int_\Sigma \exp(y'(\tau)+\lambda_{\sigma^1,\tau}x)\,d\mu(\tau)}.$$

It follows that conditionally on $(y'(\sigma))_{\sigma\in\Sigma}$ as well as the replicas $\sigma^1,\ldots,\sigma^n \in \Sigma$,

$$\mathbb{E}z(\sigma^1)\prod_{\ell=1}^n \frac{\exp y(\sigma^\ell)}{\int_\Sigma \exp y(\tau)\,d\mu(\tau)}$$
$$= \sum_{\ell=1}^n C(\sigma^1,\sigma^\ell)\mathbb{E}F(z(\sigma^1)) - n\mathbb{E}F(z(\sigma^1))\langle C(\sigma^1,\sigma^{n+1})\rangle',$$

where $\langle\cdot\rangle'$ denotes the Gibbs average (4.38) conditionally on $\sigma^1,\ldots,\sigma^n$, and therefore with the average only taken with respect to the randomness of $\sigma^{n+1}$. Since $(y'(\sigma))_{\sigma\in\Sigma}$ is independent of $z(\sigma^1)$, this equality also holds unconditionally on the randomness of $(y'(\sigma))_{\sigma\in\Sigma}$. Substituting it into the right side of (4.41) yields

$$\mathbb{E}\int_{\Sigma^n} \Phi(\sigma^1,\ldots,\sigma^n)\left(\sum_{\ell=1}^n C(\sigma^1,\sigma^\ell) - n\langle C(\sigma^1,\sigma^{n+1})\rangle'\right)$$
$$\prod_{\ell=1}^n \frac{\exp y(\sigma^\ell)}{\int_\Sigma \exp y(\tau)\,d\mu(\tau)}\,d\mu(\sigma^\ell).$$

Remembering that the Gibbs average (4.38) is a product measure allows us to absorb the Gibbs average $\langle \cdot \rangle'$ into the Gibbs average $\langle \cdot \rangle$ to obtain (4.40). This completes the proof. ∎

We now turn to the Gaussian concentration inequality. This inequality states that a Lipschitz function of a standard Gaussian random vector will concentrate around its mean, with a standard deviation of the order of the Lipschitz constant. To motivate this statement, let us take $g$ a centred Gaussian random variable with variance $v^2$, and observe that by Chebyshev's inequality, we have for every $t, \lambda > 0$ that

$$\mathbb{P}\{|g| > t\} \leq 2\mathbb{P}\{g > t\} \leq 2\exp(-\lambda t)\mathbb{E}\exp(\lambda g) = 2\exp\left(\frac{\lambda^2 v^2}{2} - \lambda t\right). \quad (4.42)$$

Optimizing over $\lambda$ to find $\lambda = \frac{t}{v^2}$ reveals that

$$\mathbb{P}\{|g| > t\} \leq 2\exp\left(-\frac{t^2}{2v^2}\right). \quad (4.43)$$

Notice that $v$ is the Lipschitz constant of the function $F_v(g) := vg$ which maps a standard Gaussian random variable to a centred Gaussian random variable with variance $v^2$. The elementary concentration inequality (4.43) can therefore be written in terms of a standard Gaussian random variable $g$ as

$$\mathbb{P}\{|F_v(g) - \mathbb{E}F_v(g)| > t\} \leq 2\exp\left(-\frac{t^2}{2\|F_v\|_{\text{Lip}}^2}\right), \quad (4.44)$$

where the Lipschitz semi-norm $\|\cdot\|_{\text{Lip}}$ is defined in (3.38). It turns out that this formula generalizes to all Lipschitz continuous functions of a Gaussian random vector.

**Theorem 4.7** (Gaussian concentration inequality). *If $g$ is a standard Gaussian vector in $\mathbb{R}^d$ and $F : \mathbb{R}^d \to \mathbb{R}$ is a Lipschitz continuous function, then for any $t \geq 0$,*

$$\mathbb{P}\{|F(g) - \mathbb{E}F(g)| \geq t\} \leq 2\exp\left(-\frac{t^2}{2\|F\|_{\text{Lip}}^2}\right). \quad (4.45)$$

*Proof.* Replacing $F$ by $F - \mathbb{E}F(g)$ if necessary and using the union bound, it suffices to show that for every Lipschitz continuous function $F : \mathbb{R}^d \to \mathbb{R}$ with $\mathbb{E}F(g) = 0$ and every $t \geq 0$, we have

$$\mathbb{P}\{F(g) \geq t\} \leq \exp\left(-\frac{t^2}{2\|F\|_{\text{Lip}}^2}\right). \quad (4.46)$$

The proof proceeds in two steps. First we prove (4.46) when $F$ is continuously differentiable, and then we use an approximation argument to establish it in general.

*Step 1: F continuously differentiable.* By Chebyshev's inequality, we have for every $\lambda > 0$ that
$$\mathbb{P}\{F(g) \geq t\} \leq \exp(-\lambda t)\mathbb{E}\exp(\lambda F(g)). \tag{4.47}$$

To bound the function $h(\lambda) := \mathbb{E}\exp(\lambda F(g))$, we use a symmetrization trick. Let $\widetilde{g}$ be an independent copy of $g$, and for each $s \in [0,1]$, introduce the Gaussian random variable $Z_s := \sqrt{s}g + \sqrt{1-s}\widetilde{g}$. Define the function $G_\lambda(x) := \exp(\lambda F(x))$, and observe that by the fundamental theorem of calculus,
$$h'(\lambda) = \mathbb{E}F(g)G_\lambda(g) = \int_0^1 \frac{d}{ds}\mathbb{E}F(g)G_\lambda(Z_s)\,ds.$$

The independence of $g$ and $\widetilde{g}$ as well as the assumption $\mathbb{E}F(g) = 0$ have played their part. This implies that
$$h'(\lambda) = \int_0^1 \mathbb{E}F(g)\nabla G_\lambda(Z_s) \cdot \left(\frac{g}{2\sqrt{s}} - \frac{\widetilde{g}}{2\sqrt{1-s}}\right)ds.$$

For each $1 \leq i \leq d$, the Gaussian integration by parts formula reveals that
$$\mathbb{E}g_i F(g)\partial_{x_i}G(Z_s) = \mathbb{E}\partial_{x_i}\big(F(g)\partial_{x_i}G_\lambda(Z_s)\big)$$
$$= \lambda\mathbb{E}|\partial_{x_i}F(g)|^2 G_\lambda(Z_s) + \sqrt{s}\mathbb{E}F(g)\partial_{x_i}^2 G_\lambda(Z_s),$$
while
$$\mathbb{E}\widetilde{g}_i F(g)\partial_{x_i}G_\lambda(Z_s) = \sqrt{1-s}\mathbb{E}F(g)\partial_{x_i}^2 G_\lambda(Z_s).$$

It follows by the Cauchy-Schwarz inequality and the fact that $Z_s$ is equal in distribution to $g$ for every $s \in [0,1]$ that
$$h'(\lambda) \leq \frac{1}{2}\lambda\|F\|_{\text{Lip}}^2 \mathbb{E}G_\lambda(g)\int_0^1 \frac{1}{\sqrt{s}}\,ds \leq \lambda\|F\|_{\text{Lip}}^2 h(\lambda).$$

Dividing both sides of this equation by $h(\lambda)$ and integrating the resulting expression shows that
$$\log h(\lambda) \leq \log h(0) + \frac{\lambda^2\|F\|_{\text{Lip}}^2}{2} = \frac{\lambda^2\|F\|_{\text{Lip}}^2}{2}.$$

Substituting this into (4.47) and optimizing over $\lambda$ to find $\lambda = \frac{t}{\|F\|_{\text{Lip}}^2}$ establishes (4.46).

*Step 2: F general.* For each $\varepsilon > 0$, we introduce the function $F_\varepsilon(x) := \mathbb{E}F(x + \varepsilon\widetilde{g})$, where $\widetilde{g}$ is an independent copy of $g$. Since the random vector $x + \varepsilon\widetilde{g}$ is standard

Gaussian, the function $F_\varepsilon$ is smooth. Moreover, it is such that $\|F_\varepsilon\|_{\mathrm{Lip}} \leqslant \|F\|_{\mathrm{Lip}}$. It follows by (4.46) applied to $F_\varepsilon$ that for every integer $n \geqslant 1$,

$$\mathbb{P}\{F_\varepsilon(g) > t + n^{-1}\} \leqslant \exp\left(-\frac{(t+n^{-1})^2}{2\|F\|_{\mathrm{Lip}}^2}\right).$$

Since $(F_\varepsilon)_{\varepsilon>0}$ converges to $F$ uniformly, the sequence $(F_\varepsilon(g))_{\varepsilon>0}$ of random variables converges in law to $F(g)$. Invoking the Portmanteau theorem (Theorem A.17) reveals that

$$\mathbb{P}\{F(g) > t + n^{-1}\} \leqslant \liminf_{\varepsilon \to 0} \mathbb{P}\{F_\varepsilon(g) > t + n^{-1}\} \leqslant \exp\left(-\frac{(t+n^{-1})^2}{2\|F\|_{\mathrm{Lip}}^2}\right),$$

and using the continuity of measure to let $n$ tend to infinity completes the proof. ∎

The Gaussian concentration inequality has many important consequences in the context of statistical mechanics. For instance, as shown in Exercise 4.4, it gives the concentration of the maximum of the coordinates of a Gaussian vector around its average. As we will see in Exercise 6.2 of Chapter 6, this implies that the maximum of the Hamiltonian in the Sherrington-Kirkpatrick model concentrates around its average. In the language of physics, one could say that the ground-state energy is self-averaging. Another quantity that is generally self-averaging is the free energy. In the context of the symmetric rank-one matrix estimation problem, we can think of the free energy (4.8) as a function of the Gaussian noise $W$ and the signal $\bar{x}$. To make this dependence explicit, let us temporarily change notation and write the free energy as

$$F_N^\circ(t, W, \bar{x}) \coloneqq F_N^\circ(t). \tag{4.48}$$

Letting $a > 0$ be such that the support of the measure $P_1$ is contained in the interval $[-\sqrt{a}, \sqrt{a}]$, we see that the measure $P_N$ is supported in the ball $\overline{B}_{\sqrt{aN}}(0)$ of radius $\sqrt{aN}$ centred at the origin. We thus deduce that for $t \geqslant 0$, $\bar{x} \in \mathbb{R}^N$ and $W^1, W^2 \in \mathbb{R}^{N \times N}$,

$$F_N^\circ(t, W^1, \bar{x}) \leqslant \sqrt{\frac{2t}{N^3}} \sup_{x \in \overline{B}_{\sqrt{aN}}(0)} (x \cdot W^1 x - x \cdot W^2 x) + F_N^\circ(t, W^2, \bar{x}). \tag{4.49}$$

Using that $x \cdot (W^1 - W^2) x = (W^1 - W^2) \cdot (x^* x)$ and the Cauchy-Schwarz inequality, we obtain that

$$F_N^\circ(t, W^1, \bar{x}) - F_N^\circ(t, W^2, \bar{x}) \leqslant \sqrt{\frac{2t}{N^3}} \sup_{x \in \overline{B}_{\sqrt{aN}}(0)} |x|^2 |W^1 - W^2| \leqslant a\sqrt{\frac{2t}{N}} |W^1 - W^2|. \tag{4.50}$$

Since the right side of (4.50) is symmetric in the pair $(W^1, W^2)$, this estimate gives us an upper bound on the Lipschitz semi-norm of the map $W \mapsto F_N^\circ(t, W, \bar{x})$. It thus follows by the Gaussian concentration inequality applied conditionally on the randomness of the signal $\bar{x}$ that, for any $\lambda > 0$,

$$\mathbb{P}\{|F_N^\circ(t,W,\bar{x}) - \mathbb{E}_W F_N^\circ(t,W,\bar{x})| \geq \lambda\} \leq 2\exp\left(-\frac{N\lambda^2}{4a^2 t}\right), \quad (4.51)$$

where $\mathbb{E}_W$ denotes the average only with respect to the randomness of the Gaussian noise $W$. To obtain the full concentration of the free energy, it remains to establish the concentration of $\mathbb{E}_W F_N^\circ(t,W,\bar{x})$ about its average with respect to the signal $\bar{x}$. This can be done using the McDiarmid inequality.

**Theorem 4.8** (McDiarmid inequality). *Let $X = (X_1, \ldots, X_n)$ be a vector of independent random variables taking values in the measurable space $S_1 \times \cdots \times S_n$, let $c_1, \ldots, c_n \geq 0$, and let $F : S_1 \times \cdots \times S_n \to \mathbb{R}$ be a measurable function such that, for every $i \in \{1, \ldots, n\}$, $x_1 \in S_1, \ldots, x_n \in S_n$, and $x_i' \in S_i$,*

$$\left|F(x_1, \ldots, x_{i-1}, x_i', x_{i+1}, \ldots, x_n) - F(x_1, \ldots, x_{i-1}, x_i, x_{i+1}, \ldots, x_n)\right| \leq c_i. \quad (4.52)$$

*For every $t \geq 0$, we have*

$$\mathbb{P}\{|F(X) - \mathbb{E} F(X)| \geq t\} \leq 2\exp\left(-\frac{2t^2}{\sum_{i=1}^n c_i^2}\right). \quad (4.53)$$

*Proof.* We decompose the proof into two steps.

*Step 1: Laplace transform estimate.* Let $a < b \in \mathbb{R}$. In this step, we show that if $X$ is a random variable taking values in the interval $[a,b]$, then

$$\log \mathbb{E} \exp \lambda (X - \mathbb{E} X) \leq \frac{\lambda^2 (b-a)^2}{8}. \quad (4.54)$$

As a preliminary observation, any random variable $X$ taking values in the interval $[a,b]$ must be such that

$$\mathbb{E}(X - \mathbb{E} X)^2 = \inf_{c \in \mathbb{R}} \mathbb{E}(X - c)^2 \leq \mathbb{E}\left(X - \frac{a+b}{2}\right)^2 \leq \left(\frac{b-a}{2}\right)^2.$$

We denote by $\psi(\lambda)$ the left side of (4.54). Since $X$ is bounded, the function $\psi$ is infinitely differentiable, and a direct calculation gives that

$$\psi'(\lambda) = \langle X - \mathbb{E} X \rangle \quad \text{and} \quad \psi''(\lambda) = \langle (X - \mathbb{E} X)^2 \rangle - \left(\langle X - \mathbb{E} X \rangle\right)^2,$$

where we introduced the notation

$$\langle Z \rangle := \frac{\mathbb{E}(Z \exp \lambda (X - \mathbb{E} X))}{\mathbb{E} \exp \lambda (X - \mathbb{E} X)}$$

to denote the Gibbs measure associated with the Hamiltonian $\lambda(X - \mathbb{E}X)$. Although this is kept implicit in the notation, we point out that this Gibbs measure depends on $\lambda$; in particular, the quantities $\langle X \rangle$ and $\mathbb{E}X$ are different in general. We see that $\psi''$ is expressed as the variance of a random variable that takes values in an interval of length $b - a$. By the preliminary observation, we deduce that $\psi''$ is uniformly bounded by $(b-a)^2/4$. Since $\psi(0) = \psi'(0) = 0$, integrating this bound yields (4.54).

*Step 2: Chebyshev inequality.* We now take up the notation and assumptions of the theorem, and introduce, for every $i \in \{1, \ldots, n\}$, the difference

$$\Delta_i := \mathbb{E}\big[f(X) \mid X_1, \ldots, X_i\big] - \mathbb{E}\big[f(X) \mid X_1, \ldots, X_{i-1}\big],$$

with the understanding that $\Delta_1 = \mathbb{E}[f(X) \mid X_1] - \mathbb{E}f(X)$. We clearly have

$$f(X) - \mathbb{E}f(X) = \sum_{i=1}^{n} \Delta_i.$$

For every $\lambda \geq 0$, we now look for an upper bound on

$$\mathbb{E}\exp\Big(\lambda \sum_{i=1}^{n} \Delta_i\Big) = \mathbb{E}\prod_{i=1}^{n} \exp \lambda \Delta_i.$$

Conditionally on $X_1, \ldots, X_{n-1}$, the random variable $\Delta_n$ is centred, and by (4.52), it takes values in an interval of length $c_n$. Since the random variables $\Delta_1, \ldots, \Delta_{n-1}$ are measurable with respect to $X_1, \ldots, X_{n-1}$, we can apply the result of the previous step conditionally on these variables, so that

$$\mathbb{E}\exp\Big(\lambda \sum_{i=1}^{n} \Delta_i\Big) \leq \mathbb{E}\exp\Big(\lambda \sum_{i=1}^{n-1} \Delta_i\Big) \exp\Big(\frac{\lambda^2 c_n^2}{8}\Big).$$

Continuing inductively from $i = n - 1$ down to $i = 1$, we conclude that

$$\mathbb{E}\exp\Big(\lambda \sum_{i=1}^{n} \Delta_i\Big) \leq \exp\Big(\frac{\lambda^2}{8} \sum_{i=1}^{n} c_i^2\Big).$$

By Chebyshev's inequality, for every $t \geq 0$ and $\lambda \geq 0$, we have

$$\mathbb{P}\{F(X) - \mathbb{E}F(X) \geq t\} \leq \exp\Big(-\lambda t + \frac{\lambda^2}{8} \sum_{i=1}^{n} c_i^2\Big).$$

Optimizing over $\lambda$ yields that

$$\mathbb{P}\{F(X) - \mathbb{E}F(X) \geq t\} \leq \exp\Big(-\frac{2t^2}{\sum_{i=1}^{n} c_i^2}\Big).$$

By symmetry and the union bound, we thus obtain (4.53). This completes the proof. ∎

To control the deviations of $\mathbb{E}_W F_N^\circ(t,W,\bar{x})$ about its average using the McDiarmid inequality, fix $t \geq 0, W \in \mathbb{R}^{N \times N}$ and $\bar{x} \in \mathbb{R}^N$. Recalling (4.6), a direct computation reveals that for any $1 \leq i \leq N$,

$$\partial_{\bar{x}_i} F_N^\circ(t,W,\bar{x}) = \frac{4t}{N^2} \langle x_i (x \cdot \bar{x}) \rangle. \tag{4.55}$$

Letting $a > 0$ be such that the support of the measure $P_1$ is contained in the interval $[-\sqrt{a}, \sqrt{a}]$, this implies that

$$\left| \partial_{\bar{x}_i} F_N^\circ(t,W,\bar{x}) \right| \leq \frac{4t a^{3/2}}{N}. \tag{4.56}$$

Averaging over the randomness of the Gaussian noise $W$ and using Jensen's inequality shows that this upper bound also holds for the averaged free energy $\mathbb{E}_W F_N^\circ(t,W,\bar{x})$. Combining this with the mean value theorem, we find that for any $i \in \{1,\ldots,N\}$ and $\bar{x}_1,\ldots,\bar{x}_N,\bar{x}_i' \in [-\sqrt{a}, \sqrt{a}]$,

$$\left| \mathbb{E}_W F_N^\circ(t,W,\bar{x}_1,\ldots,\bar{x}_{i-1},\bar{x}_i',\bar{x}_{i+1},\ldots,\bar{x}_n) \right.$$
$$\left. - \mathbb{E}_W F_N^\circ(\bar{x}_1,\ldots,\bar{x}_{i-1},\bar{x}_i,\bar{x}_{i+1},\bar{x}_n) \right| \leq \frac{4t a^{3/2}}{N} |\bar{x}_i - \bar{x}_i'| = \frac{8t a^2}{N}. \tag{4.57}$$

It follows by the McDiarmid inequality that, for every $\lambda \geq 0$,

$$\mathbb{P}\left\{ |\mathbb{E}_W F_N^\circ(t,W,\bar{x}) - \mathbb{E} F_N^\circ(t,W,\bar{x})| \geq \lambda \right\} \leq 2 \exp\left( -\frac{N\lambda^2}{32 t^2 a^4} \right). \tag{4.58}$$

Together with the triangle inequality and the upper bound (4.51), we find, for every $T < +\infty$, a constant $C < +\infty$ such that for every $t \in [0,T]$ and $\lambda \geq 0$,

$$\mathbb{P}\left\{ |F_N^\circ(t) - \overline{F}_N^\circ(t)| \geq \lambda \right\} \leq 2 \exp\left( -\frac{N\lambda^2}{C} \right). \tag{4.59}$$

This establishes the exponential concentration of the free energy about its average.

There are many other concentration inequalities for Gaussian and non-Gaussian random variables which are often used in statistical mechanics. For instance, alternatively to the McDiarmid inequality, the Efron-Stein inequality discussed in Exercise 4.7 can be used to establish the concentration of the free energy in non-Gaussian contexts, such as in Appendix B of [105] for community detection on sparse graphs. Exercise 4.8 presents the Gaussian Poincaré inequality, which can be used to obtain concentration results in $L^2$ under a weaker assumption than Lipschitz continuity of the functional. For much more on concentration inequalities we refer the interested reader to [59].

**Exercise 4.4.** Let $a \in \mathbb{R}_{\geqslant 0}$ and let $(g_1, \ldots, g_d)$ be a $d$-dimensional centred Gaussian vector with $\mathbb{E} g_i^2 \leqslant a$ for every $i \in \{1, \ldots, d\}$. Show that for every $t \geqslant 0$,

$$\mathbb{P}\{|\max_{1 \leqslant i \leqslant d} g_i - \mathbb{E} \max_{1 \leqslant i \leqslant d} g_i| \geqslant t\} \leqslant 2\exp\left(-\frac{t^2}{2a}\right). \tag{4.60}$$

**Exercise 4.5.** For $W = (W_{ij})_{1 \leqslant i,j \leqslant N} \in \mathbb{R}^{N \times N}$ a matrix of independent standard Gaussian random variables, we introduce the norm $\|W\|_* := \sup_{|x| \leqslant 1} |Wx|$. The purpose of this exercise is to show that there exists a constant $C < +\infty$ such that for every $a \geqslant C$,

$$\mathbb{P}\{\|W\|_*^2 \geqslant aN\} \leqslant \exp\left(-\frac{aN}{C}\right). \tag{4.61}$$

(i) Show that there exists a constant $C < +\infty$ such that for every $x \in \mathbb{R}^N$ of unit norm and all $a \geqslant 0$,

$$\mathbb{P}\{|Wx|^2 \geqslant aN\} \leqslant \exp\left(\left(C - \frac{a}{C}\right)N\right). \tag{4.62}$$

(ii) Consider a set $A \subseteq \mathbb{R}^N$ with the property that for every $x \in \mathbb{R}^N$ with $|x| \leqslant 1$, there exists $y \in A$ such that $|x - y| \leqslant \frac{1}{2}$. Show that

$$\|W\|_* \leqslant 2 \sup_{y \in A} |Wy|. \tag{4.63}$$

(iii) Show that one can find such a set $A$ of size exponential in $N$, and conclude that (4.61) holds.

(iv) Deduce that for every $q \geqslant 1$, there exists a constant $C < +\infty$ with

$$\mathbb{E}\|W\|_*^q \leqslant CN^{q/2}. \tag{4.64}$$

**Exercise 4.6** (Approximate Gaussian integration by parts). Let $X$ be a centred random variable with finite third moment $\mathbb{E}|X|^3 < +\infty$. Show that for every $F \in C^2(\mathbb{R};\mathbb{R})$ with $\|F''\|_\infty < +\infty$,

$$\left|\mathbb{E} X F(X) - \mathbb{E} X^2 \mathbb{E} F'(X)\right| \leqslant \frac{3}{2}\|F''\|_\infty \mathbb{E}|X|^3. \tag{4.65}$$

**Exercise 4.7** (Efron-Stein). Let $X \in \mathbb{R}^d$ be a vector of independent random variables, and let $X' \in \mathbb{R}^d$ be an independent copy of $X$. Show that for every bounded and measurable $f : \mathbb{R}^d \to \mathbb{R}$,

$$\operatorname{Var} f(X) \leqslant \frac{1}{2} \sum_{i=1}^d \mathbb{E}\big(f(X) - f(X_1, \ldots, X_{i-1}, X_i', X_{i+1}, \ldots, X_N)\big)^2. \tag{4.66}$$

**Exercise 4.8** (Gaussian Poincaré inequality). Let $Z$ be a standard Gaussian random variable, and let $f \in C_c^\infty(\mathbb{R};\mathbb{R})$ be a smooth and compactly supported function. Show that $\operatorname{Var} f(Z) \leqslant \mathbb{E}|f'(Z)|^2$.

## 4.3   A Hamilton-Jacobi approach to rank-one matrix estimation

To determine the limit free energy in the symmetric rank-one matrix estimation problem using the Hamilton-Jacobi approach, we need to find a partial differential equation satisfied by the free energy (4.9) up to an error that vanishes with $N$. Notice that the limits of the random free energy (4.8) and its average (4.9), provided that they exist, are the same. Indeed, by the free energy concentration inequality (4.59), there exists a constant $C < +\infty$ that depends only on an upper bound on $t$ and on the support of $P_1$ such that for every $\lambda > 0$,

$$\mathbb{P}\{|F_N^\circ(t) - \overline{F}_N^\circ(t)| \geq \lambda\} \leq 2\exp\left(-\frac{N\lambda^2}{C}\right). \tag{4.67}$$

Since the right side of this expression is summable, the Borel-Cantelli lemma implies that almost surely,

$$\limsup_{N \to +\infty} |F_N^\circ(t) - \overline{F}_N^\circ(t)| = 0. \tag{4.68}$$

Using the Gaussian integration by parts formula, we saw in (4.33) that

$$\mathbb{E}\langle x \cdot Wx \rangle = \sqrt{\frac{2t}{N}}\left(\mathbb{E}\langle |x|^4 \rangle - \mathbb{E}\langle (x \cdot x')^2 \rangle\right). \tag{4.69}$$

Together with the derivative computation in (4.20) this implies that

$$\partial_t \overline{F}_N^\circ(t) = \frac{1}{N^2}\mathbb{E}\langle (x \cdot \bar{x})^2 \rangle. \tag{4.70}$$

Unfortunately, there is no way of closing this equation if we can only compute derivatives in $t$ of the free energy. Indeed, the situation is analogous to that encountered if we were studying a Curie-Weiss model without any magnetization part, where there would be no parameter $h$ with respect to which to differentiate and close the equation. To overcome this issue, we will define an "enriched" free energy that also depends on an additional parameter $h$. The enriched free energy $\overline{F}_N(t, h)$ should extend the free energy $\overline{F}_N^\circ(t)$, in the sense that there is some $h_0$ with $\overline{F}_N(\cdot, h_0) = \overline{F}_N^\circ(\cdot)$; it should be simple enough that we can explicitly compute the limit of its initial condition $\overline{F}_N(0, h)$; and it should be rich enough that it allows us to close the equation, up to a small error term. In the context of statistical inference models, we also want to ensure that the enrichment does not destroy the fact that the Gibbs measure is a conditional expectation. Indeed, this property gives us access to the Nishimori identity, which plays a fundamental role in simplifying statistical inference models and distinguishing them from the more complicated spin glass models we will discuss in Chapter 6.

In the context of the symmetric rank-one matrix estimation problem, the appropriate enrichment of the free energy is obtained by assuming that, in addition to observing the noisy rank-one matrix $Y$ in (4.1), we also observe a noisy version $\widetilde{Y}$ of the signal vector $\bar{x}$,

$$\widetilde{Y} := \sqrt{2h}\,\bar{x} + z. \tag{4.71}$$

The noise vector $z = (z_i)_{i \leqslant N} \in \mathbb{R}^N$ is made of independent standard Gaussian random variables independent of the vector $\bar{x}$ and the noise matrix $W$, and the parameter $h \geqslant 0$ is a signal-to-noise ratio. The enriched symmetric rank-one matrix estimation problem is to infer the rank-one matrix $\bar{x}\bar{x}^*$ from the observation of $\mathcal{Y} := (Y, \widetilde{Y})$. Applying Bayes' formula as in Exercise 4.2 shows that the law of the signal $\bar{x}$ given the observation of $\mathcal{Y}$ is the Gibbs measure whose Hamiltonian on $\mathbb{R}^N$ is

$$H_N(t,h,x) := H_N^\circ(t,x) + \sqrt{2h}\,\widetilde{Y} \cdot x - h|x|^2. \tag{4.72}$$

In other words, for any bounded and measurable function $f : \mathbb{R}^N \to \mathbb{R}$, we have

$$\mathbb{E}[f(\bar{x}) \mid \mathcal{Y}] = \frac{\int_{\mathbb{R}^N} f(x) \exp H_N(t,h,x)\,dP_N(x)}{\int_{\mathbb{R}^N} \exp H_N(t,h,x)\,dP_N(x)}. \tag{4.73}$$

The free energy

$$F_N(t,h) := \frac{1}{N} \log \int_{\mathbb{R}^N} \exp H_N(t,h,x)\,dP_N(x) \tag{4.74}$$

of this model is again random, as it depends on $\bar{x}$, $W$ and $z$, and we denote its average by

$$\overline{F}_N(t,h) := \mathbb{E} F_N(t,h) = \frac{1}{N} \mathbb{E} \log \int_{\mathbb{R}^N} \exp H_N(t,h,x)\,dP_N(x). \tag{4.75}$$

We will again often refer to the average free energy (4.75) as simply the free energy. Through a slight abuse of notation, we will, as before, write $\langle \cdot \rangle$ for the average with respect to the Gibbs measure (4.73), and write $x, x', x''$, and so on, for independent random variables sampled according to this probability measure. That is, for every bounded and measurable function $f : \mathbb{R}^N \to \mathbb{R}$, we write

$$\langle f(x) \rangle := \frac{\int_{\mathbb{R}^N} f(x) \exp H_N(t,h,x)\,dP_N(x)}{\int_{\mathbb{R}^N} \exp H_N(t,h,x)\,dP_N(x)}, \tag{4.76}$$

and so on as in (4.12) with more replicas. The derivation of the identity (4.70) is unchanged for this more general Gibbs measure: we have for all $t, h \geqslant 0$ that

$$\partial_t \overline{F}_N(t,h) = \frac{1}{N^2} \mathbb{E}\langle (x \cdot \bar{x})^2 \rangle. \tag{4.77}$$

Recalling from (4.71) and (4.72) that

$$H_N(t,h,x) = H_N^\circ(t,x) + 2hx\cdot\bar{x} + \sqrt{2h}z\cdot x - h|x|^2, \qquad (4.78)$$

we next observe that for all $t, h \geq 0$,

$$\partial_h \overline{F}_N(t,h) = \frac{1}{N\sqrt{2h}} \mathbb{E}\langle z\cdot x\rangle + \frac{2}{N}\mathbb{E}\langle x\cdot\bar{x}\rangle - \frac{1}{N}\mathbb{E}\langle |x|^2\rangle. \qquad (4.79)$$

To integrate out the Gaussian noise $z$, we fix $1 \leq i \leq N$, we write $F(z, W, \bar{x}) := \langle x_i\rangle$, and we observe that by the Gaussian integration by parts formula in (4.29),

$$\mathbb{E}\langle z_i x_i\rangle = \mathbb{E}z_i F(z, W, \bar{x}) = \mathbb{E}\partial_{z_i} F(z, W, \bar{x}) = \sqrt{2h}\big(\mathbb{E}\langle x_i^2\rangle - \mathbb{E}\langle x_i\rangle^2\big). \qquad (4.80)$$

It follows that

$$\mathbb{E}\langle z\cdot x\rangle = \sqrt{2h}\big(\mathbb{E}\langle |x|^2\rangle - \mathbb{E}\langle x\cdot x'\rangle\big), \qquad (4.81)$$

which together with the Nishimori identity implies that

$$\partial_h \overline{F}_N(t,h) = \frac{1}{N}\mathbb{E}\langle x\cdot\bar{x}\rangle. \qquad (4.82)$$

It is reasonable to expect the variance of the inner product, or overlap, $N^{-1}x\cdot\bar{x}$ between a sample $x$ from the Gibbs measure (4.73) and the ground-truth signal $\bar{x}$ to be small, simply because it is the average of a large number of variables. If this is so, then the difference between the time derivative (4.77) and the square of the spatial derivative (4.82) would also be small, since

$$\partial_t \overline{F}_N(t,h) - \big(\partial_h \overline{F}_N(t,h)\big)^2 = \mathrm{Var}\left(\frac{x\cdot\bar{x}}{N}\right). \qquad (4.83)$$

This suggests that the enriched free energy (4.75) might converge to the function $f$ solving the Hamilton-Jacobi equation

$$\partial_t f - (\partial_h f)^2 = 0 \quad \text{on} \quad \mathbb{R}_{>0}\times\mathbb{R}_{>0} \qquad (4.84)$$

with the initial condition

$$\psi(h) := \lim_{N\to+\infty} \overline{F}_N(0,h) = \overline{F}_1(0,h). \qquad (4.85)$$

We have used the assumption that $P_N$ is a product measure to assert that the initial condition is independent of $N$. Notice that the Hamilton-Jacobi equation (4.84) is the same as the Hamilton-Jacobi equation (3.9) that appeared in the context of the Curie-Weiss model, with the sole difference that it is posed on $\mathbb{R}_{>0}\times\mathbb{R}_{>0}$ as opposed to $\mathbb{R}_{>0}\times\mathbb{R}$.

Because of this difference, the Hamilton-Jacobi equation (4.84) that we would like to use to describe the limit free energy in the symmetric rank-one matrix estimation problem does not exactly fall within the scope of the Hamilton-Jacobi equations studied in Chapter 3. Indeed, now the spatial variable is restricted to the space $\mathbb{R}_{\geq 0}$, which unlike $\mathbb{R}$ is a domain with a boundary. The right way forward would be to extend the analysis in Chapter 3 to equations posed on more general domains. How to proceed along these lines and conclude for the identification of the limit free energy in the symmetric rank-one matrix estimation problem is explained in details in [103], in a setting that closely matches our current one. This approach is also taken up in [74, 75] in great generality in finite dimension, in [73, 195] for the infinite-dimensional equations arising in Chapter 6 in relation with spin glasses, and in [104] for community detection on sparse graphs. The upshot of this analysis is that, for all equations of relevance to us, we do not have to prescribe any boundary condition for the solution, and in effect, we can proceed by simply ignoring the boundary. Intuitively, this is possible because the characteristic lines discussed below (3.90) always go towards the boundary rather than away from it.

Instead of explaining this, we prefer to use a simpler but less general workaround here, by taking advantage of the simple geometry of the domain $\mathbb{R}_{\geq 0}$ through a symmetrization trick. Although this will be convenient in our context, the reader should not surmise that this is a standard way to handle equations posed on domains with boundaries, and there are indeed more complicated models of statistical inference or spin glasses, where $\mathbb{R}_{\geq 0}$ is replaced by the space of positive semi-definite matrices, for which this trick does noes seem to be applicable.

So we will content ourselves with saying that a function $f : \mathbb{R}_{\geq 0} \times \mathbb{R}_{\geq 0} \to \mathbb{R}$ is a viscosity solution to the Hamilton-Jacobi equation (4.84) with initial condition $\psi$ if its symmetrization $\widetilde{f} : \mathbb{R}_{\geq 0} \times \mathbb{R} \to \mathbb{R}$ defined by

$$\widetilde{f}(t,h) := f(t,|h|) \tag{4.86}$$

is a viscosity solution the Hamilton-Jacobi equation

$$\partial_t \widetilde{f} - \left(\partial_h \widetilde{f}\right)^2 = 0 \quad \text{on} \quad \mathbb{R}_{>0} \times \mathbb{R} \tag{4.87}$$

subject to the symmetrized initial condition $\widetilde{\psi}(h) := \psi(|h|)$. As will be seen, the main reason why this approach makes sense is that the functions we care about are non-decreasing in the variable $h$, when this variable varies in $\mathbb{R}_{\geq 0}$. If the function $\psi$ was decreasing, then the symmetrization would instead create an unphysical singularity at the origin $h = 0$ — by this we mean that we would be able to contradict the definition of viscosity solution for (4.87) by using any point in $\mathbb{R}_{\geq 0} \times \{0\}$ as a contact point.

Our main result will be that the limit free energy (4.75) is indeed a solution to the Hamilton-Jacobi equation (4.84). Together with the Hopf-Lax formula in Theorem 3.8 and the observation that the free energy (4.9) can be recovered from the enriched free energy (4.75),

$$\overline{F}_N^\circ(t) = \overline{F}_N(t,0), \qquad (4.88)$$

this will give us the following result for the limit free energy in the symmetric rank-one matrix estimation problem.

**Theorem 4.9** (Identification of limit free energy). *For every $N \geq 1$, we denote by $\overline{F}_N : \mathbb{R}_{\geq 0} \times \mathbb{R}_{\geq 0} \to \mathbb{R}$ the enriched free energy (4.75) in the symmetric rank-one matrix estimation problem. For every $t, h \geq 0$, we have that $\overline{F}_N(t,h)$ converges to $f(t,h)$ as $N$ tends to infinity, where $f : \mathbb{R}_{\geq 0} \times \mathbb{R}_{\geq 0} \to \mathbb{R}$ is the unique viscosity solution to the Hamilton-Jacobi equation (4.84) subject to the initial condition*

$$\psi(h) := \overline{F}_1(0,h) = \mathbb{E} \log \int_\mathbb{R} \exp\left(\sqrt{2h} x z_1 + 2h x \overline{x}_1 - h x^2\right) dP_1(x). \qquad (4.89)$$

*The limit free energy $f$ admits the Hopf-Lax representation, for every $t, h \geq 0$,*

$$f(t,h) = \sup_{h' \geq 0} \left(\psi(h+h') - \frac{(h')^2}{4t}\right). \qquad (4.90)$$

*In particular, the limit of the free energy $\overline{F}_N^\circ$ defined in (4.9) is given by*

$$\lim_{N \to +\infty} \overline{F}_N^\circ(t) = \sup_{h \geq 0} \left(\psi(h) - \frac{h^2}{4t}\right). \qquad (4.91)$$

This result will be proved by applying the convex selection principle to the symmetrization $\widetilde{f}$ of any subsequential limit $f$ of the sequence $(\overline{F}_N)_{N \geq 1}$ of enriched free energies. To begin with, we will verify that the enriched free energy $\overline{F}_N$ is jointly convex, Lipschitz continuous and non-decreasing, and that the first two of these three properties are inherited by $\widetilde{f}$. We say that a function $g : \mathbb{R}_{\geq 0} \times \mathbb{R}_{\geq 0} \to \mathbb{R}$ is *non-decreasing* if, for all $t, h, h' \geq 0$, we have

$$h \leq h' \implies g(t,h) \leq g(t,h'). \qquad (4.92)$$

We will then show that each free energy $\overline{F}_N$ satisfies the Hamilton-Jacobi equation (4.84) up to some explicit error term, and use this to deduce that any smooth function that touches $f$ from above must satisfy the Hamilton-Jacobi equation (4.84) at the contact point. This will then allow us to conclude that $\widetilde{f}$ is the unique viscosity solution to the Hamilton-Jacobi equation (4.87) using the convex selection principle in Lemma 3.23. We will rely on this more general version of the convex selection principle, as opposed to the cleaner form in Theorem 3.21, because the symmetrized initial condition $\widetilde{\psi}$ is not necessarily continuously differentiable at the origin. As can easily be seen, e.g. from (S.27), the symmetrized initial condition $\widetilde{\psi}$ does belong to $C^1(\mathbb{R};\mathbb{R})$ when $\mathbb{E}\bar{x}_1 = 0$, and the reader may consider working out a simpler proof based on Theorem 3.21 under this additional assumption.

We start as announced by showing that the enriched free energy $\overline{F}_N$ is jointly convex, non-decreasing in the sense of (4.92), and Lipschitz continuous uniformly over $N$.

**Proposition 4.10.** *The enriched free energy* (4.75) *in the symmetric rank-one matrix estimation problem is jointly convex, non-decreasing, and Lipschitz continuous uniformly over N.*

The proof that the free energy $\overline{F}_N$ is convex involves a somewhat lengthy calculation. This may come as a surprise, since for the simpler models studied in previous chapters, the convexity of the free energy was a consequence of the very general observation presented in Exercise 2.6 that the log-Laplace transform of a random variable is convex. This argument cannot be applied here however. One may wonder whether this comes from the fact that we placed some square roots acting on the parameters $t$ and $h$ in the definition of the problem; and in the context of spin glasses discussed in Chapter 6, it will be immediate to see that the function $\widetilde{F}_N : (t,h) \mapsto \overline{F}_N(t^2,h^2)$ is indeed convex. So are we making the question uselessly complicated here? The main reason we want to insist on showing the convexity property is that this is a requirement for the validity of the convex selection principle. The proof of this result relies on the fact that the underlying Hamilton-Jacobi equation does not explicitly depend on the parameters $t$ and $h$. If we were to write down a Hamilton-Jacobi equation for $\widetilde{F}_N$ instead, then the parameters $t$ and $h$ would have to be explicitly present in the equation; and it turns out that the convex selection principle is actually false in this more general setting.

In the context of statistical inference, there is however a fundamental information-theoretic reason to expect that the free energy $\overline{F}_N$ as defined in (4.75) is indeed convex. For the interested reader, we briefly sketch this argument and why it does not actually allow us to show the joint convexity of $\overline{F}_N$. Recall from (4.24) that we denote by $\mathsf{I}_N(t)$ the mutual information between the signal $\bar{x}$ and the observation $Y$. By Proposition 4.4, we see that the convexity of $\overline{F}_N(\cdot,0)$ is equivalent to the concavity of the mutual information $\mathsf{I}_N$. Moreover, one can show that when we observe $Y$, we gain exactly as much information about the signal as if we were observing, for two independent copies $W^1$ and $W^2$ of $W$, the quantities

$$\sqrt{\frac{t}{N}}\bar{x}\bar{x}^* + W^1 \quad \text{and} \quad \sqrt{\frac{t}{N}}\bar{x}\bar{x}^* + W^2. \tag{4.93}$$

Notice that compared with the definition of $Y$ in (4.1), the variable $t$ has been replaced by $t/2$ in each of the two quantities above. Finally, once we have observed the first of the two quantities in (4.93), one can verify that we are going to gain at most as much information when we subsequently observe the second quantity in (4.93). In other words, the mutual information $I_N$ satisfies a sort of subadditivity property, and elementary properties of the mutual information allow us to upgrade this to the fact that $I_N$ is a concave function. The argument we just sketched therefore leads to the conclusion that the function $\overline{F}_N(\cdot,0)$ is convex. Minor variants of this argument yield that the function $\overline{F}_N$ is convex in each of the variables separately.

It may seem plausible that this argument can be generalized and lead to a conceptual information-theoretic proof that the function $\overline{F}_N$ is convex *jointly* in $(t,h)$. However, pushing this argument to a setting with multiple variables yields instead that every entry of the Hessian of $\overline{F}_N$ is non-negative. To our surprise, for the problem of community detection on sparse graphs, which unlike the situation considered in this chapter is not reducible to a problem with Gaussian noise, one can indeed show that the relevant free energy in this context is *not* convex in general [153].

To sum up this informal discussion, the proof that $\overline{F}_N$ is jointly convex does have to use some aspects of the particular structure of the class of problems we consider, and thus at least some calculations do need to be made. We also refer to Proposition 3.1 of [153] for a somewhat more general view on such calculations.

*Proof of Proposition 4.10.* The derivative computations (4.77) and (4.82) as well as the boundedness of the support of $P_1$ imply that the first order derivatives of the enriched free energy $\overline{F}_N$ are uniformly bounded. This establishes the uniform Lipschitz continuity of $\overline{F}_N$. That $\overline{F}_N$ is non-decreasing follows from the derivative computation (4.82) and the Nishimori identity, $\partial_h \overline{F}_N(t,h) = N^{-1}\mathbb{E}\langle x \rangle^2 \geqslant 0$. To prove the convexity of $\overline{F}_N$, we show that its Hessian is non-negative definite. Differentiating the expression (4.82) in $h$ and recalling (4.78), we can write

$$N\partial_h^2 \overline{F}_N(t,h)$$
$$= \mathbb{E}\langle (x\cdot\bar{x})\partial_h H_N(t,h,x)\rangle - \mathbb{E}\langle (x\cdot\bar{x})\partial_h H_N(t,h,x')\rangle$$
$$= \frac{1}{\sqrt{2h}}\mathbb{E}\langle (x\cdot\bar{x})(z\cdot x - z\cdot x')\rangle + 2\mathbb{E}\langle (x\cdot\bar{x})(x\cdot\bar{x} - x'\cdot\bar{x})\rangle - \mathbb{E}\langle (x\cdot\bar{x})(|x|^2 - |x'|^2)\rangle.$$

We now use the Gaussian integration by parts formula in (4.29) to integrate out the noise $z$, as in the calculations leading to (4.33) or to (4.81). Recalling also that when the Gibbs measure is an average over the two variables $x$ and $x'$, the underlying Hamiltonian is $H_N(t,h,x) + H_N(t,h,x')$, we obtain

$$\mathbb{E}\langle (x\cdot\bar{x})(z\cdot x)\rangle = \sqrt{2h}\mathbb{E}\langle (x\cdot\bar{x})(|x|^2 - x\cdot x')\rangle,$$
$$\mathbb{E}\langle (x\cdot\bar{x})(z\cdot x')\rangle = \sqrt{2h}\mathbb{E}\langle (x\cdot\bar{x})(|x'|^2 + x\cdot x' - 2x'\cdot x'')\rangle.$$

It follows by the Nishimori identity that

$$N\partial_h^2 \overline{F}_N(t,h) = 2\mathbb{E}\langle (x\cdot x')^2\rangle - 4\mathbb{E}\langle (x\cdot x')(x\cdot x'')\rangle + 2\mathbb{E}\langle (x\cdot x')(x''\cdot x''')\rangle. \quad (4.94)$$

We introduce the re-scaled and centred variable

$$y := \frac{1}{\sqrt{N}}(x - \langle x \rangle),$$

with $y', y'', y'''$ denoting independent copies of $y$ under the measure $\langle \cdot \rangle$, so that

$$\begin{aligned} N^2 \mathbb{E}\langle (y \cdot y')^2 \rangle &= \mathbb{E}\langle ((x - \langle x \rangle) \cdot (x' - \langle x \rangle))^2 \rangle \\ &= \mathbb{E}\langle (x \cdot x' - x \cdot \langle x \rangle - x' \cdot \langle x \rangle + \langle x \rangle \cdot \langle x \rangle)^2 \rangle \\ &= \mathbb{E}\langle (x \cdot x')^2 \rangle - 2\mathbb{E}\langle (x \cdot x')(x \cdot x'') \rangle + \mathbb{E}\langle (x \cdot x')(x'' \cdot x''') \rangle \\ &= \frac{N}{2} \partial_h^2 \overline{F}_N(t,h). \end{aligned} \quad (4.95)$$

This already shows that $\overline{F}_N$ is convex in the $h$ variable. To compute the second derivative in $t$, we can rewrite (4.77) in the form of

$$N^2 \partial_t \overline{F}_N(t,h) = \mathbb{E}\langle xx^* \cdot \bar{x}\bar{x}^* \rangle, \quad (4.96)$$

and follow through the same calculation as for the second derivative in $h$; the only difference is that each occurrence of $x$ is replaced by $xx^*$, each occurrence of $\bar{x}$ is replaced by $\bar{x}\bar{x}^*$, and so on, with $z$ being replaced by $W$. This leads to

$$N^3 \partial_t^2 \overline{F}_N(t,h) = 2\mathbb{E}\langle (xx^* \cdot x'x'^*)^2 \rangle - 4\mathbb{E}\langle (xx^* \cdot x'x'^*)(xx^* \cdot x''x''^*) \rangle \\ + 2\mathbb{E}\langle (xx^* \cdot x'x'^*)(x''x''^* \cdot x'''x'''^*) \rangle,$$

so in terms of the re-scaled and centred variable

$$\xi := \frac{1}{N}(xx^* - \langle xx^* \rangle),$$

we get

$$\frac{1}{2N} \partial_t^2 \overline{F}_N(t,h) = \mathbb{E}\langle (\xi \cdot \xi')^2 \rangle.$$

For the cross-derivative, we start from (4.96) and differentiate in $h$ to obtain

$$N^2 \partial_h \partial_t \overline{F}_N(t,h) = \frac{1}{\sqrt{2h}} \mathbb{E}\langle (xx^* \cdot \bar{x}\bar{x}^*)(z \cdot x - z \cdot x') \rangle + 2\mathbb{E}\langle (xx^* \cdot \bar{x}\bar{x}^*)(x \cdot \bar{x} - x' \cdot \bar{x}) \rangle \\ - \mathbb{E}\langle (xx^* \cdot \bar{x}\bar{x}^*)(|x|^2 - |x'|^2) \rangle.$$

A Gaussian integration by parts allows us to rewrite the first term on the right side above as

$$\mathbb{E}\langle (xx^* \cdot \bar{x}\bar{x}^*)(|x|^2 - x \cdot x' - (|x'|^2 + x \cdot x' - 2x' \cdot x'')) \rangle.$$

An application of the Nishimori identity therefore yields that

$$N^2 \partial_h \partial_t \overline{F}_N(t,h) = 2\mathbb{E}\langle (xx^* \cdot x'x'^*)(x \cdot x') \rangle - 4\mathbb{E}\langle (xx^* \cdot x'x'^*)(x \cdot x'') \rangle \\ - 2\mathbb{E}\langle (xx^* \cdot x'x'^*)(x'' \cdot x''') \rangle,$$

which can be rewritten in terms of the re-scaled and centred variables as

$$\frac{1}{2N} \partial_h \partial_t \overline{F}_N(t,h) = \mathbb{E}\langle (\xi \cdot \xi')(y \cdot y') \rangle.$$

To see that the Hessian of $\overline{F}_N$ is non-negative definite, we take $w = (a,b) \in \mathbb{R}^2$ and write

$$\frac{1}{2N} w \cdot \nabla^2 \overline{F}_N(t,h) w = a^2 \mathbb{E}\langle (\xi \cdot \xi')^2 \rangle + b^2 \mathbb{E}\langle (y \cdot y')^2 \rangle + 2ab \mathbb{E}\langle (\xi \cdot \xi')(y \cdot y') \rangle \\ = \mathbb{E}\langle (a\xi \cdot \xi' + b y \cdot y')^2 \rangle \geq 0,$$

as desired. ∎

**Lemma 4.11.** *If $f : \mathbb{R}_{\geq 0} \times \mathbb{R}_{\geq 0} \to \mathbb{R}$ is a jointly convex and Lipschitz continuous function that is non-decreasing, then its symmetrization $\widetilde{f} : \mathbb{R}_{\geq 0} \times \mathbb{R} \to \mathbb{R}$ defined by $\widetilde{f}(t,h) := f(t,|h|)$ is also jointly convex and Lipschitz continuous, with the same Lipschitz constant.*

*Proof.* The uniform Lipschitz continuity of the symmetrization $\widetilde{f}$ is an immediate consequence of the uniform Lipschitz continuity of $f$ and the reverse triangle inequality. To show that the symmetrization is also jointly convex, fix $\alpha \in [0,1]$, $t,t' \geq 0$ and $h,h' \in \mathbb{R}$. By non-decreasingness of $f$,

$$\widetilde{f}(\alpha t + (1-\alpha)t', \alpha h + (1-\alpha)h') = f(\alpha t + (1-\alpha)t', |\alpha h + (1-\alpha)h'|) \\ \leq f(\alpha t + (1-\alpha)t', \alpha |h| + (1-\alpha)|h'|).$$

It follows by convexity of $f$ that

$$\widetilde{f}(\alpha t + (1-\alpha)t', \alpha h + (1-\alpha)h') \leq \alpha f(t,|h|) + (1-\alpha) f(t',|h'|) \\ = \alpha \widetilde{f}(t,h) + (1-\alpha) \widetilde{f}(t',h').$$

This completes the proof. ∎

We now proceed analogously to our treatment of the Curie-Weiss model and aim to estimate the error term in the approximate Hamilton-Jacobi equation (4.83) in terms of quantities involving the free energy itself. Compared with (3.7) or (3.17), a new term appears that relates to the fluctuations of the free energy itself.

**Proposition 4.12.** *The enriched free energy* (4.75) *in the symmetric rank-one matrix estimation problem satisfies the approximate Hamilton-Jacobi equation*

$$0 \leq \partial_t \overline{F}_N(t,h) - \left(\partial_h \overline{F}_N(t,h)\right)^2 \leq \frac{1}{N}\partial_h^2 \overline{F}_N(t,h) + \mathbb{E}\left(\partial_h F_N - \partial_h \overline{F}_N\right)^2(t,h). \quad (4.97)$$

*Proof.* The lower bound in (4.97) is immediate from (4.83) and the non-negativity of the variance. To establish the upper bound, by (4.83), it suffices to show that

$$\operatorname{Var}\left(\frac{x \cdot \bar{x}}{N}\right) \leq \frac{1}{N}\partial_h^2 \overline{F}_N(t,h) + \mathbb{E}\left(\partial_h F_N - \partial_h \overline{F}_N\right)^2(t,h). \quad (4.98)$$

As in the Curie-Weiss model, it is reasonable to expect this variance term to be related to $\partial_h^2 \overline{F}_N(t,h)$. Although the value of this derivative was stated in (4.94), we will compute it slightly differently here for convenience. We write

$$H'_N(h,x) := \frac{1}{\sqrt{2h}} z \cdot x + 2x \cdot \bar{x} - |x|^2, \quad (4.99)$$

so that the derivative of the free energy (4.74) prior to averaging may be expressed concisely as

$$\partial_h F_N(t,h) = \frac{1}{N}\langle H'_N \rangle.$$

Differentiating this expression reveals that

$$\partial_h^2 F_N(t,h) = \frac{1}{N}\langle \partial_h H'_N \rangle + \frac{1}{N}\langle (H'_N)^2 \rangle - \frac{1}{N}\langle H'_N \rangle^2$$

$$= \frac{1}{N}\langle (H'_N)^2 \rangle - \frac{1}{N}\langle H'_N \rangle^2 - \frac{1}{N(2h)^{3/2}}\langle z \cdot x \rangle.$$

Together with (4.81) and the Nishimori identity, this implies that

$$\partial_h^2 \overline{F}_N(t,h) = \frac{1}{N}\mathbb{E}\langle (H'_N)^2 \rangle - \frac{1}{N}\mathbb{E}\langle H'_N \rangle^2 - \frac{1}{2hN}\left(\mathbb{E}\langle |x|^2 \rangle - \mathbb{E}\langle x \cdot \bar{x} \rangle\right)$$

$$= \frac{1}{N}\mathbb{E}\langle (H'_N)^2 \rangle - \frac{1}{N}\mathbb{E}\langle H'_N \rangle^2 - \frac{1}{2hN}\left(\mathbb{E}\langle |x|^2 \rangle - \mathbb{E}|\langle x \rangle|^2\right).$$

Notice that the variance of $H'_N$ may be written as

$$\operatorname{Var}(H'_N) = \mathbb{E}\left\langle \left(H'_N - \langle H'_N \rangle\right)^2 \right\rangle + \mathbb{E}\left(\langle H'_N \rangle - \mathbb{E}\langle H'_N \rangle\right)^2$$

$$= \mathbb{E}\langle (H'_N)^2 \rangle - \mathbb{E}\langle H'_N \rangle^2 + N^2 \mathbb{E}\left(\partial_h F_N - \partial_h \overline{F}_N\right)^2(t,h).$$

Up to lower-order terms, the proof consists in showing that the variance of $x \cdot \bar{x}$ is bounded from above by the variance of $H'_N$, which in turn is essentially the right side of (4.97) up to scaling. To justify this precisely, we write

$$\frac{1}{N^2}\operatorname{Var}(H'_N) = \frac{1}{N}\partial_h^2 \overline{F}_N(t,h) + \mathbb{E}\left(\partial_h F_N - \partial_h \overline{F}_N\right)^2 + \frac{1}{2hN^2}\left(\mathbb{E}\langle |x|^2 \rangle - \mathbb{E}|\langle x \rangle|^2\right).$$

To relate this back to the variance term in (4.98), observe that by the derivative computation (4.82),

$$\frac{1}{N^2}\mathrm{Var}(H'_N) - \frac{1}{N^2}\mathrm{Var}(x\cdot\bar{x}) = \frac{1}{N^2}\mathbb{E}\langle (H'_N)^2\rangle - \frac{1}{N^2}\mathbb{E}\langle (x\cdot\bar{x})^2\rangle.$$

It follows that

$$\mathrm{Var}\left(\frac{x\cdot\bar{x}}{N}\right) \leq \frac{1}{N}\partial_h^2 \overline{F}_N(t,h) + \mathbb{E}\left(\partial_h F_N - \partial_h \overline{F}_N\right)^2(t,h)$$

$$+ \frac{1}{2hN^2}\mathbb{E}\langle |x|^2\rangle - \frac{1}{N^2}\mathbb{E}\langle (H'_N)^2\rangle + \frac{1}{N^2}\mathbb{E}\langle (x\cdot\bar{x})^2\rangle. \quad (4.100)$$

What is important for the sequel is to verify that the second line in this expression is of lower order in $N$ due to a cancellation between $\mathbb{E}\langle (H'_N)^2\rangle$ and $\mathbb{E}\langle (x\cdot\bar{x})^2\rangle$. But in fact we can show that this second line in (4.100) is non-positive. We observe that

$$\mathbb{E}\langle (H'_N)^2\rangle = \frac{1}{2h}\mathbb{E}\langle (z\cdot x)^2\rangle + \frac{2}{\sqrt{2h}}\mathbb{E}\langle z\cdot x(2x\cdot\bar{x} - |x|^2)\rangle$$

$$+ 4\mathbb{E}\langle x\cdot\bar{x}(x\cdot\bar{x} - |x|^2)\rangle + \mathbb{E}\langle |x|^4\rangle.$$

Fixing $i,j \in \{1,\ldots,N\}$, two applications of the Gaussian integration by parts formula in (4.29) reveal that for $i \neq j$,

$$\mathbb{E}\langle z_i z_j x_i x_j\rangle = \sqrt{2h}\mathbb{E}\langle z_j x_j(x_i^2 - x_i x_i')\rangle = 2h\mathbb{E}\langle (x_i^2 - x_i x_i')(x_j^2 + x_j x_j' - 2x_j x_j'')\rangle.$$

while for $i = j$,

$$\mathbb{E}\langle z_i^2 x_i^2\rangle = \sqrt{2h}\mathbb{E}\langle z_i x_i(x_i^2 - x_i x_i')\rangle + \mathbb{E}\langle x_i^2\rangle$$
$$= 2h\mathbb{E}\langle (x_i^2 - x_i x_i')(x_i^2 + x_i x_i' - 2x_i x_i'')\rangle + \mathbb{E}\langle x_i^2\rangle.$$

Together with the Nishimori identity, this shows that

$$\frac{1}{2h}\mathbb{E}\langle (z\cdot x)^2\rangle = \mathbb{E}\langle |x|^4\rangle - 2\mathbb{E}\langle |x|^2(x\cdot\bar{x})\rangle - \mathbb{E}\langle (x\cdot\bar{x})^2\rangle$$

$$+ 2\mathbb{E}\langle (x\cdot\bar{x})(x\cdot x')\rangle + \frac{1}{2h}\mathbb{E}\langle |x|^2\rangle.$$

A similar computation using the Gaussian integration by parts formula gives

$$\frac{2}{\sqrt{2h}}\mathbb{E}\langle z\cdot x(2x\cdot\bar{x} - |x|^2)\rangle = 2\mathbb{E}\langle (2x\cdot\bar{x} - |x|^2)(|x|^2 - (x\cdot x'))\rangle$$

$$= 6\mathbb{E}\langle |x|^2(x\cdot\bar{x})\rangle - 4\mathbb{E}\langle (x\cdot\bar{x})(x\cdot x')\rangle - 2\mathbb{E}\langle |x|^4\rangle.$$

It follows by the Cauchy-Schwarz inequality that

$$\mathbb{E}\langle (H'_N)^2 \rangle = 3\mathbb{E}\langle (x\cdot\bar{x})^2 \rangle - 2\mathbb{E}\langle (x\cdot\bar{x})(x\cdot x') \rangle + \frac{1}{2h}\mathbb{E}\langle |x|^2 \rangle$$

$$\geq \mathbb{E}\langle (x\cdot\bar{x})^2 \rangle + \frac{1}{2h}\mathbb{E}\langle |x|^2 \rangle.$$

This shows that the second line in (4.100) is non-positive and establishes the upper bound (4.98), thereby completing the proof. ∎

To control the error term on the right side of the approximate Hamilton-Jacobi equation (4.97), at first glance, it seems like we need to establish the concentration of the derivative of the enriched free energy (4.74) about its average. However, even in the simpler setting of the Curie-Weiss model with a random external field, we do not expect this concentration of the derivative to be valid at every single point. But in any case, we do not expect the term $N^{-1}\partial_h \overline{F}_N(t,h)$ in (4.97) to tend to zero at every single point $(t,h)$ either. The whole point of the notion of viscosity solution is that it suffices to control the size of the error term in (4.97) at appropriate contact points. In particular, we see from the proof of Theorem 3.25 that it suffices to establish the smallness of the error term at every point at which the free energy is touched from above by a smooth function. Leveraging convexity again, we will show that at such points, the variance of the derivative of the free energy can be controlled in terms of

$$\mathbb{E} \sup_{(t,h)\in [0,M]^2} |F_N(t,h) - \overline{F}_N(t,h)|^2, \qquad (4.101)$$

for an adequate choice of $M < +\infty$. Leveraging the Gaussian concentration inequality and a covering argument, we will show that (4.101) is essentially of order $N^{-1}$. We start by establishing a Gaussian-type estimate on the tail of the random variable in (4.101).

**Lemma 4.13.** *For each $M \in \mathbb{R}_{\geq 0}$, there exists a constant $C < +\infty$ such that for all $\lambda > CN^{-1/2}\sqrt{\log(N)}$, we have*

$$\mathbb{P}\left\{ \sup_{(t,h)\in [0,M]^2} |F_N(t,h) - \overline{F}_N(t,h)| \geq \lambda \right\} \leq \exp\left(-\frac{N\lambda^2}{C}\right). \qquad (4.102)$$

*Proof.* We write $C < +\infty$ for a constant that does not depend on $N$ or $\lambda$, but may depend on $M$, and whose value may change as we proceed through the proof. We introduce the random norm

$$\|W\|_* := \sup_{|x|\leq 1} |Wx|.$$

To run a covering argument, we begin by establishing the Hölder continuity of the free energy $\overline{F}_N$. We fix $t,t',h,h' \in [0,M]$ as well as $x$ in the support of $P_N$, and observe that by the Cauchy-Schwarz inequality and the boundedness of the support of $P_1$, we have

$$|H_N(t,h,x) - H_N(t',h',x)| \leq \left|\sqrt{\frac{2t}{N}} - \sqrt{\frac{2t'}{N}}\right| |x \cdot Wx| + \frac{1}{N}|t-t'|\left(2(x\cdot\bar{x})^2 + |x|^4\right)$$
$$+ \left|\sqrt{2h} - \sqrt{2h'}\right| |x \cdot z| + |h-h'|\left(2|x\cdot\bar{x}| + |x|^2\right)$$
$$\leq C\sqrt{N}\left(\left|\sqrt{t} - \sqrt{t'}\right|\|W\|_* + \left|\sqrt{h} - \sqrt{h'}\right||z|\right)$$
$$+ CN\left(|t-t'| + |h-h'|\right).$$

Together with the observation that for $y \geq y'$, we have $\sqrt{y} - \sqrt{y'} \leq |y-y'|^{1/2}$ as well as $|y-y'| \leq 2|y||y-y'|^{1/2}$, this bound on the Hamiltonian implies that the free energy (4.74) is Hölder continuous on $[0,M]^2$ with

$$|F_N(t,h) - F_N(t',h')| \leq C\left(|t-t'|^{1/2} + |h-h'|^{1/2}\right)X$$

for the random variable

$$X := 1 + \frac{\|W\|_*}{\sqrt{N}} + \frac{|z|}{\sqrt{N}}.$$

We also recall from Proposition 4.10 that the averaged free energy $\overline{F}_N$ is Lipschitz continuous, uniformly over $N$. We thus deduce that for every $\lambda > 0$ and $\varepsilon > 0$,

$$\mathbb{P}\left\{\sup_{[0,M]^2}|F_N(t,h) - \overline{F}_N(t,h)| \geq \lambda\right\} \leq \mathbb{P}\left\{\sup_{A_\varepsilon}|F_N(t,h) - \overline{F}_N(t,h)| \geq \lambda/2\right\}$$
$$+ \mathbb{P}\left\{X \geq \varepsilon^{-1/2}\lambda/C\right\}$$

for the set $A_\varepsilon := \varepsilon\mathbb{Z}^2 \cap [0,M]^2$. Indeed, every point $(t,h) \in [0,M]^2$ is at distance at most $\varepsilon$ from a point in $A_\varepsilon$. A union bound and the fact that the cardinality of $A_\varepsilon$ is bounded by $C\varepsilon^{-2}$ yield that this is further bounded by

$$C\varepsilon^{-2}\sup_{A_\varepsilon}\mathbb{P}\left\{|F_N(t,h) - \overline{F}_N(t,h)| \geq \lambda/2\right\} + \mathbb{P}\left\{X \geq \varepsilon^{-1/2}\lambda/C\right\}.$$

A simple extension of the proof of the free energy concentration inequality (4.59) yields that for every $\lambda \geq 0$ and $(t,h) \in [0,M]^2$,

$$\mathbb{P}\left\{|F_N(t,h) - \overline{F}_N(t,h)| \geq \lambda\right\} \leq 2\exp\left(-\frac{N\lambda^2}{C}\right).$$

On the other hand, Exercise 4.5 and (4.43) imply that for some constant $C' < +\infty$, and any $\varepsilon > 0$ and $\lambda > C'\sqrt{\varepsilon}$,

$$\mathbb{P}\{X \geq \varepsilon^{-1/2}\lambda/C\} \leq \exp\left(-\frac{N\lambda^2}{\varepsilon C}\right).$$

Combining these two bounds and choosing $\varepsilon = N^{-1}$ reveals that for any $\lambda > C'N^{-1/2}$,

$$\mathbb{P}\left\{\sup_{[0,M]^2} |F_N(t,h) - \overline{F}_N(t,h)| \geq \lambda\right\} \leq CN^2 \exp\left(-\frac{N\lambda^2}{C}\right).$$

Whenever $\lambda > C'N^{-1/2}\sqrt{\log(N)}$ for some sufficiently large constant $C' < +\infty$, the term $N^2$ can be absorbed into the exponential to obtain (4.102) and complete the proof. ∎

We are now ready to verify the key assumption of the convex selection principle besides the convexity property itself: that any subsequential limit of $\overline{F}_N$ must satisfy the equation on a dense subset. It is convenient to phrase this in terms of contact points, just like we did in Section 3.6, see in particular Corollary 3.24.

**Lemma 4.14.** *Let $f$ be any subsequential limit of the sequence $(\overline{F}_N)_{N \geq 1}$ of enriched free energies in the symmetric rank-one matrix estimation problem, and fix $t^*, h^* > 0$. If $\phi \in C^\infty(\mathbb{R}_{>0} \times \mathbb{R}; \mathbb{R})$ is a smooth function with the property that $f - \phi$ has a strict local maximum at $(t^*, h^*) \in \mathbb{R}_{>0} \times \mathbb{R}_{>0}$, then*

$$\left(\partial_t \phi - (\partial_h \phi)^2\right)(t^*, h^*) = 0 \tag{4.103}$$

*Proof.* Abusing notation, we do not denote the subsequence along which the convergence of $\overline{F}_N$ to $f$ occurs explicitly, in effect pretending that the convergence occurs along the whole sequence. Using Exercise 3.1, we find a sequence $(t_N, h_N)_{N \geq 1}$ converging to $(t^*, h^*)$ such that $\overline{F}_N - \phi$ has a local maximum at $(t_N, h_N)$; this exercise also guarantees that the neighbourhood over which $(t_N, h_N)$ is a local maximum can be chosen independently of $N$. To control the right side of the approximate Hamilton-Jacobi equation (4.97) at the contact point $(t_N, h_N)$, the idea will be to argue similarly to the proof of Proposition 2.15 where it is shown that whenever a sequence of convex functions converges, the sequence of derivatives also converges to the derivative of the limit at every point of differentiability of the limit. The subsequential limit $f$ of the sequence $(\overline{F}_N)_{N \geq 1}$ of convex functions will be differentiable at the contact point $(t^*, h^*)$. The sequence $(F_N)_{N \geq 1}$ also converges to $f$ in the sense of second moments by Lemma 4.13, and although $F_N$ itself is not necessarily convex, its second derivative is bounded from below by a term of order one which will suffice to control the error term in (4.97). The proof therefore proceeds in three steps. First we control the upper deviation of $\overline{F}_N$ from its tangent at $(t_N, h_N)$ by a parabola, then we control the lower deviation of $F_N$ from its tangent at $(t_N, h_N)$ by a random parabola, and finally we combine these two ingredients to control the right side of the approximate Hamilton-Jacobi equation (4.97). Throughout the proof, we understand that the constant $C < +\infty$ may change from one occurrence to the next, only making sure that it does not depend on $N$.

*Step 1: Hessian of $\overline{F}_N$ upper bound.* Since $\overline{F}_N - \phi$ has a local maximum at $(t_N, h_N)$ and $\phi$ is smooth, there exists $C < +\infty$ such that for every $h' \in \mathbb{R}$ with $|h'| \leq C^{-1}$,

$$\overline{F}_N(t_N, h_N + h') - \overline{F}_N(t_N, h_N) \leq \phi(t_N, h_N + h') - \phi(t_N, h_N)$$
$$\leq h' \partial_h \phi(t_N, h_N) + C|h'|^2.$$

For every $N$ sufficiently large, we have $h_N > 0$, and therefore $\partial_h(\overline{F}_N - \phi)(t_N, h_N) = 0$. It follows that for every $h' \in \mathbb{R}$ with $|h'| \leq C^{-1}$,

$$\overline{F}_N(t_N, h_N + h') - \overline{F}_N(t_N, h_N) \leq h' \partial_h \overline{F}_N(t_N, h_N) + C|h'|^2. \tag{4.104}$$

In particular, we have
$$\partial_h^2 \overline{F}_N(t_N, h_N) \leq C, \tag{4.105}$$
which is also immediate from the observation that $\partial_h^2(\overline{F}_N - \phi)(t_N, h_N) \leq 0$.

*Step 2: Hessian of $F_N$ lower bound.* Recalling the definition of the Hamiltonian $H'_N$ in (4.99), we can write

$$\partial_h^2 F_N(t_N, h_N + h') = \frac{1}{N}\langle (H'_N)^2 \rangle - \frac{1}{N}\langle H'_N \rangle^2 - \frac{1}{N(2(h_N + h'))^{3/2}}\langle z \cdot x \rangle.$$

Since $h_N$ converges to $h^* > 0$, it remains bounded away from zero for $N$ sufficiently large. Using also the non-negativity of the variance, the Cauchy-Schwarz inequality and the fact that the measure $P_1$ has bounded support, we deduce that for every $|h'| \leq C^{-1}$, we have

$$\partial_h^2 F_N(t_N, h_N + h') \geq -\frac{C|z|}{\sqrt{N}}.$$

It follows by Taylor's theorem that for every $|h'| \leq C^{-1}$,

$$F_N(t_N, h_N + h') - F_N(t_N, h_N) - h' \partial_h F_N(t_N, h_N) \geq -\frac{C|z|}{\sqrt{N}}|h'|^2. \tag{4.106}$$

*Step 3: controlling the right side of* (4.97). Combining (4.104) and (4.106) with the fact that $(t_N, h_N)_{N \geq 1}$ converges to $(t^*, h^*)$ gives a neighbourhood $V$ of $(t^*, h^*)$ with

$$h'\left(\partial_h F_N - \partial_h \overline{F}_N\right)(t_N, h_N) \leq 2\sup_V |F_N - \overline{F}_N| + C|h'|^2\left(1 + \frac{|z|}{\sqrt{N}}\right). \tag{4.107}$$

In particular, given a deterministic $\lambda \in [0, C^{-1}]$, the bound (4.107) for

$$h' = \lambda \frac{\partial_h F_N - \partial_h \overline{F}_N}{|\partial_h F_N - \partial_h \overline{F}_N|}(t_N, h_N)$$

reads
$$\lambda|\partial_h F_N - \partial_h \overline{F}_N|(t_N, h_N) \leq 2\sup_V |F_N - \overline{F}_N| + C\lambda^2\left(1 + \frac{|z|}{\sqrt{N}}\right).$$

Squaring both sides of this inequality, taking expectations and leveraging the concentration inequality in Lemma 4.13 yields
$$\lambda^2 \mathbb{E}(\partial_h F_N - \partial_h \overline{F}_N)^2(t_N, h_N) \leq 8\mathbb{E}\sup_V (F_N - \overline{F}_N)^2 + C\lambda^4 \mathbb{E}\left(1 + \frac{|z|}{\sqrt{N}}\right)^2$$
$$\leq \frac{C}{N^{1/2}} + C\lambda^4,$$

where we have used that $\mathbb{E}|z|^2 = N\mathbb{E}z_1^2 = N$. Choosing $\lambda := N^{-1/8}$ leads to the upper bound
$$\mathbb{E}(\partial_h F_N - \partial_h \overline{F}_N)^2(t_N, h_N) \leq \frac{C}{N^{1/4}}.$$

Substituting this into the approximate Hamilton-Jacobi equation in Proposition 4.12 and remembering the Hessian bound (4.105) gives
$$0 \leq \left(\partial_t \phi - (\partial_h \phi)^2\right)(t_N, h_N) = \left(\partial_t \overline{F}_N - (\partial_h \overline{F}_N)^2\right)(t_N, h_N) \leq \frac{C}{N} + \frac{C}{N^{1/4}}.$$

Letting $N$ tend to infinity and using the smoothness of $\phi$ completes the proof. ∎

We are now ready to prove Theorem 4.9 using the convex selection principle. We will appeal to the version of this result given in Lemma 3.23. We point out that under the additional assumption that $\mathbb{E}\overline{x}_1 = 0$, this proof can be simplified, as Step 3 becomes unnecessary, and we can instead rely on the simpler version of the convex selection principle given by Theorem 3.21. The key observation is that under this additional assumption, the symmetrized initial condition $\widetilde{\psi}$ is in $C^1(\mathbb{R}; \mathbb{R})$, by (S.27).

*Proof of Theorem 4.9.* The Arzelà-Ascoli theorem together with the derivative computations (4.77) and (4.82) and the boundedness of the support of $P_1$ imply that the sequence $(\overline{F}_N)_{N \geq 1}$ is precompact. Denoting by $f$ a subsequential limit, the idea is to apply the convex selection principle in Lemma 3.23 to the symmetrization $\widetilde{f}$ of $f$ and show that it is the unique viscosity solution to the Hamilton-Jacobi equation (4.87). The proof therefore proceeds in four steps. First, we show that $\widetilde{f}$ is jointly convex and Lipschitz continuous; then we prove that it satisfies the Hamilton-Jacobi equation (4.87) on a dense subset of $\mathbb{R}_{\geq 0} \times \mathbb{R}$, and that the initial condition $\widetilde{\psi} := \widetilde{f}(0, \cdot)$ satisfies the assumption in Lemma 3.23; and finally we conclude using the uniqueness of viscosity solutions and the Hopf-Lax formula.

*Step 1: $\widetilde{f}$ is jointly convex and Lipschitz continuous.* Each free energy $\overline{F}_N$ is jointly convex, Lipschitz continuous and non-decreasing by Proposition 4.10. Moreover, the Lipschitz constant of $\overline{F}_N$ is bounded uniformly over $N$. It follows by Lemma 4.11 that the symmetrization of each free energy $\overline{F}_N$ is also jointly convex and Lipschitz continuous uniformly over $N$. Since this sequence of symmetrizations converges to $\widetilde{f}$, the function $\widetilde{f}$ must also be jointly convex and Lipschitz continuous.

*Step 2: $\widetilde{f}$ satisfies* (4.87) *on a dense set.* Following the proof of Corollary 3.24, we deduce from Lemma 4.14 that the set

$$\mathcal{A} := \left\{ (t,h) \in \mathbb{R}_{>0} \times \mathbb{R}_{>0} \mid f \text{ is differentiable at } (t,h) \right.$$
$$\left. \text{and } \left( \partial_t f - (\partial_h f)^2 \right)(t,h) = 0 \right\}$$

is dense in $\mathbb{R}_{>0} \times \mathbb{R}_{>0}$. We also observe that if $t > 0$ and $h < 0$ are such that $(t, |h|)$ is a point of differentiability of $f$, then $\widetilde{f}$ is differentiable at $(t,h)$ and we have

$$\partial_t \widetilde{f}(t,h) = \partial_t f(t,|h|) \quad \text{and} \quad \left( \partial_h \widetilde{f}(t,h) \right)^2 = \left( \partial_h f(t,|h|) \right)^2.$$

This implies that $\widetilde{f}$ satisfies the Hamilton-Jacobi equation (4.87) on a dense subset of $\mathbb{R}_{\geqslant 0} \times \mathbb{R}$.

*Step 3: $\widetilde{\psi}$ satisfies the assumption of Lemma 3.23.* We first verify that the initial condition $\widetilde{\psi}$ satisfies the assumption of Lemma 3.23 at any point $h \neq 0$. That is, for every $h \neq 0$ and $p \in \partial \widetilde{\psi}(h)$, we find $b \in \mathbb{R}$ such that $(b,p) \in \partial \widetilde{f}(0,h)$ and $b - p^2 \geqslant 0$. By the results of the previous steps and Lemma 3.22, there exists $(a,q) \in \partial \widetilde{f}(0,h)$ such that $a - q^2 = 0$. We must have in particular that $q \in \partial \widetilde{\psi}(h)$, so by Theorem 2.13, we have $q = p$, and the choice of $b = a = p^2$ satisfies our requirements. For future reference, we observe that for every $h > 0$, we thus have that $p = \partial_h \psi(h)$ and $(\partial_h \psi(h)^2, \partial_h \psi(h)) \in \partial \widetilde{f}(0,h)$. Letting $h > 0$ tend to zero and using Proposition 2.14, we deduce that

$$\left( \partial_h \psi(0)^2, \partial_h \psi(0) \right) \in \partial \widetilde{f}(0,0). \tag{4.108}$$

We now consider the case $h = 0$. Letting $p \in \partial \widetilde{\psi}(0)$, we look for $b \in \mathbb{R}$ such that $(b,p) \in \partial \widetilde{f}(0,0)$ and $b - p^2 \geqslant 0$. Since $p \in \partial \widetilde{\psi}(0)$, we have that for every $y \in \mathbb{R}$,

$$\widetilde{\psi}(y) - \widetilde{\psi}(0) = \psi(|y|) - \psi(0) \geqslant py.$$

Letting $y$ tend to zero from the left and from the right of zero, we see that this implies that $|p| \leqslant \partial_h \psi(0)$. Letting $b := \partial_h \psi(0)^2$, we therefore have that

$$b - p^2 \geqslant b - |\partial_h \psi(0)|^2 = 0.$$

It remains to verify that $(b,p) \in \partial \widetilde{f}(0,0)$. By (4.108), for every $(t,h) \in \mathbb{R}_{\geqslant 0} \times \mathbb{R}$, we have that

$$\widetilde{f}(t,h) \geqslant \widetilde{f}(0,0) + tb + h \partial_h \psi(0).$$

Using that $\widetilde{f}(t,-h) = \widetilde{f}(t,h)$, and then that $|p| \leq \partial_h \psi(0)$, we deduce that for every $(t,h) \in \mathbb{R}_{\geq 0} \times \mathbb{R}$,

$$\widetilde{f}(t,h) \geq \widetilde{f}(0,0) + tb + |h|\partial_h \psi(0) \geq \widetilde{f}(0,0) + tb + hp.$$

This shows that $(b,p) \in \partial \widetilde{f}(0,0)$, as desired.

*Step 4: conclusion.* The results of the previous steps ensure that we can invoke the convex selection principle in Lemma 3.23 to conclude that $\widetilde{f}$ is a viscosity solution to the Hamilton-Jacobi equation (4.87) with initial condition $\widetilde{\psi}$. Together with the uniqueness result in Corollary 3.7, this implies that the symmetrization of the sequence $(\overline{F}_N)_{N \geq 1}$ converges to the unique viscosity solution $\widetilde{f}$ to the Hamilton-Jacobi equation (4.87) with initial condition $\widetilde{\psi}$. In particular, the sequence $(\overline{F}_N)_{N \geq 1}$ of functions from $\mathbb{R}_{\geq 0} \times \mathbb{R}_{\geq 0}$ to $\mathbb{R}$ converges to some limit $f$, which is the restriction of $\widetilde{f}$ to $\mathbb{R}_{\geq 0} \times \mathbb{R}_{\geq 0}$. The Hopf-Lax formula and the computation for the convex dual of the square function in Exercise 2.11 imply that for all $t, h \geq 0$,

$$f(t,h) = \widetilde{f}(t,h) = \sup_{h' \in \mathbb{R}} \left( \psi(|h'|) - \frac{(h'-h)^2}{4t} \right).$$

Using that $h \geq 0$, we see that it is always better to choose $h' \geq 0$ in the optimization problem above, since for every $h' \geq 0$,

$$\psi(|-h'|) - \frac{(-h'-h)^2}{4t} \leq \psi(|h'|) - \frac{(h'-h)^2}{4t},$$

and we thus obtain that

$$f(t,h) = \sup_{h' \geq 0} \left( \psi(h') - \frac{(h'-h)^2}{4t} \right) = \sup_{h' \geq -h} \left( \psi(h+h') - \frac{(h')^2}{4t} \right). \quad (4.109)$$

Recall also from Proposition 4.10 that the function $\psi$ is non-decreasing. In the supremum on the right side of (4.109), we may therefore restrict the supremum to $h' \geq 0$, since when $h' \in [-h, 0]$, we have

$$\psi(h+h') - \frac{(h')^2}{4t} \leq \psi(h).$$

We have thus shown (4.90). Remembering (4.88) and evaluating this at $h = 0$ gives the representation (4.91) and completes the proof. ∎

We can now combine the formula (4.91) for the limit free energy in the symmetric rank-one matrix estimation problem obtained in Theorem 4.9 with Propositions 4.3 and 4.4 to determine the two information-theoretic quantities introduced at the beginning of this chapter, the minimal mean-square error (4.2) and the mutual information (4.4). It will be convenient to denote by

$$f(t) := \sup_{h \geqslant 0} \left( \psi(h) - \frac{h^2}{4t} \right) \qquad (4.110)$$

the limit free energy in the symmetric rank-one matrix estimation problem, and to write $\mathcal{D}$ for its points of differentiability,

$$\mathcal{D} := \{ t \in \mathbb{R}_{\geqslant 0} \mid f \text{ is differentiable at } t \}. \qquad (4.111)$$

By Rademacher's theorem (Theorem 2.10), the set $\mathcal{D}$ is dense in $\mathbb{R}_{\geqslant 0}$, and by the envelope theorem (Theorem 2.21), if $h_*(t) \geqslant 0$ is any maximizer of the right side of (4.110), then we have for every $t \in \mathcal{D}$ that

$$\partial_t f(t) = \frac{h_*^2(t)}{4t^2}. \qquad (4.112)$$

Together with the derivative computation (4.77), this allows us to express the minimal mean-square error in terms of any maximizer $h_*(t)$. We also recall from the envelope theorem that for every $t \geqslant 0$, if there exists a unique maximizer to the right side of (4.110), then $t \in \mathcal{D}$.

**Proposition 4.15.** *At every point of differentiability $t \in \mathcal{D}$ of the limit free energy* (4.110), *the limit minimal mean-square error is given by*

$$\mathsf{mmse}(t) := \lim_{N \to +\infty} \mathsf{mmse}_N(t) = \left( \mathbb{E}|\bar{x}_1|^2 \right)^2 - \frac{h_*^2(t)}{4t^2}, \qquad (4.113)$$

*where $h_*(t)$ denotes any maximizer of the right side of* (4.110).

*Proof.* Proposition 4.3, the asymptotic expansion (4.19), and the derivative computation (4.77) imply that for every $t \geqslant 0$,

$$\mathsf{mmse}_N(t) = \left( \mathbb{E}|\bar{x}_1|^2 \right)^2 - \frac{1}{N^2} \mathbb{E} \langle (x \cdot \bar{x})^2 \rangle + o(1) = \left( \mathbb{E}|\bar{x}_1|^2 \right)^2 - \partial_t \overline{F}_N(t, 0) + o(1).$$

Together with Proposition 2.15, Theorem 2.13 and the consequence (4.112) of the envelope theorem, this implies that at any point of differentiability $t \in \mathcal{D}$, we have

$$\mathsf{mmse}(t) = \left( \mathbb{E}|\bar{x}_1|^2 \right)^2 - \partial_t f(t, 0) = \left( \mathbb{E}|\bar{x}_1|^2 \right)^2 - \frac{h_*^2(t)}{4t^2},$$

as required. ∎

A variational formula for the limit mutual information can be obtained directly from the variational representation (4.110) of the limit free energy.

**Proposition 4.16.** *For every $t \geqslant 0$, the limit mutual information is given by*

$$\mathsf{I}(t) := \lim_{N \to +\infty} \mathsf{I}_N(t) = t\left(\mathbb{E}|\bar{x}_1|^2\right)^2 - \sup_{h \geqslant 0}\left(\psi(h) - \frac{h^2}{4t}\right). \qquad (4.114)$$

*Proof.* Proposition 4.4 and the asymptotic expansion (4.19) imply that

$$\mathsf{I}_N(t) = t\left(\mathbb{E}|x_1|^2\right)^2 - \overline{F}_N^\circ(t) + o(1).$$

Letting $N$ tend to infinity and using the variational formula (4.110) for the limit free energy completes the proof. ∎

Before moving to the next section, we mention more general models that can also be handled by the techniques we used to analyze the symmetric rank-one matrix estimation problem.

A first natural generalization consists in relaxing the symmetry assumption, and in considering matrices of higher rank. For illustration, we may imagine a recommendation system collecting ratings for each customer and product. One may model the fundamental characteristics of a customer by a vector with $K$ coordinates, for some $K$ relatively small, say $K \simeq 10$, and represent the list of customers as an $M$-by-$K$ matrix $u$, and where we can postulate that the rows representing the customers are independent and identically distributed. We also assume that the products on sale can be similarly represented by an $N$-by-$K$ matrix $v$, with i.i.d. rows representing each product, in such a way that the ratings that users attribute to the products are essentially a noisy version of the $M$-by-$N$ matrix $uv^*$. We would like to retrieve information about the matrix $uv^*$ given the noisy observations we have of it, in the regime in which $M$ and $N$ are both large and of the same order of magnitude.

We can generalize this model further, so that we can represent, for instance, a situation in which we have multiple different observations of the interactions between the rows of $u$ and $v$, or we have higher-order interactions. To encode this larger class of problems, we fix integers $p, K, L \geqslant 1$, and a matrix $A \in \mathbb{R}^{K^p \times L}$ that are deterministic (and known to the observer). Let $X$ be an $N$-by-$K$ random matrix, and denote by $X^{\otimes p}$ the $p$-fold tensor product of $X$, which we think of as an $N^p$-by-$K^p$ matrix. The more general model consists in observing a noisy version of the matrix $X^{\otimes p}A$, and asking again whether meaningful information can be recovered about $X^{\otimes p}A$ or $X$ from this observation, in the regime of large $N$. Under mild assumptions on the law of $X$ (e.g. the rows are i.i.d. and almost surely bounded in norm by a fixed constant), the techniques that we have used to analyze the symmetric rank-one matrix estimation problem allow us to answer this more general question precisely [75].

We have seen in Chapter 3, and in particular in Proposition 3.20, that the standard Curie-Weiss model can be analyzed more directly, without appealing to the convex selection principle. The inequality on the left side of (4.97) can be used in a similar way to analyze the symmetric rank-one matrix estimation problem. We chose to present a proof based on the convex selection principle because this is the only proof we are aware of that is capable of handling the more general inference problems discussed in the previous paragraph. To obtain this generalization, the main difficulty is that the variable $h \geqslant 0$ must be replaced by a positive semi-definite matrix, so the symmetrization trick we have used to handle the domain $\mathbb{R}_{\geqslant 0}$ no longer works. In this case, we must therefore make sense of Hamilton-Jacobi equations posed on domains with boundaries; we refer to [74, 75, 103] for more on this. The non-linearity appearing in the Hamilton-Jacobi equation is the mapping $q \mapsto (AA^*) \cdot q^{\otimes p}$. This mapping is not necessarily convex, so the Hopf-Lax formula no longer applies in this case. But the free energy is still jointly convex, and the Hopf representation of solutions can be adapted to this context with a boundary. We thus end up with a variational representation of the limit free energy, and can thereby conveniently analyze the minimal mean-square error and mutual information of the model.

**Exercise 4.9.** The goal of this exercise is to prove the existence of a "phase transition" for the limit of the minimal mean-square error (4.2) as we vary the signal-to-noise ratio, as was announced in the introduction to this chapter. Notice that the mean-square error of the null estimator (0.9) converges to $(\mathbb{E}|\bar{x}_1|^2)^2$ as $N$ tends to infinity. Assuming that $\mathbb{E}\bar{x}_1 = 0$ and that $\bar{x}_1$ is not constant, show that there exists a critical parameter $t_c \in \mathbb{R}_{>0}$ such that

$$\text{for every } t < t_c, \quad \lim_{N \to +\infty} \operatorname{mmse}_N(t) = \left(\mathbb{E}|\bar{x}_1|^2\right)^2, \quad \text{while} \quad (4.115)$$

$$\text{for every } t > t_c, \quad \limsup_{N \to +\infty} \operatorname{mmse}_N(t) < \left(\mathbb{E}|\bar{x}_1|^2\right)^2. \quad (4.116)$$

(What happens at $t = t_c$ depends on the law of $\bar{x}_1$.)

**Exercise 4.10.** Let $(t, h) \in \mathbb{R}_{>0} \times \mathbb{R}_{>0}$ be a point of differentiability of the limit free energy $f : \mathbb{R}_{\geqslant 0} \times \mathbb{R}_{\geqslant 0} \to \mathbb{R}$ in the symmetric rank-one matrix estimation problem. Denote by $h_*(t, h)$ any maximizer of the right side of (4.90). Prove that

$$\lim_{N \to +\infty} \frac{1}{N} \mathbb{E}|\bar{x} - \mathbb{E}[\bar{x}|\mathcal{Y}]|^2 = \mathbb{E}\bar{x}_1^2 + \frac{h - h_*(t, h)}{2t} - \partial_h \psi(h_*(t, h)). \quad (4.117)$$

## 4.4  Comparison with a concrete algorithm

In this section, we discuss a classical and computationally efficient method of constructing an estimator for the symmetric rank-one matrix $\bar{x}\bar{x}^*$ from its noisy observation $Y$ known as principal component analysis (PCA). The PCA estimator $\widehat{\bar{x}\bar{x}^*}$ of the rank-one matrix $\bar{x}\bar{x}^*$ is constructed by finding an eigenvector $v \in \mathbb{R}^N$ with $|v|^2 = N$ associated with the top eigenvalue of the symmetric matrix $Y + Y^* \in \mathbb{R}^{N \times N}$, and setting

$$\widehat{\bar{x}\bar{x}^*} := \lambda_0 v v^* \qquad (4.118)$$

for the choice of $\lambda_0 \geq 0$ that minimizes the mean-square error

$$\mathsf{mse}_N := \frac{1}{N^2} \inf_{\lambda \geq 0} \mathbb{E}|\bar{x}\bar{x}^* - \lambda v v^*|^2. \qquad (4.119)$$

By comparing the limit mean-square error of the PCA estimator to the limit minimal mean-square error computed in Proposition 4.15 for a Gaussian prior, a Rademacher prior and a sparse prior, we will see that the performance of PCA depends greatly on the identity of this prior. By "prior", we mean the law of the signal $P_N = (P_1)^{\otimes N}$. In this section, we always normalize this measure so that $\mathbb{E}\bar{x}_1 = 0$ and $\mathbb{E}\bar{x}_1^2 = 1$.

To study the performance of the PCA estimator, we start by determining its asymptotic mean-square error. We will rely on a result in random matrix theory which states that, in probability,

$$\lim_{N \to +\infty} \frac{|\bar{x} \cdot v|^2}{|\bar{x}|^2 |v|^2} = \max\left(0, 1 - \frac{1}{4t}\right). \qquad (4.120)$$

We will take this convergence for granted and refer the interested reader to Section 3.1 in [52] for its proof in a more general setting (see also [28]).

**Proposition 4.17.** *For every $t \geq 0$, the asymptotic mean-square error associated with the PCA estimator is given by*

$$\mathsf{mse}(t) := \lim_{N \to +\infty} \mathsf{mse}_N(t) = \begin{cases} 1 & \text{if } t \leq \frac{1}{4}, \\ \frac{1}{4t}\left(2 - \frac{1}{4t}\right) & \text{if } t \geq \frac{1}{4}. \end{cases} \qquad (4.121)$$

*Proof.* We start by writing, for any $\lambda \geq 0$,

$$\mathbb{E}|\bar{x}\bar{x}^* - \lambda v v^*|^2 = \mathbb{E}|\bar{x}|^4 - 2\lambda \mathbb{E}(\bar{x} \cdot v)^2 + \lambda^2 \mathbb{E}|v|^4.$$

Optimizing over $\lambda$ and recalling that $|v|^2 = N$ reveals that the optimal $\lambda$ is

$$\lambda_0 := \frac{\mathbb{E}(\bar{x} \cdot v)^2}{\mathbb{E}|v|^4} = \frac{\mathbb{E}(\bar{x} \cdot v)^2}{N^2}.$$

It follows by the asymptotic expansion (4.19) and the assumption $\mathbb{E}\bar{x}_1^2 = 1$ that

$$\mathsf{mse}_N(t) = \frac{\mathbb{E}|\bar{x}|^4}{N^2} - \left(\frac{\mathbb{E}(\bar{x} \cdot v)^2}{N^2}\right)^2 = 1 - \left(\mathbb{E}\frac{|\bar{x} \cdot v|^2}{|\bar{x}|^2 |v|^2} \frac{|\bar{x}|^2}{N}\right)^2 + o(1).$$

By the law of large numbers, we have that $|\bar{x}|^2/N$ converges to 1 almost surely as $N$ tends to infinity. Together with (4.120), this gives the convergence in probability

$$\lim_{N\to+\infty} \frac{|x\cdot v|^2}{|\bar{x}|^2|v|^2} \frac{|\bar{x}|^2}{N} = \max\left(0, 1 - \frac{1}{4t}\right).$$

Remembering that the support of $P_1$ is bounded and invoking the dominated convergence theorem for convergence in probability reveals that

$$\lim_{N\to+\infty} \mathsf{mse}_N(t) = 1 - \max\left(0, 1 - \frac{1}{4t}\right)^2.$$

Expanding and distinguishing cases completes the proof.   ∎

This result implies that the critical signal-to-noise ratio of the PCA estimator is always $t^*_{\text{PCA}} = \frac{1}{4}$. It does not depend on the prior, except to the extent that we have imposed $\mathbb{E}\bar{x}_1 = 0$ and $\mathbb{E}\bar{x}_1^2 = 1$. By considering three different priors, we will see that the PCA estimator can sometimes be optimal, sometimes be sub-optimal but with a critical signal-to-noise ratio coinciding with the optimal one, and sometimes be sub-optimal and with a critical signal-to-noise ratio differing from the optimal one.

In the context of the Gaussian prior

$$dP_1(x) := \frac{1}{\sqrt{2\pi}} \exp\left(-\frac{x^2}{2}\right) dx, \tag{4.122}$$

the PCA estimator is optimal.

**Proposition 4.18.** *For the Gaussian prior* (4.122), *the limit minimal mean-square error is given by*

$$\mathsf{mmse}(t) = \begin{cases} 1 & \text{if } t \leq \frac{1}{4}, \\ \frac{1}{4t}\left(2 - \frac{1}{4t}\right) & \text{if } t \geq \frac{1}{4}. \end{cases} \tag{4.123}$$

*Partial proof.* Although the Gaussian prior (4.122) does not have bounded support, the formula (4.110) for the limit free energy still holds in this setting. This can be shown by approximation, but we will simply admit it. In the Gaussian setting, the initial condition $\psi$ can be computed explicitly. Recall that

$$\psi(h) = \mathbb{E}\log \int_{-\infty}^{+\infty} \exp\left(\sqrt{2h}xz_1 + 2hx\bar{x}_1 - hx^2\right) dP_1(x),$$

and let $a := \sqrt{2h}z_1 + 2h\bar{x}_1$ and $b := 1 + 2h$, in such a way that the innermost integral becomes

$$\frac{1}{\sqrt{2\pi}} \int_{-\infty}^{+\infty} e^{ax - \frac{bx^2}{2}} dx = \frac{e^{\frac{a^2}{2b}}}{\sqrt{2\pi}} \int_{-\infty}^{+\infty} e^{-\frac{(x-\frac{a}{b})^2}{2/b}} dx = \frac{e^{\frac{a^2}{2b}}}{\sqrt{b}},$$

where in the last equality we have identified the density of a Gaussian random variable with mean $a/b$ and variance $1/b$. It follows by independence of $z_1$ and $\bar{x}_1$ that

$$\psi(h) = \frac{1}{2(1+2h)}\mathbb{E}\left(2hz_1^2 + 4h\sqrt{2h}z_1\bar{x}_1 + 4h^2\bar{x}_1^2\right) - \frac{1}{2}\log(1+2h)$$

$$= h - \frac{1}{2}\log(1+2h).$$

Any maximizer $h_*(t)$ of (4.110) must therefore satisfy the critical point equation

$$1 - \frac{1}{1+2h_*(t)} = \frac{h_*(t)}{2t}.$$

Rearranging reveals that $h_*(t) = 2t - \frac{1}{2}$ or $h_*(t) = 0$. If $t \leq \frac{1}{4}$, the only non-negative solution is $h_*(t) = 0$. This implies that the limit free energy is differentiable at $t$ and, by Proposition 4.15,

$$\mathsf{mmse}(t) = 1 - \frac{h_*^2(t)}{4t^2} = 1.$$

On the other hand, if $t > \frac{1}{4}$, then

$$\psi\left(2t - \frac{1}{2}\right) - \frac{(2t-\frac{1}{2})^2}{4t} = t - \frac{1}{2}\log(4t) - \frac{1}{16t} > \frac{1}{4} - \frac{1}{2}\log(1) - \frac{1}{4} = 0$$

so (4.110) is maximized by $h_*(t) = 2t - \frac{1}{2}$. Together with the envelope theorem (Theorem 2.21), this shows that for $t > \frac{1}{4}$, the free energy is differentiable at $t$, and

$$\mathsf{mmse}(t) = 1 - \frac{h_*^2(t)}{4t^2} = 1 - \frac{(2t-\frac{1}{2})^2}{4t^2} = \frac{2t}{4t^2} - \frac{1}{16t^2} = \frac{1}{4t}\left(2 - \frac{1}{4t}\right).$$

This completes the proof. ∎

In the context of a Rademacher prior

$$P_1 := \frac{1}{2}\delta_1 + \frac{1}{2}\delta_{-1} \tag{4.124}$$

the PCA estimator is sub-optimal but identifies the correct critical signal-to-noise ratio. In other words, the PCA estimator recovers part of the signal exactly in the region where this is feasible, but it does not recover as much of the signal as the optimal estimator. This is illustrated in Figure 4.1. While we expect the PCA estimator to be sub-optimal in the entire regime $t > \frac{1}{4}$, we show this rigorously only for $t$ sufficiently large.

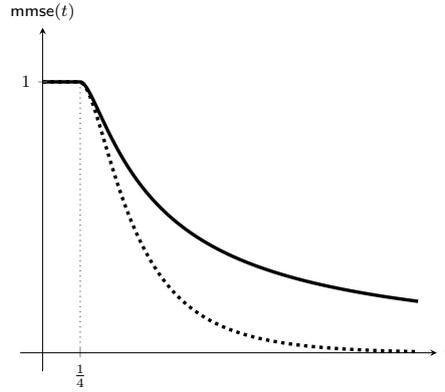

**Figure 4.1** The dotted line is the minimal mean-square error for the Rademacher prior (4.124) while the solid line is the mean-square error achieved using PCA.

**Proposition 4.19.** *Fix the Rademacher prior* (4.124). *If* $t < \frac{1}{4}$, *then* $\mathsf{mmse}(t) = 1$. *If* $t$ *is large enough and is a point of differentiability* $t \in \mathcal{D}$ *of the limit free energy* (4.110), *then*

$$\mathsf{mmse}(t) < \frac{1}{4t}\left(2 - \frac{1}{4t}\right). \tag{4.125}$$

**Remark 4.20.** For $t$ sufficiently large, a minor modification of the proof of Proposition 4.19 also shows that

$$\limsup_{N \to +\infty} \mathsf{mmse}_N(t) < \frac{1}{4t}\left(2 - \frac{1}{4t}\right),$$

without having to impose that $t \in \mathcal{D}$.

*Proof of Proposition 4.19.* The proof proceeds in two steps. Writing $h_*(t)$ for any maximizer of (4.110), first we show that $h_*(t) = 0$ for $t \leq \frac{1}{4}$ while $h_*(t) > 0$ for $t > \frac{1}{4}$, and then we prove (4.125) by showing that $h_*(t)$ is essentially $2t$ for $t$ large enough.

*Step 1: $h_*(t) = 0$ for $t \leq \frac{1}{4}$ while $h_*(t) > 0$ for $t > \frac{1}{4}$.* The derivative computation (4.82) for $N = 1$, the Nishimori identity, and the explicit form of the Rademacher prior (4.124) imply that

$$\partial_h \psi(h) = \mathbb{E}\langle x \cdot \bar{x}_1 \rangle = \mathbb{E}\langle x \rangle^2 = \mathbb{E}\tanh^2\left(\sqrt{2h}z_1 + 2h\bar{x}_1\right).$$

Averaging with respect to the randomness of $\bar{x}_1$, and using the symmetry of $\tanh^2$ and $z_1$ reveals that

$$\partial_h \psi(h) = \frac{1}{2}\mathbb{E}\tanh^2\left(\sqrt{2h}z_1 + 2h\right) + \frac{1}{2}\mathbb{E}\tanh^2\left(\sqrt{2h}z_1 - 2h\right) = \mathbb{E}\tanh^2\left(\sqrt{2h}z_1 + 2h\right).$$

Any maximizer $h_*(t)$ of (4.110) must therefore satisfy the equation

$$h_*(t) = 2t\mathbb{E}\tanh^2\left(\sqrt{2h_*(t)}z_1 + 2h_*(t)\right). \tag{4.126}$$

The second derivative computation (4.94), the Nishimori identity, and the fact that $x^2 = 1$ for any $x$ in the support of $P_1$ imply that

$$\partial_h^2 \psi(h) = 2\mathbb{E}\langle(x\bar{x})^2\rangle - 4\mathbb{E}\langle x\bar{x}xx'\rangle + 2\mathbb{E}\langle x\bar{x}\rangle^2 = 2 - 4\mathbb{E}\langle x\rangle^2 + 2\mathbb{E}\langle x\rangle^4 = 2\mathbb{E}\left(1 - \langle x\rangle^2\right)^2.$$

This means that $0 \leqslant \partial_h^2 \psi(h) \leqslant 2$. Together with the fact that $\partial_h \psi(0) = \psi(0) = 0$, this implies that $0 \leqslant \psi(h) \leqslant h^2$. It follows that for $t \leqslant \frac{1}{4}$ and $h > 0$,

$$\psi(h) - \frac{h^2}{4t} \leqslant h^2\left(1 - \frac{1}{4t}\right) \leqslant 0.$$

This means that for $t < \frac{1}{4}$, we have $h_*(t) = 0$. On the other hand, if $t > \frac{1}{4}$ and $h > 0$ is small enough, then

$$\psi(h) - \frac{h^2}{4t} = h^2\left(1 - \frac{1}{4t}\right) + o(h^2) > 0,$$

so $h_*(t) > 0$.

*Step 2: $h_*(t)$ is essentially $2t$ for $t$ large.* For any $x \in \mathbb{R}$,

$$\left|1 - \tanh^2(x)\right| = 1 - \left(\frac{e^{2x} - 1}{e^{2x} + 1}\right)^2 = \frac{4e^{2x}}{e^{4x} + 2e^{2x} + 1} \leqslant \frac{4}{e^{2x} + e^{-2x}} \leqslant 4e^{-2|x|}.$$

It follows by the critical point equation (4.126) that

$$\left|h_*(t) - 2t\right| \leqslant 2t\mathbb{E}\left|\tanh^2\left(\sqrt{2h_*(t)}z_1 + 2h_*(t)\right) - 1\right|$$
$$\leqslant 8t\mathbb{E}\exp\left(-2\left|\sqrt{2h_*(t)}z_1 + 2h_*(t)\right|\right).$$

To bound this further, introduce the set $A_{h_*(t)} := \{|\sqrt{2h_*(t)}z_1| > h_*(t)\}$, and observe that

$$\left|h_*(t) - 2t\right| \leqslant 8t\mathbb{P}\left(A_{h_*(t)}\right) + 8t\mathbb{E}\exp\left(-2\left|\sqrt{2h_*(t)}z_1 + 2h_*(t)\right|\right)\mathbf{1}_{A_{h_*(t)}^c}$$
$$\leqslant 8t\mathbb{P}\left(A_{h_*(t)}\right) + 8t\exp\left(-2h_*(t)\right).$$

Since $z_1$ is a Gaussian random variable, for any $c > 0$, we have

$$\mathbb{P}\{|z_1| > c\} = \frac{2}{\sqrt{2\pi}}\int_c^{+\infty}\exp\left(-\frac{x^2}{2}\right)dx \leqslant \sqrt{\frac{2}{\pi}}\int_c^{+\infty}\frac{x}{c}\exp\left(-\frac{x^2}{2}\right)dx$$
$$\leqslant \frac{\sqrt{2}}{c\sqrt{\pi}}\exp\left(-\frac{c^2}{2}\right).$$

Applying this with $c = \sqrt{\frac{h_*(t)}{2}}$ to bound $\mathbb{P}(A_{h_*(t)})$ reveals that

$$|h_*(t) - 2t| \leq \frac{16t}{\sqrt{\pi h_*(t)}} \exp\left(-\frac{h_*(t)}{4}\right) + 8t \exp(-2h_*(t)).$$

One can easily check that any choice of maximizer $h_*(t)$ of (4.110) must diverge to infinity as $t$ tends to infinity. Using that $h_*(t)$ exceeds a large constant for $t$ sufficiently large, we deduce from the previous display that $h_*(t)$ in fact grows linearly with $t$, and appealing again to the previous display, we obtain that for $t$ sufficiently large, we have

$$|h_*(t) - 2t| \leq \exp(-t/4).$$

It follows by Proposition 4.15 that for $t$ large enough,

$$\mathsf{mmse}(t) = 1 - \frac{h_*^2(t)}{4t^2} \leq 1 - \frac{(2t - e^{-t/4})^2}{4t^2} = \frac{e^{-t/4}}{2t} - \frac{e^{-t/2}}{4t^2}.$$

Increasing $t$ if necessary to ensure that $e^{-t/2} \leq \frac{1}{4} e^{-t/4}$ gives the strict upper bound

$$\mathsf{mmse}(t) \leq \frac{e^{-t/4}}{4t}\left(2 - \frac{1}{4t}\right) < \frac{1}{4t}\left(2 - \frac{1}{4t}\right)$$

and completes the proof. ∎

In the context of a sparse prior defined for some $p \in (0,1)$ by

$$P_1 := (1-p)\delta_0 + \frac{p}{2}\delta_{1/\sqrt{p}} + \frac{p}{2}\delta_{-1/\sqrt{p}}, \qquad (4.127)$$

the PCA estimator can be sub-optimal in a stronger sense: for $p > 0$ sufficiently small, there is a region of signal-to-noise ratios for which a non-trivial recovery of the signal is possible but for which the PCA estimator is trivial.

**Proposition 4.21.** *For the sparse prior* (4.124) *with $p > 0$ small enough, there exists a signal-to-noise ratio $t^*(p) < \frac{1}{4}$ such that for every point of differentiability $t \in \mathcal{D}$ of the limit free energy* (4.110) *with $t > t^*(p)$, we have*

$$\mathsf{mmse}(t) < 1. \qquad (4.128)$$

**Remark 4.22.** One can also show that for every $t > t^*(p)$, we have

$$\limsup_{N \to +\infty} \mathsf{mmse}_N(t) < 1,$$

without having to impose that $t \in \mathcal{D}$.

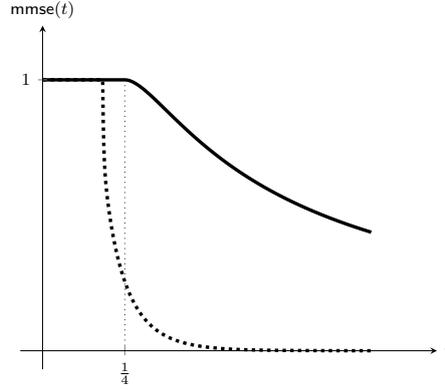

**Figure 4.2** The dotted line is the minimal mean-square error for the sparse prior (4.127) with $p = 0.05$ while the solid line is the mean-square error achieved using PCA.

*Proof of Proposition 4.21.* A direct computation using the explicit form of the sparse prior (4.124) shows that

$$\psi(h) = \mathbb{E}\log\left((1-p) + \frac{p}{2}\exp\left(\sqrt{\frac{2h}{p}}z_1 + \frac{2h}{\sqrt{p}}\bar{x}_1 - \frac{h}{p}\right)\right.$$
$$\left. + \frac{p}{2}\exp\left(-\sqrt{\frac{2h}{p}}z_1 - \frac{2h}{\sqrt{p}}\bar{x}_1 - \frac{h}{p}\right)\right).$$

The monotonicity of the logarithm together with the independence of $z_1$ and $\bar{x}_1$ give the lower bound

$$\mathbb{E}\log\left(\frac{p}{2}\exp\left(\sqrt{\frac{2h}{p}}z_1 + \frac{2h}{\sqrt{p}}\bar{x}_1 - \frac{h}{p}\right)\right)\mathbf{1}_{\{\bar{x}_1=1/\sqrt{p}\}} = \frac{p}{2}\log\left(\frac{p}{2}\right) + \frac{h}{2}.$$

Similarly,

$$\mathbb{E}\log\left(\frac{p}{2}\exp\left(-\sqrt{\frac{2h}{p}}z_1 - \frac{2h}{\sqrt{p}}\bar{x}_1 - \frac{h}{p}\right)\right)\mathbf{1}_{\{\bar{x}_1=-1/\sqrt{p}\}} = \frac{p}{2}\log\left(\frac{p}{2}\right) + \frac{h}{2}.$$

Combining these two bounds reveals that

$$\psi(h) \geqslant p\log\left(\frac{p}{2}\right) + h.$$

If we set $h = 1/2$ and choose $p^*$ small enough so $-\frac{1}{16} \leqslant p\log(p/2) \leqslant 0$ whenever $0 \leqslant p \leqslant p^*$, we have that

$$\psi\left(\frac{1}{2}\right) - \frac{(1/2)^2}{4t} \geqslant \frac{1}{2} - \frac{1}{16} - \frac{1}{16t} = \frac{1}{16}\left(7 - \frac{1}{t}\right).$$

It follows that any maximizer $h_*(t)$ of (4.110) must be strictly positive for any $t > \frac{1}{7}$, and therefore $\mathsf{mmse}(t) < 1$. Setting $t^*(p) = \frac{1}{7}$ completes the proof. ∎

## 4.5 The community detection problem

In this section, we analyze the information-theoretic properties of the stochastic block model. We do so by leveraging a universality property of the free energy (4.9) that will allow us to map the stochastic block model onto the symmetric rank-one matrix estimation problem. The stochastic block model is a simple model for networks with a community structure. It was first introduced in the machine learning and statistics literature [116, 139, 261, 264], and also emerged independently in a variety of other scientific disciplines. In the theoretical computer science community it is often called the planted partition model [58, 65, 106], while the mathematics literature also refers to it as the inhomogeneous random graph model [54]. This model has been used as a test bed for clustering and community detection algorithms used in a variety of contexts including social networks [201], protein-to-protein interaction networks [76], recommendation systems [164], medical prognosis [248], DNA folding [66], image segmentation [238], and natural language processing [30]. In this section we focus on the dense stochastic block model with two communities, which we now define precisely.

Consider $N$ individuals, each belonging to exactly one of two communities. We encode the individuals as elements of $\{1,\ldots,N\}$, and represent the community structure using a vector

$$\overline{\sigma} := (\overline{\sigma}_1,\ldots,\overline{\sigma}_N) \in \Sigma_N := \{-1,+1\}^N. \tag{4.129}$$

For convenience, we have changed our notation of the hypercube $\{-1,+1\}^N$ from $\{\pm 1\}^N$ to $\Sigma_N$; we will keep this new notation when we discuss spin glasses in Chapter 6. We understand that individuals $i$ and $j$ belong to the same community if and only if $\overline{\sigma}_i = \overline{\sigma}_j$. We sample the labels $(\overline{\sigma}_i)_{i \leqslant N}$ independently from a Bernoulli distribution $P_1$ with probability of success $p \in (0,1)$ and expectation $\overline{m}$,

$$p := P_1\{1\} = \mathbb{P}\{\overline{\sigma}_i = 1\} \quad \text{and} \quad \overline{m} := \mathbb{E}\overline{\sigma}_1 = 2p - 1. \tag{4.130}$$

The assignment vector $\overline{\sigma}$ is thus distributed according to the product law

$$P_N := P_1^{\otimes N}, \tag{4.131}$$

and the expected sizes of the communities are $Np$ and $N(1-p)$. Using the assignment vector $\overline{\sigma}$, we construct a random undirected graph $\mathbf{G}_N = (G_{ij})_{i,j \leqslant N}$ with vertex set $\{1,\ldots,N\}$ by stipulating that an edge between node $i$ and node $j$ is present with conditional probability

$$\mathbb{P}\{G_{ij} = 1 \mid \overline{\sigma}\} := \begin{cases} a_N & \text{if } \overline{\sigma}_i = \overline{\sigma}_j, \\ b_N & \text{if } \overline{\sigma}_i \neq \overline{\sigma}_j, \end{cases} \tag{4.132}$$

for some $a_N, b_N \in (0,1)$, independently of all other edges. In other words, the probability that an edge is present between node $i$ and node $j$ depends only on whether or not the individuals $i$ and $j$ belong to the same community. To express (4.132) more succinctly, it is convenient to introduce the average and the gap of $a_N$ and $b_N$,

$$c_N := \frac{a_N + b_N}{2} \quad \text{and} \quad \Delta_N := \frac{a_N - b_N}{2} \in (-c_N, c_N), \tag{4.133}$$

in such a way that

$$\mathbb{P}\{G_{ij} = 1 \mid \overline{\sigma}\} = c_N + \Delta_N \overline{\sigma}_i \overline{\sigma}_j. \tag{4.134}$$

Our inference task is to try to reconstruct the community structure $\overline{\sigma}$ as best we can, given the observation of the network of interactions $\mathbf{G}_N = (G_{ij})_{i,j \leqslant N}$. This is illustrated in Figure 4.3.

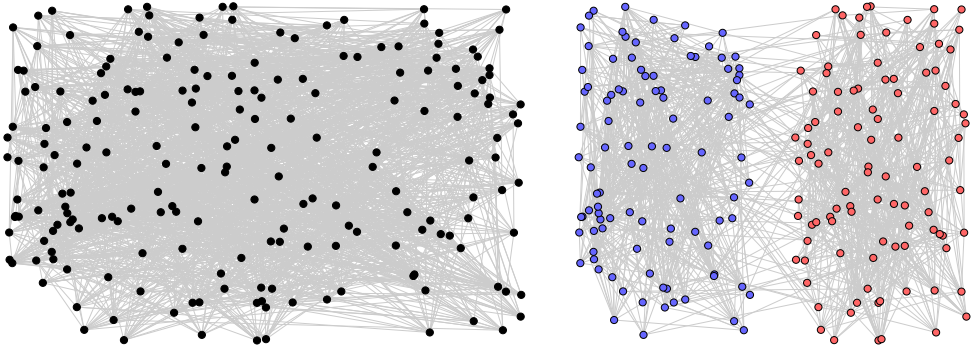

**Figure 4.3** The two figures display the same graph of interactions with 100 red nodes and 100 blue nodes, and with each node having on average 10 links with nodes of the same color and 1 link with a node of a different color (all edges are independent). On the left figure, the blue and the red nodes have been placed uniformly at random over the entire area, and the colors have been concealed. On the right figure, the blue and the red nodes have been classified by placing them randomly to the left and to the right of the area respectively. In the community detection problem, the left figure is shown to the statistician, whose goal is to infer the colouring of the nodes.

When $p = 1/2$, the symmetry between the two communities makes it clear that we can only hope to recover $\overline{\sigma}$ up to a change of sign. We may also consider the problem of recovering the matrix $\overline{\sigma}\overline{\sigma}^*$, which encodes whether or not two individuals belong to the same community or not for each pair of individuals.

In the case when $\Delta_N \leqslant 0$, it is more likely for an edge to be present between nodes in different communities, and the model is called *disassortative*. When $\Delta_N > 0$, connections are more likely between individuals in the same community, and the model is termed *assortative*.

The recovery task in the stochastic block model can be interpreted in at least two different ways, often called exact recovery and detection. The exact recovery task aims to determine the regimes of $a_N$ and $b_N$, or equivalently of $c_N$ and $\Delta_N$, for which there exists an algorithm that completely recovers the two communities with high probability, up to a global change of sign. A necessary condition for exact recovery is that the random graph $\mathbf{G}_N$ be connected; this makes exact recovery impossible in the sparse regime when the average degree of a node remains bounded as $N$ diverges. The sharp threshold for exact recovery was obtained in [2, 190], where it was shown that in the symmetric dense regime, $p = 1/2$, $a_N = a \log(N)/N$ and $b_N = b \log(N)/N$, exact recovery is possible, and efficiently so, if and only if $\sqrt{a} - \sqrt{b} \geq 2$. On the other hand, the detection task is to construct a partition of the graph $\mathbf{G}_N$ that is positively correlated with the assignment vector $\overline{\sigma}$ with high probability, possibly up to a global change of sign. This is the problem that we will focus on, in the dense regime where the average degree of a node diverges with $N$. More specifically, we will work under the two following assumptions on $c_N$ and the quantity

$$\lambda_N := \frac{N\Delta_N^2}{c_N(1-c_N)}. \tag{4.135}$$

**A1** The sequence $(\lambda_N)_{N \geq 1}$ converges to some value $\lambda \geq 0$.

**A2** The sequence $(Nc_N(1-c_N))_{N \geq 1}$ diverges to infinity.

The second of these assumptions implies that the average degree of a node $i$,

$$\mathbb{E}\deg(i) = (N-1)(c_N + \overline{m}\Delta_N), \tag{4.136}$$

diverges with $N$. The sparse regime, when the average degree of a node remains bounded even in the limit of large system size, is more challenging. In this case, identifying the regimes of $c_N$ and $\Delta_N$ for which the detection of the communities is possible has been achieved in [171, 188, 191] (in a setting with two communities). Concerning the precise identification of the limit mutual information between the community structure and the random graph $\mathbf{G}_N$, we refer the reader interested in the sparse disassortative regime to [4, 80, 91], and the reader interested in the harder sparse assortative regime to [3, 91, 105, 129, 149, 189, 192, 267].

As in the symmetric rank-one matrix estimation problem, we may be interested in studying the large-$N$ limit of the minimal mean-square error for the recovery of the matrix $\overline{\sigma}\,\overline{\sigma}^*$, that is,

$$\mathsf{mmse}_N := \frac{1}{N}\mathbb{E}|\overline{\sigma}\,\overline{\sigma}^* - \mathbb{E}[\overline{\sigma}\,\overline{\sigma}^* \mid \mathbf{G}_N]|^2, \tag{4.137}$$

or similarly for the vector $\overline{\sigma}$ itself. We will here focus on studying the asymptotic behaviour of the mutual information

$$\mathsf{I}_N := \frac{1}{N}\mathbb{E}\int_{\Sigma_N} \log\left(\frac{dP_{\overline{\sigma}|\mathbf{G}_N}}{dP_N}(\sigma)\right) dP_{\overline{\sigma}|\mathbf{G}_N}(\sigma) \tag{4.138}$$

between the graph $\mathbf{G}_N$ and the vector $\overline{\sigma}$. The asymptotic behaviour of the minimal mean-square error in (4.137) can also be identified using similar ideas, in analogy with the results of Propositions 4.3 and 4.4 and (4.21); see also [95, 161]. Observing that

$$\mathbb{P}\{\mathbf{G}_N = (G_{ij})|\overline{\sigma} = \sigma\} = \prod_{i<j}(c_N + \Delta_N \overline{\sigma}_i \overline{\sigma}_j)^{G_{ij}}(1 - c_N - \Delta_N \overline{\sigma}_i \overline{\sigma}_j)^{1-G_{ij}}, \quad (4.139)$$

we can use Bayes' formula to obtain the law of the signal $\overline{\sigma}$ conditionally on the observation of $\mathbf{G}_N$. We write it in the form of a Gibbs measure,

$$\mathbb{P}\{\overline{\sigma} = \sigma \mid \mathbf{G}_N = (G_{ij})\} = \frac{\exp H_N^{\mathrm{SBM}}(\sigma) P_N(\sigma)}{\int_{\Sigma_N} \exp H_N^{\mathrm{SBM}}(\tau) \, dP_N(\tau)}, \quad (4.140)$$

for the Hamiltonian

$$H_N^{\mathrm{SBM}}(\sigma) := \sum_{i<j} \log\left[\left(1 + \frac{\Delta_N}{c_N}\sigma_i\sigma_j\right)^{G_{ij}}\left(1 - \frac{\Delta_N}{1-c_N}\sigma_i\sigma_j\right)^{1-G_{ij}}\right]. \quad (4.141)$$

We denote the associated free energy by

$$\overline{F}_N^{\mathrm{SBM}} := \frac{1}{N} \mathbb{E} \log \int_{\Sigma_N} \exp H_N^{\mathrm{SBM}}(\sigma) \, dP_N(\sigma). \quad (4.142)$$

In the limit of large $N$, this free energy coincides with the mutual information (4.138), up to an additive constant.

**Proposition 4.23.** *Under assumptions* (**A1**) *and* (**A2**), *the limits of the free energy* (4.142) *and of the mutual information* (4.138) *differ by an additive constant,*

$$I_N = \frac{\lambda}{4} - \overline{F}_N^{\mathrm{SBM}} + o(1). \quad (4.143)$$

*Proof.* The explicit form of the likelihood in (4.139) and the definition of the Hamiltonian in (4.141) imply that

$$I_N = \frac{1}{N}\mathbb{E}H_N^{\mathrm{SBM}}(\overline{\sigma}) - \overline{F}_N^{\mathrm{SBM}}. \quad (4.144)$$

Since the coordinates of the assignment vector $\overline{\sigma}$ are i.i.d. the first term simplifies to

$$\frac{1}{N}\mathbb{E}H_N^{\mathrm{SBM}}(\overline{\sigma}) = \frac{1}{N} \cdot \binom{N}{2} \cdot \left[\mathbb{E}G_{12}\log\left(1 + \frac{\Delta_N}{c_N}\overline{\sigma}_1\overline{\sigma}_2\right)\right.$$
$$\left. + \mathbb{E}(1-G_{12})\log\left(1 - \frac{\Delta_N}{1-c_N}\overline{\sigma}_1\overline{\sigma}_2\right)\right]. \quad (4.145)$$

Averaging with respect to the randomness of $G_{12}$ conditionally on the randomness of the assignment vector $\overline{\sigma}$ reveals that this is equal to

$$\frac{N-1}{2}\left[c_N \mathbb{E}\left(1+\frac{\Delta_N}{c_N}\overline{\sigma}_1\overline{\sigma}_2\right)\log\left(1+\frac{\Delta_N}{c_N}\overline{\sigma}_1\overline{\sigma}_2\right)\right.$$
$$\left.+(1-c_N)\mathbb{E}\left(1-\frac{\Delta_N}{1-c_N}\overline{\sigma}_1\overline{\sigma}_2\right)\log\left(1-\frac{\Delta_N}{1-c_N}\overline{\sigma}_1\overline{\sigma}_2\right)\right], \quad (4.146)$$

and Taylor expanding the logarithm gives

$$\frac{1}{N}\mathbb{E}H_N^{\text{SBM}}(\overline{\sigma}) = \frac{N-1}{4}\cdot\left[\frac{\Delta_N^2}{c_N(1-c_N)}+\mathcal{O}\left(\frac{\Delta_N^3}{c_N^2}+\frac{\Delta_N^3}{(1-c_N)^2}\right)\right]. \quad (4.147)$$

Recalling the definition of the constant $\lambda_N$ in (4.135), observing that $c_N \in (0,1)$ and remembering (4.144) reveals that

$$\mathsf{I}_N = \frac{N-1}{N}\cdot\frac{\lambda_N}{4}-\overline{F}_N^{\text{SBM}}+\mathcal{O}\left(\frac{N-1}{N}\cdot\frac{\lambda_N^{3/2}}{\sqrt{Nc_N(1-c_N)}}\right). \quad (4.148)$$

Invoking assumptions (**A1**) and (**A2**) completes the proof. ∎

To determine the limit of the free energy $\overline{F}_N^{\text{SBM}}$, we will follow [95, 161] and argue that it is asymptotically equivalent to the free energy of a rank-one estimation problem of just the form we have been studying in earlier sections of this chapter. To see this, we start by simplifying the expression of the Hamiltonian $H_N^{\text{SBM}}$ and discarding some lower-order terms. Taylor expanding the logarithm in the definition of the Hamiltonian (4.141) reveals that

$$H_N^{\text{SBM}}(\sigma) = \widetilde{H}_N^{\text{SBM}}(\sigma) + \sum_{i<j}\left[G_{ij}\left(\frac{\Delta_N^2}{2(1-c_N)^2}-\frac{\Delta_N^2}{2c_N^2}\right)-\frac{\Delta_N^2}{2(1-c_N)^2}\right]$$
$$+\sum_{i<j}\mathcal{O}\left(G_{ij}\frac{\Delta_N^3}{c_N^3}+(1-G_{ij})\frac{\Delta_N^3}{(1-c_N)^3}\right) \quad (4.149)$$

for the Hamiltonian

$$\widetilde{H}_N^{\text{SBM}}(\sigma) := \sum_{i<j}\frac{\Delta_N}{c_N(1-c_N)}(G_{ij}-c_N)\sigma_i\sigma_j. \quad (4.150)$$

If we introduce the free energy

$$\widetilde{F}_N^{\text{SBM}} := \frac{1}{N}\mathbb{E}\int_{\Sigma_N}\exp\widetilde{H}_N^{\text{SBM}}(\sigma)\,dP_N(\sigma) \quad (4.151)$$

associated with the Hamiltonian (4.150), then the equality (4.149) implies that

$$\overline{F}_N^{\text{SBM}} = \widetilde{F}_N^{\text{SBM}} + \frac{1}{N}\sum_{i<j}\mathbb{E}\left[G_{ij}\left(\frac{\Delta_N^2}{2(1-c_N)^2} - \frac{\Delta_N^2}{2c_N^2}\right) - \frac{\Delta_N^2}{2(1-c_N)^2}\right]$$
$$+ \sum_{i<j}\mathcal{O}\left(\mathbb{E}G_{ij}\frac{\Delta_N^3}{c_N^3} + (1-\mathbb{E}G_{ij})\frac{\Delta_N^3}{(1-c_N)^3}\right). \quad (4.152)$$

Noticing that $\mathbb{E}G_{ij} = c_N + \Delta_N \overline{m}^2$ and remembering the definition of $\lambda_N$ in (4.135), this becomes

$$\overline{F}_N^{\text{SBM}} = \widetilde{F}_N^{\text{SBM}} - \frac{N-1}{N}\cdot\frac{\lambda_N}{4} + \mathcal{O}\left(\frac{N-1}{N}\cdot\frac{\lambda_N^{3/2}}{\sqrt{Nc_N(1-c_N)}}\right). \quad (4.153)$$

Leveraging the assumptions (**A1**) and (**A2**) shows that

$$\overline{F}_N^{\text{SBM}} = \widetilde{F}_N^{\text{SBM}} - \frac{\lambda}{4} + o(1) \quad (4.154)$$

so the free energies (4.142) and (4.151) are equal up to an additive constant. The free energy (4.151) starts to look more like the free energy (4.9) in the symmetric rank-one matrix estimation problem. Indeed, if we introduce the centred random variables

$$\widetilde{G}_{ij} := \frac{\Delta_N}{c_N(1-c_N)}\left(G_{ij} - c_N - \Delta_N\overline{\sigma}_i\overline{\sigma}_j\right), \quad (4.155)$$

then the Hamiltonian (4.150) may be written as

$$\widetilde{H}_N^{\text{SBM}}(\sigma) = \sum_{i<j}\left(\widetilde{G}_{ij}\sigma_i\sigma_j + \frac{\lambda_N}{N}\overline{\sigma}_i\overline{\sigma}_j\sigma_i\sigma_j\right). \quad (4.156)$$

A direct computation shows that the random variables $\widetilde{G}_{ij}$ have variance

$$\mathbb{E}\widetilde{G}_{ij}^2 = \frac{\lambda_N}{Nc_N(1-c_N)}\mathbb{E}\left(1 - c_N - \Delta_N\overline{\sigma}_i\overline{\sigma}_j\right)\left(c_N + \Delta_N\overline{\sigma}_i\overline{\sigma}_j\right) \quad (4.157)$$

$$= \frac{\lambda_N}{N} + \mathcal{O}\left(\frac{\lambda_N}{N}\cdot\left(\sqrt{\frac{\lambda_N}{Nc_N(1-c_N)}} + \frac{\lambda_N}{N}\right)\right). \quad (4.158)$$

Recall that under our assumption (**A2**), the average degree of a node diverges as $N$ tends to infinity. This gives some credence to the idea that in the limit of large system size, a sort of central limit theorem takes place and we may as well substitute the random variables $(\widetilde{G}_{ij})_{ij}$ by centred Gaussian random variables with the same variance. In other words, we may expect the free energy (4.151) to be asymptotically equivalent to the Gaussian free energy

$$\overline{F}_N^{\text{gauss}}(\lambda) := \frac{1}{N}\mathbb{E}\log\int_{\Sigma_N}\exp H_N^{\text{gauss}}(\lambda,\sigma)\,\mathrm{d}P_N(\sigma) \quad (4.159)$$

associated with the Hamiltonian

$$H_N^{\text{gauss}}(\lambda,\sigma) := \sum_{i<j}\left(\sqrt{\frac{\lambda}{N}}W_{ij}\sigma_i\sigma_j + \frac{\lambda}{N}\overline{\sigma}_i\overline{\sigma}_j\sigma_i\sigma_j\right).$$

As shown in Exercise 4.11, up to an additive constant, this Gaussian free energy is the free energy (4.9) in the symmetric rank-one matrix estimation problem,

$$\overline{F}_N^{\text{gauss}}(\lambda) = \overline{F}_N^\circ\left(\frac{\lambda}{4}\right) + \frac{\lambda}{4} - \frac{\lambda}{2N}. \tag{4.160}$$

To show rigorously that the free energies (4.151) and (4.159) are asymptotically equivalent up to an additive constant, and therefore that the free energy (4.142) in the stochastic block model is asymptotically equivalent to the free energy (4.9) in the symmetric rank-one matrix estimation problem for an appropriate choice of prior and signal-to-noise ratio $t$, we use an interpolation argument.

**Theorem 4.24.** *Under assumptions* (**A1**) *and* (**A2**), *the free energy* (4.142) *in the stochastic block model is asymptotically equivalent to the free energy* (4.9) *in the symmetric rank-one matrix estimation problem with Bernoulli prior,*

$$\overline{F}_N^{\text{SBM}} = \overline{F}_N^\circ\left(\frac{\lambda}{4}\right) + o(1). \tag{4.161}$$

*Proof.* By the asymptotic equivalences (4.154) and (4.160), it suffices to show that

$$\widetilde{F}_N^{\text{SBM}} = \overline{F}_N^{\text{gauss}}(\lambda) + o(1). \tag{4.162}$$

To alleviate notation, we will keep the dependence on $\lambda$ implicit. We proceed by interpolation. For each $t \in [0,1]$, we define the interpolating Hamiltonian

$$H_{N,t}(\sigma) := \sum_{i<j}\left[\left(\sqrt{t}\widetilde{G}_{ij} + \sqrt{1-t}\sqrt{\frac{\lambda}{N}}W_{ij}\right)\sigma_i\sigma_j + \left(t\cdot\frac{\lambda_N}{N} + (1-t)\cdot\frac{\lambda}{N}\right)\overline{\sigma}_i\overline{\sigma}_j\sigma_i\sigma_j\right]$$

and the interpolating free energy

$$\overline{F}_N(t) := \frac{1}{N}\mathbb{E}\log\int_{\Sigma_N}\exp H_{N,t}(\sigma)\,\mathrm{d}P_N(\sigma)$$

in such a way that by the fundamental theorem of calculus

$$\left|\widetilde{F}_N^{\text{SBM}} - \overline{F}_N^{\text{gauss}}\right| = \left|\overline{F}_N(1) - \overline{F}_N(0)\right| \leqslant \sup_{t\in[0,1]}\left|\overline{F}_N'(t)\right|. \tag{4.163}$$

To compute this derivative, write $\langle \cdot \rangle_t$ for the average with respect to the Gibbs measure associated with the interpolating Hamiltonian $H_{N,t}$, and observe that

$$\overline{F}'_N(t) = \frac{1}{2N\sqrt{t}} \sum_{i<j} \mathbb{E}\widetilde{G}_{ij}\langle \sigma_i\sigma_j\rangle_t - \frac{\sqrt{\lambda}}{2N^{3/2}\sqrt{1-t}} \sum_{i<j} \mathbb{E}W_{ij}\langle \sigma_i\sigma_j\rangle_t$$
$$+ (\lambda_N - \lambda) \cdot \frac{1}{N^2} \sum_{i<j} \mathbb{E}\langle \overline{\sigma}_i\overline{\sigma}_j\sigma_i\sigma_j\rangle_t. \qquad (4.164)$$

At this point, fix indices $i < j$ and introduce the function $F(\widetilde{G}_{ij}) = \langle \sigma_i\sigma_j\rangle_t$. A direct computation reveals that

$$\partial_{\widetilde{G}_{ij}} F = \langle \sigma_i\sigma_j \partial_{\widetilde{G}_{ij}} H_{N,t}(\sigma)\rangle_t - \langle \sigma_i\sigma_j\rangle_t \langle \partial_{\widetilde{G}_{ij}} H_{N,t}(\sigma)\rangle_t = \sqrt{t}\bigl(1 - \langle \sigma_i\sigma_j\rangle_t^2\bigr),$$
$$\partial^2_{\widetilde{G}_{ij}} F = -2\sqrt{t}\langle \sigma_i\sigma_j\rangle_t \cdot \partial_{\widetilde{G}_{ij}} F = -2t\bigl(\langle \sigma_i\sigma_j\rangle_t - \langle \sigma_i\sigma_j\rangle_t^3\bigr),$$

so the approximate Gaussian integration by parts formula in Exercise 4.6 gives

$$\mathbb{E}\widetilde{G}_{ij}\langle \sigma_i\sigma_j\rangle_t = \sqrt{t}\mathbb{E}\widetilde{G}_{ij}^2\bigl(1 - \mathbb{E}\langle \sigma_i\sigma_j\rangle_t^2\bigr) + \mathcal{O}\bigl(\mathbb{E}|\widetilde{G}_{ij}|^3\bigr).$$

On the other hand, the Gaussian integration by parts formula as in (4.29) shows that

$$\mathbb{E}W_{ij}\langle \sigma_i\sigma_j\rangle_t = \sqrt{1-t} \cdot \sqrt{\frac{\lambda}{N}} \cdot \bigl(1 - \mathbb{E}\langle \sigma_i\sigma_j\rangle_t^2\bigr).$$

It follows by (4.164) that

$$\overline{F}'_N(t) = \frac{1}{2N} \sum_{i<j} \mathbb{E}\widetilde{G}_{ij}^2\bigl(1 - \mathbb{E}\langle \sigma_i\sigma_j\rangle_t^2\bigr) - \frac{\lambda}{2N^2} \sum_{i<j} \bigl(1 - \mathbb{E}\langle \sigma_i\sigma_j\rangle_t^2\bigr)$$
$$+ \mathcal{O}\bigl(N\mathbb{E}|\widetilde{G}_{12}|^3 + |\lambda_N - \lambda|\bigr).$$

Together with (4.158) and (4.163), this implies that

$$\widetilde{F}_N^{\mathrm{SBM}} = \overline{F}_N^{\mathrm{gauss}} + \mathcal{O}\left(\lambda_N \cdot \left(\sqrt{\frac{\lambda_N}{Nc_N(1-c_N)}} + \frac{\lambda_N}{N}\right) + N\mathbb{E}|\widetilde{G}_{12}|^3 + |\lambda_N - \lambda|\right).$$

Observing that

$$\mathbb{E}|\widetilde{G}_{12}|^3 \leq \left|\frac{\Delta_N^3}{c_N^3(1-c_N)^3}\right| \mathbb{E}\bigl|(1 - c_N - \Delta_N\overline{\sigma}_i\overline{\sigma}_j)(c_N + \Delta_N\overline{\sigma}_i\overline{\sigma}_j)\bigr|$$
$$\leq \frac{4\lambda_N^{3/2}}{N\sqrt{Nc_N(1-c_N)}},$$

and remembering the assumptions (**A1**) and (**A2**) establishes (4.162) and completes the proof. ∎

Together with Theorem 4.9 on the limit free energy in the symmetric rank-one matrix estimation problem and Proposition 4.23 on the relationship between the mutual information and the free energy in the stochastic block model, this allows us to deduce the following variational representation for the limit free energy and limit mutual information in the stochastic block model.

**Theorem 4.25.** *Under assumptions* (**A1**) *and* (**A2**), *the limit of the free energy in the stochastic block model admits the variational representation*

$$\lim_{N\to+\infty} \overline{F}_N^{\text{SBM}} = \sup_{h\geqslant 0}\left(\psi(h) - \frac{h^2}{\lambda}\right),$$

*where*

$$\psi(h) := \mathbb{E}\log \int_{\Sigma_1} \exp\left(\sqrt{2h}\sigma z_1 + 2h\sigma\overline{\sigma}_1 - h\right) dP_1(\sigma), \tag{4.165}$$

*with $z_1$ a standard Gaussian random variable independent of any other randomness. The limit mutual information in the stochastic block model admits the variational representation*

$$\lim_{N\to+\infty} I_N = \frac{\lambda}{4} - \sup_{h\geqslant 0}\left(\psi(h) - \frac{h^2}{\lambda}\right). \tag{4.166}$$

**Exercise 4.11.** Show that the Gaussian free energy (4.159) and the free energy (4.9) are related by

$$\overline{F}_N^{\text{gauss}}(\lambda) = \overline{F}_N^{\circ}\left(\frac{\lambda}{4}\right) + \frac{\lambda}{4} - \frac{\lambda}{2N}. \tag{4.167}$$

**Exercise 4.12.** In the symmetric dense stochastic block model with $p = \frac{1}{2}$, show that

$$\lim_{N\to+\infty} I_N = \inf_{h\geqslant 0}\left(\frac{\lambda}{4} + \frac{h^2}{\lambda} + h - \mathbb{E}\log\cosh\left(\sqrt{2h}z_1 + 2h\overline{\sigma}_1\right)\right).$$

# Chapter 5
# Poisson point processes and extreme values

The main goal of this chapter is to lay important groundwork for our study of mean-field spin glass models in Chapter 6. The material we wish to present involves the consideration of Poisson point processes. In the first four sections of this chapter, we therefore give a self-contained introduction to this topic. This part culminates with limit theorems for extreme values of independent and identically distributed random variables, which in our view should be of interest to a very broad audience. We then turn to aspects that are specifically geared towards the analysis of mean-field spin glasses in the last three sections.

In Section 5.1, we study the topological properties of the space of point measures on a locally compact and separable metric space. We then define point processes, which are random variables on the space of point measures, in Section 5.2. In Section 5.3, we introduce Poisson point processes and study their fundamental properties. We use these objects to determine the asymptotic behaviour of extremes of i.i.d. random variables in Section 5.4. We then proceed in Section 5.5 with a thorough analysis of a special Poisson point process known as the Poisson-Dirichlet point process. This process plays a fundamental role in the analysis of mean-field spin glasses, as a building block in the hierarchical construction of the Poisson-Dirichlet cascades discussed in Section 5.6. Finally, in Section 5.7, we discuss a characterization of Poisson-Dirichlet cascades in terms of certain distributional identities, and the fundamental role played by the notion of ultrametricity.

## 5.1 The space of point measures

A Poisson point process is a special type of point process, which in turn is a random variable taking values in a space of point measures. In this section, we study the basic topological properties of the space of point measures on a locally compact and separable metric space $(S, d)$, endowed with its Borel $\sigma$-algebra $\mathcal{B}(S)$. In practice, we will mostly be interested in the case when $S$ is $\mathbb{R}^d$ or an open subset of $\mathbb{R}^d$ equipped with the Euclidean distance, but this level of generality does not cause much additional difficulty, and is at times handy.

A measure $\lambda$ on $S$ is said to be *locally finite* if for every compact set $K \subseteq S$, we have $\lambda(K) < +\infty$. A *point measure* on $S$ is a measure $\lambda$ on $S$ that is locally finite and for which there exists a collection of points $(x_i)_{i \in I}$ with $x_i \in S$ satisfying

$$\lambda = \sum_{i \in I} \delta_{x_i}, \tag{5.1}$$

where $\delta_{x_i}$ denotes the Dirac measure at $x_i$. The space of point measures on $S$ is denoted by $\mathcal{M}_\delta(S)$. Notice that a point $x_i$ in the representation (5.1) could be repeated. In particular, measures of the form $k\delta_x$ for some $k \in \mathbb{N}$ and $x \in S$ belong to $\mathcal{M}_\delta(S)$. By Exercise 5.1, the index set $I$ is always countable. We endow $\mathcal{M}_\delta(S)$ with the topology of vague convergence. This is the coarsest topology that makes the mappings

$$\lambda \mapsto \int_S f \, d\lambda \tag{5.2}$$

continuous for every $f \in C_c(S; \mathbb{R})$. Here $C_c(S; \mathbb{R})$ is the space of real-valued continuous functions with compact support in $S$. Specifically, this means that a sequence $(\lambda_n)_{n \geqslant 1} \subseteq \mathcal{M}_\delta(S)$ converges to a measure $\lambda \in \mathcal{M}_\delta(S)$ if and only if for every continuous function $f \in C_c(S; \mathbb{R})$ of compact support,

$$\lim_{n \to +\infty} \int_S f \, d\lambda_n = \int_S f \, d\lambda. \tag{5.3}$$

Notice that this form of convergence does not preserve the total mass of a point measure. For instance, the Dirac mass $\delta_n$, seen as an element of $\mathcal{M}_\delta(\mathbb{R})$, converges to the null measure as $n$ tends to infinity. We will mainly be concerned with random elements taking values in $\mathcal{M}_\delta(S)$ and their convergence. For the classical convergence theory discussed in Sections A.5 and A.6 of the appendix to apply, it will be useful to know that the topology of vague convergence on $\mathcal{M}_\delta(S)$ is metrizable. We will deduce this from the following result regarding the separability of the space of compactly supported continuous functions on $S$.

**Lemma 5.1.** *There exists a countable family $\mathcal{F} \subseteq C_c(S; \mathbb{R})$ with the property that, for any $f \in C_c(S; \mathbb{R})$ and $\varepsilon > 0$, there exists $g \in \mathcal{F}$ with*

$$\sup_{x \in S} |f(x) - g(x)| \leqslant \varepsilon \quad \text{and} \quad \operatorname{supp} g \subseteq \{x \in S \mid d(x, \operatorname{supp} f) \leqslant \varepsilon\}. \tag{5.4}$$

*In particular, the space $C_c(S; \mathbb{R})$ is separable.*

*Proof.* Let $\mathcal{D}$ be a countable and dense subset of $S$. We denote by $\mathcal{F}^+$ the set of functions in $C_c(S; \mathbb{R}_{\geqslant 0})$ that can be written in the form

$$x \mapsto \max_{y \in Y} \bigl(\lambda_y (1 - k d(x, y))_+\bigr), \tag{5.5}$$

where $Y$ is a finite subset of $\mathcal{D}$, $(\lambda_y)_{y \in Y}$ are non-negative rational numbers, $k$ is a positive integer, and we recall that $d$ is the distance function on $S$ and that $r_+ := \max(r, 0)$. Up to decomposing a function $f \in C_c(S; \mathbb{R})$ into its positive and negative parts, and setting $\mathcal{F}$ to be the set of all functions that can be written as the difference between two functions in $\mathcal{F}^+$, it suffices to fix a non-negative function $f \in C_c(S; \mathbb{R}_{\geq 0})$ as well as $\varepsilon > 0$ and to construct $g \in \mathcal{F}_+$ satisfying (5.4). Moreover, up to replacing $\varepsilon$ by $\varepsilon \|f\|_\infty$ in (5.4), we can assume without loss of generality that $f$ is bounded by one. Since $f$ is supported on a compact set $K$ and is continuous, it is uniformly continuous. Hence, there exists $\delta > 0$ such that for every $x, y \in S$, we have $|f(x) - f(y)| \leq \varepsilon$ whenever $d(x, y) \leq \delta$. Up to reducing $\delta$ if necessary, we can also assume that $\delta \leq \varepsilon$. We let $k$ be the smallest integer larger than $\delta^{-1}$. For every $y \in \mathcal{D}$, we pick $\lambda_y$ a rational number in the interval $[f(y) - \varepsilon, f(y) + \varepsilon]$. For every $\eta \in (0, \delta)$, we denote by $Y_\eta$ a finite subset of $\mathcal{D}$, made of points that are at distance at most $\eta$ from $K$, such that every $x \in K$ is at distance at most $\eta$ from a point in $Y_\eta$. The set $Y_\eta$ can be constructed by extracting a finite sub-cover from the open cover of the compact set $K$ formed by open balls of radius $\eta$ centred at the points in $\mathcal{D}$. We write $g_\eta$ for the function (5.5) with these choices of parameters, and with $Y$ replaced by $Y_\eta$. The function $g_\eta$ is Lipschitz continuous, with Lipschitz constant at most $k + k\varepsilon$. For every $x \in K$, we can find $y \in Y_\eta$ such that $d(x, y) \leq \eta \leq \delta$, and thus

$$|f(x) - g_\eta(x)| \leq |f(y) - g_\eta(y)| + \varepsilon + k\eta + k\varepsilon\eta \leq 2\varepsilon + k\eta + k\varepsilon\eta.$$

On the other hand, for $x \notin K$ but in the support of $g_\eta$, there is $y \in Y_\eta$ with $d(x, y) \leq k^{-1} \leq \delta$ and

$$g_\eta(x) = \lambda_y (1 - kd(x, y))_+.$$

It follows that

$$|f(x) - g_\eta(x)| = |g_\eta(x)| \leq |\lambda_y| \leq |f(y)| + \varepsilon \leq |f(x)| + 2\varepsilon = 2\varepsilon.$$

Choosing $\eta$ to be smaller than the minimum between $\varepsilon/k$ and $1/k$, we have shown that for every $x \in S$,

$$|f(x) - g_\eta(x)| \leq 4\varepsilon.$$

Notice also that every point in the support of $g_\eta$ is at distance at most $\eta + k^{-1} \leq \eta + \delta \leq 2\delta \leq 2\varepsilon$ from a point in $K$. We have therefore constructed a function satisfying (5.5), up to replacing $\varepsilon$ by $4\varepsilon$. In order to conclude, it remains to verify that the support of the function $g_\eta$ is compact. In order to do so, it suffices to check that, for $\varepsilon > 0$ sufficiently small, the set

$$K_\varepsilon := \{x \in S \mid \text{there exists } y \in K \text{ with } d(x, y) \leq \varepsilon\} \tag{5.6}$$

is compact. Since $S$ is locally compact, we can cover $K$ with open balls having compact closure, simply by picking an open ball having compact closure for each point in $K$. Since $K$ is compact, we can extract a finite sub-cover of $K$. As soon as $\varepsilon > 0$ is smaller than the radius of the smallest ball in this finite sub-cover, we see that the set $K_\varepsilon$ is compact. This completes the proof. ∎

This result suggests that vague convergence on $\mathcal{M}_\delta(S)$ should be determined by the set $\mathcal{F}$. If we denote by $(f_n)_{n \geq 1}$ an enumeration of the countable set $\mathcal{F}$, this would ensure that the metric $\rho : \mathcal{M}_\delta(S) \times \mathcal{M}_\delta(S) \to \mathbb{R}_{\geq 0}$ defined by

$$\rho(\lambda, \lambda') := \sum_{n=0}^{+\infty} 2^{-n} \frac{\left| \int_S f_n \, d\lambda - \int_S f_n \, d\lambda' \right|}{1 + \left| \int_S f_n \, d\lambda - \int_S f_n \, d\lambda' \right|} \tag{5.7}$$

metrizes the vague topology on $\mathcal{M}_\delta(S)$. The fact that (5.7) defines a metric is the content of Exercise 5.4.

**Theorem 5.2.** *A sequence $(\lambda_n)_{n \geq 1} \subseteq \mathcal{M}_\delta(S)$ converges vaguely to a point measure $\lambda \in \mathcal{M}_\delta(S)$ if and only if (5.3) holds for all functions $f \in \mathcal{F}$. In particular, the topology induced by the metric in (5.7) is the vague topology on $\mathcal{M}_\delta(S)$.*

*Proof.* The direct implication is immediate. Conversely, we assume that (5.3) holds for every $f \in \mathcal{F}$, and fix $f \in C_c(S; \mathbb{R})$ as well as $\varepsilon \in (0, 1)$. Letting $K$ be the support of $f$, we invoke Lemma 5.1 to find $g \in \mathcal{F}$ with

$$\sup_{x \in S} |f(x) - g(x)| \leq \varepsilon \quad \text{and} \quad \operatorname{supp} g \subseteq K_\varepsilon.$$

Here $K_\varepsilon$ denotes the $\varepsilon$-neighbourhood of $K$ defined in (5.6). By the triangle inequality,

$$\left| \int_S f \, d\lambda - \int_S f \, d\lambda_n \right| \leq \varepsilon \big( \lambda(K_\varepsilon) + \sup_{n \geq 1} \lambda_n(K_\varepsilon) \big) + \left| \int_S g \, d\lambda_n - \int_S g \, d\lambda \right|,$$

so letting $n$ tend to infinity and leveraging the assumption that (5.3) holds for $g \in \mathcal{F}$ reveals that

$$\limsup_{n \to +\infty} \left| \int_S f \, d\lambda - \int_S f \, d\lambda_n \right| \leq \varepsilon \big( \lambda(K_\varepsilon) + \sup_{n \geq 1} \lambda_n(K_\varepsilon) \big).$$

To bound this further, recall from the end of the proof of Lemma 5.1 that, decreasing $\varepsilon$ if necessary, the set $K_\varepsilon$ is compact. Decreasing $\varepsilon$ further if necessary, and invoking Lemmas A.7 and 5.1, we can find $\phi \in \mathcal{F}$ with $\phi \geq 1$ on $K_\varepsilon$. Since $(\int_S \phi \, d\lambda_n)_{n \geq 1}$ converges by assumption, it is uniformly bounded, so

$$M := \sup_{n \geq 1} \lambda_n(K_\varepsilon) \leq \sup_{n \geq 1} \int_S \phi \, d\lambda_n < +\infty. \tag{5.8}$$

It follows that

$$\limsup_{n \to +\infty} \left| \int_S f \, d\lambda - \int_S f \, d\lambda_n \right| \leq \varepsilon \lambda(K_\varepsilon) + \varepsilon M.$$

Recalling that $K_\varepsilon$ is compact for $\varepsilon > 0$ sufficiently small, and that $\lambda(K_\varepsilon)$ is therefore finite, and observing that $\lambda(K_\varepsilon)$ is non-increasing in $\varepsilon$, we obtain the desired result. ∎

Theorem 5.2 allows us to appeal to the basic convergence results recalled in Sections A.5 and A.6 of the appendix for random elements taking values in $\mathcal{M}_\delta(S)$ endowed with its Borel $\sigma$-algebra $\mathcal{B}(\mathcal{M}_\delta(S))$. In Exercise 5.2, we give a characterization of the Borel $\sigma$-algebra $\mathcal{B}(\mathcal{M}_\delta(S))$ as the smallest $\sigma$-algebra on $\mathcal{M}_\delta(S)$ which makes the evaluation map $\lambda \mapsto \lambda(A)$ measurable for every measurable set $A \in \mathcal{B}(S)$. In the next section, we will want to characterize the convergence of random elements on $\mathcal{M}_\delta(S)$ through their Laplace transform. To do this, it will be convenient to know what relatively compact sets in $\mathcal{M}_\delta(S)$ look like. Notice that the notions of compactness and sequential compactness in $\mathcal{M}_\delta(S)$ are equivalent by Theorem 5.2.

**Lemma 5.3.** *A set $M \subseteq \mathcal{M}_\delta(S)$ is relatively compact if and only if for all $f \in C_c(S; \mathbb{R})$,*

$$\sup_{\lambda \in M} \left| \int_S f \, d\lambda \right| < +\infty. \tag{5.9}$$

*Proof.* For the direct implication, we take a relatively compact set $M \subseteq \mathcal{M}_\delta(S)$, and for each $f \in C_c(S; \mathbb{R})$ we define the continuous functional $T_f : \mathcal{M}_\delta(S) \to \mathbb{R}$ by

$$T_f(\lambda) := \int_S f \, d\lambda.$$

Since $\overline{M}$ is compact, the set $T_f(\overline{M}) \subseteq \mathbb{R}$ is compact, and therefore bounded. This means that

$$\sup_{\lambda \in M} \left| \int_S f \, d\lambda \right| \leq \sup_{\lambda \in \overline{M}} \left| \int_S f \, d\lambda \right| < +\infty,$$

as desired.

For the converse implication, we assume that (5.9) holds for all $f \in C_c(S; \mathbb{R})$, and we fix a sequence $(\lambda_n)_{n \geq 1} \subseteq M$. We aim to identify an element $\lambda \in \mathcal{M}_\delta(S)$ and a subsequence of $(\lambda_n)_{n \geq 1}$ that converges vaguely to $\lambda$. We decompose the rest of the proof into three steps.

*Step 1: convergence on a given compact subset.* In this step, we show that for any given compact set $K \subseteq S$, one can find $\lambda \in \mathcal{M}_\delta(S)$ and a subsequence $(n(k))_{k \geq 1}$ such that for every $f \in C_c(S; \mathbb{R})$ with compact support in $K$,

$$\lim_{k \to +\infty} \int_S f \, d\lambda_{n(k)} = \int_S f \, d\lambda. \tag{5.10}$$

We denote by $\lambda_n|_K$ the restriction of $\lambda_n$ to $K$. Arguing as in the proof of Theorem 5.2, we see that the sequence $(\lambda_n(K))_{n \geq 1}$ must remain bounded. Up to the extraction of a subsequence, we can thus assume that it converges to some $j \geq 0$, and since the sequence $(\lambda_n(K))_{n \geq 1}$ only takes integer values, the subsequence is ultimately constant. That is, for $n$ large enough and along the subsequence, we can write

$$\lambda_n|_K = \sum_{i=1}^{j} \delta_{x_{n,i}},$$

for some $x_{n,1}, \ldots, x_{n,j} \in K$. Up to the extraction of a further subsequence, we can also assume that each of these points converges as $n$ tends to infinity, and we denote by $x_1, \ldots, x_j \in K$ their respective limits. We define $\lambda := \sum_{i=1}^{j} \delta_{x_i}$. Let $f \in C_c(S; \mathbb{R})$ be a continuous function with compact support in $K$. For every $n$ sufficiently large,

$$\int_S f \, d\lambda_n = \sum_{i=1}^{j} f(x_{n,i}),$$

so letting $n$ tend to infinity along the subsequence reveals that

$$\lim_{n \to +\infty} \int_S f \, d\lambda_n = \sum_{i=1}^{j} f(x_i) = \int_S f \, d\lambda.$$

This is (5.10).

*Step 2: covering $S$ with compact sets.* In order to patch together the different measures obtained in the previous step for different choices of the compact set $K$, we now build a countable collection $\mathcal{U}$ of open sets with compact closure which covers $S$. Since $S$ is assumed to be separable, there exists a countable set $\chi \subseteq S$ that is dense in $S$. We denote by $\mathcal{U}$ the collection of all open balls with centre a point in $\chi$, with radius a rational number, and whose closure is compact. The collection $\mathcal{U}$ is clearly countable; we claim that it covers $S$. To see this, let $x \in S$. Since $S$ is locally compact, for $r > 0$ sufficiently small, the closed ball $\overline{B}_r(x)$ is compact. Since $\chi$ is dense, one can find $y \in \chi$ whose distance to $x$ is less than $r/4$. Any open ball centred at $x$ and with radius some rational number in $(r/4, r/2)$ contains $x$, and it belongs to $\mathcal{U}$ since its closure still lies in the compact ball $\overline{B}_r(x)$.

*Step 3: constructing $\lambda$.* By the two previous steps and a diagonal extraction argument, we can find a subsequence $(n(k))_{k \geq 1}$ and, for every $U \in \mathcal{U}$, a measure $\lambda_U \in \mathcal{M}_\delta(S)$ such that for every $f \in C_c(S; \mathbb{R})$ with compact support in $\overline{U}$, we have

$$\lim_{k \to +\infty} \int_S f \, d\lambda_{n(k)} = \int_S f \, d\lambda_U.$$

It is clear that the restrictions of the measures $\lambda_U$ and $\lambda_V$ to $U \cap V$ must coincide. We can thus build a point measure $\lambda \in \mathcal{M}_\delta(S)$ such that for every $U \in \mathcal{U}$, the restrictions of $\lambda$ and $\lambda_U$ to $U$ coincide. Since $\mathcal{U}$ covers $S$, for every $f \in C_c(S; \mathbb{R})$, we can find a finite cover of the support of $f$ made of elements of $\mathcal{U}$, say $U_1,\ldots,U_n$. Letting $h_1,\ldots,h_n$ be an associated partition of unity as in Lemma A.8, we can write

$$\int_S f \, d\lambda_{n(k)} = \sum_{i=1}^n \int_S h_i f \, d\lambda_{n(k)},$$

and since each $h_i f$ has compact support in $\overline{U}_i$, we conclude that

$$\lim_{k\to+\infty} \int_S f \, d\lambda_{n(k)} = \sum_{i=1}^n \int_S h_i f \, d\lambda = \int_S f \, d\lambda.$$

The proof is thus complete. ∎

**Exercise 5.1.** Show that in the expression (5.1), the index set $I$ must be countable.

**Exercise 5.2.** Show that, for every $A \in \mathcal{B}(S)$, the mapping

$$\begin{cases} \mathcal{M}_\delta(S) & \to & \mathbb{R} \cup \{+\infty\} \\ \lambda & \mapsto & \lambda(A) \end{cases} \tag{5.11}$$

is measurable. Show that the smallest $\sigma$-algebra such that each of these mappings is measurable is the Borel $\sigma$-algebra of $\mathcal{M}_\delta(S)$.

**Exercise 5.3.** Let $(\lambda_n)_{n\geq 1}$ be a sequence of measures in $\mathcal{M}_\delta(S)$, and let $\lambda$ be a locally finite measure on $S$ such that, for every $f \in C_c(S;\mathbb{R})$, we have

$$\lim_{n\to+\infty} \int_S f \, d\lambda_n = \int_S f \, d\lambda. \tag{5.12}$$

Show that $\lambda$ must be a point measure.

**Exercise 5.4.** Show that the function in (5.7) defines a metric on $\mathcal{M}_\delta(S)$.

**Exercise 5.5.** Show that the space $\mathcal{M}_\delta(S)$ endowed with the metric (5.7) is complete.

## 5.2 Point processes

A *point process on $S$* is a random variable taking values in $\mathcal{M}_\delta(S)$, and it serves as a model for a random distribution of points in $S$. That is, a point process $\Lambda$ on $S$ is a measurable map from some abstract probability space $(\Omega, \mathcal{A}, \mathbb{P})$ to $\mathcal{M}_\delta(S)$ equipped with its Borel $\sigma$-algebra $\mathcal{B}(\mathcal{M}_\delta(S))$. For each $\omega \in \Omega$, the weight assigned by the point measure $\Lambda(\omega) \in \mathcal{M}_\delta(S)$ to any measurable set $A \subseteq S$ gives the number of points $\Lambda(\omega)(A) \subseteq [0,+\infty]$ in the support of $\Lambda(\omega)$ that belong to $A$. The law of the random process $\Lambda$ is the probability measure $\mathbb{P} \circ \Lambda^{-1}$ on $\mathcal{M}_\delta(S)$. Two point processes $\Lambda$ and $\Lambda'$ have the same law if, for every bounded and measurable function $F : \mathcal{M}_\delta(S) \to \mathbb{R}$,

$$\mathbb{E}F(\Lambda) = \mathbb{E}F(\Lambda'). \tag{5.13}$$

A sequence $(\Lambda_n)_{n \geqslant 1}$ of point processes on $S$ converges in law to a point process $\Lambda$ on $S$ if, for every bounded and continuous function $F : \mathcal{M}_\delta(S) \to \mathbb{R}$,

$$\lim_{n \to +\infty} \mathbb{E}F(\Lambda_n) = \mathbb{E}F(\Lambda). \tag{5.14}$$

In the same spirit as Theorem A.23, the next two results simplify these criteria by showing that it suffices to consider Laplace transform functions of the form

$$F(\lambda) := \exp\left(-\int_S f \, d\lambda\right) \tag{5.15}$$

for functions $f$ ranging in the set $C_c(S; \mathbb{R}_{\geqslant 0})$. Notice that the quantity $\int_S f \, d\lambda$ is well-defined and finite for every $f \in C_c(S; \mathbb{R}_{\geqslant 0})$ and $\lambda \in \mathcal{M}_\delta(S)$, since elements of $\mathcal{M}_\delta(S)$ are assumed to be locally finite.

**Proposition 5.4.** *Two point processes $\Lambda, \Lambda'$ on $S$ have the same law if and only if, for every $f \in C_c(S; \mathbb{R}_{\geqslant 0})$,*

$$\mathbb{E}\exp\left(-\int_S f \, d\Lambda\right) = \mathbb{E}\exp\left(-\int_S f \, d\Lambda'\right). \tag{5.16}$$

*Proof.* The direct implication is immediate. Conversely, we assume that (5.16) holds for every $f \in C_c(S; \mathbb{R}_{\geqslant 0})$. Consider the collection of sets of the form

$$\left\{ \lambda \in \mathcal{M}_\delta(S) \mid \left| \int_S f_i \, d(\lambda - \lambda_0) \right| < r_i \text{ for all } 1 \leqslant i \leqslant k \right\},$$

where $k \in \mathbb{N}$, $f_1, \ldots, f_k \in C_c(S; \mathbb{R})$, $r_1, \ldots, r_k > 0$, and $\lambda_0 \in \mathcal{M}_\delta(S)$ are arbitrary. This collection is closed under finite intersections, and the smallest $\sigma$-algebra containing it is the Borel $\sigma$-algebra. By Dynkin's $\pi$-$\lambda$ theorem (Theorem A.5), it therefore suffices to show that, for every $f_1, \ldots, f_k \in C_c(S; \mathbb{R})$, the random vectors

$$Y := \left(\int_S f_i \, d\Lambda\right)_{1 \leqslant i \leqslant k} \quad \text{and} \quad Z := \left(\int_S f_i \, d\Lambda'\right)_{1 \leqslant i \leqslant k}$$

have the same law. Decomposing each function into its positive and negative parts, and thus considering a vector of length $2k$ instead of $k$, we can without loss of generality assume that $f_1, \ldots, f_k$ take non-negative values. Fix $t \in \mathbb{R}^k_{\geqslant 0}$, and observe that applying the assumption (5.16) to the function $f := \sum_{i=1}^k t_i f_i$ reveals that $\mathbb{E}e^{-t \cdot Y} = \mathbb{E}e^{-t \cdot Z}$. Invoking the uniqueness theorem for Laplace transforms in Theorem A.23 completes the proof. ∎

**Remark 5.5.** If $S$ is an open subset of $\mathbb{R}^d$ equipped with the Euclidean topology, then Proposition 5.4 is also valid if we only assume that (5.16) holds for every $f \in C_c^\infty(S; \mathbb{R}_{\geqslant 0})$. Here $C_c^\infty(S; \mathbb{R}_{\geqslant 0})$ denotes the space of non-negative smooth functions with compact support in $S$. Indeed, this follows by a density argument similar to that in Theorem 5.2. In short, denoting by $K$ the support of $f$, we first find $K_\varepsilon$ a compact enlargement of $K$, and argue as in (5.8) that we can make $\sup_{n\geqslant 1} \mathbb{P}\{\Lambda_n(K_\varepsilon) \geqslant M\}$ arbitrarily small by taking $M$ sufficiently large. We then conclude by observing that we can find a sequence of functions in $C_c^\infty(S; \mathbb{R}_{\geqslant 0})$ with support in $K_\varepsilon$ that converges to $f$ in the uniform norm. More generally, in the statement of Proposition 5.4, we can replace $C_c(S; \mathbb{R}_{\geqslant 0})$ by any subset $\mathcal{F} \subseteq C_c(S; \mathbb{R}_{\geqslant 0})$ that satisfies the approximation property stated in Lemma 5.1.

We next aim to use Proposition 5.4 to argue that the convergence in law of a sequence of point processes can be assessed by looking only at the Laplace transform of the point processes. We start with a preliminary observation clarifying the notion of convergence in law of point processes with respect to the vague topology.

**Lemma 5.6.** *A sequence $(\Lambda_n)_{n\geqslant 1}$ of point processes on $S$ converges in law to a point process $\Lambda$ on $S$ if and only if, for every $f \in C_c(S; \mathbb{R}_{\geqslant 0})$, the sequence $\left(\int_S f\, d\Lambda_n\right)_{n\geqslant 1}$ converges in law to $\int_S f\, d\Lambda$.*

*Proof.* The direct implication is immediate. To show the converse implication we first decompose every $f \in C_c(S; \mathbb{R})$ into its positive and negative parts, $f_+$ and $f_-$, and assert the joint convergence in law of $\left(\int_S f_+\, d\Lambda_n, \int_S f_-\, d\Lambda_n\right)_{n\geqslant 1}$ towards $\left(\int_S f_+\, d\Lambda, \int_S f_-\, d\Lambda\right)$ by appealing to Proposition A.24. In particular, we have that for every $f \in C_c(S; \mathbb{R})$, the sequence $\left(\int_S f\, d\Lambda_n\right)_{n\geqslant 1}$ converges in law to $\int_S f\, d\Lambda$. We now proceed in two steps to show that $(\Lambda_n)_{n\geqslant 1}$ converges to $\Lambda$. First we argue that $(\Lambda_n)_{n\geqslant 1}$ is tight using Prokhorov's theorem, and then we prove that $\Lambda$ is its only subsequential limit using Proposition 5.4. The result then follows from Lemma A.19.

*Step 1: $(\Lambda_n)_{n\geqslant 1}$ is tight.* We fix $\varepsilon > 0$ and denote by $(f_k)_{k\geqslant 1}$ an enumeration of the countable set $\mathcal{F}$ in Lemma 5.1. Since $\left(\int_S f_k\, d\Lambda_n\right)_{n\geqslant 1}$ converges in law to $\int_S f_k\, d\Lambda$ as $n$ tends to infinity, for every $k \geqslant 1$, it is possible to find a constant $c_k > 0$ such that, for every $n \geqslant 1$,

$$\mathbb{P}\left\{\left|\int_S f_k\, d\Lambda_n\right| > c_k\right\} \leqslant \frac{\varepsilon}{2^k}.$$

This can be shown directly, or by appealing to the converse of Prokhorov's theorem (Theorem A.21) applied to the complete and separable space $\mathbb{R}$. With the sequence $(c_k)_{k \geqslant 1}$ at hand, we introduce the set of point measures

$$M := \left\{ \lambda \in \mathcal{M}_\delta(S) \mid \left| \int_S f_k \, d\lambda \right| \leqslant c_k \text{ for all } k \geqslant 1 \right\}.$$

An identical argument to that in Theorem 5.2 reveals that for every $g \in C_c(S; \mathbb{R})$,

$$\sup_{\lambda \in M} \left| \int_S g \, d\lambda \right| < +\infty.$$

It follows by the characterization of relatively compact sets in $\mathcal{M}_\delta(S)$ provided in Lemma 5.3 that $M$ is relatively compact. Moreover, a union bound and the choice of the sequence $(c_k)_{k \geqslant 1}$ imply that, for every $n \geqslant 1$,

$$\mathbb{P}\{\Lambda_n \notin \overline{M}\} \leqslant \mathbb{P}\{\Lambda_n \notin M\} \leqslant \sum_{k \geqslant 1} \mathbb{P}\left\{ \left| \int_S f_k \, d\Lambda_n \right| > c_k \right\} \leqslant \varepsilon.$$

This shows that the sequence $(\Lambda_n)_{n \geqslant 1}$ is tight. By Prokhorov's theorem (Theorem A.20), every subsequence of $(\Lambda_n)_{n \geqslant 1}$ admits a further subsequence that converges in law.

*Step 2: $(\Lambda_n)_{n \geqslant 1}$ admits a unique subsequential limit.* Let $(\Lambda_{n(k)})_{k \geqslant 1}$ be a subsequence of $(\Lambda_n)_{n \geqslant 1}$ converging in law to some point process $\Lambda'$. For every $f \in C_c(S; \mathbb{R}_{\geqslant 0})$, we have

$$\mathbb{E} \exp\left( -\int_S f \, d\Lambda \right) = \lim_{k \to +\infty} \mathbb{E} \exp\left( -\int_S f \, d\Lambda_{n(k)} \right) = \mathbb{E} \exp\left( -\int_S f \, d\Lambda' \right).$$

Proposition 5.4 implies that $\Lambda'$ and $\Lambda$ must have the same law, thereby completing the proof. ∎

**Proposition 5.7** (Convergence via Laplace transforms). *A sequence $(\Lambda_n)_{n \geqslant 1}$ of point processes on $S$ converges in law to a point process $\Lambda$ on $S$ if and only if, for every $f \in C_c(S; \mathbb{R}_{\geqslant 0})$, we have*

$$\lim_{n \to +\infty} \mathbb{E} \exp\left( -\int_S f \, d\Lambda_n \right) = \mathbb{E} \exp\left( -\int_S f \, d\Lambda \right). \tag{5.17}$$

*Proof.* The direct implication is immediate. For the converse, fix $f \in C_c(S; \mathbb{R}_{\geqslant 0})$, and notice that assumption (5.17) implies that for every $t \geqslant 0$,

$$\lim_{n \to +\infty} \mathbb{E} \exp\left( -t \int_S f \, d\Lambda_n \right) = \mathbb{E} \exp\left( -t \int_S f \, d\Lambda \right).$$

It follows by Proposition A.24 that the sequence $\left( \int_S f \, d\Lambda_n \right)_{n \geqslant 1}$ converges in law to $\int_S f \, d\Lambda$. Invoking Lemma 5.6 completes the proof. ∎

**Remark 5.8.** If $S$ is an open subset of $\mathbb{R}^d$ equipped with the Euclidean topology, then, as in Remark 5.5, we see that in the statement of Proposition 5.17, we can replace the condition "$f \in C_c(S; \mathbb{R}_{\geqslant 0})$" by the condition "$f \in C_c^\infty(S; \mathbb{R}_{\geqslant 0})$".

## 5.3 Poisson point processes

A special class of point processes known as Poisson point processes play a particularly important role in the study of extremes of i.i.d. random variables, spin glasses, and many other areas. Given a locally finite measure $\mu$ on $S$, we say that a point process $\Lambda$ on $S$ is a *Poisson point process* on $S$ with *intensity measure* $\mu$ if, for every integer $n \geq 1$ and pairwise disjoint measurable sets $(A_k)_{1 \leq k \leq n} \subseteq S$, the random variables $(\Lambda(A_k))_{1 \leq k \leq n}$ are independent Poisson random variables, with $\Lambda(A_k)$ having mean $\mu(A_k)$ for $1 \leq k \leq n$. In this definition, we understand that $\Lambda(A)$ is equal to $+\infty$ almost surely if $\mu(A) = +\infty$. A Poisson point process is therefore a point process $\Lambda$ for which we stipulate that the number of points in the support of $\Lambda$ that fall in any measurable set $A \subseteq S$ follows a Poisson distribution with mean $\mu(A)$. We also require that these Poisson counts be independent for pairwise disjoint sets.

It is not clear from the description of a Poisson point process that such an object should exist or be uniquely defined. We now describe how to generate a Poisson point process with intensity measure $\mu$. For simplicity, we will assume that $\mu(S)$ is finite; otherwise, we can decompose $\mu$ into a sum of measures of finite total mass, construct independent Poisson point processes for each of these measures, and verify that the sum of these point processes is a Poisson point process with intensity measure $\mu$. Let $(X_n)_{n \geq 1}$ be independent random variables with law

$$\overline{\mu} := \frac{\mu}{\mu(S)}, \tag{5.18}$$

let $N$ be a Poisson random variable with mean $\mu(S)$, and define

$$\Lambda := \sum_{i=1}^{N} \delta_{X_i}. \tag{5.19}$$

The following proposition states that this is the unique Poisson point process with intensity measure $\mu$.

**Proposition 5.9** (Poisson point process). *Let $\mu$ be a locally finite measure on $S$. There exists a Poisson point process with intensity measure $\mu$; when $\mu$ is finite, the point process* (5.19) *is one such. Moreover, given a Poisson point process $\Lambda$ with intensity measure $\mu$ and a measurable function $f : S \to \mathbb{R}_{\geq 0}$,*

$$\mathbb{E} \exp\left(-\int_S f \, \mathrm{d}\Lambda\right) = \exp\left(-\int_S \left(1 - e^{-f}\right) \mathrm{d}\mu\right). \tag{5.20}$$

*In particular, any two Poisson point processes with intensity measure $\mu$ have the same law.*

*Proof.* We start by proving the existence of a Poisson point process with intensity measure $\mu$. We restrict our attention to the case when $\mu$ is finite; otherwise, we can decompose $\mu$ into a sum of measures of finite total mass, construct independent Poisson point processes for each of these measures, and verify that the sum of these point processes is a Poisson point process with intensity measure $\mu$ by following the argument we now present, with slightly heavier notation. Let $\Lambda$ be the point process defined in (5.19). For every measurable function $f : S \to \mathbb{R}_{\geq 0}$, we have

$$\mathbb{E}\exp\left(-\int_S f\,d\Lambda\right) = \mathbb{E}\exp\left(-\sum_{i=1}^N f(X_i)\right)$$
$$= e^{-\mu(S)} \sum_{n=0}^{+\infty} \frac{(\mu(S))^n}{n!} \left(\int_S e^{-f}\,d\overline{\mu}\right)^n$$
$$= \exp\left(-\int_S (1-e^{-f})\,d\mu\right).$$

In particular, if $A_1, \ldots, A_k \subseteq S$ are disjoint measurable sets and $t_1, \ldots, t_k \geq 0$, then

$$\mathbb{E}\exp\left(-\sum_{i=1}^k t_i \Lambda(A_i)\right) = \prod_{i=1}^k \exp\left(-(1-e^{-t_i})\mu(A_i)\right).$$

Together with the injectivity of the Laplace transform in Theorem A.23 and the explicit form of the moment generating function of a Poisson random variable, this shows that $\bigl(\Lambda(A_i)\bigr)_{1 \leq i \leq k}$ are independent Poisson random variables with respective means $\bigl(\mu(A_i)\bigr)_{1 \leq i \leq k}$. We now turn to the proof of (5.20) and establish the uniqueness of Poisson point processes. Through a slight abuse of notation, denote by $\Lambda$ any Poisson point process with intensity measure $\mu$, and consider a simple function $f : S \to \mathbb{R}_{\geq 0}$ of the form

$$f := \sum_{i=1}^k t_i \mathbf{1}_{A_i} \qquad (5.21)$$

for some pairwise disjoint measurable sets $A_i$, and some non-negative constants $t_i \geq 0$. The definition of a Poisson point process implies that

$$\mathbb{E}\exp\left(-\int_S f\,d\Lambda\right) = \prod_{i=1}^k \exp\left((e^{-t_i}-1)\mu(A_i)\right) = \exp\left(-\int_S (1-e^{-f})\,d\mu\right),$$

where we have used the explicit form of the moment generating function of a Poisson random variable. For a general measurable function $f$, we can find a sequence of functions of the form (5.21) that converge to $f$ monotonically, see for instance Theorem 1.17 in [233], so (5.20) follows from the monotone convergence theorem. Invoking Proposition 5.4 establishes the uniqueness in law of the Poisson point process with intensity measure $\mu$ and completes the proof. ∎

The set of Poisson point processes is preserved by a number of transformations. To begin with, let us fix another locally compact separable metric space $(S', d')$ endowed with its Borel $\sigma$-algebra $\mathcal{B}(S')$, and show that mapping the points in the support of a Poisson point process on $S$ to the space $S'$ again yields a Poisson point process. Given a Poisson point process $\Lambda$ of the form (5.19) and a measurable map $f : S \to S'$, we write

$$f(\Lambda) := \sum_{i=1}^{N} \delta_{f(X_i)} \tag{5.22}$$

for the image point process on $S'$.

**Proposition 5.10** (Mapping theorem)**.** *We fix a locally finite measure $\mu$ on $S$ and a measurable map $f : S \to S'$ such that $\mu \circ f^{-1}$ is a locally finite measure on $S'$. If $\Lambda$ is a Poisson point process on $S$ with intensity measure $\mu$, then $f(\Lambda)$ is a Poisson point process on $S'$ with intensity measure $\mu \circ f^{-1}$.*

*Proof.* Proposition 5.9 implies that for any measurable function $g \in C_c(S'; \mathbb{R}_{\geq 0})$,

$$\mathbb{E}\exp\left(-\int_S g\,\mathrm{d}(f(\Lambda))\right) = \mathbb{E}\exp\left(-\int_S g \circ f\,\mathrm{d}\Lambda\right) = \exp\left(-\int_S \left(1 - e^{-g \circ f}\right)\mathrm{d}\mu\right).$$

Changing variables and combining Propositions 5.4 and 5.9 completes the proof. ∎

In addition to using a Poisson point process $\Lambda$ on $S$ to construct a Poisson point process on $S'$, we can also use it to construct a Poisson point process on the product space $S \times S'$ by appending a random coordinate to each point in its support. The simplest instance of this construction is to fix a sequence $(Y_i)_{i \geq 1}$ of i.i.d. random variables independent of the random variables $(X_i)_{i \geq 1}$ and $N$ appearing in the representation (5.19) of $\Lambda$, and to define the point process

$$\Lambda^* := \sum_{i=1}^{N} \delta_{(X_i, Y_i)}. \tag{5.23}$$

Here we assume again that the measure $\mu$ is finite to simplify a bit the discussion, but this restriction can easily be lifted. If we write $\mu$ for the intensity measure of $\Lambda$ and $\mathbb{P}_Y$ for the law of $Y_1$, then the explicit construction of a Poisson point process shows that $\Lambda^*$ is a Poisson point process on $S \times S'$ with intensity measure $\mu \otimes \mathbb{P}_Y$. More generally, if we use a transition function $K : S \times \mathcal{B}(S') \to [0,1]$ to generate, independently for each $i \in I$, a point $m(X_i) \in S'$ from the distribution $K(X_i, \cdot)$, and define the *marked point process*

$$\Lambda^* := \sum_{i=1}^{N} \delta_{(X_i, m(X_i))} \tag{5.24}$$

on $S \times S'$, then the *marking theorem* ensures that this is a Poisson point process with intensity measure $\mu^*(s, s') := \mu(\mathrm{d}s) \otimes K(s, \mathrm{d}s')$. We call a function $K : S \times \mathcal{B}(S') \to [0,1]$ a *transition function* if, for each $s \in S$, the section $K(s, \cdot)$ is a probability measure on $\mathcal{B}(S')$, while, for each $A \in \mathcal{B}(S')$, the section $K(\cdot, A)$ is measurable on $S$. The marking theorem is therefore a generalization of the i.i.d. construction in (5.23).

**Proposition 5.11** (Marking theorem). *If $\Lambda$ is a Poisson point process on $S$ with intensity measure $\mu$, then the point process $\Lambda^*$ in (5.24) is a Poisson point process on $S \times S'$ with intensity measure $\mu^*(s, s') := \mu(ds) \otimes K(s, ds')$.*

*Proof.* Fix a measurable function $f \in C_c(S \times S'; \mathbb{R}_{\geq 0})$, and observe that

$$\mathbb{E}\exp\left(-\int_S f \, d\Lambda^*\right) = \mathbb{E}\exp\left(-\sum_{i=1}^N f(X_i, m(X_i))\right)$$

$$= e^{-\mu(S)} \sum_{n=0}^{+\infty} \frac{(\mu(S))^n}{n!} \left(\int_S \mathbb{E}e^{-f(s,m(s))} \, d\overline{\mu}(s)\right)^n$$

$$= e^{-\mu(S)} \sum_{n=0}^{+\infty} \frac{(\mu(S))^n}{n!} \left(\int_S \int_{S'} e^{-f(s,s')} K(s, ds') \, d\overline{\mu}(s)\right)^n$$

$$= \exp\left(-\int_S \int_{S'} \left(1 - e^{-f(s,s')}\right) K(s, ds') \, d\mu(s)\right),$$

where the last equality uses the fact that $K(s, \cdot)$ is a probability measure. Invoking Propositions 5.4 and 5.9 completes the proof. ∎

For later usage, we now record the following explicit formula translating averages involving Poisson point processes into averages involving their intensity measure.

**Proposition 5.12** (Poisson average). *If $\Lambda$ is a Poisson point process on $S$ with intensity measure $\mu$, then for any measurable function $f : S \to \mathbb{R}_{\geq 0}$,*

$$\mathbb{E} \int_S f \, d\Lambda = \int_S f \, d\mu. \tag{5.25}$$

*Proof.* Approximating $f$ by a sequence of simple functions as in Proposition 5.9, it suffices to consider the case when $f : S \to \mathbb{R}_{\geq 0}$ is of the form

$$f := \sum_{i=1}^k t_i \mathbf{1}_{A_i}$$

for some pairwise disjoint measurable sets $A_i$ and non-negative constants $t_i \in \mathbb{R}_{\geq 0}$. The definition of a Poisson point process implies that

$$\mathbb{E} \int_S f \, d\Lambda = \sum_{i=1}^k t_i \mathbb{E}\Lambda(A_i) = \sum_{i=1}^k t_i \mu(A_i) = \int_S f \, d\mu,$$

as required. ∎

**Exercise 5.6.** Let $(\Lambda_n)_{n \geq 1}$ be a sequence of point processes that converges in law to a Poisson point process $\Lambda$ with intensity measure $\mu$, and fix a measurable set $A \subseteq S$ with $\mu(\partial A) = 0$. Show that the sequence of random variables $(\Lambda_n(A))_{n \geq 1}$ converges in law to $\Lambda(A)$.

## 5.4 Extremes of i.i.d. random variables

We now analyze the asymptotic behaviour of extreme values of i.i.d. random variables, which is most conveniently phrased in terms of the convergence to a Poisson point process. We view this as a fundamental result in probability theory, essentially on par with the central limit theorem, and we hope that the reader will appreciate to find a brief and self-contained presentation of this result. It will also inform our analysis of spin glasses, as will be discussed in the next sections and chapter.

Given a sequence $(X_n)_{n \geq 1}$ of i.i.d. random variables, we denote their running maximum by

$$M_n := \max_{1 \leq k \leq n} X_k. \qquad (5.26)$$

The most basic question in the study of extreme values is whether there exist sequences $(a_n)_{n \geq 1} \subseteq \mathbb{R}_{>0}$ and $(b_n)_{n \geq 1} \subseteq \mathbb{R}$ such that the normalized running maximum,

$$a_n^{-1}(M_n - b_n), \qquad (5.27)$$

converges in law to a non-trivial limit. Under suitable assumptions on the law of the random variables, it is elementary to obtain the answer to this question. Poisson point processes allow us to extend this result and capture the joint convergence in law of the largest values taken by the sequence $(X_k)_{k \leq n}$, in the limit of large $n$. The next proposition makes this precise in the case of random variables with a polynomial tail distribution.

**Proposition 5.13** (Polynomial tail). *Let $(X_n)_{n \geq 1}$ be i.i.d. random variables, and suppose that there exist $c > 0, \zeta > 0$ such that, as $x$ tends to infinity,*

$$\mathbb{P}\{X_1 \geq x\} \sim \frac{c}{x^\zeta}. \qquad (5.28)$$

*Then the point process*

$$\Lambda_n := \sum_{k=1}^n \delta_{n^{-1/\zeta} X_k} \qquad (5.29)$$

*converges in law on $\mathbb{R}_{>0}$ to the Poisson point process with intensity $\mathrm{d}\mu(x) := \frac{c\zeta}{x^{\zeta+1}} \mathrm{d}x$.*

In Proposition 5.13, and similarly throughout, we say for convenience that a sequence of point processes converges in law on $\mathbb{R}_{>0}$ to mean that it converges in law in the sense of the vague topology on $\mathcal{M}_\delta(\mathbb{R}_{>0})$. We implicitly understand that we have equipped $\mathbb{R}_{>0}$ with the Euclidean topology. One may take issue with the fact that the definition of $\Lambda_n$ may involve Dirac masses at locations outside of $\mathbb{R}_{>0}$ if $X_1$ can take values in $(-\infty, 0]$. These Dirac masses can be discarded if desired; in any case, they are not seen by the topology of $\mathcal{M}_\delta(\mathbb{R}_{>0})$, so their presence is irrelevant to the question at hand. The display in (5.28) means that $x^\zeta \mathbb{P}\{X_1 \geq x\}$ converges to $c$ as $x$ tends to infinity.

*Proof.* By Proposition 5.7, it suffices to investigate the convergence of the Laplace transform of $\Lambda_n$, and by Remark 5.8 it suffices to verify (5.17) for a given smooth and compactly supported function $f \in C_c^\infty(\mathbb{R}_{>0}; \mathbb{R}_{\geq 0})$. The independence of $(X_k)_{1 \leq k \leq n}$ implies that

$$\mathbb{E} \exp\left(-\int_0^{+\infty} f \, d\Lambda_n\right) = \mathbb{E} \exp\left(-\sum_{k=1}^n f(n^{-1/\zeta} X_k)\right) = \left(\mathbb{E} \exp\left(-f(n^{-1/\zeta} X_1)\right)\right)^n.$$

The slight difficulty of this proof is that the assumption (5.28) concerns the cumulative distribution function of $X_1$, while it would, for instance, be easier to directly have access to an assumption on the density of $X_1$. We resolve this through integration by parts,

$$\mathbb{E}\left(1 - \exp\left(-f(n^{-1/\zeta} X_1)\right)\right) = \mathbb{E} \int_0^{+\infty} f'(y) e^{-f(y)} \mathbf{1}_{\{y \leq n^{-1/\zeta} X_1\}} \, dy$$

$$= \int_0^{+\infty} f'(y) e^{-f(y)} \mathbb{P}\{X_1 \geq n^{1/\zeta} y\} \, dy.$$

The first equality leverages the fact that $f(0) = 0$. We now fix $\varepsilon > 0$, and observe that by assumption (5.28), there exists $x_\varepsilon > 0$ such that for every $x \geq x_\varepsilon$,

$$\left|x^\zeta \mathbb{P}\{X_1 \geq x\} - c\right| \leq \varepsilon. \tag{5.30}$$

Recalling also that $f$ has compact support in $\mathbb{R}_{>0}$, we see that, for $n = n(\varepsilon)$ sufficiently large, we have

$$\left|\int_0^{+\infty} f'(y) e^{-f(y)} \mathbb{P}\{X_1 \geq n^{1/\zeta} y\} \, dy - \frac{c}{n} \int_0^{+\infty} f'(y) e^{-f(y)} y^{-\zeta} \, dy\right|$$

$$\leq \frac{\varepsilon}{n} \int_0^{+\infty} |f'(y)| e^{-f(y)} y^{-\zeta} \, dy.$$

The last integral is finite, and can be absorbed up to a reparametrization of $\varepsilon > 0$. We can now "undo" the integration by parts,

$$\int_0^{+\infty} f'(y) e^{-f(y)} y^{-\zeta} \, dy = \int_0^{+\infty} f'(y) e^{-f(y)} \int_0^{+\infty} \frac{\zeta}{x^{\zeta+1}} \mathbf{1}_{\{x \geq y\}} \, dx \, dy$$

$$= \int_0^{+\infty} \int_0^x f'(y) e^{-f(y)} \, dy \, \frac{\zeta}{x^{\zeta+1}} \, dx$$

$$= \int_0^{+\infty} (1 - e^{-f(x)}) \frac{\zeta}{x^{\zeta+1}} \, dx.$$

We have thus shown that, for every $\varepsilon > 0$ and every $n$ sufficiently large,

$$\left|\mathbb{E} \exp\left(-f(n^{-1/\zeta} X_1)\right) - \left(1 - \frac{c}{n} \int_0^{+\infty} (1 - e^{-f(x)}) \frac{\zeta}{x^{\zeta+1}} \, dx\right)\right| \leq \frac{\varepsilon}{n}.$$

In other words, we have the asymptotic expansion as $n$ tends to infinity

$$\mathbb{E}\exp\bigl(-f(n^{-1/\zeta}X_1)\bigr) = 1 - \frac{c}{n}\int_0^{+\infty}\bigl(1 - e^{-f(x)}\bigr)\frac{\zeta}{x^{\zeta+1}}\,dx + o\Bigl(\frac{1}{n}\Bigr). \tag{5.31}$$

It thus follows that

$$\lim_{n\to+\infty}\Bigl(\mathbb{E}\exp\bigl(-f(n^{-1/\zeta}X_1)\bigr)\Bigr)^n = \exp\Bigl(-\int_0^{+\infty}\bigl(1 - e^{-f(x)}\bigr)\frac{c\zeta}{x^{\zeta+1}}\,dx\Bigr).$$

By Proposition 5.9, the latter expression is the Laplace transform of a Poisson point process on $\mathbb{R}_{>0}$ with intensity $d\mu(x) = \frac{c\zeta}{x^{\zeta+1}}\,dx$. This completes the proof.  ∎

Since $\int_1^{+\infty} x^{-(\zeta+1)}\,dx < +\infty$, the limit Poisson point process appearing in Proposition 5.13 contains only a finite number of points in the interval $[1,+\infty)$. Heuristically, it thus makes sense to think of Proposition 5.13 as capturing the law of the largest values of $(X_k)_{k\leq n}$. However, Proposition 5.13 actually falls short of allowing us to assert this; in particular, one cannot infer the asymptotic behaviour of the running maximum $M_n$ from Proposition 5.13 as stated. The reason is that the topology of $\mathcal{M}_\delta(\mathbb{R}_{>0})$ does not allow us to test for properties that require "observing" the point process arbitrarily far away towards infinity; for instance, if we add a Dirac mass at position $n$ to $\Lambda_n$, then this does not affect the convergence of $\Lambda_n$ in $\mathcal{M}_\delta(\mathbb{R}_{>0})$. Fortunately, it is easy to amend for this shortcoming by choosing a topology more adapted to our needs. Indeed, we can compactify the interval $\mathbb{R}_{>0}$ at $+\infty$ to obtain the set $(0,+\infty]$, equipped for instance with the distance

$$d(x,y) := |e^{-x} - e^{-y}|, \tag{5.32}$$

with the understanding that $e^{-\infty} = 0$. Writing $d\mu(x) := \frac{c\zeta}{x^{\zeta+1}}\,dx$, we note that the Poisson point process with intensity $\mu$ appearing in Proposition 5.13 is a valid point process in $\mathcal{M}_\delta((0,+\infty])$. Indeed, since $\mu([1,+\infty]) = \mu([1,+\infty)) < +\infty$, the measure $\mu$ is locally finite also on $(0,+\infty]$, and likewise for the Poisson point process. We can strengthen the statement of Proposition 5.13 so that it is now stated with respect to the more demanding topology of $\mathcal{M}_\delta((0,+\infty])$.

**Proposition 5.14.** *Let $(X_n)_{n\geq 1}$ be i.i.d. random variables, and suppose that there exist $c > 0$, $\zeta > 0$ such that (5.28) holds as $x$ tends to infinity. The point process $\Lambda_n$ in (5.29) converges in law on $(0,+\infty]$ to the Poisson point process with intensity $d\mu(x) := \frac{c\zeta}{x^{\zeta+1}}\,dx$.*

*Proof.* The proof is almost identical to that of Proposition 5.13. The only difference is that the set $C_c((0,+\infty];\mathbb{R}_{\geq 0})$ contains functions that may not have compact support in $\mathbb{R}_{>0}$, but instead converge to some limit at $+\infty$. As before, it suffices to verify (5.17) for functions $f \in C_c^\infty((0,+\infty];\mathbb{R}_{\geq 0})$, and one can also impose that $f$ is constant over $[M,+\infty]$ for $M$ is sufficiently large, by density. The proof of Proposition 5.13 then applies without further modification.  ∎

Propositions 5.13 and 5.14 can be refined so that we also learn about the index at which a given extreme value is occurring; see Exercise 5.7.

To further clarify that Proposition 5.14 indeed captures the asymptotic behaviour of the extremes of the sequence $(X_n)_{n \geqslant 1}$, we now deduce from it the asymptotic behaviour of the $k$ largest values of $(X_\ell)_{\ell \leqslant n}$, for each fixed $k$. To state this result, we first make a simple observation concerning the Poisson point process $\Lambda$ appearing in Proposition 5.14. We recall that the Poisson point process $\Lambda$ has only a finite number of points falling in the interval $[1, +\infty)$, so we can enumerate its support decreasingly. Moreover, since $\int_0^{+\infty} x^{-(\zeta+1)} \, dx$ is infinite, the support of $\Lambda$ is infinite almost surely. So with probability one, there exist $u_1 \geqslant u_2 \geqslant \cdots$ such that $\Lambda = \sum_{n=1}^{+\infty} \delta_{u_n}$.

**Proposition 5.15.** *Let $(X_n)_{n \geqslant 1}$ be i.i.d. random variables, and suppose that there exist $c > 0$, $\zeta > 0$ such that (5.28) holds as x tends to infinity. Denote by $X_{1,n} \geqslant X_{2,n} \geqslant \cdots \geqslant X_{n,n}$ the decreasing ordering of the sequence $(X_k)_{1 \leqslant k \leqslant n}$. Let $\Lambda$ be a Poisson point process with intensity $d\mu(x) := \frac{c\zeta}{x^{\zeta+1}} dx$, and let $u_1 \geqslant u_2 \geqslant \cdots$ be such that $\Lambda = \sum_{n=1}^{+\infty} \delta_{u_n}$. For each integer $k \geqslant 1$, the vector*

$$n^{-1/\zeta}(X_{1,n}, X_{2,n}, \ldots, X_{k,n}) \tag{5.33}$$

*converges in law to $(u_1, \ldots, u_k)$ as n tends to infinity.*

*Proof.* We introduce the set

$$\mathcal{C} := \Big\{ \lambda \in \mathcal{M}_\delta((0, +\infty]) \mid$$
$$\lambda = \sum_{n=1}^{+\infty} \delta_{y_n} \text{ for a strictly decreasing sequence } (y_n)_{n \geqslant 1} \subseteq \mathbb{R}_{>0} \Big\}$$

of point processes with no repeated points in their support. Recall that every measure in $\mathcal{M}_\delta((0, +\infty])$ must be locally finite, and in particular, every $\lambda \in \mathcal{M}_\delta((0, +\infty])$ has a finite number of points in $[1, +\infty]$. Fix a positive integer $k$, and denote by $F : \mathcal{M}_\delta((0, +\infty]) \to \mathbb{R}^k$ the mapping such that $F(\lambda)$ outputs the vector of the $k$ largest elements in the support of $\lambda$; while this will be irrelevant for our purposes, we can decide to complete $F(\lambda)$ with zeros in case $\lambda$ has fewer than $k$ elements in its support. The key step of the proof is to observe that the mapping $F$ is continuous at every $\lambda \in \mathcal{C}$. In other words, if $(\lambda_n)_{n \geqslant 1} \subseteq \mathcal{M}_\delta((0, +\infty])$ converges to some $\lambda \in \mathcal{C}$, then it must be that $F(\lambda_n)$ converges to $F(\lambda)$. Indeed, let $(y_n)_{n \geqslant 1} \subseteq \mathbb{R}_{>0}$ be a strictly decreasing sequence with $\lambda = \sum_{n=1}^{+\infty} \delta_{y_n}$, and observe that every function of the form

$$f(y) := \min((a(y-b))_+, 1),$$

for arbitrary $a, b > 0$, belongs to $C_c((0, +\infty]; \mathbb{R})$. Choosing $a = k$ and $b = y_1 \pm k^{-1}$, with $k > 0$ sufficiently large, it is not difficult to conclude that the largest point in $\lambda_n$ converges to the largest point in $\lambda$. We can repeat the procedure recursively, choosing $b = y_2 \pm k^{-1}$, then $b = y_3 \pm k^{-1}$, and so on, thereby concluding that $(F(\lambda_n))_{n \geqslant 1}$ converges to $F(\lambda)$ as claimed. At this point, recall the notation $\Lambda_n$ in (5.29). Since the intensity measure $d\mu(x) := \frac{c\zeta}{x^{\zeta+1}} dx$ does not contain any atoms, there are no repeated points in $\Lambda$, and the point process $\Lambda$ thus belongs to $\mathcal{C}$ with probability one. By Proposition 5.14 and Exercise A.10, we deduce that $F(\Lambda_n)$ converges in law to $F(\Lambda)$, as desired. ∎

Recall the definition of the random variable $M_n$ in (5.26). As a simple consequence of Proposition 5.15, which can also be established by a very direct and elementary computation, we obtain the convergence in law of $n^{-1/\zeta} M_n$. The limit law is characterized by the fact that, for every $x > 0$,

$$\lim_{n \to +\infty} \mathbb{P}\{M_n \leqslant n^{1/\zeta} x\} = \exp\left(-c \int_x^{+\infty} \frac{\zeta}{y^{\zeta+1}} dy\right) = \exp\left(-cx^{-\zeta}\right). \tag{5.34}$$

A random variable whose cumulative distribution function is given, for $x > 0$, by $x \mapsto \exp(-x^{-\zeta})$ is called a *Fréchet law* of parameter $\zeta > 0$. For i.i.d. random variables $(X_n)_{n \geqslant 1}$, in order that there exist sequences $(a_n)_{n \geqslant 1} \subseteq \mathbb{R}_{>0}$ and $(b_n)_{n \geqslant 1} \subseteq \mathbb{R}$ such that the normalized running maximum $a_n^{-1}(M_n - b_n)$ converges in law to a Fréchet law of parameter $\zeta > 0$, it is necessary and sufficient that, for every $t > 0$, we have

$$\lim_{x \to +\infty} \frac{\mathbb{P}\{X_1 \geqslant tx\}}{\mathbb{P}\{X_1 \geqslant x\}} = t^{-\zeta}, \tag{5.35}$$

see for instance Theorem 1.6.2 in [160]. In fact, under this assumption, we can impose that $b_n = 0$ and $\zeta \log(a_n) \sim \log(n)$, and show analogues of Propositions 5.13, 5.14, and 5.15, provided that we re-scale the random variables by $a_n^{-1}$ in place of $n^{-1/\zeta}$. The proof can be obtained by combining the arguments we have seen with those in the proof of Theorem 1.6.2 in [160].

Besides Fréchet laws, there are two other possible classes of limit laws for the re-scaled maximum of i.i.d. random variables. The first class is associated with random variables with a light tail, such as exponential or Gaussian random variables, or the logarithm of a random variable satisfying (5.28). For simplicity, we only give the example of Gaussian random variables here, but the reader can work out other examples as well.

**Proposition 5.16** (Light tail, example of Gaussians). *Let $(X_n)_{n \geqslant 1}$ be independent standard Gaussian random variables, and let*

$$a_n := \left(2\log n - \log\log n - \log(4\pi)\right)^{\frac{1}{2}}. \tag{5.36}$$

*The point process*

$$\Lambda_n := \sum_{k=1}^n \delta_{a_n(X_k - a_n)} \tag{5.37}$$

*converges in law on $\mathbb{R}$ to a Poisson point process with intensity* $d\mu(x) := e^{-x} dx$.

*Proof.* Fix a function $f \in C_c(\mathbb{R}; \mathbb{R}_{\geq 0})$ of compact support, and observe that

$$\mathbb{E}\exp\left(-\int_{\mathbb{R}} f \, d\Lambda_n\right) = \mathbb{E}\exp\left(-\sum_{k=1}^{n} f(a_n(X_k - a_n))\right)$$

$$= \left(\frac{1}{\sqrt{2\pi a_n^2}} \int_{\mathbb{R}} e^{-f(x)} \exp\left(-\frac{(x+a_n^2)^2}{2a_n^2}\right) dx\right)^n$$

$$= \exp\left[n \log\left(1 - \frac{1}{\sqrt{2\pi a_n^2}} \int_{\mathbb{R}} (1 - e^{-f(x)}) \exp\left(-\frac{(x+a_n^2)^2}{2a_n^2}\right) dx\right)\right].$$

Using the definition of the constant $a_n$ in (5.36) and the fact that $f$ is of compact support reveals that for $x$ in the support of $f$,

$$\frac{n}{\sqrt{2\pi a_n^2}} \exp\left(-\frac{(x+a_n^2)^2}{2a_n^2}\right) = \frac{n}{\sqrt{4\pi \log(n)}} \exp\left(-x - \frac{1}{2}a_n^2\right) + o(1) = e^{-x} + o(1).$$

It follows by a Taylor expansion of the logarithm that

$$\mathbb{E}\exp\left(-\int_{\mathbb{R}} f \, d\Lambda_n\right) = \exp\left(-\int_{\mathbb{R}} (1 - e^{-f(x)}) e^{-x} dx\right) + o(1).$$

Invoking Propositions 5.4 and 5.9 completes the proof. ∎

As before, we can upgrade this statement by compactifying $\mathbb{R}$ at $+\infty$, and then work out the corresponding version of Proposition 5.15.

The last possible class of limit laws concerns random variables whose support is bounded to the right, and with a probability to fall close to the top of the support that decays in a power-law fashion. Examples include the uniform law on an interval, or the law of $-1/X_1$ for $X_1$ satisfying (5.28).

**Proposition 5.17** (Polynomial tail near a point). *Let $(X_n)_{n \geq 1}$ be i.i.d. random variables, and suppose that there exist $c > 0$, $\alpha > 0$ such that, as $x < 0$ tends to $0$,*

$$\mathbb{P}\{X_1 \geq x\} \sim c|x|^\alpha. \tag{5.38}$$

*Then the point process*

$$\Lambda_n := \sum_{k=1}^{n} \delta_{n^{1/\alpha} X_k} \tag{5.39}$$

*converges in law on $\mathbb{R}_{\leq 0}$ to a Poisson point process with intensity* $d\mu(x) := c\alpha|x|^{\alpha-1} dx$.

*Proof.* For a change, we prove this by appealing to Proposition 5.14 and the mapping theorem (Proposition 5.10). Letting $x$ increase to zero in (5.38) shows that $X_1$ takes values in $\mathbb{R}_{<0}$, so the random variable $Y_1 := -1/X_1$ is well-defined and takes values in $\mathbb{R}_{>0}$ with probability one. Moreover, as $x$ tends to infinity,

$$\mathbb{P}\{Y_1 \geq x\} = \mathbb{P}\{X_1 \geq -1/x\} \sim \frac{c}{x^\alpha}.$$

It follows by Proposition 5.14 that the point process

$$\sum_{k=1}^{n} \delta_{n^{-1/\alpha} Y_k}$$

converges in law on $(0, +\infty]$ to a Poisson point process with intensity $d\nu(x) := \frac{c\alpha}{x^{\alpha+1}} dx$. Together with the mapping theorem applied to the function $f(x) := -\frac{1}{x}$, this implies that the point process

$$\Lambda_n := \sum_{k=1}^{n} \delta_{f(n^{-1/\alpha} Y_k)} = \sum_{k=1}^{n} \delta_{n^{1/\alpha} X_k}$$

converges in law on $\mathbb{R}_{\leq 0}$ to a Poisson point process with intensity $\mu := \nu \circ f^{-1}$. Observe that for any $0 < a < b$,

$$\mu([a,b]) = \int_{-1/a}^{-1/b} \frac{c\alpha}{x^{\alpha+1}} dx = \int_a^b c\alpha(-y)^{\alpha+1} \frac{dy}{y^2} = \int_a^b c\alpha |y|^{\alpha-1} dy,$$

where we have used the change of variables $y = -1/x$. Invoking Propositions 5.4 and 5.9 completes the proof. ∎

As in Proposition 5.15, one can also show that the $k$ largest values of $\Lambda_n$ converge in law to the $k$ largest values of the limit point process.

**Exercise 5.7.** Under the assumptions of Proposition 5.13, show that the point process

$$\Lambda_n := \sum_{k=1}^{+\infty} \delta_{(k/n, n^{-1/\zeta} X_k)} \tag{5.40}$$

converges in law on $\mathbb{R}_{\geq 0} \times \mathbb{R}_{>0}$ to the Poisson point process with intensity $d\mu(t,x) := \frac{c\zeta}{x^{\zeta+1}} dt\, dx$.

**Exercise 5.8.** Fix $\zeta \in (0,1)$, and let $(X_n)_{n \geq 1}$ be i.i.d. random variables taking values in $\mathbb{R}_{\geq 0}$ such that, as $x$ tends to infinity,

$$\mathbb{P}\{X_1 \geq x\} \sim x^{-\zeta}.$$

Let $\Lambda$ be a Poisson point process with intensity $d\mu(x) := \frac{\zeta}{x^{\zeta+1}} dx$, and recall that there exist $u_1 \geq u_2 \geq \cdots$ such that $\Lambda = \sum_{n=1}^{+\infty} \delta_{u_n}$.

(i) Using that $\zeta < 1$, show that $\int x \, d\Lambda(x) = \sum_{n=1}^{+\infty} u_n$ is finite almost surely. For every $n \geq 1$, we set

$$v_n := \frac{u_n}{\sum_{k=1}^{+\infty} u_k}. \tag{5.41}$$

The goal of this exercise is to show that for each fixed $k$, after normalization by $\sum_{\ell=1}^{n} X_\ell$, the ordered $k$ largest values among $(X_\ell)_{\ell \leq n}$ converge to $(v_1, \ldots, v_k)$. Although we will not show it, we mention that this is equivalent to the statement that the point process

$$\sum_{\ell=1}^{n} \delta_{X_\ell/(\sum_{i=1}^{n} X_i)}$$

converges in law on $\mathbb{R}_{>0}$ to $\sum_{n=1}^{+\infty} \delta_{v_n} = \int \delta_{x/\int y \, d\Lambda(y)} \, d\Lambda(x)$.

(ii) Show that, for every $a > 0$ sufficiently large, we have

$$\mathbb{E}(X_1 \mathbf{1}_{\{X_1 \in [a,2a]\}}) \leq 2a^{1-\zeta}. \tag{5.42}$$

(iii) Deduce that

$$\lim_{\varepsilon \to 0} \limsup_{n \to +\infty} \frac{1}{n^{1/\zeta}} \sum_{k=1}^{n} \mathbb{E}(X_k \mathbf{1}_{\{X_k \leq \varepsilon n^{1/\zeta}\}}) = 0. \tag{5.43}$$

(iv) Use this to show that for every $\varepsilon > 0$,

$$\lim_{K \to +\infty} \limsup_{n \to +\infty} \mathbb{P}\left\{ \sum_{k=K}^{n} X_{k,n} \geq \varepsilon n^{1/\zeta} \right\} = 0. \tag{5.44}$$

(v) Denoting by $X_{1,n} \geq X_{2,n} \geq \cdots \geq X_{n,n}$ the decreasing ordering of the sequence $(X_k)_{1 \leq k \leq n}$, conclude that for every integer $k \geq 1$, the vector

$$\left( \frac{X_{1,n}}{\sum_{\ell=1}^{n} X_\ell}, \ldots, \frac{X_{k,n}}{\sum_{\ell=1}^{n} X_\ell} \right) \tag{5.45}$$

converges in law to $(v_1, \ldots, v_k)$. It may be helpful to refer to Exercises A.10 and A.11.

## 5.5 Poisson-Dirichlet processes

In the context of mean-field spin glasses, the energies attributed to different configurations of the system are not independent random variables. However, one might imagine that, possibly after clumping together neighbouring configurations in a suitable way, we ultimately end up with Gibbs weights whose extreme values resemble the sequence $(v_n)_{n \geq 1}$ appearing in Exercise 5.7. This intuition turns out to be valid. For a toy model of spin glasses called the random energy model, in which we assume that the energy values of different configurations are independent, one can more readily justify this intuition, as will be explained in more detail in Section 6.3.

In fact, we will ultimately aim to capture not only the asymptotic behaviour of the largest weights of the Gibbs measure, but also how the configurations associated with these extreme values are organized in space. Remarkably, this richer structure is still described by a canonical object called a Poisson-Dirichlet cascade (or also a Ruelle probability cascade, after [234]). This object is built from a hierarchy of Poisson point processes as those appearing in Proposition 5.13. In this section, in order to prepare for the study of Poisson-Dirichlet cascades and their universality, we delve more deeply into the analysis of this Poisson point process.

The *Poisson-Dirichlet point process* with parameter $\zeta > 0$ is the Poisson point process $\Lambda$ on $\mathbb{R}_{>0}$ with intensity measure

$$\mathrm{d}\mu(x) := \frac{\zeta}{x^{\zeta+1}} \, \mathrm{d}x. \tag{5.46}$$

For each $\varepsilon > 0$, the Poisson random variable $\Lambda[\varepsilon, +\infty)$ has mean $\mu[\varepsilon, +\infty)$, so the support of $\Lambda$ has finitely many points in the interval $[\varepsilon, +\infty)$. Since $\mu(\mathbb{R}_{>0})$ is infinite, the support of $\Lambda$ has infinitely many points. We can thus order the random support of $\Lambda$ into an infinite decreasing sequence $(u_n)_{n \geq 1}$, so that this Poisson point process can be represented as

$$\Lambda := \sum_{n=1}^{+\infty} \delta_{u_n}. \tag{5.47}$$

We will refer to the process $(u_n)_{n \geq 1}$ as the *ordered Poisson-Dirichlet point process* with parameter $\zeta > 0$. In order to be able to normalize this sequence of positive weights into a probability measure, we need to restrict ourselves to the range $\zeta \in (0, 1)$ for which the sequence $(u_n)_{n \geq 1}$ is summable, as shown in the next proposition. For future reference, it will also be convenient to record that for $\zeta$ in this range, the moments of order less than $\zeta$ of the sum are finite.

**Proposition 5.18.** *The ordered Poisson-Dirichlet point process* $(u_n)_{n \geq 1}$ *with parameter* $\zeta > 0$ *is almost surely summable if and only if* $\zeta \in (0, 1)$. *Moreover, if* $\zeta \in (0, 1)$, *then for every* $0 < a < \zeta$,

$$\mathbb{E}\Big(\sum_{n=1}^{+\infty} u_n\Big)^a < +\infty. \tag{5.48}$$

*Proof.* To see that $(u_n)_{n \geq 1}$ is summable if and only if $\zeta \in (0, 1)$, recall the representation (5.47) of the Poisson point process $\Lambda$ with intensity measure (5.46), and observe that

$$\sum_{n=1}^{+\infty} u_n = \int_{\mathbb{R}_{>0}} x \, \mathrm{d}\Lambda(x) = \int_{(0,1]} x \, \mathrm{d}\Lambda(x) + \int_{(1,+\infty)} x \, \mathrm{d}\Lambda(x).$$

The second of these integrals is finite since the support of $\Lambda$ only has finitely many points in the interval $(1,+\infty)$. If $\zeta < 1$, then the second integral is almost surely finite since, by Proposition 5.12, its expectation is finite, as

$$\mathbb{E}\int_{(0,1]} x\,\mathrm{d}\Lambda(x) = \int_{(0,1]} x\,\mathrm{d}\mu(x) = \int_0^1 \frac{\zeta}{x^\zeta}\,\mathrm{d}x.$$

For $\zeta \geqslant 1$, we write

$$\int_{(0,1]} x\,\mathrm{d}\Lambda(x) = \sum_{n=0}^{+\infty} \int_{(2^{-(n+1)},2^{-n}]} x\,\mathrm{d}\Lambda(x) \geqslant \sum_{n=0}^{+\infty} 2^{-(n+1)}\Lambda(2^{-(n+1)},2^{-n}]. \quad (5.49)$$

The random variables $\left(\Lambda(2^{-(n+1)},2^{-n}]\right)_{n\geqslant 0}$ are independent Poisson random variables, and the mean of $\Lambda(2^{-(n+1)},2^{-n}]$ is

$$\int_{2^{-(n+1)}}^{2^{-n}} \frac{\zeta}{x^{\zeta+1}}\,\mathrm{d}x = 2^{\zeta n}(2^\zeta - 1).$$

In particular, when $\zeta \geqslant 1$, we can find a constant $c > 0$ such that for every $n \geqslant 0$,

$$\mathbb{P}\{2^{-(n+1)}\Lambda(2^{-(n+1)},2^{-n}] \geqslant c\} \geqslant c.$$

A Borel-Cantelli argument thus ensures that the integral on the left side of (5.49) diverges almost surely in this case.

Assuming that $\zeta \in (0,1)$, a similar argument can be used to bound the moments (5.48). Indeed, using that $x^a \leqslant 1 + x$ and $(x+y)^a \leqslant x^a + y^a$ for $x,y \geqslant 0$, we see that

$$\mathbb{E}\left(\sum_{n=1}^{+\infty} u_n\right)^a \leqslant \mathbb{E}\left(\int_{(0,1]} x\,\mathrm{d}\Lambda(x)\right)^a + \mathbb{E}\left(\int_{(1,+\infty)} x\,\mathrm{d}\Lambda(x)\right)^a$$

$$\leqslant 1 + \mathbb{E}\int_{(0,1]} x\,\mathrm{d}\Lambda(x) + \mathbb{E}\int_{(1,+\infty)} x^a\,\mathrm{d}\Lambda(x)$$

$$\leqslant 1 + \int_0^1 \frac{\zeta}{x^\zeta}\,\mathrm{d}x + \int_1^{+\infty} \frac{\zeta}{x^{1+\zeta-a}}\,\mathrm{d}x.$$

Leveraging the assumption that $a < \zeta$ completes the proof. ∎

Henceforth, unless otherwise stated, whenever we speak of a Poisson-Dirichlet point process, we implicitly assume that $\zeta \in (0,1)$. Relying on the summability of $(u_n)_{n\geqslant 1}$, we can define the *Poisson-Dirichlet process* $(v_n)_{n\geqslant 1}$ with parameter $\zeta \in (0,1)$ by normalizing the sequence $(u_n)_{n\geqslant 1}$,

$$v_n := \frac{u_n}{\sum_{k=1}^{+\infty} u_k}. \quad (5.50)$$

We think of this process as giving the weights of a Gibbs measure supported on the natural numbers, and we write $\langle \cdot \rangle$ for its associated average, with $\alpha$ denoting the canonical random variable on this space. We also denote by $(\alpha^\ell)_{\ell \geq 1}$ independent copies, or replicas, of the random variable $\alpha$ under this measure. This means that for any $k \geq 1$ and any bounded function $f : \mathbb{N}^k \to \mathbb{R}$,

$$\langle f(\alpha^1, \ldots, \alpha^k) \rangle := \sum_{j_1, \ldots, j_k \geq 1} f(j_1, \ldots, j_k) v_{j_1} \cdots v_{j_k}. \tag{5.51}$$

We now prove a fundamental invariance property of the ordered Poisson-Dirichlet point process that will allow us to explicitly compute free energy functionals associated with the Gibbs average (5.51). This invariance property will also lead to deep distributional identities that characterize the Gibbs measures whose ordered weights are given by a Poisson-Dirichlet process.

**Theorem 5.19.** *Let $(X_n, Y_n)_{n \geq 1}$ be i.i.d. pairs of random variables on $\mathbb{R}_{>0} \times \mathbb{R}^d$ independent of the ordered Poisson-Dirichlet point process $(u_n)_{n \geq 1}$, with the property that $\mathbb{E} X_1^\zeta < +\infty$. Let $\nu_\zeta$ be the probability measure on $\mathbb{R}^d$ such that*

$$\nu_\zeta(B) := \frac{\mathbb{E} X_1^\zeta \mathbf{1}_{\{Y_1 \in B\}}}{\mathbb{E} X_1^\zeta}, \tag{5.52}$$

*and let $(Y_n')_{n \geq 1}$ be a sequence of i.i.d. random variables with law $\nu_\zeta$ independent of $(u_n)_{n \geq 1}$. The Poisson point processes*

$$\sum_{n=1}^{+\infty} \delta_{(u_n X_n, Y_n)} \quad \text{and} \quad \sum_{n=1}^{+\infty} \delta_{((\mathbb{E} X_1^\zeta)^{1/\zeta} u_n, Y_n')} \tag{5.53}$$

*have the same law.*

*Proof.* Let $\mu$ be as in (5.46), and let $\nu$ be the law of $(X_1, Y_1)$ on $\mathbb{R}_{>0} \times \mathbb{R}^d$. Since $(X_n, Y_n)_{n \geq 1}$ are i.i.d. and independent of $(u_n)_{n \geq 1}$, the marking theorem (Proposition 5.11) and the mapping theorem (Proposition 5.10) applied with the function $f(u, x, y) := (ux, y)$ imply that

$$\sum_{n=1}^{+\infty} \delta_{(u_n X_n, Y_n)}$$

is a Poisson point process on $\mathbb{R}_{>0} \times \mathbb{R}^d$ with intensity measure $(\mu \otimes \nu) \circ f^{-1}$. Given two measurable sets $A \subseteq \mathbb{R}_{>0}$ and $B \subseteq \mathbb{R}^d$, Fubini's theorem and a change of variables

reveal that

$$(\mu \otimes \nu) \circ f^{-1}(A \times B) = \int_{\mathbb{R}_{>0} \times \mathbb{R}^d} \int_0^{+\infty} \mathbf{1}_{\{ux \in A\}} \mathbf{1}_{\{y \in B\}} \zeta u^{-1-\zeta} \, du \, d\nu(x,y)$$

$$= \int_{\mathbb{R}_{>0} \times \mathbb{R}^d} x^\zeta \mathbf{1}_{\{y \in B\}} \int_0^{+\infty} \mathbf{1}_{\{v \in A\}} \zeta v^{-1-\zeta} \, dv \, d\nu(x,y)$$

$$= \mu(A) \int_{\mathbb{R}_{>0} \times \mathbb{R}^d} x^\zeta \mathbf{1}_{\{y \in B\}} \, d\nu(x,y)$$

$$= \mathbb{E} X_1^\zeta (\mu \otimes \nu_\zeta)(A \times B).$$

Similarly, the marking theorem and the mapping theorem applied with the function $g(u) := (\mathbb{E} X_1^\zeta)^{1/\zeta} u$ imply that

$$\sum_{n=1}^{+\infty} \delta_{((\mathbb{E} X_1^\zeta)^{1/\zeta} u_n, Y_n')}$$

is a Poisson point process with intensity measure $(\mu \circ g^{-1}) \otimes \nu_\zeta$. Given a measurable set $A \subseteq \mathbb{R}_{>0}$, a change of variables reveals that

$$\mu \circ g^{-1}(A) = \int_0^{+\infty} \mathbf{1}_{\{(\mathbb{E} X_1^\zeta)^{1/\zeta} u \in A\}} \zeta u^{-1-\zeta} \, du = \mathbb{E} X_1^\zeta \int_0^{+\infty} \mathbf{1}_{\{v \in A\}} \zeta v^{-1-\zeta} \, dv$$

$$= \mathbb{E} X_1^\zeta \mu(A).$$

Since the Poisson point processes in (5.53) have the same intensity measure, they have the same law by Proposition 5.9. This completes the proof. ∎

**Proposition 5.20.** *Let $(X_n)_{n \geq 1}$ be i.i.d. random variables on $\mathbb{R}_{>0}$ that are independent of the ordered Poisson-Dirichlet point process $(u_n)_{n \geq 1}$ with parameter $\zeta \in (0,1)$. If $\mathbb{E} X_1^\zeta < +\infty$, then*

$$\mathbb{E} \log \langle X_\alpha \rangle = \frac{1}{\zeta} \log \mathbb{E} X_1^\zeta. \tag{5.54}$$

*Proof.* Theorem 5.19 implies that

$$\mathbb{E} \log \sum_{n=1}^{+\infty} u_n X_n = \mathbb{E} \log \sum_{n=1}^{+\infty} (\mathbb{E} X_1^\zeta)^{1/\zeta} u_n = \frac{1}{\zeta} \log \mathbb{E} X_1^\zeta + \mathbb{E} \log \sum_{n=1}^{+\infty} u_n. \tag{5.55}$$

Remembering the definition of the Poisson-Dirichlet process in (5.50) and of the Gibbs measure in (5.51) completes the proof. ∎

This result is slightly generalized in Exercise 5.9. Its counterpart for the Poisson-Dirichlet cascades discussed in the next section will be used to determine the initial condition for the infinite-dimensional Hamilton-Jacobi equation arising in the context of the spin-glass models studied in the next chapter. Another consequence of the invariance property in Theorem 5.19 is the Bolthausen-Sznitman invariance for the ordered Poisson-Dirichlet point process.

**Lemma 5.21** (Bolthausen-Sznitman invariance for PDP [57]). *Let $(g_n)_{n \geq 1}$ be a sequence of independent centred Gaussian random variables with variance $\mathsf{v} > 0$, independent of the ordered Poisson-Dirichlet point process $(u_n)_{n \geq 1}$ with parameter $\zeta \in (0,1)$. For every $t \in \mathbb{R}$, the Poisson point processes*

$$\sum_{n=1}^{+\infty} \delta_{(u_n \exp(t(g_n - t\mathsf{v}\zeta/2)), g_n - t\mathsf{v}\zeta)} \quad \text{and} \quad \sum_{n=1}^{+\infty} \delta_{(u_n, g_n)}$$

*have the same law.*

In the context of mean-field spin glasses, we think of each term $u_n$ as representing the asymptotic value of the exponential of the energy of a particular configuration; and since this energy is a Gaussian random variable, it is relatively natural to perturb $u_n$ by multiplying it with the exponential of a Gaussian random variable as done in Lemma 5.21.

*Proof of Lemma 5.21.* Consider the sequence of random variables $(X_n, Y_n)_{n \geq 1}$ with values in $\mathbb{R}_{>0} \times \mathbb{R}$ defined by

$$(X_n, Y_n) := \big(\exp(t(g_n - t\mathsf{v}\zeta/2)), g_n - t\mathsf{v}\zeta\big).$$

Using the formula for the moment generating function of a standard Gaussian random variable reveals that

$$\mathbb{E} X_1^\zeta = \exp(-t^2 \zeta^2 \mathsf{v}/2) \mathbb{E} \exp(t\zeta g_1) = 1.$$

It follows by a change of variables that the measure $\nu_\zeta$ defined by (5.52) is given by

$$\begin{aligned}
\nu_\zeta(B) &= \frac{1}{\sqrt{2\pi\mathsf{v}}} \int_\mathbb{R} e^{t\zeta(x - t\mathsf{v}\zeta/2)} \mathbf{1}_{\{x - t\mathsf{v}\zeta \in B\}} e^{-\frac{x^2}{2\mathsf{v}}} \, \mathrm{d}x \\
&= \frac{1}{\sqrt{2\pi\mathsf{v}}} \int_\mathbb{R} \mathbf{1}_{\{x - t\mathsf{v}\zeta \in B\}} e^{-\frac{(x - t\zeta\mathsf{v})^2}{2\mathsf{v}}} \, \mathrm{d}x \\
&= \mathbb{P}\{g_1 \in B\}.
\end{aligned}$$

This means that $\nu_\zeta$ is the Gaussian distribution on $\mathbb{R}$ with variance $\mathsf{v}$. Invoking Theorem 5.19 completes the proof. ∎

The Bolthausen-Sznitman invariance for the ordered Poisson-Dirichlet point process can be used to establish the *Ghirlanda-Guerra identities* for the Gibbs measure (5.51). These identities were introduced in [128] (see also [11, 130]) in the context of mean-field spin glasses. It turns out that these identities characterize the law of the weights of this Gibbs measure, as shown in Exercise 5.10. To state them, it will be convenient to introduce the "overlap"

$$R_{\ell,\ell'} := \mathbf{1}_{\{\alpha^\ell = \alpha^{\ell'}\}} \tag{5.56}$$

between two independent samples $\alpha^\ell \in \mathbb{N}$ and $\alpha^{\ell'} \in \mathbb{N}$ from the Gibbs measure (5.51). At this stage it might sound a bit strange to call the indicator function in (5.56) an overlap, but the choice of terminology will become gradually clear. For now, we note that one can represent $\alpha^\ell$ and $\alpha^{\ell'}$ as elements of a Hilbert space and interpret the overlap $R_{\ell,\ell'}$ as a scalar product between these elements. To do so, we can choose the Hilbert space to be the space of square-integrable sequences, denote by $(e_n)_{n \in \mathbb{N}}$ the canonical basis in this space, and identify each integer $n$ with the vector $e_n$, so that $R_{\ell,\ell'} = e_{\alpha^\ell} \cdot e_{\alpha^{\ell'}}$. The Ghirlanda-Guerra identities specify the law of the overlap $R_{1,n+1}$ given the overlap array $R^n := (R_{\ell,\ell'})_{\ell,\ell' \leq n}$.

**Theorem 5.22** (Ghirlanda-Guerra identities for PDP). *Let $\langle \cdot \rangle$ denote the Gibbs measure (5.51) associated with the Poisson-Dirichlet process $(v_n)_{n \geq 1}$ with parameter $\zeta \in (0,1)$. For every $n \geq 1$ and every bounded and measurable function $f$ of the overlaps $R^n := (R_{\ell,\ell'})_{\ell,\ell' \leq n}$ as defined in (5.56), we have*

$$\mathbb{E}\langle f(R^n) R_{1,n+1}\rangle = \frac{1}{n}\mathbb{E}\langle f(R^n)\rangle \mathbb{E}\langle R_{1,2}\rangle + \frac{1}{n}\sum_{\ell=2}^{n}\mathbb{E}\langle f(R^n)R_{1,\ell}\rangle, \tag{5.57}$$

*and moreover,*

$$\mathbb{E}\langle R_{1,2}\rangle = \mathbb{E}\langle \mathbf{1}_{\{\alpha^1 = \alpha^2\}}\rangle = 1 - \zeta. \tag{5.58}$$

*Proof.* Let $(g_n)_{n \geq 1}$ be a sequence of independent standard Gaussian random variables, and for each $t \in \mathbb{R}$, define the random weights

$$v_n^t := \frac{u_n \exp t\left(g_n - \frac{\zeta t}{2}\right)}{\sum_{k=1}^{+\infty} u_k \exp t\left(g_k - \frac{\zeta t}{2}\right)} = \frac{u_n \exp(tg_n)}{\sum_{k=1}^{+\infty} u_k \exp(tg_k)}.$$

We denote by $\langle \cdot \rangle_t$ the average with respect to the random measure on $\mathbb{N}$ with weights $(v_n^t)_{n \geq 1}$. By the Bolthausen-Sznitman invariance (Lemma 5.21), the Poisson point processes

$$\sum_{n=1}^{+\infty}\delta_{(v_n^t, g_n - t\zeta)} \quad \text{and} \quad \sum_{n=1}^{+\infty}\delta_{(v_n, g_n)} \tag{5.59}$$

have the same law. This implies that for any bounded measurable function $f$ of the overlaps $R^n$,

$$0 = \mathbb{E}\langle f(R^n) g_{\alpha^1}\rangle = \mathbb{E}\langle f(R^n)(g_{\alpha^1} - t\zeta)\rangle_t.$$

The first equality is obtained by averaging with respect to the Gaussian randomness first. Invoking the Gibbs Gaussian integration by parts formula (Theorem 4.6) reveals that

$$t\mathbb{E}\Big\langle f(R^n)\Big(\sum_{\ell=1}^n \mathbf{1}_{\{\alpha^1=\alpha^\ell\}} - n\mathbf{1}_{\{\alpha^1=\alpha^{n+1}\}} - \zeta\Big)\Big\rangle_t = 0.$$

Using once again that the point processes in (5.59) have the same law, recalling (5.56) and rearranging yields

$$\mathbb{E}\langle f(R^n)R_{1,n+1}\rangle = \frac{1}{n}\mathbb{E}\Big\langle f(R^n)\Big(\sum_{\ell=1}^n R_{1,\ell} - \zeta\Big)\Big\rangle = \frac{1}{n}\sum_{\ell=2}^n \mathbb{E}\langle f(R^n)R_{1,\ell}\rangle + \frac{1-\zeta}{n}\mathbb{E}\langle f(R^n)\rangle.$$

Taking $f = n = 1$ shows that $1 - \zeta = \mathbb{E}\langle R_{1,2}\rangle$, and completes the proof. ∎

Denoting by $\nu$ the law of a Bernoulli random variable with probability of success $1 - \zeta$, we deduce from Theorem 5.22 that conditionally on $R^n$, the distribution of the overlap $R_{1,n+1}$ is given by the mixture

$$\frac{1}{n}\nu + \frac{1}{n}\sum_{\ell=2}^n \delta_{R_{1,\ell}}. \tag{5.60}$$

This means that with probability $\frac{1}{n}$, the overlap between $\alpha^{n+1}$ and $\alpha^1$ is sampled independently according to $\nu$, and with probability $\frac{n-1}{n}$ it takes one of the existing values from $R_{1,2},\ldots,R_{1,n}$ uniformly at random. It turns out that one can in fact infer the entire law of the weights of the Gibbs measure (5.51) from the Ghirlanda-Guerra identities, not just this conditional marginal distribution. This point is the subject of Exercise 5.10.

**Exercise 5.9.** Under the setting of Proposition 5.20, let $(Z_n)_{n\geq 1}$ be i.i.d. random variables on $\mathbb{R}_{>0}$, independent of $(X_n)_{n\geq 1}$ and of the ordered Poisson-Dirichlet point process $(u_n)_{n\geq 1}$ with parameter $\zeta \in (0,1)$. Assume that $\mathbb{E}X_1^\zeta < +\infty$ and $\mathbb{E}Z_1^\zeta < +\infty$. Show that

$$\mathbb{E}\log\sum_{n=1}^{+\infty} u_n X_n Z_n = \mathbb{E}\log\sum_{n=1}^{+\infty} u_n Z_n + \frac{1}{\zeta}\log\mathbb{E}X_1^\zeta. \tag{5.61}$$

**Exercise 5.10.** Let $(v_j)_{j\geq 1}$ denote the decreasing enumeration of the weights of a random probability measure $\langle\cdot\rangle$ on the natural numbers. We define the overlaps as in (5.56), suppose that the Ghirlanda-Guerra identities (5.57) hold, and define $\zeta \in (0,1)$ so that (5.58) holds. The main purpose of this exercise is to show that the law of $(v_j)_{j\geq 1}$ must be the law of the Poisson-Dirichlet process with parameter $\zeta$.

(i) Consider the closed set

$$V := \Big\{x = (x_j)_{j\geq 1} \mid \sum_{j=1}^{+\infty} x_j \leq 1 \text{ and } (x_j)_{j\geq 1} \text{ is non-increasing}\Big\} \tag{5.62}$$

of the product space $[0,1]^{\mathbb{N}}$ equipped with the product topology, and for each integer $n \geq 1$, define the function $p_n : V \to \mathbb{R}$ by $p_n(x) := \sum_{j \geq 1} x_j^n$. Prove that $p_n$ is continuous for $n > 1$. Is $p_1$ continuous?

(ii) Prove that the set $\mathcal{P}$ of finite linear combinations of constants and products $p_{n_1} \cdots p_{n_k}$ for $n_1, \ldots, n_k \geq 2$ is dense in the space $C(V; \mathbb{R})$ of continuous functions from $V$ to $\mathbb{R}$ endowed with the uniform norm.

(iii) For each $k \geq 1$, $n_1, \ldots, n_k \geq 2$, and $x \in V$, define $f(x; n_1, \ldots, n_k) := \prod_{\ell=1}^{k} p_{n_\ell}(x)$, let $n := n_1 + \cdots + n_k$, and set $S(n_1, \ldots, n_k) := \mathbb{E} f\big((v_j)_{j \geq 1}, n_1, \ldots, n_k\big)$. Establish *Talagrand's identities*,

$$S(n_1+1, n_2, \ldots, n_k) = \frac{n_1 - \zeta}{n} S(n_1, \ldots, n_k) + \sum_{j=2}^{k} \frac{n_j}{n} S(n_2, \ldots, n_j + n_1, \ldots, n_k).$$

(iv) Deduce that the law of $(v_j)_{j \geq 1}$ is entirely determined by the parameter $\zeta$, and conclude that the random probability measure $\langle \cdot \rangle$ must be the Gibbs measure (5.51) associated with the Poisson-Dirichlet process with parameter $\zeta$.

## 5.6 Poisson-Dirichlet cascades

Using Poisson-Dirichlet point processes, we now iteratively construct Poisson-Dirichlet cascades (also called Ruelle probability cascades [234]). As was discussed at the opening of the previous section, these richer objects will serve as "canonical models" for describing the geometry of the Gibbs measure of mean-field spin glasses, in the limit of large system size.

For each fixed choice of integer $K \geq 1$, we build Poisson-Dirichlet cascades as random probability measures on the index set $\mathbb{N}^K$. The case $K = 1$ will correspond to the Poisson-Dirichlet process. We think of $\mathbb{N}^K$ as the set of leaves of a rooted tree with vertex set

$$\mathcal{A} := \mathbb{N}^0 \cup \mathbb{N} \cup \mathbb{N}^2 \cup \ldots \cup \mathbb{N}^K, \tag{5.63}$$

where $\mathbb{N}^0 := \{\varnothing\}$. In this identification, the root of the tree is the element $\varnothing$, and each vertex $\alpha = (n_1, \ldots, n_k) \in \mathbb{N}^k$ for $k \leq K - 1$ has children

$$\alpha n := (n_1, \ldots, n_k, n) \in \mathbb{N}^{k+1} \tag{5.64}$$

for all $n \in \mathbb{N}$. For each vertex $\alpha = (n_1, \ldots, n_k) \in \mathcal{A}$, we write $|\alpha| = k$ for its distance from the root $\varnothing$, or equivalently for the number of coordinates in $\alpha$, and for each $0 \leq \ell \leq k$ we denote by

$$\alpha_{|\ell} := (n_1, \ldots, n_\ell) \tag{5.65}$$

the vertex on the path from the root $\varnothing$ to $\alpha$ at distance $\ell$ from the root. We understand that $\alpha_{|0} = \varnothing$. For any two leaves $\alpha, \beta \in \mathbb{N}^K$ we write

$$\alpha \wedge \beta := \max\{\ell \leq K \mid \alpha_{|\ell} = \beta_{|\ell}\} \tag{5.66}$$

for the number of common vertices in the paths from the root $\varnothing$ to the leaves $\alpha$ and $\beta$.

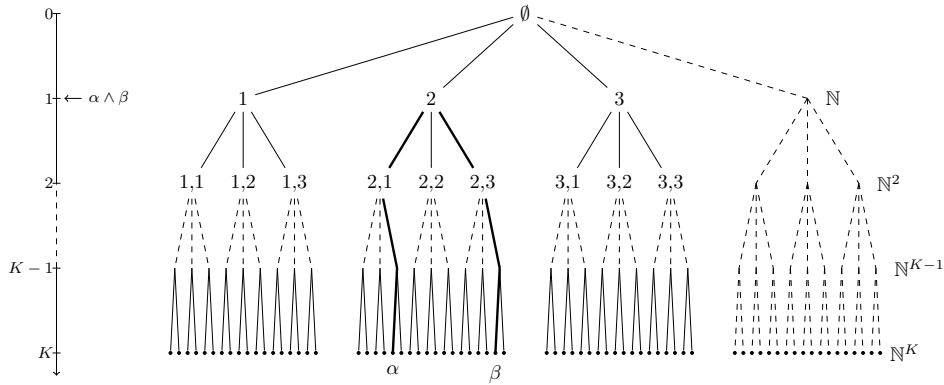

**Figure 5.1** The ancestry lines of the leaves $\alpha$ and $\beta$ meet at the first level of this regular tree of infinite degree and depth $K$.

Although the Poisson-Dirichlet cascade will be indexed by $\alpha \in \mathbb{N}^K$, its construction will involve random variables indexed by vertices of the entire tree $\mathcal{A}$. Given a sequence of parameters

$$0 = \zeta_0 < \zeta_1 < \ldots < \zeta_K < \zeta_{K+1} = 1, \tag{5.67}$$

for each non-leaf vertex $\alpha \in \mathcal{A} \setminus \mathbb{N}^K$, let $\Lambda_\alpha$ be the Poisson-Dirichlet point process with parameter $\zeta_{|\alpha|+1}$; we construct these point processes independently from one another as $\alpha$ varies. We denote by $(u_{\alpha n})_{n \geq 1}$ the ordered Poisson-Dirichlet point process corresponding to $\Lambda_\alpha$. In this way, parent vertices $\alpha \in \mathcal{A} \setminus \mathbb{N}^K$ enumerate independent Poisson point processes $\Lambda_\alpha$, while child vertices $\alpha n \in \mathcal{A} \setminus \mathbb{N}^0$ enumerate individual points $u_{\alpha n}$. To each node $\alpha \in \mathcal{A}$, we associate the quantity

$$w_\alpha := \prod_{1 \leq \ell \leq K} u_{\alpha_{|\ell}}. \tag{5.68}$$

We would now like to define the Poisson-Dirichlet cascade as the random measure on $\mathbb{N}^K$ with weights

$$v_\alpha := \frac{w_\alpha}{\sum_{\beta \in \mathbb{N}^K} w_\beta}. \tag{5.69}$$

This is possible due to an analogue of Proposition 5.18.

**Lemma 5.23.** *The sum $\sum_{\alpha \in \mathbb{N}^K} w_\alpha$ is finite with probability one.*

*Proof.* The proof proceeds by induction on $K$. The case $K = 1$ is the content of Proposition 5.18, so let us assume that the result holds when we sum over $K - 1 \geqslant 1$ levels, and prove it for sums over $K$ levels. Notice that for every $\alpha \in \mathbb{N}^{K-1}$ and $n \in \mathbb{N}$, we have $w_{\alpha n} = w_\alpha u_{\alpha n}$, and therefore

$$\sum_{\alpha \in \mathbb{N}^K} w_\alpha = \sum_{\alpha \in \mathbb{N}^{K-1}} \sum_{n \in \mathbb{N}} w_{\alpha n} = \sum_{\alpha \in \mathbb{N}^{K-1}} w_\alpha \sum_{n \in \mathbb{N}} u_{\alpha n}.$$

For each $\alpha \in \mathcal{A} \smallsetminus \mathbb{N}^K$, let $U_\alpha := \sum_{n \in \mathbb{N}} u_{\alpha n}$, so that

$$\sum_{\alpha \in \mathbb{N}^K} w_\alpha = \sum_{\alpha \in \mathbb{N}^{K-1}} w_\alpha U_\alpha = \sum_{\alpha \in \mathbb{N}^{K-2}} \sum_{n \in \mathbb{N}} w_{\alpha n} U_{\alpha n} = \sum_{\alpha \in \mathbb{N}^{K-2}} w_\alpha \sum_{n \in \mathbb{N}} u_{\alpha n} U_{\alpha n}.$$

Observe that for each $\alpha \in \mathbb{N}^{K-2}$, the random variables $(U_{\alpha n})_{n \geqslant 1}$ are i.i.d. since the Poisson point processes $(\Lambda_{\alpha n})_{n \geqslant 1}$ are. Moreover, as $\zeta_{K-2} < \zeta_{K-1}$, Proposition 5.18 implies that

$$\mathbb{E}(U_{\alpha n})^{\zeta_{K-2}} = \mathbb{E}\left( \int_0^{+\infty} x \, d\Lambda_{\alpha n}(x) \right)^{\zeta_{K-2}} < +\infty.$$

It follows by Theorem 5.19 that $(u_{\alpha n} U_{\alpha n})_{n \geqslant 1}$ and $(c u_{\alpha n})_{n \geqslant 1}$ for $c := (\mathbb{E} U_{\alpha 1}^{\zeta_{K-2}})^{1/\zeta_{K-2}}$ are Poisson point processes with the same intensity measure. In particular,

$$\sum_{n \in \mathbb{N}} u_{\alpha n} U_{\alpha n} \stackrel{d}{=} c \sum_{n \in \mathbb{N}} u_{\alpha n} = c U_\alpha, \qquad (5.70)$$

where $\stackrel{d}{=}$ denotes equality in distribution. Since the quantities on the right side of (5.70) are i.i.d. over $\alpha \in \mathbb{N}^{K-2}$ and independent of $(w_\alpha)_{\alpha \in \mathbb{N}^{K-2}}$, it follows that $\sum_{\alpha \in \mathbb{N}^K} w_\alpha$ is finite with probability one if and only if

$$\sum_{\alpha \in \mathbb{N}^{K-2}} w_\alpha U_\alpha = \sum_{\alpha \in \mathbb{N}^{K-2}} \sum_{n \in \mathbb{N}} w_\alpha u_{\alpha n} = \sum_{\alpha \in \mathbb{N}^{K-1}} w_\alpha$$

is. Invoking the induction hypothesis completes the proof. ■

**Definition 5.24** (Poisson-Dirichlet cascade). The *Poisson-Dirichlet cascade* with parameters (5.67) is the family $(v_\alpha)_{\alpha \in \mathbb{N}^K}$ obtained by normalizing the $(w_\alpha)_{\alpha \in \mathbb{N}^K}$, that is, for every $\alpha \in \mathbb{N}^K$,

$$v_\alpha := \frac{w_\alpha}{\sum_{\beta \in \mathbb{N}^K} w_\beta}. \qquad (5.71)$$

We stress that the Poisson-Dirichlet cascade is defined as the *family* of weights $(v_\alpha)_{\alpha \in \mathbb{N}^K}$, as opposed to the *set* of these weights. While the precise labelling of the weights will not be important, we do care to keep track of the underlying tree structure. In other words, we are interested in the knowledge of $(v_\alpha)_{\alpha \in \mathbb{N}^K}$ up to the action of the bijections of $\mathbb{N}^K$ that preserve the tree structure.

We think of this cascade as giving the weights of a Gibbs measure supported on the set $\mathbb{N}^K$, and we write $\langle \cdot \rangle$ for its associated average, with $\alpha$ denoting the canonical random variable. We write $(\alpha^\ell)_{\ell \geq 1}$ to denote independent copies, or replicas, of the random variable $\alpha$ under this measure. This means that for any $k \geq 1$ and any bounded function $f : (\mathbb{N}^K)^k \to \mathbb{R}$,

$$\langle f(\alpha^1, \ldots, \alpha^k) \rangle := \sum_{\alpha_1, \ldots, \alpha_k \in \mathbb{N}^K} f(\alpha_1, \ldots, \alpha_k) v_{\alpha_1} \cdots v_{\alpha_k}. \tag{5.72}$$

We now generalize the explicit formula obtained in Proposition 5.20 and the Ghirlanda-Guerra identities in Theorem 5.22 from the setting of the Poisson-Dirichlet processes to that of Poisson-Dirichlet cascades. It will be convenient to introduce a sequence $(\omega_k)_{0 \leq k \leq K}$ of i.i.d. uniform random variables on $[0,1]$. These uniform random variables can be used to generate any random variable taking values in a complete separable metric space, as explained in Exercise 5.11 in the simple case of real-valued random variables. We give ourselves a measurable function

$$X_K := X_K(\omega_0, \ldots, \omega_K) \tag{5.73}$$

from $[0,1]^K$ to $\mathbb{R}$. Assuming that $\mathbb{E} \exp \zeta_K X_K < +\infty$, we now show that certain functionals associated with the Gibbs measure (5.72) can be computed recursively by considering the sequence

$$X_k := X_k(\omega_0, \ldots, \omega_k) := \frac{1}{\zeta_{k+1}} \log \mathbb{E}_{\omega_{k+1}} \exp \zeta_{k+1} X_{k+1} \tag{5.74}$$

defined for $0 \leq k \leq K-1$. By continuity, see Exercise 5.12, we take $X_{-1} := \mathbb{E}_{\omega_0} X_0$. These functionals will come up in the definition of enriched free energies in the next chapter.

**Theorem 5.25.** *Let $(\omega_\beta)_{\beta \in \mathcal{A}}$ be a family of independent uniform random variables on $[0,1]$, independent of the Poisson-Dirichlet cascade $(v_\alpha)_{\alpha \in \mathbb{N}^K}$, and let $X_{-1}$ be defined recursively via (5.74). We have*

$$\mathbb{E} \log \langle \exp X_K(\omega_\varnothing, \omega_{\alpha_{|1}}, \ldots, \omega_{\alpha_{|K-1}}, \omega_\alpha) \rangle = X_{-1}. \tag{5.75}$$

In order to avoid ambiguities, we emphasize that on the left side of (5.75), the random variable $\alpha$ takes values in $\mathbb{N}^K$, is sampled according to the Gibbs measure defined in (5.72), and that $\varnothing, \alpha_{|1}, \ldots, \alpha_{|K-1}, \alpha$ denote the sequence of ancestors of $\alpha$, starting from the root, as defined in (5.65). The expectation $\langle \cdot \rangle$ only averages over the random variable $\alpha$, while the expectation $\mathbb{E}$ averages the weights (5.71) of the Gibbs measure as well as the random variables $(\omega_\beta)_{\beta \in \mathcal{A}}$. On the other hand, the recursive calculation defining $X_{-1}$ only involves much simpler recursive averages over the random variables $\omega_K, \ldots, \omega_0$.

*Proof of Theorem 5.25.* First, we observe that by Jensen's inequality and a simple induction, for $-1 \leqslant k \leqslant K$,

$$\mathbb{E}\exp \zeta_K X_k = \mathbb{E}(\mathbb{E}_{\omega_{k+1}} \exp \zeta_{k+1} X_{k+1})^{\zeta_K/\zeta_{k+1}} \leqslant \mathbb{E}\exp \zeta_K X_{k+1} \leqslant \mathbb{E}\exp \zeta_K X_K < +\infty,$$

so the sequence $(X_k)_{-1 \leqslant k \leqslant K}$ is well-defined. To simplify notation, for each leaf $\alpha \in \mathbb{N}^K$, write $\Omega_\alpha := (\omega_\varnothing, \omega_{\alpha_{|1}}, \ldots, \omega_{\alpha_{|K-1}}, \omega_\alpha)$ for the sequence of uniform random variables along the path joining $\alpha$ to the root of the tree. In this notation, equality (5.75) may be written as

$$\mathbb{E}\log \sum_{\alpha \in \mathbb{N}^K} w_\alpha \exp X_K(\Omega_\alpha) = \mathbb{E}\log \sum_{\alpha \in \mathbb{N}^K} w_\alpha + X_{-1}.$$

It turns out to be more convenient to establish a slightly more general equality. Let $Z > 0$ be a random variable with $\mathbb{E}Z^{\zeta_K} < +\infty$ and let $(Z_\alpha)_{\alpha \in \mathbb{N}^K}$ be a sequence of i.i.d. copies of $Z$ independent of all other random variables. We will prove that

$$\mathbb{E}\log \sum_{\alpha \in \mathbb{N}^K} w_\alpha Z_\alpha \exp X_K(\Omega_\alpha) = \mathbb{E}\log \sum_{\alpha \in \mathbb{N}^K} w_\alpha Z_\alpha + X_{-1}, \qquad (5.76)$$

and that both sides are well-defined by induction on $K$. The case $K = 1$ is a consequence of Proposition 5.20 in the form (5.55), or more directly, of the rephrasing given in Exercise 5.9. Indeed, conditionally on the randomness of $\omega_\varnothing$,

$$\mathbb{E}\log \sum_{n \in \mathbb{N}} u_n Z_n \exp X_1(\omega_\varnothing, \omega_n) = \mathbb{E}\log \sum_{n \in \mathbb{N}} u_n Z_n + \frac{1}{\zeta_1}\log \mathbb{E}\exp X_1(\omega_\varnothing, \omega_1)$$

$$= \mathbb{E}\log \sum_{\alpha \in \mathbb{N}} w_\alpha Z_\alpha + X_0.$$

Averaging with respect to the randomness of $\omega_\varnothing$ establishes (5.76) for $K = 1$. Let us therefore suppose that the result holds for $K - 1 \geqslant 1$ and prove it for $K$. For each $\alpha \in \mathbb{N}^{K-1}$ we write $\Omega_{\alpha n} := (\Omega_\alpha, \omega_{\alpha n})$, so that

$$\sum_{\alpha \in \mathbb{N}^K} w_\alpha Z_\alpha \exp X_K(\Omega_\alpha) = \sum_{\alpha \in \mathbb{N}^{K-1}} \sum_{n \in \mathbb{N}} w_{\alpha n} Z_{\alpha n} \exp X_K(\Omega_{\alpha n})$$

$$= \sum_{\alpha \in \mathbb{N}^{K-1}} w_\alpha \sum_{n \in \mathbb{N}} u_{\alpha n} Z_{\alpha n} \exp X_K(\Omega_\alpha, \omega_{\alpha n}). \qquad (5.77)$$

Denote by $\mathcal{F}_{K-1}$ the $\sigma$-algebra generated by $(w_\alpha)_{\alpha \in \mathbb{N}^{K-1}}$ and $(\Omega_\alpha)_{\alpha \in \mathbb{N}^{K-1}}$, and observe that for each $\alpha \in \mathbb{N}^{K-1}$, the random variables $(u_{\alpha n})_{n \geqslant 1}$, $(Z_{\alpha n})_{n \geqslant 1}$ and $(\omega_{\alpha n})_{n \geqslant 1}$ are independent of $\mathcal{F}_{K-1}$. It follows that

$$\mathbb{E}[(Z_{\alpha n}\exp X_K(\Omega_\alpha, \omega_{\alpha n}))^{\zeta_K} \mid \mathcal{F}_{K-1}]^{1/\zeta_K} = c(\mathbb{E}_{\omega_{\alpha n}}\exp \zeta_K X_K(\Omega_\alpha, \omega_{\alpha n}))^{1/\zeta_K}$$

$$= c\exp X_{K-1}(\Omega_\alpha)$$

for $c := (\mathbb{E}Z^{\zeta_K})^{1/\zeta_K}$. Theorem 5.19 now implies that conditionally on $\mathcal{F}_{K-1}$,

$$\sum_{n\in\mathbb{N}} u_{\alpha n} Z_{\alpha n} \exp X_K(\Omega_\alpha, \omega_{\alpha n}) \stackrel{d}{=} c U_\alpha \exp X_{K-1}(\Omega_\alpha)$$

for $U_\alpha := \sum_{n\in\mathbb{N}} u_{\alpha n}$, and where we recall that $\stackrel{d}{=}$ denotes equality in distribution. Conditionally on $\mathcal{F}_{K-1}$, both sides of this distributional equality are independent over $\alpha \in \mathbb{N}^{K-1}$, so substituting into (5.77) reveals that conditionally on $\mathcal{F}_{K-1}$,

$$\sum_{\alpha\in\mathbb{N}^K} w_\alpha Z_\alpha \exp X_K(\Omega_\alpha) \stackrel{d}{=} c \sum_{\alpha\in\mathbb{N}^{K-1}} w_\alpha U_\alpha \exp X_{K-1}(\Omega_\alpha).$$

Averaging over the randomness of $\mathcal{F}_{K-1}$ gives this equality in distribution unconditionally. The particular choice of $X_K = 0$ yields the distributional equality

$$\sum_{\alpha\in\mathbb{N}^K} w_\alpha Z_\alpha \stackrel{d}{=} c \sum_{\alpha\in\mathbb{N}^{K-1}} w_\alpha U_\alpha,$$

so proving (5.76) is reduced to establishing the equality

$$\mathbb{E}\log \sum_{\alpha\in\mathbb{N}^{K-1}} w_\alpha U_\alpha \exp X_{K-1}(\Omega_\alpha) = \mathbb{E}\log \sum_{\alpha\in\mathbb{N}^{K-1}} w_\alpha U_\alpha + X_{-1}. \tag{5.78}$$

By Proposition 5.18, the random variables $(U_\alpha)_{\alpha\in\mathbb{N}^{K-1}}$ satisfy

$$\mathbb{E}U_\alpha^{\zeta_{K-1}} = \mathbb{E}\left(\int x\,d\Lambda_\alpha(x)\right)^{\zeta_{K-1}} < +\infty$$

and are i.i.d. and independent of all other random variables in (5.78). Since the equation (5.78) is of the same type as (5.76) with $K$ replaced by $K-1$ and $(Z_\alpha)_{\alpha\in\mathbb{N}^K}$ replaced by $(U_\alpha)_{\alpha\in\mathbb{N}^{K-1}}$, invoking the induction hypothesis completes the proof. ∎

In Chapter 6, we will use this result to compute the initial condition for the Hamilton-Jacobi equation associated with spin-glass models. For this purpose, it will be convenient to establish a slightly more general version. Instead of just considering $K+1$ sources of randomness $(\omega_k)_{0\leq k\leq K}$, we generate $N(K+1)$ independent random variables $(\omega_{k,i})_{0\leq k\leq K, 1\leq i\leq N}$ sampled from the uniform distribution on $[0,1]$, and instead of just considering one function $X_K$, we give ourselves $N$ measurable functions

$$X_{K,i} := X_{K,i}(\omega_{0,i},\ldots,\omega_{K,i}) \tag{5.79}$$

for $i \in \{1,\ldots,N\}$. Assuming that $\mathbb{E}\exp \zeta_K X_{K,i} < +\infty$, we recursively define, for every $k \in \{-1,\ldots,K-1\}$ and $i \in \{1,\ldots,N\}$,

$$X_{k,i} := X_{k,i}(\omega_{0,i},\ldots,\omega_{k,i}) := \frac{1}{\zeta_{k+1}} \log \mathbb{E}_{\omega_{k+1,i}} \exp \zeta_{k+1} X_{k+1,i}. \tag{5.80}$$

Once again, by continuity, we understand that $X_{-1,i} := \mathbb{E}_{\omega_{0,i}} X_{0,i}$. The next result states that this more general iterative procedure can be used to calculate more complicated averages with respect to the Poisson-Dirichlet cascade $\langle\cdot\rangle$.

**Corollary 5.26.** *Let $N \geq 1$ be an integer, let $(\omega_{\beta,i})_{\beta \in \mathcal{A}, 1 \leq i \leq N}$ be a family of independent uniform random variables on $[0,1]$ independent of the Poisson-Dirichlet cascade $(v_\alpha)_{\alpha \in \mathbb{N}^K}$, and for each $i \in \{1,\ldots,N\}$, let $X_{-1,i}$ be defined recursively via* (5.80). *We have*

$$\frac{1}{N} \mathbb{E} \log \left\langle \exp \sum_{i=1}^N X_{K,i}(\omega_{\varnothing,i}, \omega_{\alpha_{|1},i}, \ldots, \omega_{\alpha_{|K-1},i}, \omega_{\alpha,i}) \right\rangle = \frac{1}{N} \sum_{i=1}^N X_{-1,i}. \qquad (5.81)$$

*Proof.* The main point is to realize that we can apply Theorem 5.25 to compute the term on the left side of (5.81). Indeed, we can find a measurable mapping $f = (f_1, \ldots, f_N) : [0,1] \to [0,1]^N$ that sends Lebesgue measure on $[0,1]$ to Lebesgue measure on $[0,1]^N$ (a simple splitting of the binary expansion of a real number will do). Letting $(\omega_k)_{0 \leq k \leq K}$ be independent uniform random variables on $[0,1]$, the random variables $(f_i(\omega_k))_{0 \leq k \leq K, 1 \leq i \leq N}$ have the same law as the $(\omega_{k,i})_{0 \leq k \leq K, 1 \leq i \leq N}$ in the definition of the variables $X_{k,i}$; and integrating each $\omega_k$ from $k = K$ downwards to $k = 0$ corresponds to integrating $(\omega_{i,k})_{1 \leq i \leq N}$ in the exact same way.

Let $Y_K := \sum_{i=1}^N X_{K,i}$. By independence,

$$\mathbb{E} \exp \zeta_K Y_K = \prod_{i=1}^N \mathbb{E} \exp \zeta_K X_{K,i} < +\infty.$$

Using again independence and an induction from $k = K - 1$ to $k = 0$, we verify that

$$Y_k := \frac{1}{\zeta_{k+1}} \log \mathbb{E}_{\omega_{k+1}} \exp \zeta_{k+1} Y_{k+1} = \sum_{i=1}^N X_{k,i},$$

and similarly for $k = -1$. The result is therefore a consequence of Theorem 5.25. ∎

We now turn our attention to extending the Ghirlanda-Guerra identities (Theorem 5.22) to the setting of Poisson-Dirichlet cascades. To define the analogue of the overlap (5.56), we fix a sequence

$$0 = q_{-1} \leq q_0 < q_1 < \ldots < q_K < q_{K+1} = +\infty, \qquad (5.82)$$

and define the overlap between two independent samples $\alpha^\ell$ and $\alpha^{\ell'}$ from the Poisson-Dirichlet cascade $\langle \cdot \rangle$ by

$$R_{\ell,\ell'} := q_{\alpha^\ell \wedge \alpha^{\ell'}}. \qquad (5.83)$$

In the context of mean-field spin glasses, we will think of this overlap $R_{\ell,\ell'}$ as representing the scalar product between the configurations whose Gibbs weights are represented by $v_{\alpha_\ell}$ and $v_{\alpha_{\ell'}}$.

In order to interpret the overlap (5.83) as a scalar product, we fix a separable Hilbert space $H$ and let $(e_\alpha)_{\alpha \in \mathcal{A}}$ be an orthonormal family in $H$. For each $\alpha \in \mathbb{N}^K$, we define $h_\alpha \in H$ by

$$h_\alpha := \sum_{k=0}^{K} (q_k - q_{k-1})^{1/2} e_{\alpha_{|k}}, \tag{5.84}$$

and observe that for every $\alpha^1, \alpha^2 \in \mathbb{N}^K$,

$$h_{\alpha^1} \cdot h_{\alpha^2} = \sum_{k=0}^{\alpha^1 \wedge \alpha^2} (q_k - q_{k-1}) = q_{\alpha^1 \wedge \alpha^2} = R_{1,2}. \tag{5.85}$$

In other words, if we identify each leaf $\alpha \in \mathbb{N}^K$ with the element $h_\alpha \in H$, then the overlap (5.83) between two leaves $\alpha^1$ and $\alpha^2$ becomes the scalar product between their representatives $h_{\alpha^1}$ and $h_{\alpha^2}$ in $H$.

To prove the Ghirlanda-Guerra identities (5.57) for the Poisson-Dirichlet cascades, we will rely on a generalization of the Bolthausen-Sznitman invariance property (Lemma 5.21). For technical reasons, to extend the Bolthausen-Sznitman invariance property to the setting of the Poisson-Dirichlet cascade, it will be useful to first assume that $q_0 = 0$ in the sequence (5.82). This assumption will be lifted with ease when we prove the Ghirlanda-Guerra identities, since these identities possess a natural invariance property as we vary the sequence $(q_k)_{0 \leq k \leq K}$.

Given a family $(z_\alpha)_{\alpha \in \mathcal{A} \setminus \mathbb{N}^0}$ of independent standard Gaussian random variables, we introduce the family $(Z_q(\alpha))_{\alpha \in \mathbb{N}^K}$ of Gaussian random variables defined by

$$Z_q(\alpha) := \sum_{k=0}^{K} (q_k - q_{k-1})^{1/2} z_{\alpha_{|k}} = \sum_{k=1}^{K} (q_k - q_{k-1})^{1/2} z_{\alpha_{|k}}. \tag{5.86}$$

Notice that the covariance structure

$$\mathbb{E} Z_q(\alpha^1) Z_q(\alpha^2) = q_{\alpha^1 \wedge \alpha^2} = R_{1,2} \tag{5.87}$$

of these random variables matches the definition of the overlaps of the Poisson-Dirichlet cascade introduced in (5.83). As a side remark, we also point out that if we view the $Z_q(\alpha)$'s as elements of the Hilbert space of square-integrable random variables on the underlying probability space, then the mapping $\alpha \mapsto Z_q(\alpha)$ provides us with a second way to realize the overlap (5.83) as a scalar product.

As was mentioned below Definition 5.24, we want to keep track of the information on the weights $(v_\alpha)_{\alpha \in \mathbb{N}^K}$ up to relabellings of $\mathbb{N}^K$ that preserve the tree structure. It would therefore not suffice for our purposes to obtain a version of the Bolthausen-Sznitman invariance that only concerns the *set* of weights $(w_\alpha)_{\alpha \in \mathbb{N}^K}$, or equivalently the point process $\sum_{\alpha \in \mathbb{N}^K} \delta_{w_\alpha}$ (whether or not each $w_\alpha$ is paired with some Gaussian random variable is not the point we want to stress here). To make this precise, we say that a bijection $\pi : \mathbb{N}^K \to \mathbb{N}^K$ *preserves the tree structure* if for every $\alpha, \beta \in \mathbb{N}^K$, we have that $\pi(\alpha) \wedge \pi(\beta) = \alpha \wedge \beta$.

**Lemma 5.27** (Bolthausen-Sznitman invariance for PDC [57]). *Suppose that $q_0 = 0$ in the sequence (5.82), let $(Z_q(\alpha))_{\alpha \in \mathbb{N}^K}$ denote the family of Gaussian random variables defined by (5.86) and taken independently of the random weights $(w_\alpha)_{\alpha \in \mathbb{N}^K}$ defined in (5.68), and let $b := \sum_{k=1}^{K}(q_k - q_{k-1})\zeta_k$. For each $t \in \mathbb{R}$, there exists a random permutation $\pi: \mathbb{N}^K \to \mathbb{N}^K$ that preserves the tree structure such that*

$$\left(w_{\pi(\alpha)}\exp\left(t\left(Z_q(\pi(\alpha))-\frac{tb}{2}\right)\right), Z_q(\pi(\alpha))-tb\right)_{\alpha \in \mathbb{N}^K} \tag{5.88}$$

*has the same law as $(w_\alpha, Z_q(\alpha))_{\alpha \in \mathbb{N}^K}$.*

*Proof.* Given a family $(z_\alpha)_{\alpha \in \mathcal{A} \smallsetminus \mathbb{N}^0}$ of independent standard Gaussian random variables, we define the family $(g_\alpha)_{\alpha \in \mathcal{A} \smallsetminus \mathbb{N}^0}$ of independent Gaussian random variables with variances $\mathsf{v}_{|\alpha|} := q_{|\alpha|} - q_{|\alpha|-1}$ by

$$g_\alpha := (q_{|\alpha|} - q_{|\alpha|-1})^{1/2} z_\alpha.$$

For each $\alpha \in \mathcal{A} \smallsetminus \mathbb{N}^K$, we recall that $\Lambda_\alpha = (u_{\alpha n})_{n \geq 1}$ is a Poisson-Dirichlet point process with parameter $\zeta_{|\alpha|+1} \in (0,1)$, and that these processes are independent as $\alpha$ varies in $\mathcal{A} \smallsetminus \mathbb{N}^K$. It follows by the Bolthausen-Sznitman invariance in Lemma 5.21 that the Poisson point processes

$$\Lambda_\alpha^+ := \sum_{n=1}^{+\infty} \delta_{(u_{\alpha n}, g_{\alpha n})} \quad \text{and} \quad \Lambda_\alpha' := \sum_{n=1}^{+\infty} \delta_{(u_{\alpha n}\exp(t(g_{\alpha n}-t\mathsf{v}_{|\alpha|+1}\zeta_{|\alpha|+1}/2)), g_{\alpha n}-t\mathsf{v}_{|\alpha|+1}\zeta_{|\alpha|+1})}$$

have the same law. We also observe that the random variables $(\Lambda_\alpha^+)_{\alpha \in \mathcal{A} \smallsetminus \mathbb{N}^K}$ are independent, and this also holds for $(\Lambda_\alpha')_{\alpha \in \mathcal{A} \smallsetminus \mathbb{N}^K}$. We denote by $\pi_\alpha: \mathbb{N} \to \mathbb{N}$ the permutation that sorts the random variables $(u_{\alpha n}\exp(t(g_{\alpha n}-t\mathsf{v}_{|\alpha|+1}\zeta_{|\alpha|+1}/2)))_{n \geq 1}$ decreasingly, and set

$$U_\alpha^+ := (u_{\alpha n}, g_{\alpha n})_{n \geq 1},$$

$$U_\alpha' := \left(u_{\alpha \pi_\alpha(n)}\exp\left(t\left(g_{\alpha \pi_\alpha(n)}-\frac{t\mathsf{v}_{|\alpha|+1}\zeta_{|\alpha|+1}}{2}\right)\right), g_{\alpha \pi_\alpha(n)}-t\mathsf{v}_{|\alpha|+1}\zeta_{|\alpha|+1}\right)_{n \geq 1}.$$

Since $(\Lambda_\alpha^+)_{\alpha \in \mathcal{A} \smallsetminus \mathbb{N}^K}$ and $(\Lambda_\alpha')_{\alpha \in \mathcal{A} \smallsetminus \mathbb{N}^K}$ have the same law, we have that $(U_\alpha^+)_{\alpha \in \mathcal{A} \smallsetminus \mathbb{N}^K}$ and $(U_\alpha')_{\alpha \in \mathcal{A} \smallsetminus \mathbb{N}^K}$ have the same law as well. We next define recursively the permutation $\pi: \mathcal{A} \to \mathcal{A}$ such that $\pi(\varnothing) = \varnothing$ and, for every $\alpha \in \mathcal{A} \smallsetminus \mathbb{N}^K$ and $n \in \mathbb{N}$,

$$\pi(\alpha n) := \pi(\alpha)\pi_{\pi(\alpha)}(n). \tag{5.89}$$

The fact that $\pi$ is indeed a permutation can be shown by induction on the depth of the tree. Notice also that $\pi$ preserves the parent-child relationship: by (5.89), for every $\alpha \in \mathcal{A} \smallsetminus \mathbb{N}^K$ and $n \in \mathbb{N}$, the permutation $\pi$ sends $\alpha n$ to a child of $\pi(\alpha)$. In particular, for each $k \in \{0,\ldots,K\}$, the restriction of $\pi$ to $\mathbb{N}^k$ is a permutation that preserves the tree structure. We now show by induction that, for every $k \in \{0,\ldots,K-1\}$,

$$(U_\alpha^+)_{|\alpha|\leq k} \stackrel{d}{=} (U_{\pi(\alpha)}')_{|\alpha|\leq k}, \tag{5.90}$$

where we recall that $\stackrel{d}{=}$ denotes equality in distribution. We have already observed that the identity (5.90) is valid when $k = 0$. Proceeding inductively, suppose that the identity is valid with $k \in \{0, \ldots, K-2\}$. Recall that the two families $(U_\alpha^+)_{|\alpha| \leq k}$ and $(U_\alpha^+)_{|\alpha|=k+1}$ are independent, and that the random variables $(U_\alpha^+)_{|\alpha|=k+1}$ are independent as $\alpha$ varies. We denote by $\mathcal{F}_k$ the $\sigma$-algebra generated by $(U_\alpha^+)_{|\alpha| \leq k}$. The family $(U'_{\pi(\alpha)})_{|\alpha| \leq k}$ is $\mathcal{F}_k$-measurable, and so is $(\pi(\alpha))_{|\alpha| \leq k+1}$. Conditionally on $\mathcal{F}_k$, the random variables $(U'_\alpha)_{|\alpha|=k+1}$ are independent, and the family has the same law as $(U_\alpha^+)_{|\alpha|=k+1}$. Since $(\pi(\alpha))_{|\alpha|=k+1}$ is $\mathcal{F}_k$-measurable, this remains valid for the family $(U'_{\pi(\alpha)})_{|\alpha|=k+1}$. Combining these observations yields (5.90) with $k$ replaced by $k+1$. We can now deduce from (5.90) that

$$\Big( \prod_{k=1}^K u_{\alpha_{|k}}, \sum_{k=1}^K g_{\alpha_{|k}} \Big)_{\alpha \in \mathbb{N}^K}$$

has the same law as

$$\Big( \prod_{k=1}^K u_{\pi(\alpha)_{|k}} \exp\Big(t\Big(g_{\pi(\alpha)_{|k}} - \frac{t v_k \zeta_k}{2}\Big)\Big), \sum_{k=1}^K (g_{\pi(\alpha)_{|k}} - t v_k \zeta_k) \Big)_{\alpha \in \mathbb{N}^K}.$$

Recalling the definition of the family of Gaussian random variables $(Z_q(\alpha))_{\alpha \in \mathbb{N}^K}$, of the weights $(w_\alpha)_{\alpha \in \mathbb{N}^K}$ and of the constant $b = \sum_{k=1}^K v_k \zeta_k$ completes the proof. ∎

**Theorem 5.28** (Ghirlanda-Guerra identities for PDC [62])**.** *For every $n \geq 1$, every bounded and measurable function $f$ of the overlaps $R^n = (R_{\ell,\ell'})_{\ell,\ell' \leq n}$, and every bounded and measurable function $\psi : \mathbb{R} \to \mathbb{R}$, we have*

$$\mathbb{E}\langle f(R^n) \psi(R_{1,n+1}) \rangle = \frac{1}{n} \mathbb{E}\langle f(R^n) \rangle \mathbb{E}\langle \psi(R_{1,2}) \rangle + \frac{1}{n} \sum_{\ell=2}^n \mathbb{E}\langle f \psi(R_{1,\ell}) \rangle \quad (5.91)$$

*and*

$$\mathbb{E}\langle \psi(R_{1,2}) \rangle = \sum_{k=0}^K \psi(q_k)(\zeta_{k+1} - \zeta_k). \quad (5.92)$$

*In particular, for every $k \in \{0, \ldots, K\}$,*

$$\mathbb{E}\langle \mathbf{1}_{\{\alpha^1 \wedge \alpha^2 = k\}} \rangle = \zeta_{k+1} - \zeta_k. \quad (5.93)$$

*Proof.* We first explain why we can assume that $q_0 = 0$ without loss of generality. Notice that the overlap $R_{\ell,\ell'}$ is obtained by first computing $\alpha^\ell \wedge \alpha^{\ell'}$, and then mapping the result through the function

$$\begin{cases} \{0, \ldots, K\} & \to & \mathbb{R} \\ k & \mapsto & q_k. \end{cases}$$

Defining $\widehat{R}_{\ell,\ell'} := \alpha^\ell \wedge \alpha^{\ell'}$, it therefore follows that (5.91)-(5.92) can be rephrased as the statement that, for every $f : \mathbb{R}^{n \times n} \to \mathbb{R}$ and $\psi : \mathbb{R} \to \mathbb{R}$,

$$\mathbb{E}\langle f(\widehat{R}^n)\psi(\widehat{R}_{1,n+1})\rangle = \frac{1}{n}\mathbb{E}\langle f(\widehat{R}^n)\rangle \mathbb{E}\langle \psi(\widehat{R}_{1,2})\rangle + \frac{1}{n}\sum_{\ell=2}^n \mathbb{E}\langle f\psi(\widehat{R}_{1,\ell})\rangle,$$

and

$$\mathbb{E}\langle \psi(\widehat{R}_{1,2})\rangle = \sum_{k=0}^K \psi(k)(\zeta_{k+1} - \zeta_k).$$

Notice that these two statements do not depend on the choice of the values of $(q_0, \ldots, q_K)$. In other words, as soon as we can prove Theorem 5.28 for some choice of $(q_0, \ldots, q_K)$ satisfying (5.82), we can deduce its validity for all such choices. In particular, there is no loss of generality in imposing that $q_0 = 0$. (This invariance is not valid for Lemma 5.27, since the law of $Z_q$ depends on the choice of $(q_0, \ldots, q_K)$.)

So we assume from now on that $q_0 = 0$. Moreover, as the overlap takes only finitely many values $q_0, \ldots, q_K$, by polynomial interpolation we may assume without loss of generality that $\psi$ is a polynomial. By linearity it therefore suffices to prove (5.91) and (5.92) when $\psi(x) = x^m$ for some integer $m \in \mathbb{N}$. Since (5.91) and (5.92) are clear for $m = 0$, let us assume that $m \geq 1$. Given a family $(z_\alpha)_{\alpha \in \mathcal{A} \setminus \mathbb{N}^0}$ of independent standard Gaussian random variables, we introduce the family $(Z_{q^m}(\alpha))_{\alpha \in \mathbb{N}^K}$ of Gaussian random variables defined by

$$Z_{q^m}(\alpha) := \sum_{k=1}^K \left(q_k^m - q_{k-1}^m\right)^{1/2} z_{\alpha_{|k}}.$$

We write $(w_\alpha)_{\alpha \in \mathbb{N}^K}$ for the random weights in (5.68), and let

$$b_m := \sum_{k=1}^K (q_k^m - q_{k-1}^m)\zeta_k.$$

For each $t \in \mathbb{R}$ and $\alpha \in \mathbb{N}^K$, we define the random weight

$$v_\alpha^{t,m} := \frac{w_\alpha \exp t\left(Z_{q^m}(\alpha) - \frac{tb_m}{2}\right)}{\sum_{\alpha \in \mathbb{N}^K} w_\alpha \exp t\left(Z_{q^m}(\alpha) - \frac{tb_m}{2}\right)} = \frac{w_\alpha \exp(tZ_{q^m}(\alpha))}{\sum_{\alpha \in \mathbb{N}^K} w_\alpha \exp(tZ_{q^m}(\alpha))},$$

and denote by $\langle \cdot \rangle_{t,m}$ the Gibbs measure (5.72) with the weights $(v_\alpha)_{\alpha \in \mathbb{N}^K}$ replaced by the weights $(v_\alpha^{t,m})_{\alpha \in \mathbb{N}^K}$. By the Bolthausen-Sznitman invariance (Lemma 5.27), there exists a random permutation $\pi : \mathbb{N}^K \to \mathbb{N}^K$ that preserves the tree structure such that the families

$$\left(v_{\pi(\alpha)}^{t,m}, Z_{q^m}(\pi(\alpha)) - tb_m\right)_{\alpha \in \mathbb{N}^K} \quad \text{and} \quad \left(v_\alpha^{0,m}, Z_{q^m}(\alpha)\right)_{\alpha \in \mathbb{N}^K} \quad (5.94)$$

have the same law. This implies that for any $f : \mathbb{R}^{n \times n} \to \mathbb{R}$,

$$\begin{aligned}
0 &= \mathbb{E}\langle f(R^n) Z_{q^m}(\alpha^1)\rangle \\
&= \mathbb{E} \sum_{\alpha^1,\ldots,\alpha^n \in \mathbb{N}^K} f(R^n) Z_{q^m}(\alpha^1) v_{\alpha^1}^{0,m} \cdots v_{\alpha^n}^{0,m} \\
&= \mathbb{E} \sum_{\alpha^1,\ldots,\alpha^n \in \mathbb{N}^K} f(R^n) \big(Z_{q^m}(\pi(\alpha^1)) - t b_m\big) v_{\pi(\alpha^1)}^{t,m} \cdots v_{\pi(\alpha^n)}^{t,m} \\
&= \mathbb{E} \sum_{\alpha^1,\ldots,\alpha^n \in \mathbb{N}^K} f(R^n) \big(Z_{q^m}(\alpha^1) - t b_m\big) v_{\alpha^1}^{t,m} \cdots v_{\alpha^n}^{t,m} \\
&= \mathbb{E}\langle f(R^n) (Z_{q^m}(\alpha^1) - t b_m)\rangle_{t,m},
\end{aligned}$$

where in the fourth equality, we used that $\pi$ is a permutation that preserves the tree structure to ensure that the action of $\pi$ on $R^n$ is trivial. Invoking the Gibbs Gaussian integration by parts formula (Theorem 4.6) yields that

$$t \mathbb{E}\Big\langle f(R^n)\Big(\sum_{\ell=1}^n R_{1,\ell}^m - n R_{1,n+1}^m - b_m\Big)\Big\rangle_{t,m} = 0.$$

Using again that the families in (5.94) have the same law with $\pi$ preserving the tree structure, recalling that $R_{1,1} = q_K$ by (5.83), and rearranging, we obtain that

$$\mathbb{E}\langle f(R^n) R_{1,n+1}^m\rangle = \frac{1}{n}\sum_{\ell=2}^n \mathbb{E}\langle f(R^n) R_{1,\ell}^m\rangle + \frac{q_K^m - b_m}{n}\mathbb{E}\langle f(R^n)\rangle. \qquad (5.95)$$

Taking $f = n = 1$ reveals that

$$\mathbb{E}\langle R_{1,2}^m\rangle = q_K^m - b_m = \sum_{k=0}^K q_k^m(\zeta_{k+1} - \zeta_k),$$

where we have used that $\zeta_0 = 0$ and $\zeta_{K+1} = 1$. This establishes (5.92) and, together with (5.95), also gives (5.91). Taking $\psi(x) := \mathbf{1}_{\{x = q_k\}}$ for $0 \leqslant k \leqslant K$ completes the proof. ∎

As was already mentioned several times, in the notion of Poisson-Dirichlet cascade, we pay attention to keeping track of the tree structure of $\mathbb{N}^K$, as opposed to simply studying the set of weights $\{w_\alpha \mid \alpha \in \mathbb{N}^K\}$, or equivalently the point process $\sum_{\alpha \in \mathbb{N}^K} \delta_{w_\alpha}$. In fact, this point process has the same law as the point process $\sum_{n=1}^{+\infty} \delta_{v_n}$, if $(v_n)_{n \geqslant 1}$ is a Poisson-Dirichlet process of parameter $\zeta_K$. This can be shown by induction on the number of layers of the cascade. Alternatively, one can use Theorem 5.28 to observe that the overlap defined by $\widetilde{R}_{\ell,\ell'} := \mathbf{1}_{\{\alpha^\ell = \alpha^{\ell'}\}}$ satisfies the Ghirlanda-Guerra identities with $\mathbb{E}\langle \mathbf{1}_{\{\alpha^1 = \alpha^2\}}\rangle = 1 - \zeta_K$, and then use Exercise 5.10 to conclude. But the cascade $(w_\alpha)_{\alpha \in \mathbb{N}^K}$ itself is richer as it keeps track of how the indices are organized in the tree structure, an aspect which cannot be recovered from $\sum_{\alpha \in \mathbb{N}^K} \delta_{w_\alpha}$.

**Exercise 5.11.** Given a right-continuous and non-decreasing function $F : \mathbb{R} \to [0, 1]$, let $F^{-1} : (0, 1) \to \mathbb{R}$ denote its *quantile transform*,

$$F^{-1}(x) := \inf\{s \in \mathbb{R} \mid F(s) \geq x\}. \tag{5.96}$$

(i) Show that for all $t \in \mathbb{R}$ and $x \in (0, 1)$, we have $F(t) < x$ if and only if $t < F^{-1}(x)$.

(ii) Let $X$ be a real-valued random variable, and denote by $F$ its distribution function, $F(x) := \mathbb{P}\{X \leq x\}$. Letting $U$ be a uniform random variable on the interval $[0, 1]$, show that $F^{-1}(U)$ and $X$ have the same law.

**Exercise 5.12.** Fix $\varepsilon > 0$ and a random variable $X$ such that $\mathbb{E} \exp \varepsilon |X| < +\infty$. Prove that

$$\lim_{\zeta \searrow 0} \frac{1}{\zeta} \log \mathbb{E} \exp \zeta X = \mathbb{E} X. \tag{5.97}$$

## 5.7  Ultrametricity as a universal property

We have already discussed informally that Poisson-Dirichlet cascades provide us with "canonical models" for the asymptotic behaviour of Gibbs measures of mean-field spin glasses. The goal of this section is to give some substance to this claim. We start by discussing background motivation in more detail.

In the context of mean-field spin glasses, as in statistical inference, we will be concerned with the study of a random Gibbs measure on $\mathbb{R}^N$, and we would like to find ways to capture its asymptotic behaviour. To fix notation, we denote by $\langle \cdot \rangle_N$ the expectation with respect to this random Gibbs measure over $\mathbb{R}^N$, with canonical random variable $\sigma$, and by $(\sigma^\ell)_{\ell \geq 1}$ independent copies of the random variable $\sigma$ under $\langle \cdot \rangle_N$. Recalling that $\langle \cdot \rangle_N$ is itself random, we denote by $\mathbb{E}$ the expectation with respect to this additional source of randomness.

It is not at all clear how to capture asymptotic properties of a sequence of probability measures defined on spaces of increasingly large dimension. Besides, the spin-glass models we aim to consider have a very high degree of symmetry, and there is no axis or other subspace of small dimension that we would particularly want to concentrate upon. A very fruitful approach is to keep track of the array of overlaps $R := (R_{\ell, \ell'})_{\ell, \ell' \geq 1}$ between the different replicas. For one thing, it is shown in Lemma 1.7 of [211] that the law of this overlap array characterizes the underlying Gibbs measure up to a random orthogonal transformation. Since the models we will study have a very large group of symmetries, due to their mean-field character, there is hopefully not much loss of information if one only understands the Gibbs measure up to an orthogonal transformation. Moreover, under modest assumptions on the model, we know a priori that each coordinate of $R$ takes values in a bounded set, say $[-1, 1]$ if the support of the Gibbs measure on $\mathbb{R}^N$ is contained in $B_{\sqrt{N}}(0)$.

So we can at least extract convergent subsequences, and wonder about properties of the limit overlap array. From now on we denote by $R$ any such limit overlap array.

The array $R$ must be symmetric positive semi-definite. Notice also that if $\pi$ is a permutation of finitely many indices in $\mathbb{N}$, then the array $(R_{\pi(\ell),\pi(\ell')})_{\ell,\ell'\geq 1}$ has the same law as the array $R$. Under these assumptions, the Dovbysh-Sudakov theorem (Theorem 1.7 in [211]) asserts that the off-diagonal part of the array $R$ can still be represented as an array of overlaps. More precisely, there exists a random probability measure $G$ on some Hilbert space, with associated expectation $\langle\cdot\rangle$, such that, using again the notation $(\sigma^\ell)_{\ell\geq 1}$ for the replicas under this new probability measure, the off-diagonal elements $(R_{\ell,\ell'})_{\ell\neq\ell'}$ of the limit array have the same joint distribution as $(\sigma^\ell\cdot\sigma^{\ell'})_{\ell\neq\ell'}$. This result can be thought of as a variant of the de Finetti representation theorem for exchangeable sequences.

Notice that in the statement of the Dovbysh-Sudakov theorem, we excluded the diagonal elements of the array $R$ from the representation. The result does not extend to on-diagonal elements as stated, as their might be some "loss of norm" as we pass to the limit. We illustrate that this phenomenon of loss of norm can indeed happen on a simple example. We take the uniform measure on the sphere of radius $\sqrt{N}$ in $\mathbb{R}^N$ as the sequence of Gibbs measures (so these are non-random in this example). In this case, the limit overlap array is $R_{\ell,\ell'} = \mathbf{1}_{\{\ell=\ell'\}}$. The only choice then to represent the off-diagonal elements of this array is to take the limit Gibbs measure to be concentrated at the origin, so $\sigma^\ell\cdot\sigma^\ell = 0$, and this does not match the on-diagonal elements of the limit overlap array that are equal to 1.

For convenience, until the end of this section, *we redefine $R$ to be $(R_{\ell,\ell'})_{\ell\neq\ell'}$*, i.e. we leave aside the on-diagonal elements.

A crucial observation for mean-field spin-glass models is that, up to the addition of a small perturbation to the Hamiltonian, we can ensure that any possible limit overlap array $R$ must satisfy the Ghirlanda-Guerra identities. The perturbation added to the Hamiltonian is small in the sense that it does not alter the value of the limit free energy. We have seen in Exercise 5.10 that when the overlap can only takes values in $\{0,1\}$, an overlap array that satisfies the Ghirlanda-Guerra identities must have the same law as that of a Poisson-Dirichlet process. It turns out that this observation admits a powerful generalization: if an overlap array satisfies the Ghirlanda-Guerra identities, then it must have the same law as that of a Poisson-Dirichlet cascade. This is the operationalization of the idea that Poisson-Dirichlet cascades provide us with "canonical models" for what the Gibbs measure of a mean-field spin glass asymptotically looks like.

To state this characterization of Poisson-Dirichlet cascades precisely, let us start by recapitulating its context. We study a random probability measure $G$ supported on the unit ball of a Hilbert space $H$, whose expectation we denote by $\langle\cdot\rangle$, and we denote by $\mathbb{E}$ the expectation with respect to the randomness of the random measure $G$ itself; we write $(\sigma^\ell)_{\ell\geq 1}$ for the replicas under $G$, and $R = (\sigma^\ell\cdot\sigma^{\ell'})_{\ell\neq\ell'}$

stands for their overlap array. We denote by $\zeta$ the law of one overlap under $\mathbb{E}\langle\cdot\rangle$, so that for every measurable set $A \subseteq \mathbb{R}$,

$$\zeta(A) = \mathbb{E}\langle \mathbf{1}_{\{R_{1,2} \in A\}}\rangle. \tag{5.98}$$

We say that the measure $G$ satisfies the Ghirlanda-Guerra identities if for every $n \geq 1$, every bounded measurable function $f$ of the overlaps $R^n = (R_{\ell,\ell'})_{\ell \neq \ell' \leq n}$, and every bounded measurable function $\psi : \mathbb{R} \to \mathbb{R}$, we have

$$\mathbb{E}\langle f(R^n)\psi(R_{1,n+1})\rangle = \frac{1}{n}\mathbb{E}\langle f(R^n)\rangle\mathbb{E}\langle\psi(R_{1,2})\rangle + \frac{1}{n}\sum_{\ell=2}^{n}\mathbb{E}\langle f(R^n)\psi(R_{1,\ell})\rangle. \tag{5.99}$$

Equivalently, the measure $G$ satisfies the Ghirlanda-Guerra identities if, for every $n \geq 1$, the law of $R_{1,n+1}$ under $\mathbb{E}\langle\cdot\rangle$ conditionally on $R^n$ is given by

$$\frac{1}{n}\zeta + \frac{1}{n}\sum_{\ell=2}^{n}\delta_{R_{1,\ell}}. \tag{5.100}$$

**Theorem 5.29** (Characterization of PDC [210]). *If the random probability measure $G$ satisfies the Ghirlanda-Guerra identities (5.99), then the law of the infinite overlap array $R = (R_{\ell,\ell'})_{\ell \neq \ell'}$ under $\mathbb{E}\langle\cdot\rangle$ is uniquely determined by the probability measure $\zeta$ in (5.98). Moreover, the measure $\zeta$ is supported on $\mathbb{R}_{\geq 0}$.*

Concerning the on-diagonal elements $(\sigma^\ell \cdot \sigma^\ell)_{\ell \geq 1}$, one can show that the Ghirlanda-Guerra identities impose that $G$ be supported on the sphere in $H$ whose radius is the square root of the top of the support of $\zeta$. In other words, denoting by $q^*$ the top of the support of $\zeta$, we have, with probability one under $\mathbb{E}\langle\cdot\rangle$, that

$$\|\sigma^\ell\|^2 = q^*. \tag{5.101}$$

We refer to Theorem 2.15 of [211] for a (relatively simple) proof of this fact.

When the support of $\zeta$ is finite, we already know from Theorem 5.28 that the overlap array associated with a Poisson-Dirichlet cascade satisfies the Ghirlanda-Guerra identities, with $\zeta$ being the measure given by

$$\zeta = \sum_{k=0}^{K}(\zeta_{k+1} - \zeta_k)\delta_{q_k}. \tag{5.102}$$

Recall also from the discussion around (5.84)-(5.85) that the overlap array associated with a Poisson-Dirichlet cascade can indeed be interpreted as a matrix of scalar products for a suitable choice of random probability measure. Theorem 5.29 thus tells us that, when the support of $\zeta$ is finite, the law of the overlap array coming from $G$ is the same as that coming from the Poisson-Dirichlet cascade with parameters $(\zeta_k)_{1 \leq k \leq K}$ and $(q_k)_{0 \leq k \leq K}$ chosen so that (5.102) holds.

For more general $\zeta$, we have not yet built a canonical example of an overlap array satisfying (5.98) and the Ghirlanda-Guerra identities. This can be achieved by approximation, as will be explained in Corollary 5.32.

In short, up to extending the notion of Poisson-Dirichlet cascade and its associated overlap array by continuity, we see that if a random probability measure $G$ satisfies the Ghirlanda-Guerra identities, then the law of its overlap array must be that of a Poisson-Dirichlet cascade.

The fundamental step in the proof of Theorem 5.29 involves the concept of ultrametricity. We say that the support of the measure $G$ is *ultrametric* if, with probability one over $\mathbb{E}\langle \cdot \rangle$, we have

$$R_{1,3} \geqslant \min(R_{1,2}, R_{2,3}). \tag{5.103}$$

In such a circumstance, we may also say that the overlap array $R$ is ultrametric. Recalling from (5.101) that $G$ is supported on a sphere, we can rephrase the ultrametricity property as the statement that, with probability one over $\mathbb{E}\langle \cdot \rangle$, we have

$$\|\sigma^1 - \sigma^3\| \leqslant \max\left(\|\sigma^1 - \sigma^2\|, \|\sigma^2 - \sigma^3\|\right). \tag{5.104}$$

Equivalently, this says that with probability one over $\mathbb{E}$, every choice of three points $\sigma^1$, $\sigma^2$, and $\sigma^3$ in the support of $G$ must be such that (5.104) holds. In words, this stronger form of the triangle inequality means that every triangle has its two longest sides of equal length. Indeed, for any three points $\sigma^1$, $\sigma^2$, and $\sigma^3$ in the support of $G$, if the segment joining $\sigma^1$ to $\sigma^2$ happens to be the shortest of the three sides of the triangle, then (5.104) ensures that the segment joining $\sigma^1$ to $\sigma^3$ is at most as long as that joining $\sigma^2$ to $\sigma^3$; and by symmetry, we conclude that these two segments have the same length. In terms of overlaps, this means that among $R_{1,2}$, $R_{1,3}$, and $R_{2,3}$, the two smallest overlaps must be equal.

The notion of ultrametricity is intimately tied with tree structures. For convenience we will assume that the support of the overlap distribution $\zeta$ is finite, and argue that the ultrametricity property allows us to organize the support of $G$ into the leaves of a tree, so that the distance (or equivalently for us: the overlap) between any two leaves $\alpha$ and $\beta$ is only a function of $\alpha \wedge \beta$, the depth of their most recent common ancestor. Indeed, the ultrametricity property implies that, for every $p$ in the support of $\zeta$, the relation

$$\sigma^1 \sim_p \sigma^2 \iff R_{1,2} \geqslant p \tag{5.105}$$

defines an equivalence relation. We then draw the tree as in Figure 5.1, so that the nodes at depth $k$ represent the equivalence classes of the relation $\sim_p$ for $p$ the $k^{\text{th}}$ atom in the support of $\zeta$, i.e. for $p = q_k$ if $\zeta$ is as in (5.102). We refer to Lemma 4.2 of [70] for a more precise construction.

In particular, one can easily see that Poisson-Dirichlet cascades satisfy the ultrametricity property (5.103), using that for any three leaves $\alpha$, $\beta$, $\gamma$ we have

$$\alpha \wedge \gamma \geqslant \min(\alpha \wedge \beta, \beta \wedge \gamma). \tag{5.106}$$

The crucial step towards establishing Theorem 5.29 is to show that the Ghirlanda-Guerra identities imply ultrametricity. The setting for this result is as discussed in the paragraph above (5.98) and we do not restate it here.

**Theorem 5.30** (Panchenko ultrametricity [210])**.** *If the random probability measure G satisfies the Ghirlanda-Guerra identities* (5.99)*, then its support is ultrametric.*

We refer the interested reader to Theorem 2.14 in [211] for a proof of this result, and to [208] for a simpler argument in the special case when the support of $\zeta$ is finite.

Once Theorem 5.30 is established, the proof of Theorem 5.29 is comparatively less demanding. We first illustrate its core mechanics by deriving the joint law of the three overlaps $(R_{1,2}, R_{1,3}, R_{2,3})$, assuming the validity of the Ghirlanda-Guerra identities and the ultrametricity property. We will then generalize this argument and thereby deduce Theorem 5.29 from Theorem 5.30.

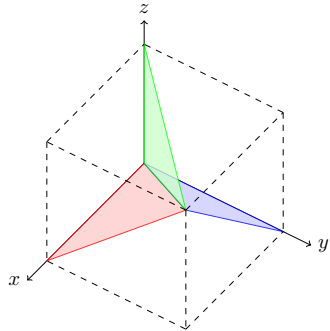

**Figure 5.2** When ultrametricity holds, the law of the triple of overlaps $(R_{1,2}, R_{1,3}, R_{2,3})$ is supported on the set of points $(x, y, z) \in [0,1]^3$ with two equal coordinates smaller than the third.

Theorem 5.30 states that the joint law of the triple of overlaps $(R_{1,2}, R_{1,3}, R_{2,3})$ is supported on the set of points with two equal coordinates smaller than the third. This set is displayed in Figure 5.2. We first write down the law of this triple away from the diagonal $\{x = y = z\}$. Starting with the off-diagonal part of the red isosceles triangle in Figure 5.2, which corresponds to the set $\{(x, y, z) \in [0,1]^3 \mid z = y < x\}$, we have that for every bounded measurable function $\psi : \mathbb{R}^3 \to \mathbb{R}$,

$$\mathbb{E}\langle \psi(R_{1,2}, R_{1,3}, R_{2,3}) \mathbf{1}_{\{R_{1,2} > R_{1,3} = R_{2,3}\}} \rangle = \frac{1}{2} \int_{\{x > y\}} \psi(x, y, y) \, d\zeta(x) \, d\zeta(y). \tag{5.107}$$

Indeed, recall that ultrametricity forces the two smallest overlaps among $(R_{1,2}, R_{1,3}, R_{2,3})$ to always be equal, so

$$\mathbb{E}\langle \psi(R_{1,2},R_{1,3},R_{2,3})\mathbf{1}_{\{R_{1,2}>R_{1,3}=R_{2,3}\}}\rangle = \mathbb{E}\langle \psi(R_{1,2},R_{1,3},R_{1,3})\mathbf{1}_{\{R_{1,2}>R_{1,3}\}}\rangle, \quad (5.108)$$

and then the Ghirlanda-Guerra identities and (5.98) yield (5.107). By symmetry, we can also write down the joint law of the three overlaps on the off-diagonal parts of the blue and green triangles in Figure 5.2.

To determine the joint law on the diagonal $\{(x,y,z) \in [0,1]^3 \mid x = y = z\}$, we can use ultrametricity to write

$$\mathbb{E}\langle \psi(R_{1,2},R_{1,3},R_{2,3})\mathbf{1}_{\{R_{1,2}=R_{1,3}=R_{2,3}\}}\rangle$$
$$= \mathbb{E}\langle \psi(R_{1,2},R_{1,2},R_{1,2})\mathbf{1}_{\{R_{1,2}=R_{1,3}\}}\rangle - \mathbb{E}\langle \psi(R_{1,2},R_{1,2},R_{1,2})\mathbf{1}_{\{R_{1,2}=R_{1,3}<R_{2,3}\}}\rangle.$$

Using ultrametricity once again, we can drop the condition $R_{1,2} = R_{1,3}$ in the second term since $R_{1,2} < R_{2,3}$, and then use the Ghirlanda-Guerra identities to find that

$$\mathbb{E}\langle \psi(R_{1,2},R_{1,3},R_{2,3})\mathbf{1}_{\{R_{1,2}=R_{1,3}=R_{2,3}\}}\rangle$$
$$= \frac{1}{2}\int \psi(x,x,x)\,\mathrm{d}\zeta(x) + \frac{1}{2}\int_{\{x=y\}} \psi(x,x,x)\,\mathrm{d}\zeta(x)\,\mathrm{d}\zeta(y)$$
$$- \frac{1}{2}\int_{\{x<z\}} \psi(x,x,x)\,\mathrm{d}\zeta(x)\,\mathrm{d}\zeta(z)$$
$$= \frac{1}{2}\int_{\{x>y\}} \psi(x,x,x)\,\mathrm{d}\zeta(x)\,\mathrm{d}\zeta(y) + \int_{\{x=y\}} \psi(x,x,x)\,\mathrm{d}\zeta(x)\,\mathrm{d}\zeta(y).$$

Bringing all terms together, we conclude that for any bounded measurable function $\psi : \mathbb{R}^3 \to \mathbb{R}$, we have

$$\mathbb{E}\langle \psi(R_{1,2},R_{1,3},R_{2,3})\rangle$$
$$= \frac{1}{2}\int_{\{x>y\}} \left(\psi(x,y,y) + \psi(y,x,y) + \psi(y,y,x) + \psi(x,x,x)\right) \mathrm{d}\zeta(x)\,\mathrm{d}\zeta(y)$$
$$+ \int_{\{x=y\}} \psi(x,x,x)\,\mathrm{d}\zeta(x)\,\mathrm{d}\zeta(y). \quad (5.109)$$

This shows in particular that the law of the triple $(R_{1,2},R_{1,3},R_{2,3})$ is uniquely determined by the probability measure $\zeta$, as stated in Theorem 5.29. We now generalize this argument and show that, under the assumptions of ultrametricity and of the validity of the Ghirlanda-Guerra identities, the law of the full overlap array is completely determined in terms of the law of $R_{1,2}$. We recall that the setting is as explained in the paragraph above (5.98).

**Proposition 5.31.** *If the random probability measure $G$ satisfies the Ghirlanda-Guerra identities and its support is ultrametric, then the law of the infinite overlap array $R = (R_{\ell,\ell'})_{\ell \neq \ell'}$ under $\mathbb{E}\langle \cdot \rangle$ is uniquely determined by the probability measure $\zeta$ in* (5.98).

*Proof.* Let $n \geq 1$ be an integer. We show the following statement by induction on $n$. Letting $(\ell_1, \ell'_1), \ldots, (\ell_n, \ell'_n)$ be $n$ distinct pairs of integers such that $\ell_k < \ell'_k$ for every $k \leq n$, letting $\preccurlyeq_1, \ldots, \preccurlyeq_{n-1} \in \{=, <\}$, and letting $\psi : \mathbb{R}^n \to \mathbb{R}$ be a bounded and measurable function, the quantity

$$\mathbb{E}\langle \psi(R_{\ell_1,\ell'_1}, \ldots, R_{\ell_n,\ell'_n}) \mathbf{1}_{\{R_{\ell_1,\ell'_1} \preccurlyeq_1 \cdots \preccurlyeq_{n-1} R_{\ell_n,\ell'_n}\}} \rangle$$

is uniquely determined from the knowledge of the law $\zeta$ of the overlap $R_{1,2}$.

The claim is immediate for $n = 1$. We now assume that the claim is valid for some integer $n \geq 1$, and give ourselves $(\ell_1, \ell'_1), \ldots, (\ell_{n+1}, \ell'_{n+1})$ as above, order relations $\preccurlyeq_1, \ldots, \preccurlyeq_n \in \{=, <\}$, and a bounded and measurable function $\psi : \mathbb{R}^{n+1} \to \mathbb{R}$. By symmetry between the replicas, we may assume that $\ell_{n+1} = 1$ and $\ell'_{n+1} = L$, where $L$ is the largest element of $\{\ell'_1, \ldots, \ell'_{n+1}\}$. We distinguish two cases.

*Case 1.* Suppose that the order relation $\preccurlyeq_n$ is $<$. If

$$L \notin \{\ell'_1, \ldots, \ell'_n\}, \tag{5.110}$$

then we can readily apply the Ghirlanda-Guerra identities to simplify the expectation we aim to compute, since we recall that conditionally on $R^{L-1}$, the law of $R_{1,L}$ under $\mathbb{E}\langle \cdot \rangle$ is given by (5.100) (substituting $n$ by $L-1$ there). Otherwise, let $k \leq n$ be such that $\ell'_k = L$. By ultrametricity, the events $R_{k,L} < R_{1,L}$ and $R_{1,k} < R_{1,L}$ are identical, and on this event, we have that $R_{1,k} = R_{k,L}$. In the expression

$$\mathbb{E}\langle \psi(R_{\ell_1,\ell'_1}, \ldots, R_{\ell_n,\ell'_n}, R_{1,n+1}) \mathbf{1}_{\{R_{\ell_1,\ell'_1} \preccurlyeq_1 \cdots \preccurlyeq_{n-1} R_{\ell_n,\ell'_n} < R_{1,n+1}\}} \rangle,$$

we can therefore substitute every occurrence of $R_{k,L}$ among $R_{\ell_1,\ell'_1}, \ldots, R_{\ell_n,\ell'_n}$ by $R_{1,k}$. We can do so for every $k$ such that $\ell'_k = L$, until we arrive at an expression that involves no term of the form $R_{k,L}$ for any $k$ among the overlaps $R_{\ell_1,\ell'_1}, \ldots, R_{\ell_n,\ell'_n}$. Once this is done, we can simplify the expression by applying the Ghirlanda-Guerra identities, as already observed.

*Case 2.* We now assume that the order relation $\preccurlyeq_n$ is the equality $=$. There might be several consecutive occurrences of the equality sign suffixing the sequence $\preccurlyeq_1, \ldots, \preccurlyeq_n$, and our goal is to reduce this number of suffixing equality signs by at least one, so that by induction, we can ultimately reduce all computations to occurrences of Case 1. We write

$$\mathbb{E}\langle \psi(R_{\ell_1,\ell'_1},\ldots,R_{\ell_n,\ell'_n},R_{1,L})\mathbf{1}_{\{R_{\ell_1,\ell'_1} \preccurlyeq_1 \cdots \preccurlyeq_{n-1} R_{\ell_n,\ell'_n} = R_{1,L}\}}\rangle$$
$$= \mathbb{E}\langle \psi(R_{\ell_1,\ell'_1},\ldots,R_{\ell_n,\ell'_n},R_{\ell_n,\ell'_n})\mathbf{1}_{\{R_{\ell_1,\ell'_1} \preccurlyeq_1 \cdots \preccurlyeq_{n-1} R_{\ell_n,\ell'_n}\}}\rangle$$
$$- \mathbb{E}\langle \psi(R_{\ell_1,\ell'_1},\ldots,R_{\ell_n,\ell'_n},R_{\ell_n,\ell'_n})\mathbf{1}_{\{R_{\ell_1,\ell'_1} \preccurlyeq_1 \cdots \preccurlyeq_{n-1} R_{\ell_n,\ell'_n} < R_{1,L}\}}\rangle$$
$$- \mathbb{E}\langle \psi(R_{\ell_1,\ell'_1},\ldots,R_{\ell_n,\ell'_n},R_{\ell_n,\ell'_n})\mathbf{1}_{\{R_{\ell_1,\ell'_1} \preccurlyeq_1 \cdots \preccurlyeq_{n-1} R_{\ell_n,\ell'_n}\}}\mathbf{1}_{\{R_{\ell_n,\ell'_n} > R_{1,L}\}}\rangle.$$

The first term on the right side of this equality can be computed by the induction hypothesis, while the second one is an instance of Case 1. The last term can be decomposed into a sum, by partitioning the event according to the position of $R_{1,L}$ relatively to the overlaps $R_{\ell_1,\ell'_1},\ldots,R_{\ell_n,\ell'_n}$. That is, we rewrite the last term of the previous display as a sum of terms of the form

$$\mathbb{E}\langle \psi(R_{\ell_1,\ell'_1},\ldots,R_{\ell_n,\ell'_n},R_{\ell_n,\ell'_n})\mathbf{1}_{\{R_{\ell_1,\ell'_1} \preccurlyeq_1 \cdots \preccurlyeq_{k-1} R_{\ell_k,\ell'_k} \preccurlyeq'_k R_{1,L} \preccurlyeq''_k R_{\ell_{k+1},\ell'_{k+1}} \preccurlyeq_{k+1} \cdots \preccurlyeq_{n-1} R_{\ell_n,\ell'_n}\}}\rangle,$$

for various choices of $k$ and of the order relations $\preccurlyeq'_k, \preccurlyeq''_k \in \{=,<\}$. Since this decomposition comes from a partitioning of the event $R_{1,L} < R_{\ell_n,\ell'_n}$, the number of equality signs suffixing the sequence

$$\preccurlyeq_1,\ldots,\preccurlyeq_{k-1},\preccurlyeq'_k,\preccurlyeq''_k,\preccurlyeq_{k+1},\ldots,\preccurlyeq_{n-1}$$

is the same as the number of equality signs suffixing the sequence $\preccurlyeq_1,\ldots,\preccurlyeq_{n-1}$, which is indeed one less than in the sequence $\preccurlyeq_1,\ldots,\preccurlyeq_n$ we started with. This shows that we can progressively bring this number to zero, while otherwise only summing other terms that we can already compute. This completes the induction argument. ∎

The fact that when $G$ satisfies the Ghirlanda-Guerra identities, the law of $R_{1,2}$ must be supported on $\mathbb{R}_{\geqslant 0}$ is shown in Exercise 5.13. Theorem 5.29 is now an immediate consequence of Theorem 5.30 and Proposition 5.31.

As was discussed below the statement of Theorem 5.29, when the measure $\zeta$ in (5.98) has finite support, the overlap arrays associated with Poisson-Dirichlet cascades are the canonical representatives of random overlap arrays that satisfy the Ghirlanda-Guerra identities. It is desirable to extend this construction and be able to also discuss such canonical representatives for more general choices of the measure $\zeta$. This can easily be achieved by approximation, up to the extraction of a subsequence. In fact, it follows from Proposition 5.31 that no subsequence extraction is necessary. In particular, the law of the overlap array associated with a Poisson-Dirichlet cascade is a continuous function of $\zeta$.

**Corollary 5.32.** *Let $M > 0$, let $(\zeta_n)_{n \geqslant 1}$ be a sequence of probability measures with finite support in $[0,M]$ that converge weakly to a probability measure $\zeta$. For each* $n \geq 1$, let $R^{(n)} = (R^{(n)}_{\ell,\ell'})_{\ell \neq \ell' \geq 1}$ denote the overlap array associated with a Poisson-Dirichlet cascade such that the law of of $R^{(n)}_{1,2}$ is $\zeta_n$. Then $(R^{(n)})_{n \geq 1}$ converges in law in the sense of finite-dimensional distributions.

*Proof.* Since each entry of $R^{(n)}$ takes values in $[0,M]$, it is clear that we can extract subsequentially convergent limits. Let $R$ denote one such limit array. By construction, the overlap $R$ satisfies the Ghirlanda-Guerra identities (5.99) for continuous functions $f$ and $\psi$. Using for instance Exercise A.8, we see that this property guarantees that the conditional law of $R_{1,n+1}$ given $(R_{\ell,\ell'})_{\ell \neq \ell' \leq n}$ is the law displayed in (5.100), and so $R$ satisfies the Ghirlanda-Guerra identities (5.99) for arbitrary bounded measurable functions $f$ and $\psi$. By the Portmanteau theorem (Theorem A.17), we also see that $R$ is ultrametric, and that the law of $R_{1,2}$ is $\zeta$. By Proposition 5.31, the law of $R$ is uniquely determined by these properties, so the proof is complete. ∎

If desired, we can apply the Dovbysh-Sudakov theorem (Theorem 1.7 in [211]) to obtain that the limit overlap array in Corollary 5.32 can be represented in the form of $(\sigma^\ell \cdot \sigma^{\ell'})_{\ell \neq \ell'}$, where $(\sigma^\ell)_{\ell \geq 1}$ are independent samples from a random probability measure $G$ on some Hilbert space.

**Exercise 5.13** (Talagrand's positivity principle). Under the setting explained in the paragraph above (5.98), we assume that the random probability measure $G$ satisfies the Ghirlanda-Guerra identities. Our goal is to show that the law $\zeta$ of one overlap, see (5.98), must be supported in $\mathbb{R}_{\geq 0}$.

(i) Let $A$ be a measurable subset of $\mathbb{R}$, let $n \geq 1$ be an integer, and let $A_n$ denote the event that for every $\ell \neq \ell' \leq n$, we have $R_{\ell,\ell'} \in A$. Show that

$$\mathbf{1}_{A_{n+1}} \geq \mathbf{1}_{A_n}\left(1 - \sum_{\ell=1}^n \mathbf{1}_{\{R_{\ell,n+1} \notin A\}}\right). \tag{5.111}$$

(ii) Show that for every $n \geq 1$, we have

$$\mathbb{E}\langle \mathbf{1}_{A_n}\rangle \geq (\zeta(A))^{n-1}. \tag{5.112}$$

(iii) Conclude that the support of $\zeta$ must be a subset of $\mathbb{R}_{\geq 0}$.

# Chapter 6
# Mean-field spin glasses

In this chapter, we enter more firmly into the realm of mean-field spin glasses. One of the most basic models in this class is the Sherrington-Kirkpatrick (SK) model. Like the Curie-Weiss model and the models of statistical inference studied in previous chapters, but unlike the Ising model on $\mathbb{Z}^d$, it is a mean-field model, in the sense that the variables that compose it are mutually exchangeable. And like the models of statistical inference, but unlike the Ising and Curie-Weiss models, the Hamiltonian defining the Gibbs measure of the SK model is random; in other words, the interactions between the variables are disordered. In the late 1970's, Giorgio Parisi proposed a self-contained description for the limit free energy of the SK model, now known as the Parisi formula, using sophisticated non-rigorous arguments [217, 218, 219]. This progressively led to important progress on the theoretical understanding of these models [177, 178, 179, 220], at the physics level of rigour. A mathematical proof of the Parisi formula was then given in the early 2000's in [132, 250]. A more robust proof was later found in [210, 212] using the idea of ultrametricity discussed in the previous section, and also inspired by [11, 13, 21, 57, 128, 207, 208, 209].

As was already touched upon briefly in the overview to this book, there are many seemingly modest generalizations of the SK model, such as the bipartite model described in (0.4) and represented in Figure 0.2, that remain poorly understood mathematically. An important motivation for writing this book is to lay some groundwork for making progress on this problem.

We will not provide a proof of the Parisi formula in this book. Rather, our goal is to explain how the Parisi formula relates to the Hamilton-Jacobi approach, and more specifically, how it can be interpreted as the solution to an infinite-dimensional Hamilton-Jacobi equation evaluated at a point. In Section 6.1, we describe the SK model and state the Parisi formula. In Section 6.2, we explore a naive attempt at enriching the free energy in the SK model by mimicking what we did in the context of statistical inference. This enriched free energy satisfies an approximate

Hamilton-Jacobi equation, with an error term that is the variance of the overlap between two replicas. We encountered a similar situation in statistical inference, and in this context, we could leverage the Nishimori identity to ensure that the error term vanishes in the limit of large system size. Unfortunately, this is not the case at low temperature in spin-glass models such as the SK model, and we will need to refine our approach. In Section 6.3, we briefly discuss a simplification of the SK model known as the random energy model, in the hope of gaining some insight into the appropriate way of enriching the free energy. This will hint at a potential role for the Poisson-Dirichlet cascades already discussed in Section 5.6. In Section 6.4, we use this object to define a more sophisticated enriched free energy, and argue informally that this object should satisfy an infinite-dimensional Hamilton-Jacobi equation, up to an error term that is expected to vanish in the limit of large system size. Moreover, the non-linearity in the equation associated with the SK model is convex, so the Hopf-Lax formula suggests a variational representation for the limit free energy. In Section 6.5, we connect this Hopf-Lax formula with the Parisi formula. Taking the validity of the Parisi formula for granted, this gives some substance to the heuristic derivations performed earlier. One of the advantages of the Hamilton-Jacobi approach is that it leads to the formulation of a natural conjecture for the limit free energy of more difficult models such as the bipartite model, for which variational formulas break down. In Section 6.6, we conclude this chapter by discussing the state of current research on this point.

## 6.1   The Sherrington-Kirkpatrick model and the Parisi formula

The SK model is often motivated as a pure optimization problem known as the *Dean's problem*. We imagine that the Dean of a university is tasked with dividing a group of $N$ students indexed by the elements of $\{1,\ldots,N\}$ into two dormitories indexed by $-1$ and $+1$. An allocation of the $N$ students into the two dormitories can be identified with a *configuration vector*

$$\sigma = (\sigma_1,\ldots,\sigma_N) \in \Sigma_N := \{-1,+1\}^N. \tag{6.1}$$

To decide on dormitory allocations, the Dean is provided with a collection of parameters $(g_{ij})_{i,j \leqslant N}$ called *interaction parameters* which describe how much student $i$ likes or dislikes student $j$. A positive parameter means that student $i$ likes student $j$, while a negative parameter means that student $i$ dislikes student $j$. The Dean strives to maximize the "comfort function"

$$c_N(\sigma) := \sum_{i,j=1}^{N} g_{ij}\sigma_i\sigma_j. \tag{6.2}$$

We define the comfort function in this way so that it most resembles an Ising-type model; the possibly more natural choice of $\sum_{i,j} g_{ij}\mathbf{1}_{\{\sigma_i=\sigma_j\}} = \frac{1}{2}c_N(\sigma) + \frac{1}{2}\sum_{i,j}g_{ij}$ would not change the analysis much. The story of the Dean is a simple illustration that "complex" optimization problems of this form arise quite naturally. For arbitrary choices of the coefficients $(g_{ij})_{i,j \leqslant N}$, finding near-optimizers to the function $c_N$ is known to be NP-hard [22]. Instead of this worst-case analysis, we would like to consider the problem of optimizing the function $c_N$ for a "typical" instantiation. We encode this by positing that the interaction parameters $(g_{ij})_{i,j \leqslant N}$ are independent standard Gaussian random variables. The key assumption here is that the random variables are independent and identically distributed; arguing as in the proof of Theorem 4.24 allows us to extend many results to a much broader class of distributions, see for instance Exercise 6.4. The Hamiltonian of the SK model is defined as the comfort function $c_N(\sigma)$, suitably normalized, and for this choice of independent standard Gaussian interaction parameters $(g_{ij})_{i,j \leqslant N}$. Precisely, for every $\sigma \in \Sigma_N$, we set

$$H_N(\sigma) := \frac{1}{\sqrt{N}} \sum_{i,j=1}^{N} g_{ij} \sigma_i \sigma_j. \tag{6.3}$$

We may as well think of $H_N$ as a function defined on $\mathbb{R}^N$. The family $(H_N(\sigma))_{\sigma \in \mathbb{R}^N}$ is a centred Gaussian process with covariance structure

$$\mathbb{E} H_N(\sigma^1) H_N(\sigma^2) = \frac{1}{N} \sum_{i,j=1}^{N} \sigma_i^1 \sigma_j^1 \sigma_i^2 \sigma_j^2 = N \left( \frac{1}{N} \sum_{i=1}^{N} \sigma_i^1 \sigma_i^2 \right)^2 = N R_{1,2}^2 \tag{6.4}$$

that depends on the spin configurations $\sigma^1$ and $\sigma^2$ only through their normalized scalar product, or *overlap*,

$$R_{1,2} := \frac{1}{N} \sigma^1 \cdot \sigma^2 = \frac{1}{N} \sum_{i=1}^{N} \sigma_i^1 \sigma_i^2. \tag{6.5}$$

Notice that we could also have directly defined $(H_N(\sigma))_{\sigma \in \mathbb{R}^N}$ by specifying that it is a centred Gaussian process whose covariance is given by (6.4); the main utility of the explicit formula (6.3), besides its possibly more intuitive appeal, is that it guarantees the existence of such a process. The law of $(H_N(\sigma))_{\sigma \in \mathbb{R}^N}$ is invariant under permutations of the coordinates of $\sigma$. This is why the SK model is said to be *mean-field*. Instead of focusing on analyzing the maximum

$$\frac{1}{N} \max_{\sigma \in \Sigma_N} H_N(\sigma), \tag{6.6}$$

and consistently with the rest of the book, we prefer to analyze the free energy

$$\overline{F}_N(\beta) := \frac{1}{N} \mathbb{E} \log \sum_{\sigma \in \Sigma_N} \exp(\beta H_N(\sigma)) \tag{6.7}$$

at each inverse temperature $\beta > 0$. The reason for the normalization by $\sqrt{N}$ in (6.3) is that the configurations $\sigma$ that contribute meaningfully to the sum in (6.7) are such that $H_N(\sigma)$ is of order $N$ (see for instance Exercise 6.1), so that the exponential can be in balance with the cardinality of $\Sigma_N$. As we vary the parameter $\beta$, the sum will therefore interpolate between a small-$\beta$ situation in which the main contribution consists of the very numerous terms that have $H_N(\sigma)$ not very large, and a large-$\beta$ situation where the sum is dominated by the few terms for which $H_N(\sigma)$ is very close to the maximum. For large $\beta$, this is made more precise in Exercise 6.3.

A first formula for the limit of the free energy $\overline{F}_N(\beta)$ was proposed by Sherrington and Kirkpatrick in [237]. However, in this same work, they observed that their formula predicted that the entropy of the limit Gibbs measure would become negative for large $\beta$, but we recall from (1.7) that the entropy is always non-negative, so their proposed solution could not be valid for large $\beta$. A few years later, Giorgio Parisi [217, 218, 219] proposed a more sophisticated formula now known as the *Parisi formula*, which was then established rigorously in [132, 250] (see also [211]). It states that

$$\lim_{N\to+\infty} \overline{F}_N(\beta) = \inf_{\zeta \in \mathcal{D}[0,1]} \left( \Phi_\zeta(0,0) - \beta^2 \int_0^1 t\zeta(t)\,\mathrm{d}t + \log(2) \right), \tag{6.8}$$

where $\Phi_\zeta : [0,1] \times \mathbb{R} \to \mathbb{R}$ is the solution to the parabolic equation

$$\begin{cases} -\partial_t \Phi_\zeta(t,x) = \beta^2 \left( \partial_x^2 \Phi_\zeta(t,x) + \zeta(t)\left( \partial_x \Phi_\zeta(t,x) \right)^2 \right) & \text{on } [0,1] \times \mathbb{R}^d \\ \Phi_\zeta(1,x) = \log\cosh(x) & \text{for } x \in \mathbb{R} \end{cases} \tag{6.9}$$

solved backwards in time from $t = 1$ to $t = 0$, and where the infimum in (6.8) is taken over the space

$$\mathcal{D}[0,1] := \{\zeta : [0,1] \to [0,1] \mid \zeta \text{ is right-continuous and}$$

$$\text{non-decreasing with } \zeta(1) = 1\} \tag{6.10}$$

of probability distribution functions on $[0,1]$. The equation (6.9) is sometimes called the *Parisi PDE*. One of the main goals of this chapter is to recast this statement using the point of view provided by the Hamilton-Jacobi approach. Before doing this, we mention two directions in which the SK model can be generalized.

On the one hand, the covariance structure (6.4) can be generalized in the following way. Given a sequence $(\beta_p)_{p\geq 1}$ of non-negative numbers that decrease sufficiently rapidly, the Hamiltonian

$$H_N(\sigma) := \sum_{p \geq 1} \beta_p^{1/2} H_{N,p}(\sigma) \tag{6.11}$$

of the *mixed p-spin model* on $\Sigma_N$ is given by a linear combination of *pure p-spin* Hamiltonians on $\Sigma_N$,

$$H_{N,p}(\sigma) := \frac{1}{N^{(p-1)/2}} \sum_{i_1,\ldots,i_p=1}^N g_{i_1\ldots i_p} \sigma_{i_1} \cdots \sigma_{i_p}, \tag{6.12}$$

where $(g_{i_1\ldots i_p})_{i_1,\ldots,i_p \leqslant N}$ are independent standard Gaussian random variables, independent of each other as we vary $p$. The Hamiltonian (6.11) is a centred Gaussian process with covariance structure

$$\mathbb{E} H_N(\sigma^1) H_N(\sigma^2) = N\xi(R_{1,2}) \tag{6.13}$$

that depends on the spin configurations $\sigma^1$ and $\sigma^2$ only through their overlap $R_{1,2}$, and where we have set

$$\xi(x) := \sum_{p \geqslant 1} \beta_p x^p. \tag{6.14}$$

The SK model corresponds to the choice $\beta_p = \mathbf{1}_{\{p=2\}}$, and so $\xi(x) = x^2$. More than the explicit form of $H_N(\sigma)$ given in (6.11) and (6.12), it is the function $\xi$ characterizing the covariance structure of the Hamiltonian that will play a central role. In particular, this function is essentially the non-linearity in the Hamilton-Jacobi equation we will need to consider. We refer to the function $\xi$ as the *covariance function*. The limit of the free energy

$$\overline{F}_N(\beta) := \frac{1}{N} \mathbb{E} \log \sum_{\sigma \in \Sigma_N} \exp(\beta H_N(\sigma)) \tag{6.15}$$

in the mixed $p$-spin models admits a variational representation analogous to (6.8). The connection between the SK model and the Hamilton-Jacobi approach that will be presented in this chapter extends to this more general setting [196].

Another possible way to generalize the SK model consists in considering that each coordinate of the configuration vector (6.1) can itself be a vector; the bipartite model is an example of such a setup. For illustration, we briefly motivate the bipartite model through a slightly modified version of the Dean's problem.

Suppose the Dean of a university actually needs to split students and tutors into two groups. We assume that there is an equal number $N$ of students and tutors for convenience, but this is inessential. An allocation of the $2N$ students and tutors into the two groups $+1$ and $-1$ can be identified with a vector

$$\sigma = (\sigma_1, \sigma_2) = \left((\sigma_{1,1}, \sigma_{1,2}, \ldots, \sigma_{1,N}), (\sigma_{2,1}, \sigma_{2,2}, \ldots, \sigma_{2,N})\right) \in \Sigma_N^2, \tag{6.16}$$

where $\sigma_1 \in \Sigma_N$ encodes the group allocations for the students and $\sigma_2 \in \Sigma_N$ the group allocations for the tutors. We postulate that the Dean now only has access to a collection of parameters $(g_{ij})_{i,j \leqslant N}$, where $g_{ij}$ encodes the quality of the interaction between student $i$ and tutor $j$. The Dean aims to maximize the comfort function

$$c_N(\sigma) := \sum_{i,j=1}^{N} g_{ij} \sigma_{1,i} \sigma_{2,j} \tag{6.17}$$

over all $\sigma \in \Sigma_N^2$. Choosing the interaction parameters $(g_{ij})_{i,j \leqslant N}$ to be independent standard Gaussian random variables and re-scaling $c_N$, we define the *bipartite model*

[119, 120, 152, 156] whose Hamiltonian on $\Sigma_N^2$ is given by

$$H_N(\sigma) := \frac{1}{\sqrt{N}} \sum_{i,j=1}^N g_{ij}\sigma_{1,i}\sigma_{2,j}. \tag{6.18}$$

The Hamiltonian (6.18) is a centred Gaussian process with covariance structure

$$\mathbb{E}H_N(\sigma^1)H_N(\sigma^2) = N\left(\frac{\sigma_1^1 \cdot \sigma_1^2}{N}\right)\left(\frac{\sigma_2^1 \cdot \sigma_2^2}{N}\right). \tag{6.19}$$

So in this model, the relevant covariance function takes the two "intra-community" overlaps as arguments, and is defined on $\mathbb{R}^2$ by

$$\xi(x,y) := xy. \tag{6.20}$$

The fundamental difference between the bipartite model and the mixed $p$-spin models is that the covariance function $\xi$ is always convex for the latter, but not for the former. To be precise, for scalar (as opposed to vector) models, for which the function $\xi$ is defined on $\mathbb{R}$, the relevant question is whether the function is convex on $\mathbb{R}_{\geq 0}$, and the function in (6.14) satisfies this convexity property for every choice of the non-negative parameters $(\beta_p)_{p\geq 1}$. Due to the absence of a comparable convexity property for the bipartite model, the limit of the free energy

$$\frac{1}{N}\mathbb{E}\log\sum_{\sigma\in\Sigma_N^2}\exp(\beta H_N(\sigma)) \tag{6.21}$$

with $H_N$ as in (6.18) has not yet been identified. As will be explained below, the Hamilton-Jacobi approach suggests a natural conjecture for what this limit could be. Similar considerations apply to essentially arbitrary vector spin-glass models, provided that the covariance of the Hamiltonian can be expressed as a function of the matrix of overlaps between the different types of coordinates. The dividing line between the models for which the limit free energy has or has not been identified is exactly whether or not the underlying covariance function $\xi$ is convex. In the most general models in which each spin configuration $\sigma$ can be decomposed into a fixed number $D$ of different components $\sigma = (\sigma_1, \ldots, \sigma_D)$, the covariance function takes the matrix of all possible overlaps $N^{-1}(\sigma_d^1 \cdot \sigma_{d'}^2)_{d,d'\leq D}$ as input, and the precise question is whether the function $\xi$ is convex over the space of $D$-by-$D$ positive semi-definite matrices.

**Exercise 6.1.** For the Hamiltonian $H_N$ defined in (6.3), show that there exists a constant $C < +\infty$ such that

$$C^{-1}N \leq \mathbb{E}\max_{\sigma\in\Sigma_N} H_N(\sigma) \leq CN. \tag{6.22}$$

**Exercise 6.2.** For the Hamiltonian $H_N$ defined in (6.3), prove that almost surely,

$$\lim_{N\to+\infty}\left|\frac{1}{N}\max_{\sigma\in\Sigma_N}H_N(\sigma) - \frac{1}{N}\mathbb{E}\max_{\sigma\in\Sigma_N}H_N(\sigma)\right| = 0. \tag{6.23}$$

**Exercise 6.3.** For the Hamiltonian $H_N$ defined in (6.3) and the free energy $\overline{F}_N(\beta)$ defined in (6.7), show that

$$\frac{1}{N}\mathbb{E}\max_{\sigma\in\Sigma_N}H_N(\sigma) \leqslant \overline{F}_N(\beta) \leqslant \frac{1}{N}\mathbb{E}\max_{\sigma\in\Sigma_N}H_N(\sigma) + \frac{\log(2)}{\beta}. \tag{6.24}$$

**Exercise 6.4.** Fix a collection $(x_{ij})_{i,j\leqslant N}$ of i.i.d. centred random variables with variance one, and write

$$H_N^x(\sigma) := \frac{1}{\sqrt{N}}\sum_{i,j\leqslant N}x_{ij}\sigma_i\sigma_j \quad\text{and}\quad \overline{F}_N^x(\beta) := \frac{1}{N}\mathbb{E}\log\sum_{\sigma\in\Sigma_N}\exp\bigl(\beta H_N^x(\sigma)\bigr) \tag{6.25}$$

for their associated SK Hamiltonian and free energy. Denote by $\overline{F}_N(\beta)$ the Gaussian SK free energy (6.7). Assuming that $\mathbb{E}|x_{11}|^3 < +\infty$, show that for every inverse temperature parameter $\beta > 0$,

$$\limsup_{N\to+\infty}\left|\overline{F}_N^x(\beta) - \overline{F}_N(\beta)\right| = 0. \tag{6.26}$$

## 6.2 A first attempt at a Hamilton-Jacobi approach to the SK model

We now return to the setting of the SK model. We prefer to slightly change notation and thereby slightly generalize the model in the following way. We fix a probability measure $P_1$ on $\Sigma_1$, and for each integer $N \geqslant 1$, let

$$P_N := (P_1)^{\otimes N} \tag{6.27}$$

denote its $N$-fold tensor product. Given a "time" parameter $t \geqslant 0$, we introduce the free energy

$$\overline{F}_N^\circ(t) := -\frac{1}{N}\mathbb{E}\log\int_{\Sigma_N}\exp\bigl(\sqrt{2t}H_N(\sigma) - Nt\bigr)\mathrm{d}P_N(\sigma), \tag{6.28}$$

where the Hamiltonian $H_N$ is that of the SK model, see (6.3). When $P_1$ is the uniform measure on $\Sigma_1$, the free energy (6.7) we had defined earlier is related to this new quantity according to

$$\overline{F}_N(\beta) = -\overline{F}_N^\circ\left(\frac{\beta^2}{2}\right) + \log(2) + \frac{\beta^2}{2}. \tag{6.29}$$

So we may as well focus on trying to compute $\overline{F}_N^\circ(t)$ for each value of $t \geq 0$. By (6.4), for every $\sigma \in \Sigma_N$, the random variable $H_N(\sigma)$ is a centred Gaussian with variance $N$, so

$$\mathbb{E}\exp\left(\sqrt{2t}H_N(\sigma) - Nt\right) = 1. \tag{6.30}$$

In the language of statistical physics, one could say that we have normalized the Hamiltonian so that the annealed free energy is always zero. Combining (6.30) with Jensen's inequality shows that $\overline{F}_N^\circ(t) \geq 0$; this is why we chose to include a negative sign in front of the free energy (6.28). Choosing to parameterize the inverse temperature by $\sqrt{2t}$ in place of $\beta$ is convenient because then the variance of the Gaussian field $\sqrt{2t}H_N(\sigma)$ scales linearly, as a Brownian motion would. As was the case in Chapter 4, see in particular the discussion around (4.93), this parametrization then satisfies a certain semigroup property. While we will not dwell much on this point, this choice of parametrization will allow us in particular to obtain expressions for the derivatives of the free energy that do not explicitly depend on $t$.

In the expression (6.28) for the free energy, we chose to use an integral symbol for what turns out to be just a sum over $\Sigma_N$, possibly with some suitable weights if the measure $P_1$ is not centred. The main reason we chose to denote it in this way is that the results we will discuss concerning $\overline{F}_N^\circ(t)$ can be generalized straightforwardly to the case where the measure $P_1$ is arbitrary with compact support, provided that we replace the term $-Nt$ in (6.28) by $-t|\sigma|^4/N$, or more generally by $-Nt\xi(N^{-1}|\sigma|^2)$ for the mixed $p$-spin models in (6.13). In other words, this term should be half the variance of $\sqrt{2t}H_N(\sigma)$, so that the identity in (6.30) generalizes nicely. If the limit of $\overline{F}_N^\circ(t)$ can be identified, then the possibly unwelcome extra term $-Nt\xi(N^{-1}|\sigma|^2)$ can later be removed by essentially quoting Proposition 3.20 on the generalized Curie-Weiss model. This is explained using a slightly different argument in Section 5 of [198]. Here we will stick to the case in which $P_1$ is supported on $\Sigma_1$, in which case this additional term is constant and therefore does not cause any difficulty.

As in Section 4.3 on statistical inference, there is no hope of finding a partial differential equation solved by the function $\overline{F}_N^\circ(t)$ in (6.28) on its own. Rather, we need to enrich the Hamiltonian to contain more terms, so that the free energy then depends on additional parameters. The goal then is to find additional terms that are sufficiently simple that we can analyze them when setting $t = 0$; but that are sufficiently expressive that we can compensate any small variation in $t$ by small variations in the other parameters, thereby obtaining a Hamilton-Jacobi equation. To get a first rough intuition for what one could do, we start by rewriting the Hamiltonian $H_N(\sigma)$ for the SK model as

$$\sum_{i=1}^N \left(\frac{1}{\sqrt{N}} \sum_{j=1}^N g_{i,j}\sigma_j\right)\sigma_i. \tag{6.31}$$

Inspired by the analysis of the previous chapters, we would like to try to only add terms to the Hamiltonian that do not encode any interaction between the variables;

in other words, linear functions of $\sigma$ would be best. We would like to add a term that would at least somewhat resemble (6.31); since the term $\frac{1}{\sqrt{N}} \sum_{j=1}^{N} g_{i,j}\sigma_j$ is random, we may try adding a new term in the Hamiltonian proportional to $\sum_{i=1}^{N} z_i\sigma_i$, where $z = (z_1, \ldots, z_N) \in \mathbb{R}^N$ is a standard Gaussian vector.

The brief informal argument we just gave is at best a first hint for what one could try to do. We will not attempt to make it more precise here, and instead just proceed with the suggested idea and see how well it fares. So we fix a standard Gaussian vector $z = (z_1, \ldots, z_N) \in \mathbb{R}^N$ independent of $(H_N(\sigma))_{\sigma \in \mathbb{R}^N}$, and for each $h \geqslant 0$, we consider the enriched free energy

$$\overline{F}_N(t,h) := -\frac{1}{N}\mathbb{E}\log\int_{\Sigma_N} \exp\left(\sqrt{2t}H_N(\sigma) - Nt + \sqrt{2h}z\cdot\sigma - Nh\right) dP_N(\sigma). \quad (6.32)$$

We have kept the practice of adding the compensating term $-Nh$, i.e. half the variance of $\sqrt{2h}z\cdot\sigma$, in the exponential. For measures $P_N$ that are not necessarily supported on $\Sigma_N$, replacing this term by $-h|\sigma|^2$ (and $-Nt$ by $-Nt\xi(N^{-1}|\sigma|^2)$) ensures that all the calculations below remain valid. As usual, we write $\langle\cdot\rangle$ for the Gibbs average associated with the Hamiltonian

$$H_N(t,h,\sigma) := \sqrt{2t}H_N(\sigma) - Nt + \sqrt{2h}z\cdot\sigma - Nh, \quad (6.33)$$

we denote by $\sigma$ the canonical random variable under $\langle\cdot\rangle$, and by $(\sigma^\ell)_{\ell\geqslant 1}$ independent copies of $\sigma$ under $\langle\cdot\rangle$. We stress that the probability measure $\langle\cdot\rangle$ is random, since it depends on the realization of $(H_N(\sigma))_{\sigma\in\mathbb{R}^N}$ (or equivalently, of $(g_{i,j})_{1\leqslant i,j\leqslant N}$) and $z$. Leveraging the Gibbs Gaussian integration by parts formula (Theorem 4.6), it is readily verified that

$$\partial_h \overline{F}_N(t,h) = -\frac{1}{N}\mathbb{E}\left\langle\frac{1}{\sqrt{2h}}z\cdot\sigma - N\right\rangle = \mathbb{E}\langle R_{1,2}\rangle, \quad (6.34)$$

while

$$\partial_t \overline{F}_N(t,h) = -\frac{1}{N}\mathbb{E}\left\langle\frac{1}{\sqrt{2t}}H_N(\sigma) - N\right\rangle = \mathbb{E}\langle R_{1,2}^2\rangle. \quad (6.35)$$

This implies that

$$\partial_t \overline{F}_N(t,h) - \left(\partial_h \overline{F}_N(t,h)\right)^2 = \mathrm{Var}(R_{1,2}), \quad (6.36)$$

where we write Var to denote the variance over all sources of randomness, i.e. with respect to the measure $\mathbb{E}\langle\cdot\rangle$. The formula (6.36) looks very similar to what we found in (3.5) for the Curie-Weiss model and in (4.83) for symmetric rank-one matrix estimation, and suggests that the limit $f$ of the enriched free energy $\overline{F}_N$ should satisfy the Hamilton-Jacobi equation

$$\partial_t f - \left(\partial_h f\right)^2 = 0 \quad \text{on} \quad \mathbb{R}_{>0} \times \mathbb{R}_{>0}. \quad (6.37)$$

Unfortunately, unlike previous models, in the spin-glass setting this guess is too naive. More precisely, this guess can be shown to be correct for small values of $t \geqslant 0$, where it can be used to recover the original formula of Sherrington and Kirkpatrick; but it fails for large values of $t \geqslant 0$, where the variance of the overlap $R_{1,2}$ no longer vanishes in the large $N$ limit. In the language of statistical physics, the SK model exhibits "replica symmetry breaking" at low temperature. Physicists use the phrase "replica symmetry breaking" because if we sample four replicas $\sigma^1, \ldots, \sigma^4$ from the Gibbs measure, one could expect by symmetry to witness the same overlap value between $\sigma^1$ and $\sigma^2$ and between $\sigma^3$ and $\sigma^4$. The symmetry between the replicas is not literally broken, but indeed these overlap values will in general be different because the overlap between two replicas is random, and remains so even in the limit of large system size.

To sum up, the consideration of the enriched free energy in (6.32) seems to be pointing in the right direction, but is insufficient, and we therefore need a more refined way of enriching the free energy.

## 6.3   Some insight from the random energy model

To gain a deeper understanding of the structure of the Gibbs measure and determine how to appropriately enrich the free energy (6.28), we take a small detour and analyze a simplified model called the random energy model, or REM [94]. The Hamiltonian $\sqrt{2t}H_N(\sigma)$ appearing in the free energy (6.28) is a family of Gaussian random variables indexed by $\sigma \in \Sigma_N$ with variance $2tN$. The main difficulty in understanding the behaviour of this Hamiltonian is that these Gaussian random variables are correlated, by (6.4). If we consider pure $p$-spin models instead, then the covariance between $H_N(\sigma^1)$ and $H_N(\sigma^2)$ becomes $NR_{1,2}^p$. For $\sigma^1, \sigma^2 \in \Sigma_N$ and in the limit of large $p$, this would simplify to $N\mathbf{1}_{\{R_{1,2}=1\}}$. That is to say, in this limit, the random variables $(H_N(\sigma))_{\sigma \in \Sigma_N}$ have become independent. The random energy model corresponds to setting exactly this as the definition: for a family $(E_\sigma)_{\sigma \in \Sigma_N}$ of independent standard Gaussian random variables, $t \geqslant 0$ and $\sigma \in \Sigma_N$, we set

$$H_N^{\text{REM}}(t, \sigma) := \sqrt{2tN}E_\sigma. \tag{6.38}$$

Our goal in this section is to inquire about the structure of the associated Gibbs measure,

$$G_N^{\text{REM}}(t, \sigma) = \frac{\exp H_N^{\text{REM}}(t, \sigma)}{\sum_{\tau \in \Sigma_N} \exp H_N^{\text{REM}}(t, \tau)}, \tag{6.39}$$

especially for large values of $t$. The hope is that the Gibbs measure of "true" spin glasses behaves similarly, possibly after clumping together neighbouring configurations in some way so as to tame the effect of correlations.

Since we are dealing with a family of independent standard Gaussian random variables, a natural starting point is to appeal to Proposition 5.16. Substituting $n$ by $2^N$ there, we introduce the constant

$$a_N := \left(2N\log(2) - \log(N) - \log\log(2) - \log(4\pi)\right)^{\frac{1}{2}}, \tag{6.40}$$

so that Proposition 5.16 ensures that the point process

$$\sum_{\sigma \in \Sigma_N} \delta_{a_N(E_\sigma - a_N)} = \sum_{\sigma \in \Sigma_N} \delta_{\frac{a_N}{\sqrt{2tN}}(H_N^{\mathrm{REM}}(t,\sigma) - \sqrt{2tN}a_N)} \tag{6.41}$$

converges in law on $\mathbb{R}$ to a Poisson point process with intensity $d\nu(x) := e^{-x}dx$. If we introduce the parameter

$$\zeta := \sqrt{\frac{\log(2)}{t}} \tag{6.42}$$

and observe that $\frac{a_N}{\sqrt{2tN}} = \zeta + o(1)$, this suggests that the point process

$$\Lambda_N := \sum_{\sigma \in \Sigma_N} \delta_{\zeta H_N^{\mathrm{REM}}(t,\sigma) - a_N^2} \tag{6.43}$$

should also converge in law on $\mathbb{R}$ to a Poisson point process with intensity $\nu$. Indeed, the lower-order terms should not impact the convergence. Perhaps the simplest way to justify this is to essentially repeat the proof of Proposition 5.16.

**Proposition 6.1.** *The point process $\Lambda_N$ in (6.43) converges in law on $\mathbb{R}$ to a Poisson point process with intensity $d\nu(x) := e^{-x}dx$.*

*Proof.* Fix a function $f \in C_c(\mathbb{R}; \mathbb{R}_{\geq 0})$ of compact support, and observe that

$$\mathbb{E}\exp\left(-\int_\mathbb{R} f\,d\Lambda_N\right)$$
$$= \exp\left[2^N \log\left(1 - \frac{1}{\sqrt{4\pi tN\zeta^2}} \int_\mathbb{R} (1 - e^{-f(x)}) \exp\left(-\frac{(x+a_N^2)^2}{4tN\zeta^2}\right)dx\right)\right].$$

Using the definition of the constant $a_N$ in (6.40), the asymptotic expansion

$$\frac{a_N^2}{2tN\zeta^2} = 1 + o(1),$$

and the fact that $f$ is of compact support reveals that uniformly over $x$ in the support of $f$,

$$\frac{2^N}{\sqrt{4\pi tN\zeta^2}} \exp\left(-\frac{(x+a_N^2)^2}{4tN\zeta^2}\right) = \frac{2^N}{\sqrt{4\pi tN\zeta^2}} \exp\left(-x - \frac{1}{2}a_N^2\right) + o(1) = e^{-x} + o(1).$$

It follows by a Taylor expansion of the logarithm that

$$\mathbb{E}\exp\left(-\int_\mathbb{R} f\,d\Lambda_N\right) = \exp\left(-\int_\mathbb{R} (1 - e^{-f(x)}) e^{-x}\,dx\right) + o(1).$$

Invoking Propositions 5.4 and 5.9 completes the proof. ∎

Together with the mapping theorem (Proposition 5.10) applied to the function $f(x) := e^{x/\zeta}$, this result implies that the point process

$$\sum_{\sigma \in \Sigma_N} \delta_{\exp(H_N^{\text{REM}}(t,\sigma) - a_N^2/\zeta)} \tag{6.44}$$

converges in law on $\mathbb{R}_{>0}$ to the Poisson point process with intensity measure

$$d\mu(x) := \frac{\zeta}{x^{\zeta+1}} dx. \tag{6.45}$$

Notice that the Gibbs measure (6.39) is invariant under translations of the Hamiltonian (6.38) by a constant, since

$$G_N^{\text{REM}}(t,\sigma) = \frac{\exp H_N^{\text{REM}}(t,\sigma)}{\sum_{\tau \in \Sigma_N} \exp H_N^{\text{REM}}(t,\tau)} = \frac{\exp\left(H_N^{\text{REM}}(t,\sigma) - a_N^2/\zeta\right)}{\sum_{\tau \in \Sigma_N} \exp\left(H_N^{\text{REM}}(t,\tau) - a_N^2/\zeta\right)}. \tag{6.46}$$

In view also of Exercise 5.8, this suggests that the point process supported on the Gibbs weights of the REM should converge in law to the re-normalization of the Poisson point process with intensity measure (6.45). More precisely, if we assume that $\zeta \in (0,1)$ so that the Poisson-Dirichlet process $(v_n)_{n \geq 1}$ is well-defined, then as $N$ tends to infinity we should expect that

$$\sum_{\sigma \in \Sigma_N} \delta_{G_N^{\text{REM}}(t,\sigma)} \text{ converges in law to } \sum_{n=1}^{+\infty} \delta_{v_n}. \tag{6.47}$$

Observe that $\zeta \in (0,1)$ if and only if $t > \log(2)$. It is possible to show that in this regime, the convergence in law (6.47) does indeed take place. We therefore expect the Poisson-Dirichlet process to play an important role in the enrichment of the free energy (6.28).

## 6.4  A Hamilton-Jacobi approach to the SK model

In the previous section, we suggested that the weights of the Gibbs measure at low temperature, perhaps after some suitable grouping-together of neighbouring sites, might resemble the weights of the Poisson-Dirichlet process introduced in Section 5.5. Although we will not aim to justify this precisely, this intuition turns out to be valid. Recall that our ultimate goal is to find a more refined random magnetic field to replace $\sqrt{2h}z$ in (6.32) so that we can find relationships between derivatives with respect to the various parameters. Recall also from (6.35) that the $t$-derivative of the free energy is the second moment of the overlap $R_{1,2}$, and that using the naive random magnetic field $\sqrt{2h}z$, we could create the relation (6.36) which replaces the second moment of $R_{1,2}$ by its variance. In other words, in this case we would be done if $R_{1,2}$ was constant. Now, if we were to enrich the free energy with a random magnetic field that has the overlap of a Poisson-Dirichlet process, then this would still not suffice in general. Indeed, we see from (5.56) that the overlap of a Poisson-Dirichlet process can only take two values, so armed with this, all we could possibly hope for would be to compensate for circumstances in which $R_{1,2}$ itself is known to take only two values. Compared with the naive attempt in (6.32), this seems to mark some further progress, but it is still insufficient in general. Aspects of the geometry of the Gibbs measure are important, beyond the simple recording of the statistics of the Gibbs weights.

Our work in Sections 5.6 and 5.7 suggests a clear way forward: Poisson-Dirichlet cascades will provide us with the richest possible class of overlap structures we can expect to encounter.

This is admittedly still only a relatively vague intuition. We will not push much towards giving it more substance, and mostly stick again to the attitude of just trying it out and seeing that things now actually work well. For the curious and motivated reader, we point out however the following possibly interesting route to explore. Below (6.31), we gestured at the idea that, since the term $\frac{1}{\sqrt{N}}\sum_{j=1}^{N} g_{i,j}\sigma_j$ is random, it might make sense to add a new term in the Hamiltonian proportional to $\sum_{i=1}^{N} z_i\sigma_i$, where $z = (z_1, \ldots, z_N) \in \mathbb{R}^N$ is a standard Gaussian vector. To be more precise, what we are really concerned about is the law of $\frac{1}{\sqrt{N}}\sum_{j=1}^{N} g_{i,j}\sigma_j$ under the Gibbs measure for a "typical" but fixed realization of the random coefficients $(g_{i,j})_{i,j \leqslant N}$. Exercise 6.5 explores this point in a simple setting, and suggests that aspects of the structure of the Gibbs measure such as whether the overlap is concentrated must be taken into consideration here.

We now proceed to construct the refined magnetic field that will be used as a replacement for the simple field $\sqrt{2h}z$ appearing in (6.32). Given an integer $K \geqslant 1$ and parameters

$$0 = \zeta_0 < \zeta_1 < \cdots < \zeta_K < \zeta_{K+1} = 1, \tag{6.48}$$

$$0 = q_{-1} \leqslant q_0 < q_1 < \cdots < q_K < q_{K+1} = +\infty, \tag{6.49}$$

as in (5.67) and (5.82), we define the right-continuous non-decreasing path $\mathsf{q} : [0, 1) \to \mathbb{R}_{\geqslant 0}$ by

$$\mathsf{q} := \sum_{k=0}^{K} q_k \mathbf{1}_{[\zeta_k, \zeta_{k+1})}. \tag{6.50}$$

The reason for storing the families of parameters $(\zeta_k)_{1 \leqslant k \leqslant K}$ and $(q_k)_{0 \leqslant k \leqslant K}$ in the form of this path will become clear as we proceed. We recall that we denote by $\mathcal{A}$ the tree in (5.63), whose set of leaves is $\mathbb{N}^K$, that $\alpha_{|\ell}$ denotes the ancestor of $\alpha$ at depth $\ell$ as in (5.65), and that $\alpha \wedge \beta$ denotes the depth of the most common ancestor of the leaves $\alpha$ and $\beta$ as in (5.66). We give ourselves a family $(z_\alpha)_{\alpha \in \mathcal{A}}$ of independent standard Gaussian random vectors in $\mathbb{R}^N$ independent of all other sources of randomness, and set, for every $\alpha \in \mathbb{N}^K$,

$$Z_{\mathsf{q}}(\alpha) := \sum_{k=0}^{K} (2q_k - 2q_{k-1})^{1/2} z_{\alpha_{|k}}. \tag{6.51}$$

The family of random variables $(Z_{\mathsf{q}}(\alpha))_{\alpha \in \mathbb{N}^K}$ is Gaussian, and its covariance structure is given, for every $\alpha^1, \alpha^2 \in \mathbb{N}^K$, by

$$\mathbb{E} Z_{\mathsf{q}}(\alpha^1) Z_{\mathsf{q}}(\alpha^2) = q_{\alpha^1 \wedge \alpha^2}. \tag{6.52}$$

This will be our new random external field. We let $(v_\alpha)_{\alpha \in \mathbb{N}^K}$ be a Poisson-Dirichlet cascade with parameters $(\zeta_k)_{1 \leqslant k \leqslant K}$, independent of all other sources of randomness, and we pair these two objects in our definition of the enriched free energy as

$$\overline{F}_N(t, \mathsf{q}) := -\frac{1}{N} \mathbb{E} \log \int_{\Sigma_N} \sum_{\alpha \in \mathbb{N}^K} \exp\big(H_N(t, \mathsf{q}, \sigma, \alpha)\big) v_\alpha \, dP_N(\sigma), \tag{6.53}$$

for the enriched Hamiltonian

$$H_N(t, \mathsf{q}, \sigma, \alpha) := \sqrt{2t} H_N(\sigma) - Nt + Z_{\mathsf{q}}(\alpha) \cdot \sigma - N q_K. \tag{6.54}$$

We denote by $\langle \cdot \rangle$ the associated Gibbs measure, with canonical random variable $(\sigma, \alpha)$. Explicitly, this means that for every bounded measurable function $f : \mathbb{R}^N \times \mathbb{N}^K \to \mathbb{R}$, we have

$$\langle f(\sigma, \alpha) \rangle := \frac{\int_{\Sigma_N} \sum_{\alpha \in \mathbb{N}^K} f(\sigma, \alpha) \exp\big(H_N(t, \mathsf{q}, \sigma, \alpha)\big) v_\alpha \, dP_N(\sigma)}{\int_{\Sigma_N} \sum_{\alpha \in \mathbb{N}^K} \exp\big(H_N(t, \mathsf{q}, \sigma, \alpha)\big) v_\alpha \, dP_N(\sigma)}. \tag{6.55}$$

We denote by $(\sigma^\ell, \alpha^\ell)_{\ell \geqslant 1}$ independent copies of $(\sigma, \alpha)$ under $\langle \cdot \rangle$, also called replicas. We may now discuss two types of overlaps: those that involve the $\sigma$ variables, which we may call the $\sigma$-overlaps and keep writing as

$$R_{1,2} := \frac{1}{N} \sum_{i \leqslant N} \sigma_i^1 \sigma_i^2; \tag{6.56}$$

and those that involve the $\alpha$ variables, which we may call the $\alpha$-overlaps and will simply denote by $\alpha^1 \wedge \alpha^2$.

Choosing $P_1$ to be the uniform measure on $\Sigma_1$, the free energy (6.7) in the SK model is related to this enriched free energy (6.53) by

$$\overline{F}_N(\beta) = -\overline{F}_N\left(\frac{\beta^2}{2}, 0\right) + \log(2) + \frac{\beta^2}{2}. \tag{6.57}$$

Using the Gaussian integration by parts formula, we now compute the derivatives of the enriched free energy with respect to $t$ and $(q_k)_{0 \leqslant k \leqslant K}$. We will find a relationship between these quantities, up to an error expressed in terms of the variance of the $\sigma$-overlap given the $\alpha$-overlap. The conditional expectation at play here is with respect to the measure $\mathbb{E}\langle\cdot\rangle$. For instance, we may write, for any function $f : \Sigma_N^2 \to \mathbb{R}$,

$$\mathbb{E}\langle f(\sigma^1,\sigma^2) \mid \alpha^1 \wedge \alpha^2 \rangle := \sum_{k=0}^{K} \mathbf{1}_{\{\alpha^1 \wedge \alpha^2 = k\}} \frac{\mathbb{E}\langle \mathbf{1}_{\{\alpha^1 \wedge \alpha^2 = k\}} f(\sigma^1,\sigma^2)\rangle}{\mathbb{E}\langle \mathbf{1}_{\{\alpha^1 \wedge \alpha^2 = k\}}\rangle}. \tag{6.58}$$

Before we proceed, we recall that in Theorem 5.28, specifically in (5.93), we managed to identify the law of the $\alpha$-overlap when the underlying Gibbs measure is simply the measure with weights $(v_\alpha)_{\alpha \in \mathbb{N}^K}$ on the leaves of the tree. One key observation is that the proof of this result carries verbatim to the present setting with a more complicated Gibbs measure. The point is that the Bolthausen-Sznitman invariance property of Poisson-Dirichlet cascades stated in Lemma 5.27 remains valid even if we multiply each $w_\alpha$ by some fixed function of $Z_{\mathsf{q}}(\alpha)$, since the family $(Z_{\mathsf{q}}(\alpha))_{\alpha \in \mathbb{N}^K}$ is independent of $(w_\alpha)$ and its law is invariant under permutations that preserve the tree structure. The other random field denoted by $(Z_q(\alpha))_{\alpha \in \mathbb{N}^K}$ there should be taken as independent from the random field $(Z_{\mathsf{q}}(\alpha))_{\alpha \in \mathbb{N}^K}$ we have defined here. In particular, we may replace each $w_\alpha$ by

$$w_\alpha \int_{\Sigma_N} \exp H_N(t,\mathsf{q},\sigma,\alpha)\,\mathrm{d}P_N(\sigma), \tag{6.59}$$

and then proceed through the proof of Theorem 5.28 to obtain that, with the Gibbs measure as in (6.55), we still have for every $k \in \{0,\ldots,K\}$ that

$$\mathbb{E}\langle \mathbf{1}_{\{\alpha^1 \wedge \alpha^2 = k\}}\rangle = \zeta_{k+1} - \zeta_k. \tag{6.60}$$

**Lemma 6.2.** *For every $t \geqslant 0$ and every path $\mathsf{q} : [0,1] \to \mathbb{R}$ of the form (6.50), we have*

$$\partial_t \overline{F}_N(t,\mathsf{q}) = \mathbb{E}\langle R_{1,2}^2 \rangle, \tag{6.61}$$

*and for every $k \in \{0,\ldots,K\}$,*

$$\partial_{q_k} \overline{F}_N(t,\mathsf{q}) = \mathbb{E}\langle \mathbf{1}_{\{\alpha^1 \wedge \alpha^2 = k\}} R_{1,2}\rangle. \tag{6.62}$$

*In particular,*

$$\partial_t \overline{F}_N(t,\mathsf{q}) - \sum_{k=0}^{K} (\zeta_{k+1} - \zeta_k) \left(\frac{\partial_{q_k} \overline{F}_N(t,\mathsf{q})}{\zeta_{k+1} - \zeta_k}\right)^2$$
$$= \mathbb{E}\langle (R_{1,2} - \mathbb{E}\langle R_{1,2} \mid \alpha^1 \wedge \alpha^2\rangle)^2\rangle. \tag{6.63}$$

*Proof.* The same calculation as for (6.35) still gives

$$\partial_t \overline{F}_N(t,\mathrm{q}) = -\frac{1}{N}\mathbb{E}\left\langle \frac{1}{\sqrt{2t}} H_N(\sigma) - N \right\rangle = \mathbb{E}\langle R_{1,2}^2 \rangle.$$

Concerning the derivative with respect to $q_k$, the Gibbs Gaussian integration by parts formula (Theorem 4.6) implies that for every $k \in \{0, \ldots, K-1\}$,

$$\partial_{q_k} \overline{F}_N(t,\mathrm{q}) = -\frac{1}{N}\mathbb{E}\left\langle (2q_k - 2q_{k-1})^{-1/2} z_{\alpha_{|k}} \cdot \sigma - (2q_{k+1} - 2q_k)^{-1/2} z_{\alpha_{|k+1}} \cdot \sigma \right\rangle$$

$$= \frac{1}{N}\mathbb{E}\left\langle \left(\mathbf{1}_{\{\alpha^1_{|k} = \alpha^2_{|k}\}} - \mathbf{1}_{\{\alpha^1_{|k+1} = \alpha^2_{|k+1}\}}\right) \sigma^1 \cdot \sigma^2 \right\rangle$$

$$= \mathbb{E}\left\langle \mathbf{1}_{\{\alpha^1 \wedge \alpha^2 = k\}} R_{1,2} \right\rangle.$$

For $k = K$, the last term on the first two lines above disappears, and we have

$$\partial_{q_K} \overline{F}_N(t,\mathrm{q}) = \frac{1}{N}\mathbb{E}\left\langle \mathbf{1}_{\{\alpha^1_{|K} = \alpha^2_{|K}\}} \sigma^1 \cdot \sigma^2 \right\rangle = \mathbb{E}\left\langle \mathbf{1}_{\{\alpha^1 \wedge \alpha^2 = K\}} R_{1,2} \right\rangle.$$

So the relation (6.62) is indeed valid for every $k \in \{0, \ldots, K\}$. Combining this with (6.60) shows that

$$\frac{\partial_{q_k} \overline{F}_N(t,\mathrm{q})}{\zeta_{k+1} - \zeta_k} = \frac{\mathbb{E}\langle \mathbf{1}_{\{\alpha^1 \wedge \alpha^2 = k\}} R_{1,2} \rangle}{\mathbb{E}\langle \mathbf{1}_{\{\alpha^1 \wedge \alpha^2 = k\}} \rangle} = \mathbb{E}\langle R_{1,2} \mid \alpha^1 \wedge \alpha^2 = k \rangle. \tag{6.64}$$

We next observe that the right side of (6.63) can be rewritten as

$$\mathbb{E}\langle R_{1,2}^2 \rangle - \mathbb{E}\left\langle \left(\mathbb{E}\langle R_{1,2} \mid \alpha^1 \wedge \alpha^2 \rangle\right)^2 \right\rangle$$

$$= \partial_t \overline{F}_N(t,\mathrm{q}) - \sum_{k=0}^{K} \mathbb{E}\left\langle \mathbf{1}_{\{\alpha^1 \wedge \alpha^2 = k\}} \left(\mathbb{E}\langle R_{1,2} \mid \alpha^1 \wedge \alpha^2 = k \rangle\right)^2 \right\rangle$$

$$= \partial_t \overline{F}_N(t,\mathrm{q}) - \sum_{k=0}^{K} (\zeta_{k+1} - \zeta_k) \left(\frac{\partial_{q_k} \overline{F}_N(t,\mathrm{q})}{\zeta_{k+1} - \zeta_k}\right)^2.$$

This completes the proof. ∎

By Jensen's inequality, we have

$$\mathbb{E}\left\langle \left(\mathbb{E}\langle R_{1,2} \mid \alpha^1 \wedge \alpha^2 \rangle\right)^2 \right\rangle \geq \mathbb{E}\langle R_{1,2}^2 \rangle, \tag{6.65}$$

and thus

$$\mathbb{E}\left\langle \left(R_{1,2} - \mathbb{E}\langle R_{1,2} \mid \alpha^1 \wedge \alpha^2 \rangle\right)^2 \right\rangle \leq \mathbb{E}\left\langle \left(R_{1,2} - \mathbb{E}\langle R_{1,2}^2 \rangle\right)^2 \right\rangle. \tag{6.66}$$

So the right side of (6.63), which we think of as an error term, is smaller than or equal to the error term we found in (6.36) in our naive attempt at enriching the free energy.

The right side of (6.63) measures the size of the conditional variance of the $\sigma$-overlap given the $\alpha$-overlap. Can we ensure that this is small? In other words, can we assert that we can essentially infer the $\sigma$-overlap from the observation of the $\alpha$-overlap? Recall also that, as for the Curie-Weiss model and the models of statistical inference, we do not aim for the error term to be small for absolutely all choices of the parameters, but only for "most" such choices.

As was discovered in [213, 214, 215], the Panchenko ultrametricity theorem (Theorem 5.30) implies the following synchronization phenomenon: up to a small perturbation of the Hamiltonian (to ensure a family of Ghirlanda-Guerra identities), the $\sigma$- and the $\alpha$-overlaps must be asymptotically monotonically coupled. A pair of real random variables $(X,Y)$ is said to be *monotonically coupled* if, denoting by $(X',Y')$ an independent copy of $(X,Y)$, we have

$$\mathbb{P}\{X < X' \text{ and } Y' < Y\} = 0. \tag{6.67}$$

Denoting by $F_X^{-1}$ and $F_Y^{-1}$ the inverse cumulative distribution functions of $X$ and $Y$ respectively, one can show (see Proposition 2.1 in [197]) that the pair $(X,Y)$ is monotonically coupled if and only if its law is that of $(F_X^{-1}(U), F_Y^{-1}(U))$, where $U$ is a uniform random variable on the interval $[0,1]$.

In the monotonically coupled pair $(X,Y)$, we think of $X$ as the $\sigma$-overlap and of $Y$ as the $\alpha$-overlap, and wonder whether we can recover $X$ from the observation of $Y$. Naturally, if $Y$ is deterministic, then the observation of $Y$ cannot help much in predicting $X$. On the other hand, if the law of $Y$ is nicely spread out, for instance if it is a uniform random variable on $[0,1]$, then the observation of $Y$ allows us to recover the "secret" random variable $U$ in the representation above, and thus to perfectly infer the value of $X$.

Coming back to the analysis of the overlaps, we thus see that with the synchronization phenomenon in place, the error term on the right side of (6.63) will be small provided that the law of the $\alpha$-overlap is sufficiently spread out. Recall also from (6.60) that we know the law of the $\alpha$-overlap exactly. Combining all this, it was shown in Proposition 5.5 of [195] that when a family of Ghirlanda-Guerra identities all hold asymptotically as $N$ tends to infinity, and when $\zeta_{k+1} - \zeta_k = 1/(K+1)$ for every $k \in \{0,\ldots,K\}$, we have that

$$\limsup_{N\to+\infty} \mathbb{E}\left\langle \left(R_{1,2} - \mathbb{E}\langle R_{1,2} \mid \alpha^1 \wedge \alpha^2\rangle\right)^2\right\rangle \leqslant \frac{12}{K}. \tag{6.68}$$

This informal discussion motivates us to give more substance to the idea that the free energy $\overline{F}_N(t,\mathsf{q})$ should really be thought of as a function of the path $\mathsf{q}$, as was already suggested by the notation, and to seek to make sense of it for more general paths than the piecewise-constant paths we have so far considered. To start with, we give an informal justification of the fact that we decided to bundle together the parameters $(\zeta_k)_{1 \leqslant k \leqslant K}$ and $(q_k)_{0 \leqslant k \leqslant K}$ into the path $\mathsf{q}$ in (6.50). We have already discussed extensively, in particular in Section 5.7, that it makes sense to focus on understanding the behaviour of overlap arrays. While we have chosen to define the $\alpha$-overlap simply as $\alpha^1 \wedge \alpha^2$, in view of (6.52), it would have been more faithful to the definition of our random field $Z_q$ to have defined this same quantity as seen through the sequence of parameters $(q_k)_{0 \leq k \leq K}$, namely $q_{\alpha^1 \wedge \alpha^2}$. And we recall from (6.60) that for every $k \in \{0, \ldots, K\}$, we have

$$\mathbb{E} \langle \mathbf{1}_{\{q_{\alpha^1 \wedge \alpha^2} = q_k\}} \rangle = \zeta_{k+1} - \zeta_k. \tag{6.69}$$

Recalling the definition of the path q in (6.50), we thus see that the law of $q_{\alpha^1 \wedge \alpha^2}$, which one might argue is the "true" $\alpha$-overlap, is the law of $q(U)$, where $U$ is a uniform random variable over $[0,1]$. With this in mind, it is perhaps not surprising that the function $\overline{F}_N(t, q)$ satisfies a natural continuity property as a function of the path q. We write

$$Q(\mathbb{R}_{\geq 0}) := \{q : [0,1) \to \mathbb{R}_{\geq 0} \mid q \text{ is right-continuous and non-decreasing}\}, \tag{6.70}$$

and for every $p \in [1, +\infty]$, we let $Q_p(\mathbb{R}_{\geq 0}) := Q(\mathbb{R}_{\geq 0}) \cap L^p([0,1];\mathbb{R})$. We also note that when we say that a path q is piecewise-constant, we mean that we can partition the interval $[0,1)$ into a *finite* number of intervals in the interior of which q is constant. In other words, a piecewise-constant path $q \in Q(\mathbb{R}_{\geq 0})$ is any path that can be represented in the form of (6.50) for a suitable choice of the parameters.

**Proposition 6.3.** *For any pair of piecewise-constant paths* $q, q' \in Q(\mathbb{R}_{\geq 0})$,

$$\left| \overline{F}_N(t, q) - \overline{F}_N(t, q') \right| \leq \int_0^1 |q(u) - q'(u)| \, du. \tag{6.71}$$

*In particular, the enriched free energy* (6.53) *admits a unique extension to the space* $\mathbb{R}_{\geq 0} \times Q_1(\mathbb{R}_{\geq 0})$ *satisfying* (6.71) *for all times* $t \geq 0$ *and paths* $q, q' \in Q_1(\mathbb{R}_{\geq 0})$.

*Proof.* The key step of the proof is to realize that we can allow for repetitions in the parameters $(q_k)_{0 \leq k \leq K}$ without encountering contradictions. The bound (6.71) then follows by simple estimates on the derivatives obtained in Lemma 6.2. We decompose the proof into three steps.

*Step 1: integrating the Poisson-Dirichlet cascade.* We use Theorem 5.25 to describe an explicit procedure for integrating out the Poisson-Dirichlet cascade. We fix parameters $(\zeta_k)_{1 \leq k \leq K}$ and $(q_k)_{0 \leq k \leq K}$ as in (6.48)-(6.49), and let q be as in (6.50). For every $y_0, \ldots, y_K \in \mathbb{R}^N$, we define

$$X_K(y_0, \ldots, y_K) :=$$
$$\log \int_{\Sigma_N} \exp\left( \sqrt{2t} H_N(\sigma) - Nt + \sum_{k=0}^K (2q_k - 2q_{k-1})^{1/2} y_k \cdot \sigma - N q_K \right) dP_N(\sigma), \tag{6.72}$$

and then recursively, for every $k \in \{0, \ldots, K-1\}$,

$$X_k = X_k(y_0, \ldots, y_k) := \frac{1}{\zeta_{k+1}} \log \mathbb{E}_{y_{k+1}} \exp(\zeta_{k+1} X_{k+1}),$$

where we denote by $\mathbb{E}_{y_{k+1}}$ the integration of the variable $y_{k+1}$ according to the $N$-dimensional standard Gaussian probability measure. We also define $X_{-1} := \mathbb{E}_{y_0} X_0$, with the same interpretation. Using a measurable mapping from $[0,1]$ to $\mathbb{R}^N$ that sends the Lebesgue measure on $[0,1]$ to the standard Gaussian probability measure on $\mathbb{R}^N$ (which can be constructed by first mapping the Lebesgue measure on $[0,1]$ to the Lebesgue measure on $[0,1]^N$ by splitting the binary expansion of the argument, and then using Exercise 5.11 coordinate by coordinate), we see that we can apply Theorem 5.25 and obtain that

$$-N\mathbb{E}\big[F_N(t, \mathsf{q}) \mid (H_N(\sigma))_{\sigma \in \mathbb{R}^N}\big] = X_{-1}, \tag{6.73}$$

where we write

$$F_N(t, \mathsf{q}) := -\frac{1}{N} \log \int_{\Sigma_N} \sum_{\alpha \in \mathbb{N}^K} \exp\big(H_N(t, \mathsf{q}, \sigma, \alpha)\big) v_\alpha \, \mathrm{d}P_N(\sigma).$$

In other words, $F_N(t, \mathsf{q})$ is the quantity on the right side of (6.53) before taking the expectation, and we have

$$\overline{F}_N(t, \mathsf{q}) = \mathbb{E} F_N(t, \mathsf{q}).$$

In the derivation of (6.73), we used that that the conditional expectation appearing there consists in averaging over the randomness coming from the Poisson-Dirichlet cascade and the random field $Z_\mathsf{q}$, since these are the only other sources of randomness and they have been chosen independently of $(H_N(\sigma))_{\sigma \in \mathbb{R}^N}$.

*Step 2: allowing for possible repetitions of parameters.* We fix an integer $K \geq 1$ as well as sequences of parameters

$$0 = \zeta_0 < \zeta_1 < \cdots < \zeta_K < \zeta_{K+1} = 1,$$
$$0 = q_{-1} \leq q_0 \leq q_1 \leq \cdots \leq q_K < q_{K+1} = +\infty,$$

and consider the piecewise-constant path

$$\mathsf{q} := \sum_{k=0}^{K} q_k \mathbf{1}_{[\zeta_k, \zeta_{k+1})}.$$

Since $\mathsf{q}$ is piecewise-constant, we have already given a definition of $\overline{F}_N(t, \mathsf{q})$, but this definition would require that we rewrite $\mathsf{q}$ into a form that suppresses possible repetitions in the parameters. We claim that if we instead naively follow the procedure with the repeated parameters, we obtain the same value for $\overline{F}_N(t, \mathsf{q})$.

If there is a repetition in the parameters $(q_k)_{0\leqslant k\leqslant K}$, say $q_{k-1} = q_k$ for some $k \in \{0,\ldots,K\}$, then the term indexed by $k$ in the sum on the right side of (6.72) vanishes, and thus for every $y_0,\ldots,y_k \in \mathbb{R}^N$, we have that

$$X_{k-1}(y_0,\ldots,y_{k-1}) = X_k(y_0,\ldots,y_k).$$

It is therefore clear that removing this repetition in the sequence $(q_k)_{0\leqslant k\leqslant K}$ does not change the value of the result.

*Step 3: extending Lemma 6.2 to the case of repeated indices.* We now justify that the identity (6.62) from Lemma 6.2 remains valid even if we allow for repetitions in the parameters $(q_k)_{0\leqslant k\leqslant K}$. We know that (6.62) holds when the $(q_k)_{0\leqslant k\leqslant K}$ satisfy (6.49). Moreover, the free energy is clearly a continuous function of the parameters $(q_k)_{0\leqslant k\leqslant K}$, and the right side of (6.62) is also continuous with respect to these parameters. A classical analysis argument thus suffices to conclude that the free energy $\overline{F}_N$ is differentiable also at the boundary of the set of non-decreasing sequences $(q_k)_{0\leqslant k\leqslant K}$, and that the relation (6.62) still holds there.

*Step 4: extending to $Q_1(\mathbb{R}_{\geqslant 0})$.* Now that we can allow for repetitions in the parameters, for any two piecewise-constant paths $\mathsf{q}$ and $\mathsf{q}'$, we can always find an integer $K \geqslant 1$ and sequences of parameters $0 = \zeta_0 < \zeta_1 < \cdots < \zeta_K < \zeta_{K+1} = 1$, $0 = q_{-1} \leqslant q_0 \leqslant \cdots \leqslant q_K$ and $0 = p_{-1} \leqslant p_0 \leqslant \cdots \leqslant p_K$ with

$$\mathsf{q} = \sum_{k=0}^{K} q_k \mathbf{1}_{[\zeta_k,\zeta_{k+1})} \quad \text{and} \quad \mathsf{q}' = \sum_{k=0}^{K} p_k \mathbf{1}_{[\zeta_k,\zeta_{k+1})}.$$

For each $s \in [0,1]$, we introduce the piecewise-constant path

$$\mathsf{q}^s := \sum_{k=0}^{K} (sq_k + (1-s)p_k) \mathbf{1}_{[\zeta_k,\zeta_{k+1})},$$

and define the interpolating free energy $\varphi(s) := \overline{F}_N(t,\mathsf{q}^s)$. The fundamental theorem of calculus and the chain rule imply that

$$\left|\overline{F}_N(t,\mathsf{q}) - \overline{F}_N(t,\mathsf{q}')\right| \leqslant \sup_{s\in[0,1]} |\varphi'(s)| \leqslant \sum_{k=0}^{K} \sup_{s\in[0,1]} |\partial_{q_k}\overline{F}_N(t,\mathsf{q}^s)||q_k - p_k|.$$

It follows from the previous step and (6.60) that

$$\left|\overline{F}_N(t,\mathsf{q}) - \overline{F}_N(t,\mathsf{q}')\right| \leqslant \sum_{k=0}^{K} (\zeta_{k+1} - \zeta_k)|q_k - p_k| = \int_0^1 |\mathsf{q}(u) - \mathsf{q}'(u)|\,\mathrm{d}u.$$

A density argument then allows us to uniquely extend the enriched free energy to the space $Q_1(\mathbb{R}_{\geqslant 0})$ while preserving the bound (6.71). This completes the proof. ∎

To write equation (6.63) in Lemma 6.2 in a way that depends directly on the path $q \in Q(\mathbb{R}_{\geq 0})$ as opposed to the parameters $(\zeta_k)_{1 \leq k \leq K}$ and $(q_k)_{0 \leq k \leq K}$, we introduce additional notation. Given a function $h : \mathbb{R}_{\geq 0} \times Q_2(\mathbb{R}_{\geq 0}) \to \mathbb{R}$, a time $t \geq 0$ and a path $q \in Q_2(\mathbb{R}_{\geq 0})$, we say that $h$ admits a *Gateaux derivative* at $(t, q) \in \mathbb{R}_{\geq 0} \times Q_2(\mathbb{R}_{\geq 0})$ if there exists a function $\partial_q h(t, q, \cdot) \in L^2([0,1]; \mathbb{R})$ such that for every $q' \in L^2([0,1]; \mathbb{R})$ with $q + \varepsilon q' \in Q_2(\mathbb{R}_{\geq 0})$ for $\varepsilon > 0$ sufficiently small, we have as $\varepsilon > 0$ tends to zero that

$$h(t, q + \varepsilon q') - h(t, q) = \varepsilon \int_0^1 \partial_q h(t, q, u) q'(u) \, du + o(\varepsilon). \tag{6.74}$$

Let $q$ be as in (6.50), and assume that the enriched free energy $\overline{F}_N$ admits a Gateaux derivative at the pair $(t, q)$. For every other path $q' = \sum_{k=0}^K p_k \mathbf{1}_{[\zeta_k, \zeta_{k+1})}$ associated with a non-decreasing sequence $0 = p_{-1} \leq p_0 \leq \cdots \leq p_K$, we have

$$\overline{F}_N(t, q + \varepsilon q') - \overline{F}_N(t, q) = \varepsilon \sum_{k=0}^K p_k \int_{\zeta_k}^{\zeta_{k+1}} \partial_q \overline{F}_N(t, q, u) \, du + o(\varepsilon). \tag{6.75}$$

On the other hand, a Taylor expansion reveals that

$$\overline{F}_N(t, q + \varepsilon q') - \overline{F}_N(t, q) = \varepsilon \sum_{k=0}^K p_k \partial_{q_k} \overline{F}_N(t, q) \mathbf{1}_{[\zeta_k, \zeta_{k+1})} + o(\varepsilon). \tag{6.76}$$

Since these equalities must hold for all piecewise-constant paths $q'$, we must have

$$\partial_{q_k} \overline{F}_N(t, q) = \int_{\zeta_k}^{\zeta_{k+1}} \partial_q \overline{F}_N(t, q, u) \, du. \tag{6.77}$$

Assuming that $\max_k(\zeta_{k+1} - \zeta_k)$ is small, we have the approximate equality

$$\int_0^1 \partial_q \overline{F}_N(t, q, u)^2 \, du \simeq \sum_{k=0}^K (\zeta_{k+1} - \zeta_k) \left( \frac{1}{\zeta_{k+1} - \zeta_k} \int_{\zeta_k}^{\zeta_{k+1}} \partial_q \overline{F}_N(t, q, u) \, du \right)^2, \tag{6.78}$$

see also Exercise 6.6 for a more precise argument. Up to an error that tends to zero as $\max_k(\zeta_{k+1} - \zeta_k)$ tends to zero, we can therefore rewrite the identity (6.63) from Lemma 6.2 as

$$\partial_t \overline{F}_N(t, q) - \int_0^1 \partial_q \overline{F}_N(t, q, u)^2 \, du \simeq \mathrm{Var}(R_{1,2} \mid \alpha^1 \wedge \alpha^2), \tag{6.79}$$

where we used the notation

$$\mathrm{Var}(R_{1,2} \mid \alpha^1 \wedge \alpha^2) := \mathbb{E}\left\langle \left(R_{1,2} - \mathbb{E}\langle R_{1,2} \mid \alpha^1 \wedge \alpha^2 \rangle\right)^2 \right\rangle. \tag{6.80}$$

Combined with the synchronization mechanism discussed earlier, see in particular (6.68), this suggests that the limit free energy $f : \mathbb{R}_{\geq 0} \times \mathcal{Q}_2(\mathbb{R}_{\geq 0}) \to \mathbb{R}$ should satisfy the infinite-dimensional Hamilton-Jacobi equation

$$\partial_t f(t, \mathsf{q}) - \int_0^1 \partial_{\mathsf{q}} f(t, \mathsf{q}, u)^2 \, du = 0 \quad \text{on} \quad \mathbb{R}_{>0} \times \mathcal{Q}_2(\mathbb{R}_{\geq 0}). \tag{6.81}$$

Our next goal will be to explain that this guess allows us to correctly predict the Parisi formula presented in (6.8).

**Exercise 6.5.** We consider the SK model on $\Sigma_N$ with $P_N := (P_1)^{\otimes N}$ for a possibly non-centred probability measure $P_1$ on $\Sigma_1$, and we assume that there exists $q \in \mathbb{R}$ with

$$\lim_{N \to +\infty} \mathbb{E}\langle (R_{1,2} - q)^2 \rangle = 0.$$

Let $g = (g_1, \ldots, g_N)$ be a standard Gaussian vector and $Z$ be a standard Gaussian random variable independent of all other sources of randomness and of each other. We understand that the probability $\mathbb{P}$ and associated expectation $\mathbb{E}$ also takes the average over $g$, and we denote by $\mathbb{E}_Z$ the average over the randomness of $Z$. Show that for every $F \in C_c(\mathbb{R}; \mathbb{R})$ and $\varepsilon > 0$, we have

$$\lim_{N \to +\infty} \mathbb{P}\left\{ \left| \left\langle F\left(\frac{g \cdot \sigma}{\sqrt{N}}\right) \right\rangle - \mathbb{E}_Z F\left(\frac{g \cdot \langle \sigma \rangle}{\sqrt{N}} + \sqrt{1-q}Z\right) \right| > \varepsilon \right\} = 0. \tag{6.82}$$

In words, for most realizations of $g$, the law of the random variable $N^{-1/2} g \cdot \sigma$ under $\langle \cdot \rangle$ is essentially Gaussian, with mean $N^{-1/2} g \cdot \langle \sigma \rangle$ and variance $1 - q$.

**Exercise 6.6.** Let $f \in L^2([0,1]; \mathbb{R})$. Show that for every $\varepsilon > 0$, there exists $\delta > 0$ such that for every $0 = \zeta_0 < \zeta_1 < \cdots < \zeta_K < \zeta_{K+1} = 1$, if $\max_k(\zeta_{k+1} - \zeta_k) \leq \delta$, then

$$\left| \int_0^1 f^2 - \sum_{k=0}^K (\zeta_{k+1} - \zeta_k) \left( \frac{1}{\zeta_{k+1} - \zeta_k} \int_{\zeta_k}^{\zeta_{k+1}} f \right)^2 \right| \leq \varepsilon. \tag{6.83}$$

## 6.5 Connecting the Hamilton-Jacobi approach with the Parisi formula

The goal of this section is to demonstrate that the Hamilton-Jacobi equation (6.81) that we derived heuristically allows us to correctly predict the Parisi formula (6.8).

Our first task will be to identify the initial condition to this Hamilton-Jacobi equation. Just as in the Curie-Weiss model and in the models of statistical inference we considered in Chapter 4, we will be able to leverage the product structure of the reference measure $P_N$ to show that the quantity $\overline{F}_N(0, \mathsf{q})$ in fact does not depend on $N$. In order to connect to the Parisi formula as presented in (6.8), we will also relate this quantity to the second-order partial differential equation appearing in (6.9).

Once this is done, we need to make sense of the infinite-dimensional Hamilton-Jacobi equation (6.81). There are at least three possible ways to approach this problem. A first possibility would be to seek to adapt the notion of viscosity solution to infinite-dimensional spaces. Alternatively, we could devise finite-dimensional approximations of the equation, and then try to show the convergence of these approximations as the dimension is sent to infinity. Finally, we could decide that since the non-linearity in (6.104) is convex, whichever notion of solution we come up with will admit a variational representation analogous to the Hopf-Lax formula from Theorem 3.8, and so we may as well define the solution according to this variational formula directly. We will start by exploring the last of these options, since it is the easiest. As will be stressed in Section 6.6, this approach is no longer applicable for non-convex models such as the bipartite model, and we will therefore also outline the workings of the second of these three options there. That all these approaches yield the same function when they are applicable is shown in [73].

We start by studying the initial condition associated with the Hamilton-Jacobi equation (6.81). In the next lemma, we will show in particular that

$$\psi(\mathsf{q}) := \lim_{N \to +\infty} \overline{F}_N(0,\mathsf{q}) = \overline{F}_1(0,\mathsf{q}). \tag{6.84}$$

For each piecewise-constant path $\mathsf{q}$, Theorem 5.25 describes an explicit procedure for computing this quantity in terms of recursive averages over the Gaussian measure of quantities of the form

$$\frac{1}{\zeta_k} \log \mathbb{E} \exp \zeta_k X_k.$$

In order to encode this recursive procedure even for general paths, one can view this calculation in terms of a second-order partial differential equation. To illustrate the idea, let us observe that if $(B_t)_{t \geq 0}$ is a standard Brownian motion and $f \in C_c(\mathbb{R}; \mathbb{R}_{\geq 0})$, then the function

$$u(t,x) := \mathbb{E} f(x + B_{2t}) \tag{6.85}$$

solves the heat equation $\partial_t u = \partial_x^2 u$, and a simple calculation yields that the function $\widetilde{u} := \zeta_k^{-1} \log u$ is then a solution to $\partial_t \widetilde{u} = \partial_x^2 \widetilde{u} + \zeta_k (\partial_x \widetilde{u})^2$. The representation of the initial condition in the next lemma will be an iteration of this observation.

To state this precisely, we start by recalling that we denote by $\mathsf{Q}_\infty(\mathbb{R}_{\geq 0})$ the space of paths in $\mathsf{Q}(\mathbb{R}_{\geq 0})$ that remain bounded. For every such path $\mathsf{q} \in \mathsf{Q}_\infty(\mathbb{R}_{\geq 0})$, we define

$$\mathsf{q}(1) := \lim_{u \nearrow 1} \mathsf{q}(u). \tag{6.86}$$

In particular, for $\mathsf{q}$ as in (6.50), we have $\mathsf{q}(1) = q_K$. For every $\mathsf{q} \in \mathsf{Q}_\infty(\mathbb{R}_{\geq 0})$, we also define its right-continuous inverse $\mathsf{q}^{-1} : \mathbb{R} \to [0,1]$ according to

$$\mathsf{q}^{-1}(t) := \sup\{u \in [0,1] \mid \mathsf{q}(u) \leq t\}. \tag{6.87}$$

We understand the supremum above to be zero in case the set is empty. Finally, for every $t \geq 0$ and $x \in \mathbb{R}$, we write

$$\phi(t,x) := \log \int_{\Sigma_1} \exp(x\sigma - t) \, dP_1(\sigma). \tag{6.88}$$

In case $P_1$ were not supported on $\Sigma_1$, we would replace the term $t$ in the exponential by $t\sigma^2$.

**Lemma 6.4** (Parisi PDE). *For every piecewise-constant path $\mathsf{q} \in \mathcal{Q}(\mathbb{R}_{\geq 0})$ and integer $N \geq 1$, we have*

$$\overline{F}_N(0,\mathsf{q}) = \overline{F}_1(0,\mathsf{q}) = -\Phi^\mathsf{q}(0,0), \tag{6.89}$$

*where the function $\Phi^\mathsf{q} : [0,\mathsf{q}(1)] \times \mathbb{R} \to \mathbb{R}$ satisfies*

$$\begin{cases} -\partial_t \Phi^\mathsf{q}(t,x) = \partial_x^2 \Phi^\mathsf{q}(t,x) + \mathsf{q}^{-1}(t)(\partial_x \Phi^\mathsf{q}(t,x))^2 & \text{on } [0,\mathsf{q}(1)] \times \mathbb{R}, \\ \Phi^\mathsf{q}(\mathsf{q}(1),x) = \phi(\mathsf{q}(1),x) & \text{for } x \in \mathbb{R}. \end{cases} \tag{6.90}$$

Before proceeding to the proof, we discuss the nature of the partial differential equation (6.90) and the associated notion of solution. At small scales, the first-order term $(\partial_x \Phi^\mathsf{q})^2$ becomes less and less important compared with the second-order term $\partial_x^2 \Phi^\mathsf{q}$, which is regularizing. As a consequence, this equation is much easier to study than the first-order Hamilton-Jacobi equations that have accompanied us throughout the book, and the only slight difficulty is that the function $\mathsf{q}^{-1}$ has some jump discontinuities. For each $\mathsf{q} \in \mathcal{Q}_\infty(\mathbb{R}_{\geq 0})$, we can choose among several equivalent definitions of the notion of being a solution to (6.90).

Perhaps the most intuitive notion of solution, which we could call the notion of classical solution, would be to ask that $\Phi^\mathsf{q}$ be a continuous function such that, at every $(t,x) \in [0,\mathsf{q}(1)] \times \mathbb{R}$ except those $t$'s that are points of discontinuity of $\mathsf{q}^{-1}$, the function $\Phi^\mathsf{q}$ is continuously differentiable in $t$ and twice continuously differentiable in $x$, and the first line of (6.90) holds for each such $(t,x)$. To guarantee uniqueness, this must be supplemented with a mild growth condition, in our case we can simply require $\Phi^\mathsf{q}$ to have at most linear growth in $x$ at infinity, uniformly over $t$.

An alternative and analytically very convenient notion of solution is that of mild solution. We could for instance look for a function in $C([0,\mathsf{q}(1)];C^2(\mathbb{R};\mathbb{R}))$ with at most linear growth and such that, for every $t \in (0,\mathsf{q}(1))$, we have

$$\Phi^\mathsf{q}(t,x) = e^{(\mathsf{q}(1)-t)\Delta}\phi(\mathsf{q}(1),\cdot)(x)$$
$$+ \int_t^{\mathsf{q}(1)} \mathsf{q}^{-1}(t) e^{(s-t)\Delta}(\partial_x \Phi^\mathsf{q}(s,\cdot))^2(x) \, ds, \tag{6.91}$$

where we denote by $(e^{t\Delta})_{t \geq 0}$ the heat semigroup; explicitly, the heat semigroup is such that for every measurable function $f : \mathbb{R} \to \mathbb{R}$ for which the integral is well-defined,
$$e^{t\Delta}f(x) := \int_{\mathbb{R}} \frac{1}{\sqrt{4\pi t}} \exp\left(-\frac{(x-y)^2}{4t}\right) f(y)\, dy.$$

This formulation allows one to show the well-posedness of equation (6.90) over a short time interval using a fixed-point argument, and this can then be iterated using estimates from the maximum principle; see also [147]. If desired, one can then show the equivalence between the notions of classical and mild solutions.

We will not go into details on these points. At least, the arguments we now present yield in particular that a solution to (6.90) indeed exists, starting with the proof below covering the case of piecewise-constant q.

*Proof of Lemma 6.4.* We fix a piecewise-constant path $\mathsf{q} \in \mathcal{Q}(\mathbb{R}_{\geqslant 0})$ of the form (6.50). Applying Corollary 5.26 to the random variables
$$X_{K,i} := \log \int_{\Sigma_1} \exp(Z_{\mathsf{q}}^i(\alpha)\sigma - q_K)\, dP_1(\sigma)$$
shows that $\overline{F}_N(0,\mathsf{q}) = \overline{F}_1(0,\mathsf{q})$, so the initial condition (6.84) is well-defined and is given by
$$\psi(\mathsf{q}) = \overline{F}_1(0,\mathsf{q}) = -\log \mathbb{E} \int_{\Sigma_1} \sum_{\alpha \in \mathbb{N}^K} \exp(Z_{\mathsf{q}}(\alpha)\sigma - q_K) v_\alpha\, dP_1(\sigma). \qquad (6.92)$$

To express $\psi(\mathsf{q})$ as the value at $(0,0)$ of the solution to (6.90), we denote by $(B_t)_{t\geqslant 0}$ a standard Brownian motion, and recursively define the function $\Phi^{\mathsf{q}} : [0,\mathsf{q}(1)] \times \mathbb{R} \to \mathbb{R}$ by $\Phi^{\mathsf{q}}(\mathsf{q}(1),x) := \phi(q_K,x)$ and
$$\Phi^{\mathsf{q}}(t,x) := \zeta_{k+1}^{-1} \log \mathbb{E} \exp \zeta_{k+1} \Phi^{\mathsf{q}}(q_{k+1}, x + B_{2q_{k+1}-2t}) \qquad (6.93)$$
for $t \in [q_k, q_{k+1})$ and $-1 \leqslant k \leqslant K-1$. It then follows from Theorem 5.25 that
$$\psi(\mathsf{q}) = -\Phi^{\mathsf{q}}(0,0).$$

We now show that $\Phi^{\mathsf{q}}$ satisfies the equation (6.90). It will be important to observe that
$$\mathsf{q}^{-1} = \sum_{k=0}^{K+1} \zeta_k \mathbf{1}_{[q_{k-1},q_k)}. \qquad (6.94)$$

For each $k \in \{-1,\ldots,K-1\}$, we consider the function defined for $t \in [q_k, q_{k+1})$ by
$$\widetilde{\Phi}^{\mathsf{q}}(t,x) := \exp \zeta_{k+1} \Phi^{\mathsf{q}}(t,x).$$

As for (6.85), we see that the function $\widetilde{\Phi}^{\mathsf{q}}$ satisfies
$$-\partial_t \widetilde{\Phi}^{\mathsf{q}} = \partial_x^2 \widetilde{\Phi}^{\mathsf{q}} \quad \text{on} \quad [q_k, q_{k+1}) \times \mathbb{R}.$$

Taking the logarithm, we obtain that

$$-\partial_t \Phi^q = \partial_x^2 \Phi^q + \zeta_{k+1}(\partial_x \Phi^q)^2 \quad \text{on} \quad [q_k, q_{k+1}) \times \mathbb{R}.$$

Remembering from (6.94) that $q^{-1} = \zeta_{k+1}$ on $[q_k, q_{k+1})$ completes the proof. ∎

In order to extend Lemma 6.4 to arbitrary paths $q \in Q_\infty(\mathbb{R}_{\geq 0})$, one possible route would be to show the continuity of $\Phi^q(0,0)$ as a function of q, since we already know the continuity of $\overline{F}_1(0,\cdot)$ from Proposition 6.3. Since we have not been very precise in our discussion on the well-posedness of (6.90), we prefer to proceed in a slightly different way: by leveraging Lemma 6.4 and slightly upgrading Proposition 6.3, we will justify that taking solutions to (6.90) for a sequence of piecewise-constant paths that converge to some arbitrary $q \in Q_\infty(\mathbb{R}_{\geq 0})$, we can pass to the limit and obtain a solution to (6.90) for this limit path. In particular, this shows that a solution to (6.90) exists for arbitrary $q \in Q_\infty(\mathbb{R}_{\geq 0})$, and the only point of interest that we did not justify with full precision is that this solution is unique; this uniqueness property can be obtained by a classical Picard fixed-point argument.

Before starting the proof, we make a simple observation that is convenient when wanting to compare $\Phi^q$ for different choices of q. The observation is that, if we were to define $\Phi^q$ over a larger time interval than $[0, q(1)]$, then we would end up with a function that coincides with $\Phi^q$ on $[0, q(1)]$. In other words, for each $T \geq q(1)$ and for $u$ the solution to

$$\begin{cases} -\partial_t u(t,x) = \partial_x^2 u(t,x) + q^{-1}(t)(\partial_x u(t,x))^2 & \text{on } [0,T] \times \mathbb{R}, \\ u(T,x) = \phi(T,x) & \text{for } x \in \mathbb{R}, \end{cases} \qquad (6.95)$$

we have that $u$ and $\Phi^q$ coincide on $[0, q(1)] \times \mathbb{R}$. This comes from the fact that the function $\phi$ in (6.88) is such that

$$-\partial_t \phi = \partial_x^2 \phi + (\partial_x \phi)^2.$$

Indeed, for the natural Gibbs measure associated with (6.88), the derivative $-\partial_t \phi$ computes the second moment of $\sigma$, while $\partial_x^2 \phi$ computes its variance and $\partial_x \phi$ its first moment. With this observation in mind, we can without loss of generality think of the function $\Phi^q$ as defined on $\mathbb{R}_{\geq 0} \times \mathbb{R}$, by simply extending it to be equal to $\phi$ outside of $[0, q(1)] \times \mathbb{R}$. For every $T \geq q(1)$, this extension satisfies the partial differential equation in (6.95).

**Lemma 6.5.** *(i) For every piecewise-constant path $q \in Q(\mathbb{R}_{\geq 0})$, and for $\Phi^q$ the solution to (6.90), we have for every $t \geq 0$ and $x \in \mathbb{R}$ that*

$$|\partial_x \Phi^q(t,x)| \leq 1 \qquad \text{and} \qquad 0 \leq \partial_x^2 \Phi^q(t,x) \leq 1. \qquad (6.96)$$

(ii) *For every piecewise-constant* $\mathsf{q}, \mathsf{q}' \in Q(\mathbb{R}_{\geqslant 0})$, $t \geqslant 0$ *and* $x \in \mathbb{R}$,

$$|\Phi^{\mathsf{q}}(t,x) - \Phi^{\mathsf{q}'}(t,x)| \leqslant \int_0^1 |\mathsf{q}(u) - \mathsf{q}'(u)| \, du. \tag{6.97}$$

(iii) *For every piecewise-constant* $\mathsf{q}, \mathsf{q}' \in Q(\mathbb{R}_{\geqslant 0})$, $t \geqslant 0$ *and* $x \in \mathbb{R}$,

$$|\partial_x \Phi^{\mathsf{q}}(t,x) - \partial_x \Phi^{\mathsf{q}'}(t,x)| \leqslant 4 \int_0^1 |\mathsf{q}(u) - \mathsf{q}'(u)| \, du. \tag{6.98}$$

*Proof.* We decompose the proof into four steps.

*Step 1: reduction to $t = 0$.* We first show that it suffices to prove the estimates in the statement of Lemma 6.5 for $t = 0$. Recall from (6.94) that for $\mathsf{q}$ as in (6.50), we have

$$\mathsf{q}^{-1} = \sum_{k=0}^{K+1} \zeta_k \mathbf{1}_{[q_{k-1}, q_k)}.$$

We claim that for each fixed $t_0 \in [0, \mathsf{q}(1)]$, say with $t_0 \in [q_{k_0-1}, q_{k_0})$ for some $k_0 \in \{0, \ldots, K\}$, there is a simple relationship between $\Phi^{\mathsf{q}}(t_0 + t, x)$ and $\Phi^{\widetilde{\mathsf{q}}}(t, x)$, where $\widetilde{\mathsf{q}} \in Q(\mathbb{R}_{\geqslant 0})$ is such that

$$\widetilde{\mathsf{q}}^{-1} = \zeta_{k_0} \mathbf{1}_{[0, q_{k_0} - t_0)} + \sum_{k=k_0+1}^{K+1} \zeta_k \mathbf{1}_{[q_{k-1} - t_0, q_k - t_0)}.$$

Since when a constant is added to the initial condition in (6.90), this only changes the solution by this same constant, we see indeed that for every $t \in [0, \mathsf{q}(1) - t_0]$ and $x \in \mathbb{R}$, we have

$$\Phi^{\mathsf{q}}(t_0 + t, x) + t_0 = \Phi^{\widetilde{\mathsf{q}}}(t, x).$$

This justifies that we may as well focus on the case $t = 0$. If $P_1$ were not supported on $\Sigma_1$, we would need to argue slightly differently, seeking instead to interpret $\Phi^{\widetilde{\mathsf{q}}}$ as a version of $\Phi^{\mathsf{q}}$ for a modified reference measure $P_1$.

*Step 2: proof of (i).* The representation (6.93) of the function $\Phi^{\mathsf{q}}$ for a piecewise-constant path $\mathsf{q} \in Q(\mathbb{R}_{\geqslant 0})$ together with Theorem 5.25 imply that

$$\Phi^{\mathsf{q}}(0,x) = \mathbb{E} \log \int_{\Sigma_1} \sum_{\alpha \in \mathbb{N}^K} \exp(\sigma(x + Z_{\mathsf{q}}(\alpha)) - q_K) v_\alpha \, dP_1(\sigma). \tag{6.99}$$

Denoting by $\langle \cdot \rangle$ the Gibbs average for $N = 1$ associated with the Hamiltonian $(\sigma, \alpha) \mapsto \sigma(x + Z_{\mathsf{q}}(\alpha))$, we have

$$\partial_x \Phi^{\mathsf{q}}(0,x) = \mathbb{E}\langle \sigma \rangle \quad \text{and} \quad \partial_x^2 \Phi^{\mathsf{q}}(0,x) = \mathbb{E}\langle 1 - \sigma^1 \sigma^2 \rangle = 1 - \mathbb{E}\langle \sigma \rangle^2.$$

From these equalities, it is immediate that $|\partial_x \Phi^{\mathsf{q}}(0,x)| \leqslant 1$ and $0 \leqslant \partial_x^2 \Phi^{\mathsf{q}}(0,x) \leqslant 1$.

*Step 3: proof of (ii).* The result follows from the observation that using (6.99), the proof of Proposition 6.3 applies verbatim.

*Step 4: proof of (iii).* We fix two piecewise-constant paths $\mathsf{q},\mathsf{q}' \in Q(\mathbb{R}_{\geq 0})$. We can assume without loss of generality that there exist an integer $K \geq 1$ and sequences of parameters $0 = \zeta_0 < \zeta_1 < \cdots < \zeta_K < \zeta_{K+1} = 1$, $0 = q_{-1} \leq q_0 \leq \cdots \leq q_K$ and $0 = p_{-1} \leq p_0 \leq \cdots \leq p_K$ such that

$$\mathsf{q} = \sum_{k=0}^{K} q_k \mathbf{1}_{[\zeta_k, \zeta_{k+1})} \quad \text{and} \quad \mathsf{q}' = \sum_{k=0}^{K} p_k \mathbf{1}_{[\zeta_k, \zeta_{k+1})}.$$

For each $s \in [0,1]$ and $k \in \{0,\ldots,K\}$, we define $q_k^s := sq_k + (1-s)p_k$, and the associated piecewise-constant path

$$\mathsf{q}^s := \sum_{k=0}^{K} q_k^s \mathbf{1}_{[\zeta_k, \zeta_{k+1})}.$$

Notice that we can choose the parameters $(q_k)_{0 \leq k \leq K}$ and $(p_k)_{0 \leq k \leq K}$ in such a way that for every $k \in \{0,\ldots,K\}$, we have $q_{k-1} \neq q_k$ or $p_{k-1} \neq p_k$. This ensures that for every $s \in (0,1)$, the sequence $(q_k^s)_{0 \leq k \leq K}$ is strictly increasing. We set

$$\varphi(s) := \partial_x \Phi^{\mathsf{q}^s}(0,x) = \mathbb{E}\langle \sigma \rangle_s,$$

where $\langle \cdot \rangle_s$ denotes the Gibbs average associated with the Hamiltonian $(\sigma, \alpha) \mapsto \sigma(x + Z_{\mathsf{q}^s}(\alpha))$. A direct computation and the chain rule give

$$\varphi'(s) = \mathbb{E}\langle \partial_s Z_{\mathsf{q}^s}(\alpha^1) - \sigma^1 \sigma^2 \partial_s Z_{\mathsf{q}^s}(\alpha^2) \rangle_s$$

$$= \sum_{k=0}^{K} (2q_k^s - 2q_{k-1}^s)^{-1/2}\big((q_k - p_k) - (q_{k-1} - p_{k-1})\big)\mathbb{E}\langle z_{\alpha^1_{|k}} - \sigma^1 \sigma^2 z_{\alpha^2_{|k}} \rangle_s.$$

The Gibbs Gaussian integration by parts formula (Theorem 4.6) and the symmetry between replicas reveal that

$$\mathbb{E}\langle z_{\alpha^1_{|k}} \rangle_s = (2q_k^s - 2q_{k-1}^s)^{1/2}\mathbb{E}\langle \sigma^1 - \sigma^2 \mathbf{1}_{\{\alpha^1_{|k}=\alpha^2_{|k}\}} \rangle_s$$

$$\mathbb{E}\langle \sigma^1 \sigma^2 z_{\alpha^2_{|k}} \rangle_s = (2q_k^s - 2q_{k-1}^s)^{1/2}\mathbb{E}\langle \sigma^1 + \sigma^2 \mathbf{1}_{\{\alpha^1_{|k}=\alpha^2_{|k}\}} - 2\sigma^1 \sigma^2 \sigma^3 \mathbf{1}_{\{\alpha^1_{|k}=\alpha^2_{|k}\}} \rangle_s.$$

It follows that

$$\varphi'(s) = 2\sum_{k=0}^{K} \big((q_k - p_k) - (q_{k-1} - p_{k-1})\big)\mathbb{E}\langle(\sigma^1\sigma^2\sigma^3 - \sigma^1)\mathbf{1}_{\{\alpha^1_{|k}=\alpha^2_{|k}\}}\rangle_s$$

$$= 2\sum_{k=0}^{K-1}(q_k - p_k)\mathbb{E}\langle(\sigma^1\sigma^2\sigma^3 - \sigma^1)(\mathbf{1}_{\{\alpha^1_{|k}=\alpha^2_{|k}\}} - \mathbf{1}_{\{\alpha^1_{|k+1}=\alpha^2_{|k+1}\}})\rangle_s$$

$$\phantom{= 2\sum_{k=0}^{K-1}} + 2(q_K - p_K)\mathbb{E}\langle(\sigma^1\sigma^2\sigma^3 - \sigma^1)\mathbf{1}_{\{\alpha^1_{|K}=\alpha^2_{|K}\}}\rangle_s$$

$$= 2\sum_{k=0}^{K}(q_k - p_k)\mathbb{E}\langle(\sigma^1\sigma^2\sigma^3 - \sigma^1)\mathbf{1}_{\{\alpha^1 \wedge \alpha^2 = k\}}\rangle_s,$$

where we used that $q_{-1} = p_{-1} = 0$ in the second equality. Remembering that the property (6.60) is valid under very general conditions including our present one, we can rewrite this as

$$\varphi'(s) = 2 \sum_{k=0}^{K} (q_k - p_k)(\zeta_{k+1} - \zeta_k) \mathbb{E} \langle \sigma^1 \sigma^2 \sigma^3 - \sigma^1 \mid \alpha^1 \wedge \alpha^2 = k \rangle_s.$$

We can finally use the fundamental theorem of calculus to conclude that

$$\left| \partial_x \Phi^{\mathsf{q}}(0,x) - \partial_x \Phi^{\mathsf{q}'}(0,x) \right| \leqslant 4 \sum_{k=0}^{K} |q_k - p_k|(\zeta_{k+1} - \zeta_k) = 4 \int_0^1 |\mathsf{q}(u) - \mathsf{q}'(u)| \, du.$$

This completes the proof of (iii). ∎

Although we will not make use of this fact, we note in passing that Lemma 6.5 can be generalized, with essentially no change to the proof, to yield that for every integer $k \geqslant 1$, the derivative $\partial_x^k \Phi^{\mathsf{q}}(t,x)$ is bounded uniformly over the choice of piecewise-constant $\mathsf{q} \in Q(\mathbb{R}_{\geqslant 0})$, $t \geqslant 0$ and $x \in \mathbb{R}$, and the mapping of $\mathsf{q} \mapsto \partial_x^k \Phi^{\mathsf{q}}(t,x)$ is Lipschitz continuous with respect to the $L^1$ norm, uniformly over $t \geqslant 0$ and $x \in \mathbb{R}$.

We can now extend the identity (6.89) from Lemma 6.4 to arbitrary bounded paths.

**Proposition 6.6.** *For every path $\mathsf{q} \in Q_\infty(\mathbb{R}_{\geqslant 0})$, the initial condition (6.84) is given by*

$$\psi(\mathsf{q}) = \overline{F}_1(0, \mathsf{q}) = -\Phi^{\mathsf{q}}(0,0), \tag{6.100}$$

*where the function $\Phi^{\mathsf{q}} : [0, \mathsf{q}(1)] \times \mathbb{R} \to \mathbb{R}$ satisfies the equation (6.90).*

*Proof.* We recall that, as explained around (6.95), for any $\mathsf{q} \in Q_\infty(\mathbb{R}_{\geqslant 0})$, we may as well think of $\Phi^{\mathsf{q}}$ as being defined on $\mathbb{R}_{\geqslant 0} \times \mathbb{R}$. We fix an arbitrary path $\mathsf{q} \in Q_\infty(\mathbb{R}_{\geqslant 0})$, and take a sequence of piecewise-constant paths $(\mathsf{q}_n)_{n \geqslant 1}$ that converge to $\mathsf{q}$ in the $L^1$ norm. For convenience we may also assume that for every $n \geqslant 1$, we have $\mathsf{q}_n(1) = \mathsf{q}(1)$. By (ii) in Lemma 6.5, the functions $\Phi^{\mathsf{q}_n}$ converge pointwise to some limit; we define $\Phi^{\mathsf{q}}$ to be this limit, and proceed to show that this limit solves (6.90) in the sense of satisfying (6.91). Since $\mathsf{q}_n(1) = \mathsf{q}(1)$, we have for every $n \geqslant 1$ and $t \in [0, \mathsf{q}(1)]$ that

$$\Phi^{\mathsf{q}_n}(t,x) = e^{(\mathsf{q}(1)-t)\Delta} \phi(\mathsf{q}(1), \cdot)(x)$$
$$+ \int_t^{\mathsf{q}(1)} \mathsf{q}_n^{-1}(s) e^{(s-t)\Delta} (\partial_x \Phi^{\mathsf{q}_n}(s, \cdot))^2 (x) \, ds. \tag{6.101}$$

By (iii) in Lemma 6.5, the sequence $(\partial_x \Phi^{\mathsf{q}_n})_{n \geqslant 1}$ converges to some limit, uniformly over $\mathbb{R}_{\geqslant 0} \times \mathbb{R}$. Using the uniform bound on $\partial_x^2 \Phi^{\mathsf{q}_n}$ from (i) of Lemma 6.5, we see that the functions $\partial_x \Phi^{\mathsf{q}_n}$ are also Lipschitz continuous in $x$, uniformly over $\mathbb{R}_{\geqslant 0} \times \mathbb{R}$ and over $n \geqslant 1$. It is therefore clear that $\Phi^q$ is differentiable in $x$, that the limit of $\partial_x \Phi^{q_n}$ as $n$ tends to infinity is $\partial_x \Phi^q$, and since $q_n$ converges to $q$ in the $L^1$ norm, that

$$\lim_{n \to +\infty} \|\partial_x \Phi^{q_n} - \partial_x \Phi^q\|_{L^\infty(\mathbb{R}_{\geqslant 0} \times \mathbb{R}; \mathbb{R})} = 0. \tag{6.102}$$

Observing also that, by Fubini's theorem (see e.g. Proposition 2.17 of [235] for a detailed argument),

$$\int_0^1 |q_n(u) - q(u)| \, du = \int_0^{+\infty} |q_n^{-1}(s) - q^{-1}(s)| \, ds, \tag{6.103}$$

we also have that the right side of (6.103) converges to 0 as $n$ tends to infinity. Combining this with (6.102), we conclude from (6.101) that $\Phi^q$ is indeed a mild solution to (6.90). Recalling from Lemma 6.4 that for each $n \geqslant 1$, we have

$$\overline{F}_1(0, q_n) = -\Phi^{q_n}(0, 0),$$

and that $\overline{F}_1(0, \cdot)$ is Lipschitz continuous by Proposition 6.3, this completes the proof. ∎

We recall from (6.81) that, through informal arguments, we suggested that the limit $f$ of the free energy $\overline{F}_N$ might satisfy the Hamilton-Jacobi equation

$$\partial_t f(t, q) - \int_0^1 \partial_q f(t, q, u)^2 \, du = 0 \quad \text{on} \quad \mathbb{R}_{>0} \times Q_2(\mathbb{R}_{\geqslant 0}), \tag{6.104}$$

with an initial condition which we have now firmly identified as

$$f(0, q) = \psi(q) = -\Phi^q(0, 0). \tag{6.105}$$

While we have only really made sense of $\Phi^q(0, 0)$ for $q \in Q_\infty(\mathbb{R}_{\geqslant 0})$, via the partial differential equation in (6.90), we recall that $-\Phi^q(0, 0) = \overline{F}_1(0, q)$ and that $\overline{F}_1(0, \cdot)$ is Lipschitz continuous with respect to the $L^1$ norm, by Proposition 6.3, so we can extend $\Phi^q(0, 0)$ to every $q \in Q_1(\mathbb{R}_{\geqslant 0})$ by continuity. The main reason we prefer to think of (6.104) as being posed on $\mathbb{R}_{>0} \times Q_2(\mathbb{R}_{\geqslant 0})$ as opposed to, for instance, being posed on $\mathbb{R}_{>0} \times Q_\infty(\mathbb{R}_{\geqslant 0})$, is that the space $Q_2(\mathbb{R}_{\geqslant 0})$ is naturally embedded in the Hilbert space $L^2([0, 1]; \mathbb{R})$, and that the Hilbert-space geometry of $L^2([0, 1]; \mathbb{R})$ is directly tied with our notion of derivative in (6.74), and is also closely related to the Euclidean-space geometry of our finite-dimensional approximations to (6.104). This is not fundamental however, and we could as well decide to try to make sense of (6.104) as being posed on $\mathbb{R}_{>0} \times Q_\infty(\mathbb{R}_{\geqslant 0})$, which we might still think of as a subset of $L^2([0, 1]; \mathbb{R})$ if we want to.

As was discussed at the opening of this section, there are several different options for us to try to make sense of the equation (6.104). Since the non-linearity in the equation is convex, we may expect that whichever notion of solution we come up with will also admit a variational representation analogous to the Hopf-Lax formula. There are however two potential sources of problem for the validity of such a variational representation. The first one is the infinite-dimensional nature of the equation. The second one is that the space $Q_2(\mathbb{R}_{\geq 0})$ of non-decreasing paths from $[0,1)$ to $\mathbb{R}_{\geq 0}$ has a boundary: if we restrict to piecewise-constant paths of the form (6.50), we see that we must impose that the sequence $(q_k)_{0 \leq k \leq K}$ be a sequence of non-negative numbers, and be non-decreasing.

If we simplify this further and only consider the space of constant paths, which corresponds to the function introduced in (6.32), then this reduces to the constraint that the constant path must be non-negative, that is, $h \geq 0$ in the notation from (6.32). In this simplest case of a domain with a boundary, we have seen in Chapter 4 that, under some assumptions on the initial condition $\psi_1 : \mathbb{R} \to \mathbb{R}$, we can make sense of the solution to

$$\partial_t f_1 - (\partial_h f_1)^2 = 0 \quad \text{on} \quad \mathbb{R}_{>0} \times \mathbb{R}_{\geq 0} \tag{6.106}$$

subject to the initial condition $f_1(0, \cdot) = \psi_1$, and that the solution is given, for every $t, h \geq 0$, by

$$f_1(t, h) = \sup_{h' \geq 0} \left( \psi_1(h + h') - \frac{(h')^2}{4t} \right); \tag{6.107}$$

see in particular (4.84) and (4.90). While we exploited the convexity of $\psi_1$ to make use of the convex selection principle there, the fact that we could essentially disregard the contribution of the boundary was only really using that the function $\psi_1$ was non-decreasing. Using the language of characteristics discussed in Section 3.5, the point is to ensure that the characteristic lines always go from the inside of the domain towards the boundary, never the other way around.

In our present setting, a similar monotonicity property also holds. If we only consider constant paths, with the free energy as in (6.32), then we found in (6.34) that the derivative with respect to the value $h$ of the constant path is given by

$$N^{-1}\mathbb{E}\langle \sigma^1 \cdot \sigma^2 \rangle = N^{-1}\mathbb{E}\langle \sigma \rangle^2 \geq 0. \tag{6.108}$$

We fix an integer $K \geq 1$, denote

$$U_K := \left\{ q = (q_0, \ldots, q_K) \in \mathbb{R}_{\geq 0}^{K+1} \mid 0 \leq q_0 \leq q_1 \leq \cdots \leq q_K \right\}, \tag{6.109}$$

and for each sequence $q = (q_k)_{0 \leq k \leq K} \in U_K$, we consider the piecewise-constant path

$$\mathsf{q}_q := \sum_{k=0}^{K} q_k \mathbf{1}_{\left[\frac{k}{K+1}, \frac{k+1}{K+1}\right)}. \tag{6.110}$$

While so far we had kept the dependence of $\mathsf{q}$ on its parameters implicit, it will be useful to make this dependence explicit for a short while, as we do in the notation $\mathsf{q}_q$ instead of just q. One can show that

$$0 \leqslant \partial_{q_0}\overline{F}_N(t,\mathsf{q}_q) \leqslant \cdots \leqslant \partial_{q_K}\overline{F}_N(t,\mathsf{q}_q). \tag{6.111}$$

The proof of this result combines the Gaussian integration by parts used for the proof of Lemma 6.2 with the explicit procedure for integrating the Poisson-Dirichlet cascade presented in Theorem 5.25 and Corollary 5.26. We will not give a detailed proof of (6.111), and only refer to Lemma 2.4 of [196] or Proposition 14.3.2 of [254]. (An alternative proof can be devised by analyzing the properties of the equation in (6.90).) In analogy with (4.84) and (4.90), we can therefore hope that there is a nice variational formula for equations posed on $U_K$, even though this space has a non-trivial boundary.

Recall from Lemma 6.2 and the discussion that followed that if $K$ is chosen sufficiently large, then we expect that

$$\partial_t \overline{F}_N(t,\mathsf{q}_q) - (K+1)\sum_{k=0}^{K}(\partial_{q_k}\overline{F}_N(t,\mathsf{q}_q))^2 \simeq 0. \tag{6.112}$$

Since the convex dual of the mapping $(p_k)_{0\leqslant k \leqslant K} \mapsto (K+1)\sum_{k=0}^{K} p_k^2$ is the mapping $(q_k)_{0\leqslant k \leqslant K} \mapsto \frac{1}{4(K+1)}\sum q_k^2$, and since we argued that boundary aspects should hopefully behave as nicely as they have in the context of (4.90), we are led to expect that, for $K$ sufficiently large and $\mathsf{q}_q$ as in (6.110),

$$f(t,\mathsf{q}_q) \simeq \sup_{q' \in U_K}\left(\psi(\mathsf{q}_{q+q'}) - \frac{1}{4t(K+1)}\sum_{k=0}^{K}(q'_k)^2\right). \tag{6.113}$$

Using also that $\mathsf{q}_{q+q'} = \mathsf{q}_q + \mathsf{q}_{q'}$, this can be rewritten as

$$f(t,\mathsf{q}_q) \simeq \sup_{q' \in U_K}\left(\psi(\mathsf{q}_q + \mathsf{q}_{q'}) - \frac{1}{4t}\int_0^1 (\mathsf{q}_{q'}(u))^2 \, du\right). \tag{6.114}$$

Sending $K$ to infinity, we see that the natural candidate for a Hopf-Lax variational representation of the solution to (6.104)-(6.105) is, for every $t \geqslant 0$ and $\mathsf{q} \in \mathcal{Q}_2(\mathbb{R}_{\geqslant 0})$,

$$f(t,\mathsf{q}) = \sup_{\mathsf{q}' \in \mathcal{Q}_2(\mathbb{R}_{\geqslant 0})}\left(\psi(\mathsf{q}+\mathsf{q}') - \frac{1}{4t}\int_0^1 (\mathsf{q}'(u))^2 \, du\right). \tag{6.115}$$

We note that for the more general mixed $p$-spin models defined in (6.13)-(6.14), the relation (6.61) would be replaced by

$$\partial_t \overline{F}_N(t,\mathsf{q}) = \mathbb{E}\langle \xi(R_{1,2})\rangle, \tag{6.116}$$

while the derivatives with respect to each $q_k$ are still as in (6.62). In this context, we would therefore be led to suspect that the limit free energy $f$ might satisfy[1]

$$\partial_t f(t,\mathsf{q}) - \int_0^1 \xi(\partial_\mathsf{q} f(t,\mathsf{q},u))\,\mathrm{d}u = 0 \quad \text{on} \quad \mathbb{R}_{>0} \times \mathsf{Q}_2(\mathbb{R}_{\geqslant 0}), \tag{6.117}$$

with the same initial condition (6.105). The inequality (6.111) remains valid in this case, and in particular, the argument of the function $\xi$ in (6.117) is always non-negative. This is why the convexity of $\xi$ on $\mathbb{R}_{\geqslant 0}$ is sufficient for the validity of a Hopf-Lax variational representation for $f$. Following along the heuristic argument above would thus suggest that

$$f(t,\mathsf{q}) = \sup_{\mathsf{q}' \in \mathsf{Q}_\infty(\mathbb{R}_{\geqslant 0})} \left( \psi(\mathsf{q}+\mathsf{q}') - t \int_0^1 \xi^*\!\left(\frac{\mathsf{q}'(u)}{t}\right) \mathrm{d}u \right), \tag{6.118}$$

where $\xi^* : \mathbb{R} \to \mathbb{R}$ is the convex dual of $\xi : \mathbb{R}_{\geqslant 0} \to \mathbb{R}$, so that for every $s \in \mathbb{R}$,

$$\xi^*(s) := \sup_{r \geqslant 0}(rs - \xi(r)). \tag{6.119}$$

In this section we take the (rather convenient!) point of view of simply *defining* the solution to (6.104)-(6.105) as being the function given in (6.115). The main goal of this section is to show that our guessed formula (6.115) is indeed consistent with the Parisi formula for the free energy of the SK model that was presented in (6.8)-(6.9). We recall from (6.57) that the free energy (6.7) is equal to

$$-\overline{F}_N\!\left(\frac{\beta^2}{2},0\right) + \log(2) + \frac{\beta^2}{2}, \tag{6.120}$$

with $\overline{F}_N(t,\mathsf{q})$ being the enriched free energy defined in (6.53) with the choice of $P_N = (\frac{1}{2}\delta_{-1} + \frac{1}{2}\delta_1)^{\otimes N}$. If we believe that $\overline{F}_N(t,\mathsf{q})$ converges to $f(t,\mathsf{q})$ as defined in (6.115) for arbitrary choices of $t \geqslant 0$ and $\mathsf{q} \in \mathsf{Q}_2(\mathbb{R}_{\geqslant 0})$, then in particular, it should be that the limit of the free energy in (6.7) is

$$-f\!\left(\frac{\beta^2}{2},0\right) + \log(2) + \frac{\beta^2}{2}. \tag{6.121}$$

We now show that this guess indeed coincides with the Parisi formula (6.8).

**Theorem 6.7** ([196]). *Let $f : \mathbb{R}_{\geqslant 0} \times \mathsf{Q}_2(\mathbb{R}_{\geqslant 0}) \to \mathbb{R}$ be the function defined by (6.115), which we interpret as the solution to the Hamilton-Jacobi equation (6.104)-(6.105). For every $\beta \geqslant 0$, we have*

$$-f\!\left(\frac{\beta^2}{2},0\right) + \log(2) + \frac{\beta^2}{2} = \inf_{\zeta \in \mathcal{D}[0,1]} \left( \Phi_\zeta(0,0) - \beta^2 \int_0^1 t\zeta(t)\,\mathrm{d}t + \log(2) \right), \tag{6.122}$$

---

[1] The reader who wonders whether the integral in (6.117) is finite can derive from (6.62) that, assuming $\partial_\mathsf{q} f(t,\mathsf{q},\cdot)$ is well-defined, it must be uniformly bounded.

*where the space* $\mathcal{D}[0,1]$ *is defined in* (6.10) *and the function* $\Phi_\zeta : [0,1] \times \mathbb{R} \to \mathbb{R}$ *denotes the solution to the Parisi PDE* (6.9). *In particular, the limit of the free energy* (6.7) *in the SK model is given by*

$$\lim_{N \to +\infty} \overline{F}_N(\beta) = -f\left(\frac{\beta^2}{2}, 0\right) + \log(2) + \frac{\beta^2}{2}. \tag{6.123}$$

*Proof.* The Hopf-Lax representation (6.115) implies that

$$-f\left(\frac{\beta^2}{2}, 0\right) + \frac{\beta^2}{2} = \inf_{\mathsf{q} \in Q_2(\mathbb{R}_{\geqslant 0})} \left(\Phi^{\mathsf{q}}(0,0) + \frac{1}{2\beta^2} \int_0^1 |\mathsf{q}(u)|^2 \,\mathrm{d}u + \frac{\beta^2}{2}\right). \tag{6.124}$$

To establish (6.122), for each $a > 0$ it will be convenient to write

$$Q_{\leqslant a} := \{\mathsf{q} \in Q_\infty(\mathbb{R}_{\geqslant 0}) \mid \mathsf{q}(1) \leqslant a\} \tag{6.125}$$

for the paths in $Q_\infty(\mathbb{R}_{\geqslant 0})$ that are bounded by $a$. The proof now proceeds in three steps. First we show that the infimum on the right side of (6.124) can be restricted to paths in $Q_{\leqslant \beta^2}$, then we use a change of variables to normalize to paths in $Q_{\leqslant 1}$, and finally we perform another change of variables to replace a path $\mathsf{q} \in Q_{\leqslant 1}$ by its right-continuous inverse $\mathsf{q}^{-1} \in \mathcal{D}[0,1]$.

*Step 1: restricting to* $Q_{\leqslant \beta^2}$. We fix a path $\mathsf{q} \in Q_2(\mathbb{R}_{\geqslant 0})$, and define the path $\widetilde{\mathsf{q}} \in Q_{\leqslant \beta^2}$ by $\widetilde{\mathsf{q}} := \min(\mathsf{q}, \beta^2)$. We denote by $u^* := \inf\{u \in [0,1] \mid \mathsf{q}(u) \geqslant \beta^2\}$ the first point at which $\mathsf{q}$ exceeds $\beta^2$, and observe that by the Lipschitz property of $\Phi^\mathsf{q}$ in Proposition 6.3,

$$\left(\Phi^{\widetilde{\mathsf{q}}}(0,0) + \frac{1}{2\beta^2} \int_0^1 |\widetilde{\mathsf{q}}(u)|^2 \,\mathrm{d}u\right) - \left(\Phi^{\mathsf{q}}(0,0) + \frac{1}{2\beta^2} \int_0^1 |\mathsf{q}(u)|^2 \,\mathrm{d}u\right)$$

$$\leqslant \int_0^1 |\mathsf{q}(u) - \widetilde{\mathsf{q}}(u)| \,\mathrm{d}u - \frac{1}{2\beta^2} \int_0^1 \left(|\mathsf{q}(u)|^2 - |\widetilde{\mathsf{q}}(u)|^2\right) \mathrm{d}u$$

$$= \int_{u^*}^1 \left(\mathsf{q}(u) - \beta^2\right) \mathrm{d}u - \int_{u^*}^1 \left(\mathsf{q}(u) - \beta^2\right) \frac{\mathsf{q}(u) + \beta^2}{2\beta^2} \,\mathrm{d}u$$

$$\leqslant 0,$$

where we used that $\mathsf{q}(u) \geqslant \beta^2$ for $u \geqslant u^*$ in the final inequality. Together with (6.124), this implies that

$$-f\left(\frac{\beta^2}{2}, 0\right) + \frac{\beta^2}{2} = \inf_{\mathsf{q} \in Q_{\leqslant \beta^2}} \left(\Phi^{\mathsf{q}}(0,0) + \frac{1}{2\beta^2} \int_0^1 |\mathsf{q}(u)|^2 \,\mathrm{d}u + \frac{\beta^2}{2}\right). \tag{6.126}$$

*Step 2: re-normalizing to* $Q_{\leqslant 1}$. Substituting $\mathsf{q}$ with $\beta^2 \mathsf{q}$ in the previous display, we obtain that

$$-f\left(\frac{\beta^2}{2}, 0\right) + \frac{\beta^2}{2} = \inf_{\mathsf{q} \in Q_{\leqslant 1}} \left(\Phi^{\beta^2 \mathsf{q}}(0,0) + \frac{\beta^2}{2} \int_0^1 |\mathsf{q}(u)|^2 \,\mathrm{d}u + \frac{\beta^2}{2}\right). \tag{6.127}$$

Our goal now is to relate, for each fixed $q \in Q_{\leq 1}$, the function $\Phi^{\beta^2 q}$ from (6.90) with the function $\Phi_\zeta$ solution to (6.9) for a suitable choice of $\zeta \in \mathcal{D}[0,1]$. We choose $\zeta := q^{-1}$, the right-continuous inverse of $q$. By definition, the function $\Phi_{q^{-1}} : [0,1] \times \mathbb{R} \to \mathbb{R}$ solves

$$\begin{cases} -\partial_t \Phi_{q^{-1}}(t,x) = \beta^2 \left( \partial_x^2 \Phi_{q^{-1}}(t,x) + q^{-1}(t) \left( \partial_x \Phi_{q^{-1}}(t,x) \right)^2 \right) & \text{on } [0,1] \times \mathbb{R}^d \\ \Phi_{q^{-1}}(1,x) = \log \cosh(x) & \text{for } x \in \mathbb{R}. \end{cases}$$

We write $\widetilde{q} := \beta^2 q$, and let $\widetilde{\Phi} : [0,\beta^2] \times \mathbb{R} \to \mathbb{R}$ be defined by $\widetilde{\Phi}(t,x) := \Phi_{q^{-1}}(t/\beta^2, x)$. Using that $\widetilde{q}^{-1}(t) = q^{-1}(t/\beta^2)$, we obtain that $\widetilde{\Phi}$ satisfies

$$\begin{cases} -\partial_t \widetilde{\Phi}(t,x) = \partial_x^2 \widetilde{\Phi}(t,x) + \widetilde{q}^{-1}(t) \left( \partial_x \widetilde{\Phi}(t,x) \right)^2 & \text{on } [0,\beta^2] \times \mathbb{R}^d \\ \widetilde{\Phi}(\beta^2, x) = \log \cosh(x) & \text{for } x \in \mathbb{R}. \end{cases}$$

We recall from the discussion around (6.95) that we may as well think of the function $\Phi^{\beta^2 q} = \Phi^{\widetilde{q}}$ as being defined on $[0,\beta^2] \times \mathbb{R}$, since $\beta^2 \geq \beta^2 q(1)$. We also recall that the addition of a constant to the initial condition in (6.90) simply changes the solution by the addition of this constant. Using also (6.88), we therefore obtain that, for every $t \in [0,\beta^2]$ and $x \in \mathbb{R}$,

$$\Phi^{\beta^2 q}(t,x) = \widetilde{\Phi}(t,x) - \beta^2.$$

Combining this with (6.127) yields that

$$-f\left(\frac{\beta^2}{2}, 0\right) + \frac{\beta^2}{2} = \inf_{q \in Q_{\leq 1}} \left( \Phi_{q^{-1}}(0,0) + \frac{\beta^2}{2} \int_0^1 |q(u)|^2 \, du - \frac{\beta^2}{2} \right). \tag{6.128}$$

*Step 3: re-parametrizing by $\mathcal{D}[0,1]$.* Our next step consists in replacing the optimization variable $q \in Q_{\leq 1}$ in (6.128) by its right-continuous inverse. So we fix $q \in Q_{\leq 1}$ and denote by $\zeta : [0,1] \to [0,1]$ its right-continuous inverse, that is, for every $t \in [0,1]$,

$$\zeta(t) := \sup\{u \in [0,1] \mid q(u) \leq t\}. \tag{6.129}$$

Since $q$ is non-decreasing, we have that for every $u, t \in [0,1]$,

$$u < \zeta(t) \implies q(u) \leq t \quad \text{and} \quad u > \zeta(t) \implies q(u) > t.$$

Since $\zeta \in \mathcal{D}[0,1]$, we can write this function as $t \mapsto \mu[0,t]$ for some probability measure $\mu$ on $[0,1]$. We infer from the previous display that, for every $t \in [0,1]$,

$$\mu[0,t] = \int_0^1 \mathbf{1}_{\{q(u) \leq t\}} \, du. \tag{6.130}$$

The mapping $A \mapsto \int_0^1 \mathbf{1}_{\{q(u) \in A\}} \, du$, where $A$ is any measurable subset of $[0,1]$, defines a probability measure on $[0,1]$; it is the image of the Lebesgue measure on $[0,1]$ through the mapping q. By (6.130), this probability measure has the same cumulative distribution function as $\mu$. By Dynkin's $\pi$-$\lambda$ theorem (see Theorem A.5 and Exercise A.3), we deduce that this measure is equal to the measure $\mu$. That is, for every bounded measurable function $g : [0,1] \to \mathbb{R}$,

$$\int_0^1 g(q(u)) \, du = \int_0^1 g(t) \, d\mu(t). \tag{6.131}$$

This corresponds to the fact that if we sample a uniform random variable $U$ on $[0,1]$, then the law of $q(U)$ is $\mu$. We use (6.131) with the function $g$ replaced by the square function to obtain that

$$\int_0^1 |q(u)|^2 \, du = \int_0^1 t^2 \, d\mu(t).$$

We then re-express the latter integral in terms of $\zeta$ by integrating by parts, that is,

$$\int_0^1 |q(u)|^2 \, du = 1 - 2 \int_0^1 \int_0^1 s \mathbf{1}_{\{t \leq s\}} \, ds \, d\mu(t)$$

$$= 1 - 2 \int_0^1 s \int_0^1 \mathbf{1}_{\{t \leq s\}} \, d\mu(t) \, ds$$

$$= 1 - 2 \int_0^1 s \zeta(s) \, ds.$$

Since the correspondence between $q \in \mathcal{Q}_{\leq 1}$ and $\zeta \in \mathcal{D}[0,1]$ is bijective, we can perform the change of variables $q \to \zeta$ in (6.128). Using also the identity in the previous display, we obtain that

$$-f\left(\frac{\beta^2}{2}, 0\right) + \frac{\beta^2}{2} = \inf_{\zeta \in \mathcal{D}[0,1]} \left( \Phi_\zeta(0,0) - \beta^2 \int_0^1 t \zeta(t) \, dt \right), \tag{6.132}$$

as announced. ∎

## 6.6  Towards non-convex models

Let us start by summarizing what we have done so far. For the SK model, we argued for the introduction of an enriched free energy $\overline{F}_N(t,q)$ which depends on a variable $t \geq 0$ and a path $q \in \mathcal{Q}_1(\mathbb{R}_{\geq 0})$. This free energy is defined precisely in (6.53) for piecewise-constant paths, and then extended by continuity. The quantity we primarily want to compute is the large-$N$ limit of $\overline{F}_N(t,0)$, but the introduction of the additional variable is likely to help us to discover this limit. Indeed, we argued informally that if $f = f(t, \mathsf{q}) : \mathbb{R}_{\geqslant 0} \times \mathcal{Q}_2(\mathbb{R}_{\geqslant 0}) \to \mathbb{R}$ denotes the limit of the free energy, we expect $f$ to be a solution to

$$\partial_t f(t,\mathsf{q}) - \int_0^1 \partial_\mathsf{q} f(t,\mathsf{q},u)^2 \, \mathrm{d}u = 0 \quad \text{on} \quad \mathbb{R}_{>0} \times \mathcal{Q}_2(\mathbb{R}_{\geqslant 0}). \tag{6.133}$$

Moreover, thanks to the fact that the reference measure $P_N$ is a product measure, we have that $\overline{F}_N(0,\mathsf{q}) = \overline{F}_1(0,\mathsf{q})$ for every $N$, so the initial condition for $f$ is

$$f(0,\mathsf{q}) = \psi(\mathsf{q}) := \overline{F}_1(0,\mathsf{q}). \tag{6.134}$$

This object is computable relatively explicitly, since it only involves averages with respect to a single spin variable. For piecewise-constant paths q, we can compute this quantity recursively using Theorem 5.25, and then use the continuity of $\psi$ for more general paths. Alternatively, we can also express this calculation in the form of the parabolic equation in (6.90).

Since the non-linearity in (6.133) is convex, we expect that the solution to this equation admits a Hopf-Lax variational representation, which should take the form

$$f(t,\mathsf{q}) = \sup_{\mathsf{q}' \in \mathcal{Q}_2(\mathbb{R}_{\geqslant 0})} \left( \psi(\mathsf{q}+\mathsf{q}') - \frac{1}{4t} \int_0^1 (\mathsf{q}'(u))^2 \, \mathrm{d}u \right). \tag{6.135}$$

And indeed, when setting $\mathsf{q} = 0$ in this formula, we verified that this guess yields the answer predicted by Parisi in [217, 218, 219] and proved rigorously in [132, 250]. In fact, the convergence of the enriched free energy $\overline{F}_N(t,\mathsf{q})$ to the function $f(t,\mathsf{q})$ defined in (6.133) can be shown rigorously for arbitrary choices of $\mathsf{q} \in \mathcal{Q}_2(\mathbb{R}_{\geqslant 0})$.

**Theorem 6.8** ([198]). *Let $\overline{F}_N(t,\mathsf{q})$ be the enriched free energy defined in (6.53) and extended by continuity in Proposition 6.3, and let $f(t,\mathsf{q})$ be as in (6.135). For every $t \geqslant 0$ and $\mathsf{q} \in \mathcal{Q}_2(\mathbb{R}_{\geqslant 0})$, we have*

$$\lim_{N \to +\infty} \overline{F}_N(t,\mathsf{q}) = f(t,\mathsf{q}). \tag{6.136}$$

In fact, this result was shown to be valid for arbitrary choices of the measure $P_1$ of compact support and for general mixed $p$-spin models, provided that we replace the definition of the Hamiltonian in (6.54) by

$$H_N(t,\mathsf{q},\sigma,\alpha) := \sqrt{2t} H_N(\sigma) - Nt\xi\left(\frac{|\sigma|^2}{N}\right) + Z_\mathsf{q}(\alpha) \cdot \sigma - q_K |\sigma|^2, \tag{6.137}$$

and with $f$ as in (6.118).

The main motivation for exploring this chain of ideas is to make progress towards the understanding of non-convex models. A typical representative in this class is the bipartite model defined in (6.18). We recall that for this model, the configuration vector $\sigma$ is split into two parts, $\sigma = (\sigma_1, \sigma_2) \in \Sigma_N^2$. While we keep the length of each of these two components the same for convenience of notation, we point out that we do not plan to exploit any symmetry between these two components; to make this clearest, one can think that under the reference measure $P_N$, the vectors $(\sigma_{1,i})_{1 \leq i \leq N}$ and $(\sigma_{2,i})_{1 \leq i \leq N}$ are independent with i.i.d. coordinates, but that the mean of $\sigma_{1,i}$ may be different from the mean of $\sigma_{2,i}$.

For this model, the naive extension of the free energy we explored in Section 6.2 would already involve the choice of two additional parameters instead of one, so that we can act on each of the two coordinates of $\sigma$ separately. Explicitly, the analogue of the free energy in (6.32) would be, for every $t \geq 0$ and $h = (h_1, h_2) \in \mathbb{R}_{\geq 0}^2$,

$$\overline{F}_N(t,h) := -\frac{1}{N} \mathbb{E} \log \int_{\Sigma_N^2} \exp\Big(\sqrt{2t} H_N(\sigma) - Nt \\ + \sqrt{2h_1} z_1 \cdot \sigma_1 - Nh_1 + \sqrt{2h_2} z_2 \cdot \sigma_2 - Nh_2 \Big) dP_N(\sigma), \quad (6.138)$$

where $z_1$ and $z_2$ are independent standard $N$-dimensional Gaussian vectors, independent from $(H_N(\sigma))_{\sigma \in \Sigma_N^2}$. The same Gaussian integration by parts calculations as those in (6.34) and (6.35) yield that, for each $a \in \{1, 2\}$,

$$\partial_{h_a} \overline{F}_N(t,h) = N^{-1} \mathbb{E} \langle \sigma_a^1 \cdot \sigma_a^2 \rangle, \quad (6.139)$$

while

$$\partial_t \overline{F}_N(t,h) = N^{-2} \mathbb{E} \langle (\sigma_1^1 \cdot \sigma_1^2)(\sigma_2^1 \cdot \sigma_2^2) \rangle. \quad (6.140)$$

This suggests that we should consider the relation

$$\big(\partial_t \overline{F}_N - \partial_{h_1} \overline{F}_N \, \partial_{h_2} \overline{F}_N\big)(t,h) \\ = N^{-2} \mathbb{E} \big\langle \big(\sigma_1^1 \cdot \sigma_1^2 - \mathbb{E} \langle \sigma_1^1 \cdot \sigma_1^2 \rangle\big)\big(\sigma_2^1 \cdot \sigma_2^2 - \mathbb{E} \langle \sigma_2^1 \cdot \sigma_2^2 \rangle\big) \big\rangle. \quad (6.141)$$

This term would be small if the overlaps were concentrated. Indeed, using that $|xy| \leq (x^2 + y^2)/2$, we see that

$$\big|\big(\partial_t \overline{F}_N - \partial_{h_1} \overline{F}_N \, \partial_{h_2} \overline{F}_N\big)(t,h)\big| \leq \frac{1}{2N^2} \sum_{a=1}^{2} \mathbb{E} \big\langle \big(\sigma_a^1 \cdot \sigma_a^2 - \mathbb{E} \langle \sigma_a^1 \cdot \sigma_a^2 \rangle\big)^2 \big\rangle. \quad (6.142)$$

However, as for the SK model, these overlaps will typically not have a small variance when $t$ is large. As in Section 6.4, we therefore need to define a more sophisticated enriched free energy. To do so, we give ourselves an integer $K \geq 1$ and parameters

$$0 = \zeta_0 < \zeta_1 < \cdots < \zeta_K < \zeta_{K+1} = 1, \quad (6.143)$$

$$0 = q_{-1,a} \leq q_{0,a} \leq q_{1,a} \leq \cdots \leq q_{K,a}, \quad (6.144)$$

for $a \in \{1, 2\}$. For each $k \in \{-1, \ldots, K\}$, we write $q_k := (q_{k,1}, q_{k,2}) \in \mathbb{R}^2_{\geq 0}$, and we also impose that $q_k \neq q_{k+1}$. Notice that the sequence $(q_k)_{0 \leq k \leq K}$ is increasing for the partial order induced by $\mathbb{R}^2_{\geq 0}$. We encode all these parameters into the path $\mathsf{q} : [0,1) \to \mathbb{R}^2_{\geq 0}$ defined by

$$\mathsf{q} := \sum_{k=0}^{K} q_k \mathbf{1}_{[\zeta_k, \zeta_{k+1})}. \tag{6.145}$$

We write

$$\mathsf{Q}(\mathbb{R}^2_{\geq 0}) := \{\mathsf{q} : [0,1) \to \mathbb{R}^2_{\geq 0} \mid \mathsf{q} \text{ is right-continuous and non-decreasing}\}, \tag{6.146}$$

where we understand that the notion of non-decreasingness is for the partial order induced by $\mathbb{R}^2_{\geq 0}$, that is, a path $\mathsf{q} : [0,1) \to \mathbb{R}^2_{\geq 0}$ is non-decreasing if for every $u \leq v \in [0,1)$, we have $\mathsf{q}(v) - \mathsf{q}(u) \in \mathbb{R}^2_{\geq 0}$.

We give ourselves a Poisson-Dirichlet cascade $(v_\alpha)_{\alpha \in \mathbb{N}^K}$ built over the tree $\mathcal{A}$ with the parameters in (6.143), and two families $(z_{\alpha,a})_{\alpha \in \mathcal{A}, a \in \{1,2\}}$ of independent standard $N$-dimensional random variables, independent of each other and of every other source of randomness. We then set, for each $a \in \{1, 2\}$,

$$Z_{\mathsf{q},a}(\alpha) := \sum_{k=0}^{K} (2q_{k,a} - 2q_{k-1,a})^{1/2} z_{\alpha_{|k}, a}. \tag{6.147}$$

Finally, we define the enriched free energy

$$\overline{F}_N(t, \mathsf{q}) := -\frac{1}{N} \mathbb{E} \log \int_{\Sigma_N^2} \sum_{\alpha \in \mathbb{N}^K} \exp(H_N(t, \mathsf{q}, \sigma, \alpha)) v_\alpha \, \mathrm{d}P_N(\sigma) \tag{6.148}$$

associated with the Hamiltonian

$$H_N(t, \mathsf{q}, \sigma, \alpha) := \sqrt{2t} H_N(\sigma) - Nt + \sum_{a=1}^{2} (Z_{\mathsf{q},a}(\alpha) \cdot \sigma_a - N q_{K,a}). \tag{6.149}$$

The proof of Proposition 6.3 applies to this setting as well and allows us to extend $\overline{F}_N$ to $\mathbb{R}_{\geq 0} \times \mathsf{Q}_1(\mathbb{R}^2_{\geq 0})$, where we write $\mathsf{Q}_p(\mathbb{R}^2_{\geq 0}) := \mathsf{Q}(\mathbb{R}^2_{\geq 0}) \cap L^p([0,1]; \mathbb{R}^2)$ for $p \in [1, +\infty]$. The derivative calculations from Lemma 6.2 can be performed similarly here: the relation (6.140) is still valid, while for every $k \in \{0, \ldots, K\}$ and $a \in \{1, 2\}$,

$$\partial_{q_{k,a}} \overline{F}_N(t, \mathsf{q}) = N^{-1} \mathbb{E} \langle \mathbf{1}_{\{\alpha^1 \wedge \alpha^2 = k\}} \sigma_a^1 \cdot \sigma_a^2 \rangle. \tag{6.150}$$

In particular,

$$\partial_t \overline{F}_N(t, \mathsf{q}) - \sum_{k=0}^{K} (\zeta_{k+1} - \zeta_k) \left( \frac{\partial_{q_{k,1}} \overline{F}_N(t, \mathsf{q})}{\zeta_{k+1} - \zeta_k} \right) \left( \frac{\partial_{q_{k,2}} \overline{F}_N(t, \mathsf{q})}{\zeta_{k+1} - \zeta_k} \right)$$
$$= N^{-2} \mathbb{E} \langle (\sigma_1^1 \cdot \sigma_1^2 - \mathbb{E} \langle \sigma_1^1 \cdot \sigma_1^2 \mid \alpha^1 \wedge \alpha^2 \rangle)(\sigma_2^1 \cdot \sigma_2^2 - \mathbb{E} \langle \sigma_2^1 \cdot \sigma_2^2 \mid \alpha^1 \wedge \alpha^2 \rangle) \rangle. \tag{6.151}$$

As for (6.142), this quantity can be bounded in absolute value by

$$\frac{1}{2N^2} \sum_{a=1}^{2} \mathbb{E}\left\langle \left(\sigma_a^1 \cdot \sigma_a^2 - \mathbb{E}\langle \sigma_a^1 \cdot \sigma_a^2 \mid \alpha^1 \wedge \alpha^2 \rangle\right)^2 \right\rangle. \quad (6.152)$$

Up to the addition of a small perturbation to the Hamiltonian, one can enforce the asymptotic validity of a family of Ghirlanda-Guerra identities involving different linear combinations of all the overlaps at play. This implies the ultrametricity of these overlaps, which in turn implies that they must be synchronized. In short, we can make sure that the quantity in (6.152) is small "most of the time", provided that the integer $K$ is sufficiently large. This suggests that the limit $f : \mathbb{R}_{\geq 0} \times Q_2(\mathbb{R}_{\geq 0}^2) \to \mathbb{R}$ of the free energy should satisfy the infinite-dimensional Hamilton-Jacobi equation

$$\partial_t f(t, \mathsf{q}) - \int_0^1 \partial_{\mathsf{q}_1} f(t, \mathsf{q}, u) \, \partial_{\mathsf{q}_2} f(t, \mathsf{q}, u) \, \mathrm{d}u = 0 \quad \text{on} \quad \mathbb{R}_{>0} \times Q_2(\mathbb{R}_{\geq 0}^2), \quad (6.153)$$

where we write $\mathsf{q} = (\mathsf{q}_1, \mathsf{q}_2)$ to denote the two components of the path $\mathsf{q}$, or equivalently, we write $\partial_{\mathsf{q}} f(t, \mathsf{q}, u) = (\partial_{\mathsf{q}_1} f(t, \mathsf{q}, u), \partial_{\mathsf{q}_2} f(t, \mathsf{q}, u))$. The initial condition, which we still denote by $\psi : Q_1(\mathbb{R}_{\geq 0}^2) \to \mathbb{R}$, is such that $\psi(\mathsf{q}) = \overline{F}_1(0, \mathsf{q})$, and it can be determined more explicitly using the same procedure as for the SK model.

It is at this stage that a new problem occurs. Indeed, we can no longer expect a Hopf-Lax variational representation to be valid for the solution to (6.153), since the non-linearity, which in essence is the mapping $\xi(x, y) := xy$ featuring in the covariance calculation (6.19), is neither convex nor concave. If we were to ignore this and blindly try to write down a Hopf-Lax variational formula anyway, we would need to consider the convex dual of $\xi$, which would be defined, for $a, b \in \mathbb{R}$, analogously to (6.119), as

$$\xi^*(a, b) := \sup_{x, y \geq 0} (ax + by - xy). \quad (6.154)$$

Since this function takes the value $+\infty$ on $\mathbb{R}_{\geq 0}^2 \setminus \{0\}$, pursuing this direction does not look promising. A more refined attempt in the same direction is also shown to be invalid in Subsection 6.2 of [195].

In the context of the models of statistical inference explored in Chapter 4, we could alternatively leverage the fact that the initial condition is always convex. This convexity property is valid even for the very general models discussed at the end of Section 4.3, and some of those models would be associated with a non-convex non-linearity in the equation. Also in view of the convex selection principle, the next question we need to ask is therefore whether the initial condition $\psi$ in our present context is convex or concave, or whether it is neither of those. Unfortunately, the answer is negative in general. In Exercise 6.7, a simple example with only one type of spin taking values in $\Sigma_1$ is given for which the initial condition is neither convex nor concave, even if we restrict it to the space of constant paths[2]. We also recall that, in order for a Hopf variational representation to be valid, the convexity or concavity of the initial condition is not only sufficient, but also necessary, as is clear from the formula (3.68) with $t = 0$.

In short, our usual candidate variational representations for the solution to (6.153), namely the Hopf-Lax and the Hopf formulas, are invalid in general. This motivates to try to make sense of the solution to the Hamilton-Jacobi equation (6.153) directly. Perhaps the easiest way to proceed is to construct solutions to finite-dimensional approximations to (6.153), and then pass to the limit. That is, in the context of the SK model and inspired by (6.112), we would take a large integer $K$, define $U_K$ as in (6.109), and seek a function $f^{(K)} : \mathbb{R}_{\geq 0} \times U_K \to \mathbb{R}$ that solves the finite-dimensional Hamilton-Jacobi equation

$$\partial_t f^{(K)} - (K+1) \sum_{k=0}^{K} \left(\partial_{q_k} f^{(K)}\right)^2 = 0 \quad \text{on} \quad \mathbb{R}_{>0} \times U_K \tag{6.155}$$

subject to the initial condition $f^{(K)}(0, q) = \psi(\mathsf{q}_q)$, with $\mathsf{q}_q$ as in (6.50). For the bipartite model, we would set

$$U_K^2 := \{q = (q_0, \ldots, q_K) \in (\mathbb{R}_{\geq 0}^2)^{K+1} \mid 0 \leq q_0 \leq q_1 \leq \cdots \leq q_K\}, \tag{6.156}$$

where we understand that $q_k \leq q_{k+1}$ stands for $q_{k+1} - q_k \in \mathbb{R}_{\geq 0}^2$, and we would aim to find a function $f^{(K)} : \mathbb{R}_{\geq 0} \times U_K^2 \to \mathbb{R}$ that solves the finite-dimensional Hamilton-Jacobi equation

$$\partial_t f^{(K)} - (K+1) \sum_{k=0}^{K} \partial_{q_{k,1}} f^{(K)} \partial_{q_{k,2}} f^{(K)} = 0 \quad \text{on} \quad \mathbb{R}_{\geq 0} \times U_K^2, \tag{6.157}$$

subject to the appropriate initial condition. Continuing with the bipartite model, for each $\mathsf{q} \in \mathcal{Q}_2(\mathbb{R}_{\geq 0}^2)$, we would then want to use piecewise-constant approximations to $\mathsf{q}$ such as

$$\mathsf{q} \simeq \sum_{k=0}^{K} \mathbf{1}_{[\frac{k}{K+1}, \frac{k+1}{K+1})} \frac{1}{K+1} \int_{\frac{k}{K+1}}^{\frac{k+1}{K+1}} \mathsf{q}(u) \, du. \tag{6.158}$$

Together with (6.110), this motivates the introduction of the sequence $q^{(K)}(\mathsf{q}) = (q_k^{(K)}(\mathsf{q}))_{0 \leq k \leq K}$ defined, for each $k \in \{0, \ldots, K\}$, by

$$q_k^{(K)}(\mathsf{q}) := \frac{1}{K+1} \int_{\frac{k}{K+1}}^{\frac{k+1}{K+1}} \mathsf{q}(u) \, du. \tag{6.159}$$

---

[2] For the expert reader who may think that this contradicts the main result of [25], we point out that the convexity shown there is with respect to the *inverse function* to our paths q. It is the convexity or concavity with respect to q that is relevant in the study of the variational representations to (6.153).

We also add that Exercise 6.7 shows the invalidity of a Hopf variational representation to (6.153) for some choices of the reference measure, but we do not exclude the possibility that it is valid for other choices.

It is shown in [195] that the functions $f^{(K)}$ are well-defined, and that for each $q \in Q_2(\mathbb{R}_{\geq 0}^2)$, the limit

$$f(t,q) := \lim_{K \to +\infty} f^{(K)}(t, q^{(K)}(q)) \qquad (6.160)$$

exists and is finite. One may then decide that this limit object $f$ is what we understand to be the solution to (6.153) with initial condition $\psi$. A more intrinsic notion, directly making sense of what it means to be a viscosity solution to (6.153) in infinite dimensions, was introduced in [73] and shown to be equivalent to this definition by successive approximations. In the context of the SK model, one can write Hopf-Lax formulas for each of the finite-dimensional approximations and then pass to the limit in the resulting expressions, thereby verifying that the notion of solution just introduced here coincides with the variational definition in (6.135) (see also [73, 74]). From now on, whenever we refer to the solution to (6.153), we refer to the function defined through the procedure ending with (6.160).

We are thus led to the following question.

**Question 6.9.** Let $\overline{F}_N(t,q)$ be the enriched free energy for the bipartite model, as defined in (6.148), and let $f(t,q)$ be the solution to (6.153) with initial condition $\psi$. Is it true that for every $t \geq 0$ and $q \in Q_2(\mathbb{R}_{\geq 0}^2)$, we have

$$\lim_{N \to +\infty} \overline{F}_N(t,q) = f(t,q) \quad ? \qquad (6.161)$$

We do not know the answer to this question. The following result gives an inequality between the two terms in (6.161).

**Theorem 6.10** ([195]). *With the same notation as in Question 6.9, we have for every $t \geq 0$ and $q \in Q_2(\mathbb{R}_{\geq 0}^2)$ that*

$$\liminf_{N \to +\infty} \overline{F}_N(t,q) \geq f(t,q). \qquad (6.162)$$

One can also show that any subsequential limit of $\overline{F}_N$ must satisfy the equation (6.153) "almost everywhere". This is shown in a sense involving finite-dimensional approximations in [195], and also directly in infinite dimensions in [70] for a precise infinite-dimensional generalization of the notion of "almost everywhere". As was explained in Example 3.1, the verification of such a property does not suffice to identify the solution uniquely, even in finite dimensions.

We now describe how Question 6.9 would be phrased for more general mean-field spin glasses with a finite number of types. For a fixed integer $D \geq 1$, we take our configuration vector to be of the form $\sigma = (\sigma_1, \ldots, \sigma_D) \in (\mathbb{R}^N)^D$. We take $P_N$ to be a probability measure on $(\mathbb{R}^N)^D \simeq \mathbb{R}^{D \times N}$ such that if $\sigma$ is sampled according to $P_N$, and if we think of $\sigma$ as a $D$-by-$N$ matrix, then the $N$ columns of $\sigma$ are independent and identically distributed with a fixed distribution $P_1$ of bounded support in $\mathbb{R}^D$. We give ourselves a centred Gaussian field $(H_N(\sigma))_{\sigma \in \mathbb{R}^{D \times N}}$ such that, for some smooth function $\xi : \mathbb{R}^{D \times D} \to \mathbb{R}$ and for every $\sigma^1, \sigma^2 \in (\mathbb{R}^N)^D$,

$$\mathbb{E} H_N(\sigma^1) H_N(\sigma^2) = N\xi\left(\left(\frac{\sigma_d^1 \cdot \sigma_{d'}^2}{N}\right)_{d,d' \leqslant D}\right). \tag{6.163}$$

The mean-field character of the model is encoded by the fact that this covariance depends only on the $D$-by-$D$ matrix of scalar products between the different "types" of spins. The set of functions $\xi$ for which such a Gaussian random field exists is described precisely in Proposition 6.6 of [197]. Writing $S^D_{\geqslant 0}$ to denote the set of $D$-by-$D$ positive semi-definite matrices, we give ourselves a family of real parameters $(\zeta_k)_{1 \leqslant k \leqslant K}$ as in (6.143) and

$$0 = q_{-1} \leqslant q_0 < q_1 < \cdots < q_K \in S^D_{\geqslant 0}, \tag{6.164}$$

where we understand that $q_k \leqslant q_{k+1}$ means that $q_{k+1} - q_k \in S^D_{\geqslant 0}$, and $q_k < q_{k+1}$ means that $q_k \leqslant q_{k+1}$ and $q_k \neq q_{k+1}$. We encode this into a path q as in (6.145), an element of the space of paths $Q(S^D_{\geqslant 0})$ defined analogously to (6.146). We take $(v_\alpha)_{\alpha \in \mathbb{N}^K}$ to be the Poisson-Dirichlet cascade associated with the parameters $(\zeta_k)_{1 \leqslant k \leqslant K}$, and let $(z_\alpha)_{\alpha \in \mathcal{A}}$ be independent $D$-by-$N$ random matrices whose entries are independent standard Gaussians. We then define

$$Z_{\mathsf{q}}(\alpha) := \sum_{k=0}^{K} (2q_k - 2q_{k-1})^{1/2} z_{\alpha_{|k}}, \tag{6.165}$$

which is a random $D$-by-$N$ matrix. Recall that for any two matrices $a$ and $b$ of the same size, we denote by $a \cdot b$ their entry-wise scalar product, that is, $a \cdot b := \operatorname{tr}(ab^*)$. We define the enriched free energy as

$$\overline{F}_N(t,\mathsf{q}) := -\frac{1}{N} \mathbb{E} \log \int_{\Sigma_N^2} \sum_{\alpha \in \mathbb{N}^K} \exp\left(H_N(t,\mathsf{q},\sigma,\alpha)\right) v_\alpha \, dP_N(\sigma) \tag{6.166}$$

for the Hamiltonian

$$H_N(t,\mathsf{q},\sigma,\alpha)$$
$$:= \sqrt{2t} H_N(\sigma) - Nt\xi\left(\left(\frac{\sigma_d \cdot \sigma_{d'}}{N}\right)_{d,d' \leqslant D}\right) + Z_{\mathsf{q}}(\alpha) \cdot \sigma - \sigma \cdot q_K \sigma. \tag{6.167}$$

Similar arguments as for the bipartite model lead us to expect that, as $N$ tends to infinity, the free energy $\overline{F}_N$ in (6.148) converges to the function $f : \mathbb{R}_{\geqslant 0} \times Q_2(S^D_{\geqslant 0}) \to \mathbb{R}$ that solves the infinite-dimensional Hamilton-Jacobi equation

$$\partial_t f(t,\mathsf{q}) - \int_0^1 \xi\left(\partial_{\mathsf{q}} f(t,\mathsf{q},u)\right) du = 0 \quad \text{on} \quad \mathbb{R}_{>0} \times Q_2(S^D_{\geqslant 0}) \tag{6.168}$$

subject to the initial condition $f(0,\mathsf{q}) = \psi(\mathsf{q}) := \overline{F}_1(0,\mathsf{q})$. One can construct a solution to this equation exactly as in the procedure leading to (6.160).

**Question 6.11.** Let $\overline{F}_N(t,\mathsf{q})$ be the enriched free energy (6.166), and let $f(t,\mathsf{q})$ be the solution to (6.168) with initial condition $\psi$. Is it true that for every $t \geq 0$ and $\mathsf{q} \in \mathcal{Q}_2(S^D_{\geq 0})$, we have

$$\lim_{N \to +\infty} \overline{F}_N(t,\mathsf{q}) = f(t,\mathsf{q}) \quad ? \tag{6.169}$$

The results presented in Theorem 6.10 and in the paragraph below have been extended to this setting [70, 197]. We also mention that, for certain problems of community detection which, unlike the situation explored in Section 4.5, involve random graphs whose average degree remains bounded, very similar difficulties emerge [105, 153]. When the covariance function $\xi$ is convex over $S^D_{\geq 0}$, the function $f$ in Question 6.11 admits a Hopf-Lax variational representation [73], and the identity (6.169) is indeed valid [70].

In Theorem 3.19, we observed that whenever the non-linearity or the initial condition of the finite-dimensional Hamilton-Jacobi equation (3.20) is convex or concave, the solution to this equation must be a critical value of the function $\mathcal{J}_{t,x}$ defined in (3.53). Remarkably, the limit of the free energy $\overline{F}_N(t,\mathsf{q})$ defined in (6.166), assuming it exists, must also satisfy this property for the natural analogue of $\mathcal{J}_{t,x}$ associated with the infinite-dimensional Hamilton-Jacobi equation (6.168). What is truly remarkable is that this property is valid even in situations in which neither the non-linearity nor the initial condition possess any global convexity or concavity property. To state this precisely, we define, for every $t \geq 0$ and $\mathsf{p}, \mathsf{q}, \mathsf{q}' \in \mathcal{Q}_2(S^D_{\geq 0})$,

$$\mathcal{J}_{t,\mathsf{q}}(\mathsf{q}',\mathsf{p}) := \psi(\mathsf{q}') + \int_0^1 \mathsf{p}(u) \cdot (\mathsf{q} - \mathsf{q}')(u)\,du + t \int_0^1 \xi(\mathsf{p}(u))\,du. \tag{6.170}$$

**Theorem 6.12** ([70]). *Suppose that the free energy $\overline{F}_N$ defined in (6.166) converges pointwise to a function $f : \mathbb{R}_{\geq 0} \times \mathcal{Q}_2(S^D_{\geq 0}) \to \mathbb{R}$. For every $t \geq 0$ and $\mathsf{q} \in \mathcal{Q}_2(S^D_{\geq 0})$, there exists $\mathsf{q}', \mathsf{p} \in \mathcal{Q}_2(S^D_{\geq 0})$ such that*

$$\mathsf{q} = \mathsf{q}' - t\nabla\xi(\mathsf{p}), \qquad \mathsf{p} = \partial_q\psi(\mathsf{q}'), \tag{6.171}$$

*and*

$$f(t,\mathsf{q}) = \mathcal{J}_{t,\mathsf{q}}(\mathsf{q}',\mathsf{p}). \tag{6.172}$$

The conditions in (6.171) correspond to the requirement that the pair $(\mathsf{q}',\mathsf{p})$ be a critical point of the functional $\mathcal{J}_{t,\mathsf{q}}$. They can be written more explicitly as

$$\mathsf{q}(u) = \mathsf{q}'(u) - t\nabla\xi(\mathsf{p}(u)), \qquad \mathsf{p}(u) = \partial_q\psi(\mathsf{q}',u) \qquad du\text{-a.e. in }[0,1). \tag{6.173}$$

As was observed in Section 3.5, these conditions can be interpreted as saying that the characteristic line starting at $\mathsf{q}'$ passes through the point $(t,\mathsf{q})$. Roughly speaking, it is also shown in [70] that, up to a small perturbation of the parameters, one can choose the critical point $(q',p)$ in Theorem 6.12 in such a way that the overlap matrix $N^{-1}(\sigma_d \cdot \sigma'_{d'})_{d,d' \leqslant D}$ converges in law to $p(U)$, where $U$ is a uniform random variable on $[0,1]$.

**Exercise 6.7.** Let $P_1$ be a probability measure on $\mathbb{R}$ with compact support, let $z$ be a standard Gaussian random variable, and consider the free energy

$$f(t) := -\mathbb{E}\log \int_{\mathbb{R}} \exp\left(\sqrt{2t}z\sigma - t\sigma^2\right) dP_1(\sigma)$$

associated with the Hamiltonian $H(t,\sigma) := \sqrt{2t}z\sigma - t\sigma^2$. Denoting by $\langle \cdot \rangle$ its associated Gibbs measure, show that

$$\partial_t^2 f = 2\mathbb{E}\big(\langle \sigma^2 \rangle - \langle \sigma \rangle^2\big)\big(\langle \sigma^2 \rangle - 3\langle \sigma \rangle^2\big).$$

Deduce that for a suitable choice of measure $P_1$, the free energy $f$ is neither convex nor concave.

# Appendix A
# Basic results in analysis and probability

In this appendix, we give detailed proofs of the basic results in analysis and probability theory that are used in the main text. The topics covered are classical, so readers already familiar with the material can safely skip it; we have merely included it to keep the book as self-contained as possible. The key results that we prove are the Caratheodory extension theorem, the Dynkin $\pi$-$\lambda$ theorem, the Riesz representation theorem, the Stone-Weierstrass theorem, the Lebesgue differentiation theorem, the Portmanteau theorem, the Prokhorov theorem, and the injectivity of the Fourier and Laplace transforms. Although some of these results can be proved in greater generality, we will focus on the setting in which we give ourselves a metric space $S$, whose metric we denote by $d$.

## A.1 Constructing measures

One of the most classical problems in geometry is to determine the length, area, or volume of a curve, surface or solid. This problem can be generalized to the setting of a metric space $S$ as the task of assigning an adequate notion of size or measure to a subset of $S$. Intuitively, any reasonable notion of measure should be countably additive, monotone, and assign no mass to the empty set. Unfortunately, examples such as that due to Vitali [260] show that it is not in general possible to assign such a notion of measure to all sets in $S$. To overcome this issue, we will only define the measure of a collection $\mathcal{S}$ of subsets of $S$ which make the pair $(S, \mathcal{S})$ a measurable space.

A pair $(S, \mathcal{S})$ is a *measurable space* if $\mathcal{S}$ is a $\sigma$-algebra of subsets of $S$. A collection $\mathcal{A}$ of subsets of $S$ is called an *algebra* if it contains $S$, it is closed under finite unions and it is closed under complements,

(i) $S \in \mathcal{A}$,

(ii) for every $A, B \in \mathcal{A}$, we have $A \cup B \in \mathcal{A}$,

(iii) if $A \in \mathcal{A}$, then $A^c \in \mathcal{A}$.

A collection $\mathcal{S}$ of subsets of $S$ is called a *σ-algebra* if it is an algebra which is closed under countable unions,

(iv) if $(A_n)_{n \geqslant 1}$ is a sequence of sets in $\mathcal{S}$, then $\bigcup_{n \geqslant 1} A_n \in \mathcal{S}$.

To verify that an algebra $\mathcal{A}$ is a σ-algebra, it suffices to verify that it is countably additive on disjoint sets,

(iv') if $(A_n)_{n \geqslant 1}$ is a sequence of disjoint sets in $\mathcal{S}$, then $\bigcup_{n \geqslant 1} A_n \in \mathcal{S}$.

Indeed, if $(A_n)_{n \geqslant 1}$ is a sequence of sets in $\mathcal{S}$, then the sequence

$$F_n := A_n \setminus \bigcup_{i=1}^{n-1} A_i \tag{A.1}$$

of disjoint sets is also in $\mathcal{S}$ and has the same union as $(A_n)_{n \geqslant 1}$. The simplest example of a σ-algebra is the *σ-algebra generated* by a collection $\mathcal{E}$ of subsets of $S$,

$$\sigma(\mathcal{E}) := \bigcap \{\mathcal{M} \mid \mathcal{M} \text{ is a σ-algebra containing } \mathcal{E}\}. \tag{A.2}$$

This is the smallest σ-algebra containing $\mathcal{E}$ and it is well-defined as the intersection of a collection of σ-algebras is again a σ-algebra. We will typically endow a metric space $S$ with its Borel σ-algebra. The *Borel σ-algebra* on a metric space $S$ is the σ-algebra generated by all open sets,

$$\mathcal{B}(S) := \sigma(\{U \mid U \text{ is open in } S\}). \tag{A.3}$$

A *measure* $\mu$ on a measurable space $(S, \mathcal{S})$ is a non-negative set function $\mu : \mathcal{S} \to [0, +\infty]$ which assigns no mass to the empty set and is *countably additive*,

(i) $\mu(\emptyset) = 0$,

(ii) for every sequence $(A_n)_{n \geqslant 1}$ of disjoint sets in $\mathcal{S}$,

$$\mu\left(\bigcup_{n=1}^{\infty} A_n\right) = \sum_{n=1}^{+\infty} \mu(A_n). \tag{A.4}$$

A measure $\mu$ is said to be a *finite measure* if $\mu(S) < +\infty$, and a finite measure $\mathbb{P}$ is said to be a *probability measure* if $\mathbb{P}(S) = 1$. The simplest example of a measure is the Dirac measure at a point $x \in S$ defined on the measurable space $(S, 2^S)$ by

$$\delta_x(A) := \mathbf{1}_A(x). \tag{A.5}$$

Thinking back to the geometric problem of determining length, we would like our first non-trivial example of a measure to be the measure $m$ on the real line which gives the length of a set. Although there are many sets whose length we do not know a priori, we know the value of the measure $m$ on the algebra

$$\mathcal{A} := \left\{ \bigcup_{i=1}^{n}(a_i, b_i] \mid n \geq 1 \text{ and } (a_i, b_i] \text{ are disjoint intervals} \right\} \tag{A.6}$$

of finite unions of disjoint half-open intervals,

$$m\left(\bigcup_{i=1}^{n}(a_i, b_i]\right) := \sum_{i=1}^{n}(b_i - a_i). \tag{A.7}$$

We will first construct a measure on the Borel $\sigma$-algebra $\mathcal{B}(\mathbb{R})$ which satisfies (A.7) using the Caratheodory extension theorem, and then we will prove that such a measure must be unique using the Dynkin $\pi$-$\lambda$ theorem. This strategy of constructing a measure by first specifying its value on a large enough but simple collection of sets such as an algebra, then establishing its existence on the $\sigma$-algebra generated by this collection of sets using the Caratheodory extension theorem, and finally obtaining its uniqueness using the Dynkin $\pi$-$\lambda$ theorem is very general.

The idea behind the Caratheodory extension theorem is to first define an outer measure on the collection $2^S$ of all subsets of $S$, and then restrict this outer measure to the $\sigma$-algebra of Caratheodory measurable sets where it becomes a measure. An *outer measure* $\mu^*$ on $S$ is a non-negative set function $\mu^* : 2^S \to [0, +\infty]$ which assigns no mass to the empty set, is *monotone* and *countably sub-additive*,

(i) $\mu^*(\emptyset) = 0$,

(ii) for all subsets $A, B$ of $S$ with $A \subseteq B$,

$$\mu^*(A) \leq \mu^*(B), \tag{A.8}$$

(iii) for every sequence $(A_n)_{n \geq 1}$ of subsets of $S$,

$$\mu^*\left(\bigcup_{n=1}^{\infty} A_n\right) \leq \sum_{n=1}^{+\infty} \mu^*(A_n). \tag{A.9}$$

The most common way to define an outer measure is to start with a collection $\mathcal{E} \subseteq 2^S$ of elementary sets and a non-negative set function $\rho : \mathcal{E} \to [0, +\infty]$, and define the set function $\mu^* : 2^S \to [0, +\infty]$ by

$$\mu^*(A) := \inf\left\{\sum_{n=1}^{+\infty} \rho(E_n) \mid E_n \in \mathcal{E} \text{ and } A \subseteq \bigcup_{n=1}^{\infty} E_n\right\}. \tag{A.10}$$

Provided that $\emptyset, S \in \mathcal{E}$ and $\rho(\emptyset) = 0$, this set function $\mu^*$ is an outer measure. The assumption $S \in \mathcal{E}$ ensures that the infimum in (A.10) is never taken over the empty set, so $\mu^*$ is well-defined.

**Proposition A.1.** *If $\mathcal{E} \subseteq 2^S$ is a collection of sets with $\varnothing, S \in \mathcal{E}$ and $\rho : \mathcal{E} \to [0, +\infty]$ is a non-negative set function with $\rho(\varnothing) = 0$, then (A.10) is an outer measure.*

*Proof.* Taking $E_n = \varnothing \in \mathcal{E}$ for all $n \geq 1$ shows that $\mu^*(\varnothing) = \rho(\varnothing) = 0$. If $A, B$ are subsets of $S$ with $A \subseteq B$, then the infimum defining $\mu^*(A)$ is taken over a larger set than that defining $\mu^*(B)$ which means that $\mu^*(A) \leq \mu^*(B)$. To establish countable sub-additivity, we fix a sequence $(A_n)_{n \geq 1}$ of subsets of $S$ as well as $\varepsilon > 0$. For each $n \geq 1$, let $(E_{n,i})_{i \geq 1} \subseteq \mathcal{E}$ be such that for every $n \geq 1$,

$$A_n \subseteq \bigcup_{j=1}^{\infty} E_{n,i} \quad \text{and} \quad \sum_{i=1}^{+\infty} \rho(E_{n,i}) \leq \mu^*(A_n) + \frac{\varepsilon}{2^n}.$$

Since $(E_{n,i})_{n,i \geq 1}$ forms a cover of $\bigcup_{n=1}^{\infty} A_n$, the definition of $\mu^*$ implies that

$$\mu^*\left(\bigcup_{n=1}^{\infty} A_n\right) \leq \sum_{n,i=1}^{+\infty} \rho(E_{n,i}) \leq \sum_{n=1}^{+\infty} \mu^*(A_n) + \varepsilon.$$

Letting $\varepsilon$ tend to zero completes the proof. ∎

The key insight behind the Caratheodory extension theorem is that an outer measure restricted to the $\sigma$-algebra of Caratheodory measurable sets becomes a measure. A subset $A$ of $S$ is *Caratheodory measurable* with respect to an outer measure $\mu^*$ if, for every other subset $E$ of $S$,

$$\mu^*(E) = \mu^*(E \cap A) + \mu^*(E \cap A^c). \tag{A.11}$$

In other words, $A$ is Caratheodory measurable if it can be used to partition any other set $E$ in an additive fashion. We will write

$$\mathcal{S}_* := \{A \subseteq S \mid A \text{ is Caratheodory measurable}\} \tag{A.12}$$

for the collection of Caratheodory measurable sets.

**Theorem A.2** (Caratheodory extension). *If $\mu^*$ is an outer measure on $S$, then the collection $\mathcal{S}_*$ of Caratheodory measurable sets is a $\sigma$-algebra, and the restriction of $\mu^*$ to $\mathcal{S}_*$ is a measure.*

*Proof.* Since $\mu^*(\varnothing) = 0$ and the definition of Caratheodory measurable is symmetric in $A$ and $A^c$, the set $\mathcal{S}_*$ contains $S$ and is closed under complements. To show that it is an algebra, fix two Caratheodory measurable sets $A, B \in \mathcal{S}_*$ as well as an arbitrary subset $E$ of $S$. By Caratheodory measurability of $A$ and $B$,

$$\mu^*(E) = \mu^*(E \cap A \cap B) + \mu^*(E \cap A \cap B^c) + \mu^*(E \cap A^c \cap B) + \mu^*(E \cap A^c \cap B^c).$$

Combining the fact that $A \cup B = (A \cap B) \cup (A \cap B^c) \cup (A^c \cap B)$ with the sub-additivity of $\mu^*$ reveals that

$$\mu^*(E \cap (A \cup B)) \leq \mu^*(E \cap A \cap B) + \mu^*(E \cap A \cap B^c) + \mu^*(E \cap A^c \cap B).$$

It follows that

$$\mu^*(E \cap (A \cup B)) + \mu^*(E \cap (A \cup B)^c) \leq \mu^*(E),$$

where we have used the fact that $(A \cup B)^c = A^c \cap B^c$. By sub-additivity of $\mu^*$ this inequality is in fact an equality so $A \cup B \in \mathcal{S}_*$ and $\mathcal{S}_*$ is an algebra. To prove that $\mathcal{S}_*$ is a $\sigma$-algebra, it suffices to fix a sequence $(A_n)_{n \geq 1}$ of disjoint sets in $\mathcal{S}_*$ and show that their union also lies in $\mathcal{S}_*$. Fix an integer $n \geq 1$, and to simplify notation, let $B_n := \bigcup_{i=1}^n A_i \in \mathcal{S}_*$ and $B := \bigcup_{n=1}^\infty A_n$. For any subset $E$ of $S$, the Caratheodory measurability of $A_n$ and the assumption that the sets $A_i$ are disjoint imply that

$$\mu^*(E \cap B_n) = \mu^*(E \cap B_n \cap A_n) + \mu^*(E \cap B_n \cap A_n^c)$$
$$= \mu^*(E \cap A_n) + \mu^*(E \cap B_{n-1}).$$

It follows by Caratheodory measurability of $B_n$ and an induction that

$$\mu^*(E) = \mu^*(E \cap B_n) + \mu^*(E \cap B_n^c) = \sum_{i=1}^n \mu^*(E \cap A_i) + \mu^*(E \cap B_n^c).$$

Noticing that $B^c \subseteq B_n^c$, leveraging the monotonicity of $\mu^*$, and letting $n$ tend to infinity shows that

$$\mu^*(E) \geq \sum_{n=1}^{+\infty} \mu^*(E \cap A_n) + \mu^*(E \cap B^c). \tag{A.13}$$

Invoking the sub-additivity of $\mu^*$ reveals that

$$\mu^*(E) \geq \mu^*\left(\bigcup_{n=1}^\infty (E \cap A_n)\right) + \mu^*(E \cap B^c) = \mu^*(E \cap B) + \mu^*(E \cap B^c) \geq \mu^*(E),$$

so in fact these inequalities are equalities and $B \in \mathcal{S}_*$ as required for $\mathcal{S}_*$ to be a $\sigma$-algebra. That $\mu^*$ restricts to a measure on $\mathcal{S}_*$ is immediate from (A.13) with $E = B$ and the sub-additivity of outer measure. This completes the proof. ∎

The Caratheodory extension theorem is often stated as the possibility to extend a pre-measure from an algebra to the $\sigma$-algebra it generates. A *pre-measure* $\mu_0$ on an algebra $\mathcal{A}$ is a non-negative set function $\mu_0 : \mathcal{A} \to [0, +\infty]$ which assigns no mass to the empty set and is countably additive whenever this makes sense,

(i) $\mu(\emptyset) = 0$,

(ii) if $(A_n)_{n \geq 1}$ is a sequence of disjoint sets in $\mathcal{A}$ with $\bigcup_{n=1}^{\infty} A_n \in \mathcal{A}$, then

$$\mu_0\left(\bigcup_{n=1}^{\infty} A_n\right) = \sum_{n=1}^{+\infty} \mu_0(A_n). \tag{A.14}$$

We have opted for the more general version of the Caratheodory extension theorem as we will prove the Riesz representation theorem by extending a set function from the collection of open sets, which does not form an algebra, to the Borel $\sigma$-algebra $\mathcal{B}(S)$. Nonetheless, let us establish this simplified version of the Caratheodory extension theorem and use it to extend the pre-measure $m$ defined by (A.7) to a measure on the Borel $\sigma$-algebra $\mathcal{B}(\mathbb{R})$. The fact that (A.7) is a well-defined pre-measure on the algebra (A.6) is the content of Exercise A.2.

**Corollary A.3** (Simplified Caratheodory extension). *If $\mu_0 : \mathcal{A} \to \mathbb{R}$ is a pre-measure on an algebra $\mathcal{A}$, then $\mu_0$ can be extended to a measure $\mu$ on the $\sigma$-algebra $\sigma(\mathcal{A})$.*

*Proof.* Denote by $\mu^* : 2^S \to [0, +\infty]$ the set function defined by (A.10) with $\rho = \mu_0$, and recall that it is an outer measure by Proposition A.1. It follows by the Caratheodory extension theorem that the restriction $\mu$ of $\mu^*$ to the $\sigma$-algebra $\mathcal{S}_*$ of Caratheodory measurable sets is a measure. To see the $\mu$ is an extension of $\mu_0$, fix $A \in \mathcal{A}$, and let $(E_n)_{n \geq 1} \subseteq \mathcal{A}$ be a covering of $A$. Introduce the disjoint sets

$$F_n := A \cap \left(E_n \setminus \bigcup_{i=1}^{n-1} E_i\right) \in \mathcal{A}$$

whose union is $A$. The definition of a pre-measure implies that

$$\mu_0(A) = \sum_{n=1}^{+\infty} \mu_0(F_n) \leq \sum_{n=1}^{+\infty} \mu_0(E_n),$$

and taking the infimum over all coverings $(E_n)_{n \geq 1}$ shows that $\mu_0(A) \leq \mu^*(A)$. Taking the covering $E_1 = A$ and $E_n = \emptyset$ for $n \geq 2$ shows that this is in fact an equality, so $\mu_0 = \mu^* = \mu$ on $\mathcal{A}$. To show that the extension $\mu$ is defined on $\sigma(\mathcal{A})$, it suffices to verify that $\mathcal{A} \subseteq \mathcal{S}_*$. Fix a set $A \in \mathcal{A}$ as well as an arbitrary subset $E$ of $S$ and $\varepsilon > 0$. Let $(E_n)_{n \geq 1} \subseteq \mathcal{A}$ be a sequence of elementary sets with

$$E \subseteq \bigcup_{n=1}^{\infty} E_n \quad \text{and} \quad \mu^*(E) \leq \sum_{n=1}^{+\infty} \mu_0(E_n) + \varepsilon.$$

Since $\mu_0$ is additive on $\mathcal{A}$,

$$\mu^*(E) + \varepsilon \geq \sum_{n=1}^{+\infty} \mu_0(E_n \cap A) + \sum_{n=1}^{+\infty} \mu_0(E_n \cap A^c) \geq \mu^*(E \cap A) + \mu^*(E \cap A^c).$$

Since $\varepsilon$ is arbitrary, together with the sub-additivity of outer measure, this shows that $A \in \mathcal{S}_*$ and completes the proof. ∎

A direct application of this result gives *a* measure $m$ on the Borel $\sigma$-algebra $\mathcal{B}(\mathbb{R})$ which satisfies (A.7) on the algebra (A.6). We would like to say that this is *the* measure on the Borel $\sigma$-algebra $\mathcal{B}(\mathbb{R})$ which gives the length of a set. We therefore need to establish the uniqueness of such a measure. This will be done using the Dynkin $\pi$-$\lambda$ theorem. A *$\pi$-system* is a collection $\mathcal{P}$ of sets that is closed under intersections,

(i) if $A, B \in \mathcal{P}$, then $A \cap B \in \mathcal{P}$.

A *$\lambda$-system* is a collection $\mathcal{L}$ of sets that contains the set $S$, is closed under complements and is closed under increasing unions,

(i) $S \in \mathcal{L}$,

(ii) if $A, B \in \mathcal{L}$ are such that $A \subseteq B$, then $B \setminus A \in \mathcal{L}$,

(iii) if $(A_n)_{n \geq 1}$ is an increasing sequence of sets in $\mathcal{L}$ with $A_1 \subseteq A_2 \subseteq \cdots$, then $\bigcup_{n=1}^\infty A_n \in \mathcal{L}$.

**Lemma A.4.** *If $\mathcal{L}$ is both a $\pi$-system and a $\lambda$-system, then it is a $\sigma$-algebra.*

*Proof.* Since a $\lambda$-system is closed under complements, it suffices to show that $\mathcal{L}$ is closed under countable unions. Fix a sequence $(A_n)_{n \geq 1}$ of sets in $\mathcal{L}$, and consider the increasing sequence
$$B_n := \bigcup_{i=1}^n A_i.$$
An induction relying on the fact that $B_1 \cup B_2 = (B_1^c \cap B_2^c)^c$ and the assumption that $\mathcal{L}$ is a $\pi$-system shows that $(B_n)_{n \geq 1} \subseteq \mathcal{L}$. Observing that $(B_n)_{n \geq 1}$ increases to $\bigcup_{n=1}^\infty A_n$ and remembering that $\mathcal{L}$ is a $\lambda$-system completes the proof. ∎

**Theorem A.5** (Dynkin $\pi$-$\lambda$)**.** *If $\mathcal{P}$ is a $\pi$-system and $\mathcal{L}$ is a $\lambda$-system with $\mathcal{P} \subseteq \mathcal{L}$, then*
$$\sigma(\mathcal{P}) \subseteq \mathcal{L}. \tag{A.15}$$

*Proof.* Denote by $\ell(\mathcal{P})$ the smallest $\lambda$-system containing $\mathcal{P}$. By minimality of $\ell(\mathcal{P})$ and $\sigma(\mathcal{P})$, it suffices to prove that $\ell(\mathcal{P})$ is a $\sigma$-algebra. Indeed, if this is the case, then
$$\sigma(\mathcal{P}) \subseteq \ell(\mathcal{P}) \subseteq \mathcal{L}.$$
By Lemma A.4 it is in fact sufficient to show that $\ell(\mathcal{P})$ is a $\pi$-system. For each $A \in \ell(\mathcal{P})$ define the $\lambda$-system
$$\mathcal{G}_A := \{B \mid A \cap B \in \ell(\mathcal{P})\}.$$
Since $\mathcal{P}$ is a $\pi$-system, for every $A \in \mathcal{P}$, we have $\mathcal{P} \subseteq \mathcal{G}_A$. It follows by minimality of $\ell(\mathcal{P})$ that $\ell(\mathcal{P}) \subseteq \mathcal{G}_A$. This means that for every $A \in \mathcal{P}$ and $B \in \ell(\mathcal{P})$, we have $A \cap B \in \ell(\mathcal{P})$. In other words, for every $B \in \ell(\mathcal{P})$, we have $\mathcal{P} \subseteq \mathcal{G}_B$. It follows by minimality of $\ell(\mathcal{P})$ that for every $B \in \ell(\mathcal{P})$, we have $\ell(\mathcal{P}) \subseteq \mathcal{G}_B$. This shows that $\ell(\mathcal{P})$ is a $\pi$-system and completes the proof. ∎

As a direct application of the Dynkin $\pi$-$\lambda$ theorem, we can now show that there is a unique $\sigma$-finite measure on the Borel $\sigma$-algebra $\mathcal{B}(\mathbb{R})$ which corresponds to our intuition about the length of a set. A measure $\mu$ on a measurable space $(S, \mathcal{S})$ is $\sigma$-finite if there exists a sequence of sets $(A_n)_{n \geq 1} \subseteq \mathcal{S}$ which cover $S$ and have finite measure,

$$S = \bigcup_{n=1}^{\infty} A_n \quad \text{and} \quad \mu(A_n) < +\infty \text{ for all } n \geq 1. \tag{A.16}$$

**Corollary A.6.** *There is a unique measure m on the Borel $\sigma$-algebra $\mathcal{B}(\mathbb{R})$ with*

$$m\left(\bigcup_{i=1}^{n}(a_i, b_i]\right) := \sum_{i=1}^{n}(b_i - a_i) \tag{A.17}$$

*on the algebra* (A.6) *of finite unions of disjoint intervals. This measure is the Lebesgue measure.*

*Proof.* Notice that any measure satisfying (A.17) must be $\sigma$-finite. The existence of a measure $m$ satisfying (A.17) on the algebra (A.6) has already been established by combining the Caratheodory extension theorem in Corollary A.3 with Exercise A.2. To obtain uniqueness, fix two $\sigma$-finite measures $\mu$ and $\nu$ on $\mathcal{B}(\mathbb{R})$ which coincide on the algebra (A.6). Fix an increasing sequence $(I_n)_{n \geq 1} \subseteq \mathcal{A}$ of intervals $I_1 \subseteq I_2 \subseteq \cdots$ which cover $\mathbb{R}$ and have finite $\mu$ and $\nu$ measures. By Exercise A.3, for each integer $n \geq 1$, the collection of sets

$$\mathcal{L}_n := \{B \in \sigma(\mathcal{A}) \mid \mu(B \cap I_n) = \nu(B \cap I_n)\} \tag{A.18}$$

forms a $\lambda$-system. Since this $\lambda$-system contains the $\pi$-system (A.6) and this $\pi$-system generates the Borel $\sigma$-algebra, the Dynkin $\pi$-$\lambda$ theorem implies that $\mathcal{L}_n = \mathcal{B}(\mathbb{R})$. It follows that for every $n \geq 1$ and Borel set $B \in \mathcal{B}(\mathbb{R})$,

$$\mu(B \cap I_n) = \nu(B \cap I_n). \tag{A.19}$$

Using the continuity from below of measure established in Exercise A.1 shows that $\mu = \nu$ on $\mathcal{B}(\mathbb{R})$ and completes the proof. ∎

The Caratheodory extension theorem and the Dynkin $\pi$-$\lambda$ theorem can be combined to uniquely define many measures by first defining them on a simpler class of sets. For instance, product measures can be constructed by defining them on the algebra of finite disjoint unions of rectangles. In the next section, we will combine these two results to establish the Riesz representation theorem by extending a set function from the collection of open sets to the Borel $\sigma$-algebra. Another important consequence of the Dynkin $\pi$-$\lambda$ theorem is that every finite measure on the Borel $\sigma$-algebra of a metric space is closed regular. A measure $\mu$ is *closed regular* if for every measurable set $A \subseteq S$,

$$\mu(A) = \sup\{\mu(F) \mid F \subseteq A \text{ and } F \text{ is closed}\}. \tag{A.20}$$

For a finite measure, taking complements shows that being closed regular is equivalent to being open regular. A measure $\mu$ is *open regular* if for every measurable set $A \subseteq S$,

$$\mu(A) = \inf\{\mu(U) \mid A \subseteq U \text{ and } U \text{ is open}\}. \tag{A.21}$$

The fact that every finite measure on the Borel $\sigma$-algebra of a metric space is closed regular is the content of Exercise A.4.

**Exercise A.1.** Fix a measure $\mu$ on a measurable space $(S, \mathcal{S})$.

(i) Show that $\mu$ is monotone in the sense that for all measurable sets $A, B \in \mathcal{S}$ with $A \subseteq B$, we have
$$\mu(A) \leq \mu(B). \tag{A.22}$$

(ii) Show that $\mu$ is continuous from below in the sense that for any increasing sequence of measurable sets $(A_n)_{n \geq 1} \subseteq \mathcal{S}$ with $A_1 \subseteq A_2 \subseteq \cdots$, we have
$$\mu\left(\bigcup_{n=1}^{\infty} A_n\right) = \lim_{n \to +\infty} \mu(A_n). \tag{A.23}$$

**Exercise A.2.** Let $F : \mathbb{R} \to \mathbb{R}$ be an increasing and right-continuous function, and write $\mathcal{A}$ for the algebra (A.6) of finite unions of disjoint intervals on $\mathbb{R}$. Show that the non-negative set function $m_F : \mathcal{A} \to [0, +\infty]$ defined by

$$m_F\left(\bigcup_{i=1}^n (a_i, b_i]\right) := \sum_{i=1}^n \left(F(b_i) - F(a_i)\right) \tag{A.24}$$

is a well-defined pre-measure.

**Exercise A.3.** Let $\mathcal{P}$ be a $\pi$-system of sets in $S$ and let $\mu$ and $\nu$ be measures on the measurable space $(S, \sigma(\mathcal{P}))$ which coincide on $\mathcal{P}$. Show that for any $A \in \mathcal{P}$ with $\mu(A) < +\infty$, the collection of sets

$$\mathcal{L} := \{B \in \sigma(\mathcal{P}) \mid \mu(A \cap B) = \nu(A \cap B)\} \tag{A.25}$$

forms a $\lambda$-system. Deduce that a finite measure on $(S, \sigma(\mathcal{P}))$ is determined by its values on $\mathcal{P}$.

**Exercise A.4.** Fix a finite measure $\mu$ on the Borel $\sigma$-algebra a metric space $S$. The purpose of this exercise is to show that $\mu$ is closed regular in the sense that (A.20) holds for any measurable set $A \subseteq S$.

(i) Show that the set

$$\mathcal{L} := \left\{ A \in \mathcal{S} \mid \mu(A) = \sup_{\substack{F \subseteq A \\ F \text{ closed}}} \mu(F) \text{ and } \mu(A^c) = \sup_{\substack{F \subseteq A^c \\ F \text{ closed}}} \mu(F) \right\} \quad (A.26)$$

is a $\lambda$-system.

(ii) Prove that $\mathcal{L}$ contains all open sets.

(iii) Conclude that $\mu$ is closed regular.

## A.2   The Riesz representation theorem

The Riesz representation theorem allows us to identify a positive linear functional on the space of continuous functions on a compact metric space $S$ with a unique measure on the Borel $\sigma$-algebra $\mathcal{B}(S)$. A *positive linear functional* on the space of continuous functions on $S$ is a linear functional $T : C(S;\mathbb{R}) \to \mathbb{R}$ that is non-negative on non-negative functions,

(i) for all $a, b \in \mathbb{R}$ and $f, g \in C(S;\mathbb{R})$, we have $T(af + bg) = aT(f) + bT(g)$,

(ii) for all $f \in C(S;\mathbb{R})$ with $f \geq 0$, we have $T(f) \geq 0$.

Notice that we do not require the continuity of the linear functional $T$; however, as shown in Exercise A.5, this continuity is implied by the positivity of $T$. To prove the Riesz representation theorem, we will rely on two classical topological results whose proofs are particularly simple in the metric setting.

**Lemma A.7** (Urysohn)**.** *If $F$ is a closed subset of a metric space $S$ that is contained in some open set $U \supseteq F$, then there exists a continuous function $f : S \to [0,1]$ which is equal to one on $F$ and vanishes outside $U$,*

$$\mathbf{1}_F \leq f \leq \mathbf{1}_U. \quad (A.27)$$

*Proof.* Denote by $d$ the metric on $S$, and define the function $f : S \to [0,1]$ by

$$f(x) := \frac{d(x, U^c)}{d(x, F) + d(x, U^c)}.$$

Since the distance to a closed set is continuous and $F \cap U^c = \emptyset$, the function $f$ is continuous. Moreover, we have $f(x) = 0$ if and only if $x \in U^c$ and $f(x) = 1$ if and only if $x \in F$. This completes the proof.  ∎

**Lemma A.8** (Partition of unity). *If $K$ is a compact subset of a metric space $S$ and $(U_i)_{i \leq n}$ is an open cover of $K$, then there exist continuous functions $(h_i)_{i \leq n}$ with*

$$\sum_{i=1}^{n} h_i = 1 \text{ on } K \quad \text{and} \quad 0 \leq h_i \leq \mathbf{1}_{U_i} \text{ for } 1 \leq i \leq n. \tag{A.28}$$

*The collection of functions $(h_i)_{i \leq n}$ is called a* partition of unity *on $K$ subordinate to the open cover $(U_i)_{i \leq n}$.*

*Proof.* Fix $x \in K$ and let $1 \leq i \leq n$ be such that $x \in U_i$. Denote by $B_{r(x)}(x)$ an open ball centred at $x$ of small enough radius $r(x) > 0$ so that $x \in \overline{B}_{r(x)}(x) \subseteq U_i$. The collection of sets $(B_{r(x)}(x))_{x \in K}$ forms an open cover of the compact set $K$ from which it is possible to extract a finite sub-cover $(B_{r(x_j)}(x_j))_{j \leq m}$. Let $V_i$ be the union of the closed balls $\overline{B}_{r(x_j)}(x_j)$ which lie in the open set $U_i$, and invoke Urysohn's lemma to find a collection of continuous functions $(g_i)_{i \leq n}$ with $\mathbf{1}_{V_i} \leq g_i \leq \mathbf{1}_{U_i}$ for $1 \leq i \leq n$. Since $(V_i)_{i \leq n}$ covers $K$, we have $\sum_{i=1}^{n} g_i \geq 1$ on $K$. To modify the functions $(g_i)_{i \leq n}$ so their sum becomes one on $K$, apply Urysohn's lemma once again to find a continuous functions $g$ with

$$\mathbf{1}_K \leq g \leq \mathbf{1}_{\left\{\sum_{i=1}^{n} g_i > 0\right\}}.$$

Define $g_{n+1} := 1 - g$ so that $\sum_{i=1}^{n+1} g_i > 0$ on $S$, and define the family of functions $(h_i)_{i \leq n}$ by

$$h_i := \frac{g_i}{\sum_{i=1}^{n+1} g_i}.$$

Since the support of $h_i$ coincides with that of $g_i$, we have that $0 \leq h_i \leq \mathbf{1}_{U_i}$. Moreover, as $g_{n+1} = 0$ on $K$, we have that $\sum_{i=1}^{n} h_i = 1$ on $K$. This completes the proof. ∎

**Theorem A.9** (Riesz representation). *If $S$ is a compact metric space and $T : C(S; \mathbb{R}) \to \mathbb{R}$ is a positive linear functional on the space of continuous functions on $S$, then there exists a unique measure $\mu$ on the Borel $\sigma$-algebra $\mathcal{B}(S)$ such that for all continuous functions $f \in C(S; \mathbb{R})$,*

$$T(f) = \int_S f \, d\mu. \tag{A.29}$$

*Proof.* We proceed in four steps. First, we will establish the uniqueness of $\mu$ by showing that it is entirely determined by $T$ on open sets. We will then define an outer measure $\mu^*$ on $2^S$ by leveraging the value that $\mu$ should attribute to open sets. The Caratheodory extension theorem will then guarantee that the restriction $\mu$ of this outer measure $\mu^*$ to the Borel $\sigma$-algebra is a measure. Finally we will show that $T$ can be identified with $\mu$ through equation (A.29).

*Step 1: uniqueness of $\mu$.* Fix a measure $\mu$ on $\mathcal{B}(S)$ for which (A.29) holds, and fix an open set $U \subseteq S$. For each integer $n \geqslant 1$, define the closed and therefore compact set
$$K_n := \left\{ x \in S \mid d(x, U^c) \geqslant \frac{1}{n} \right\}$$
in such a way that $U = \bigcup_{n=1}^\infty K_n$. Given an integer $n \geqslant 1$, invoke Urysohn's lemma to find a continuous function $f_n$ with $\mathbf{1}_{K_n} \leqslant f_n \leqslant \mathbf{1}_U$. Since $\mu$ satisfies (A.29), we have that $\mu(K_n) \leqslant T(f_n)$. It follows by the continuity of measure in Exercise A.1 that
$$\mu(U) = \lim_{n \to +\infty} \mu(K_n) \leqslant \limsup_{n \to +\infty} T(f_n) \leqslant \sup \{ T(f) \mid f \in C(S; \mathbb{R}) \text{ and } 0 \leqslant f \leqslant \mathbf{1}_U \}.$$

On the other hand, if $f$ is a continuous function with $0 \leqslant f \leqslant \mathbf{1}_U$, then (A.29) implies that $T(f) \leqslant \mu(U)$. Taking the supremum over all such functions $f$ shows that
$$\mu(U) = \sup \{ T(f) \mid f \in C(S; \mathbb{R}) \text{ and } 0 \leqslant f \leqslant \mathbf{1}_U \}. \tag{A.30}$$

This determines the measure $\mu$ on the $\pi$-system of open sets. Since $\mu$ is a finite measure, this in fact determines $\mu$ on the Borel $\sigma$-algebra by Exercise A.3.

*Step 2: defining the outer measure $\mu^*$.* Inspired by the value (A.30) that $\mu$ must take on open sets, given an open set $U \subseteq S$, let
$$\rho(U) := \sup \{ T(f) \mid f \in C(S; \mathbb{R}) \text{ and } 0 \leqslant f \leqslant \mathbf{1}_U \},$$
and define the function $\mu^* : 2^S \to [0, +\infty]$ by
$$\mu^*(A) := \inf \{ \rho(U) \mid A \subseteq U \text{ and } U \text{ is open} \}. \tag{A.31}$$

To prove that this defines an outer measure, it suffices to show that $\rho$ is sub-additive on the collection of open sets. Indeed, if this is the case, then
$$\mu^*(A) = \inf \left\{ \sum_{n=1}^{+\infty} \rho(U_n) \,\Big|\, A \subseteq \bigcup_{n=1}^\infty U_n \text{ and } U_n \text{ is open for every } n \geqslant 1 \right\}$$
which is an outer measure by Proposition A.1. Given a sequence $(U_n)_{n \geqslant 1}$ of open sets, let $U := \bigcup_{n=1}^\infty U_n$, and fix a continuous function $f$ with $0 \leqslant f \leqslant \mathbf{1}_U$. Denote by $K$ the closed support of $f$. Since $K$ is compact as a closed subset of the compact space $S$ and it is contained in $U$, it is possible to extract a finite sub-cover $(U_i)_{i \leqslant n}$ of $K$. By Lemma A.8 there exists a partition of unity $(h_i)_{i \leqslant n}$ on $K$ subordinate to the open cover $(U_i)_{i \leqslant n}$. Since $\sum_{i=1}^n h_i = 1$ on the support of $K$, we have that $f = \sum_{i=1}^n g h_i$. It follows by linearity of $T$ that
$$T(f) = \sum_{i=1}^n T(f h_i) \leqslant \sum_{i=1}^n \rho(U_i) \leqslant \sum_{n=1}^{+\infty} \mu(U_n),$$
where we have used that $0 \leqslant f h_i \leqslant \mathbf{1}_{U_i}$. Taking the supremum over all continuous functions $f$ with $0 \leqslant f \leqslant \mathbf{1}_U$ shows that $\rho$ is sub-additive on the collection of open sets and that $\mu^*$ is therefore an outer measure. The monotonicity of $\rho$ ensures that this outer measures agrees with $\rho$ on the collection of open sets.

*Step 3: $\mu^*$ restricts to a measure $\mu$ on $\mathcal{B}(S)$.* Since open sets generate the Borel $\sigma$-algebra, by the Caratheodory extension theorem, to show that the restriction $\mu$ of the outer measure $\mu^*$ to the Borel $\sigma$-algebra is a measure, it suffices to prove that every open set is Caratheodory measurable. By sub-additivity of outer measure this comes down to showing that for every open set $U$ and every subset $E$ of $S$ with $\mu^*(E) < +\infty$, we have

$$\mu^*(E) \geq \mu^*(E \cap U) + \mu^*(E \cap U^c). \tag{A.32}$$

Given $\varepsilon > 0$, let $V$ be an open set with $E \subseteq V$ and $\rho(V) \leq \mu^*(E) + \varepsilon$. By monotonicity of outer measure and the fact that $\mu^*$ coincides with $\rho$ on the collection of open sets,

$$\mu^*(E \cap U) + \mu^*(E \cap U^c) \leq \mu^*(V \cap U) + \mu^*(V \cap U^c) = \rho(V \cap U) + \mu^*(V \cap U^c).$$

To bound this further, find $f \in C(S;\mathbb{R})$ with $f \leq \mathbf{1}_{V \cap U}$ and $\rho(V \cap U) \leq T(f) + \varepsilon$. Denote by $K \subseteq U$ the compact support of $f$, and observe that by monotonicity of outer measure,

$$\mu^*(E \cap U) + \mu^*(E \cap U^c) \leq T(f) + \rho(V \cap K^c) + \varepsilon,$$

where we have used that $V \cap K^c$ is an open set. At this point, let $g \in C(S;\mathbb{R})$ be such that $g \leq \mathbf{1}_{V \cap K^c}$ and $\rho(V \cap K^c) \leq T(g) + \varepsilon$ so that

$$\mu^*(E \cap U) + \mu^*(E \cap U^c) \leq T(f + g) + 2\varepsilon.$$

As the supports of $f$ and $g$ are disjoint and both contained in $V$, we have that $f + g \leq \mathbf{1}_V$. It follows by definition of $\rho$ that

$$\mu^*(E \cap U) + \mu^*(E \cap U^c) \leq \rho(V) + 2\varepsilon \leq \mu^*(E) + 3\varepsilon.$$

Letting $\varepsilon$ tend to zero establishes (A.32) and shows that $\mu := \mu^*|_{\mathcal{B}(S)}$ defines a measure on the Borel $\sigma$-algebra $\mathcal{B}(S)$.

*Step 4: $T$ can be identified with $\mu$.* To show that $T$ can be identified with the measure $\mu$ through equation (A.29), up to replacing $f$ by $-f$, it suffices to prove that for every $f \in C(S;\mathbb{R})$,

$$T(f) \leq \int_S f \, d\mu. \tag{A.33}$$

Moreover, up to replacing $f$ by $f + \|f\|_\infty$, it suffices to establish (A.33) for non-negative functions $f \in C(S;\mathbb{R})$. Fix a non-negative $f \in C(S;\mathbb{R})$, and write $K$ for its compact support. Since $f$ is bounded, its range belongs to some bounded interval $[a,b] \subseteq \mathbb{R}$. Given $\varepsilon > 0$, let $a = y_0 < y_1 < \ldots < y_{n-1} < y_n = b$ be a partition of the interval $[a,b]$ with $y_{i+1} - y_i < \varepsilon$ for $0 \leqslant i \leqslant n-1$. Introduce the collection of disjoint Borel sets $E_i := f^{-1}((y_{i-1}, y_i]) \cap K$ which partition $K = \bigcup_{i=1}^n E_i$. For each $1 \leqslant i \leqslant n$, leveraging the continuity of $f$ and the definition (A.31) of the outer measure $\mu^*$, find an open set $U_i$ containing $E_i$ with

$$\rho(U_i) \leqslant \mu(E_i) + \frac{\varepsilon}{n} \quad \text{and} \quad f(x) < y_i + \varepsilon \text{ on } U_i. \tag{A.34}$$

Invoke Lemma A.8 to find a partition of unity $(h_i)_{i \leqslant n}$ on $K$ subordinate to the open cover $(U_i)_{i \leqslant n}$. Since $f = \sum_{i=1}^n f h_i$ and $f h_i \leqslant (y_i + \varepsilon) h_i$, the linearity of $T$ implies that

$$T(f) = \sum_{i=1}^n T(f h_i) \leqslant \sum_{i=1}^n (y_i + \varepsilon) T(h_i) \leqslant \sum_{i=1}^n (y_i + \varepsilon) \rho(U_i).$$

It follows by the properties (A.34) of the open sets $(U_i)_{i \leqslant n}$ that

$$T(f) \leqslant \sum_{i=1}^n (y_i + \varepsilon)\left(\mu(E_i) + \frac{\varepsilon}{n}\right) \leqslant \sum_{i=1}^n (y_i - \varepsilon) \mu(E_i) + \varepsilon(2\mu(K) + \|f\|_\infty + \varepsilon).$$

Combining this with the fact that $f(x) > y_{i-1} > y_i - \varepsilon$ on $E_i$ and letting $\varepsilon$ tend to zero establishes (A.33) and completes the proof. ∎

**Exercise A.5.** Show that any positive linear functional $T : C(S; \mathbb{R}) \to \mathbb{R}$ on the space of continuous functions on a compact metric space $S$ is continuous with respect to the uniform norm.

## A.3  The Stone-Weierstrass theorem

In this section, we prove a generalization of the well-known theorem of Weierstrass [262] which states that any continuous function on a compact interval $[a,b]$ is the uniform limit of polynomials on $[a,b]$. A probabilistic proof of this classical result leveraging Bernstein polynomials is outlined in Exercise A.6. Following the work of Stone [240], we will generalize this result in two main directions. We will consider continuous functions on a compact metric space $S$ as opposed to compact sub-intervals of the real line, and we will be concerned with the density of collections of functions which form an algebra in the space of continuous functions on $S$ as opposed to focusing on polynomials. We will prove a version of the Stone-Weierstrass theorem for real-valued functions and another for complex-valued functions. The Stone-Weierstrass theorem is a very helpful tool as it often allows us to deduce an identity for all continuous functions once it has been proved for a smaller family of functions such as monomials. In Section A.6, we will use the real version of the Stone-Weierstrass theorem to show that a non-negative random variable is entirely determined by its Laplace transform, and we will leverage the complex version of the Stone-Weierstrass theorem to prove that a random variable is entirely determined by its characteristic function.

The real version of the Stone-Weierstrass theorem states that a collection of functions $\mathcal{F}$ which forms an algebra in the space of continuous real-valued functions on a compact metric space $S$ is dense in $C(S;\mathbb{R})$ endowed with the uniform norm,

$$\|f-g\|_\infty = \sup_{x \in S} |f(x) - g(x)|, \tag{A.35}$$

provided it contains constants and separates points. An *algebra* in $C(S;\mathbb{R})$ is a real vector subspace $\mathcal{F}$ of $C(S;\mathbb{R})$ which is closed under multiplication,

(i) for all $\lambda, \mu \in \mathbb{R}$ and $f, g \in \mathcal{F}$, we have $\lambda f + \mu g \in \mathcal{F}$,

(ii) for all $f, g \in \mathcal{F}$, we have $fg \in \mathcal{F}$.

A collection $\mathcal{F}$ of real-valued continuous functions on $S$ *separates points* if, for all distinct $x \neq y \in S$, there exists a function $f \in \mathcal{F}$ with $f(x) \neq f(y)$.

**Theorem A.10** (Stone-Weierstrass real version). *If $(S,d)$ is a compact metric space and $\mathcal{F}$ is an algebra in $C(S;\mathbb{R})$ which separates points and contains constants, then $\mathcal{F}$ is dense in $(C(S;\mathbb{R}), \|\cdot\|_\infty)$.*

*Proof.* The proof proceeds in two steps. First, we will show that $\overline{\mathcal{F}}$ is a lattice, in the sense that for all $f, g \in \overline{\mathcal{F}}$, we have $\min(f,g) \in \overline{\mathcal{F}}$ and $\max(f,g) \in \overline{\mathcal{F}}$. We will then fix $f \in C(S;\mathbb{R})$ as well as $\varepsilon > 0$, and we will use the lattice property of $\overline{\mathcal{F}}$ to construct an explicit function $g \in \overline{\mathcal{F}}$ with $\|f - g\|_\infty \leq \varepsilon$. This will establish the density of $\overline{\mathcal{F}}$ in $C(S;\mathbb{R})$. The density of $\mathcal{F}$ in $C(S;\mathbb{R})$ is then immediate from its density in $\overline{\mathcal{F}}$.

*Step 1: $\overline{\mathcal{F}}$ is a lattice.* To show that $\overline{\mathcal{F}}$ is a lattice, it suffices to show that for any $f \in \overline{\mathcal{F}}$, we have $|f| \in \overline{\mathcal{F}}$. Indeed, this implies that for any $f, g \in \overline{\mathcal{F}}$, we have

$$\min(f,g) = \frac{f+g-|f-g|}{2} \in \overline{\mathcal{F}} \quad \text{and} \quad \max(f,g) = \frac{f+g+|f-g|}{2} \in \overline{\mathcal{F}}.$$

Fix $f \in \overline{\mathcal{F}}$, and let $M > 0$ be large enough so the range of $f$ is contained in the interval $[-M, M]$. By the Weierstrass approximation theorem in Exercise A.6, the absolute value function $x \mapsto |x|$ defined on the compact interval $[-M, M]$ can be approximated by polynomials of $x$. This means that $|f(x)|$ can be approximated by polynomials of $f(x)$. Since $\mathcal{F}$ is an algebra which contains constants, any polynomial of $f(x)$ belongs to $\overline{\mathcal{F}}$, so $|f| \in \overline{\mathcal{F}}$.

*Step 2: constructing an approximating function.* Observe that for any pair of distinct points $s \neq t \in S$ and all values $c, d \in \mathbb{R}$, there exists $g \in \mathcal{F}$ with $g(s) = c$ and $g(t) = d$.

Indeed, as $\mathcal{F}$ separates points, there exists $h \in \mathcal{F}$ with $h(s) \neq h(t)$, and the function $g \in \mathcal{F}$ defined by

$$g(x) := c + (d-c)\frac{h(x) - h(s)}{h(t) - h(s)}$$

satisfies the desired properties. We now fix a continuous function $f \in C(S;\mathbb{R})$ as well as $\varepsilon > 0$. For all pairs of distinct points $s \neq t \in S$, we write $g_{s,t} \in \mathcal{F}$ for a function with $g_{s,t}(s) = f(s)$ and $g_{s,t}(t) = f(t)$. By continuity, the function $g_{s,t}$ approximates $f$ in small enough neighbourhoods of both $s$ and $t$. To construct a function $g \in \overline{\mathcal{F}}$ which approximates $f$ everywhere, fix $s \in S$, and consider the open neighbourhood of $t$,

$$U_t := \{x \in S \mid g_{s,t}(x) < f(x) + \varepsilon\}.$$

Since $(U_t)_{t \in S}$ defines an open cover of the compact set $S$, there exists a finite sub-cover $(U_{t_i})_{i \leqslant n}$. The function $g_s : S \to \mathbb{R}$ defined by

$$g_s(x) := \min_{1 \leqslant i \leqslant n} g_{s,t}(x)$$

is such that $g_s(s) = f(s)$ and $g_s(x) < f(x) + \varepsilon$ for all $x \in S$. Moreover, we have that $g_s \in \overline{\mathcal{F}}$ as $\overline{\mathcal{F}}$ is a lattice. Similarly, fix $t \in S$ and consider the open neighbourhood of $s$,

$$V_s := \{x \in S \mid g_s(x) > f(x) - \varepsilon\}.$$

Since $(V_s)_{s \in S}$ defines an open cover of the compact set $S$, there exists a finite sub-cover $(U_{s_j})_{j \leqslant m}$. The function $g : S \to \mathbb{R}$ defined by

$$g(x) := \max_{1 \leqslant j \leqslant m} g_s(x)$$

belongs to the lattice $\overline{\mathcal{F}}$ and is such that $\|f - g\|_\infty \leqslant \varepsilon$. This completes the proof. ∎

As shown in Exercise A.7, the real version of the Stone-Weierstrass theorem does not generalize to the setting of complex-valued continuous functions on $S$. It turns out that in addition to containing constants and separating points, for an algebra $\mathcal{F}$ in $C(S;\mathbb{C})$ to be dense with respect to the uniform norm,

$$\|f - g\|_\infty = \|\mathrm{Re}(f) - \mathrm{Re}(g)\|_\infty + \|\mathrm{Im}(f) - \mathrm{Im}(g)\|_\infty, \qquad (A.36)$$

it must also be closed under conjugation. This means that for all $f \in \mathcal{F}$, we need to have $\overline{f} \in \mathcal{F}$. It is worth emphasizing that an algebra in $C(S;\mathbb{C})$ is a complex vector subspace of $C(S;\mathbb{C})$ which is closed under multiplication. In particular, the first condition (i) in the definition of an algebra must hold for all $\lambda, \mu \in \mathbb{C}$.

**Theorem A.11** (Stone-Weierstrass complex version). *If $(S,d)$ is a compact metric space and $\mathcal{F}$ is an algebra in $C(S;\mathbb{C})$ which separates points, contains constants and is closed under complex conjugation, then $\mathcal{F}$ is dense in $(C(S;\mathbb{C}), \|\cdot\|_\infty)$.*

*Proof.* For each $f \in \mathcal{F}$, we have

$$\mathrm{Re}(f) = \frac{f + \overline{f}}{2} \in \mathcal{F} \quad \text{and} \quad \mathrm{Im}(f) = \frac{f - \overline{f}}{2i} \in \mathcal{F},$$

so the set $\mathcal{F}_\mathbb{R}$ of real and imaginary parts of functions in $\mathcal{F}$ is contained in $\mathcal{F}$, and it is an algebra of real-valued functions which separates points and contains constants. It follows by the real version of the Stone-Weierstrass theorem that $\mathcal{F}_\mathbb{R}$ is dense in $(C(S;\mathbb{R}), \|\cdot\|_\infty)$. Since every function $f \in C(S;\mathbb{C})$ may be written as $\mathrm{Re}(f) + i\mathrm{Im}(f)$ for the real-valued functions $\mathrm{Re}(f), \mathrm{Im}(f) \in C(S;\mathbb{R})$, and a sequence of complex-valued functions converges with respect to $\|\cdot\|_\infty$ if and only if its associated sequences of real and imaginary parts converge with respect to $\|\cdot\|_\infty$, we have that $\mathcal{F}$ is dense in $(C(S;\mathbb{C}), \|\cdot\|_\infty)$. This completes the proof. ∎

**Exercise A.6** (Weierstrass approximation). The purpose of this exercise is to give a probabilistic proof of the Stone-Weierstrass theorem on a compact interval $[a,b] \subseteq \mathbb{R}$. That is, we will show that for every $\varepsilon > 0$ and continuous function $f \in ([a,b];\mathbb{R})$, there exists a polynomial $P$ on $[a,b]$ with $\|f - P\|_\infty \leq \varepsilon$.

 (i) Argue that there is no loss in generality in considering the interval $[0,1]$.

(ii) Fix $\varepsilon > 0$ and $f \in C([0,1];\mathbb{R})$. Given $x \in [0,1]$, let $(X_n)_{n \geq 1}$ be a sequence of independent and identically distributed random variables with $X_1 \sim \mathrm{Ber}(x)$. For each $n \geq 1$, consider the sample average $S_n := \frac{1}{n}\sum_{i=1}^n X_i$, and define the function $P_n(x) := \mathbb{E}f(S_n)$. Prove that $P_n$ is a polynomial and that for $n$ large enough, we have $\|f - P\|_\infty \leq \varepsilon$.

**Exercise A.7.** Denote by $\mathcal{C} := \{e^{it} \mid t \in [0,2\pi]\} \subseteq \mathbb{C}$ the unit circle, and consider the complex-valued function $f : \mathcal{C} \to \mathbb{C}$ defined by $f(z) := \overline{z}$. Prove that $f$ cannot be approximated in $(C(\mathcal{C},\mathbb{C}), \|\cdot\|_\infty)$ by polynomials (of the $z$ variable).

## A.4   The Lebesgue differentiation theorem

Before moving from the realm of analysis to that of probability theory, we prove the Lebesgue differentiation theorem on Euclidean space $S = \mathbb{R}^d$. In the one-dimensional setting, $d = 1$, this result generalizes the fundamental theorem of calculus.

**Proposition A.12** (Fundamental theorem of calculus). *If $f : \mathbb{R} \to \mathbb{R}$ is continuous and $a \in \mathbb{R}$, then the integral function $F : \mathbb{R} \to \mathbb{R}$ defined by*

$$F(x) := \int_a^x f(t)\, dt \tag{A.37}$$

*is differentiable and $F' = f$.*

*Proof.* Fix $x \in \mathbb{R}$ as well as $\varepsilon > 0$. The continuity of $f$ gives $\delta > 0$ with $|f(t) - f(x)| < \varepsilon$ whenever $|x - t| < \delta$. If $h \in \mathbb{R}$ is such that $|h| < \delta$, then

$$\left| \frac{F(x+h) - F(x)}{h} - f(x) \right| \leq \frac{1}{h} \int_x^{x+h} |f(t) - f(x)| \, dt \leq \varepsilon.$$

This completes the proof. ∎

The Lebesgue differentiation theorem will replace the continuity assumption of the integrand by local integrability and will weaken the conclusion to almost everywhere differentiability of the integral function. In higher dimensions, we will show in Theorem A.16 that if $f \in L^1_{\mathrm{loc}}(\mathbb{R}^d; \mathbb{R})$, then for almost every $x \in \mathbb{R}^d$,

$$\lim_{r \searrow 0} \frac{1}{m(B_r(x))} \int_{B_r(x)} |f(y) - f(x)| \, dy = 0, \tag{A.38}$$

where we denote by $B_r(x)$ the open Euclidean ball of radius $r > 0$ centred at $x \in \mathbb{R}^d$, and by $m(B_r(x))$ its Lebesgue measure. When $d = 1$, this readily implies that, for almost every $x \in \mathbb{R}$,

$$\lim_{r \to 0} \frac{1}{r} \int_x^{x+r} f(t) \, dt = f(x). \tag{A.39}$$

To prove (A.38), we will apply the theory of maximal operators to the Hardy-Littlewood maximal operator defined on the space $L^1_{\mathrm{loc}}(\mathbb{R}^d; \mathbb{R})$ of locally integrable functions on $\mathbb{R}^d$ by

$$H^* f(x) := \sup_{r > 0} \frac{1}{m(B_r(x))} \left| \int_{B_r(x)} f(y) \, dy \right|, \tag{A.40}$$

Although the basic result that we will use from the theory of maximal operators holds in far greater generality, we will only state it on Euclidean space $\mathbb{R}^d$ endowed with $d$-dimensional Lebesgue measure $m$ in the context of operators on $L^1(\mathbb{R}^d; \mathbb{R})$.

The *maximal operator* associated with a family $(T_t)_{t>0}$ of linear operators from $L^1(\mathbb{R}^d; \mathbb{R})$ to the space of measurable functions on $\mathbb{R}^d$ is the operator $T^*$ from $L^1(\mathbb{R}^d; \mathbb{R})$ to the space of measurable functions on $\mathbb{R}^d$ defined by

$$T^* f(x) = \sup_{t > 0} |T_t f(x)| \tag{A.41}$$

The Hardy-Littlewood maximal operator (A.40) is the maximal operator associated with the family $(H_r)_{r>0}$ of linear operators defined on $L^1(\mathbb{R}^d; \mathbb{R})$ by

$$H_r f(x) := \frac{1}{m(B_r(x))} \int_{B_r(x)} f(y) \, dy. \tag{A.42}$$

To use the Hardy-Littlewood maximal operator to establish the almost everywhere convergence (A.38), we will leverage the fact that it is weak-$(1,1)$. An operator $T^*$ from $L^1(\mathbb{R}^d;\mathbb{R})$ to the space of measurable functions on $\mathbb{R}^d$ is *weak-$(1,1)$* if there exists a constant $C > 0$ such that for all $f \in L^1(\mathbb{R}^d;\mathbb{R})$,

$$m\{|T^*f(x)| > \lambda\} \leq \frac{C}{\lambda} \int_{\mathbb{R}^d} |f(x)| \, dx. \tag{A.43}$$

It turns out that whenever the maximal operator of a family of linear operators is weak-$(1,1)$, the space of functions where the sequence of operators converges pointwise to the identity is closed in $L^1(\mathbb{R}^d;\mathbb{R})$. In particular, if it contains a dense subspace such as the space of continuous functions, it must be all of $L^1(\mathbb{R}^d;\mathbb{R})$.

**Proposition A.13.** *If $(T_t)_{t>0}$ is a family of linear operators from $L^1(\mathbb{R}^d;\mathbb{R})$ to the space of measurable functions on $\mathbb{R}^d$ whose maximal operator $T^*$ is weak-$(1,1)$, then*

$$\left\{ f \in L^1(\mathbb{R}^d;\mathbb{R}) \;\middle|\; \lim_{t\to 0} T_t f(x) = f(x) \text{ for almost every } x \in \mathbb{R}^d \right\} \tag{A.44}$$

*is closed in $L^1(\mathbb{R}^d;\mathbb{R})$.*

*Proof.* Let $(f_n)_{n \geq 1}$ be a sequence of functions converging to some $f \in L^1(\mathbb{R}^d;\mathbb{R})$ with the property that $\lim_{t \to 0} T_t f_n(x) = f_n(x)$ almost everywhere. The triangle inequality and Chebyshev's inequality imply that for any $\lambda > 0$,

$$\begin{aligned} m\left\{ \limsup_{t \to 0} |T_t f(x) - f(x)| > \lambda \right\} &\leq m\left\{ \limsup_{t \to 0} |T_t f(x) - T_t f_n(x)| > \lambda/2 \right\} \\ &\quad + m\left\{ \limsup_{t \to 0} |f(x) - f_n(x)| > \lambda/2 \right\} \\ &\leq m\{ T^*(f - f_n)(x) > \lambda/2 \} + \frac{2}{\lambda} \int_{\mathbb{R}^d} |f - f_n| \, dx \\ &\leq \frac{2(C+1)}{\lambda} \int_{\mathbb{R}^d} |f - f_n| \, dx. \end{aligned}$$

Letting $n$ tend to infinity and leveraging the sub-additivity of measure shows that the measure of the set where $\limsup_{t \to 0} |T_t f(x) - f(x)| > 0$ is equal to zero. This is equivalent to the statement that $\lim_{t \to 0} T_t f(x) = f(x)$ for almost every $x \in \mathbb{R}^d$, which means that the set (A.44) is closed as required. ∎

Proving the Lebesgue differentiation theorem therefore comes down to showing that the Hardy-Littlewood maximal operator (A.40) is weak-$(1,1)$. This will follow from the Vitali covering lemma, which essentially says that the measure of an arbitrary union of balls is concentrated on a finite disjoint sub-collection of these balls.

**Lemma A.14** (Vitali covering). *If $\mathcal{C}$ is a collection of open balls in $\mathbb{R}^d$ and $U$ denotes their union, then for any $c < m(U)$, there exist disjoint balls $(B_i)_{i \leqslant n} \subseteq \mathcal{C}$ with*

$$c < 3^d \sum_{i=1}^{n} m(B_i). \tag{A.45}$$

*Proof.* Invoking Exercise A.4 allows us to find a compact set $K \subseteq U$ with $m(K) > c$. Since $\mathcal{C}$ is an open cover of the compact set $K$, it is possible to extract a finite sub-cover $(C_j)_{j \leqslant m} \subseteq \mathcal{C}$. We now define the sequence $(B_i)_{i \leqslant n}$ iteratively by taking $B_1$ to be the ball $C_j$ of largest radius, and by taking $B_i$ to be the ball $C_j$ of largest radius disjoint from $\bigcup_{k=1}^{i-1} B_k$. This construction ensures that if a ball $C_j$ is not selected, then it must intersect one of the selected balls $B_i$, and if $i$ is the smallest index such that this intersection takes place, then the radius of $C_j$ is at most the radius of $B_i$. This means that $K \subseteq \bigcup_{i=1}^{n} 3B_i$, where $3B_i$ denotes the ball with the same centre as $B_i$ but with three times the radius. It follows by sub-additivity of measure that

$$c < m(K) \leqslant \sum_{i=1}^{n} m(3B_i) = 3^d \sum_{i=1}^{n} m(B_i).$$

The last equality uses the homogeneity of $d$-dimensional Lebesgue measure which is immediate from the Dynkin $\pi$-$\lambda$ theorem and the fact that for every $a > 0$, we have $m(aR) = a^d m(R)$ for any finite union of disjoint rectangles $R$. This completes the proof. ∎

**Lemma A.15** (Hardy-Littlewood). *The Hardy-Littlewood maximal operator (A.40) is weak-$(1,1)$.*

*Proof.* Fix $\lambda > 0$, and let $E_\lambda := \{H^* f(x) > \lambda\}$. Given $x \in E_\lambda$, let $r_x > 0$ be such that

$$\frac{1}{m(B_{r_x}(x))} \int_{B_{r_x}(x)} |f(y)| \, dy \geqslant \frac{1}{m(B_r(x))} \left| \int_{B_r(x)} f(y) \, dy \right| > \lambda.$$

The collection $(B_{r_x}(x))_{x \in E_\lambda}$ covers $E_\lambda$, so, given any $c < m(E_\lambda)$, the Vitali covering theorem gives a finite sub-collection $(B_{r_{x_i}}(x_i))_{i \leqslant n}$ of disjoint balls with

$$c \leqslant 3^d \sum_{i=1}^{n} m(B_{r_{x_i}}(x_i)) \leqslant \frac{3^d}{\lambda} \sum_{i=1}^{n} \int_{B_{r_{x_i}}(x_i)} |f(y)| \, dy \leqslant \frac{3^d}{\lambda} \int_{\mathbb{R}^d} |f(y)| \, dy.$$

Letting $c$ increase to $m(E_\lambda)$ completes the proof. ∎

**Theorem A.16** (Lebesgue differentiation). *If $f \in L^1_{\text{loc}}(\mathbb{R}^d; \mathbb{R})$, then for almost every $x \in \mathbb{R}^d$,*

$$\lim_{r \to 0} \frac{1}{m(B_r(x))} \int_{B_r(x)} |f(y) - f(x)| \, dy = 0. \tag{A.46}$$

*Proof.* The proof proceeds in two steps. First we show that for almost every $x \in \mathbb{R}^d$,

$$\lim_{r \to 0} \frac{1}{m(B_r(x))} \int_{B_r(x)} f(y) \, dy = f(x), \tag{A.47}$$

and then we leverage this to prove (A.46).

*Step 1: proving* (A.47). Combining Lemma A.15 with Proposition A.13 shows that the set of functions where (A.47) holds is closed in $L^1(\mathbb{R}^d; \mathbb{R})$. Since it contains the dense set of continuous functions by the higher-dimensional version of Proposition A.12, it must be all of $L^1(\mathbb{R}^d; \mathbb{R})$. Noticing that any locally integrable function on $\mathbb{R}^d$ can be extended to an integrable function on $\mathbb{R}^d$ that agrees with it on a small enough ball around a given point $x \in \mathbb{R}^d$ establishes (A.47) for every locally integrable function.

*Step 2: proving* (A.46). Fix $c \in \mathbb{R}$, and observe that (A.47) applied to the function $f_c(x) := |f(x) - c|$ implies that, for almost every $x \in \mathbb{R}^d$,

$$\lim_{r \to 0} \frac{1}{m(B_r(x))} \int_{B_r(x)} |f(y) - c| \, dy = |f(x) - c|. \tag{A.48}$$

By additivity of measure and countability of the rationals, we can in fact make this convergence uniform for $c \in \mathbb{Q}$. More precisely, it is possible to find a null set $N \subseteq \mathbb{R}^d$ such that for every $x \in \mathbb{R}^d \setminus N$ and $c \in \mathbb{Q}$ the convergence in (A.48) holds. If we now fix $x \in \mathbb{R}^d \setminus N$ as well as $\varepsilon > 0$, then, by density of the rationals, it is possible to find $c \in \mathbb{Q}$ with $|f(x) - c| \leq \varepsilon$. It follows that

$$\frac{1}{m(B_r(x))} \int_{B_r(x)} |f(y) - f(x)| \, dy \leq \frac{1}{m(B_r(x))} \int_{B_r(x)} |f(y) - c| \, dy + |f(x) - c|.$$

Taking the limsup as $r$ tends to zero and then letting $\varepsilon$ tend to zero completes the proof. ∎

## A.5   A topological characterization of weak convergence

For the remainder of this chapter we will focus on the weak convergence of probability measures on a metric space $S$ endowed with its Borel $\sigma$-algebra $\mathcal{B}(S)$. We write $C_b(S; \mathbb{R})$ for the space of bounded and continuous functions on $S$. A sequence $(\mathbb{P}_n)_{n \geq 1}$ of probability measures on $S$ *converges weakly* to a probability measure $\mathbb{P}$ on $S$ if, for every bounded and continuous function $f \in C_b(S; \mathbb{R})$,

$$\lim_{n \to +\infty} \int_S f \, d\mathbb{P}_n = \int_S f \, d\mathbb{P}. \tag{A.49}$$

A sequence $(X_n)_{n \geq 1}$ of random elements on $S$ *converges in law* to a random element $X$ on $S$ if the sequence of laws of the $X_n$ converges weakly to the law of $X$.

The *law* of the random element $X$ on $S$ is the probability measure $\mathbb{P}_X$ on $S$ defined by
$$\mathbb{P}_X(A) := \mathbb{P}\{X \in A\}. \tag{A.50}$$

This notion of convergence is motivated by the fact that, as shown in Exercise A.8, a probability measure is determined by its action on the space $C_b(S;\mathbb{R})$ of continuous and bounded functions. In the setting of the real line, $S = \mathbb{R}$, Exercise A.9 outlines a proof of the fact that convergence in law coincides with the familiar convergence in distribution seen in some introductory probability courses. In the setting of general metric spaces, the weak convergence of probability measures can be characterized topologically through the *Portmanteau theorem*.

**Theorem A.17** (Portmanteau). *For a sequence $(\mathbb{P}_n)_{n \geq 1}$ of probability measures on S, the following statements are equivalent.*

 (i) *$(\mathbb{P}_n)_{n \geq 1}$ converges weakly to the probability measure $\mathbb{P}$.*

 (ii) *(A.49) holds for all bounded Lipschitz functions $f : S \to \mathbb{R}$.*

 (iii) *For any open set $U \subseteq S$, we have $\liminf_{n \to +\infty} \mathbb{P}_n(U) \geq \mathbb{P}(U)$.*

 (iv) *For any closed set $F \subseteq S$, we have $\limsup_{n \to +\infty} \mathbb{P}_n(F) \leq \mathbb{P}(F)$.*

 (v) *For any set $A \in \mathcal{B}(S)$ with $\mathbb{P}(\partial A) = 0$, we have $\lim_{n \to +\infty} \mathbb{P}_n(A) = \mathbb{P}(A)$.*

*Proof.* (i) $\Rightarrow$ (ii): This is immediate from the fact that any bounded Lipschitz function on $S$ is continuous.

(ii) $\Rightarrow$ (iii): Let $U$ be an open set and define the sequence of bounded Lipschitz functions $(f_m)_{m \geq 1}$ on $S$ by $f_m(s) := \min(1, m d(s, U^c))$. Since $U^c$ is closed, the sequence $(f_m)_{m \geq 1}$ increases to $\mathbf{1}_U$. It follows that
$$\int_S f_m \, d\mathbb{P}_n \leq \mathbb{P}_n(U).$$

Using (ii) to let $n$ tend to infinity and the monotone convergence theorem to let $m$ tend to infinity establishes (iii).

(iii) $\Leftrightarrow$ (iv): This is immediate by taking complements.

(iv) $\Rightarrow$ (v): Let $A \subseteq S$ be a measurable set with $\mathbb{P}(\partial A) = 0$. Since $\text{int}(A)$ is open, $\overline{A}$ is closed and $\text{int}(A) \subseteq \overline{A}$, assumption (iv) and its equivalence to (iii) imply that
$$\mathbb{P}(\text{int}(A)) \leq \liminf_{n \to +\infty} \mathbb{P}_n(\text{int} A) \leq \limsup_{n \to +\infty} \mathbb{P}_n(\overline{A}) \leq \mathbb{P}(\overline{A}).$$

Combining this with the fact that $\mathbb{P}(\overline{A}) = \mathbb{P}(\text{int} A) = \mathbb{P}(A)$ by the assumption that $\mathbb{P}(\partial A) = 0$ establishes (v).

(v) $\Rightarrow$ (i): Fix $f \in C_b(S;\mathbb{R})$, and for $y \in \mathbb{R}$ let $F_y := \{s \in S : f(s) = y\}$ be the level set

of $f$. Since the sets $(F^y)_{y\in\mathbb{R}}$ are disjoint, there exist at most countably many values of $y$ for which $\mathbb{P}(F_y) > 0$. Given $\varepsilon > 0$ it is therefore possible to find an integer $N \geqslant 1$ and a sequence $a_0 \leqslant \cdots \leqslant a_N$ with $\max_k(a_k - a_{k-1}) \leqslant \varepsilon$ and $\mathbb{P}(F_{a_k}) = 0$ such that $f$ maps into the interval $(a_0, a_N)$. Introduce the function $f_\varepsilon : S \to \mathbb{R}$ defined by

$$f_\varepsilon := \sum_{k=0}^{N-1} a_k \mathbf{1}_{B_k},$$

where $B_k := \{s \in S : a_k \leqslant f(s) < a_{k+1}\}$. Since $f$ is continuous, for $0 \leqslant k \leqslant N-1$, we have $\partial B_k \subseteq F_{a_k} \cup F_{a_{k+1}}$, and therefore $\mathbb{P}(\partial B_k) = 0$. It follows by (v) that

$$\lim_{n \to +\infty} \int_S f_\varepsilon \, d\mathbb{P}_n = \sum_{k=0}^{N-1} a_k \mathbb{P}(B_k) = \int_S f_\varepsilon \, d\mathbb{P}.$$

Using that $\|f_\varepsilon - f\|_\infty \leqslant \varepsilon$ to let $\varepsilon$ tend to zero completes the proof. ∎

**Remark A.18.** A similar argument as that used to show that (ii) ⇒ (iii) can be used to prove that (ii) ⇒ (iv). For future reference, let us present this argument in the more general setting where we have a sequence of finite measures $(\mu_n)_{n\geqslant 1}$ and a finite measure $\mu$ with the property that for every $f \in C_b(S;\mathbb{R})$,

$$\lim_{n \to +\infty} \int_S f \, d\mu_n = \int_S f \, d\mu. \tag{A.51}$$

Given a closed set $F \subseteq S$, define the sequence of bounded and Lipschitz functions $(g_m)_{m\geqslant 1}$ on $S$ by $g_m(s) := \max(0, 1 - md(x, F))$. Since $F$ is closed, the sequence $(g_m)_{m\geqslant 1}$ decreases to $\mathbf{1}_F$. It follows that

$$\int_S g_m \, d\mu_n \geqslant \mu_n(F). \tag{A.52}$$

Using (A.51) to let $n$ tend to infinity and then the dominated convergence theorem to let $m$ tend to infinity shows that $\limsup_{n \to +\infty} \mu_n(F) \leqslant \mu(F)$ as required.

**Exercise A.8.** Consider two probability measures $\mathbb{P}$ and $\mathbb{Q}$ on $S$ such that, for every $f \in C_b(S;\mathbb{R})$, we have

$$\int_S f \, d\mathbb{P} = \int_S f \, d\mathbb{Q}. \tag{A.53}$$

Show that $\mathbb{P} = \mathbb{Q}$.

**Exercise A.9.** Let $(\mathbb{P}_n)_{n\geqslant 1}$ be a sequence of probability measures on $\mathbb{R}$, and let $\mathbb{P}$ be a probability measure on $\mathbb{R}$. Denote by $F_n : \mathbb{R} \to [0,1]$ and $F : \mathbb{R} \to [0,1]$ the distribution functions associated with $\mathbb{P}_n$ and $\mathbb{P}$,

$$F_n(t) := \mathbb{P}((-\infty, t]) \quad \text{and} \quad F(t) := \mathbb{P}((-\infty, t]). \tag{A.54}$$

(i) Prove that $(\mathbb{P}_n)_{n\geqslant 1}$ converges weakly to $\mathbb{P}$ if and only if at every point $t$ of continuity of $F$, the numerical sequence $(F_n(t))_{n\geqslant 1}$ converges to $F(t)$.

(ii) Is restricting to points of continuity of $F$ necessary?

**Exercise A.10.** The goal of this exercise is to show the continuous mapping theorem. Let $S, S'$ be two metric spaces, let $f : S \to S'$ be a measurable function, and let

$$\mathscr{C}_f := \{x \in S \mid f \text{ is continuous at } x\}. \tag{A.55}$$

Let $(X_n)_{n\geqslant 1}$ and $X_\infty$ be random variables taking values in $S$ with $(X_n)_{n\geqslant 1}$ converging in law to $X_\infty$ as $n$ tends to infinity.

(i) Show that the set $\mathscr{C}_f$ is measurable.

(ii) Assuming that $\mathbb{P}\{X_\infty \in \mathscr{C}_f\} = 1$, show that $f(X_n)$ converges in law to $f(X_\infty)$ as $n$ tends to infinity.

**Exercise A.11.** Let $(Y_{k,n})_{k,n\geqslant 1}$ and $(Y_n)_{n\geqslant 1}$ be random elements with values in $S$ such that for every $\varepsilon > 0$,

$$\lim_{k\to+\infty} \limsup_{n\to+\infty} \mathbb{P}\{d(Y_{k,n}, Y_n) \geqslant \varepsilon\} = 0. \tag{A.56}$$

Suppose that there exist a sequence of random elements $(Z_k)_{k\geqslant 1}$ and a random element $Y$ such that, for every $k \geqslant 1$, the sequence $(Y_{k,n})_{n\geqslant 1}$ converges in law to $Z_k$ as $n$ tends to infinity, and such that the sequence $(Z_k)_{k\geqslant 1}$ converges in law to $Y$ as $k$ tends to infinity. Show that $(Y_n)_{n\geqslant 1}$ converges in law to $Y$.

## A.6   Weak convergence through tightness and uniqueness

The topological characterization of weak convergence in the Portmanteau theorem is very useful to prove abstract results; however, it is difficult to use to establish the weak convergence of a specific sequence of probability measures. Typically, showing that a sequence $(\mathbb{P}_n)_{n\geqslant 1}$ of probability measures converges weakly is a two-step process. First, we show that the sequence $(\mathbb{P}_n)_{n\geqslant 1}$ is relatively compact, in the sense that any subsequence of $(\mathbb{P}_n)_{n\geqslant 1}$ admits a weakly convergent subsequence, and then we show that there is only one possible subsequential limit.

**Lemma A.19.** *Fix a sequence of probability measures $(\mathbb{P}_n)_{n\geqslant 1}$. If there exists a probability measure $\mathbb{P}$ with the property that any sequence $(n(k))_{k\geqslant 1}$ admits a further subsequence $(n(k(r)))_{r\geqslant 1}$ such that $(\mathbb{P}_{n(k(r))})_{r\geqslant 1}$ converges weakly to $\mathbb{P}$, then $(\mathbb{P}_n)_{n\geqslant 1}$ itself converges weakly to $\mathbb{P}$.*

*Proof.* Suppose for the sake of contradiction that the sequence $(\mathbb{P}_n)_{n \geqslant 1}$ does not converge weakly to $\mathbb{P}$. This means that there exist a continuous and bounded function $f \in C_b(S; \mathbb{R})$, some $\varepsilon > 0$ and a sequence $(n(k))_{k \geqslant 1}$ such that for all $k \geqslant 1$,

$$\left| \int_S f \, d\mathbb{P}_{n(k)} - \int_S f \, d\mathbb{P} \right| > \varepsilon.$$

The sequence $(\mathbb{P}_{n(k)})_{k \geqslant 1}$ cannot have a subsequence that converges weakly to $\mathbb{P}$. This contradiction completes the proof. ∎

In general, relative compactness of the sequence $(\mathbb{P}_n)_{n \geqslant 1}$ is obtained through the Prokhorov theorem, while uniqueness of the limiting measure is obtained through case-by-case considerations using tools such as characteristic functions or Laplace transforms. The Prokhorov theorem states that the relative compactness of the sequence $(\mathbb{P}_n)_{n \geqslant 1}$ is implied by its uniform tightness. A sequence $(\mathbb{P}_n)_{n \geqslant 1}$ of probability measures on $S$ is *uniformly tight* if, for any $\varepsilon > 0$, there exists a compact set $K \subseteq S$ such that, for every $n \geqslant 1$,

$$\mathbb{P}_n(K) \geqslant 1 - \varepsilon. \tag{A.57}$$

There are many different proofs of the Prokhorov theorem, some of which require heavy machinery from measure theory. We will present a relatively simple proof based on the Riesz representation theorem and the upshot of Exercise A.4 that every finite measure on a metric space is closed regular.

**Theorem A.20** (Prokhorov). *If the sequence $(\mathbb{P}_n)_{n \geqslant 1}$ of probability measures is uniformly tight, then it admits a weakly convergent subsequence.*

*Proof.* The proof proceeds in four steps. We start by extracting a subsequence $(n(k))_{k \geqslant 1}$ along which the integral $\int_S f \, d\mathbb{P}_{n(k)}$ converges for every $f \in C_b(S; \mathbb{R})$, we then use the Riesz representation theorem to identify a candidate limit measure $\mathbb{P}$ for the sequence $(\mathbb{P}_{n(k)})_{k \geqslant 1}$. Finally, we prove that $\mathbb{P}$ is a probability measure and that $(\mathbb{P}_{n(k)})_{k \geqslant 1}$ converges weakly to $\mathbb{P}$.

*Step 1: extracting the subsequence.* Since $(\mathbb{P}_n)_{n \geqslant 1}$ is uniformly tight, for every $r \geqslant 1$ there exists a compact set $K_r \subseteq S$ such that, for every $n \geqslant 1$,

$$\mathbb{P}_n(K_r) \geqslant 1 - \frac{1}{r}.$$

Observe that each of the spaces $C(K_r; \mathbb{R})$ is separable with respect to the uniform topology. Indeed, for each $r \geqslant 1$, the compactness of $K_r$ implies its separability, so, denoting by $(x_n)_{n \geqslant 1}$ a dense subset of $K_r$, and applying the Stone-Weierstrass theorem to the algebra of rational linear combinations of products of the continuous functions 1 and $d_n(x) \coloneqq d(x, x_n)$ establishes the separability of $C(K_r, \mathbb{R})$. It is therefore possible to find a countable set $C_r$ that is dense in $C(K_r; \mathbb{R})$. A diagonalization argument yields a subsequence $(n(k))_{k \geq 1}$ such that the functional $T : \bigcup_{r \geq 1} C_r \to \mathbb{R}$ defined by

$$T(f) := \lim_{k \to +\infty} \int_S f \, d\mathbb{P}_{n(k)} \tag{A.58}$$

is well-defined. By density of $C_r$ in $C(K_r; \mathbb{R})$, this functional is actually well-defined on $\bigcup_{r \geq 1} C(K_r; \mathbb{R})$. In fact, it is also well-defined on $C_b(S; \mathbb{R})$. Indeed, for any $f \in C_b(S; \mathbb{R})$ and $r \geq 1$, we have $f|_{K_r} \in C_b(K_r; S)$ and

$$\left| \int_S f \, d\mathbb{P}_{n(k)} - \int_{K_r} f|_{K_r} \, d\mathbb{P}_{n(k)} \right| \leq \|f\|_\infty \mathbb{P}_{n(k)}(K_r^c) \leq \frac{\|f\|_\infty}{r}.$$

*Step 2: identifying the candidate limit.* For each $r \geq 1$ introduce the compact set $\widetilde{K}_r := \bigcup_{i=1}^r K_i$. Since $T$ is a positive linear functional on $C(\widetilde{K}_r; \mathbb{R})$, the Riesz representation theorem gives a finite measure $\mu_r$ on $\widetilde{K}_r$ such that

$$T(f) = \int_{\widetilde{K}_r} f \, d\mu_r \tag{A.59}$$

for all $f \in C(\widetilde{K}_r; \mathbb{R})$. The measure $\mu_r$ on $\widetilde{K}_r$ induces the measure $\overline{\mu}_r(\cdot) := \mu_r(\cdot \cap \widetilde{K}_r)$ on $S$. We would now like to define the candidate measure $\mathbb{P} : \mathcal{B}(S) \to \mathbb{R}_{\geq 0}$ by

$$\mathbb{P}(A) := \lim_{r \to +\infty} \overline{\mu}_r(A). \tag{A.60}$$

To show that this limit is well-defined, fix $r_1 < r_2$ as well as $A \in \mathcal{B}(S)$. Given $\varepsilon > 0$ use the closed regularity (A.20) of finite measures established in Exercise A.4 to find a closed set $F \subseteq A \cap \widetilde{K}_{r_1}$ with

$$\overline{\mu}_{r_1}(A) \leq \overline{\mu}_{r_1}(F) + \varepsilon.$$

Similarly, use the open regularity (A.21) of finite measures to find an open set $U \subseteq \widetilde{K}_{r_2}$ with $A \cap \widetilde{K}_{r_2} \subseteq U$ and

$$\overline{\mu}_{r_2}(U) - \varepsilon \leq \overline{\mu}_{r_2}(A).$$

Using Urysohn's lemma find a function $f \in C(\widetilde{K}_{r_2}; \mathbb{R})$ with $\mathbf{1}_F \leq f \leq \mathbf{1}_U$. Since $F \subseteq \widetilde{K}_{r_1} \subseteq \widetilde{K}_{r_2}$ and $T$ is given by the limit (A.58),

$$\overline{\mu}_{r_1}(A) \leq \overline{\mu}_{r_1}(F) + \varepsilon \leq \int_{\widetilde{K}_{r_1}} f \, d\mu_{r_1} + \varepsilon = T(f|_{\widetilde{K}_{r_1}}) \leq T(f|_{\widetilde{K}_{r_2}}) = \int_{\widetilde{K}_{r_2}} f \, d\mu_{r_2} + \varepsilon.$$

Together with the fact that $f \leq \mathbf{1}_U$, this implies that

$$\overline{\mu}_{r_1}(A) \leq \overline{\mu}_{r_2}(U) + \varepsilon \leq \overline{\mu}_{r_2}(A) + 2\varepsilon.$$

Letting $\varepsilon$ tend to zero shows that $\overline{\mu}_{r_1}(A) \leqslant \overline{\mu}_{r_2}(A)$, and proves that the set function $\mathbb{P} : \mathcal{B}(S) \to \mathbb{R}_{\geqslant 0}$ in (A.60) is well-defined as a monotone limit.

*Step 3: showing $\mathbb{P}$ is a probability measure.* To show that $\mathbb{P}$ is a measure, fix a collection $(A_i)_{i \geqslant 1}$ of pairwise disjoint sets in $\mathcal{B}(S)$, and let $A = \bigcup_{i \geqslant 1} A_i$. On the one hand, for any $n \geqslant 1$,

$$\mathbb{P}(A) \geqslant \lim_{r \to +\infty} \overline{\mu}_r\left(\bigcup_{i=1}^n A_i\right) = \lim_{r \to +\infty} \sum_{i=1}^n \overline{\mu}_r(A_i) = \sum_{i=1}^n \mathbb{P}(A_i),$$

so letting $n$ tend to infinity gives the lower bound $\mathbb{P}(A) \geqslant \sum_{i=1}^{+\infty} \mathbb{P}(A_i)$. On the other hand, given $\varepsilon > 0$ there is $r$ large enough so

$$\mathbb{P}(A) \leqslant \overline{\mu}_r(A) + \varepsilon = \sum_{i=1}^{+\infty} \overline{\mu}_r(A_i) + \varepsilon \leqslant \sum_{i=1}^{+\infty} \mathbb{P}(A_i) + \varepsilon,$$

where the final inequality uses that $\mathbb{P}(A_i)$ is the monotonically increasing limit of $\overline{\mu}_r(A_i)$. Letting $\varepsilon$ tend to zero establishes the matching upper bound and shows that $\mathbb{P}$ is a measure. To see that it is a probability measure, observe that for every $r \geqslant 1$,

$$1 = T(1) \geqslant \lim_{r \to +\infty} \overline{\mu}_r(S) = \mathbb{P}(S) \geqslant \mu_r(\widetilde{K}_r) = \lim_{k \to +\infty} \mathbb{P}_{n(k)}(\widetilde{K}_r) \geqslant 1 - \frac{1}{r}.$$

Letting $r$ tend to infinity shows that $\mathbb{P}(S) = 1$ as required.

*Step 4: establishing weak convergence.* To conclude that the sequence $(\mathbb{P}_{n(k)})_{k \geqslant 1}$ of probability measures converges weakly to $\mathbb{P}$, fix a closed set $F \subseteq S$ as well as $r \geqslant 1$. Letting $k$ tend to infinity in the bound

$$\mathbb{P}_{n(k)}(F) \leqslant \mathbb{P}_{n(k)}(F \cap \widetilde{K}_r) + \frac{1}{r}$$

reveals that

$$\limsup_{k \to +\infty} \mathbb{P}_{n(k)}(F) \leqslant \limsup_{k \to +\infty} \mathbb{P}_{n(k)}(F \cap \widetilde{K}_r) + \frac{1}{r} \leqslant \overline{\mu}_r(F) + \frac{1}{r} \leqslant \mathbb{P}(F) + \frac{1}{r}.$$

The second inequality combines the fact that $(\mathbb{P}_{n(k)})_{k \geqslant 1}$ converges weakly to $\overline{\mu}_r$ on $\widetilde{K}_r$ by (A.58) and (A.59) with Remark A.18. Letting $r$ tend to infinity and invoking the Portmanteau theorem completes the proof. ∎

It turns out that the converse of the Prokhorov theorem is also true for any complete and separable metric space.

**Theorem A.21.** *Suppose $S$ is a complete and separable metric space. If the sequence $(\mathbb{P}_n)_{n \geqslant 1}$ of probability measures on $S$ converges weakly, then it is uniformly tight.*

*Proof.* Suppose for the sake of contradiction that $(\mathbb{P}_n)_{n \geq 1}$ converges weakly to some probability measure $\mathbb{P}$ but that it is not uniformly tight. This means that there exist $\varepsilon > 0$ and $r > 0$ such that for any finite family of open balls $B_1, \ldots, B_m$ of radius $r > 0$ there is some $n \geq 1$ with $\mathbb{P}_n\bigl(\bigcup_{i=1}^m B_i\bigr) < 1 - \varepsilon$. Indeed, if this were not the case, then given $\varepsilon > 0$ and $k > 0$ it would be possible to find a finite family of open balls $B_1^{(k)}, \ldots, B_{n(k)}^{(k)}$ of radius $k^{-1}$ with

$$\mathbb{P}_n\left(\bigcup_{i=1}^{n(k)} B_i^{(k)}\right) \geq 1 - \varepsilon 2^{-k}$$

for every $n \geq 1$. The set

$$K := \bigcap_{k=1}^{\infty} \bigcup_{i=1}^{n(k)} \overline{B}_i^{(k)}$$

would then satisfy $\mathbb{P}_n(K) \geq 1 - \varepsilon$ for every $n \geq 1$. Since $S$ is complete, the closed and totally bounded set $K$ is compact. This would contradict the assumption that $(\mathbb{P}_n)_{n \geq 1}$ is not uniformly tight. With this $r > 0$ at hand, use the separability of $S$ to find a collection $(B_i)_{i \geq 1}$ of open balls of radius $r > 0$ with $S = \bigcup_{i=1}^{\infty} B_i$. For each $k \geq 1$, let $A_k := \bigcup_{i=1}^k B_i$, and find $n(k) \geq 1$ with $\mathbb{P}_{n(k)}(A_k) < 1 - \varepsilon$. Since $A_m$ is open for each $m \geq 1$ and $A_m \subseteq A_k$ for every $k$ large enough, the Portmanteau theorem implies that

$$\mathbb{P}(A_m) \leq \liminf_{k \to +\infty} \mathbb{P}_{n(k)}(A_m) \leq \liminf_{k \to +\infty} \mathbb{P}_{n(k)}(A_k) < 1 - \varepsilon.$$

By the continuity of measure established in Exercise A.1 this contradicts the fact that $S = \bigcup_{m=1}^{\infty} A_m$ and completes the proof. ∎

In the context of our two-step strategy to establish weak convergence, we now turn our attention to the uniqueness of the limiting measure in the setting of Euclidean space $\mathbb{R}^d$. We will discuss two closely related transforms that entirely determine a probability measure, the characteristic function and the Laplace transform. The key difference between these two transforms is that the former is well-defined for any probability measure while the latter is slightly simpler but only well-defined for a smaller class of probability measures; here we will consider probability measures supported on the positive quadrant $\mathbb{R}^d_{\geq 0}$. The *characteristic function* of a probability measure $\mathbb{P}$ on $\mathbb{R}^d$ is the complex-valued function $\varphi : \mathbb{R}^d \to \mathbb{C}$ defined by

$$\varphi_{\mathbb{P}}(t) := \int_{\mathbb{R}^d} e^{it \cdot x} \, d\mathbb{P}(x). \tag{A.61}$$

In the case when $\mathbb{P}$ admits a density with respect to Lebesgue measure, the characteristic function is the Fourier transform of the density. The *Laplace transform* of a probability measure $\mathbb{P}$ on $\mathbb{R}^d_{\geq 0}$ is the real-valued function $\mathcal{L}_{\mathbb{P}} : \mathbb{R}^d_{\geq 0} \to \mathbb{R}$ defined by

$$\mathcal{L}_{\mathbb{P}}(\lambda) := \int_{\mathbb{R}^d_{\geq 0}} e^{-\lambda \cdot x} \, d\mathbb{P}(x). \tag{A.62}$$

**Theorem A.22.** *If two probability measures $\mathbb{P}$ and $\mathbb{Q}$ on $\mathbb{R}^d$ have the same characteristic function, $\varphi_{\mathbb{P}} = \varphi_{\mathbb{Q}}$, then they must be equal, $\mathbb{P} = \mathbb{Q}$.*

*Proof.* By a straightforward approximation argument, it suffices to show that for every continuous function of compact support $f \in C_c(\mathbb{R}^d; \mathbb{R})$, we have

$$\int_{\mathbb{R}^d} f \, d\mathbb{P} = \int_{\mathbb{R}^d} f \, d\mathbb{Q}. \tag{A.63}$$

Fix a continuous function of compact support $f \in C_c(\mathbb{R}^d; \mathbb{R})$ as well as $\varepsilon > 0$. Using Exercise A.12 on the tightness of probability measures on complete and separable metric spaces, find $K$ large enough so that $f$ vanishes outside the cube $[-K/2, K/2]^d$, and so that $\mathbb{P}((-K,K)^d) \geq 1-\varepsilon$ and $\mathbb{Q}((-K,K)^d) \geq 1-\varepsilon$. For each $m \in \mathbb{Z}^d$ denote by $g_m : [-K, K]^d \to \mathbb{C}$ the function

$$g_m(x) := \exp\left(\frac{i\pi}{2K} m \cdot x\right), \tag{A.64}$$

and write $\mathcal{A}$ for the algebra of finite linear combinations of the $(g_m)_{m \in \mathbb{Z}^d}$. It is readily verified that $\mathcal{A}$ is an algebra in $C(S; \mathbb{C})$ that separates points in the compact set $[-K, K]^d$. Since $\mathcal{A}$ is closed under complex conjugation, the complex version of the Stone-Weierstrass theorem gives $g \in \mathcal{A}$ such that, for all $x \in [-K, K]^d$,

$$|f(x) - g(x)| \leq \frac{\varepsilon}{2}.$$

Remembering that $f$ is supported on $[-K/2, K/2]^d$, up to multiplying $g$ by a nonnegative function which is equal to one on $[-K/2, K/2]^d$, vanishes outside $[-K, K]^d$ and is bounded by one, assume without loss of generality that $g$ is supported on $[-K, K]^d$. This ensures that $f(x) = g(x)$ for $x \in \mathbb{R}^d \setminus [-K, K]^d$. Combining this with the fact that $\int_{\mathbb{R}^d} g \, d\mathbb{P} = \int_{\mathbb{R}^d} g \, d\mathbb{Q}$ by the equality of characteristic functions reveals that

$$\left| \int_{\mathbb{R}^d} f \, d\mathbb{P} - \int_{\mathbb{R}^d} f \, d\mathbb{Q} \right| \leq \int_{\mathbb{R}^d} |f - g| \, d\mathbb{P} + \int_{\mathbb{R}^d} |f - g| \, d\mathbb{Q} \leq \varepsilon.$$

Letting $\varepsilon$ tend to zero establishes (A.63) and completes the proof. ∎

**Theorem A.23.** *If two probability measures $\mathbb{P}$ and $\mathbb{Q}$ on $\mathbb{R}^d_{\geq 0}$ have the same Laplace transform, $\mathcal{L}_{\mathbb{P}} = \mathcal{L}_{\mathbb{Q}}$, then they must be equal, $\mathbb{P} = \mathbb{Q}$.*

*Proof.* Up to replacing the function (A.64) by the function $g_m : [-K, K]^d \to \mathbb{R}$ defined by $g_m(x) = e^{m \cdot x}$, the proof is identical to that of Theorem A.22. ∎

We now provide an example that shows how the ideas developed in this section can be used to establish the weak convergence of a sequence of probability measures on the positive half-space by looking only at their Laplace transforms.

**Proposition A.24.** *A sequence $(\mathbb{P}_n)_{n\geq 1}$ of probability measures on $\mathbb{R}^d_{\geq 0}$ converges weakly to a probability measure $\mathbb{P}$ on $\mathbb{R}^d_{\geq 0}$ if and only if, for every $\lambda \in \mathbb{R}^d_{\geq 0}$,*

$$\lim_{n\to+\infty} \int_{\mathbb{R}^d_{\geq 0}} e^{-\lambda\cdot x}\, d\mathbb{P}_n(x) = \int_{\mathbb{R}^d_{\geq 0}} e^{-\lambda\cdot x}\, d\mathbb{P}(x). \tag{A.65}$$

*Proof.* The direct implication is immediate. To show the converse implication, we proceed in two steps. First we argue that $(\mathbb{P}_n)_{n\geq 1}$ is uniformly tight using the Prokhorov theorem, and then we prove that $\mathbb{P}$ is its only subsequential limit using Theorem A.23. The result then follows from Lemma A.19.

*Step 1: $(\mathbb{P}_n)_{n\geq 1}$ is uniformly tight.* Fix $\varepsilon > 0$ and $M > 0$. Denote by $K_M := [0,M]^d$ the box of side-length $M$ in $\mathbb{R}^d_{\geq 0}$, and for a choice of $\lambda^* > 0$ to be decided shortly, consider the constant vector

$$\lambda := (\lambda^*,\ldots,\lambda^*) \in \mathbb{R}^d_{\geq 0}.$$

Since $\mathcal{L}_\mathbb{P}$ is continuous and $\mathcal{L}_\mathbb{P}(0) = 1$, we can fix $\lambda^* > 0$ small enough such that for every $n$ sufficiently large,

$$\int_{\mathbb{R}^d_{\geq 0}} e^{-\lambda\cdot x}\, d\mathbb{P}_n(x) = \mathcal{L}_{\mathbb{P}_n}(\lambda) \geq 1 - \frac{\varepsilon}{2}.$$

It follows that for every $n$ sufficiently large,

$$1 - \frac{\varepsilon}{2} \leq \mathbb{P}_n(K_M) + \int_{\mathbb{R}^d_{\geq 0}\setminus K_M} e^{-\lambda\cdot x}\, d\mathbb{P}_n(x) \leq \mathbb{P}_n(K_M) + e^{-M\lambda^*}\mathbb{P}_n(K_M^c).$$

Rearranging reveals that

$$\mathbb{P}_n(K_M^c) \leq \frac{\varepsilon}{2(1-e^{-M\lambda^*})},$$

so choosing $M = \frac{\log(2)}{\lambda^*}$ ensures that $\mathbb{P}_n(K_M^c) \leq \varepsilon$. This shows that $(\mathbb{P}_n)_{n\geq 1}$ is uniformly tight, so any subsequence of $(\mathbb{P}_n)_{n\geq 1}$ admits a weakly convergent subsequence by the Prokhorov theorem.

*Step 2: $(\mathbb{P}_n)_{n\geq 1}$ admits a unique subsequential limit.* Let $(\mathbb{P}_{n(k)})_{k\geq 1}$ be a subsequence of $(\mathbb{P}_n)_{n\geq 1}$ converging weakly to some probability measure $\mathbb{Q}$ on $\mathbb{R}^d_{\geq 0}$. On the one hand, the definition of weak convergence implies that for every $\lambda \in \mathbb{R}^d_{\geq 0}$,

$$\lim_{k\to+\infty} \int_S e^{-\lambda\cdot x}\, d\mathbb{P}_{n(k)} = \int_S e^{-\lambda\cdot x}\, d\mathbb{Q}.$$

On the other hand, by assumption, for every $\lambda \in \mathbb{R}^d_{\geq 0}$,

$$\lim_{k\to+\infty} \int_S e^{-\lambda\cdot x}\, d\mathbb{P}_{n(k)} = \int_S e^{-\lambda\cdot x}\, d\mathbb{P}.$$

Combining the uniqueness of weak limits with Theorem A.23 implies that $\mathbb{P} = \mathbb{Q}$. This completes the proof. ∎

**Exercise A.12.** Suppose $S$ is a complete and separable metric space. Show that any probability measure $\mathbb{P}$ on $S$ is tight in the sense that for every $\varepsilon > 0$ there exists a compact set $K \subseteq S$ with $\mathbb{P}(K) \geq 1 - \varepsilon$.

# Appendix S
# Solutions to exercises

## S.1 Introduction to statistical mechanics

**Exercise 1.1.** We treat each question separately.

(i) Bounding the sum of the non-negative terms from below by a single term in the sum, and bounding each term in the sum by its maximum shows that

$$\log \max_{1 \leqslant k \leqslant K} \exp(Na_k) \leqslant \log \sum_{k=1}^{K} \exp(Na_k) \leqslant \log(K) + \log \max_{1 \leqslant k \leqslant K} \exp(Na_k).$$

It follows that

$$\frac{1}{N}\left|\log \sum_{k=1}^{K} \exp(Na_k) - \log \max_{1 \leqslant k \leqslant K} \exp(Na_k)\right| \leqslant \frac{\log(K)}{N}. \quad \text{(S.1)}$$

Letting $N \to +\infty$ establishes (1.15).

(ii) If $K = \exp(o(N))$, then (S.1) ensures that (1.15) remains valid. On the other hand, if $\log(K_N)$ grows faster than $N$, then there exists a constant $c > 0$ with $\log(K_N) \geqslant cN$ for all $N$ large enough. The constant sequence with $a_k = 1$ for every $k \geqslant 1$ is such that

$$\frac{1}{N}\left|\log \sum_{k=1}^{K_N} \exp(Na_k) - \log \max_{1 \leqslant k \leqslant K} \exp(Na_k)\right| = \frac{1}{N}\left|\log(K_N) + N - N\right| = \frac{\log(K_N)}{N} \geqslant c$$

so (1.15) cannot hold. This means that (1.15) holds for arbitrary choices of the sequence $(a_k)_{k \geqslant 1}$ if and only if $K = \exp(o(N))$.

**Exercise 1.2.** Consider the two-by-two matrix $M = (M_{i,j})_{i,j \in \{\pm 1\}}$ defined by

$$M_{ij} := e^{\beta ij + hj},$$

and denote by $\lambda_1 \leq \lambda_2$ its eigenvalues,

$$\lambda_1 = e^\beta \cosh(h) - \left(e^{2\beta} \cosh^2(h) - 2\sinh(2\beta)\right)^{\frac{1}{2}},$$
$$\lambda_2 = e^\beta \cosh(h) + \left(e^{2\beta} \cosh^2(h) - 2\sinh(2\beta)\right)^{\frac{1}{2}}.$$

Notice that $e^{2\beta}\cosh^2(h) \geq e^{2\beta} \geq 2\sinh(2\beta)$, so the square root is well-defined. Given a spin configuration $\sigma \in \{\pm 1\}^N$, let $\sigma_{N+1} = \sigma_1$ in such a way that

$$F_N(\beta, h) = \frac{1}{N} \log \sum_{\sigma \in \{\pm 1\}^N} \prod_{i=1}^N M_{\sigma_i \sigma_{i+1}} = \frac{1}{N} \log \sum_{\sigma_1 \in \Sigma_1} M^N_{\sigma_1 \sigma_1} = \frac{1}{N} \log \operatorname{tr}(M^N)$$
$$= \frac{1}{N} \log\left(\lambda_1^N + \lambda_2^N\right).$$

Since $\lambda_2 > 0$ and $\lambda_2 > \lambda_1$, the limit free energy is given by

$$f(\beta, h) := \lim_{N \to +\infty} F_N(\beta, h) = \log \lambda_2 + \lim_{N \to +\infty} \frac{1}{N} \log\left(1 + \left(\frac{\lambda_1}{\lambda_2}\right)^N\right) = \log \lambda_2.$$

The quantity $\partial_h F_N(\beta, h)$ is the mean magnetization (1.23), and is therefore bounded uniformly by one. It follows by smoothness of $f$ in $h$, Propositions 2.11 and 2.15, and Exercise 2.6 that for every $h \in \mathbb{R}$, the limit (1.25) of the mean magnetization is given by

$$m(\beta, h) := \lim_{N \to +\infty} \partial_h F_N(\beta, h) = \partial_h f(\beta, h) = \frac{e^\beta \sinh(h)}{(e^{2\beta} \cosh^2(h) - 2\sinh(2\beta))^{1/2}},$$

and is therefore a continuous function of $h$.

**Exercise 1.3.** We treat each question separately.

(i) The Hamiltonian (1.28) may be written as

$$H_N(\sigma) = \sum_{\substack{i,j \in B_{N+1} \\ i \sim j}} \mathbf{1}_{\{\sigma_i \sigma_j = 1\}} - \sum_{\substack{i,j \in B_{N+1} \\ i \sim j}} \mathbf{1}_{\{\sigma_i \sigma_j = -1\}}$$
$$= \left|\{\{i,j\} \in B_{N+1}^2 \mid i \sim j\}\right| - 2 \sum_{\substack{i,j \in B_{N+1} \\ i \sim j}} \mathbf{1}_{\{\sigma_i \sigma_j = -1\}}.$$

The final sum counts the number of pairs of nearest neighbours $(i,j) \in B_{N+1}^2$ with $\sigma_i \sigma_j = -1$. This corresponds to the number of edges $|\Gamma(\sigma)|$ in $\Gamma(\sigma)$. It follows that $H_N(\sigma) + 2|\Gamma(\sigma)|$ gives the total number of nearest neighbour pairs in the lattice $B_{N+1}^2$, so it depends only on the geometry of the lattice and not the specific configuration $\sigma \in \{\pm 1\}^{B_N}$.

(ii) Given a configuration $\sigma \in \{\pm 1\}^{B_N}$ with $\sigma_0 = -1$, let

$$j(\sigma) := \min\{j \in \{1,\ldots,N+1\} \mid \sigma_{(0,j)} = 1\}$$

be the coordinate of the first spin to right of $\sigma_0$ that is unaligned with $\sigma_0$. Notice that such a spin must exist as all spins in $B_{N+1} \setminus B_N$ are unaligned with $\sigma_0$. Since the contours in $\Gamma(\sigma)$ delimit regions in which the sign of the spins is constant, there must be a contour which separates the origin from the point $(0, j(\sigma))$ and therefore surrounds the origin.

(iii) Given a spin configuration $\sigma \in \{\pm 1\}^{B_N}$, consider a partition $\mathcal{P}(\sigma)$ of the set of contours $\Gamma(\sigma)$,

$$\Gamma(\sigma) = \bigcup_{\eta \in \mathcal{P}(\sigma)} \eta.$$

Since $H_N(\sigma) + 2|\Gamma(\sigma)|$ does not depend on the configuration $\sigma$ by (i), we have

$$\langle \mathbf{1}_{\{\gamma \subseteq \Gamma(\sigma)\}} \rangle = \frac{\sum_{\sigma \in \{\pm 1\}^{B_N}} \mathbf{1}_{\{\gamma \subseteq \Gamma(\sigma)\}} \exp(-2\beta|\Gamma(\sigma)|)}{\sum_{\sigma \in \{\pm 1\}^{B_N}} \exp(-2\beta|\Gamma(\sigma)|)}$$

$$= \exp(-2\beta|\gamma|) \cdot \frac{\sum_{\sigma \in \{\pm 1\}^{B_N}} \mathbf{1}_{\{\gamma \subseteq \Gamma(\sigma)\}} \prod_{\eta \in \mathcal{P}(\sigma) \setminus \gamma} \exp(-2\beta|\eta|)}{\sum_{\sigma \in \{\pm 1\}^{B_N}} \prod_{\eta \in \mathcal{P}(\sigma)} \exp(-2\beta|\eta|)}.$$

It therefore suffices to prove that

$$\sum_{\sigma \in \{\pm 1\}^{B_N}} \mathbf{1}_{\{\gamma \subseteq \Gamma(\sigma)\}} \prod_{\eta \in \mathcal{P}(\sigma) \setminus \gamma} \exp(-2\beta|\eta|) \leq \sum_{\sigma \in \{\pm 1\}^{B_N}} \prod_{\eta \in \mathcal{P}(\sigma)} \exp(-2\beta|\eta|). \quad \text{(S.2)}$$

Denote by $\Sigma_\gamma$ the set of configurations $\sigma \in \{\pm 1\}^{B_N}$ with $\gamma \subseteq \Gamma(\sigma)$, and define the flipping map $\mathcal{F} : \Sigma_\gamma \to \{\pm 1\}^{B_N}$ by

$$\mathcal{F}(\sigma)_i := \begin{cases} -\sigma_i & \text{if } i \in \text{int}(\gamma), \\ \sigma_i & \text{otherwise.} \end{cases}$$

Observe that $\mathcal{F}$ is injective and satisfies $\Gamma(\mathcal{F}(\sigma)) = \mathcal{P}(\sigma) \setminus \gamma$. It follows that the left side of (S.2) is given by

$$\sum_{\sigma \in \{\pm 1\}^{B_N}} \mathbf{1}_{\{\gamma \in \Gamma(\sigma)\}} \prod_{\eta \in \mathcal{P}(\mathcal{F}(\sigma))} \exp(-2\beta|\eta|)$$

which is in turn equal to

$$\sum_{\sigma \in \{\pm 1\}^{B_N}} \mathbf{1}_{\{\sigma \subseteq \mathcal{F}(\Sigma_\gamma)\}} \prod_{\eta \in \mathcal{P}(\sigma)} \exp(-2\beta|\eta|) \leq \sum_{\sigma \in \{\pm 1\}^{B_N}} \prod_{\eta \in \mathcal{P}(\sigma)} \exp(-2\beta|\eta|).$$

This establishes (S.2).

(iv) To bound $\langle \sigma_0 \rangle$ from below away from zero we will use that

$$\langle \sigma_0 \rangle = \langle \mathbf{1}_{\{\sigma_0=1\}} \rangle - \langle \mathbf{1}_{\{\sigma_0=-1\}} \rangle = 1 - 2\langle \mathbf{1}_{\{\sigma_0=-1\}} \rangle, \tag{S.3}$$

and instead bound $\langle \mathbf{1}_{\{\sigma_0=-1\}} \rangle$ from above away from 1/2. Denote by $\Gamma_0$ the set of all contours surrounding the origin, and write $\Gamma_0^k$ for the set of contours in $\Gamma_0$ of length $k$. Combining (ii) with the union bound and (iii) reveals that

$$\langle \mathbf{1}_{\{\sigma_0=-1\}} \rangle \leq \sum_{\gamma \in \Gamma_0} \langle \mathbf{1}_{\{\gamma \subseteq \Gamma(\sigma)\}} \rangle \leq \sum_{\gamma \in \Gamma_0} \exp(-2\beta|\gamma|) = \sum_{k=1}^{+\infty} |\Gamma_0^k| \exp(-2\beta k).$$

To bound the size of $\Gamma_0^k$, observe that any path enclosing the origin must have length at least 4. Moreover, it must cross the point $(0, j)$ for some $j \in \{1, \ldots, k/2\}$, and at each of its $k$ nodes it can go in at most 4 directions. This means that $k$ must be larger than 4 for $\Gamma_0^k$ not to be empty, and that $\Gamma_0^k$ can contain at most $k4^k$ edges. It follows that

$$\langle \mathbf{1}_{\{\sigma_0=-1\}} \rangle \leq \sum_{k=4}^{+\infty} k\left(4e^{-2\beta}\right)^k.$$

Taking $\beta > \frac{3}{2}\log 2$ and invoking the dominated convergence theorem shows that the right side of this expression tends to zero as $\beta \to +\infty$. In particular, for any $\delta > 0$ there exists $\beta$ large enough with $\limsup_{N \to +\infty} \langle \mathbf{1}_{\{\sigma_0=-1\}} \rangle \leq \frac{\delta}{2}$. Taking $\delta < 1$ and substituting this into (S.3) gives $\beta$ large enough with $\liminf_{N \to +\infty} \langle \sigma_0 \rangle > 0$ as required.

(v) The symmetry between spin sites when $h = 0$ implies that the mean magnetization (1.23) simplifies to

$$m_N(\beta, 0) = \left\langle \frac{1}{|B_2|} \sum_{i \in B_2} \sigma_i \right\rangle = \langle \sigma_0 \rangle.$$

It follows by (iv) that for $\beta$ large enough,

$$\liminf_{N \to +\infty} m_N(\beta, 0) > 0$$

which means that the two-dimensional Ising model with no magnetic field and with $+1$ boundary condition has a non-zero asymptotic mean magnetization. This suggests that it has a "memory" of the slight tilt it has been exposed to via the boundary condition, and therefore that it exhibits ferromagnetic behaviour. Although we will not pursue this further, one can show by subadditivity arguments that the limit free energy associated with this model is well-defined and does not depend on the boundary condition we choose; and use Proposition 2.15 to assert that any subsequential limit of $m_N(\beta, 0)$ must belong to the subdifferential (in $h$) of this limit free energy. By symmetry, we conclude that for $\beta > 0$ large enough, the limit free energy of the two-dimensional Ising model is not differentiable in $h$ at $h = 0$.

## S.2   Convex analysis and large deviation principles

**Exercise 2.1.** Fix $\alpha \in (0,1)$ and $x, y \in C_2 - C_1$. Let $x_1, y_1 \in C_1$ and $x_2, y_2 \in C_2$ be such that $x = x_2 - x_1$ and $y = y_2 - y_1$. Observe that

$$\alpha x + (1-\alpha)y = (\alpha x_2 + (1-\alpha)y_2) - (\alpha x_1 + (1-\alpha)y_1) \in C_2 - C_1. \qquad (S.4)$$

This establishes the convexity of $C_2 - C_1$.

**Exercise 2.2.** It is clear that $\mathrm{int}(C) \subseteq \mathrm{int}(\overline{C})$. To prove the converse containment, fix $z \in \mathrm{int}(\overline{C})$ and $x \in \mathrm{int}(C)$. Let $\gamma > 0$ be small enough so that $y := z + \gamma(z - x) \in \overline{C}$, and observe that

$$z = (1-\alpha)x + \alpha y$$

for $\alpha = \frac{1}{\gamma + 1} \in (0, 1)$. If $\varepsilon > 0$ is small enough so $B_\varepsilon(x) \subseteq C$, then $B_{(1-\alpha)\varepsilon}(z) \subseteq C$. This shows that $z \in \mathrm{int}(C)$, and therefore $\mathrm{int}(\overline{C}) = \mathrm{int}(C)$. The equality of the boundaries follows from the fact that $\partial A = \overline{A} \setminus \mathrm{int}(A)$ for every set $A \subseteq \mathbb{R}^d$.

**Exercise 2.3.** Since $\mathrm{int}(C) \subseteq C$ and $C$ is closed, we have $\overline{\mathrm{int}(C)} \subseteq C$. To prove the converse inclusion, fix $x \in C$ and $x' \in \mathrm{int}(C)$. Let $\varepsilon > 0$ be small enough so $B_\varepsilon(x') \subseteq C$. By convexity of $C$, for all $t \in [0,1]$ and $y \in B_\varepsilon(x')$, we have $(1-t)x + ty \in C$. This means that, for all $t \in [0,1]$,

$$B_{t\varepsilon}((1-t)x + tx') \subseteq C,$$

and thus $(1-t)x + tx' \in \mathrm{int}(C)$. It follows that

$$x = \lim_{t \searrow 0} ((1-t)x + tx') \in \overline{\mathrm{int}(C)}.$$

This completes the proof.

**Exercise 2.4.** We say that points $x_1, \ldots, x_k \in \mathbb{R}^d$ are *affinely independent* if

$$\sum_{i=1}^{k} \alpha_i x_i = 0 \quad \text{with } \alpha_i \in \mathbb{R} \text{ and } \sum_{i=1}^{k} \alpha_i = 0$$

implies that $\alpha_1 = \cdots = \alpha_k = 0$. This is equivalent to the linear independence of the vectors $x_2 - x_1, \ldots, x_k - x_1$. In particular, if $x_1, \ldots, x_k$ are affinely independent in $\mathbb{R}^d$, then $k \leq d + 1$. It therefore suffices to prove that any point $x \in \mathrm{conv}(A)$ can be written as a convex combination of affinely independent points in $A$. Since $x \in \mathrm{conv}(A)$, it may be written as

$$x = \sum_{i=1}^{k} \lambda_i x_i \quad \text{for some } x_i \in A, \lambda_i \geq 0 \text{ with } \sum_{i=1}^{k} \lambda_i = 1 \text{ and } k \geq 1. \qquad (S.5)$$

We take $k$ minimal among all such representations. Suppose for the sake of contradiction that $x_1, \ldots, x_k$ are not affinely independent. This means that there exist $\alpha_1, \ldots, \alpha_k \in \mathbb{R}$, not all zero, with

$$\sum_{i=1}^{k} \alpha_i x_i = 0 \quad \text{and} \quad \sum_{i=1}^{k} \alpha_i = 0.$$

Denote by $I$ the set of indices $i \in \{1, \ldots, k\}$ such that $\alpha_i > 0$. Since the set $I$ is not empty, we can choose $m \in I$ to be such that

$$\frac{\lambda_m}{\alpha_m} = \inf_{i \in I} \frac{\lambda_i}{\alpha_i}.$$

In the representation

$$x = \sum_{i=1}^{k} \left( \lambda_i - \frac{\lambda_m}{\alpha_m} \alpha_i \right) x_i,$$

the coefficients sum to one and are all non-negative by choice of $m$. Since the $m$'th coefficient is zero, this contradicts the minimality of $k$ among all representations of the form (S.5). This shows that $x_1, \ldots, x_k$ are affinely independent, and completes the proof.

**Exercise 2.5.** Let $(x^k)_{k \geq 1} \subseteq \text{conv}(A)$ be a sequence in the convex hull of $A$. By Exercise 2.4, for each $k \geq 1$, there exist $x_1^k, \ldots, x_{d+1}^k \in A$ and $\alpha_1^k, \ldots, \alpha_{d+1}^k \in [0, 1]$ with

$$x^k = \sum_{i=1}^{d+1} \alpha_i^k x_i^k \quad \text{and} \quad \sum_{i=1}^{d+1} \alpha_i^k = 1.$$

Since $A$ and $[0,1]$ are compact, there are subsequences of $(x_i^k)_{k \in \mathcal{K}}$ and $(\alpha_i^k)_{k \in \mathcal{K}}$ with the same index set $\mathcal{K}$ that converge to some $x_i$ and $\alpha_i$, respectively, for every $i \in \{1, \ldots, d+1\}$. Since $A$ is compact, each $x_i$ belongs to $A$. It follows that the subsequence $(x^k)_{k \in \mathcal{K}}$ converges to the point

$$x = \sum_{i=1}^{d+1} \alpha_i x_i \in \text{conv}(A).$$

This completes the proof.

**Exercise 2.6.** Fix $\lambda, \mu \in \mathbb{R}^d$ and $\alpha \in (0,1)$. By Hölder's inequality with conjugate exponents $p = \alpha^{-1}$ and $q = (1-\alpha)^{-1}$, we have

$$\exp \psi(\alpha \lambda + (1-\alpha)\mu) = \mathbb{E} \exp(\alpha \lambda \cdot X) \exp((1-\alpha)\mu \cdot X)$$
$$\leq (\mathbb{E} \exp(\lambda \cdot X))^\alpha (\mathbb{E} \exp(\mu \cdot X))^{1-\alpha}$$
$$= (\exp \psi(\lambda))^\alpha (\exp \psi(\mu))^{1-\alpha}.$$

Taking logarithms completes the proof.

**Exercise 2.7.** First suppose that $f_\alpha$ is convex for each $\alpha \in I$, and fix $x, y \in \mathbb{R}^d$ as well as $\lambda \in (0,1)$. For every $\alpha \in I$, we have

$$f_\alpha(\lambda x + (1-\lambda)y) \leq \lambda f_\alpha(x) + (1-\lambda)f_\alpha(y).$$

Taking the supremum over $\alpha \in I$ gives

$$\sup_{\alpha \in I} f_\alpha(\lambda x + (1-\lambda)y) \leq \lambda \sup_{\alpha \in I} f_\alpha(x) + (1-\lambda) \sup_{\alpha \in I} f_\alpha(y)$$

which establishes the convexity of $\sup_{\alpha \in I} f_\alpha$. Now, suppose instead that $f_\alpha$ is lower semi-continuous for every $\alpha \in I$. We observe that

$$\{(x,\lambda) \in \mathbb{R}^d \times \mathbb{R} \mid \sup_{\alpha \in I} f_\alpha(x) \leq \lambda\} = \bigcap_{\alpha \in I} \{(x,\lambda) \in \mathbb{R}^d \times \mathbb{R} \mid f_\alpha(x) \leq \lambda\}.$$

Since the intersection of closed sets is closed, this shows that the epigraph of $\sup_{\alpha \in I} f_\alpha$ is closed, and thus that $\sup_{\alpha \in I} f_\alpha$ is lower semi-continuous. We note in passing that the first part of the exercise can also be proved by reasoning about the convexity of the epigraphs.

**Exercise 2.8.** Let $x, x' \in \mathbb{R}^d$ and $\alpha \in (0,1)$. We want to verify that

$$g(\alpha x + (1-\alpha)x') \leq \alpha g(x) + (1-\alpha)g(x'). \tag{S.6}$$

Without loss of generality, we may assume that $g(x)$ and $g(x')$ are finite. For any given $\varepsilon > 0$, one can find $y, y' \in \mathbb{R}^k$ such that $f(x,y) \leq g(x) + \varepsilon$ and $f(x',y') \leq g(x') + \varepsilon$. By the definition of $g$ as an infimum and the convexity of $f$, we have

$$g(\alpha x + (1-\alpha)x') \leq f(\alpha x + (1-\alpha)x', \alpha y + (1-\alpha)y')$$
$$\leq \alpha f(x,y) + (1-\alpha)f(x',y').$$

Letting $\varepsilon$ tend to zero, we obtain the result.

**Exercise 2.9.** We treat each question separately.

(i) Since $\underline{f}$ is a supremum of lower semi-continuous functions, it is lower semi-continuous, by Exercise 2.7.

(ii) Denote $h(x) := \liminf_{y \to x} f(y)$, and let $g$ be lower semi-continuous with $g \geq f$. By the lower semi-continuity of $g$ and the inequality $g \leq f$, we have

$$g(x) \leq \liminf_{y \to x} g(y) \leq \liminf_{y \to x} f(x) = h(x).$$

Taking the supremum over $g$, we obtain that $\underline{f} \leq h$. To show the converse inequality, it suffices to verify that $h$ is lower semi-continuous, since we clearly have $h \leq f$.

Let $(x_n)_{n \geqslant 1}$ be a sequence converging to $x \in \mathbb{R}^d$, and $\varepsilon > 0$. For each $n \geqslant 1$, we can find a point $x'_n$ with $|x'_n - x_n| \leqslant 1/n$ and

$$f(x'_n) \leqslant h(x_n) + \varepsilon.$$

In particular, the sequence $(x'_n)_{n \geqslant 1}$ converges to $x$, and

$$h(x) = \liminf_{y \to x} f(y) \leqslant \liminf_{n \to +\infty} f(x'_n) \leqslant \liminf_{n \to +\infty} h(x_n) + \varepsilon.$$

Since $\varepsilon > 0$ is arbitrary, this shows that $h$ is lower semi-continuous, and thus completes the argument.

(iii) Let $x, y \in \mathbb{R}^d$ and $\alpha \in (0,1)$. We aim to show that

$$\underline{f}(\alpha x + (1-\alpha)x') \leqslant \alpha \underline{f}(x) + (1-\alpha)\underline{f}(x').$$

Without loss of generality, we may assume that $\underline{f}(x)$ and $\underline{f}(x')$ are finite. We fix $\varepsilon > 0$. For each $r > 0$, one can find $x_r, x'_r \in \mathbb{R}^d$ with $|x_r - x| + |x'_r - x'| \leqslant r$ such that

$$f(x_r) \leqslant \underline{f}(x) + \varepsilon \quad \text{and} \quad f(x'_r) \leqslant \underline{f}(x') + \varepsilon.$$

In particular, we have

$$|\alpha x_r + (1-\alpha)x'_r - \alpha x - (1-\alpha)x'_r| \leqslant r$$

and

$$f(\alpha x_r + (1-\alpha)x'_r) \leqslant \alpha \underline{f}(x) + (1-\alpha)\underline{f}(x') + \varepsilon.$$

Since this is valid for every $r > 0$, we deduce that

$$\underline{f}(\alpha x + (1-\alpha)x') \leqslant \alpha \underline{f}(x) + (1-\alpha)\underline{f}(x') + \varepsilon.$$

Since $\varepsilon > 0$ was arbitrary, the proof is complete.

**Exercise 2.10.** We treat each question separately.

(i) Since $f^*$ is the supremum of affine functions, this is immediate from Exercise 2.7.

(ii) The definition of the convex dual implies that $f(x) + f^*(\lambda) \geqslant \lambda \cdot x$ for every $x, \lambda \in \mathbb{R}^d$. It follows that for every $x, \lambda \in \mathbb{R}^d$,

$$x \cdot \lambda - f^*(\lambda) \leqslant f(x) + f^*(\lambda) - f^*(\lambda) = f(x).$$

Taking the supremum over $x \in \mathbb{R}^d$ shows that $f^{**} \leqslant f$.

(iii) If $f \leqslant g$ and $x, \lambda \in \mathbb{R}^d$, then
$$x \cdot \lambda - g(x) \leqslant x \cdot \lambda - f(x) \leqslant f^*(\lambda).$$

Taking the supremum over $x \in \mathbb{R}^d$ shows that $g^* \leqslant f^*$ and completes the proof.

**Exercise 2.11.** Fix $\lambda \in \mathbb{R}^d$, and observe that by the Cauchy-Schwarz inequality,
$$f^*(\lambda) = \frac{1}{2} \sup_{x \in \mathbb{R}^d} \left(2\lambda \cdot x - |x|^2\right) \leqslant \frac{1}{2} \sup_{x \in \mathbb{R}^d} \left(2|\lambda||x| - |x|^2\right).$$

The parabola $y \mapsto 2|\lambda|y - y^2$ has a unique maximum at $y = |\lambda|$, so $f^*(\lambda) \leqslant f(\lambda)$. The matching lower bound is obtained by choosing $x = \lambda$.

**Exercise 2.12.** Fix $x \in \mathbb{R}^d$ with $|x| > L$. By definition of the convex dual and Lipschitz continuity of $f$,
$$f^*(x) \geqslant \sup_{p \in \mathbb{R}^d} \left(p \cdot x - L|p|\right) - f(0) \geqslant \limsup_{\alpha \to +\infty} \alpha |x|\bigl(|x| - L\bigr) - f(0) = +\infty.$$

This completes the proof.

**Exercise 2.13.** We use the notation $\underline{f}^* := (\underline{f})^*$. Since $\underline{f} \leqslant f$, part (iii) of Exercise 2.10 yields that $f^* \geqslant \underline{f}^*$. To show the converse bound, we fix $\lambda \in \mathbb{R}^d$, $\varepsilon > 0$, and $x \in \mathbb{R}^d$ such that
$$f^*(\lambda) \leqslant \lambda \cdot x - f(x) - \varepsilon.$$

By Exercise 2.9, we can find $x' \in \mathbb{R}^d$ such that $|x' - x| \leqslant \varepsilon$ and $\underline{f}(x) \geqslant f(x) - \varepsilon$. Combining this with the previous display, we obtain that
$$\underline{f}^*(\lambda) \geqslant \lambda \cdot x' - \underline{f}'(x)$$
$$\geqslant \lambda \cdot x - f(x) - (|\lambda| + 1)\varepsilon$$
$$\geqslant f^*(\lambda) - (|\lambda| + 2)\varepsilon.$$

Since $\varepsilon > 0$ is arbitrary, this completes the proof that $f^* = \underline{f}^*$. By Exercise 2.9, the function $\underline{f}$ is convex and lower semi-continuous, so the Fenchel-Moreau theorem ensures that $\underline{f}^{**} = \underline{f}$. Conjugating the relation $f^* = \underline{f}^*$ therefore yields that $f^{**} = \underline{f}$, as announced.

**Exercise 2.14.** The definition of the polar $K^\circ$ ensures that for each fixed $x \in K$ and all $v \in K^\circ$, we have $x \cdot v \leqslant 0$. This means that $K \subseteq K^{\circ\circ}$. To prove the converse inclusion, suppose for the sake of contradiction that there exists $x \in K^{\circ\circ} \smallsetminus K$. Since $K$ is a closed convex set, the supporting hyperplane theorem (Theorem 2.2) gives a non-zero $v \in \mathbb{R}^d$ with
$$\sup\{v \cdot y \mid y \in K\} < v \cdot x.$$

If there were $y \in K$ with $v \cdot y > 0$, then we would have $v \cdot \lambda y < v \cdot x$ for all $\lambda > 0$ which would lead to a contradiction upon letting $\lambda$ tend to infinity. This means that $v \in K^\circ$. Observe also that $0 = \lim_{\lambda \searrow 0} \lambda y \in K$, where $y$ is any point in $K$, so $v \cdot x \geq 0$. Since $v \in K^\circ$, this contradicts the fact the $x \in K^{\circ\circ}$ and completes the proof.

**Exercise 2.15.** We treat each question separately.

(i) Fix $x \in \text{int}(C)$ and $n \in \mathbf{n}_C(x)$. For every $\varepsilon > 0$ small enough, we have $x + \varepsilon n \in C$, and therefore
$$\varepsilon |n|^2 = n \cdot (x + \varepsilon n - x) \leq 0.$$
This implies that $\mathbf{n}_C(x) \subseteq \{0\}$. The converse inclusion is trivial. To show that $\mathbf{T}_C(x) = \mathbb{R}^d$, fix $v \in \mathbb{R}^d$, and let $\varepsilon > 0$ be small enough so $x + \varepsilon v \in C$. Observe that
$$v = \varepsilon^{-1}(x + \varepsilon v - x) \in \mathbf{T}_C(x),$$
as required.

(ii) By Exercise 2.14, it suffices to show that $\mathbf{n}_C(x) = \mathbf{T}_C(x)^\circ$. On the one hand, if $v \in \mathbf{n}_C(x)$, then for all $x' \in C$ and $\lambda \geq 0$, we have $v \cdot \lambda(x' - x) \leq 0$. Since this inequality is preserved under limits, it follows that $v \cdot y \leq 0$ for all $y \in \mathbf{T}_C(x)$, and therefore that $\mathbf{n}_C(x) \subseteq \mathbf{T}_C(x)^\circ$. On the other hand, if $v \in \mathbf{T}_C(x)^\circ$, then for all $x' \in C$, we have $v \cdot (x' - x) \leq 0$ since $x' - x \in \mathbf{T}_C(x)$. This shows that $\mathbf{T}_C(x)^\circ \subseteq \mathbf{n}_C(x)$.

(iii) If $v \in \mathbf{T}_C(x)$, there exist sequences $(x'_i)_{i \geq 1} \subseteq C$ and $(\lambda_i)_{i \geq 1} \subseteq \mathbb{R}_{\geq 0}$ with $\lambda_i(x'_i - x) \to v$. Letting
$$t_i := \min(2^{-i}, \lambda_i^{-1}) \quad \text{and} \quad x_i := x + t_i \lambda_i (x'_i - x)$$
gives a sequence $(t_i)_{i \geq 1} \subseteq \mathbb{R}_{>0}$ decreasing to 0 and a sequence $(x_i)_{i \geq 1} \subseteq \mathbb{R}^d$ converging to $x$ with $t_i^{-1}(x_i - x) \to v$. Moreover, since $t_i \lambda_i \in [0, 1]$, we have that $x_i \in C$ for every $i \geq 1$. This completes the proof of the direct implication. The converse implication is immediate.

**Exercise 2.16.** It suffices to show that for every $x \in H \setminus C$ the infimum
$$d(x, C) := \inf\{\|x - y\|^2 \mid y \in C\}$$
is achieved. Indeed, all other arguments in the proof of Lemma 2.1 extend without change to the Hilbert space setting. Fix $x \in H \setminus C$, and let $(y_n)_{n \geq 1}$ be a minimizing sequence with
$$\lim_{n \to +\infty} \|x - y_n\|^2 = d(x, C).$$
We would like to show that $(y_n)_{n \geq 1}$ is Cauchy. Fix $n, m \geq 1$, and observe that
$$\|y_n - y_m\|^2 = \|y_n - x\|^2 + \|y_m - x\|^2 - 2(y_n - x) \cdot (y_m - x).$$

Similarly,
$$\|y_n + y_m - 2x\|^2 = \|y_n - x\|^2 + \|y_m - x\|^2 + 2(y_n - x) \cdot (y_m - x).$$

Combining these two equalities reveals that
$$\|y_n - y_m\|^2 = 2\left(\|y_n - x\|^2 + \|y_m - x\|^2 - 2\left\|\frac{y_n + y_m}{2} - x\right\|^2\right).$$

By convexity of $C$, we have $\frac{y_n + y_m}{2} \in C$, and therefore
$$\|y_n - y_m\|^2 \leqslant 2\left(\|y_n - x\|^2 + \|y_m - x\|^2 - 2d(x,C)\right)$$

This shows that $(y_n)_{n \geqslant 1}$ is Cauchy. Since $C$ is a closed subspace of the Hilbert space $H$, it is a complete metric space. The sequence $(y_n)_{n \geqslant 1}$ therefore converges to some $y_0 \in C$. By continuity of the norm, we have $\|y_0 - x\|^2 = d(x,C)$ which completes the proof.

**Exercise 2.17.** We treat each question separately.

(i) To show that $H = C \oplus C^\perp$ we need to show that every $x \in H$ may be written uniquely as $x = y + z$ for some $y \in C$ and $z \in C^\perp$. Fix $x \in H$, and invoke Exercise 2.16 to define $y := P_C(x)$ and $z := x - P_C(x)$. It is clear that $x = y + z$ and that $y \in C$. To see that $z \in C^\perp$, fix $w \in C$. Since $C$ is a vector space, applying the characterization (2.29) of the projection $P_C(x)$ to the vector $w + P_C(x)$ reveals that
$$z \cdot w = (x - P_C(x)) \cdot (w + P_C(x) - P_C(x)) \leqslant 0.$$

Since $C$ is symmetric, this shows that $z \in C^\perp$, so $x$ may be written as $x = y + z$ for some $y \in C$ and $z \in C^\perp$. To show that this decomposition is unique, suppose that we may also write $x = w + v$ for some other $w \in C$ and $v \in C^\perp$. Since $C$ and $C^\perp$ are both vector spaces, we have $y - w \in C$ and $z - v \in C^\perp$, so Pythagoras' theorem implies that
$$0 = \|x - x\|^2 = \|y - w + z - v\|^2 = \|y - w\|^2 + \|z - v\|^2$$
which implies that $y = w$ and $z = v$ as required.

(ii) The uniqueness of $y \in H$ is clear. Indeed, if $f(x) = x \cdot y = x \cdot y'$ for some $y, y' \in H$, then $x \cdot (y - y') = 0$ for all $x \in H$, and choosing $x = y - y'$ shows that $y = y'$. We now establish the existence of $y$. If $f = 0$, we can choose $y = 0$, so let us assume that $f \neq 0$. Consider the closed subspace $C := \ker f$ of $H$. Since $f \neq 0$, the space $C$ is not all of $H$, so Exercise 2.17 gives $z \in C^\perp$ with $\|z\| = 1$. For each $x \in H$, we have $f(x)z - f(z)x \in C$ by linearity of $f$. It follows that
$$0 = z \cdot (f(x)z - f(z)x) = f(x) - f(z)z \cdot x.$$

Setting $y := zf(z)$ completes the proof.

**Exercise 2.18.** Fix a step function of the form $\phi := \sum_{i=1}^n \phi_i \mathbf{1}_{[a_{i-1}, a_i)}$ for some $(\phi_i)_{i \leq n} \subseteq \mathbb{R}$ and non-decreasing $(a_i)_{0 \leq i \leq n} \subseteq [a,b]$. If we denote by $L$ the Lipschitz constant of $F$, then

$$|T(\phi)| \leq \sum_{i=1}^n |\phi_i| |F(a_i) - F(a_{i-1})| \leq L \sum_{i=1}^n |\phi_i| (a_i - a_{i-1}) = L \int_a^b |\phi(x)| \, dx$$

It follows by the Cauchy-Schwarz inequality that $T$ defines a continuous linear functional on the space of step functions on $L^2([a,b];\mathbb{R})$. By density, it extends to a unique continuous linear functional on the Hilbert space $L^2([a,b];\mathbb{R})$. Invoking the Riesz representation theorem on Hilbert spaces established in Exercise 2.17 gives a unique function $f \in L^2([a,b];\mathbb{R})$ with

$$T(g) = \int_a^b f(t) g(t) \, dt$$

for all $g \in L^2([a,b];\mathbb{R})$. If we now fix $x \in [a,b]$, and apply this representation to the step function $g := \mathbf{1}_{[a,x)}$ we find that $F(x) - F(a) = T(g) = \int_a^x f(t) \, dt$. It follows that for any $h \in \mathbb{R} \setminus \{0\}$ with $x + h \in [a,b]$,

$$\left| \frac{F(x+h) - F(x)}{h} - f(x) \right| \leq \frac{1}{|h|} \int_{x-|h|}^{x+|h|} |f(t) - f(x)| \, dt.$$

Invoking the Lebesgue differentiation theorem established in Theorem A.16 completes the proof.

**Exercise 2.19.** We treat each question separately.

(i) The convex dual of $x \mapsto I_{\{x \geq 0\}}$ is the mapping

$$z \mapsto \sup_{x \in \mathbb{R}^d} \left( x \cdot z - I_{\{x \geq 0\}} \right) = I_{\{z \leq 0\}}.$$

(ii) We start by observing that

$$I_{\{x \geq 0\}} = \sup_{z \in \mathbb{R}^d} \left( x \cdot z - I_{\{z \leq 0\}} \right), \tag{S.7}$$

as can be checked directly or by appealing to the previous step and the Fenchel-Moreau theorem. We can also write

$$I_{\{Ax \geq c\}} = \sup_{y \in \mathbb{R}^k} \left( (Ax - c) \cdot y - I_{\{y \leq 0\}} \right)$$
$$= \sup_{y \in \mathbb{R}^k} \left( x \cdot A^* y - c \cdot y - I_{\{y \leq 0\}} \right).$$

Summing the two previous displays yields that

$$\phi(x) = \sup_{z \in \mathbb{R}^d,\, y \in \mathbb{R}^k} \left( x \cdot z - I_{\{z \leq 0\}} + x \cdot A^* y - c \cdot y - I_{\{y \leq 0\}} \right)$$

$$= \sup_{z \in \mathbb{R}^d,\, y \in \mathbb{R}^k} \left( x \cdot z - c \cdot y - I_{\{z \leq A^* y\}} - I_{\{y \leq 0\}} \right)$$

$$= \sup_{z \in \mathbb{R}^d} \left( x \cdot z - \psi(z) \right).$$

To go from the first to the second line, we took the supremum over $y$ first and then replaced $z$ by $z - A^* y$. We obtained the announced identity.

(iii) The infimum on the left side of (2.37) can be rewritten as

$$-\sup_{x \in \mathbb{R}^d} \left( -b \cdot x - \phi(x) \right) = -\phi^*(-b),$$

and by the previous step, this is $-\psi^{**}(-b)$, provided that $\psi^{**}$ is well-defined. We also observe that

$$-\psi(-b) = -\inf_{y \in \mathbb{R}^k} \left( c \cdot y + I_{\{y \leq 0\}} + I_{\{A^* y \geq -b\}} \right)$$

$$= \sup_{y \in \mathbb{R}^k} \left( -c \cdot y - I_{\{-y \geq 0\}} - I_{\{-A^* y \leq b\}} \right)$$

$$= \sup_{y \in \mathbb{R}^k} \left( c \cdot y - I_{\{y \geq 0\}} - I_{\{A^* y \leq b\}} \right),$$

using the change of variables $y \mapsto -y$. This is the optimization problem on the right side of (2.37), so our goal is to show that $\psi^{**}(-b) = \psi(-b)$. Since

$$\psi(z) = \inf_{y \in \mathbb{R}^k} \left( -c \cdot y + I_{\{-y \in K_z\}} \right),$$

the assumption on $K_z$ ensures that $\psi$ is finite in a neighbourhood of $-b$. Arguing as in Exercise 2.8, we see that $\psi$ must be convex, with the caveat that it might a priori take the value $-\infty$. But by Remark 2.7, the fact that $\psi$ is finite on some open set rules out this possibility, so $\psi$ takes values in $\mathbb{R} \cup \{+\infty\}$. Let $\underline{\psi}$ be the lower semi-continuous envelope of $\psi$, as defined in Exercise 2.9. By Exercise 2.13, we have that $\psi^{**} = \underline{\psi}$. By Proposition 2.9, the function $\psi$ is continuous at $-b$, and thus part (ii) of Exercise 2.9 ensures that $\underline{\psi}(-b) = \psi(-b)$. This completes the proof.

**Exercise 2.20.** To alleviate notation, we introduce the set

$$C := \mathbf{n}_{\mathrm{dom}\, f}(x) + \mathrm{conv}(S(x)).$$

The proof proceeds in three steps. First, we show that $C \subseteq \partial f(x)$, then we prove that $C$ is closed and convex, and finally we show that $\partial f(x) \subseteq C$.

*Step 1: $C \subseteq \partial f(x)$.* By Theorem 2.13, Proposition 2.14 and the continuity of $f$ on $\operatorname{dom} f$, we have $S(x) \subseteq \partial f(x)$. Together with the convexity of $\partial f(x)$, this shows that $\operatorname{conv}(S(x)) \subseteq \partial f(x)$. We now prove that $\partial f(x) + \mathbf{n}_{\operatorname{dom} f}(x) \subseteq \partial f(x)$. Fix $p \in \partial f(x)$ and $n \in \mathbf{n}_{\operatorname{dom} f}(x)$. For every $x' \in \operatorname{dom} f$, we have

$$f(x') \geqslant f(x) + p \cdot (x' - x) \geqslant f(x) + (p+n) \cdot (x' - x),$$

where the first inequality uses that $p \in \partial f(x)$, and the second uses that $n \in \mathbf{n}_{\operatorname{dom} f}(x)$. This shows that $p + n \in \partial f(x)$, and therefore that $\partial f(x) + \mathbf{n}_{\operatorname{dom} f}(x) \subseteq \partial f(x)$. In particular, we have $C \subseteq \partial f(x)$.

*Step 2: $C$ is closed and convex.* By Exercise 2.1, the set $C$ is convex. To show that it is also closed, fix a sequence $(p_i)_{i \geqslant 1} \subseteq C$ converging to some point $p \in \mathbb{R}^d$. For each $i \geqslant 1$, let $n_i \in \mathbf{n}_{\operatorname{dom} f}(x)$ and $s_i \in \operatorname{conv}(S(x))$ be such that $p_i = n_i + s_i$. The Lipschitz continuity of $f$ on $\operatorname{dom} f$ ensures that the sequence $(s_i)_{i \geqslant 1}$ is bounded. Up to passing to a subsequence, we may therefore suppose that it converges to some point $s \in \overline{\operatorname{conv}}(S(x)) = \operatorname{conv}(S(x))$. We have implicitly used Exercise 2.5 to say that $\operatorname{conv}(S(x))$ is compact, and therefore closed, as the convex hull of the compact set $S(x)$. The sequence $(n_i)_{i \geqslant 1}$ also converges along this subsequence to some point $n \in \mathbf{n}_{\operatorname{dom} f}(x)$. We have implicitly used that the normal cone is closed. Taking the limit along this subsequence shows that $p = n + s \in C$, and therefore that $C$ is closed.

*Step 3: $\partial f(x) \subseteq C$.* Suppose for the sake of contradiction that there exists $p \in \partial f(x) \smallsetminus C$. The supporting hyperplane theorem (Theorem 2.2) then gives a non-zero vector $v \in \mathbb{R}^d$ with

$$\sup \{v \cdot (n+s) \mid n \in \mathbf{n}_{\operatorname{dom} f}(x) \text{ and } s \in \operatorname{conv}(S(x))\} < v \cdot p. \tag{S.8}$$

Recall the definition of the polar in Exercise 2.14, and observe that $v \in \mathbf{n}_{\operatorname{dom} f}(x)^\circ$. Indeed, if this were not the case, there would exist $n \in \mathbf{n}_{\operatorname{dom} f}(x)$ with $v \cdot n > 0$. Given $y_0 \in \overline{\operatorname{conv}}(S(x))$, which exists by Rademacher's theorem (Theorem 2.10), we would then have, for all $\lambda > 0$,

$$\lambda v \cdot n + v \cdot y_0 = (\lambda n + y_0) \cdot v < v \cdot p.$$

Letting $\lambda$ tend to infinity would give a contradiction. It follows by Exercise 2.15 that $v \in \mathbf{n}_{\operatorname{dom} f}(x)^\circ = \mathbf{T}_{\operatorname{dom} f}(x)$, where $\mathbf{T}_{\operatorname{dom} f}(x)$ denotes the tangent cone to $\operatorname{dom} f$ at $x$. Combining Exercises 2.15 and 2.3 with Rademacher's theorem gives a sequence $(x_i)_{i \geqslant 1} \subseteq \operatorname{int}(\operatorname{dom} f)$ converging to $x$ and a sequence $(t_i)_{i \geqslant 1} \subseteq \mathbb{R}_{>0}$ decreasing to 0 such that $\nabla f(x_i)$ exists for all $i \geqslant 1$, and the sequence $(v_i)_{i \geqslant 1}$ defined by $v_i := t_i^{-1}(x_i - x)$ converges to $v$. By Lipschitz continuity of $f$ on $\operatorname{dom} f$, the sequence $(\nabla f(x_i))_{i \geqslant 1}$ is bounded, so, up to passing to a subsequence, assume without loss of generality that it converges to some $c \in C$. Since $\nabla f(x+t_i v_i) \in \partial f(x+t_i v_i)$ by Theorem 2.13 and $p \in \partial f(x)$ by assumption, for every $i \geq 1$, we have

$$-\nabla f(x+t_i v_i) \cdot v_i \leq \frac{f(x) - f(x+t_i v_i)}{t_i} \leq -v_i \cdot p.$$

Letting $i$ tend to infinity shows that $-c \cdot v \leq -v \cdot p$, and therefore $c \cdot v \geq v \cdot p$. This contradicts (S.8) and completes the proof.

**Exercise 2.21.** A direct computation shows that the rate function $I$ in (2.51) is decreasing on $(0, p)$ and increasing on $(p, 1)$. Using this, it is readily verified that the string of inequalities (2.56) implies (2.54) and (2.55). Conversely, let us assume (2.54) and (2.55), and deduce (2.56). Fix a Borel set $A \subseteq \mathbb{R}$, and observe that

$$\mathbb{P}\{S_N \in A\} \leq \mathbb{P}\{S_N \in A \cap (-\infty, p]\} + \mathbb{P}\{S_N \in A \cap [p, +\infty)\}.$$

Denote by $x \in \overline{A}$ the infimum of $A \cap [p, +\infty)$. Since $x \geq p$, applying (2.54) reveals that

$$\limsup_{N \to +\infty} \frac{1}{N} \log \mathbb{P}\{S_N \in A \cap [p, +\infty)\} \leq \limsup_{N \to +\infty} \log \mathbb{P}\{S_N \geq x\} = -I(x) \leq -\inf_{x \in \overline{A}} I(x).$$

A similar argument shows that

$$\limsup_{N \to +\infty} \frac{1}{N} \mathbb{P}\{S_N \in A \cap (-\infty, p]\} \leq -\inf_{x \in \overline{A}} I(x).$$

Combining these two bounds establishes the upper bound in (2.56). To prove the lower bound, fix $x \in \text{int}(A)$, and let $\varepsilon > 0$ be such that $[x - \varepsilon, x + \varepsilon] \subseteq A$. Observe that

$$\mathbb{P}\{S_N \in A\} \geq \mathbb{P}\{S_N \in [x - \varepsilon, x + \varepsilon)\} = \mathbb{P}\{S_N \geq x - \varepsilon\} - \mathbb{P}\{S_N \geq x + \varepsilon\}$$

If $x > p$ and $\varepsilon > 0$ is small enough so $x - \varepsilon > p$, then (2.54) implies that

$$\lim_{N \to +\infty} \frac{1}{N} \log \mathbb{P}\{S_N \geq x - \varepsilon\} = -I(x - \varepsilon) \quad \text{and} \quad \lim_{N \to +\infty} \frac{1}{N} \log \mathbb{P}\{S_N \geq x + \varepsilon\} = -I(x + \varepsilon).$$

Since $I$ is increasing on the interval $(p, 1)$, we have $-I(x - \varepsilon) \geq -I(x + \varepsilon)$, and therefore

$$\liminf_{N \to +\infty} \frac{1}{N} \log \mathbb{P}\{S_N \in A\} \geq -I(x - \varepsilon).$$

Letting $\varepsilon$ tend to zero and using the continuity of $I$ gives the lower bound

$$\liminf_{N \to +\infty} \frac{1}{N} \log \mathbb{P}\{S_N \in A\} \geq -I(x).$$

An identical bound holds if $x \leq p$. Taking the supremum over all $x \in \text{int}(A)$ establishes the lower bound in the string of inequalities (2.56) and completes the proof.

**Exercise 2.22.** The log-Laplace transform of a Bernoulli random variable $X$ is

$$\psi(\lambda) = \log \mathbb{E}\exp(\lambda X) = \log\bigl(pe^\lambda + (1-p)\bigr).$$

It follows by definition of the convex dual that

$$\psi^*(x) = \sup_{\lambda \in \mathbb{R}} \bigl(\lambda x - \log\bigl(pe^\lambda + (1-p)\bigr)\bigr). \tag{S.9}$$

Elementary calculus reveals that this supremum is achieved for $\lambda \in \mathbb{R}$ with

$$x = \frac{pe^\lambda}{pe^\lambda + (1-p)}.$$

Rearranging shows that

$$e^\lambda = \frac{x}{p} \cdot \frac{1-p}{1-x}.$$

Substituting this into (S.9) gives

$$\psi^*(x) = x\log\left(\frac{x}{p}\right) - x\log\left(\frac{1-x}{1-p}\right) - \log\left(x \cdot \frac{1-p}{1-x} + (1-p)\right)$$
$$= x\log\left(\frac{x}{p}\right) - x\log\left(\frac{1-x}{1-p}\right) - \log\left(\frac{1-p}{1-x}\right)$$
$$= x\log\left(\frac{x}{p}\right) + (1-x)\log\left(\frac{1-x}{1-p}\right)$$

as required.

**Exercise 2.23.** We recall that by Proposition 2.9, the function $f$ is continuous; we also recall that

$$f^*(\lambda) = \sup_{x \in \mathbb{R}^d} (\lambda \cdot x - f(x)). \tag{S.10}$$

We fix $\lambda \in \mathbb{R}^d$. The assumption of (2.119) and the continuity of $f$ guarantee that for some constant $C < +\infty$, we have for every $x \in \mathbb{R}^d$ that

$$f(x) \geq (|\lambda| + 2)|x| - C.$$

Plugging this estimate into (S.10) and using the Cauchy-Schwarz inequality yields that, for every $\lambda' \in \mathbb{R}^d$,

$$\lambda' \cdot x - f(x) \leq C - (|\lambda| + 2 - |\lambda'|)|x|.$$

In particular, we have that $f^*(\lambda)$ is finite. The bound above also implies that for every $\lambda' \in \mathbb{R}^d$ with $|\lambda'| \leq |\lambda| + 1$, we may as well restrict the supremum in the definition of $f^*(\lambda')$ in (S.10) to values of $x$ that range in a fixed compact set $K$. Since $f$ is also continuous, the supremum in (S.10) is achieved for those values of $\lambda'$. Since the mapping $(\lambda, x) \mapsto \lambda \cdot x - f(x)$ is continuous and continuously differentiable with respect to its first variable, we are thus in position to appeal to the envelope theorem (Theorem 2.21). The claim will follow once we show that the supremum in (S.10) is achieved at a unique point $x^*(\lambda)$, and that $\lambda \mapsto x^*(\lambda)$ is continuous.

Suppose that there are two distinct optimizers $a \neq b \in \mathbb{R}^d$ for (S.10). We have in particular that
$$\lambda \cdot a - f(a) = \lambda \cdot b - f(b).$$
The strict convexity assumption then implies that
$$\frac{1}{2}\bigl(\lambda \cdot a - f(a) + \lambda \cdot b - f(b)\bigr) < \lambda \cdot \frac{a+b}{2} - f\Bigl(\frac{a+b}{2}\Bigr).$$
Since this contradicts the optimality assumption on $a$ and $b$, this completes the proof of uniqueness. Denoting by $x^*(\lambda)$ the optimum, we now show that the mapping $\lambda \mapsto x^*(\lambda)$ is continuous. Let $(\lambda_n)_{n \geq 1}$ be a sequence of points in $\mathbb{R}^d$ converging to $\lambda$. Up to extraction of a subsequence, we may assume that $x^*(\lambda_n)$ converges to some $a \in \mathbb{R}^d$. Passing to the limit in the identity
$$f^*(\lambda_n) = \lambda_n \cdot x^*(\lambda_n) - f(x^*(\lambda_n)),$$
we conclude that $a$ must be an optimizer in (S.10). By uniqueness of the optimizer, we have that $a = x^*(\lambda)$, as desired.

**Exercise 2.24.** We treat each question separately.

(i) The envelope theorem (Theorem 2.21) implies that $\partial_t f(t, 0) = m_0(t)^2$. If $(t_n)_{n \geq 1} \subseteq \mathbb{R}_{\geq 0}$ is a sequence converging to some $t \in \mathbb{R}_{\geq 0}$, then for every $n \geq 1$ we have $\partial_t f(t_n, 0) = m_0(t_n)^2$. The sequence $(m_0(t_n))_{n \geq 1}$ is uniformly bounded by one, and by (2.115) any of its subsequential limits must satisfy the fixed point equation $m = \tanh(2tm)$. Together with (2.116) this means that it must be one of $\pm m_0(t)$. It follows that
$$\lim_{n \to +\infty} \partial_t f(t_n, 0) = \lim_{n \to +\infty} m_0(t_n)^2 = m_0(t)^2 = \partial_t f(t, 0).$$
This establishes the continuity of $t \mapsto \partial_t f(t, 0)$.

(ii) Fix $h > 0$, and observe that by Proposition 2.22,
$$f(t_c, h) = t_c m_h^+(t_c)^2 + h m_h^+(t_c) - \frac{1 + m_h^+(t_c)}{2} \log(1 + m_h^+(t_c))$$
$$- \frac{1 - m_h^+(t_c)}{2} \log(1 - m_h^+(t_c)).$$

Remembering that $t_c = \frac{1}{2}$ and Taylor expanding the logarithm reveals that

$$f(t_c,h) = \frac{1}{2}m_h^+(t_c)^2 + hm_h^+(t_c) - \frac{1}{2}m_h^+(t_c)^2 + \mathcal{O}\left(m_h^+(t_c)^3\right)$$
$$= (1+o(1))hm_h^+(t_c),$$

where the last equality uses the fact that $m_h(t_c) \to m_0(t_c) = 0$ as $h$ decreases to zero. On the other hand, Taylor expanding the hyperbolic tangent and using the fixed point equation (2.115) implies that

$$m_h(t_c) = m_h(t_c) + h - \frac{(m_h(t_c)+h)^3}{3} + \mathcal{O}\left((m_h(t_c)+h)^5\right)$$
$$= m_h(t_c) - \frac{m_h(t_c)^3}{3} + h(1+o(1)).$$

Rearranging gives $m_h(t_c) = h^{1/3}(1+o(1))$. It follows that $f(t_c,h) = h^{1+\frac{1}{3}}(1+o(1))$, and therefore

$$\lim_{h \searrow 0} \frac{\log f(t_c,h)}{\log h} = 1 + \frac{1}{3}.$$

An identical argument for $h < 0$ allows us to conclude that $\delta = 3$.

(iii) We can and do restrict our attention to the case of $t > t_c$. Taylor expanding the hyperbolic tangent and using the fixed point equation (2.115) implies that, as $t$ tends to $t_c$,

$$m_0(t) = 2tm_0(t) - \frac{8t^3 m_0(t)^3}{3} + \mathcal{O}\left((tm_0(t))^3\right)$$
$$= 2tm_0(t) - (1+o(1))\frac{8t^3 m_0(t)^3}{3}$$

where the last equality uses the fact that $m_0(t) \to m_0(t_c) = 0$ as $t$ tends to $t_c$. Rearranging shows that

$$2(t-t_c)m_0(t) = (1+o(1))\frac{8t^3 m_0(t)^3}{3},$$

and therefore

$$m_0(t) = (1+o(1))\left(\frac{3(t-t_c)}{4t^3}\right)^{1/2}.$$

It follows that

$$\lim_{t \to t_c} \frac{\log m_0(t)}{\log(t-t_c)} = \frac{1}{2},$$

so $\beta = \frac{1}{2}$.

## S.3 Hamilton-Jacobi equations

**Exercise 3.1.** Fix $r > 0$ sufficiently small that $x$ is a strict maximum of $f$ in $\overline{B}_r(x)$. By continuity of $f_N$ and compactness of $\overline{B}_r(x)$, let $x_N \in \overline{B}_r(x)$ be such that

$$f_N(x_N) \geqslant f_N(y) \tag{S.11}$$

for every $y \in \overline{B}_r(x)$. Since $(x_N)_{N \geqslant 1}$ stays in a compact set, it admits a subsequence converging to some $x^* \in \overline{B}_r(x)$. Letting $N$ tend to infinity in (S.11) and leveraging the local uniform convergence of $f_N$ to $f$ reveals that $f(x^*) \geqslant f(y)$ for every $y \in \overline{B}_r(x)$. Since $x$ is a strict maximum we must have $x^* = x$, so the only limit point of $(x_N)_{N \geqslant 1}$ is $x$. This implies that $(x_N)_{N \geqslant 1}$ converges to $x$ as required. In particular, for $N$ sufficiently large $x_N$ is in the open ball $B_r(x)$, and is thus a local maximum of $f_N$.

**Exercise 3.2.** It suffices to prove that the definition of subsolution with "strict local maximum" implies the one with "local maximum". Let $f : \mathbb{R}_{\geqslant 0} \times \mathbb{R}^d \to \mathbb{R}$ be a subsolution for the definition with "strict local maximum", and fix a smooth function $\phi \in C^\infty(\mathbb{R}_{>0} \times \mathbb{R}^d; \mathbb{R})$ with the property that $f - \phi$ has a local maximum at $(t^*, x^*) \in \mathbb{R}_{>0} \times \mathbb{R}^d$. Consider the smooth function $\eta : \mathbb{R}_{>0} \times \mathbb{R}^d \to \mathbb{R}$ defined by

$$\eta(t, x) := \phi(t, x) + \varepsilon\big((t - t^*)^2 + |x - x^*|^2\big).$$

Observe that in a small enough neighbourhood of $(t^*, x^*)$, we have

$$(f - \eta)(t, x) \leqslant (f - \phi)(t^*, x^*) - \varepsilon\big(|t - t^*|^2 + |x - x^*|^2\big) < (f - \eta)(t^*, x^*),$$

so $(t^*, x^*)$ is a strict local maximum of $f - \eta$. It follows that

$$\big(\partial_t \phi - \mathsf{H}(\nabla \phi)\big)(t^*, x^*) = \big(\partial_t \eta - \mathsf{H}(\nabla \eta)\big)(t^*, x^*) \leqslant 0$$

which is the required subsolution criterion for the definition with "local maximum".

**Exercise 3.3.** Arguing as in Exercise 3.2, it suffices to prove that the definition of subsolution with "global maximum" implies the one with "local maximum". Let $f : \mathbb{R}_{\geqslant 0} \times \mathbb{R}^d \to \mathbb{R}$ be a subsolution for the definition with "global maximum", and let $\phi \in C^\infty(\mathbb{R}_{>0} \times \mathbb{R}^d; \mathbb{R})$ be such that $f - \phi$ has a local maximum at $(t^*, x^*) \in \mathbb{R}_{>0} \times \mathbb{R}^d$. Denote by $B_r(t^*, x^*)$ an open ball centred at $(t^*, x^*)$ on which $f - \phi$ attains its maximum at $(t^*, x^*)$. Let $\theta : \mathbb{R}_{\geqslant 0} \times \mathbb{R}^d \to \mathbb{R}$ be a non-negative smooth function with $\theta = 1$ on $\overline{B}_{r/2}(t^*, x^*)$ and $\theta = 0$ outside $B_{2r/3}(t^*, x^*)$, and define the parabola

$$\chi(t, x) := (t - t^*)^2 + |x - x^*|^2.$$

Introduce the constants

$$L := \sup_{t > 0} \|f(t, \cdot)\|_{\mathrm{Lip}}, \quad M := \sup_{B_r(t^*, x^*)} |\phi|, \quad M' := \sup_{B_r(t^*, x^*)} \big(|\partial_t \theta| + \|\nabla \theta\|_2\big),$$

and define the smooth function $\eta : \mathbb{R}_{>0} \times \mathbb{R}^d \to \mathbb{R}$ by

$$\eta(t,x) := \phi(t,x)\theta(t,x) + \left(\frac{L}{r} + \frac{M}{r^2} + \frac{2MM'}{r}\right)\chi(t,x).$$

We will now verify that $(t^*, x^*)$ is a global maximum of $(f - \eta)$. On the one hand, for any $(t,x) \in \overline{B}_{r/2}(t^*, x^*)$,

$$(f - \eta)(t,x) \leq (f - \phi)(t,x) \leq (f - \phi)(t^*, x^*) = (f - \eta)(t^*, x^*).$$

On the other hand, for any $(t,x) \in \overline{B}_r(t^*, x^*) \setminus B_{r/2}(t^*, x^*)$,

$$(f - \eta)(t,x) \leq (f - \phi)(t,x) + |\phi(t,x)||\theta(t,x) - \theta(t^*, x^*)| - \frac{2MM'}{r}\chi(t,x)$$

$$\leq (f - \phi)(t^*, x^*) + \frac{MM'r}{2}\left(\frac{|t - t^*|}{r/2} + \frac{|x - x^*|}{r/2}\right) - \frac{2MM'}{r}\chi(t,x)$$

$$\leq (f - \eta)(t^*, x^*),$$

where we used the mean value theorem in the second inequality. Finally, for any $(t,x) \in \mathbb{R}^d \setminus B_r(t^*, x^*)$,

$$(f - \eta)(t,x) \leq (f - \phi)(t^*, x^*) + |f(t,x) - f(t^*, x^*)| + |\phi(t^*, x^*)| - \left(\frac{L}{r} + \frac{M}{r^2}\right)\chi(t,x)$$

$$\leq (f - \eta)(t^*, x^*) + Lr\left(\frac{|t - t^*|}{r} + \frac{|x - x^*|}{r}\right) + M - \left(\frac{L}{r} + \frac{M}{r^2}\right)\chi(t,x)$$

$$\leq (f - \eta)(t^*, x^*).$$

This shows that $(t^*, x^*)$ is indeed a global maximum of $f - \eta$. Since $f$ is a subsolution for the definition with "global maximum", it follows that

$$\left(\partial_t \phi - \mathsf{H}(\nabla\phi)\right)(t^*, x^*) = \left(\partial_t \eta - \mathsf{H}(\nabla\eta)\right)(t^*, x^*) \leq 0$$

which is the required subsolution criterion for the definition with "local maximum".

**Exercise 3.4.** Given a function $\phi \in C^1(\mathbb{R}_{>0} \times \mathbb{R}^d; \mathbb{R})$, denote by $P(\phi)$ the property that for every $(t^*, x^*) \in \mathbb{R}_{>0} \times \mathbb{R}^d$ for which $f - \phi$ has a local maximum at $(t^*, x^*)$, we have

$$\left(\partial_t \phi - \mathsf{H}(\nabla\phi)\right)(t^*, x^*) \leq 0.$$

It is clear that whenever $P(\phi)$ holds for every $\phi \in C^1(\mathbb{R}_{>0} \times \mathbb{R}^d; \mathbb{R})$, it also holds for every smooth function $\phi \in C^\infty(\mathbb{R}_{>0} \times \mathbb{R}^d; \mathbb{R})$. Let us now suppose that $P(\phi)$ holds for every $\phi \in C^\infty(\mathbb{R}_{>0} \times \mathbb{R}^d; \mathbb{R})$, and fix $\phi \in C^1(\mathbb{R}_{>0} \times \mathbb{R}^d; \mathbb{R})$ and $(t^*, x^*) \in \mathbb{R}_{>0} \times \mathbb{R}^d$ with the property that $f - \phi$ has a strict local maximum at $(t^*, x^*)$. For each $\varepsilon > 0$, denote by $\phi_\varepsilon \in C^\infty(\mathbb{R}_{>0} \times \mathbb{R}^d; \mathbb{R})$ the mollification of $\phi$. Since $(\phi_\varepsilon)_{\varepsilon>0}$ converges locally uniformly to $\phi$, Exercise 3.1 gives a sequence $(t_\varepsilon, x_\varepsilon)_{\varepsilon>0} \subseteq \mathbb{R}_{>0} \times \mathbb{R}^d$ converging to $(t^*, x^*)$ with the property that $(t_\varepsilon, x_\varepsilon)$ is a local maximum of $f - \phi_\varepsilon$ for each $\varepsilon > 0$. The property $P(\phi_\varepsilon)$ implies that for every $\varepsilon > 0$,

$$\left(\partial_t \phi_\varepsilon - \mathsf{H}(\nabla \phi_\varepsilon)\right)(t_\varepsilon, x_\varepsilon) \leq 0. \tag{S.12}$$

Since $(\partial_t \phi_\varepsilon)_{\varepsilon>0}$ and $(\nabla \phi_\varepsilon)_{\varepsilon>0}$ converge to $\partial_t \phi$ and $\nabla \phi$ locally uniformly by the properties of mollifiers, letting $\varepsilon$ tend to zero in (S.12) and using the local Lipschitz continuity of $\mathsf{H}$ shows that

$$\left(\partial_t \phi - \mathsf{H}(\nabla \phi)\right)(t^*, x^*) \leq 0.$$

Invoking Exercise 3.2 establishes $P(\phi)$ and completes the proof.

**Exercise 3.5.** Fix $\phi \in C^\infty(\mathbb{R}_{\geq 0} \times \mathbb{R}^d; \mathbb{R})$ with the property that $f - \phi$ has a local maximum at the point $(t^*, x^*) \in \mathbb{R}_{>0} \times \mathbb{R}^d$. Since $f \in C^1(\mathbb{R}_{\geq 0} \times \mathbb{R}^d; \mathbb{R})$ and $t > 0$, we have

$$(\partial_t f - \partial_t \phi)(t^*, x^*) = 0 \quad \text{and} \quad (\nabla f - \nabla \phi)(t^*, x^*) = 0.$$

Combining this with the assumption that $f$ satisfies the equation at $(t^*, x^*)$ shows that $f$ is a viscosity subsolution. An identical argument shows that it is also a viscosity supersolution. This completes the proof.

**Exercise 3.6.** Let $\phi \in C^\infty(\mathbb{R}_{>0} \times \mathbb{R}^d; \mathbb{R})$ be such that $(f - \Phi) - \phi = f - (\phi + \Phi)$ has a local maximum at $(t^*, x^*) \in \mathbb{R}_{>0} \times \mathbb{R}^d$. Since $f$ is a subsolution to the Hamilton-Jacobi equation (3.20) and $\phi + \Phi \in C^\infty(\mathbb{R}_{>0} \times \mathbb{R}^d; \mathbb{R})$, we have

$$\partial_t(\phi + \Phi)(t^*, x^*) \leq \mathsf{H}(\nabla \phi + \nabla \Phi)(t^*, x^*) \leq \mathsf{H}\left(\nabla \phi(t^*, x^*)\right) + V|\nabla \Phi(t^*, x^*)|.$$

Remembering that $\partial_t \Phi \geq V|\nabla \Phi|$ shows that $\left(\partial_t \phi - \mathsf{H}(\nabla \phi)\right)(t^*, x^*) \leq 0$ as required.

**Exercise 3.7.** Fix $x^* \in \mathbb{R}^d$, and observe that $g(t, x) \coloneqq f(t, x^* + x)$ is a viscosity solution to the Hamilton-Jacobi equation (3.20). It follows by the comparison principle that

$$\sup_{\mathbb{R}_{\geq 0} \times \mathbb{R}^d} (f - g) = \sup_{\{0\} \times \mathbb{R}^d} (f - g) \leq L|x^*|.$$

An identical argument can be used to control $g - f$ and show that for every $x \in \mathbb{R}^d$,

$$|f(t, x) - f(t, x^* + x)| \leq \sup_{\mathbb{R}_{\geq 0} \times \mathbb{R}^d} |f - g| \leq L|x^*|$$

as desired.

**Exercise 3.8.** Replacing $f(t, \cdot)$ and $g(t, \cdot)$ by $f(t, \cdot + x)$ and $g(t, \cdot + x)$ respectively, assume without loss of generality that $x = 0$. Up to renaming $f$ and $g$ if necessary, suppose also that $f(T,0) \geqslant g(T,0)$. Given $M > 2L$, the comparison principle with $R = VT$ implies that

$$f(T,0) - g(T,0) \leqslant \sup_{y \in \mathbb{R}^d} \left( f(0,y) - g(0,y) - M(|y| - VT)_+ \right). \tag{S.13}$$

If $|y| > VT$, and $z \in \partial B_{VT}(0)$ is such that $\operatorname{dist}(y, \partial B_{VT}) = |y - z|$, then

$$f(0,y) - g(0,y) - M(|y| - VT)_+ \leqslant 2L|y - z| - M|y| + MVT,$$

where we used that $f(0,z) = g(0,z)$. To bound this further, notice that by the case of equality in the triangle inequality and the choice of $z$,

$$|y| = \left| y - \frac{VT}{|y|} y \right| + \left| \frac{VT}{|y|} y \right| \geqslant |y - z| + VT.$$

It follows that for all $y \in \mathbb{R}^d$ with $|y| > VT$,

$$f(0,y) - g(0,y) - M(|y| - VT)_+ \leqslant (2L - M)(|y| - VT) \leqslant 0.$$

Since $f(0,\cdot)$ and $g(0,\cdot)$ coincide on $B_{VT}$, this inequality actually holds for every $y \in \mathbb{R}^d$. Recalling (S.13) completes the proof.

**Exercise 3.9.** The strategy will be to show that $g$ is a viscosity solution to the Hamilton-Jacobi equation (3.20). Let $\phi \in C^\infty(\mathbb{R}_{>0} \times \mathbb{R}^d; \mathbb{R})$ be such that $g - \phi$ admits a local maximum at $(t^*, x^*) \in \mathbb{R}_{>0} \times \mathbb{R}^d$. Since $g$ is a subsolution to $\partial_t g - \mathsf{K}(\nabla g) = 0$, we have

$$\left( \partial_t \phi - \mathsf{K}(\nabla \phi) \right)(t^*, x^*) \leqslant 0. \tag{S.14}$$

Moreover, for any $\varepsilon > 0$ and $z \in \mathbb{R}^d$, we have

$$\phi(t^*, x^*) - \phi(t^*, x^* + \varepsilon z) \leqslant g(t^*, x^*) - g(t^*, x^* + \varepsilon z) \leqslant L\varepsilon |z|.$$

Dividing by $\varepsilon > 0$, letting $\varepsilon$ tend to zero and choosing $z = -\nabla \phi(t^*, x^*)$ shows that $|\nabla \phi(t^*, x^*)| \leqslant L$. Since $\mathsf{K}$ and $\mathsf{H}$ coincide on the ball of radius $L$, it follows by (S.14) that

$$\left( \partial_t \phi - \mathsf{H}(\nabla \phi) \right)(t^*, x^*) \leqslant 0.$$

This shows that $g$ is a subsolution to the Hamilton-Jacobi equation (3.20). An identical argument shows that it is also a supersolution to this Hamilton-Jacobi equation. Invoking the uniqueness result in Corollary 3.7 completes the proof.

**Exercise 3.10.** We treat each question separately.

(i) The fact that $\phi_t \in C^1(\mathbb{R}^d; \mathbb{R}^d)$ is clear. Since $\nabla \psi$ is Lipschitz continuous and takes values in a bounded set, and $\nabla \mathsf{H}$ is locally Lipschitz continuous, we also see that the mapping

$$x \mapsto \nabla \mathsf{H}(\nabla \psi(x))$$

is Lipschitz continuous. We denote its Lipschitz constant by $L$, and set $T := L^{-1}$. We first show that $\phi_t$ is bijective for every $t < T$. We fix $t < T$ and $y \in \mathbb{R}^d$, and for every $x \in \mathbb{R}^d$, we set $\rho(x) := y + t \nabla \mathsf{H}(\nabla \psi(x))$. By the definition of $T$, we have for every $x, x' \in \mathbb{R}^d$ that

$$|\rho(x') - \rho(x)| \leq \frac{t}{T} |x' - x|.$$

The mapping $\rho$ is therefore a contraction. By the Banach fixed-point theorem, there exists a unique $x \in \mathbb{R}^d$ such that $x = \rho(x)$. Recalling the definition of $\rho$, we see that the property $x = \rho(x)$ is equivalent to $y = \phi_t(x)$. We have thus shown that $\phi_t$ is bijective. One can then use the implicit function theorem to guarantee that the inverse of $\phi_t$ also belongs to $C^1(\mathbb{R}^d; \mathbb{R}^d)$.

(ii) We denote the inverse mapping to $\phi_t$ by $\phi_t^{-1}$. Arguing as in the previous step, one can show that the mapping $(t, x) \mapsto \phi_t^{-1}(x)$ is in $C^1([0, T] \times \mathbb{R}^d; \mathbb{R})$. Introducing the function $\mathsf{K} : \mathbb{R}^d \to \mathbb{R}$ defined by

$$\mathsf{K}(p) := \mathsf{H}(p) - \nabla \mathsf{H}(p) \cdot p,$$

we can represent the function $u$ as

$$u(t, x) = \psi\left(\phi_t^{-1}(x)\right) + t \mathsf{K}\left(\nabla \psi\left(\phi_t^{-1}(x)\right)\right).$$

This shows that $u \in C^1([0, T] \times \mathbb{R}^d; \mathbb{R})$. To compute its derivatives, let us fix some notation. We write $\mathrm{I}_d \in \mathbb{R}^{d \times d}$ for the identity matrix on $\mathbb{R}^d$, and given a continuously differentiable function $v : \mathbb{R}^d \to \mathbb{R}^k$, we denote by $\nabla v : \mathbb{R}^d \to \mathbb{R}^{d \times k}$ its Jacobian,

$$\nabla v = \begin{bmatrix} \partial_{x_1} v_1 & \cdots & \partial_{x_1} v_k \\ \vdots & \ddots & \vdots \\ \partial_{x_d} v_1 & \cdots & \partial_{x_d} v_k \end{bmatrix}. \tag{S.15}$$

In this notation, given another continuously differentiable function $w : \mathbb{R}^m \to \mathbb{R}^d$, the chain rule reads

$$\nabla \big(v(w(x))\big) = \nabla w(x) \nabla v(w(x)).$$

Taking the gradient on both sides of the identity

$$\phi_t\left(\phi_t^{-1}(x)\right) = x \tag{S.16}$$

reveals that

$$\nabla \phi_t^{-1}(x) \left( \mathrm{I}_d - t \nabla^2 \psi\left(\phi_t^{-1}(x)\right) \nabla^2 \mathsf{H}\left(\nabla \psi\left(\phi_t^{-1}(x)\right)\right) \right) = \mathrm{I}_d. \tag{S.17}$$

Combining the fact that
$$\nabla K(p) = -\nabla^2 H(p)p \tag{S.18}$$
with a direct computation shows that
$$\nabla u(t,x) = \nabla \phi_t^{-1}(x)\left(\nabla \psi(z) - t\nabla^2 \psi(z) \nabla^2 H(\nabla \psi(z))\nabla \psi(z)\right)$$
$$= \nabla \phi_t^{-1}(x)\left(I_d - t\nabla^2 \psi(z) \nabla^2 H(\nabla \psi(z))\right)\nabla \psi(z),$$
where we have written $z := \phi_t^{-1}(x)$ to simplify notation. Together with the identity (S.17), this implies that
$$\nabla u(t,x) = \nabla \psi\left(\phi_t^{-1}(x)\right).$$
and therefore that the gradient of $u$ is constant along the characteristic line $t \mapsto (t, \phi_t(x))$, with $\nabla u(t, \phi_t(x)) = \nabla \psi(x)$.

(iii) To deduce that $u$ satisfies the Hamilton-Jacobi equation (3.20) at every point in $[0,T) \times \mathbb{R}^d$, it suffices to show that for every $t < T$ and $x \in \mathbb{R}^d$, we have
$$\partial_t u(t,x) = H\left(\nabla \psi\left(\phi_t^{-1}(x)\right)\right). \tag{S.19}$$
Consistently with the notation in (S.15), we think of $\partial_t \phi_t$ and $\partial_t \phi_t^{-1}$ as "row vectors". In particular, we have that
$$\partial_t \phi_t(x) = \nabla^* H(\nabla \psi(x)),$$
where $\nabla^* H(p)$ is the row vector obtained as the transposition of the column vector $\nabla H(p)$. Taking the derivative on both sides of (S.16) and rearranging gives the identity
$$\partial_t \phi_t^{-1}(x)\left(I_d - t\nabla^2 \psi\left(\phi_t^{-1}(x)\right)\nabla^2 H\left(\nabla \psi\left(\phi_t^{-1}(x)\right)\right)\right) = \nabla^* H\left(\nabla \psi\left(\phi_t^{-1}(x)\right)\right). \tag{S.20}$$
To simplify notation, again let $z := \phi_t^{-1}(x)$. Recalling (S.18), we can write
$$\partial_t u(t,x) = \partial_t \phi_t^{-1}(x) \nabla \psi(z) + K(\nabla(\psi(z))) + t\partial_t \phi_t^{-1}(x) \nabla^2 \psi(z) \nabla K(\nabla \psi(z))$$
$$= \partial_t \phi_t^{-1}(x) \nabla \psi(z) + H(\nabla \psi(z)) - \nabla H(\nabla \psi(z)) \cdot \nabla \psi(z)$$
$$- t\partial_t \phi_t^{-1}(x) \nabla^2 \psi(z) \nabla H(\nabla \psi(z)) \nabla \psi(z)$$
$$= \partial_t \phi_t^{-1}(x)\left(I_d - t\nabla^2 \psi(z) \nabla^2 H(\nabla \psi(z))\right)\nabla \psi(z)$$
$$+ H(\nabla \psi(z)) - \nabla H(\nabla \psi(z)) \cdot \nabla \psi(z).$$
Using also (S.20), we obtain the desired result (S.19).

## S.4  Statistical inference

**Exercise 4.1.** Since the matrix $A$ has rank one, its image is a one-dimensional subspace of $\mathbb{R}^d$. It is therefore generated by some non-zero vector $u \in \mathbb{R}^d$. In particular, for every $w \in \mathbb{R}^d$, there exists $\lambda(w) \in \mathbb{R}$ with $Aw = \lambda(w)u$. Defining the vector $v \in \mathbb{R}^d$ by $v_i = \lambda(e_i)$ for $1 \leqslant i \leqslant d$, where $e_i$ denotes the canonical basis vector in $\mathbb{R}^d$, gives the decomposition $A = uv^*$. Normalizing $u$ and $v$ so that they become unit vectors gives the decomposition $A = \lambda \bar{u}\bar{v}^*$ for some scalar $\lambda \in \mathbb{R}$ and some unit vectors $\bar{u}, \bar{v} \in \mathbb{R}^d$. Notice that $v$ is non-zero since $A$ would otherwise map $\mathbb{R}^d$ onto zero and therefore be of rank zero. The symmetry of $A$ implies that $\bar{u}\bar{v}^* = \bar{v}\bar{u}^*$. It follows that
$$(\bar{u}\cdot\bar{v})^2 = \bar{u}^*\,\bar{v}\bar{u}^*\,\bar{v} = \bar{u}^*\,\bar{u}\bar{v}^*\,\bar{v} = |u|^2|v|^2.$$
This is the case of equality in the Cauchy-Schwarz inequality, so $\bar{u}$ and $\bar{v}$ are equal up to a multiplicative constant. Since $\bar{u}$ and $\bar{v}$ are both unit vectors, they are in fact equal up to a sign. This means that, up to changing the sign of $\lambda$, we have $A = \lambda \bar{u}\bar{u}^*$. The uniqueness of $\bar{u}$ comes from the fact that the image of $A$ is spanned by the unit vector $\bar{u}$, and it is spanned by a unique unit vector up to a sign. To establish the uniqueness of $\lambda$ observe that $\lambda = \mathrm{tr}(A)$. This completes the proof.

**Exercise 4.2.** To show that the conditional law of $\bar{x}$ given $Y$ is the Gibbs measure (4.7), we need to prove that for all bounded and measurable functions $f : \mathbb{R}^N \to \mathbb{R}$ and $g : \mathbb{R}^{N \times N} \to \mathbb{R}$,

$$\mathbb{E}f(\bar{x})g(Y) = \mathbb{E}\frac{\int_{\mathbb{R}^N} f(x)\exp H_N^\circ(t,x)\,\mathrm{d}P_N(x)}{\int_{\mathbb{R}^N} \exp H_N^\circ(t,x)\,\mathrm{d}P_N(x)}\,g(Y). \tag{S.21}$$

It will be convenient to define the function $h : \mathbb{R}^N \times \mathbb{R}^{N^2}$ by

$$h(x,y) := \sqrt{\frac{2t}{N}}x\cdot yx - \frac{t}{N}|x|^4.$$

Conditionally on the randomness of $\bar{x}$, the random variable $Y$ is Gaussian with density

$$\frac{1}{(2\pi)^{N^2/2}}\exp\left(-\frac{1}{2}\left|y - \sqrt{\frac{2t}{N}}\bar{x}\bar{x}^*\right|^2\right)$$

$$= \frac{1}{(2\pi)^{N^2/2}}\exp\left(h(x,y) - \frac{1}{2}|y|^2\right). \tag{S.22}$$

It follows that the right side of (S.21) is given by

$$\int_{\mathbb{R}^N}\int_{\mathbb{R}^N}\int_{\mathbb{R}^{N^2}} \frac{f(x')g(y)\exp(h(x',y)+h(x,y))}{(2\pi)^{N^2/2}\int_{\mathbb{R}^N}\exp h(x'',y)\,\mathrm{d}P_N(x'')}\,e^{-\frac{1}{2}|y|^2}\,\mathrm{d}y\,\mathrm{d}P_N(x')\,\mathrm{d}P_N(x).$$

The change of variables $x = x'$ reveals that this is equal to

$$\int_{\mathbb{R}^N}\int_{\mathbb{R}^N}\int_{\mathbb{R}^{N^2}} \frac{f(x)g(y)\exp(h(x',y))\exp(h(x,y))}{(2\pi)^{N^2/2}\int_{\mathbb{R}^N}\exp h(x,y)\,d P_N(x)} e^{-\frac{1}{2}|y|^2}\,dy\,dP_N(x')\,dP_N(x)$$

$$= \frac{1}{(2\pi)^{N^2/2}}\int_{\mathbb{R}^N}\int_{\mathbb{R}^{N^2}} f(x)g(y)\exp\left(h(x,y) - \frac{1}{2}|y|^2\right) dP_N(x)\,dy.$$

Remembering (S.22) establishes (S.21) and completes the proof.

**Exercise 4.3.** A direct computation leveraging the equality (S.22) reveals that for any two Borel sets $A \subseteq \mathbb{R}^{N^2}$ and $B \subseteq \mathbb{R}^N$,

$$\mathbb{P}\{Y \in A, \bar{x} \in B\} = \frac{1}{(2\pi)^{N^2/2}}\int_B\int_A \exp\left(-\frac{1}{2}\left|y - \sqrt{\frac{2t}{N}}\bar{x}\bar{x}^*\right|^2\right) dy\,d\mathbb{P}_{\bar{x}}(x)$$

$$= \frac{1}{(2\pi)^{N^2/2}}\int_B\int_A \exp\left(\sqrt{\frac{2t}{N}} x\cdot yx - \frac{t}{N}|x|^4\right) e^{-\frac{1}{2}|y|^2}\,dy\,d\mathbb{P}_{\bar{x}}(x)$$

$$= \int_A\int_B \exp\left(\sqrt{\frac{2t}{N}} x\cdot yx - \frac{t}{N}|x|^4\right) d\mathbb{P}_{\bar{x}}(x)\,d\mathbb{P}_W(y). \qquad (S.23)$$

We now treat each question separately.

(i) Given two Borel sets $A \subseteq \mathbb{R}^{N^2}$ and $B \subseteq \mathbb{R}^N$, the definition of the conditional expectation implies that

$$\mathbb{P}\{(\bar{x}, Y) \in A \times B\} = \int_B\int_A d\mathbb{P}_{\bar{x}|Y}\,d\mathbb{P}_Y = \int_B\int_A \frac{d\mathbb{P}_{\bar{x}|Y}}{d\mathbb{P}_{\bar{x}}}\,d\mathbb{P}_{\bar{x}}\,d\mathbb{P}_Y.$$

This means that

$$\frac{d\mathbb{P}_{\bar{x},Y}}{d\mathbb{P}_{\bar{x}}\otimes\mathbb{P}_Y} = \frac{d\mathbb{P}_{\bar{x}|Y}}{d\mathbb{P}_{\bar{x}}},$$

and therefore

$$\mathsf{I}(\bar{x};Y) = \int_{\mathbb{R}^N\times\mathbb{R}^{N^2}} \log\left(\frac{d\mathbb{P}_{\bar{x}|Y}}{d\mathbb{P}_{\bar{x}}}\right) d\mathbb{P}_{\bar{x},Y} = \int_{\mathbb{R}^{N^2}}\int_{\mathbb{R}^N} \log\left(\frac{d\mathbb{P}_{\bar{x}|Y}}{d\mathbb{P}_{\bar{x}}}\right) d\mathbb{P}_{\bar{x}|Y}\,d\mathbb{P}_Y.$$

Remembering the definition of the mutual information (4.4) shows that

$$\mathsf{I}(\bar{x};Y) = \int_{\mathbb{R}^{N^2}} \mathsf{H}(\mathbb{P}_{\bar{x}|Y}\mid\mathbb{P}_{\bar{x}})\,d\mathbb{P}_Y = \mathbb{E}\mathsf{H}(\mathbb{P}_{\bar{x}|Y}\mid\mathbb{P}_{\bar{x}}) = N\mathsf{I}_N(t)$$

as desired.

(ii) By independence of $\bar{x}$ and $W$,

$$\frac{d\mathbb{P}_{\bar{x},Y}}{d\mathbb{P}_{\bar{x},W}} = \frac{d\mathbb{P}_{\bar{x}} \otimes \mathbb{P}_Y}{d\mathbb{P}_{\bar{x}} \otimes \mathbb{P}_W} \frac{d\mathbb{P}_{\bar{x},Y}}{d\mathbb{P}_{\bar{x}} \otimes \mathbb{P}_Y} = \frac{d\mathbb{P}_Y}{d\mathbb{P}_W} \frac{d\mathbb{P}_{\bar{x},Y}}{d\mathbb{P}_{\bar{x}} \otimes \mathbb{P}_Y}.$$

It follows that

$$H(\mathbb{P}_{\bar{x},Y} \mid \mathbb{P}_{\bar{x},W}) = \int \log\left(\frac{d\mathbb{P}_Y}{d\mathbb{P}_W}\right) d\mathbb{P}_{\bar{x},Y} + \int \log\left(\frac{d\mathbb{P}_{\bar{x},Y}}{d\mathbb{P}_{\bar{x}} \otimes \mathbb{P}_Y}\right) d\mathbb{P}_{\bar{x},Y}$$

$$= \int \log\left(\frac{d\mathbb{P}_Y}{d\mathbb{P}_W}\right) d\mathbb{P}_{\bar{x},Y} + H(\mathbb{P}_{\bar{x},Y} \mid \mathbb{P}_{\bar{x}} \otimes \mathbb{P}_Y).$$

Remembering the definition of the mutual information in terms of the relative entropy and observing that the term $\frac{d\mathbb{P}_Y}{d\mathbb{P}_W}$ in the first integral is independent of the coordinate associated with the law of $\bar{x}$ simplifies this to

$$H(\mathbb{P}_{\bar{x},Y} \mid \mathbb{P}_{\bar{x},W}) = H(\mathbb{P}_Y \mid \mathbb{P}_W) + I(\bar{x}; Y)$$

as required.

(iii) Applying (S.23) with $B = \mathbb{R}^N$ shows that

$$\frac{d\mathbb{P}_Y}{d\mathbb{P}_W}(y) = \int_{\mathbb{R}^N} \exp\left(\sqrt{\frac{2t}{N}} x \cdot yx - \frac{t}{N}|x|^4\right) d\mathbb{P}_{\bar{x}}(x).$$

It follows that

$$\overline{F}_N^\circ(t) = \frac{1}{N} \mathbb{E} \log \int_{\mathbb{R}^N} \exp\left(\sqrt{\frac{2t}{N}} x \cdot Yx - \frac{t}{N}|x|^4\right) d\mathbb{P}_{\bar{x}}(x)$$

$$= \frac{1}{N} \int_{\mathbb{R}^{N^2}} \log \int_{\mathbb{R}^N} \exp\left(\sqrt{\frac{2t}{N}} x \cdot yx - \frac{t}{N}|x|^4\right) d\mathbb{P}_{\bar{x}}(x) d\mathbb{P}_Y(y)$$

$$= \frac{1}{N} \int_{\mathbb{R}^{N^2}} \log\left(\frac{d\mathbb{P}_Y}{d\mathbb{P}_W}(y)\right) d\mathbb{P}_Y(y)$$

$$= \frac{1}{N} H(\mathbb{P}_Y \mid \mathbb{P}_W).$$

Similarly, (S.23) implies that

$$\frac{d\mathbb{P}_{\bar{x},Y}}{d\mathbb{P}_{\bar{x},W}} = \exp\left(\sqrt{\frac{2t}{N}} x \cdot yx - \frac{t}{N}|x|^4\right),$$

so we also have

$$H(\mathbb{P}_{\bar{x},Y} \mid \mathbb{P}_{\bar{x},W})$$

$$= \int_{\mathbb{R}^{N^2}} \int_{\mathbb{R}^N} \exp\left(\sqrt{\frac{2t}{N}} x \cdot yx - \frac{t}{N}|x|^4\right)\left(\sqrt{\frac{2t}{N}} x \cdot yx - \frac{t}{N}|x|^4\right) d\mathbb{P}_{\bar{x}}(x) \, d\mathbb{P}_W(y)$$

$$= \int_{\mathbb{R}^N} \int_{\mathbb{R}^{N^2}} \left(\sqrt{\frac{2t}{N}} x \cdot yx - \frac{t}{N}|x|^4\right) d\mathbb{P}_{Y|\bar{x}}(y) \, d\mathbb{P}_{\bar{x}}(x)$$

$$= \mathbb{E}_{\bar{x}} \mathbb{E}_{Y|\bar{x}} \left(\sqrt{\frac{2t}{N}} \bar{x} \cdot Y\bar{x} - \frac{t}{N}|\bar{x}|^4\right),$$

where the second equality uses (S.22). Since the conditional law of $Y$ given $\bar{x}$ is Gaussian with mean $\sqrt{\frac{2t}{N}}\bar{x}\bar{x}^*$, conditionally on the randomness of $\bar{x}$, we have

$$\mathbb{E}_{Y|\bar{x}} \sqrt{\frac{2t}{N}} \bar{x} \cdot Y\bar{x} = \frac{2t}{N} \bar{x} \cdot \bar{x}\bar{x}^* \bar{x} = \frac{2t}{N}|\bar{x}|^4,$$

and therefore

$$H(\mathbb{P}_{\bar{x},Y} \mid \mathbb{P}_{\bar{x},W}) = \frac{t}{N} \mathbb{E}|\bar{x}|^4$$

as desired.

(iv) Combining the three previous parts shows that

$$I_N(t) = \frac{1}{N} I(\bar{x}; Y) = \frac{1}{N} H(\mathbb{P}_{\bar{x},Y} \mid \mathbb{P}_{\bar{x},W}) - \frac{1}{N} H(\mathbb{P}_Y \mid \mathbb{P}_W) = \frac{t}{N^2} \mathbb{E}|x|^4 - \overline{F}_N^\circ(t)$$

which agrees with (4.18). This completes the proof.

**Exercise 4.4.** Let $A \in \mathbb{R}^{d \times d}$ be the positive semi-definite matrix such that

$$A^2 = (\mathbb{E} g_i g_j)_{1 \leq i,j \leq d}.$$

We denote the rows of $A$ by $a_1, \ldots, a_d \in \mathbb{R}^d$. Letting $z$ be a standard Gaussian vector on $\mathbb{R}^d$, we have $g \stackrel{d}{=} Az$, where $\stackrel{d}{=}$ denotes equality in distribution. In particular,

$$\max_{1 \leq i \leq d} g_i \stackrel{d}{=} \max_{1 \leq i \leq d} (Az)_i = \max_{1 \leq i \leq d} a_i \cdot z = F(z)$$

for the function $F : \mathbb{R}^d \to \mathbb{R}$ defined by $F(z) := \max_{1 \leq i \leq d} a_i \cdot z$. Observe that $F$ is Lipschitz continuous with $\|F\|_{\mathrm{Lip}} \leq \sqrt{a}$. Indeed, fix $z^1, z^2 \in \mathbb{R}^d$ with $F(z^1) \geq F(z^2)$ and let $1 \leq j \leq d$ be such that $a_j \cdot z^1 = F(z^1)$. The Cauchy-Schwarz inequality and the fact that $|a_j|^2 = \mathbb{E} g_j^2 \leq a$ imply that

$$|F(z^1) - F(z^2)| = a_j \cdot z^1 - \max_{1 \leq i \leq d} a_i \cdot z^2 \leq a_j \cdot (z^1 - z^2) \leq \sqrt{a}|z^1 - z^2|.$$

Invoking the Gaussian concentration inequality (Theorem 4.7) completes the proof.

**Exercise 4.5.** We treat each question separately.

(i) Fix $\lambda > 0$ and $x \in \mathbb{R}^N$ of unit norm. Chebyshev's inequality reveals that

$$\mathbb{P}\{|Wx|^2 \geq aN\} \leq \exp(-\lambda aN)\mathbb{E}\exp(\lambda|Wx|^2). \tag{S.24}$$

To bound this further, observe that the random variables $((Wx)_i)_{i \leq N}$ are independent centred Gaussian with variance

$$\mathbb{E}(Wx)_i^2 = \mathbb{E}\sum_{j,k=1}^N W_{ij}W_{ik}x_jx_k = \sum_{j=1}^N x_j^2 = 1.$$

It follows by independence that

$$\mathbb{E}\exp(\lambda|Wx|^2) = \mathbb{E}\prod_{i=1}^N \exp(\lambda(Wx)_i^2) = \left(\mathbb{E}\exp(\lambda Z^2)\right)^N$$

for a standard Gaussian $Z$. A change of variables shows that for $\lambda < 1/2$,

$$\mathbb{E}\exp(\lambda Z^2) = \frac{1}{\sqrt{2\pi}}\int_{\mathbb{R}} e^{-x^2(\frac{1}{2}-\lambda)}\,dx = \frac{1}{\sqrt{1-2\lambda}}\frac{1}{\sqrt{2\pi}}\int_{\mathbb{R}} e^{-\frac{y^2}{2}}\,dy = \frac{1}{\sqrt{1-2\lambda}}.$$

This means that $\mathbb{E}\exp(\lambda|Wx|^2) \leq \exp(CN)$ for some $C > 2$. Substituting this into (S.24) establishes (4.62).

(ii) Given $x$ with $|x| \leq 1$, let $y \in A$ be such that $|x - y| \leq 1/2$. By the triangle inequality,

$$|Wx| \leq \|W\|_* |x - y| + |Wy| \leq \frac{1}{2}\|W\|_* + \sup_{y \in A}|Wy|.$$

Taking the supremum in $x$ and rearranging gives the required bound.

(iii) To define the set $A$ we will leverage a standard upper bound (see [259, Corollary 4.2.13]) on the smallest number of balls of radius $\varepsilon > 0$ required to cover the unit ball. In general, the smallest number of balls of radius $\varepsilon > 0$ required to cover a compact set $K$ is known as the *covering number* of the set $K$, and it is denoted by $\mathcal{N}(K,\varepsilon)$. A related quantity is the packing number of a compact set $K$. A subset $\mathcal{N}$ of $K$ is $\varepsilon$-*separated* if $|x - y| > \varepsilon$ for all distinct points $x, y \in \mathcal{N}$. The *packing number* of $K$ is the size of the largest possible $\varepsilon$-separated subset of $K$, and it is denoted by $\mathcal{P}(K,\varepsilon)$. By maximality, it is clear that the balls of radius $\varepsilon > 0$ centred at the points in an $\varepsilon$-separated subset of $K$ cover the set $K$. In particular, the covering number is bounded by the packing number,

$$\mathcal{N}(K,\varepsilon) \leq \mathcal{P}(K,\varepsilon).$$

We now use this observation to find an upper bound on the covering number of the unit ball $B_1(0) \subseteq \mathbb{R}^N$. Denote by $m$ Lebesgue measure on $\mathbb{R}^N$, fix $\varepsilon > 0$ and let $P$ be a maximal $\varepsilon$-separated subset of $B_1(0)$. Since $|x-y| > \varepsilon$ for all distinct points $x,y \in P$, the balls $(B_{\varepsilon/2}(x))_{x \in P}$ are disjoint and contained in the expanded unit ball $B_{1+\varepsilon/2}(0)$. It follows by the additivity and scaling properties of Lebesgue measure that

$$\mathcal{P}(B_1(0), \varepsilon)\left(\frac{\varepsilon}{2}\right)^N m(B_1(0)) = \sum_{x \in P} m(B_{\varepsilon/2}(x)) \leq \left(1+\frac{\varepsilon}{2}\right)^N m(B_0(1)).$$

Rearranging reveals that

$$\mathcal{N}(B_1(0), \varepsilon) \leq \mathcal{P}(B_1(0), \varepsilon) \leq \left(1+\frac{2}{\varepsilon}\right)^N. \tag{S.25}$$

Applying (S.25) with $\varepsilon = 1/2$ and invoking (ii) allows us to choose a set $A$ of size $5^N$ with the property that

$$\|W\|_* \leq 2\sup_{y \in A} |Wy|$$

Together with (i) and the union bound this implies that

$$\mathbb{P}\{\|W\|_*^2 \geq aN\} \leq \mathbb{P}\{4\sup_{y \in A}|Wy|^2 \geq aN\} \leq 5^N \exp\left(\left(C-\frac{a}{C}\right)N\right).$$

Redefining the constant $C$ establishes (4.61).

(iv) The layer-cake representation and the change of variables $\lambda = \lambda' N^{q/2}$ give

$$\mathbb{E}\|W\|_*^q = \int_0^{+\infty} \mathbb{P}\{\|W\|_*^2 > \lambda^{2/q}\} d\lambda = N^{q/2} \int_0^{+\infty} \mathbb{P}\{\|W\|_*^2 \geq (\lambda')^{2/q} N\} d\lambda'.$$

Applying (4.61) and redefining the constant $C$ completes the proof.

**Exercise 4.6.** Taylor's theorem, the assumption that $\mathbb{E}X = 0$, and the mean value theorem imply that

$$|\mathbb{E}XF(X) - \mathbb{E}X^2 \mathbb{E}F'(X)| \leq |\mathbb{E}X(F(0)+F'(0)X) - \mathbb{E}X^2 \mathbb{E}F'(X)| + \frac{1}{2}\|F''\|_\infty \mathbb{E}|X|^3$$

$$\leq \mathbb{E}X^2 |F'(X) - F'(0)| + \frac{1}{2}\|F''\|_\infty \mathbb{E}|X|^3$$

$$\leq \mathbb{E}|X|^3 \|F''\|_\infty + \frac{1}{2}\|F''\|_\infty \mathbb{E}|X|^3.$$

This completes the proof.

**Exercise 4.7.** Denote by $\mathbb{E}_i$ the average with respect to the randomness of $X_{i+1},\ldots,X_d$ and by $\mathbb{E}^i$ the average with respect to the randomness of $X_i$. Let $f_i := \mathbb{E}_i f - \mathbb{E}_{i-1} f$ with the convention that $\mathbb{E}_0 f = \mathbb{E} f$ and $\mathbb{E}_d f = f$. Observe that for $1 \leqslant i < j \leqslant d$,

$$\mathbb{E} f_i f_j = \mathbb{E} f_i \mathbb{E}^j f_j = 0.$$

Combining this with the fact that $f - \mathbb{E} f = \sum_{i=1}^d f_i$ and Jensen's inequality reveals that

$$\mathbb{E}(f - \mathbb{E} f)^2 = \sum_{i=1}^d \mathbb{E} f_i^2 = \sum_{i=1}^d \mathbb{E}\big(\mathbb{E}_i(f - \mathbb{E}^i f)\big)^2 \leqslant \sum_{i=1}^d \mathbb{E}(f - \mathbb{E}^i f)^2. \qquad (\text{S.26})$$

At this point, introduce the function

$$f'_i := f(X_1,\ldots,X_{i-1},X'_i,X_{i+1},\ldots,X_d).$$

Since $X'_i$ is an independent copy of $X_i$,

$$\mathbb{E}^i(f - \mathbb{E}^i f)^2 = \mathbb{E}^i f^2 - (\mathbb{E}^i f)^2 = \mathbb{E}^i f^2 - \mathbb{E}^i f \mathbb{E}^i f'_i = \mathbb{E}^i f^2 - \mathbb{E}^i f f'_i = \frac{1}{2}\mathbb{E}^i (f - f'_i)^2.$$

Substituting this into (S.26) completes the proof.

**Exercise 4.8.** We first give a proof that leverages the Efron-Stein inequality obtained in the previous exercise, and then sketch an alternative argument that is more direct.

For each integer $N \geqslant 1$, let $\varepsilon_1,\ldots,\varepsilon_N$ be i.i.d. Rademacher random variables. Consider the normalized partial sum

$$S_N := \frac{1}{\sqrt{N}} \sum_{i=1}^N \varepsilon_i,$$

and recall that $(S_N)_{N \geqslant 1}$ converges to $Z$ in distribution by the central limit theorem. Denote by $\varepsilon'_1,\ldots,\varepsilon'_N$ another set of i.i.d. Rademacher random variables, also independent of $\varepsilon_1,\ldots,\varepsilon_N$, and for each $1 \leqslant i \leqslant N$ consider the normalized partial sum

$$S_N^i := \frac{1}{\sqrt{N}} \Big(\sum_{j \neq i} \varepsilon_j + \varepsilon'_i\Big) = S_N + \frac{\varepsilon'_i - \varepsilon_i}{\sqrt{N}}.$$

If we denote by $\mathbb{E}_i$ the average with respect to the randomness of the vector $(\varepsilon_i, \varepsilon'_i)$, then the Efron-Stein inequality implies that

$$\mathrm{Var} f(S_N) \leqslant \frac{1}{2} \sum_{i=1}^N \mathbb{E}(f(S_N) - f(S_N^i))^2 = \frac{1}{2} \sum_{i=1}^N \mathbb{E}\mathbb{E}_i(f(S_N) - f(S_N^i))^2.$$

Notice that $f(S_N) - f(S_N^i)$ is non-zero if and only if $\varepsilon_i \neq \varepsilon_i'$, so

$$\mathbb{E}_i(f(S_N) - f(S_N^i))^2 = \frac{1}{2}\left(f\left(S_N + \frac{1-\varepsilon_i}{\sqrt{N}}\right) - f\left(S_N + \frac{-1-\varepsilon_i}{\sqrt{N}}\right)\right)^2.$$

It follows by Taylor's theorem that

$$\begin{aligned}
\operatorname{Var} f(S_N) &\leq \frac{1}{4}\sum_{i=1}^{N}\mathbb{E}\left(f\left(S_N + \frac{1-\varepsilon_i}{\sqrt{N}}\right) - f\left(S_N + \frac{-1-\varepsilon_i}{\sqrt{N}}\right)\right)^2 \\
&= \frac{1}{4}\sum_{i=1}^{N}\mathbb{E}\left(\frac{2}{\sqrt{N}}f'(S_N) + \mathcal{O}\left(\frac{\|f''\|_\infty}{N}\right)\right)^2 \\
&= \frac{1}{4}\sum_{i=1}^{N}\left(\frac{4}{N}\mathbb{E}f'(S_N)^2 + \mathcal{O}\left(\frac{\|f''\|_\infty \|f'\|_\infty}{N^{3/2}}\right)\right) \\
&= \mathbb{E}f'(S_N)^2 + \mathcal{O}\left(\frac{\|f''\|_\infty \|f'\|_\infty}{N^{1/2}}\right).
\end{aligned}$$

Letting $N$ tend to infinity and using the fact that $(S_N)_{N\geq 1}$ converges to $Z$ in distribution completes the proof.

We now briefly sketch an alternative proof of the Gaussian Poincaré inequality. For every $x \in \mathbb{R}$ and $t > 0$, we denote the heat kernel evaluated at $(t,x)$ by

$$\Phi_t(x) := \frac{1}{\sqrt{2\pi t}}\exp\left(-\frac{x^2}{2t}\right).$$

We recall that $\partial_t \Phi = \frac{1}{2}\partial_x^2 \Phi$, and we write $\star$ for the convolution operator (acting only on the space variable). A direct calculation and an application of Jensen's inequality give that, for every $t \in (0,1]$,

$$\partial_t\left(\Phi_t \star (\Phi_{1-t} \star f)^2\right) = \Phi_t \star \left(\Phi_{1-t} \star f'\right)^2 \leq \Phi_1 \star (f')^2.$$

Integrating this inequality over $t \in (0,1]$ and evaluating the result at $x = 0$, we obtain another proof of the Gaussian Poincaré inequality.

**Exercise 4.9.** Let $f$ be the function defined in (4.110), and denote by

$$t_c := \sup\{t \geq 0 \mid f(t) = 0\} \in \mathbb{R}_{\geq 0} \cup \{+\infty\}$$

its largest zero. Notice that $\psi(0) = 0$, so $f$ is non-negative and $f(0) = 0$. By Proposition 4.10 and Theorem 4.9, the function $f$ is convex, so we have $f(t) > 0$ for every $t > t_c$. By Proposition 4.15, and since $f$ is differentiable almost everywhere, it suffices to verify that $t_c \in \mathbb{R}_{>0}$. From the definition of $\psi$ in (4.85) or (4.89), we see that $\psi$ is a smooth function and satisfies $\psi(0) = 0$, while it follows from (4.82) and the Nishimori identity that

$$\partial_h \psi(0) = \mathbb{E}\langle x \cdot x' \rangle = (\mathbb{E}\bar{x}_1)^2 = 0. \quad (S.27)$$

The Gibbs measure above is with $t = 0$, $h = 0$, and $N = 1$; in this case, this Gibbs measure is simply the measure $P_1$, and $x'$ is an independent sample from $P_1$. Recalling also from (4.82) that $\psi$ is Lipschitz continuous, we deduce that there exists $C < +\infty$ such that for every $h \in \mathbb{R}$,

$$\psi(h) \leq Ch^2.$$

For every $t < C/4$, the supremum in (4.110) is therefore achieved at $h = 0$, and thus $t_c$ is strictly positive. Since we assume that $P_1$ is not the Dirac mass at 0, we see from (4.95) that $\partial_h^2 \psi(0) > 0$, and in particular there exists $h > 0$ such that $\psi(h) > 0$. It is therefore clear from (4.110) that $f(t) > 0$ if $t$ is sufficiently large, and thus $t_c$ is finite.

**Exercise 4.10.** A direct computation together with (4.73), the Nishimori identity and the derivative computation (4.82) imply that

$$\mathbb{E}|\bar{x} - \mathbb{E}[\bar{x}|\mathcal{Y}]|^2 = \mathbb{E}|\bar{x}|^2 - 2\mathbb{E}\bar{x} \cdot \mathbb{E}[\bar{x}|\mathcal{Y}] + \mathbb{E}|\mathbb{E}[\bar{x}|\mathcal{Y}]|^2$$
$$= N\mathbb{E}\bar{x}_1^2 - 2\mathbb{E}\langle \bar{x} \cdot x \rangle + \mathbb{E}\langle x \cdot x' \rangle$$
$$= N\mathbb{E}\bar{x}_1^2 - N\partial_h \overline{F}_N(t,h).$$

It follows by Propositions 4.10, 2.11 and 2.15 that at the interior point of differentiability $(t,h) \in \mathbb{R}_{>0} \times \mathbb{R}_{>0}$,

$$\lim_{N \to +\infty} \frac{1}{N} \mathbb{E}|\bar{x} - \mathbb{E}[\bar{x}|\mathcal{Y}]|^2 = \mathbb{E}\bar{x}_1^2 - \partial_h f(t,h).$$

Invoking the envelope theorem (Theorem 2.21) completes the proof.

**Exercise 4.11.** In this context, the prior $P_1$ is the law of a Bernoulli random variable with probability of success $p \in (0,1)$, so it will be convenient to write $\bar{\sigma} \in \Sigma_N$ and $\sigma \in \Sigma_N$ in place of $\bar{x}$ and $x$ respectively. In this notation, for any $t \geq 0$,

$$H_N^\circ(t,\sigma) = \sum_{i<j} \left( \sqrt{\frac{2t}{N}}(W_{ij} + W_{ji})\sigma_i\sigma_j + \frac{4t}{N}\sigma_i\sigma_j\bar{\sigma}_i\bar{\sigma}_j \right) + \sum_{i=1}^{N} \sqrt{\frac{2t}{N}} W_{ii} + 2t - tN.$$

Noticing that $\frac{W_{ij}+W_{ji}}{\sqrt{2}}$ is again a standard Gaussian, we have the equality in distribution jointly over $\sigma$,

$$H_N^\circ(t,\sigma) \stackrel{d}{=} \sum_{i<j} \left( \sqrt{\frac{4t}{N}} W_{ij}\sigma_i\sigma_j + \frac{4t}{N}\sigma_i\sigma_j\bar{\sigma}_i\bar{\sigma}_j \right) + \sum_{i=1}^{N} \sqrt{\frac{2t}{N}} W_{ii} + 2t - tN,$$

where $\stackrel{d}{=}$ denotes equality in distribution. Applying this with $t = \frac{\lambda}{4}$ reveals that

$$\overline{F}_N^\circ\left(\frac{\lambda}{4}\right) = \overline{F}_N^{\text{gauss}}(\lambda) + \frac{1}{N}\sum_{i=1}^N \sqrt{\frac{\lambda}{2N}}\mathbb{E}W_{ii} + \frac{\lambda}{2N} - \frac{\lambda}{4}.$$

Remembering that $\mathbb{E}W_{ii} = 0$ and rearranging completes the proof.

**Exercise 4.12.** By Theorem 4.25,

$$\lim_{N\to+\infty} \mathsf{I}_N = \inf_{h\geqslant 0}\left(\frac{\lambda}{4} + \frac{h^2}{\lambda} - \psi(\lambda)\right) \tag{S.28}$$

for the initial condition $\psi : \mathbb{R}_{\geqslant 0} \to \mathbb{R}$ defined by

$$\psi(h) := \mathbb{E}\log \int_{\Sigma_1} \exp\left(\sqrt{2h}\sigma z_1 + 2h\sigma\overline{\sigma}_1 - h\right)\mathrm{d}P_1(\sigma).$$

Since $P_1$ is the law of a symmetric Bernoulli random variable, the initial condition simplifies to

$$\psi(h) = \mathbb{E}\log\frac{1}{2}\left(\exp\left(\sqrt{2h}z_1 + 2h\overline{\sigma}_1\right) - \exp\left(-(\sqrt{2h}z_1 + 2h\overline{\sigma}_1)\right)\right) - h$$
$$= \mathbb{E}\log\cosh\left(\sqrt{2h}z_1 + 2h\overline{\sigma}_1\right) - h.$$

Substituting this into (S.28) completes the proof.

## S.5 Poisson point processes and extreme values

**Exercise 5.1.** By Step 2 of the proof of Lemma 5.3, there exists a sequence $(K_n)_{n \geq 1}$ of compact subsets of $S$ that covers $S$. The assumption that $\lambda$ is locally finite implies that for each $n \geq 1$, the set of indices $\{i \in I \mid x_i \in K_n\}$ is finite. The index set $I$ is the union of these sets as $n$ ranges over the natural numbers, and is therefore countable.

**Exercise 5.2.** We start by showing that for every open set $U \subseteq S$, the mapping $\lambda \mapsto \lambda(U)$ is measurable. Since $U$ is open, the sequence $(f_n)_{n \geq 1}$ of continuous functions defined by $f_n(s) := \min(1, n\,\mathrm{dist}(s, U^c))$ increases to $\mathbf{1}_U$. It follows by the monotone convergence theorem that

$$\lambda(U) = \int_S \mathbf{1}_U \, d\lambda = \lim_{n \to +\infty} \int_S f_n \, d\lambda.$$

Each of the mappings $\lambda \mapsto \int_S f_n \, d\lambda$ is continuous, and therefore measurable with respect to the Borel $\sigma$-algebra $\mathcal{B}(\mathcal{M}_\delta(S))$. As a limit of measurable functions, the mapping $\lambda \mapsto \lambda(U)$ is itself measurable with respect to the Borel $\sigma$-algebra $\mathcal{B}(\mathcal{M}_\delta(S))$. This means that the $\lambda$-system of measurable sets $A \subseteq S$ such that $\lambda \mapsto \lambda(A)$ is measurable contains the $\pi$-system of open sets. It follows by the Dynkin $\pi$-$\lambda$ theorem that each of the mappings $\lambda \mapsto \lambda(A)$ is measurable with respect to the Borel $\sigma$-algebra $\mathcal{B}(\mathcal{M}_\delta(S))$. To see that this is the smallest $\sigma$-algebra for which this property holds, it suffices to approximate any function $f \in C_c(S; \mathbb{R})$ by a sequence of measurable functions that take a finite number of values.

**Exercise 5.3.** For a direct argument, we can reproduce the proof of the converse implication in Lemma 5.3 — alternatively, we can also appeal to this lemma with the set $M = \{\lambda_n \mid n \in \mathbb{N}\}$.

**Exercise 5.4.** It is clear that $\rho$ is non-negative and symmetric with $\rho(\lambda, \lambda) = 0$ for all $\lambda \in \mathcal{M}_\delta(S)$. It therefore suffices to verify that $\rho$ satisfies the triangle inequality. Fix $\lambda, \lambda', \lambda'' \in \mathcal{M}_\delta(S)$, and define the function $\phi : \mathbb{R}_{\geq 0} \to \mathbb{R}_{\geq 0}$ by

$$\phi(x) := \frac{x}{1+x}.$$

Observe that $\phi$ is non-decreasing and concave on $\mathbb{R}_{\geq 0}$, and that it is such that

$$\rho(\lambda, \lambda'') = \sum_{n=0}^{+\infty} 2^{-n} \phi\left(\left|\int_S f_n \, d(\lambda - \lambda'')\right|\right).$$

Since $\phi$ is non-decreasing on $\mathbb{R}_{\geq 0}$, we have

$$\rho(\lambda, \lambda'') \leq \sum_{n=0}^{+\infty} 2^{-n} \phi\left(\left|\int_S f_n \, d(\lambda - \lambda')\right| + \left|\int_S f_n \, d(\lambda' - \lambda'')\right|\right).$$

To bound this further, observe that by concavity of $\phi$, for any $a > 0$ and $b \geqslant 0$, we have
$$\frac{\phi(a+b) - \phi(b)}{a} \leqslant \frac{\phi(a) - \phi(0)}{a},$$
and therefore $\phi(a+b) \leqslant \phi(a) + \phi(b)$. This inequality holds trivially for $a = 0$. It follows that
$$\rho(\lambda, \lambda'') \leqslant \sum_{n=0}^{+\infty} 2^{-n} \phi\left(\left|\int_S f_n \, d(\lambda - \lambda')\right|\right) + \sum_{n=0}^{+\infty} 2^{-n} \phi\left(\left|\int_S f_n \, d(\lambda' - \lambda'')\right|\right)$$
$$= \rho(\lambda, \lambda') + \rho(\lambda', \lambda'').$$

This completes the proof.

**Exercise 5.5.** Let $(\lambda_n)_{n \geqslant 1}$ be a Cauchy sequence in $\mathcal{M}_\delta(S)$. Given $\varepsilon \in (0, 1)$ and $k \geqslant 1$, it is possible to find $n, m \geqslant 1$ large enough so that
$$\frac{\left|\int_S f_k \, d(\lambda_n - \lambda_m)\right|}{1 + \left|\int_S f_k \, d(\lambda_n - \lambda_m)\right|} \leqslant \varepsilon.$$

Rearranging reveals that
$$\left|\int_S f_k \, d(\lambda_n - \lambda_m)\right| \leqslant \frac{\varepsilon}{1 - \varepsilon}.$$

This implies that for every $f \in \mathcal{F}$, the sequence $\left(\int_S f \, d\lambda_n\right)_{n \geqslant 1}$ is Cauchy in $\mathbb{R}$. We now use a density argument to show that this in fact true for every $f \in C_c(S; \mathbb{R})$. Fix $f \in C_c(S; \mathbb{R})$ as well as $\varepsilon > 0$, and denote by $K$ the compact support of $f$. Invoke Lemma 5.1 to find $g \in \mathcal{F}$ with $|f(x) - g(x)| \leqslant \varepsilon$ for all $x \in S$ and $\operatorname{supp} g \subseteq K_\varepsilon$. Since the sequence $\left(\int_S g \, d\lambda_n\right)_{n \geqslant 1}$ is Cauchy, there exist $n, m \geqslant 1$ large enough with
$$\left|\int_S f \, d(\lambda_n - \lambda_m)\right| = \left|\int_S g \, d(\lambda_n - \lambda_m)\right| + \left|\int_S (f - g) \, d(\lambda_n - \lambda_m)\right| \leqslant \varepsilon + 2\varepsilon \sup_{k \geqslant 1} \lambda_k(K_\varepsilon).$$

Arguing as in the part of the proof of Theorem 5.2 leading to (5.8) reveals that for $\varepsilon$ small enough, we have $\sup_{k \geqslant 1} \lambda_k(K_\varepsilon) < +\infty$. This shows that for each $f \in C_c(S; \mathbb{R})$, the sequence $\left(\int_S f \, d\lambda_n\right)_{n \geqslant 1}$ is Cauchy in $\mathbb{R}$, and therefore convergent and uniformly bounded. Invoking Lemma 5.3 completes the proof.

**Exercise 5.6.** By Proposition A.24 it suffices to show that, given a measurable set $A \subseteq S$ with $\mu(\partial A) = 0$ and $t \geqslant 0$, we have
$$\lim_{n \to +\infty} \mathbb{E} \exp(-t\Lambda_n(A)) = \mathbb{E} \exp(-t\Lambda(A)). \tag{S.29}$$

The proof very closely resembles that of the Portmanteau theorem. To alleviate notation, let $U := \operatorname{int}(A)$ and $F := \overline{A}$. Define the sequence $(f_m)_{m \geqslant 1}$ of bounded and continuous functions on $S$ by $f_m(s) := \min\bigl(1, m\,d(s, U^c)\bigr)$. Since $U^c$ is closed, the sequence $(f_m)_{m \geqslant 1}$ increases to $\mathbf{1}_U$. It follows that for every $n, m \geqslant 1$,

$$\mathbb{E}\exp\bigl(-t\Lambda_n(U)\bigr) \leqslant \mathbb{E}\exp\Bigl(-t\int_S f_m\,\mathrm{d}\Lambda_n\Bigr).$$

Letting $n$ tend to infinity and then $m$ tend to infinity shows that

$$\limsup_{n \to +\infty} \mathbb{E}\exp\bigl(-t\Lambda_n(U)\bigr) \leqslant \liminf_{m \to +\infty} \mathbb{E}\exp\Bigl(-t\int_S f_m\,\mathrm{d}\Lambda\Bigr) = \mathbb{E}\exp\bigl(-t\Lambda(U)\bigr).$$

To obtain a matching lower bound, define the sequence $(g_m)_{m \geqslant 1}$ of bounded and continuous functions on $S$ by $g_m(s) := \max\bigl(0, 1 - m\,d(x, F)\bigr)$. Since $F$ is closed, the sequence $(g_m)_{m \geqslant 1}$ decreases to $\mathbf{1}_F$. It follows that for every $n, m \geqslant 1$,

$$\mathbb{E}\exp\bigl(-t\Lambda_n(F)\bigr) \geqslant \mathbb{E}\exp\Bigl(-t\int_S g_m\,\mathrm{d}\Lambda_n\Bigr).$$

Letting $n$ tend to infinity and then $m$ tend to infinity shows that

$$\liminf_{n \to +\infty} \mathbb{E}\exp\bigl(-t\Lambda_n(F)\bigr) \geqslant \limsup_{m \to +\infty} \mathbb{E}\exp\Bigl(-t\int_S g_m\,\mathrm{d}\Lambda\Bigr) = \mathbb{E}\exp\bigl(-t\Lambda(F)\bigr).$$

Combining these two asymptotic bounds reveals that

$$\mathbb{E}\exp\bigl(-t\Lambda(F)\bigr) \leqslant \liminf_{n \to +\infty} \mathbb{E}\exp\bigl(-t\Lambda_n(A)\bigr) \leqslant \limsup_{n \to +\infty} \mathbb{E}\exp\bigl(-t\Lambda_n(A)\bigr)$$
$$\leqslant \mathbb{E}\exp\bigl(-t\Lambda(U)\bigr).$$

Recall that $\Lambda(F)$ is a Poisson random variable with mean $\mu(F)$ while $\Lambda(U)$ is a Poisson random variable with mean $\mu(U)$. Since $\mu(\partial A) = 0$, we have that $\mu(F) = \mu(A) = \mu(U)$, so these two random variables are in fact equal in distribution. The above string of inequalities is therefore a string of equalities. This establishes (S.29) and completes the proof.

**Exercise 5.7.** By Proposition 5.7, it suffices to investigate the convergence of the Laplace transform of $\Lambda_n$, and by Remark 5.8 it suffices to verify (5.17) for a given smooth and compactly supported function $g \in C_c^\infty(\mathbb{R}_{\geqslant 0} \times \mathbb{R}_{>0}; \mathbb{R}_{\geqslant 0})$. The independence of $(X_k)_{1 \leqslant k \leqslant n}$ implies that

$$\mathbb{E}\exp\Bigl(-\int_{\mathbb{R}_{\geqslant 0} \times \mathbb{R}_{>0}} g\,\mathrm{d}\Lambda_n\Bigr) = \mathbb{E}\exp\Bigl(-\sum_{k=1}^n g\bigl(k/n, n^{-\frac{1}{\zeta}} X_k\bigr)\Bigr)$$
$$= \exp\Bigl(\sum_{k=1}^n \log \mathbb{E}\exp\bigl(-g(k/n, n^{-1/\zeta} X_1)\bigr)\Bigr).$$

Applying the argument leading to (5.31) in the proof of Proposition 5.13 to the function $f(x) := g(k/n, x)$ shows that for any $1 \leq k \leq n$,

$$\mathbb{E}\exp\left(-g(k/n, n^{-1/\zeta}X_1)\right) = 1 - \frac{c}{n}\int_0^{+\infty}\left(1 - e^{-g(k/n,x)}\right)\frac{\zeta}{x^{\zeta+1}}\,dx + o\left(\frac{1}{n}\right).$$

Moreover, the $o(1/n)$ error term can be controlled uniformly over $k$. We thus have

$$\mathbb{E}\exp\left(-\int_{\mathbb{R}_{\geq 0}\times\mathbb{R}_{>0}} g\,d\Lambda_n\right) = \exp\left(\sum_{k=1}^{n}\log\left(1 - \frac{c}{n}\int_0^{+\infty}\left(1 - e^{-g(k/n,x)}\right)\frac{\zeta}{x^{\zeta+1}}\,dx + o\left(\frac{1}{n}\right)\right)\right),$$

which together with a Taylor expansion of the logarithm reveals that

$$\mathbb{E}\exp\left(-\int_{\mathbb{R}_{\geq 0}\times\mathbb{R}_{>0}} g\,d\Lambda_n\right) = \exp\left(o(1) - \sum_{k=1}^{n}\frac{c}{n}\int_0^{+\infty}\left(1 - e^{-g(k/n,x)}\right)\frac{\zeta}{x^{\zeta+1}}\,dx\right).$$

Letting $n$ tend to infinity, identifying a Riemann sum, and invoking Propositions 5.4 and 5.9 completes the proof.

**Exercise 5.8.** We treat each question separately.

(i) By Proposition 5.12, we have

$$\mathbb{E}\int_{(0,1]} x\,d\Lambda(x) = \int_0^1 \zeta x^{-\zeta}\,dx,$$

which is finite since $\zeta < 1$. In particular we have that $\int_{(0,1]} x\,d\Lambda(x)$ is finite almost surely. Since $\int_1^{+\infty}\zeta x^{-(\zeta+1)}\,dx < +\infty$, there are finitely many points in $\Lambda$ that fall in the interval $[1,+\infty)$. Combining these two observations leads to the claim.

(ii) Since $X_1\mathbf{1}_{\{X_1\in[a,2a]\}} \leq 2a$,

$$\mathbb{E}X_1\mathbf{1}_{\{X_1\in[a,2a]\}} \leq 2a\mathbb{P}\{X_1 \in [a,2a]\} < 2a\mathbb{P}\{X_1 \geq a\}.$$

For $a$ large enough we obtain the upper bound

$$\mathbb{E}X_1\mathbf{1}_{\{X_1\in[a,2a]\}} \leq 2aa^{-\zeta} = 2a^{1-\zeta}$$

as required.

(iii) Fix $\varepsilon > 0$, let $a > 0$ be large enough so the bound in (ii) holds, and fix $n \geq 1$ large enough so $\varepsilon n^{1/\zeta} > a$. The interval $[0,\varepsilon n^{1/\zeta}]$ can be partitioned as

$$[0,\varepsilon n^{1/\zeta}] \subseteq [0,a] \cup \bigcup_{i=0}^{i^*-1}(2^i a, 2^{i+1}a]$$

for $i^* = \lceil \log_2(\varepsilon n^{1/\zeta}/a) \rceil$. It follows that

$$\frac{1}{n^{1/\zeta}} \sum_{k=1}^{n} \mathbb{E} X_k \mathbf{1}_{\{X_k \leq \varepsilon n^{1/\zeta}\}} \leq n^{1-1/\zeta}\left(\mathbb{E} X_1 \mathbf{1}_{\{X_1 \leq a\}} + \sum_{i=0}^{i^*-1} \mathbb{E} X_1 \mathbf{1}_{\{X_1 \in [2^i a, 2^{i+1} a]\}}\right)$$

$$\leq n^{1-1/\zeta}\left(a + 2a^{1-\zeta} \sum_{i=0}^{i^*-1} 2^{i(1-\zeta)}\right)$$

$$\leq n^{1-1/\zeta}\left(a + \frac{2\varepsilon^{1-\zeta} n^{1/\zeta - 1}}{2^{1-\zeta} - 1}\right),$$

and therefore,

$$\limsup_{n \to +\infty} \frac{1}{n^{1/\zeta}} \sum_{k=1}^{n} \mathbb{E} X_k \mathbf{1}_{\{X_k \leq \varepsilon n^{1/\zeta}\}} \leq \frac{2\varepsilon^{1-\zeta}}{2^{1-\zeta} - 1}.$$

Letting $\varepsilon$ tend to zero gives the desired asymptotic.

(iv) Fix $\varepsilon, \delta > 0$, and let $K$ be large enough so

$$\mathbb{P}\{u_K > \delta\} \leq \frac{\delta}{2}.$$

This is possible since

$$\mathbb{P}\{u_K > \delta\} = \mathbb{P}\{\Lambda[\delta, +\infty) > K\} = \mathbb{P}\{\Pi > K\}$$

for a Poisson random variable $\Pi$ with mean $\mu[\delta, +\infty) = \delta^{-\zeta}$. Using the convergence in law established in Proposition 5.15, find $n \geq K$ such that

$$\mathbb{P}\{X_{K,n} > \delta\} \leq \delta.$$

To alleviate notation introduce the sets

$$A := \left\{\sum_{k=K}^{n} X_{k,n} \geq \varepsilon n^{1/\zeta}\right\} \quad \text{and} \quad B := \{X_{K,n} \leq \delta n^{1/\zeta}\},$$

and observe that

$$\mathbb{P}(A) = \mathbb{P}(A \cap B) + \mathbb{P}(A \cap B^c) \leq \mathbb{P}(A \cap B) + \mathbb{P}(B^c) \leq \mathbb{P}(A \cap B) + \delta.$$

To bound this further we use a generalized version of Chebyshev's inequality,

$$\varepsilon n^{1/\zeta} \mathbb{P}(A \cap B) = \varepsilon n^{1/\zeta} \mathbb{E} \mathbf{1}_A \mathbf{1}_B \leq \sum_{k=K}^{n} \mathbb{E} X_{k,n} \mathbf{1}_A \mathbf{1}_B \leq \sum_{k=K}^{n} \mathbb{E} X_{k,n} \mathbf{1}_B,$$

which implies that

$$\mathbb{P}\left\{\sum_{k=K}^{n} X_{k,n} \geq \varepsilon n^{1/\zeta}\right\} \leq \frac{1}{\varepsilon n^{1/\zeta}} \sum_{k=K}^{n} \mathbb{E} X_{k,n} \mathbf{1}_B + \delta$$

$$\leq \frac{1}{\varepsilon n^{1/\zeta}} \sum_{k=K}^{n} \mathbb{E} X_{k,n} \mathbf{1}_{\{X_{k,n} \leq \delta n^{1/\zeta}\}} + \delta.$$

Recalling that $X_{k,n} \geq 0$ for $1 \leq k \leq n$, and using the previous part to let $n$ tend to infinity and then $\delta$ tend to zero establishes the claim.

(v) Introduce the sequences $(Y_{k',n})_{k',n \geqslant k}$ and $(Y_n)_{n \geqslant 1}$ of random elements defined by

$$Y_{k',n} := \left( \frac{X_{1,n}}{\sum_{\ell=1}^{k'} X_{\ell,n}}, \ldots, \frac{X_{k,n}}{\sum_{\ell=1}^{k'} X_{\ell,n}} \right) \quad \text{and} \quad Y_n := \left( \frac{X_{1,n}}{\sum_{\ell=1}^{n} X_{\ell,n}}, \ldots, \frac{X_{k,n}}{\sum_{\ell=1}^{n} X_{\ell,n}} \right)$$

respectively. Proposition 5.15 and the continuous mapping theorem (Exercise A.10) imply that, for each fixed $k' \geqslant k$, as $n$ tends to infinity, the sequence $(Y_{k',n})_{n \geqslant 1}$ converges in law to the random element

$$Z_{k'} := \left( \frac{u_1}{\sum_{\ell=1}^{k'} u_\ell}, \ldots, \frac{u_k}{\sum_{\ell=1}^{k'} u_\ell} \right).$$

In turn, as $k'$ tends to infinity, the sequence $(Z_{k'})_{k' \geqslant k}$ converges almost surely, and therefore in law, to the random element

$$Y := (v_1, \ldots, v_k).$$

By the union bound, for every $\varepsilon > 0$,

$$\mathbb{P}\{|Y_{k',n} - Y_n|_1 \geqslant \varepsilon\} \leqslant \sum_{j=1}^{k} \mathbb{P}\left\{ \left| \frac{X_{j,n}}{\sum_{\ell=1}^{k'} X_{\ell,n}} - \frac{X_{j,n}}{\sum_{\ell=1}^{n} X_{\ell,n}} \right| \geqslant \frac{\varepsilon}{k} \right\}$$

$$\leqslant \sum_{j=1}^{k} \mathbb{P}\left\{ \frac{n^{-1/\zeta} X_{j,n}}{n^{-1/\zeta} \sum_{\ell=1}^{k'} X_{\ell,n} \cdot n^{-1/\zeta} \sum_{\ell=1}^{n} X_{\ell,n}} \cdot n^{-1/\zeta} \sum_{\ell=k'}^{n} X_{\ell,n} \geqslant \frac{\varepsilon}{k} \right\}$$

$$\leqslant \sum_{j=1}^{k} \mathbb{P}\left\{ n^{-1/\zeta} \sum_{\ell=k'}^{n} X_{\ell,n} \geqslant \frac{\varepsilon}{k} \cdot n^{-1/\zeta} X_{1,n} \right\},$$

where the second inequality uses that $\left| \frac{1}{a} - \frac{1}{a+b} \right| = \frac{b}{a(a+b)}$ for every $a,b > 0$. To bound this further, fix $\delta > 0$, and let $\delta' > 0$ be small enough so

$$\mathbb{P}\{u_1 \geqslant \delta'\} = \mathbb{P}\{\Lambda[\delta', +\infty) > 1\} = \mathbb{P}\{\Pi > 1\} \geqslant 1 - \delta/2,$$

where $\Pi$ denotes a Poisson random variable with mean $\mu[\delta', +\infty) = (\delta')^{-\zeta}$. Using the convergence in law established in Proposition 5.15, find $n$ large enough such that

$$\mathbb{P}\{n^{-1/\zeta} X_{1,n} \geqslant \delta'\} \geqslant 1 - \delta,$$

and observe that

$$\mathbb{P}\{|Y_{k',n} - Y_n|_1 \geqslant \varepsilon\} \leqslant \sum_{j=1}^{k} \mathbb{P}\left\{ n^{-1/\zeta} \sum_{\ell=k'}^{n} X_{\ell,n} \geqslant \frac{\varepsilon \delta'}{k} \right\} + \delta.$$

It follows by (iv) that for every $\varepsilon > 0$,

$$\lim_{k' \to +\infty} \limsup_{n \to +\infty} \mathbb{P}\{|Y_{k',n} - Y_n|_1 \geqslant \varepsilon\} \leqslant \delta.$$

Letting $\delta$ tend to zero and invoking Exercise A.11 completes the proof.

**Exercise 5.9.** Applying equality (5.55) with $X_n$ replaced by $X_n Z_n$ reveals that

$$\mathbb{E}\log \sum_{n=1}^{+\infty} u_n X_n Z_n = \frac{1}{\zeta} \log \mathbb{E}(X_1 Z_1)^\zeta + \mathbb{E}\log \sum_{n=1}^{+\infty} u_n$$

$$= \frac{1}{\zeta} \log \mathbb{E} X_1^\zeta + \frac{1}{\zeta} \log \mathbb{E} Z_1^\zeta + \mathbb{E}\log \sum_{n=1}^{+\infty} u_n.$$

Applying this equality again with $X_n Z_n$ replaced by $Z_n$ completes the proof.

**Exercise 5.10.** We treat each question separately.

(i) Fix $x = (x_j)_{j \geq 1} \in V$. Since $\sum_{j=1}^{+\infty} x_j \leq 1$, there must exist $m \geq 1$ with $x_m \leq m^{-1}$. It follows by non-decreasingness of the coordinates of $x$ that

$$\left| p_n((x_j)) - \sum_{j=1}^{m-1} x_j^n \right| = \sum_{j=m}^{+\infty} x_j^n \leq x_m^{n-1} \sum_{j=m}^{+\infty} x_j \leq x_m^{n-1} \leq m^{1-n}.$$

This shows that $p_n$ is continuous for $n \geq 2$ as the uniform limit of continuous functions. To see that $p_1$ is not continuous, consider the sequence $(x^n)_{n \geq 1} \subseteq V$ defined by

$$x_j^n := (\log(n))^{-1} j^{-1} \mathbf{1}_{\{2 \leq j \leq n\}}.$$

Notice that $x^n \in V$ as

$$\frac{1}{\log(n)} \sum_{j=2}^{n} \frac{1}{j} \leq \frac{1}{\log(n)} \int_1^n \frac{1}{x} dx = 1.$$

Moreover, $(x^n)_{n \geq 1}$ converges to the zero sequence but

$$p_1(x^n) = \frac{1}{\log(n)} \sum_{j=2}^{n} \frac{1}{j} \geq \frac{1}{\log(n)} \int_2^n \frac{1}{x} dx = \frac{\log(n) - \log(2)}{\log(n)}$$

is asymptotically lower bounded by 1 and can therefore not converge to 0. This shows that $p_1$ is not continuous.

(ii) The set $\mathcal{P}$ is closed under multiplication and linear combinations so it forms an algebra of continuous functions on the compact set $V$. By the real version of the Stone-Weierstrass theorem (Theorem A.10), it therefore suffices to show that $\mathcal{P}$ separates points. Fix distinct elements $x \neq y \in V$, and introduce the measures

$$\mu_x := \sum_{j \geq 1} x_j^2 \delta_{x_j} \quad \text{and} \quad \mu_y := \sum_{j \geq 1} y_j^2 \delta_{y_j}$$

on $[0, 1]$. Since the coordinates of $x$ are decreasing and sum to less than one, each non-zero $x_j$ has finite multiplicity in the sequence $(x_j)_{j \geq 1}$, and this multiplicity is given by $\mu_x(x_j)/x_j^2$. It follows that the measures $\mu_x$ and $\mu_y$ are different, so there must exist a polynomial $P:[0,1] \to \mathbb{R}$ with

$$\sum_{j \geq 1} x_j^2 P(x_j) = \int_0^1 P \, d\mu_x \neq \int_0^1 P \, d\mu_y = \sum_{j \geq 1} y_j^2 P(y_j).$$

Noticing that the function $z \mapsto \sum_{j \geq 1} z_j^2 P(z_j)$ belongs to the algebra $\mathcal{P}$ shows that $\mathcal{P}$ separates points as desired.

(iii) Fix $k \geq 1$ as well as $n_1, \ldots, n_k \geq 2$, and, through a slight abuse of notation, consider the function $f: V \to \mathbb{R}$ defined by $f(x) := \prod_{\ell=1}^k p_{n_\ell}(x)$. Observe that

$$S(n_1, \ldots, n_k) = \mathbb{E} f\big((v_j)_{j \geq 1}\big) = \sum_{j_1, \ldots, j_k \geq 1} v_{j_1}^{n_1} \cdots v_{j_k}^{n_k}.$$

To simplify this expression, we rewrite it as the average of some function $g$ of the overlaps $R^n = (R_{\ell, \ell'})_{\ell, \ell' \leq n}$ sampled from the average (5.51) associated with the weights $(v_j)_{j \geq 1}$. Let $n := n_1 + \cdots + n_k$, and partition the set $\{1, \ldots, n\}$ into $k$ subsets $\mathcal{I} = (I_j)_{j \leq k}$ with $|I_j| = n_j$ and $1 \in I_1$. For $\ell, \ell' \leq n$ write $\ell \sim \ell'$ if and only if $\ell$ and $\ell'$ belong to the same set in $\mathcal{I}$, and introduce the function

$$g(R^n) := \prod_{1 \leq \ell \sim \ell' \leq n} R_{\ell, \ell'} = \prod_{1 \leq \ell \sim \ell' \leq n} \mathbf{1}_{\{\sigma^\ell = \sigma^{\ell'}\}}.$$

By construction,

$$\mathbb{E}\langle g(R^n) \rangle = \mathbb{E} \sum_{j_1, \ldots, j_k \geq 1} v_{j_1}^{n_1} \cdots v_{j_k}^{n_k} = S(n_1, \ldots, n_k),$$

so the Ghirlanda-Guerra identities (5.57) imply that

$$\mathbb{E}\langle R_{1, n+1} g(R^n) \rangle = \frac{1}{n} S(n_1, \ldots, n_k) \mathbb{E}\langle R_{1,2} \rangle + \frac{1}{n} \sum_{\ell=2}^n \mathbb{E}\langle g(R^n) R_{1, \ell} \rangle.$$

A direct computation shows that

$$\mathbb{E}\langle g(R^n) R_{1, n+1} \rangle = \mathbb{E}\langle g(R^n) \mathbf{1}_{\{\sigma^1 = \sigma^{n+1}\}} \rangle = S(n_1 + 1, n_2, \ldots, n_k).$$

Similarly, for any index $1 \leq \ell \leq n$ with $\ell \in I_1$,

$$\mathbb{E}\langle g(R^n) R_{1, \ell} \rangle = \mathbb{E}\langle g(R^n) \mathbf{1}_{\{\sigma^1 = \sigma^\ell\}} \rangle = \mathbb{E}\langle g(R^n) \rangle = S(n_1, \ldots, n_k),$$

while for any index in $1 \leq \ell \leq n$ with $\ell \notin I_1$

$$\mathbb{E}\langle g(R^n) R_{1, \ell} \rangle = \mathbb{E}\langle g(R^n) \mathbf{1}_{\{\sigma^1 = \sigma^\ell\}} \rangle = S(n_2, \ldots, n_j + n_1, \ldots, n_k).$$

Indeed, multiplying $g$ by $\mathbf{1}_{\{\sigma^1=\sigma^l\}}$ has the effect of merging the set $I_\ell$ with the set $I_1$ in the sense that $g(R^n)\mathbf{1}_{\{\sigma^1=\sigma^l\}}$ is defined just as $g$ with the only difference that these two subsets are merged into one to form a new partition of $\{1,\ldots,n\}$ consisting of $k-1$ sets. Putting this all together reveals that

$$S(n_1+1, n_2, \ldots, n_k) = \frac{1}{n} S(n_1, \ldots, n_k) \mathbb{E}\langle R_{1,2}\rangle + \frac{1}{n} \sum_{\ell \in I_1 \setminus \{1\}} S(n_1, \ldots, n_k)$$
$$+ \frac{1}{n} \sum_{j=2}^{k} \sum_{\ell \in I_j} S(n_2, \ldots, n_j + n_1, \ldots, n_k).$$

Combining like terms and remembering that $R_{1,2}$ is a Bernoulli random variable with probability of success $1-\zeta$ by (5.58) yields

$$S(n_1+1, n_2, \ldots, n_k) = \frac{n_1 - \zeta}{n} S(n_1, \ldots, n_k) + \sum_{j=2}^{k} \frac{n_j}{n} S(n_2, \ldots, n_j + n_1, \ldots, n_k)$$

as required.

(iv) The sum of the coordinates in each term on the right side of the Talagrand identities has been reduced by one compared to the left side. A recursive application of the Talagrand identities therefore allows us to express all the quantities $S(n_1, \ldots, n_k)$ in terms of $S(2)$, which is equal to $1-\zeta$ by (5.58). Together with the density of $\mathcal{P}$ in the space of continuous functions on $V$, this means that for any continuous function $f : V \to \mathbb{R}$, the expectation $\mathbb{E} f((v_j)_{j \geq 1})$ is entirely determined by the Ghirlanda-Guerra identities (5.57) and the parameter $\zeta \in (0,1)$ determining the law (5.58) of the overlap $R_{1,2}$. Since the Poisson-Dirichlet process satisfies the same properties, the sequence $(v_j)_{j \geq 1}$ must have the same law as a Poisson-Dirichlet process with parameter $\zeta$. This completes the proof.

**Exercise 5.11.** We treat each question separately.

(i) To begin with, observe that for any $x \in (0,1)$,

$$F(F^{-1}(x)) = \lim_{n \to +\infty} F\left(F^{-1}(x) + \frac{1}{n}\right) \geq x,$$

where we have used the right-continuity of $F$. Now, suppose for the sake of contradiction that there exist $t \in \mathbb{R}$ and $x \in (0,1)$ with $F(t) < x$ but $t \geq F^{-1}(x)$. By non-decreasingness of $F$, we have

$$F(t) \geq F(F^{-1}(x)) \geq x > F(t)$$

which is a contradiction. Conversely, suppose that $t < F^{-1}(x)$ but $F(t) \geq x$. By definition of $F^{-1}$ we have $F^{-1}(x) \leq t$ which is a contradiction.

(ii) Given $t \in \mathbb{R}$ and $x \in (0,1)$, the conclusion of (i) may equivalently be expressed as the fact that $F^{-1}(x) \leq t$ if and only if $x \leq F(t)$. It follows that

$$\mathbb{P}\{F^{-1}(U) \leq t\} = \mathbb{P}\{U \leq F(t)\} = F(t).$$

Since the distribution function characterizes the law of a random variable, this completes the proof.

**Exercise 5.12.** By the dominated convergence theorem, the mapping

$$\zeta \mapsto \mathbb{E} \exp \zeta X$$

is differentiable on $(-\varepsilon, \varepsilon)$, with derivative given by

$$\zeta \mapsto \mathbb{E}(X \exp \zeta X).$$

The mapping

$$\zeta \mapsto \log \mathbb{E} \exp \zeta X$$

is therefore also differentiable on this interval, and its derivative at zero is $\mathbb{E}X$, as announced.

**Exercise 5.13.** We treat each question separately.

(i) We write

$$\mathbf{1}_{A_{n+1}} = \mathbf{1}_{A_n} \prod_{\ell=1}^{n} \mathbf{1}_{\{R_{\ell,n+1} \in A\}} = \mathbf{1}_{A_n} \prod_{\ell=1}^{n} \left(1 - \mathbf{1}_{\{R_{\ell,n+1} \notin A\}}\right),$$

and we clearly have

$$\prod_{\ell=1}^{n} \left(1 - \mathbf{1}_{\{R_{\ell,n+1} \notin A\}}\right) \geq 1 - \sum_{\ell=1}^{n} \mathbf{1}_{\{R_{\ell,n+1} \notin A\}},$$

yielding the announced inequality (5.111).

(ii) Using the result of the previous question, we can write

$$\mathbb{E}\langle \mathbf{1}_{A_{n+1}} \rangle \geq \mathbb{E}\langle \mathbf{1}_{A_n} \rangle - \sum_{\ell=1}^{n} \mathbb{E}\langle \mathbf{1}_{A_n} \mathbf{1}_{\{R_{\ell,n+1} \notin A\}} \rangle.$$

The Ghirlanda-Guerra identities imply that the sum above is $(1 - \zeta(A))\mathbb{E}\langle \mathbf{1}_{A_n} \rangle$, so

$$\mathbb{E}\langle \mathbf{1}_{A_{n+1}} \rangle \geq \zeta(A) \mathbb{E}\langle \mathbf{1}_{A_n} \rangle.$$

An induction argument completes the proof of (5.112).

(iii) Arguing by contradiction, let us assume that the support of $\zeta$ is not contained in $\mathbb{R}_{\geq 0}$. This means that we can find $\varepsilon > 0$ such that $\zeta(-\infty, -\varepsilon] > 0$. By the result of the previous step applied to the set $A := (-\infty, -\varepsilon]$, the event

$$A_n := \{R_{\ell,\ell'} \leq -\varepsilon \text{ for all } \ell \neq \ell' \leq n\}$$

has positive probability. On this event, we have

$$\Big|\sum_{\ell=1}^{n} \sigma_\ell\Big|^2 = \sum_{\ell,\ell'=1}^{n} R_{\ell,\ell'} \leq n - n(n-1)\varepsilon,$$

which can be made negative by taking $n$ sufficiently large. Since this is absurd, the proof is complete.

## S.6 Mean-field spin glasses

**Exercise 6.1.** Fix $\beta > 0$. Jensen's inequality and the explicit formula for the moment generating function of a Gaussian random variable imply that

$$\mathbb{E} \max_{\sigma \in \Sigma_N} H_N(\sigma) = \frac{1}{\beta} \log \exp \mathbb{E} \beta \max_{\sigma \in \Sigma_N} H_N(\sigma) \leq \frac{1}{\beta} \log \mathbb{E} \exp \beta \max_{\sigma \in \Sigma_N} H_N(\sigma)$$

$$\leq \frac{1}{\beta} \log \sum_{\sigma \in \Sigma_N} \mathbb{E} \exp \beta H_N(\sigma) \leq \frac{1}{\beta} \log 2^N \exp(\beta^2 N/2)$$

$$= \frac{N \log 2}{\beta} + \frac{\beta N}{2}.$$

Optimizing over $\beta > 0$ to find $\beta^2 = 2 \log 2$ yields

$$\mathbb{E} \max_{\sigma \in \Sigma_N} H_N(\sigma) \leq \sqrt{\log 4} N.$$

This establishes the upper bound. To prove the lower bound, we optimize one coordinate at a time (some would call it a "greedy algorithm"). For each integer $1 \leq i \leq N$ and spins $\sigma_1, \ldots, \sigma_{i-1} \in \Sigma_1$, define the random variable

$$X_i(\sigma_1, \ldots, \sigma_{i-1}) := \sum_{j=1}^{i-1} (g_{i,j} + g_{j,i}) \sigma_j$$

in such a way that

$$H_N(\sigma) = \frac{1}{\sqrt{N}} \sum_{i=1}^{N} g_{i,i} + \frac{1}{\sqrt{N}} \sum_{i=1}^{N} \sigma_i X_i(\sigma_1, \ldots, \sigma_{i-1}).$$

Fix $\sigma_1, \ldots, \sigma_N \in \Sigma_1$ inductively so that, for every $i \in \{1, \ldots, N\}$, the spin $\sigma_i$ has the same sign as $X_i(\sigma_1, \ldots, \sigma_{i-1})$ (with some arbitrary decision rule in the case that $X_i(\sigma_1, \ldots, \sigma_{i-1}) = 0$). With this choice of $\sigma = (\sigma_1, \ldots, \sigma_N) \in \Sigma_N$, the random variables $(X_i(\sigma_1, \ldots, \sigma_{i-1}))_{1 \leq i \leq N}$ are independent, and $X_i(\sigma_1, \ldots, \sigma_{i-1})$ is a Gaussian random variable with variance $2(i-1)$. These random variables are also independent of the family $(g_{i,i})_{1 \leq i \leq N}$. Since

$$H_N(\sigma) = \frac{1}{\sqrt{N}} \sum_{i=1}^{N} g_{i,i} + \frac{1}{\sqrt{N}} \sum_{i=1}^{N} |X_i(\sigma_1, \ldots, \sigma_{i-1})|,$$

and the absolute value of a standard Gaussian has mean $\sqrt{\frac{2}{\pi}}$, we obtain that

$$\mathbb{E} \max_{\sigma \in \Sigma_N} H_N(\sigma) \geq \sqrt{\frac{2}{\pi N}} \sum_{i=1}^{N} \sqrt{2(i-1)} \geq \sqrt{\frac{2}{\pi N}} \int_0^{N-1} \sqrt{2x} \, dx \geq C^{-1} N$$

for some constant $C > 0$. This completes the proof.

**Exercise 6.2.** Combining the covariance structure (6.4) of the SK Hamiltonian with Exercise 4.4 reveals that for any $\varepsilon > 0$,

$$\mathbb{P}\left\{\left|\frac{1}{N}\max_{\sigma \in \Sigma_N} H_N(\sigma) - \frac{1}{N}\mathbb{E}\max_{\sigma \in \Sigma_N} H_N(\sigma)\right| \geq \varepsilon\right\} \leq 2\exp\left(-\frac{N\varepsilon^2}{2}\right).$$

Since the right side of this equation is summable in $N$, the result follows by a simple application of the Borel-Cantelli lemma.

**Exercise 6.3.** Bounding the maximum over all configurations by the sum over all configurations reveals that

$$\frac{1}{N}\mathbb{E}\max_{\sigma \in \Sigma_N} H_N(\sigma) = \frac{1}{N\beta}\mathbb{E}\log\exp\beta \max_{\sigma \in \Sigma_N} H_N(\sigma) \leq \frac{\overline{F}_N(\beta)}{\beta}.$$

On the other hand, bounding the Hamiltonian at each configuration by the maximum of the Hamiltonian over all configurations shows that

$$\frac{\overline{F}_N(\beta)}{\beta} \leq \frac{1}{N\beta}\mathbb{E}\log\sum_{\sigma \in \Sigma_N}\exp\beta \max_{\sigma \in \Sigma_N} H_N(\sigma) \leq \frac{\log(2)}{\beta} + \frac{1}{N}\mathbb{E}\max_{\sigma \in \Sigma_N} H_N(\sigma).$$

This completes the proof.

**Exercise 6.4.** For each $t \in [0,1]$, define the interpolating Hamiltonian

$$H_{N,t}(\sigma) := \frac{1}{\sqrt{N}}\sum_{i,j=1}^{N}\left(\sqrt{t}x_{ij} + \sqrt{1-t}g_{ij}\right)\sigma_i\sigma_j$$

in such a way that $H_{N,0}(\sigma) = H_N(\sigma)$ and $H_{N,1}(\sigma) = H_N^x(\sigma)$. As usual, denote by $\overline{F}_N(t)$ its corresponding interpolating free energy and by $\langle \cdot \rangle_t$ the average with respect to its associated Gibbs measure. Fix $t \in (0,1)$, and observe that

$$\overline{F}_N'(t) = \frac{\beta}{N}\mathbb{E}\langle \partial_t H_{N,t}'(\sigma)\rangle_t$$

$$= \frac{\beta}{2N^{3/2}\sqrt{t}}\sum_{i,j=1}^{N}\mathbb{E}x_{ij}\langle\sigma_i\sigma_j\rangle_t - \frac{\beta}{2N^{3/2}\sqrt{1-t}}\sum_{i,j=1}^{N}\mathbb{E}g_{ij}\langle\sigma_i\sigma_j\rangle_t$$

$$\leq \sum_{i,j=1}^{N}\left|\frac{\beta}{2N^{3/2}\sqrt{t}}\mathbb{E}x_{ij}\langle\sigma_i\sigma_j\rangle_t - \frac{\beta}{2N^{3/2}\sqrt{1-t}}\mathbb{E}g_{ij}\langle\sigma_i\sigma_j\rangle_t\right|. \tag{S.30}$$

We will now bound the terms in this expression using the Gaussian integration by parts formula (Theorem 4.5) and the approximate Gaussian integration by parts formula (Exercise 4.6). Fix indices $1 \leq i, j \leq N$ and consider $F(g) := \langle \sigma_i \sigma_j \rangle_t$ as a function of the Gaussian vector $g := (g_{ij})_{i,j \leq N}$. A direct computation reveals that

$$\partial_{g_{ij}} F = \beta \langle \sigma_i \sigma_j \partial_{g_{ij}} H_{N,t}(\sigma) \rangle_t - \langle \sigma_i \sigma_j \rangle_t \langle \beta \partial_{g_{ij}} H_{N,t}(\sigma) \rangle_t$$

$$= \frac{\beta \sqrt{1-t}}{\sqrt{N}} \langle \sigma_i \sigma_j \sigma_i \sigma_j \rangle_t - \frac{\beta \sqrt{1-t}}{\sqrt{N}} \langle \sigma_i \sigma_j \rangle_t \langle \sigma_i \sigma_j \rangle_t$$

$$= \frac{\beta \sqrt{1-t}}{\sqrt{N}} \left(1 - \langle \sigma_i \sigma_j \rangle_t^2 \right),$$

so the Gaussian integration by parts formula yields

$$\frac{\beta}{2N^{3/2}\sqrt{1-t}} \mathbb{E} g_{ij} \langle \sigma_i \sigma_j \rangle_t = \frac{\beta}{2N^{3/2}\sqrt{1-t}} \mathbb{E} \partial_{g_{ij}} F = \frac{\beta^2}{2N^2} (1 - \mathbb{E} \langle \sigma_i \sigma_j \rangle_t^2). \quad (\text{S.31})$$

If instead we consider $F(x) := \langle \sigma_i \sigma_j \rangle_t$ as a function of the random vector $x := (x_{ij})$, then

$$\partial_{x_{ij}} F = \beta \langle \sigma_i \sigma_j \partial_{x_{ij}} H_{N,t}(\sigma) \rangle_t - \langle \sigma_i \sigma_j \rangle_t \langle \beta \partial_{x_{ij}} H_{N,t}(\sigma) \rangle_t = \frac{\beta \sqrt{t}}{\sqrt{N}} \left(1 - \langle \sigma_i \sigma_j \rangle_t^2 \right)$$

$$\partial_{x_{ij}}^2 F = -2 \frac{\beta \sqrt{t}}{\sqrt{N}} \langle \sigma_i \sigma_j \rangle_t \partial_{x_{ij}} F = \frac{2\beta^2 t}{N} \left(\langle \sigma_i \sigma_j \rangle_t^3 - \langle \sigma_i \sigma_j \rangle_t\right).$$

In particular,

$$\left| \partial_{x_{ij}}^2 F \right| \leq \frac{4\beta^2 t}{N},$$

so the approximate Gaussian integration by parts formula and the assumption that $\mathbb{E} x_{ij}^2 = 1$ imply that

$$\left| \frac{\beta}{2N^{3/2}\sqrt{t}} \mathbb{E} x_{ij} \langle \sigma_i \sigma_j \rangle_t - \frac{\beta^2}{2N^2} (1 - \mathbb{E}\langle \sigma_i \sigma_j \rangle_t^2) \right| \leq \frac{\beta}{2N^{3/2}\sqrt{t}} \frac{6\beta^2 t}{N} \mathbb{E}|x_{11}|^3 = \frac{3\beta^3 \sqrt{t}}{N^{5/2}} \mathbb{E}|x_{11}|^3.$$

Combining this with (S.30), (S.31), and the fundamental theorem of calculus reveals that

$$\left| \overline{F}_N^x(\beta) - \overline{F}_N(\beta) \right| = \left| \overline{F}_N(1) - \overline{F}_N(0) \right| \leq \sup_{t \in [0,1]} \left| \overline{F}_N'(t) \right| \leq \frac{3\beta^3}{\sqrt{N}} \mathbb{E}|x_{11}|^3.$$

Letting $N$ tend to infinity completes the proof.

**Exercise 6.5.** By Chebyshev's inequality, it suffices to show that

$$\lim_{N \to +\infty} \mathbb{E}\left( \left\langle F\left(\frac{g \cdot \sigma}{\sqrt{N}}\right) \right\rangle - \mathbb{E}_Z F\left(\frac{g \cdot \langle \sigma \rangle}{\sqrt{N}} + \sqrt{1-q} Z\right) \right)^2 = 0. \quad (\text{S.32})$$

This will follow from the fact that

$$\lim_{N \to +\infty} \mathbb{E}\left\langle F\left(\frac{g \cdot \sigma}{\sqrt{N}}\right)\right\rangle^2 = \lim_{N \to +\infty} \mathbb{E}\left(\left\langle F\left(\frac{g \cdot \sigma}{\sqrt{N}}\right)\right\rangle \mathbb{E}_Z F\left(\frac{g \cdot \langle \sigma \rangle}{\sqrt{N}} + \sqrt{1-q}Z\right)\right)$$

$$= \lim_{N \to +\infty} \mathbb{E}\left(\mathbb{E}_Z F\left(\frac{g \cdot \langle \sigma \rangle}{\sqrt{N}} + \sqrt{1-q}Z\right)\right)^2$$

$$= \mathbb{E}_{Z_1,Z_2,Z_3} F\left(\sqrt{q}Z_1 + \sqrt{1-q}Z_2\right) F\left(\sqrt{q}Z_1 + \sqrt{1-q}Z_3\right)$$

for some independent standard Gaussian random variables $Z_1, Z_2, Z_3$ also independent of all other sources of randomness. We first prove rigorously that

$$\lim_{N \to +\infty} \mathbb{E}\left\langle F\left(\frac{g \cdot \sigma}{\sqrt{N}}\right)\right\rangle^2 = \mathbb{E}_{Z_1,Z_2,Z_3} F\left(\sqrt{q}Z_1 + \sqrt{1-q}Z_2\right) F\left(\sqrt{q}Z_1 + \sqrt{1-q}Z_3\right), \quad (S.33)$$

and then describe how to adapt the argument for the computation of the two other limits. Since $g$ is independent of all other sources of randomness,

$$\mathbb{E}\left\langle F\left(\frac{g \cdot \sigma}{\sqrt{N}}\right)\right\rangle^2 = \mathbb{E}\left\langle \mathbb{E}_g F\left(\frac{g \cdot \sigma^1}{\sqrt{N}}\right) F\left(\frac{g \cdot \sigma^2}{\sqrt{N}}\right)\right\rangle.$$

To simplify this quantity, observe that conditionally on the replicas $\sigma^1, \sigma^2 \in \Sigma_N$, the average

$$\mathbb{E}_g F\left(\frac{g \cdot \sigma^1}{\sqrt{N}}\right) F\left(\frac{g \cdot \sigma^2}{\sqrt{N}}\right)$$

can be expressed as some continuous function of the covariance matrix

$$\mathsf{C}_N := \begin{bmatrix} 1 & R_{1,2} \\ R_{1,2} & 1 \end{bmatrix}$$

of the random vector $G_N := N^{-1/2}(g \cdot \sigma^1, g \cdot \sigma^2)$. This means that there exists a bounded continuous function $\Gamma : \mathbb{R}^{2 \times 2} \to \mathbb{R}$ with

$$\mathbb{E}\left\langle F\left(\frac{g \cdot \sigma}{\sqrt{N}}\right)\right\rangle^2 = \mathbb{E}\langle \Gamma(\mathsf{C}_N) \rangle.$$

To determine the asymptotic behaviour of this quantity, denote by

$$\mathsf{C} := \begin{bmatrix} 1 & q \\ q & 1 \end{bmatrix}$$

the covariance matrix of $G := (\sqrt{q}Z_1 + \sqrt{1-q}Z_2, \sqrt{q}Z_1 + \sqrt{1-q}Z_3)$, and fix $\varepsilon > 0$ as well as $\delta > 0$. Introduce the set

$$A_{\delta,N} := \{\sigma^1, \sigma^2 \in \Sigma_N \mid |R_{1,2} - q| \leq \delta\},$$

and observe that

$$\mathbb{E}\langle \mathbf{1}_{\{|\Gamma(C_N)-\Gamma(C)|>\varepsilon\}}\rangle \leq \mathbb{E}\langle \mathbf{1}_{\{|\Gamma(C_N)-\Gamma(C)|>\varepsilon\}}\mathbf{1}_{A_{\delta,N}}\rangle + \mathbb{E}\langle \mathbf{1}_{A^c_{\delta,N}}\rangle. \tag{S.34}$$

Since $|C_N - C| \leq \sqrt{2}\delta$ on the set $A_{\delta,N}$, and $\Gamma$ is continuous, by choosing $\delta$ small enough relative to $\varepsilon$, the first term in (S.34) can be made to vanish. It follows by Chebyshev's inequality that

$$\mathbb{E}\langle \mathbf{1}_{\{|\Gamma(C_N)-\Gamma(C)|>\varepsilon\}}\rangle \leq \mathbb{E}\langle \mathbf{1}_{A^c_{\delta,N}}\rangle \leq \frac{\mathbb{E}\langle (R_{1,2}-q)^2\rangle}{\delta^2}.$$

Letting $N$ tend to infinity shows that $(\Gamma(C_N))_{N \geq 1}$ converges in probability to $\Gamma(C)$. Since this sequence of random variables is bounded as $F$ is, the limit (S.33) holds by the dominated convergence theorem. We now briefly describe the main differences in the computation of the limits of

$$\mathbb{E}\left(\left\langle F\left(\frac{g\cdot\sigma}{\sqrt{N}}\right)\right\rangle \mathbb{E}_Z F\left(\frac{g\cdot\langle\sigma\rangle}{\sqrt{N}}+\sqrt{1-q}Z\right)\right) \quad \text{and} \quad \mathbb{E}\left(\mathbb{E}_Z F\left(\frac{g\cdot\langle\sigma\rangle}{\sqrt{N}}+\sqrt{1-q}Z\right)\right)^2.$$

The first of these quantities can be expressed as

$$\mathbb{E}\left\langle \mathbb{E}_g \mathbb{E}_Z F\left(\frac{g\cdot\sigma}{\sqrt{N}}\right) F\left(\frac{g\cdot\langle\sigma\rangle}{\sqrt{N}}+\sqrt{1-q}Z\right)\right\rangle,$$

so its asymptotic behaviour is determined by that of the covariance matrix

$$C_N := \begin{bmatrix} 1 & \frac{\sigma\cdot\langle\sigma\rangle}{N} \\ \frac{\sigma\cdot\langle\sigma\rangle}{N} & \langle R_{1,2}\rangle + 1 - q \end{bmatrix}$$

which converges to $C$ as the asymptotic overlap distribution concentrates on the singleton $q$. Similarly, the second of these quantities can be expressed as

$$\mathbb{E}\mathbb{E}_g \mathbb{E}_{Z_1,Z_2} F\left(\frac{g\cdot\langle\sigma\rangle}{\sqrt{N}}+\sqrt{1-q}Z_1\right) F\left(\frac{g\cdot\langle\sigma\rangle}{\sqrt{N}}+\sqrt{1-q}Z_2\right),$$

so its asymptotic behaviour is determined by that of the covariance matrix

$$C_N := \begin{bmatrix} \langle R_{1,2}\rangle + 1 - q & \langle R_{1,2}\rangle \\ \langle R_{1,2}\rangle & \langle R_{1,2}\rangle + 1 - q \end{bmatrix}$$

which again converges to $C$. This completes the proof.

**Exercise 6.6.** For every $g \in L^1([0,1];\mathbb{R})$, we use the notation

$$\fint_{\zeta_k}^{\zeta_{k+1}} g := \frac{1}{\zeta_{k+1}-\zeta_k}\int_{\zeta_k}^{\zeta_{k+1}} g.$$

We first observe that, by Jensen's inequality, we have for every $g \in L^2([0,1];\mathbb{R})$ that

$$0 \leq \sum_{k=0}^{K} \int_{\zeta_k}^{\zeta_{k+1}} \left(g - \fint_{\zeta_k}^{\zeta_{k+1}} g\right)^2 = \int_0^1 g^2 - \sum_{k=0}^{K}(\zeta_{k+1}-\zeta_k)\left(\fint_{\zeta_k}^{\zeta_{k+1}} g\right)^2 \leq \int_0^1 g^2. \quad \text{(S.35)}$$

We also recall that the set of continuous functions is dense in $L^2([0,1];\mathbb{R})$, see for instance Theorem 2.19 of [6]. Let $\phi \in C([0,1];\mathbb{R})$ be a continuous function such that $\|f - \phi\|_{L^2} \leq \varepsilon$. Appealing to (S.35) with $g = f - \phi$, we see that we make an error of at most $\varepsilon$ if we substitute $f$ by $\phi$ in the expression on the left side of (6.83). Since $\phi$ is uniformly continuous, there exists $\delta > 0$ such that for every $x, y \in [0,1]$ with $|x - y| \leq \delta$, we have $|\phi(y) - \phi(x)| \leq \varepsilon$. If $\max_k(\zeta_{k+1}-\zeta_k) \leq \delta$, then

$$\sum_{k=0}^{K} \int_{\zeta_k}^{\zeta_{k+1}} \left(\phi - \fint_{\zeta_k}^{\zeta_{k+1}} \phi\right)^2 \leq \varepsilon^2.$$

Using also the identity in (S.35), we therefore obtain the desired result, up to a reparametrization of $\varepsilon > 0$.

**Exercise 6.7.** A direct computation reveals that

$$-\partial_t^2 f = \mathbb{E}\langle \partial_t^2 H(t,\sigma)\rangle + \mathbb{E}\langle(\partial_t H(t,\sigma))^2\rangle - \mathbb{E}\langle \partial_t H(t,\sigma)\rangle^2.$$

We now use Gaussian integration by parts to simplify each of these terms. On the one hand, the Gibbs Gaussian integration by parts formula (Theorem 4.6) implies that

$$\mathbb{E}\langle \partial_t^2 H(t,\sigma)\rangle = -\frac{1}{(2t)^{3/2}}\mathbb{E}\langle z\sigma\rangle = -\frac{1}{2t}\bigl(\mathbb{E}\langle\sigma^2\rangle - \mathbb{E}\langle\sigma\rangle^2\bigr).$$

On the other hand,

$$\mathbb{E}\langle(\partial_t H(t,\sigma))^2\rangle - \mathbb{E}\langle\partial_t H(t,\sigma)\rangle^2$$
$$= \frac{1}{2t}\mathbb{E}\langle z^2\sigma(\sigma-\sigma')\rangle + \frac{2}{(2t)^{1/2}}\bigl(\mathbb{E}\langle z\sigma''\sigma^2\rangle - \mathbb{E}\langle z\sigma^3\rangle\bigr) + \mathbb{E}\langle\sigma^4\rangle - \mathbb{E}\langle\sigma^2\rangle^2,$$

and the Gaussian integration by parts formula (Theorem 4.5) reveals that

$$\mathbb{E}\langle z^2\sigma(\sigma-\sigma')\rangle = \mathbb{E}z\cdot z\frac{\int \sigma(\sigma-\sigma')\exp\bigl(H(t,\sigma)+H(t,\sigma')\bigr)\,dP(\sigma)\,dP(\sigma')}{\int \exp\bigl(H(t,\sigma)+H(t,\sigma')\bigr)\,dP(\sigma)\,dP(\sigma')}$$
$$= \mathbb{E}\langle\sigma(\sigma-\sigma')\rangle + \sqrt{2t}\,\mathbb{E}\langle z\sigma(\sigma-\sigma')(\sigma+\sigma')\rangle$$
$$\qquad - 2\sqrt{2t}\,\mathbb{E}\langle z\sigma(\sigma-\sigma')\sigma''\rangle$$
$$= \mathbb{E}\langle\sigma(\sigma-\sigma')\rangle + 2t\mathbb{E}\langle(\sigma-\sigma')(\sigma+\sigma')(\sigma^2+\sigma\sigma'-2\sigma\sigma'')\rangle$$
$$\qquad - 2t\mathbb{E}\langle(\sigma-\sigma')(\sigma''+\sigma''')(\sigma^2+\sigma\sigma'+\sigma\sigma''+\sigma\sigma'''-4\sigma\sigma'''')\rangle$$
$$= \mathbb{E}\langle\sigma^2\rangle - \mathbb{E}\langle\sigma\rangle^2 + 2h\bigl(\mathbb{E}\langle\sigma^4\rangle - 4\mathbb{E}\langle\sigma^3\rangle\langle\sigma\rangle - 3\mathbb{E}\langle\sigma^2\rangle^2$$
$$\qquad + 12\mathbb{E}\langle\sigma^2\rangle\langle\sigma\rangle^2 - 6\mathbb{E}\langle\sigma\rangle^4\bigr).$$

Similarly, the Gibbs Gaussian integration by parts formula implies that

$$\mathbb{E}\langle z\sigma''\sigma^2\rangle - \mathbb{E}\langle z\sigma^3\rangle = (2t)^{1/2}\big(\mathbb{E}\langle \sigma^2(\sigma''^2 + \sigma''\sigma - 2\sigma'\sigma'')\rangle - \mathbb{E}\langle \sigma^2(\sigma^2 - \sigma\sigma')\rangle\big)$$
$$= (2t)^{1/2}\big(\mathbb{E}\langle \sigma^2\rangle^2 + 2\mathbb{E}\langle \sigma^3\rangle\langle\sigma\rangle - 2\mathbb{E}\langle \sigma^2\rangle\langle\sigma\rangle^2 - \mathbb{E}\langle \sigma^4\rangle\big).$$

It follows that

$$\partial_t^2 f = 6\mathbb{E}\langle\sigma\rangle^4 - \mathbb{E}\langle\sigma^4\rangle + 4\mathbb{E}\langle\sigma^3\rangle\langle\sigma\rangle + 3\mathbb{E}\langle\sigma^2\rangle^2 - 12\mathbb{E}\langle\sigma^2\rangle\langle\sigma\rangle^2$$
$$- 2\big(\mathbb{E}\langle\sigma^2\rangle^2 + 2\mathbb{E}\langle\sigma^3\rangle\langle\sigma\rangle - 2\mathbb{E}\langle\sigma^2\rangle\langle\sigma\rangle^2 - \mathbb{E}\langle\sigma^4\rangle\big) + \mathbb{E}\langle\sigma^2\rangle^2 - \mathbb{E}\langle\sigma^4\rangle$$
$$= 6\mathbb{E}\langle\sigma\rangle^4 + 2\mathbb{E}\langle\sigma^2\rangle^2 - 8\mathbb{E}\langle\sigma^2\rangle\langle\sigma\rangle^2$$
$$= 2\mathbb{E}\big(\langle\sigma^2\rangle - \langle\sigma\rangle^2\big)\big(\langle\sigma^2\rangle - 3\langle\sigma\rangle^2\big).$$

as required. We now consider the probability measure $P_1 := p\delta_1 + (1-p)\delta_{-1}$ for a choice of $p$ to be determined. It will be important to observe that

$$\langle\sigma\rangle = \frac{p\exp(\sqrt{2t}z) - (1-p)\exp(-\sqrt{2t}z)}{p\exp(\sqrt{2t}z) + (1-p)\exp(-\sqrt{2t}z)}.$$

The Gibbs Gaussian integration by parts formula implies that

$$\partial_t f(t) = \mathbb{E}\langle\sigma^2\rangle - \frac{1}{(2t)^{1/2}}\mathbb{E}\langle z\sigma\rangle = \mathbb{E}\langle\sigma\rangle^2$$

from which it is clear that $\partial_t f(0) = (2p-1)^2 < 1$. The dominated convergence theorem reveals that $\lim_{t\to+\infty} \partial_t f(t) = 1$, so there exists $t^* \geq 0$ with $\partial_t^2 f(t^*) > 0$. However, by the formula we obtained earlier, and observing that the Gibbs measure $\langle\cdot\rangle$ is simply the measure $P_1$ when $t = 0$, we have

$$\partial_t^2 f(0) = 2\langle(\sigma - \langle\sigma\rangle)^2\rangle(1 - 3\langle\sigma\rangle^2).$$

This quantity is strictly negative whenever $p \in (0,1)$ is such that $\langle\sigma\rangle^2 > 1/3$. Since $\langle\sigma\rangle = 2p - 1$, an explicit calculation yields that this occurs if and only if

$$p \in \left(0, \frac{3-\sqrt{3}}{6}\right) \cup \left(\frac{3+\sqrt{3}}{6}, 1\right) \supseteq \left(0, \frac{1}{5}\right) \cup \left(\frac{4}{5}, 1\right).$$

This shows that for this choice of measure $P_1$, the free energy $f$ is neither concave nor convex.

## S.7 Basic results in analysis and probability

**Exercise A.1.** Fix measurable sets $A, B \in \mathcal{S}$ with $A \subseteq B$, and observe that by additivity of measure,

$$\mu(B) = \mu(A) + \mu(B \setminus A) \geq \mu(A).$$

This establishes monotonicity. To prove continuity from below, fix an increasing sequence of measurable sets $(A_n)_{n \geq 1} \subseteq \mathcal{S}$, and observe that by countable additivity of measure,

$$\mu\left(\bigcup_{n=1}^{\infty} A_n\right) = \sum_{n=1}^{+\infty} \mu(A_n \setminus A_{n-1}) = \lim_{n \to +\infty} \sum_{i=1}^{n} \mu(A_i \setminus A_{i-1}) = \lim_{n \to +\infty} \mu(A_n).$$

This completes the proof.

**Exercise A.2.** To see that the function $m_F$ is well-defined, observe that $m_F$ is finitely additive. Indeed, if $((a_i, b_i])_{1 \leq i \leq n}$ is a finite sequence of disjoint intervals and

$$(a, b] := \bigcup_{i=1}^{n} (a_i, b_i],$$

then, up to relabelling, we may assume that $a = a_1 < b_1 \leq a_2 < \ldots < b_{n-1} \leq a_n < b_n = b$, and therefore

$$m_F\left(\bigcup_{i=1}^{n}(a_i, b_i]\right) = F(b) - F(a) = \sum_{i=1}^{n} (F(b_i) - F(a_i)) = \sum_{i=1}^{n} m_F((a_i, b_i]).$$

It follows that for any two finite sequences $((a_i, b_i])_{1 \leq i \leq n}$ and $((c_j, d_j])_{1 \leq j \leq m}$ of disjoint intervals with the same union,

$$m_F\left(\bigcup_{i=1}^{n}(a_i, b_i]\right) = \sum_{i=1}^{n} \sum_{j=1}^{m} m_F((a_i, b_i] \cap (c_j, d_j]) = m_F\left(\bigcup_{j=1}^{m}(c_j, d_j]\right).$$

This shows that $m_F$ is well-defined. To prove that it is a pre-measure, fix a sequence $(I_n)_{n \geq 1}$ of disjoint intervals with $I_n = (a_n, b_n]$ whose union $I$ belongs to $\mathcal{A}$. Since the union of $(I_n)_{n \geq 1}$ is a finite union of intervals of the form $(a, b]$, the sequence $(I_n)_{n \geq 1}$ can be partitioned into finitely many subsequences such that the union of the intervals in each subsequence is a single interval of the form $(a, b]$. Considering each subsequence separately and using the finite additivity of $m_F$, assume without loss of generality that $I = (a, b]$. On the one hand, by finite additivity,

$$m_F(I) = m_F\left(\bigcup_{i=1}^{n} I_i\right) + m_F\left(I \setminus \bigcup_{i=1}^{n} I_i\right) \geq m_F\left(\bigcup_{i=1}^{n} I_i\right) = \sum_{i=1}^{n} m_F(I_i),$$

so letting $n$ tend to infinity gives the lower bound,

$$m_F(I) \geqslant \sum_{n=1}^{+\infty} m_F(I_n). \tag{S.36}$$

To obtain the matching upper bound, we will only consider the case when $a$ and $b$ are both finite; the general case can be deduced through a simple limiting argument. Fix $\varepsilon > 0$ as well as $n \geqslant 0$, and use the continuity of $F$ to find $\delta > 0$ and $\delta_n > 0$ with $F(a+\delta) - F(a) \leqslant \varepsilon$ and $F(b_n + \delta_n) - F(b_n) \leqslant \varepsilon 2^{-n}$. The open intervals $((a_n, b_n + \delta_n))_{n \geqslant 1}$ cover the compact set $[a+\delta, b]$, so it is possible to extract a finite subcover. Up to relabelling, denote by $((a_i, b_i + \delta_i))_{1 \leqslant i \leqslant n}$ this finite subcover and assume that $b_i + \delta_i \in (a_{n+1}, b_{n+1} + \delta_{n+1})$ so that the subcover is ordered according to the right endpoint of each interval. The choice of $\delta$ and the non-decreasingness of $F$ imply that

$$m_F(I) \leqslant F(b) - F(a+\delta) + \varepsilon \leqslant F(b_n + \delta_n) - F(a_n) + \sum_{i=1}^{n-1}\bigl(F(a_{i+1}) - F(a_i)\bigr) + \varepsilon.$$

The choice $b_i + \delta_i \in (a_{n+1}, b_{n+1} + \delta_{n+1})$, the non-decreasingness of $F$ and the choice of $\delta_i$ reveal that this is bounded further by

$$F(b_n + \delta_n) - F(a_n) + \sum_{i=1}^{n-1}\bigl(F(b_i + \delta_i) - F(a_i)\bigr) + \varepsilon \leqslant \sum_{i=1}^{n}\bigl(F(b_i) - F(a_i) + \varepsilon 2^{-i}\bigr) + \varepsilon.$$

It follows by definition of $m_F$ that

$$m_F(I) \leqslant \sum_{i=1}^{n}\bigl(m_F(I_i) + \varepsilon 2^{-i}\bigr) + \varepsilon \leqslant \sum_{n=1}^{+\infty} m_F(I_n) + 2\varepsilon.$$

Letting $\varepsilon$ tend to zero and combining the resulting bound with (S.36) completes the proof.

**Exercise A.3.** Since $\mu$ and $\nu$ agree on $\mathcal{P}$ and $A \in \mathcal{P}$, it is clear that $S \in \mathcal{L}$. To verify the second property of a $\lambda$-system, fix $B, C \in \mathcal{L}$ with $B \subseteq C$. By additivity of measure,

$$\mu(A \cap C \smallsetminus B) = \mu(A \cap C) - \mu(A \cap B) = \nu(A \cap C) - \nu(A \cap B) = \nu(A \cap C \smallsetminus B),$$

where the assumption that $\mu(A) = \nu(A) < +\infty$ ensures that the subtraction is well-defined. Finally, if $(B_n)_{n \geqslant 1}$ is an increasing sequence of sets in $\mathcal{L}$ and $B := \bigcup_{n=1}^{\infty} B_n$, then the continuity of measure established in Exercise A.1 implies that

$$\mu(A \cap B) = \lim_{n \to +\infty} \mu(A \cap B_n) = \lim_{n \to +\infty} \nu(A \cap B_n) = \nu(A \cap B)$$

which means that $B \in \mathcal{L}$. This shows that $\mathcal{L}$ is a $\lambda$-system. The final claim follows from the Dynkin $\pi$-$\lambda$ theorem by taking $A = S$. This completes the proof.

**Exercise A.4.** We treat each question separately.

(i) Since $S$ and $\varnothing$ are both closed sets, we have $S \in \mathcal{L}$. Suppose $A, B \in \mathcal{L}$ with $A \subseteq B$. Fix $\varepsilon > 0$ and let $F_A \subseteq A$ and $F_{A^c} \subseteq A^c$ be closed sets with

$$\mu(F_A) + \frac{\varepsilon}{2} \geqslant \mu(A) \quad \text{and} \quad \mu(F_{A^c}) + \frac{\varepsilon}{2} \geqslant \mu(A^c).$$

Similarly, let $F_B \subseteq B$ and $F_{B^c} \subseteq B^c$ be closed sets with

$$\mu(F_B) + \frac{\varepsilon}{2} \geqslant \mu(B), \quad \mu(F_{B^c}) + \frac{\varepsilon}{2} \geqslant \mu(B^c).$$

By additivity of measure

$$\mu(B \smallsetminus A) = \mu(B) + \mu(A^c) - \mu(S) \leqslant \mu(F_B) + \mu(F_{A^c}) - \mu(S) + \varepsilon = \mu(F_B \smallsetminus F_{A^c}^c) + \varepsilon$$

and

$$\mu\big((B \smallsetminus A)^c\big) = \mu(A) + \mu(B^c) \leqslant \mu(F_A) + \mu(F_{B^c}) + \varepsilon = \mu(F_A \cup F_{B^c}) + \varepsilon.$$

Since $F_B \smallsetminus F_{A^c}^c \subseteq B \smallsetminus A$ and $F_A \cup F_{B^c} \subseteq (B \smallsetminus A)^c$ are closed sets, we have $B \smallsetminus A \in \mathcal{L}$. Finally, suppose $(A_n)_{n \geqslant 1} \subseteq \mathcal{L}$ increases to $A$ and fix $\varepsilon > 0$. For each $n \geqslant 1$, let $F_n \subseteq A_n^c$ be a closed set with

$$\mu(A_n^c \smallsetminus F_n) = \mu(A_n^c) - \mu(F_n) \leqslant \frac{\varepsilon}{2^n}.$$

Using the continuity of measure established in Exercise A.1, fix $N$ large enough so $\mathbb{P}(A \smallsetminus A_N) \leqslant \varepsilon/2$, and let $F \subseteq A_N$ be a closed set with $\mu(A_N) \leqslant \mu(F) + \varepsilon/2$. By monotonicity and subadditivity of measure,

$$\mu(A) = \mu(A \smallsetminus A_N) + \mu(A_N) \leqslant \mu(F) + \varepsilon$$

and

$$\mu\Big(A^c \smallsetminus \bigcap_{n \geqslant 1} F_n\Big) = \mu\Big(\bigcap_{n \geqslant 1} A_n^c \smallsetminus \bigcap_{n \geqslant 1} F_n\Big) \leqslant \mu\Big(\bigcup_{n \geqslant 1} A_n^c \smallsetminus F_n\Big) \leqslant \sum_{n \geqslant 1} \mu(A_n^c \smallsetminus F_n) \leqslant \varepsilon.$$

Since $F \subseteq A$ and $\bigcap_{n \geqslant 1} F_n \subseteq A^c$ are closed, we have $A \in \mathcal{L}$. This shows that $\mathcal{L}$ is a $\lambda$-system.

(ii) Fix an open set $U \subseteq S$, and define the increasing sequence $(F_n)_{n \geqslant 1}$ of closed sets by

$$F_n := \{s \in S \mid d(s, U^c) \geqslant 1/n\}.$$

Since $U$ is open, we have $U = \bigcup_{n \geqslant 1} F_n$. It follows by the continuity of measure established in Exercise A.1 that

$$\mu(U) = \lim_{n \to +\infty} \mu(F_n) = \sup_{n \geqslant 1} \mu(F_n) \leqslant \sup_{\substack{F \subseteq U \\ F \text{ closed}}} \mu(F) \leqslant \mu(U),$$

and therefore $U \in \mathcal{L}$.

(iii) Since $\mathcal{L}$ is a $\lambda$-system and the set of open sets forms a $\pi$-system that generates the Borel $\sigma$-algebra, the result follows by the Dynkin $\pi$-$\lambda$ theorem.

**Exercise A.5.** Fix a positive linear functional $T: C(S;\mathbb{R}) \to \mathbb{R}$ as well as a continuous function $f \in C(S;\mathbb{R})$. Since $\|f\|_\infty \cdot 1 \pm f \geq 0$, the positivity of $T$ implies that

$$T(\|f\|_\infty \cdot 1 + f) \geq 0 \quad \text{and} \quad T(\|f\|_\infty \cdot 1 - f) \geq 0$$

It follows by linearity of $T$ that

$$-\|f\|_\infty T(1) \leq T(f) \leq \|f\|_\infty T(1).$$

Rearranging reveals that $|T(f)| \leq \|f\|_\infty T(1)$ and completes the proof.

**Exercise A.6.** Suppose that we have established the result on the space of continuous functions $C([0,1];\mathbb{R})$, and fix $f \in C([a,b];\mathbb{R})$ as well as $\varepsilon > 0$. Define the function $g: [0,1] \to \mathbb{R}$ by

$$g(x) := f(a + x(b-a)),$$

and find a polynomial $Q$ on $[0,1]$ with $\|g - Q\|_\infty \leq \varepsilon$. Consider the polynomial $P$ on $[0,1]$ defined by

$$P(x) := Q\left(\frac{x-a}{b-a}\right),$$

and observe that for any $x \in [a,b]$, we have

$$|f(x) - P(x)| = \left|Q\left(\frac{x-a}{b-a}\right) - g\left(\frac{x-a}{b-a}\right)\right| \leq \varepsilon.$$

It therefore suffices to prove the result on $C([0,1];\mathbb{R})$. Fix $\varepsilon > 0$ and $f \in C([0,1];\mathbb{R})$. Given $x \in [0,1]$, let $(X_n)_{n \geq 1}$ be a sequence of independent and identically distributed random variables with $X_1 \sim \text{Ber}(x)$. For each $n \geq 1$, consider the sample average $S_n := \frac{1}{n}\sum_{i=1}^n X_i$, and define the polynomial

$$P_n(x) := \mathbb{E}f(S_n) = \sum_{0 \leq i \leq n} f\left(\frac{i}{n}\right)\binom{n}{i}x^i(1-x)^{n-i}.$$

To express $P_n$ as a polynomial we have used the fact that $\sum_{i=1}^n X_i$ is a Binomial random variable with $n$ trials and probability of success $x$. To show that for $n$ large enough $P_n$ is within $\varepsilon$ of $f$, use the uniform continuity of $f$ on $[0,1]$ to find $\delta > 0$ so that $|f(y) - f(x)| \leq \frac{\varepsilon}{2}$ whenever $|x - y| \leq \delta$, and observe that by Chebyshev's inequality and independence of the $(X_n)_{n \geq 1}$, for any $x \in [0,1]$,

$$|P_n(x) - f(x)| \leq \mathbb{E}|f(S_n) - f(x)|\mathbf{1}_{\{|S_n - x| \leq \delta\}} + 2\|f\|_\infty \mathbb{P}\{|S_n - x| > \delta\}$$
$$\leq \frac{\varepsilon}{2} + \frac{2\|f\|_\infty}{\delta^2}\text{Var}(S_n)$$
$$\leq \frac{\varepsilon}{2} + \frac{2\|f\|_\infty}{n\delta^2}$$

Choosing $n$ large enough ensures that $\|f - P\|_\infty \leq \varepsilon$ and completes the proof.

**Exercise A.7.** Denote by $\langle \cdot, \cdot \rangle : L^2(\mathcal{C}; \mathbb{C}) \times L^2(\mathcal{C}; \mathbb{C}) \to \mathbb{C}$ the canonical inner product on $L^2(\mathcal{C}; \mathbb{C})$,
$$\langle f, g \rangle := \int_0^{2\pi} f(e^{it}) \overline{g}(e^{it}) \, dt.$$
Fix a polynomial $P: \mathcal{C} \to \mathbb{C}$ of the form $P(z) := \sum_{j=0}^n a_j z^j$ for some complex coefficients $(a_j)_{j \leq n} \subseteq \mathbb{C}$. Observe that the functions $f$ and $P$ are orthogonal,
$$\int_0^{2\pi} \overline{f}(e^{it}) P(e^{it}) \, dt = \sum_{j=0}^n a_j \int_0^{2\pi} e^{i(j+1)t} \, dt = 0.$$
It follows that
$$2\pi = \langle f, f \rangle = \langle f, f - P \rangle \leq 2\pi \|f\|_\infty \|f - P\|_\infty \leq 2\pi \|f - P\|_\infty,$$
so $\|f - P\|_\infty \geq 1$ for any polynomial $P$. This completes the proof.

**Exercise A.8.** Fix an open set $U \subseteq S$ and let $F = U^c$. Since $d(x, F) = 0$ if and only if $x \in F$, the sequence $(f_m)_{m \geq 1} \subseteq C_b(S)$ defined by $f_m(x) := \min(1, md(x, F))$ increases pointwise to $\mathbf{1}_U$ as $m$ tends to infinity. It follows by the monotone convergence theorem that $\mathbb{P}(U) = \mathbb{Q}(U)$. Invoking the Dynkin $\pi$-$\lambda$ theorem completes the proof.

**Exercise A.9.** We treat each question separately.

(i) To begin with, suppose that $(\mathbb{P}_n)_{n \geq 1}$ converges weakly to $\mathbb{P}$, and let $t \in \mathbb{R}$ be a point of continuity of $F$. Fix $\varepsilon > 0$ and let $\phi_1$ and $\phi_2$ be piecewise linear functions with
$$\mathbf{1}_{\{x \leq t - \varepsilon\}} \leq \phi_1(x) \leq \mathbf{1}_{\{x \leq t\}} \leq \phi_2(x) \leq \mathbf{1}_{\{x \leq t + \varepsilon\}}.$$
Integrating the outer inequalities with respect to $\mathbb{P}$ and the inner inequalities with respect to $\mathbb{P}_n$ reveals that
$$F(t - \varepsilon) \leq \int_\mathbb{R} \phi_1 \, d\mathbb{P} \leq \liminf_{n \to +\infty} F_n(t) \leq \limsup_{n \to +\infty} F_n(t) \leq \int_\mathbb{R} \phi_2 \, d\mathbb{P} \leq F(t + \varepsilon).$$
We have used the fact that $\phi_1, \phi_2 \in C_b(\mathbb{R}; \mathbb{R})$ to let $n$ tend to infinity. Letting $\varepsilon$ tend to zero and remembering that $t$ is a point of continuity of $F$ shows that $(F_n(t))_{n \geq 1}$ converges to $F(t)$. Conversely, suppose that for any point $t$ of continuity of $F$, the sequence $(F_n(t))_{n \geq 1}$ converges to $F(t)$. Fix $f \in C_b(\mathbb{R}; \mathbb{R})$ as well as $\varepsilon > 0$ and $k \geq 1$. Denote by $\mathcal{C}_F$ the set of points of continuity of $F$. The set $\mathcal{C}_F$ is countable by monotonicity of $F$, so there exists $M \in \mathcal{C}_F$ with $\mathbb{P}((-M, M]^c) \leq \varepsilon$. Let $(x_i^k)_{0 \leq i \leq k} \subseteq \mathcal{C}_F$ be a partition of $[-M, M]$ with $\max_{1 \leq i \leq k-1} |x_i^k - x_{i-1}^k| \leq k^{-1}$. Using the fact that $(F_n)_{n \geq 1}$ converges to $F$ at $M$ and $-M$, find $n$ large enough so $\mathbb{P}_n((-M, M]^c) \leq 2\varepsilon$. We introduce the approximating function
$$f_k(x) := \sum_{1 \leq i \leq k} f(x_i^k) \mathbf{1}_{(x_{i-1}^k, x_i^k]}(x),$$

and we combine the uniform continuity of $f$ on $[-M,M]$ with the fact that $\max_{1\leqslant i\leqslant k-1}|x^k_{i+1} - x^k_i| \leqslant k^{-1}$ to find $k$ large enough so $|f_k(x) - f(x)| \leqslant \varepsilon$ whenever $|x| \leqslant M$. Since $(x^k_i)_{0\leqslant i\leqslant k} \subseteq C_F$, we have

$$\lim_{n\to+\infty} \int_{\mathbb{R}} f_k \, d\mathbb{P}_n = \sum_{1\leqslant i\leqslant k} f_k(x^k_i)(F(x^k_i) - F(x^k_{i-1})) = \int_{\mathbb{R}} f_k \, d\mathbb{P}.$$

It follows that for $n$ large enough,

$$\left| \int_{\mathbb{R}} f \, d\mathbb{P} - \int_{\mathbb{R}} f \, d\mathbb{P}_n \right| \leqslant 3\varepsilon \|f\|_\infty + 3\varepsilon.$$

Letting $n$ tend to infinity and then $\varepsilon$ tend to zero completes the proof.

(ii) The assumption is necessary. Consider the sequence $(\mathbb{P}_n)_{n\geqslant 1}$ of probability measures with $\mathbb{P}_n(\{n^{-1}\}) = 1$. This sequence converges weakly to the probability measure $\mathbb{P}$ defined by $\mathbb{P}(\{0\}) = 1$. However, the sequence of distribution functions of $\mathbb{P}_n$ fails to converge to the distribution function of $\mathbb{P}$ at the point of discontinuity $t = 0$.

**Exercise A.10.** We treat each question separately.

(i) For every $\varepsilon > 0$ and $\delta > 0$, introduce the set

$$U_{\varepsilon,\delta} := \{x \in S \mid \text{there exist } y, z \in S \text{ such that } d(x,y) < \delta, \, d(x,z) < \delta,$$
$$\text{and } d'(f(y), f(z)) \geqslant \varepsilon\},$$

where $d'$ denotes the metric on $S'$. The space $U_{\varepsilon,\delta}$ is open. Moreover, the complement of $\mathscr{C}_f$ is given by

$$\bigcup_{\varepsilon>0} \bigcap_{\delta>0} U_{\varepsilon,\delta},$$

where the ranges for $\varepsilon$ and $\delta$ can be restricted to cover only rational values. This shows that $\mathscr{C}_f$ is measurable.

(ii) Let $F$ be a closed subset of $S'$. We have that

$$\mathbb{P}\{f(X_n) \in F\} = \mathbb{P}\{X_n \in f^{-1}(F)\} \leqslant \mathbb{P}\{X_n \in \overline{f^{-1}(F)}\}.$$

By the Portmanteau theorem, we deduce that

$$\limsup_{n\to+\infty} \mathbb{P}\{f(X_n) \in F\} \leqslant \mathbb{P}\{X_\infty \in \overline{f^{-1}(F)}\}.$$

Since $X_\infty$ takes values in $\mathscr{C}_f$ with probability one, the last quantity can be rewritten as

$$\mathbb{P}\{X_\infty \in \mathscr{C}_f \cap \overline{f^{-1}(F)}\}.$$

To show that $\mathscr{C}_f \cap \overline{f^{-1}(F)}$ is a subset of $F$, fix $x \in \mathscr{C}_f \cap \overline{f^{-1}(F)}$, and let $(x_n)_{n \geq 1}$ be a sequence of elements of $f^{-1}(F)$ that converges to $x$. Since $f(x_n) \in F$, the sequence $(x_n)_{n \geq 1}$ converges to $x \in \mathscr{C}_f$. Combined with the fact that $F$ is a closed set, this implies that $f(x) \in F$, and thus $x \in f^{-1}(F)$. It follows that

$$\limsup_{n \to +\infty} \mathbb{P}\{f(X_n) \in F\} \leq \mathbb{P}\{X_\infty \in \mathscr{C}_f \cap \overline{f^{-1}(F)}\} \leq \mathbb{P}\{f(X_\infty) \in F\}.$$

Invoking the Portmanteau theorem again completes the proof.

**Exercise A.11.** By the Portmanteau theorem, it suffices to show that for any closed set $F \subseteq S$, we have $\limsup_{n \to +\infty} \mathbb{P}\{Y_n \in F\} \leq \mathbb{P}\{Y \in F\}$. Fix a closed set $F \subseteq S$ as well as $\varepsilon > 0$, and let $F_\varepsilon := \{x \in S \mid d(x, F) \leq \varepsilon\}$ be the $\varepsilon$-neighbourhood of $F$. Observe that

$$\mathbb{P}\{Y_n \in F\} \leq \mathbb{P}\{Y_{k,n} \in F_\varepsilon\} + \mathbb{P}\{d(Y_{k,n}, Y_n) \leq \varepsilon\}.$$

Since $(Y_{k,n})_{n \geq 1}$ converges in law to $Z_k$ and $F_\varepsilon$ is closed, the Portmanteau theorem implies that

$$\limsup_{n \to +\infty} \mathbb{P}\{Y_n \in F\} \leq \mathbb{P}\{Z_k \in F_\varepsilon\} + \limsup_{n \to +\infty} \mathbb{P}\{d(Y_{k,n}, Y_n) \leq \varepsilon\}.$$

Another application of the Portmanteau theorem to the convergence in law of $(Z_k)_{k \geq 1}$ to $Y$ reveals that $\limsup_{n \to +\infty} \mathbb{P}\{Y_n \in F\} \leq \mathbb{P}\{Y \in F_\varepsilon\}$. Leveraging the fact that $F$ is closed shows that $F_\varepsilon \searrow F$, so letting $\varepsilon$ tend to zero completes the proof.

**Exercise A.12.** This is immediate from Theorem A.21 upon realizing that the sequence $(\mathbb{P})_{n \geq 1}$ converges weakly to $\mathbb{P}$. It is also possible to give a direct proof. Since $S$ is separable there exists a countable set $(x_i)_{i \geq 1}$ such that for every $m \geq 1$, we have $S = \bigcup_{i \geq 1} \overline{B}_{1/m}(x_i)$. Given $\varepsilon > 0$, the continuity of measure established in Exercise A.1 implies the existence of $n(m) \geq 1$ with

$$\mathbb{P}\left(S \setminus \bigcup_{i=1}^{n(m)} \overline{B}_{1/m}(x_i)\right) \leq \frac{\varepsilon}{2^m}.$$

It follows that the set $K := \bigcap_{m=1}^{\infty} \bigcup_{i=1}^{n(m)} \overline{B}_{1/m}(x_i)$ is such that

$$\mathbb{P}(S \setminus K) = \mathbb{P}\left(\bigcup_{m=1}^{\infty} S \setminus \bigcup_{i=1}^{n(m)} \overline{B}_{1/m}(x_i)\right) \leq \sum_{m=1}^{+\infty} \frac{\varepsilon}{2^m} = \varepsilon.$$

Since $S$ is complete and $K$ is closed and totally bounded, $K$ is compact. This completes the proof.